\documentclass[12pt]{amsbook}
\usepackage{amssymb}
\usepackage[all]{xy}

\usepackage[colorlinks=true,linkcolor=magenta,citecolor=magenta]{hyperref}
\makeindex

\newtheorem{theorem}{Theorem}[chapter]
\newtheorem{definition}[theorem]{Definition}
\newtheorem{proposition}[theorem]{Proposition}
\newtheorem{conjecture}[theorem]{Conjecture}
\newtheorem{exercise}[theorem]{Exercise}

\begin{document}

\title{Invitation to Hadamard matrices}

\author{Teo Banica}
\address{Department of Mathematics, University of Cergy-Pontoise, F-95000 Cergy-Pontoise, France. {\tt teo.banica@gmail.com}}


\begin{abstract}
An Hadamard matrix is a square matrix $H\in M_N(\pm1)$ whose rows and pairwise orthogonal. More generally, we can talk about the complex Hadamard matrices, which are the square matrices $H\in M_N(\mathbb C)$ whose entries are on the unit circle, $|H_{ij}|=1$, and whose rows and pairwise orthogonal. The main examples are the Fourier matrices, $F_N=(w^{ij})$ with $w=e^{2\pi i/N}$, and at the level of the general theory, the complex Hadamard matrices can be thought of as being some sort of exotic, generalized Fourier matrices. We discuss here the basic theory of the Hadamard  matrices, real and complex, with emphasis on the complex matrices, and their geometric and analytic aspects.
\end{abstract}

\maketitle

\chapter*{Preface}

Linear algebra is full of mysteries, with sometimes even single matrices hiding interesting mathematics, worth a lengthy contemplation. Well-known examples include the Pauli spin matrices, which are cult objects in physics, at the core of basic quantum mechanics, then the Dirac matrices, at the core of quantum electrodynamics (QED), and the Gell-Mann matrices, at the core of quantum chromodynamics (QCD).

\bigskip

This book is about a class of matrices which are particularly beautiful, no matter your aesthetics, and whose study is fun and pleasant, bringing us into lots of interesting mathematics, coming from algebra, geometry, analysis and probability. And which are of course potentially useful for something. These are the Hadamard matrices.

\bigskip

A complex Hadamard matrix is a square matrix $H\in M_N(\mathbb C)$ whose entries are on the unit circle in the complex plane, $|H_{ij}|=1$, and whose rows are pairwise orthogonal, with respect to the usual scalar product on $\mathbb C^N$. The central example is the Fourier matrix, $F_N=(w^{ij})$ with $w=e^{2\pi iN}$, with the name coming from the fact that this is the matrix of the Fourier transform over the cyclic group $G=\mathbb Z_N$. In general, a complex Hadamard matrix can be thought of as being a kind of ``generalized Fourier matrix", and the applications of the complex Hadamard matrices come from this.

\bigskip

There has been a lot of work on the Hadamard matrices, starting with Sylvester and Hadamard, long time ago, who looked at such matrices in the real case, $H\in M_N(\mathbb R)$. Here the Hadamard matrix condition states that we must have $H\in M_N(\pm1)$, and that when comparing any two rows, the number of matchings must equal the number of mismatchings. The whole subject belongs to combinatorics, design theory and group theory, although there are some interesting analytic and probabilistic aspects as well.

\bigskip

Later on, it was realized that the general complex case, $H\in M_N(\mathbb C)$, is worth attention too, with motivation coming from discrete Fourier analysis, in a large sense. The subject here belongs to linear algebra, real algebraic geometry, combinatorics again, with plenty of constructions involving all sorts of tricky roots of unity, and with interesting analytic and probabilistic aspects as well. As for the potential applications, these belong to quantum physics, via constuctions involving operator algebras and quantum groups.

\bigskip

All in all, many things to be explained, and this book is an introduction to all this, with the aim of keeping things simple, but reasonably complete. 

\bigskip

The first half of the book, Parts I and II, deals with the real Hadamard matrices, whose basic theory is quite elementary, and then with the basic theory in the complex case, using elementary algebraic and geometric techniques. Everything here is accessible with a minimal knowledge of linear algebra, and calculus in several variables.

\bigskip

The second half of the book, Parts III and IV, contains more advanced material, erring on the graduate side. We will discuss here advanced analytic techniques for dealing with the complex Hadamard matrices, and then we will have a look into potential applications to theoretical physics, at the level of quantum groups and operator algebras.

\bigskip

Although many things will be discussed in this book, this remains an introduction to the subject. There has been a huge amount of work in the real case, and we will discuss here only the very basic ideas behind this work. The same goes for the construction and classification work in the complex case, with once again a lot of literature waiting to be consulted, by the interested reader. As in what regards the applications, both in the real and the complex case, our discussion here will be something modest too, with the main aim being that of explaining the relation between the quantum groups and the Hadamard matrices, which is where the applications to quantum physics should come from.

\bigskip

There are several books dedicated to the Hadamard matrices, including Agaian \cite{aga}, Horadam \cite{hor} and Seberry-Yamada \cite{sya}, all focusing on the real case, and by using algebraic methods. It is our hope that the present book can stand as a nice complement to these, written from a physicist's viewpoint, and as an invitation to the subject. 

\bigskip

This book is partly based on a number of research papers that I wrote, and I am grateful to Julien Bichon, Ion Nechita and Jean-Marc Schlenker, for our joint work on the subject. Many thanks go as well to my cats, for advice with hunting techniques, martial arts, and more. When doing linear algebra, all this knowledge is very useful.

\bigskip

\

{\em Cergy, July 2024}

\smallskip

{\em Teo Banica}

\baselineskip=15.95pt
\tableofcontents
\baselineskip=14pt

\part{Hadamard matrices}

\ \vskip50mm

\begin{center}
{\em And only say that you'll be mine

In no others' arms entwine

Down beside where the waters flow

Down by the banks of the Ohio}
\end{center}

\chapter{Hadamard matrices}

\section*{1a. Hadamard matrices}

We will be mainly interested in this book in the complex Hadamard matrices, but let us start with some beautiful pure mathematics, regarding the real case. The definition that we need, going back to 19th century work of Sylvester \cite{syl}, on topics such as tessellated pavements and ornamental tile-work, is as follows:

\index{Sylvester}
\index{Hadamard matrix}
\index{real Hadamard matrix}

\begin{definition}
An Hadamard matrix is a square binary matrix, 
$$H\in M_N(\pm1)$$
whose rows are pairwise orthogonal, with respect to the scalar product on $\mathbb R^N$. 
\end{definition}

There are many examples of such matrices, and we will discuss this, in what follows. To start with, here is an example, which is a particularly beautiful one:
$$K_4=\begin{pmatrix}
-1&1&1&1\\
1&-1&1&1\\
1&1&-1&1\\
1&1&1&-1
\end{pmatrix}$$

Observe that this matrix has many interesting extra features, such as being symmetric, bistochastic, and circulant. Here is another example, also at $N=4$, which is interesting too, because it reminds the combinatorics of the Klein group $\mathbb Z_2\times\mathbb Z_2$:
$$W_4=\begin{pmatrix}
1&1&1&1\\
1&-1&1&-1\\
1&1&-1&-1\\
1&-1&-1&1
\end{pmatrix}$$

Summarizing, we have examples of Hadamard matrices, usually coming from certain interesting algebraic and combinatorial properties of $\mathbb R^N$, which are waiting to be explored. In general now, as a first theoretical observation, we do not really need real numbers in order to talk about the Hadamard matrices, because we have:

\index{mismatchings}

\begin{proposition}
A binary matrix $H\in M_N(\pm1)$ is Hadamard when its rows have the property that, when comparing any two of them,
$$\begin{matrix}
e_1&\ldots&e_N\\
f_1&\ldots&f_N
\end{matrix}$$
the number of matchings $(e_i=f_i)$ equals the number of mismatchings $(e_i\neq f_i)$.
\end{proposition}

\begin{proof}
This is clear from definitions. Indeed, the scalar product on $\mathbb R^N$ is given by:
$$<x,y>=\sum_ix_iy_i$$

Thus, when computing the scalar product between two rows, the matchings contribute with $1$ factors, and the mismatchings with $-1$ factors, and this gives the result. 
\end{proof}

As a consequence of the above result, we can replace if we want the $1,-1$ entries of our matrix by any two symbols, of our choice. Here is an example of an Hadamard matrix, and to be more precise, the above matrix $W_4$, written with this convention:
$$\begin{matrix}
\heartsuit&\heartsuit&\heartsuit&\heartsuit\\
\heartsuit&\clubsuit&\heartsuit&\clubsuit\\
\heartsuit&\heartsuit&\clubsuit&\clubsuit\\
\heartsuit&\clubsuit&\clubsuit&\heartsuit
\end{matrix}$$

However, it is probably better to run away from this, and use real numbers instead, as in Definition 1.1, with the idea in mind of connecting the Hadamard matrices to the foundations of modern mathematics, namely Calculus 1 and Calculus 2. So, getting back now to the real numbers, here is our first result:

\index{orthogonal group}

\begin{proposition}
For a square matrix $H\in M_N(\pm1)$, the following are equivalent:
\begin{enumerate}
\item The rows of $H$ are pairwise orthogonal, and so $H$ is Hadamard.

\item The columns of $H$ are pairwise orthogonal, and so $H^t$ is Hadamard.

\item The rescaled matrix $U=H/\sqrt{N}$ is orthogonal, $U\in O_N$.
\end{enumerate}
\end{proposition}

\begin{proof}
The idea here is that the equivalence between (1) and (2) is not exactly obvious, but both these conditions can be shown to be equivalent to (3), as follows:

\medskip

$(1)\iff(3)$ Since the rows of $U=H/\sqrt{N}$ have norm 1, this matrix is orthogonal precisely when its rows are pairwise orthogonal. But this latter condition is equivalent to the fact that the rows of $H=\sqrt{N}U$ are pairwise orthogonal, as desired.

\medskip

$(2)\iff(3)$ The same argument as above shows that $H^t$ is Hadamard precisely when its rescaling $U^t=H^t/\sqrt{N}$ is orthogonal. But since a matrix $U\in M_N(\mathbb R)$ is orthogonal precisely when its transpose $U^t\in M_N(\mathbb R)$ is orthogonal, this gives the result.
\end{proof}

As an abstract consequence of the above result, let us record:

\begin{theorem}
The set of the $N\times N$ Hadamard matrices is
$$Y_N=M_N(\pm 1)\cap\sqrt{N}O_N$$
where $O_N$ is the orthogonal group, the intersection being taken inside $M_N(\mathbb R)$.
\end{theorem}

\begin{proof}
This follows from the equivalence $(1)\iff(3)$ in Proposition 1.3, which tells us that an arbitrary $H\in M_N(\pm1)$ belongs to $Y_N$ if and only if it belongs to $\sqrt{N}O_N$.
\end{proof}

As a conclusion to what we have so far, the set $Y_N$ that we are interested in appears as a kind of set of ``special rational points'' of the real algebraic manifold $\sqrt{N}O_N$. Thus, we are doing some kind of algebraic geometry here, of precise type to be determined. In the simplest case, $N=2$, the Hadamard matrices are elementary to compute, and the set $Y_2$ consists precisely of the rational points of $\sqrt{2}O_2$, the result being as follows:

\begin{theorem}
The binary matrices $H\in M_2(\pm1)$ are split $50$-$50$ between Hadamard and non-Hadamard, the Hadamard ones being as follows,
$$\begin{pmatrix}1&1\\1&-1\end{pmatrix}\ \ \qquad
\begin{pmatrix}1&1\\-1&1\end{pmatrix}\ \ \qquad
\begin{pmatrix}1&-1\\1&1\end{pmatrix}\ \ \qquad
\begin{pmatrix}-1&1\\1&1\end{pmatrix}$$
$$\begin{pmatrix}1&-1\\-1&-1\end{pmatrix}\qquad
\begin{pmatrix}-1&1\\-1&-1\end{pmatrix}\qquad
\begin{pmatrix}-1&-1\\1&-1\end{pmatrix}\qquad
\begin{pmatrix}-1&-1\\-1&1\end{pmatrix}$$
and the non-Hadamard ones being the remaining ones. Also, we have $Y_2=M_2(\mathbb Q)\cap\sqrt{2}O_2$, with the intersection being taken inside $M_N(\mathbb R)$.
\end{theorem}

\begin{proof}
We have two assertions to be proved, which are both elementary:

\medskip

(1) In what regards the classification, this is best done by using the Hadamard matrix criterion from Proposition 1.2, which at $N=2$ simply tells us that, once the first row is chosen, the choices for the second row, as for our matrix to be Hadamard, are exactly $50\%$. The solutions are those in the statement, listed according to the lexicographic order, with respect to the standard way of reading, left to right, and top to bottom.

\medskip

\index{rotation matrix}

(2) In order to prove the second assertion, we use the fact that $O_2$ consists of 2 types of matrices, namely rotations $R_t$ and symmetries $S_t$. To be more precise, we first have the rotation of angle $t\in\mathbb R$, which is given by the following formula:
$$R_t=\begin{pmatrix}\cos t&-\sin t\\ \sin t&\cos t\end{pmatrix}$$

We also have the symmetry with respect to the $Ox$ axis rotated by $t/2\in\mathbb R$:
$$S_t=\begin{pmatrix}\cos t&\sin t\\ \sin t&-\cos t\end{pmatrix}$$

Now by multiplying everything by $\sqrt{2}$, we are led to the following formula:
$$\sqrt{2}O_2=\left\{\begin{pmatrix}c&-s\\s&c\end{pmatrix}\,,\ \begin{pmatrix}c&s\\s&-c\end{pmatrix}\Big|c^2+s^2=2\right\}$$

In order to find now the matrices from $\sqrt{2}O_2$ having rational entries, we must solve the following equation, over the integers:
$$x^2+y^2=2z^2$$

But this is equivalent to $y^2-z^2=z^2-x^2$, which is impossible for obvious reasons, unless we have $x^2=y^2=z^2$. Thus, the rational points come from $c^2=s^2=1$, and so we have a total of $2\times2\times 2=8$ rational points, which can only be the points of $Y_2$.
\end{proof}

At higher values of $N$, we cannot expect $Y_N$ to consist of the rational points of $\sqrt{N}O_N$. As a basic counterexample, we have the following matrix, which is not Hadamard:
$$\begin{pmatrix}
2&0&0&0\\
0&2&0&0\\
0&0&2&0\\
0&0&0&2
\end{pmatrix}\in 2O_4$$

Summarizing, it is quite unclear what $Y_N$ is, geometrically speaking. We can, however, solve this question by using complex numbers, in the following way:

\index{Hadamard matrix}
\index{complex Hadamard matrix}
\index{unitary group}

\begin{theorem}
The Hadamard matrices appear as the real points,
$$Y_N=M_N(\mathbb R)\cap X_N$$
of the complex Hadamard matrix manifold, which is given by:
$$X_N=M_N(\mathbb T)\cap\sqrt{N}U_N$$
Thus, $Y_N$ is the real part of an intersection of smooth real algebraic manifolds.
\end{theorem}

\begin{proof}
This is a version of Theorem 1.4, which can be established in two ways:

\medskip

(1) We can either define a complex Hadamard matrix to be a matrix $H\in M_N(\mathbb T)$, with $\mathbb T$ standing as usual for the unit circle in the complex plane, whose rows are pairwise orthogonal, with respect to the scalar product of $\mathbb C^N$, then work out a straightforward complex analogue of Proposition 1.3, which gives the formula of $X_N$ in the statement, and then observe that the real points of $X_N$ are the real Hadamard matrices.

\medskip

(2) Or, we can directly use Theorem 1.4, which formally gives the result, as follows:
\begin{eqnarray*}
Y_N
&=&M_N(\pm1)\cap\sqrt{N}O_N\\
&=&\big[M_N(\mathbb R)\cap M_N(\mathbb T)\big]\cap\big[M_N(\mathbb R)\cap\sqrt{N}U_N\big]\\
&=&M_N(\mathbb R)\cap\big[M_N(\mathbb T)\cap\sqrt{N}U_N\big]\\
&=&M_N(\mathbb R)\cap X_N
\end{eqnarray*}

We will be back to this, and more precisely with full details regarding (1), starting from chapter 5 below, when studying the complex Hadamard matrices.
\end{proof}

Summarizing, the Hadamard matrices do belong to real algebraic geometry, but in a quite subtle way. We will be back to all this, gradually, in what follows. 

\section*{1b. Walsh matrices}

Let us discuss now the examples of Hadamard matrices, with a systematic study at $N=4,6,8,10$ and so on, continuing the study from Theorem 1.5. In order to cut a bit from complexity, we can use the following notion:

\index{equivalent matrices}
\index{Hadamard equivalence}

\begin{definition}
Two Hadamard matrices are called equivalent, and we write $H\sim K$, when it is possible to pass from $H$ to $K$ via the following operations:
\begin{enumerate}
\item Permuting the rows, or the columns.

\item Multiplying the rows or columns by $-1$.
\end{enumerate}
\end{definition}

Observe that we do not include the transposition operation $H\to H^t$ in our list of allowed operations. This is because Proposition 1.3, while looking quite elementary, rests however on a deep linear algebra fact, namely that the transpose of an orthogonal matrix is orthogonal as well, and this can produce complications later on.

\bigskip

As another comment, there is of course a certain group $G$ acting there, made of two copies of $S_N$, one for the rows and one for the columns, and of two copies of $\mathbb Z_2^N$, once again one for the rows, and one for the columns. The equivalence classes of the Hadamard matrices are then the orbits of the action $G\curvearrowright Y_N$. It is possible to be a bit more explicit here, with a formula for $G$ and so on, but we will not need this.

\bigskip

Given an Hadamard matrix $H\in M_N(\pm1)$, we can use the above two operations in order to put $H$ in a ``nice'' form. Although there is no clear definition for what ``nice'' should mean, for the Hadamard matrices, with this being actually a quite subtle problem, that we will discuss later on, here is something that we can look for:

\index{dephased matrix}

\begin{definition}
An Hadamard matrix is called dephased when it is of the form
$$H=\begin{pmatrix}
1&\ldots&1\\
\vdots&*\\
1
\end{pmatrix}$$
that is, when the first row and the first column consist of $1$ entries only.
\end{definition}

Here the terminology comes from the complex Hadamard matrices, introduced in Theorem 1.6 and its proof. Indeed, when regarding $H\in M_N(\pm1)$ as a complex matrix, $H\in M_N(\mathbb T)$, the $-1$ entries have ``phases'', equal to $\pi$, and assuming that $H$ is dephased means to assume that we have no phases, on the first row and the first column.

\bigskip

Observe that, up to the equivalence relation, any Hadamard matrix $H\in M_N(\pm1)$ can be put in dephased form. Moreover, the dephasing operation is unique, if we use only the operations (2) in Definition 1.7, namely row and column multiplications by $-1$. The point now is that, with these notions in hand, we can formulate a nice classification result:

\begin{theorem}
There is only one Hadamard matrix at $N=2$, namely
$$W_2=\begin{pmatrix}1&1\\ 1&-1\end{pmatrix}$$
up to the above equivalence relation for such matrices.
\end{theorem}

\begin{proof}
The matrix in the statement $W_2$, called Walsh matrix, is clearly Hadamard. Conversely, given $H\in M_N(\pm1)$ Hadamard, we can dephase it, as follows:
$$\begin{pmatrix}a&b\\c&d\end{pmatrix}
\to\begin{pmatrix}1&1\\ac&bd\end{pmatrix}
\to\begin{pmatrix}1&1\\1&abcd\end{pmatrix}$$

Now since the dephasing operation preserves the class of the Hadamard matrices, we must have $abcd=-1$, and so we obtain by dephasing the matrix $W_2$.
\end{proof}

At $N=3$ we cannot have examples, due to the orthogonality condition between the rows, which forces $N$ to be even, for obvious reasons. At $N=4$ now, we have several examples. In order to discuss them, let us start with:

\index{tensor product}
\index{double indices}
\index{Walsh matrix}

\begin{proposition}
If $H\in M_M(\pm1)$ and $K\in M_N(\pm1)$ are Hadamard matrices, then so is their tensor product, constructed in double index notation as follows:
$$H\otimes K\in M_{MN}(\pm1)\quad,\quad 
(H\otimes K)_{ia,jb}=H_{ij}K_{ab}$$
In particular the Walsh matrices, $W_N=W_2^{\otimes n}$ with $N=2^n$, are all Hadamard.
\end{proposition}

\begin{proof}
The matrix in the statement $H\otimes K$ has indeed $\pm1$ entries, and its rows $R_{ia}$ are pairwise orthogonal, as shown by the following computation:
\begin{eqnarray*}
<R_{ia},R_{kc}>
&=&\sum_{jb}H_{ij}K_{ab}\cdot H_{kj}K_{cb}\\
&=&\sum_jH_{ij}H_{kj}\sum_bK_{ab}K_{cb}\\
&=&M\delta_{ik}\cdot N\delta_{ac}\\
&=&MN\delta_{ia,kc}
\end{eqnarray*}

As for the second assertion, this follows from this, $W_2$ being Hadamard.
\end{proof}

Before going further, we should clarify a bit our tensor product notations. In order to write $H\in M_N(\pm1)$ the indices of $H$ must belong to $\{1,\ldots,N\}$, or at least to an ordered set $\{\alpha_1,\ldots,\alpha_N\}$. But with double indices we are indeed in this latter situation, because we can use the lexicographic order on these indices. To be more precise, by using the lexicographic order on the double indices, we have the following result:

\index{lexicographic order}
\index{tensor product}

\begin{proposition}
Given $H\in M_M(\pm1)$ and $K\in M_N(\pm1)$, we have
$$H\otimes K=\begin{pmatrix}
H_{11}K&\ldots&H_{1M}K\\ 
\vdots&&\vdots\\ 
H_{M1}K&\ldots&H_{MM}K
\end{pmatrix}$$
with respect to the lexicographic order on the double indices.
\end{proposition}

\begin{proof}
We recall that the tensor product is given by $(H\otimes K)_{ia,jb}=H_{ij}K_{ab}$. Now by using the lexicographic order on the double indices, we obtain:
\begin{eqnarray*}
H\otimes K
&=&\begin{pmatrix}
(H\otimes K)_{11,11}&(H\otimes K)_{11,12}&\ldots&(H\otimes K)_{11,MN}\\
(H\otimes K)_{12,11}&(H\otimes K)_{12,12}&\ldots&(H\otimes K)_{12,MN}\\
\vdots&\vdots&&\vdots\\ 
\vdots&\vdots&&\vdots\\ 
(H\otimes K)_{MN,11}&(H\otimes K)_{MN,12}&\ldots&(H\otimes K)_{MN,MN}
\end{pmatrix}\\
&&\\
&=&\begin{pmatrix}
H_{11}K_{11}&H_{11}K_{12}&\ldots&H_{1M}K_{MN}\\
H_{11}K_{21}&H_{11}K_{22}&\ldots&H_{1M}K_{2N}\\
\vdots&\vdots&&\vdots\\ 
\vdots&\vdots&&\vdots\\ 
H_{M1}K_{N1}&H_{M1}K_{N2}&\ldots&H_{MM}K_{NN}
\end{pmatrix}
\end{eqnarray*}

Thus, by making blocks, we are led to the formula in the statement.
\end{proof}

As a basic example for the tensor product construction, the matrix $W_4$, obtained by tensoring the matrix $W_2$ with itself, is given by:
$$W_4
=\begin{pmatrix}W_2&W_2\\ W_2&-W_2\end{pmatrix}
=\begin{pmatrix}1&1&1&1\\ 1&-1&1&-1\\ 1&1&-1&-1\\ 1&-1&-1&1\end{pmatrix}$$

Getting back now to our classification work, here is the result at $N=4$:

\begin{theorem}
There is only one Hadamard matrix at $N=4$, namely
$$W_4=W_2\otimes W_2$$
up to the standard equivalence relation for such matrices.
\end{theorem}

\begin{proof}
Consider an Hadamard matrix $H\in M_4(\pm1)$, assumed to be dephased:
$$H=\begin{pmatrix}1&1&1&1\\ 1&a&b&c\\ 1&d&e&f\\ 1&g&h&i\end{pmatrix}$$

By orthogonality of the first 2 rows, we must have $\{a,b,c\}=\{-1,-1,1\}$. Thus by permuting the last 3 columns, we can assume that our matrix is as follows:
$$H=\begin{pmatrix}1&1&1&1\\ 1&-1&1&-1\\ 1&m&n&o\\ 1&p&q&r\end{pmatrix}$$

Now by orthogonality of the first 2 columns, we must have $\{m,p\}=\{-1,1\}$. Thus by permuting the last 2 rows, we can further assume that our matrix is as follows:
$$H=\begin{pmatrix}1&1&1&1\\ 1&-1&1&-1\\ 1&1&x&y\\ 1&-1&z&t\end{pmatrix}$$

But this gives the result, because the orthogonality of the rows gives $x=y=-1$. Indeed, with these values of $x,y$ plugged in, our matrix becomes:
$$H=\begin{pmatrix}1&1&1&1\\ 1&-1&1&-1\\ 1&1&-1&-1\\ 1&-1&z&t\end{pmatrix}$$

Now from the orthogonality of the columns we obtain:
$$z=-1\ ,\ t=1$$

Thus, up to equivalence of Hadamard matrices we have $H=W_4$, as claimed.
\end{proof}

The case $N=5$ is excluded, because the orthogonality condition between the rows forces $N\in 2\mathbb N$. The point now is that $N=6$ is excluded as well, because we have:

\index{size of Hadamard matrix}

\begin{theorem}
The size of an Hadamard matrix $H\in M_N(\pm1)$ must satisfy 
$$N\in\{2\}\cup 4\mathbb N$$
with this coming from the orthogonality condition between the first $3$ rows.
\end{theorem}

\begin{proof}
By permuting the rows and columns or by multiplying them by $-1$, as to rearrange the first 3 rows, we can always assume that our matrix looks as follows:
$$H=\begin{pmatrix}
1\ldots\ldots 1&1\ldots\ldots 1&1\ldots\ldots 1&1\ldots\ldots 1\\
1\ldots\ldots 1&1\ldots\ldots 1&-1\ldots -1&-1\ldots -1\\
1\ldots\ldots 1&-1\ldots -1&1\ldots\ldots 1&-1\ldots -1\\
\underbrace{\ldots\ldots\ldots}_x&\underbrace{\ldots\ldots\ldots}_y&\underbrace{\ldots\ldots\ldots}_z&\underbrace{\ldots\ldots\ldots}_t
\end{pmatrix}$$

Now if we denote by $x,y,z,t$ the sizes of the block columns, as indicated, the orthogonality conditions between the first 3 rows give the following system of equations:
$$(1\perp 2)\quad:\quad x+y=z+t$$
$$(1\perp 3)\quad:\quad x+z=y+t$$
$$(2\perp 3)\quad:\quad x+t=y+z$$

The numbers $x,y,z,t$ being such that the average of any two equals the average of the other two, and so equals the global average, the solution of our system is:
$$x=y=z=t$$

We therefore conclude that the size of our Hadamard matrix, which is the number $N=x+y+z+t$, must be a multiple of 4, as claimed.
\end{proof}

The above result is something very interesting, and we should mention that a similar analysis with 4 rows or more does not give any further restriction on the possible values of the size $N\in\mathbb N$. In fact, the celebrated Hadamard Conjecture (HC), that we will discuss in a moment, states that there should be an Hadamard matrix at any $N\in 4\mathbb N$.

\index{HC}
\index{Hadamard Conjecture}

\bigskip

Now back to our small $N$ study, the case $N=6$ being excluded by Theorem 1.13, we have to discuss the case $N=8$. Here we have as basic example the Walsh matrix $W_8$, and we will prove that, up to equivalence, this is the only Hadamard matrix at $N=8$. In order to prove this, we will use the $3\times N$ matrix analysis from the proof of Theorem 1.13. To be more precise, we will first improve this into a $4\times N$ matrix result, and then, by assuming $N=8$, we will discuss the case where we have 5 rows or more. Let us start by giving a name to the rectangular matrices that we are interested in:

\index{PHM}
\index{partial Hadamard matrix}

\begin{definition}
A partial Hadamard matrix (PHM) is a rectangular matrix 
$$H\in M_{M\times N}(\pm1)$$
whose rows are pairwise orthogonal, with respect to the scalar product of $\mathbb R^N$. 
\end{definition}

We refer to Hall \cite{hal}, Ito \cite{ito} and Verheiden \cite{ver} for a number of results regarding the PHM. In what follows we will just develop some basic theory, useful in connection with our $N=8$ questions, but we will be back to the PHM, later. We first have:

\index{equivalent PHM}
\index{standard form}
\index{dephased PHM}

\begin{definition}
Two PHM are called equivalent when we can pass from one to the other by permuting rows or columns, or multiplying the rows or columns by $-1$. Also:
\begin{enumerate}
\item We say that a PHM is in dephased form when its first row and its first column consist of $1$ entries.

\item We say that a PHM is in standard form when it is dephased, with the $1$ entries moved to the left as much as possible, by proceeding from top to bottom. 
\end{enumerate}
\end{definition}

With these notions in hand, let us go back now to the proof of Theorem 1.13. The study there concerns the $3\times N$ case, and we can improve this, as follows:

\begin{proposition}
The standard form of the dephased PHM at $M=2,3,4$ is as follows, with $\pm$ standing respectively for various horizontal vectors filled with $\pm1$,
$$H=\begin{pmatrix}+&+\\\underbrace{+}_{N/2}&\underbrace{-}_{N/2}\end{pmatrix}$$
$$H=\begin{pmatrix}+&+&+&+\\+&+&-&-\\\underbrace{+}_{N/4}&\underbrace{-}_{N/4}&\underbrace{+}_{N/4}&\underbrace{-}_{N/4}\end{pmatrix}$$
$$H=\begin{pmatrix}
+&+&+&+&+&+&+&+\\
+&+&+&+&-&-&-&-\\
+&+&-&-&+&+&-&-\\
\underbrace{+}_a&\underbrace{-}_b&\underbrace{+}_b&\underbrace{-}_a&\underbrace{+}_b&\underbrace{-}_a&\underbrace{+}_a&\underbrace{-}_b
\end{pmatrix}$$
and with $a,b\in\mathbb N$ being subject to the condition $a+b=N/4$.
\end{proposition}

\begin{proof}
Here the $2\times N$ assertion is clear, and the $3\times N$ assertion is something that we already know. Let us pick now an arbitrary partial Hadamard matrix $H\in M_{4\times N}(\pm1)$, assumed to be in standard form, as in Definition 1.15 (2). According to the $3\times N$ result, applied to the upper $3\times N$ part of our matrix, our matrix must look as follows:
$$H=\begin{pmatrix}
+&+&+&+&+&+&+&+\\
+&+&+&+&-&-&-&-\\
+&+&-&-&+&+&-&-\\
\underbrace{+}_x&\underbrace{-}_{x'}&\underbrace{+}_{y'}&\underbrace{-}_y&\underbrace{+}_{z'}&\underbrace{-}_z&\underbrace{+}_t&\underbrace{-}_{t'}
\end{pmatrix}$$

To be more precise, our matrix must be indeed of the above form, with $x,y,z,t$ and $x',y',z',t'$ being certain integers, subject to the following relations:
$$x+x'=y+y'=z+z'=t+t'=\frac{N}{4}$$

In terms of these parameters, the missing orthogonality conditions are:
$$(1\perp 4)\quad:\quad x+y'+z'+t=x'+y+z+t'$$
$$(2\perp 4)\quad:\quad x+y'+z+t'=x'+y+z'+t$$
$$(3\perp 4)\quad:\quad x+y+z'+t'=x'+y'+z+t$$

Now observe that these orthogonality conditions can be written as follows:
$$(x-x')-(y-y')-(z-z')+(t-t')=0$$
$$(x-x')-(y-y')+(z-z')-(t-t')=0$$
$$(x-x')+(y-y')-(z-z')-(t-t')=0$$

But this latter system can be solved by using the basic averaging argument from the proof of Theorem 1.13, the solution being as follows:
$$x-x'=y-y'=z-z'=t-t'$$

Now by putting everything together, the conditions to be satisfied by the block lengths are as follows, with $a,b\in\mathbb N$ being subject to the condition $a+b=N/4$:
$$x=y=z=t=a$$
$$x'=y'=z'=t'=b$$

Thus, we are led to the conclusion in the statement.
\end{proof}

In the case $N=8$, that we are interested in here, in view of our classification program from the square matrix case, we have the following more precise result:

\begin{proposition}
There are exactly two $4\times 8$ partial Hadamard matrices, namely
$$I=(W_4\ W_4)$$ 
$$J=(W_4\ K_4)$$
us to the standard equivalence relation for such matrices. 
\end{proposition}

\begin{proof}
We use the last assertion in Proposition 1.16, regarding the $4\times N$ partial Hadamard matrices, at $N=8$. In the case $a=2,b=0$, the solution is:
$$P=\begin{pmatrix}
+&+&+&+&&+&+&+&+\\
+&+&+&+&&-&-&-&-\\
+&+&-&-&&+&+&-&-\\
+&+&-&-&&-&-&+&+
\end{pmatrix}$$

In the case $a=1,b=1$, the solution is:
$$Q=\begin{pmatrix}
+&+&+&+&&+&+&+&+\\
+&+&+&+&&-&-&-&-\\
+&+&-&-&&+&+&-&-\\
+&-&+&-&&+&-&+&-
\end{pmatrix}$$

Finally, in the case $a=0,b=2$, the solution is:
$$R=\begin{pmatrix}
+&+&+&+&&+&+&+&+\\
+&+&+&+&&-&-&-&-\\
+&+&-&-&&+&+&-&-\\
-&-&+&+&&+&+&-&-
\end{pmatrix}$$

Now observe that, by permuting the columns of $P$, we can obtain the following matrix, which is precisely the matrix $I=(W_4\ W_4)$ from the statement:
$$I=\begin{pmatrix}
+&+&+&+&&+&+&+&+\\
+&-&+&-&&+&-&+&-\\
+&+&-&-&&+&+&-&-\\
+&-&-&+&&+&-&-&+
\end{pmatrix}$$

Also, by permuting the columns of $Q$, we can obtain the following matrix, which is equivalent to the matrix $J=(W_4\ K_4)$ from the statement:
$$J'=\begin{pmatrix}
+&+&+&+&&+&+&+&+\\
+&-&+&-&&-&-&+&+\\
+&+&-&-&&-&+&-&+\\
+&-&-&+&&-&+&+&-
\end{pmatrix}$$

Finally, regarding the last solution, $R$, by switching the sign on the last row we obtain $R\sim P$, and so we have $R\sim P\sim I$, which finishes the proof.
\end{proof}

We can now go back to the classification problems for the usual, square Hadamard matrices at $N=8$, and we have here the following result:

\index{Walsh matrix}

\begin{theorem}
The third Walsh matrix, namely
$$W_8=\begin{pmatrix}W_4&W_4\\W_4&-W_4\end{pmatrix}$$
is the unique Hadamard matrix at $N=8$, up to equivalence.
\end{theorem}

\begin{proof}
We use Proposition 1.17, which splits the discussion into two cases:

\medskip

\underline{Case 1}. We must look here for completions of the following matrix $I$:
$$I=\begin{pmatrix}
1&1&1&1&&1&1&1&1\\ 
1&-1&1&-1&&1&-1&1&-1\\ 
1&1&-1&-1&&1&1&-1&-1\\ 
1&-1&-1&1&&1&-1&-1&1
\end{pmatrix}$$

This is something quite technical, which can be basically done in 3 steps, as follows:

\medskip

(1) Let us first try to complete this partial $4\times 8$ Hadamard matrix into a partial $5\times 8$ Hadamard matrix. The completion must look as follows:
$$I'=\begin{pmatrix}
1&1&1&1&&1&1&1&1\\ 
1&-1&1&-1&&1&-1&1&-1\\ 
1&1&-1&-1&&1&1&-1&-1\\ 
1&-1&-1&1&&1&-1&-1&1\\
a&b&c&d&&a'&b'&c'&d'
\end{pmatrix}$$

The system of equations for the orthogonality conditions is as follows:
$$(1\perp 5)\quad:\quad a+b+c+d+a'+b'+c'+d'=0$$
$$(2\perp 5)\quad:\quad a-b+c-d+a'-b'+c'-d'=0$$
$$(3\perp 5)\quad:\quad a+b-c-d+a'+b'-c'-d'=0$$
$$(4\perp 5)\quad:\quad a-b-c+d+a'-b'-c'+d'=0$$

Now observe that this system of equations can be written as follows:
$$(a+a')+(b+b')+(c+c')+(d+d')=0$$
$$(a+a')-(b+b')+(c+c')-(d+d')=0$$
$$(a+a')+(b+b')-(c+c')-(d+d')=0$$
$$(a+a')-(b+b')-(c+c')+(d+d')=0$$

Since the matrix of this latter system is the Walsh $W_4$, which is Hadamard, and so rescaled orthogonal, and in particular invertible, the solution is:
$$(a',b',c',d')=-(a,b,c,d)$$

Thus, in order to complete $I$ into a partial $5\times 8$ Hadamard matrix, we can pick any vector $(a,b,c,d)\in(\pm1)^4$, and then set $(a',b',c',d')=-(a,b,c,d)$.

\medskip

(2) Now let us try to complete $I$ into a full Hadamard matrix $H\in M_8(\pm1)$. By using the above observation, applied to each of the 4 lower rows of $H$, we conclude that $H$ must be of the following special form, with $L\in M_4(\pm1)$ being a certain matrix:
$$H=\begin{pmatrix}W_4&W_4\\L&-L\end{pmatrix}$$

Now observe that, in order for $H$ to be Hadamard, $L$ must be Hadamard. Thus, the solutions are those above, with $L\in M_4(\pm1)$ being Hadamard.

\medskip

(3) As a third step now, let us recall from Theorem 1.12 that we must have $L\sim W_4$. However, in relation with our problem, we cannot really use this in order to conclude directly that we have $H\sim W_8$. To be more precise, in  order not to mess up the structure of $I=(W_4\ W_4)$, we are allowed now to use only operations on the rows. And the conclusion here is that, up to equivalence, we have 2 solutions, as follows:
$$P=\begin{pmatrix}W_4&W_4\\W_4&-W_4\end{pmatrix}\quad,\quad 
Q=\begin{pmatrix}W_4&W_4\\K_4&-K_4\end{pmatrix}$$

We will see in moment that these two solutions are actually equivalent, but let us pause now our study of Case 1, after all this work done, and discuss Case 2.

\medskip

\underline{Case 2}. Here we must look for completions of the following matrix $J$:
$$J=\begin{pmatrix}
1&1&1&1&&-1&1&1&1\\ 
1&-1&1&-1&&1&-1&1&1\\ 
1&1&-1&-1&&1&1&-1&1\\ 
1&-1&-1&1&&1&1&1&-1
\end{pmatrix}$$

Let us first try to complete this partial $4\times 8$ Hadamard matrix into a partial $5\times 8$ Hadamard matrix. The completion must look as follows:
$$J'=\begin{pmatrix}
1&1&1&1&&-1&1&1&1\\ 
1&-1&1&-1&&1&-1&1&1\\ 
1&1&-1&-1&&1&1&-1&1\\ 
1&-1&-1&1&&1&1&1&-1\\
a&b&c&d&&x&y&z&t
\end{pmatrix}$$

The system of equations for the orthogonality conditions is as follows:
$$(1\perp 5)\quad:\quad a+b+c+d-x+y+z+t=0$$
$$(2\perp 5)\quad:\quad a-b+c-d+x-y+z+t=0$$
$$(3\perp 5)\quad:\quad a+b-c-d+x+y-z+t=0$$
$$(4\perp 5)\quad:\quad a-b-c+d+x+y+z-t=0$$

When regarded as a system in $x,y,z,t$, the matrix of the system is $K_4$, which is invertible. Thus, the vector $(x,y,z,t)$ is uniquely determined by the vector $(a,b,c,d)$:
$$(a,b,c,d)\to(x,y,z,t)$$

We have 16 vectors $(a,b,c,d)\in(\pm1)^4$ to be tried, and the first case, covering 8 of them, is that of the row vectors of $\pm W_4$. Here we have an obvious solution, with $(x,y,z,t)$ appearing at right of $(a,b,c,d)$ inside the following matrices, which are Hadamard: 
$$R=\begin{pmatrix}W_4&K_4\\W_4&-K_4\end{pmatrix}\quad,\quad 
S=\begin{pmatrix}W_4&K_4\\-W_4&K_4\end{pmatrix}$$

As for the second situation, this is that of the 8 binary vectors $(a,b,c,d)\in(\pm1)^4$ which are not row vectors of $\pm W_4$. But this is the same as saying that, up to permutations, we have $(a,b,c,d)=\pm(-1,1,1,1)$. In this latter case, and with $+$ sign, the system is:
$$-x+y+z+t=-2$$
$$x-y+z+t=2$$
$$x+y-z+t=2$$
$$x+y+z-t=2$$

By summing the first equation with the other ones we obtain the following system, whose solution is $y=z=t=0$, not corresponding to an Hadamard matrix:
$$y+z=y+t=z+t=0$$

Summarizing, we are done with the $5\times 8$ completion problem in Case 2, the solutions coming from the rows of the matrices $R,S$ above. Now when using this, as for getting up to full $8\times8$ completions, the $R,S$ cases obviously cannot mix, and so we are left with the Hadamard matrices $R,S$, as being the only solutions. In order to conclude now, observe that we have $R=Q^t$ and $R\sim S$. Also, we have $P\sim Q$, and this finishes the proof.
\end{proof}

The above proof was of course quite long. It is possible to improve a bit things, with various algebraic tricks, but basically this is how the situation is, with each classification result for the Hadamard matrices needing a lot of routine row-by-row study.

\section*{1c. Paley matrices}

We have seen that the Hadamard matrices can be classified up to order $N=8$, with the Walsh matrices being the only ones. We discuss now the case $N\geq12$, where new phenomena appear. At $N=12$ there is no Walsh matrix, but we can use a construction due to Paley \cite{pal}. Let $q=p^r$ be an odd prime power, consider the associated finite field $\mathbb F_q$, and then consider the quadratic character $\chi:\mathbb F_q\to\{-1,0,1\}$, given by: 
$$\chi(a)=\begin{cases}
0&{\rm if}\ a=0\\
1&{\rm if}\ a=b^2,b\neq0\\
-1&{\rm otherwise}
\end{cases}$$

\index{finite field}

We can construct then the following matrix, with indices in $\mathbb F_q$:
$$Q_{ab}=\chi(b-a)$$

With these conventions, the Paley construction of Hadamard matrices, which works at $N=12$ and at many other values of $N\in4\mathbb N$, is as follows:

\index{Paley matrix}
\index{symmetric matrix}
\index{skew-symmetric}

\begin{theorem}
Given an odd prime power $q=p^r$, construct $Q_{ab}=\chi(b-a)$ as above. We have then constructions of Hadamard matrices, as follows:
\begin{enumerate}
\item Paley $1$: if $q=3(4)$ we have a matrix of size $N=q+1$, as follows:
$$P_N^1=1+\begin{pmatrix}
0&1&\ldots&1\\
-1\\
\vdots&&Q\\
-1
\end{pmatrix}$$

\item Paley $2$: if $q=1(4)$ we have a matrix of size $N=2q+2$, as follows:
$$P_N^2=\begin{pmatrix}
0&1&\ldots&1\\
1\\
\vdots&&Q\\
1
\end{pmatrix}\quad:\quad 0\to\begin{pmatrix}1&-1\\ -1&-1\end{pmatrix}\quad,\quad\pm1\to\pm\begin{pmatrix}1&1\\1&-1\end{pmatrix}$$
\end{enumerate}
These matrices are skew-symmetric $(H+H^t=2)$, respectively symmetric $(H=H^t)$.
\end{theorem}

\begin{proof}
In order to simplify the presentation, we will denote by $1$ all the identity matrices, of any size, and by $\mathbb I$ all the rectangular all-one matrices, of any size as well. It is elementary to check that the matrix $Q_{ab}=\chi(a-b)$ has the following properties:
$$QQ^t=q1-\mathbb I$$
$$Q\mathbb I=\mathbb IQ=0$$

In addition, we have the following formulae, which are elementary as well, coming from the fact that $-1$ is a square in $\mathbb F_q$ precisely when $q=1(4)$:
$$q=1(4)\implies Q=Q^t\ \ \,$$
$$q=3(4)\implies Q=-Q^t$$

With these observations in hand, the proof goes as follows:

\medskip

(1) With our conventions for the symbols $1$ and $\mathbb I$, explained above, the matrix in the statement is as follows:
$$P_N^1=\begin{pmatrix}1&\mathbb I\\ -\mathbb I&1+Q\end{pmatrix}$$
 
With this formula in hand, the Hadamard matrix condition follows from:
\begin{eqnarray*}
P_N^1(P_N^1)^t
&=&\begin{pmatrix}1&\mathbb I\\ -\mathbb I&1+Q\end{pmatrix}\begin{pmatrix}1&-\mathbb I\\ \mathbb I&1-Q\end{pmatrix}\\
&=&\begin{pmatrix}N&0\\ 0&\mathbb I+1-Q^2\end{pmatrix}\\
&=&\begin{pmatrix}N&0\\ 0&N\end{pmatrix}
\end{eqnarray*}

(2) If we denote by $G,F$ the matrices in the statement, which replace respectively the $0,1$ entries, then we have the following formula for our matrix:
$$P_N^2=\begin{pmatrix}0&\mathbb I\\ \mathbb I&Q\end{pmatrix}\otimes F+1\otimes G$$

With this formula in hand, the Hadamard matrix condition follows from:
\begin{eqnarray*}
(P_N^2)^2
&=&\begin{pmatrix}0&\mathbb I\\ \mathbb I&Q\end{pmatrix}^2\otimes F^2+\begin{pmatrix}1&0\\ 0&1\end{pmatrix}\otimes G^2+\begin{pmatrix}0&\mathbb I\\ \mathbb I&Q\end{pmatrix}\otimes(FG+GF)\\
&=&\begin{pmatrix}q&0\\ 0&q\end{pmatrix}\otimes 2+\begin{pmatrix}1&0\\ 0&1\end{pmatrix}\otimes 2+\begin{pmatrix}0&\mathbb I\\ \mathbb I&Q\end{pmatrix}\otimes0\\
&=&\begin{pmatrix}N&0\\ 0&N\end{pmatrix}
\end{eqnarray*}

Finally, the last assertion is clear, from the above formulae relating $Q,Q^t$.
\end{proof}

As an illustration for the above result, we have:

\begin{theorem}
We have Paley $1$ and $2$ matrices at $N=12$, which are equivalent:
$$P_{12}^1\sim P_{12}^2$$
In fact, this matrix is the unique Hadamard one at $N=12$, up to equivalence.
\end{theorem}

\begin{proof}
We have $12=11+1$, with $11=3(4)$ being prime, so the Paley 1 construction applies indeed, with the first row vector of $Q$ being:
$$q=(0+-+++---+-)$$

Also, we have $12=2\times 5+2$, with $5=1(4)$ being prime, so the Paley 2 construction applies as well, with the first row vector of $Q$ being:
$$q=(0+--+)$$

It is routine then to check that we have $P_{12}^1\sim P_{12}^2$, by some computations in the spirit of those from the end of the proof of Theorem 1.18, and with the matrix $P_{12}\sim P_{12}^1\sim P_{12}^2$ being as follows, with the $\pm$ signs standing for $\pm1$ entries:
$$P_{12}=\left(
\begin{array}{ccccccccccccccc}
+&+&+&+&&-&+&+&+&&-&+&+&+\\
+&-&+&-&&+&-&+&+&&+&-&+&+\\
+&+&-&-&&+&+&-&+&&+&+&-&+\\
+&-&-&+&&+&+&+&-&&+&+&+&-\\
\\
-&+&-&-&&+&-&+&+&&-&+&+&-\\
-&+&+&+&&+&+&-&+&&+&-&+&-\\
-&-&-&+&&+&+&+&+&&-&-&-&+\\
+&+&+&-&&+&+&+&-&&-&-&-&-\\
\\
-&+&+&+&&+&-&+&-&&+&+&-&+\\
+&+&-&+&&+&-&-&-&&-&-&+&+\\
+&-&+&+&&+&-&-&+&&-&+&-&-\\
-&-&+&-&&+&+&-&-&&-&+&+&+
\end{array}\right)$$
 
As for the last assertion, regarding uniqueness, this is something quite technical, requiring some clever block decomposition techniques. Alternatively, it is possible to verify this by using a computer, although programming such things is not exactly trivial.
\end{proof}

At $N=16$ now, the situation becomes fairly complicated, as follows:

\begin{theorem}
The Hadamard matrices at $N=16$ are as follows:
\begin{enumerate}
\item We have the Walsh matrix $W_{16}$.

\item There are no Paley matrices.

\item Besides $W_{16}$, we have $4$ more matrices, up to equivalence.
\end{enumerate}
\end{theorem}

\begin{proof}
Once again, this is a mixture of elementary and more advanced results:

\medskip

(1) This is clear.

\medskip

(2) This comes from the fact that we have $16=15+1$, with $15$ not being a prime power, and from the fact that we have $16=2\times 7+2$, with $7\neq1(4)$.

\medskip

(3) This is something very technical, basically requiring a computer.
\end{proof}

At $N=20$ and bigger, the situation becomes quite complicated, and the study is usually done with a mix of advanced algebraic methods, and computer techniques. The overall conclusion is that the number of Hadamard matrices of size $N\in4\mathbb N$ grows with $N$, in exponential fashion. In particular, we are led in this way into:

\index{HC}
\index{Hadamard Conjecture}

\begin{conjecture}[Hadamard Conjecture (HC)]
There is at least one Hadamard matrix 
$$H\in M_N(\pm1)$$
for any integer $N\in 4\mathbb N$.
\end{conjecture}

This conjecture, going back to the 19th century, is one of the most beautiful statements in combinatorics, linear algebra, and mathematics in general. Quite remarkably, the numeric verification so far goes up to the number of the beast:
$$\mathfrak N=666$$

\index{number of the beast}

Our purpose now will be that of gathering some evidence for this conjecture. By using the Walsh construction, we have examples at each $N=2^n$. We can add various examples coming from the Paley 1 and Paley 2 constructions, and we are led to:

\index{Paley matrix}
\index{Walsh matrix}

\begin{theorem}
The HC is verified at least up to $N=88$, as follows:
\begin{enumerate}
\item At $N=4,8,16,32,64$ we have Walsh matrices.

\item At $N=12,20,24,28,44,48,60,68,72,80,84,88$ we have Paley $1$ matrices.

\item At $N=36,52,76$ we have Paley $2$ matrices.

\item At $N=40,56$ we have Paley $1$ matrices tensored with $W_2$.
\end{enumerate}
However, at $N=92$ these constructions (Walsh, Paley, tensoring) don't work.
\end{theorem}

\begin{proof}
First of all, the numbers in (1-4) are indeed all the multiples of 4, up to 88. As for the various assertions, the proof here goes as follows:

\medskip

(1) This is clear.

\medskip

(2) Here the number $N-1$ takes the following values:
$$q=11,19,23,27,43,47,59,67,71,79,83,87$$

These are all prime powers, so we can apply the Paley 1 construction.

\medskip

(3) Since $N=4(8)$ here, and $N/2-1$ takes the values $q=17,25,37$, all prime powers, we can indeed apply the Paley 2 construction, in these cases.

\medskip

(4) At $N=40$ we have indeed $P_{20}^1\otimes W_2$, and at $N=56$ we have $P_{28}^1\otimes W_2$.

\medskip

Finally, we have $92-1=7\times13$, so the Paley 1 construction does not work, and $92/2=46$, so the Paley 2 construction, or tensoring with $W_2$, does not work either.
\end{proof}

At $N=92$ now, the situation is considerably more complicated, and we have:

\index{Williamson matrix}
\index{circulant matrix}
\index{quaternion units}

\begin{theorem}
Assuming that $A,B,C,D\in M_K(\pm1)$ are circulant, symmetric, pairwise commute and satisfy the condition
$$A^2+B^2+C^2+D^2=4K$$
the following $4K\times4K$ matrix is Hadamard, called of Williamson type:
$$H=\begin{pmatrix}
A&B&C&D\\
-B&A&-D&C\\
-C&D&A&-B\\
-D&-C&B&A
\end{pmatrix}$$
Moreover, matrices $A,B,C,D$ as above exist at $K=23$, where $4K=92$.
\end{theorem}

\begin{proof}
We use the same method as for the Paley theorem, namely tensor calculus. Consider the following matrices $1,i,j,k\in M_4(0,1)$, called the quaternion units:
$$1=\begin{pmatrix}
1&0&0&0\\
0&1&0&0\\
0&0&1&0\\
0&0&0&1
\end{pmatrix}
\qquad,\qquad
i=\begin{pmatrix}
0&1&0&0\\
1&0&0&0\\
0&0&0&1\\
0&0&1&0
\end{pmatrix}$$
$$j=\begin{pmatrix}
0&0&1&0\\
0&0&0&1\\
1&0&0&0\\
0&1&0&0
\end{pmatrix}
\qquad,\qquad
k=\begin{pmatrix}
0&0&0&1\\
0&0&1&0\\
0&1&0&0\\
1&0&0&0
\end{pmatrix}$$

These matrices describe the positions of the $A,B,C,D$ entries in the matrix $H$ from the statement, and so this matrix can be written as follows:
$$H=A\otimes 1+B\otimes i+C\otimes j+D\otimes k$$

Assuming now that $A,B,C,D$ are symmetric, we have:
\begin{eqnarray*}
HH^t
&=&(A\otimes 1+B\otimes i+C\otimes j+D\otimes k)\\
&&(A\otimes 1-B\otimes i-C\otimes j-D\otimes k)\\
&=&(A^2+B^2+C^2+D^2)\otimes 1-([A,B]-[C,D])\otimes i\\
&&-([A,C]-[B,D])\otimes j-([A,D]-[B,C])\otimes k
\end{eqnarray*}

Now assume that our matrices $A,B,C,D$ pairwise commute, and satisfy as well the condition in the statement, namely $A^2+B^2+C^2+D^2=4K$. In this case, it follows from the above formula that we have $HH^t=4K$, so we obtain indeed an Hadamard matrix.

\medskip

In general, finding such matrices is a difficult task, and this is where Williamson's extra assumption that $A,B,C,D$ should be taken circulant comes from. Finally, regarding the $K=23$ construction, which produces an Hadamard matrix of order $N=92$, this comes via a computer search. See Williamson \cite{wil} and Baumert-Golomb-Hall \cite{bgh}. 
\end{proof}

Things get even worse at higher values of $N$, where more and more complicated constructions are needed. The whole subject is quite technical, and, as already mentioned, human knowledge here stops so far at $\mathfrak N=666$. See \cite{aga}, \cite{dda}, \cite{dgo}, \cite{hor}, \cite{kta}, \cite{sya}.

\index{number of the beast}

\section*{1d. Cocyclic matrices}

We have seen so far that the combinatorial and algebraic theory of the Hadamard matrices, while very nice at the elementary level, ultimately leads into some difficult questions. There are at least two potential exits from this, namely:

\bigskip

(1) Do analysis. There are many things that can be done here, starting with the Hadamard determinant bound \cite{had}, and we will discuss this in chapter 2, and afterwards. Whether all this can help or not in relation with the Hadamard Conjecture remains to be seen, but at least we'll have some fun, and do some interesting mathematics.

\bigskip

(2) Do geometry. When allowing the entries of $H$ to be complex numbers, we reach to geometric questions, and the Hadamard Conjecture problematics dissapears, because the Fourier matrix, namely $F_N=(w^{ij})$ with $w=e^{2\pi i/N}$, is an example of such matrix at any $N\in\mathbb N$. We will discuss this later, starting from chapter 5 below.

\bigskip

Getting back now to algebra and combinatorics, as a conceptual finding on the subject, however, we have the recent theory of the cocyclic Hadamard matrices, that we will briefly explain now. This theory is based on the following notion:

\index{cocycle}
\index{cocyclic matrix}

\begin{definition}
A cocycle on a finite group $G$ is a matrix $H\in M_G(\pm1)$ satisfying:
$$H_{gh}H_{gh,k}=H_{g,hk}H_{hk}$$
$$H_{11}=1$$
If the rows of $H$ are pairwise orthogonal, we say that $H$ is a cocyclic Hadamard matrix.
\end{definition}

Here the definition of the cocycles is the usual one, with the equations coming from the fact that $F=\mathbb Z_2\times G$ must be a group, with multiplication as follows:
$$(u,g)(v,h)=(H_{gh}\cdot uv,gh)$$

\index{Walsh matrix}

As a basic illustration for the above notion, the Walsh matrix $H=W_{2^n}$ is cocyclic, coming from the group $G=\mathbb Z_2^n$, with cocycle as follows:
$$H_{gh}=(-1)^{<g,h>}$$

As explained by de Launey, Flannery and Horadam in \cite{dfh}, and in other papers, many other known examples of Hadamard matrices are cocyclic, and this leads to:

\begin{conjecture}[Cocyclic Hadamard Conjecture]
There is at least one cocyclic Hadamard matrix $H\in M_N(\pm1)$, for any $N\in 4\mathbb N$.
\end{conjecture}

\index{Cocyclic Hadamard Conjecture}

Having such a statement formulated is certainly a big advance with respect to the HC, and this is probably the main achievement of modern Hadamard matrix theory. However, in what regards a potential proof, there is no clear strategy here, at least so far. We will be back to such questions, in relation with advanced algebra, in chapters 13-16 below, with the fact that the construction $\mathbb Z_2^n\to W_{2^n}$ can be extended as to cover all the Hadamard matrices, by replacing $\mathbb Z_2^n$ with a suitable quantum permutation group. However, in what regards the potential applications to the HC, there is no clear strategy here either.

\bigskip

Finally, as a last algebraic topic, let us discuss the Circulant Hadamard Conjecture. Besides analysis in a large sense, as explained above, another potential way of getting away from the difficult HC questions is that of looking at various special classes of Hadamard matrices. However, in practice, this often leads to quite complicated mathematics too.

\bigskip

\index{circulant matrix}
\index{CHC}
\index{Circulant Hadamard Conjecture}

Illustrating and famous here is the situation in the circulant case. Given a vector $\gamma\in(\pm 1)^N$, one can ask whether the corresponding circulant matrix $H\in M_N(\pm 1)$, defined by $H_{ij}=\gamma_{j-i}$, is Hadamard or not. Here is a solution to the problem:
$$K_4=\begin{pmatrix}-1&1&1&1\\ 1&-1&1&1\\ 1&1&-1&1\\ 1&1&1&-1\end{pmatrix}$$

More generally, any vector $\gamma\in(\pm 1)^4$ satisfying $\sum\gamma_i=\pm 1$ is a solution to the problem, with the corresponding Hadamard matrix being equivalent to $K_4$. The following conjecture, due to Ryser \cite{rys}, states that there are no other solutions:

\begin{conjecture}[Circulant Hadamard Conjecture (CHC)]
There is no circulant Hadamard matrix of size $N\times N$, for any $N\neq 4$.
\end{conjecture}

The fact that such a simple-looking problem is still open might seem quite surprising. Indeed, if we denote by $S\subset\{1,\ldots,N\}$ the set of positions of the $-1$ entries of $\gamma$, the Hadamard matrix condition is simply $|S\cap(S+k)|=|S|-N/4$, for any $k\neq 0$, taken modulo $N$. Thus, the above conjecture simply states that at $N\neq 4$, such a set $S$ cannot exist. Let us record here this latter statement, also due to Ryser \cite{rys}:

\index{Ryser Conjecture}

\begin{conjecture}[Ryser Conjecture]
Given an integer $N>4$, there is no set $S\subset\{1,\ldots,N\}$ satisfying the condition
$$|S\cap(S+k)|=|S|-N/4$$
for any $k\neq 0$, taken modulo $N$.
\end{conjecture}

There has been a lot of work on this conjecture, starting with \cite{rys}. However, as it was the case with the HC, all this leads to complicated combinatorics, design theory, algebra and number theory, and so on, and there is no clear idea here, at least so far.

\section*{1e. Exercises}

There has been a lot of linear algebra and combinatorics in this chapter, and doing some more linear algebra and combinatorics will be our purpose here. First we have: 

\begin{exercise}
Verify that we have indeed the formula
$$H\otimes(K\otimes L)=(H\otimes K)\otimes L$$
when using the lexicographic order on the triple indices.
\end{exercise}

This is a very instructive exercise, making you familiar with tensor products and multiple indices, and with this knowledge being a very useful asset.

\begin{exercise}
Write down an explicit equivalence $K_4\sim W_4$.
\end{exercise}

This is normally something quite simple, just some fun with basic matrices.

\begin{exercise}
Write down the matrix $P_4^1$, and prove that $P_4^1\sim W_4$.
\end{exercise}

Again, this is something elementary, with just a bit of thinking been needed at the beginning, in order to figure out what the Paley matrix $P_4^1$ exactly is.

\begin{exercise}
Write down the matrix $P_8^1$, and prove that $P_8^1\sim W_8$.
\end{exercise}

This is certainly more difficult than the previous two exercises, but surely can be done, either by following some ideas from our classification at $N=8$, or by doing it directly, using your intuition. In case you want to ``cheat'' by using a computer, you are of course welcome to do so, because programming such things is very instructive too.

\begin{exercise}
Prove that we have $P_{12}^1\sim P_{12}^2$.
\end{exercise}

Finally, a more advanced question is that of looking at the various examples of Hadamard matrices constructed in this chapter, and see which of them are cocyclic.

\chapter{Analytic aspects}

\section*{2a. Determinant bound}

We have seen so far that the algebraic theory of the Hadamard matrices, while very nice at the elementary level, ultimately leads to some difficult questions. So, let us step now into analytic questions. The first result here, found in 1893 by Hadamard \cite{had}, about 25 years after Sylvester's 1867 founding paper \cite{syl}, and which actually led to such matrices being called Hadamard, is a determinant bound, as follows:

\index{Hadamard theorem}
\index{Hadamard determinant bound}
\index{determinant bound}
\index{Hadamard matrix}
\index{volume of parallelepiped}
\index{Sylvester}

\begin{theorem}
Given a matrix $H\in M_N(\pm1)$, we have
$$|\det H|\leq N^{N/2}$$
with equality precisely when $H$ is Hadamard.
\end{theorem}

\begin{proof}
We use here the fact, which often tends to be forgotten, that the determinant of a system of $N$ vectors in $\mathbb R^N$ is the signed volume of the associated parallelepiped:
$$\det(H_1,\ldots,H_N)=\pm vol<H_1,\ldots,H_N>$$

This is actually the definition of the determinant, in case you have forgotten the basics, with the need for the sign coming for having good additivity properties. Now in the case where our vectors have their entries in $\{\pm1\}$, we therefore have the following inequality, with equality precisely when our vectors are pairwise orthogonal:
\begin{eqnarray*}
|\det(H_1,\ldots,H_N)|
&\leq&||H_1||\times\ldots\times||H_N||\\
&=&(\sqrt{N})^N
\end{eqnarray*}

Thus, we have obtained the result, straight from the definition of $\det$.
\end{proof}

The above result is quite interesting, philosophically speaking. Let us recall indeed from chapter 1 that the set formed by the $N\times N$ Hadamard matrices is:
$$Y_N=M_N(\pm1)\cap\sqrt{N}O_N$$

Thus, what we have in Theorem 2.1 is an analytic method for locating this Hadamard matrix set $Y_N$ inside the space of binary matrices $M_N(\pm1)$. But this suggests doing several other analytic things, as for instance looking at the maximizers $H\in M_N(\pm1)$ of the quantity $|\det H|$, at values $N\in\mathbb N$ which are not multiples of 4. Things here are quite tricky, and as a basic result on the subject, at $N=3$ the situation is as follows:

\begin{proposition}
For a matrix $H\in M_3(\pm1)$ we have $|\det H|\leq4$, and this estimate is sharp, with the equality case being attained by the matrix
$$Q_3=\begin{pmatrix}
1&1&1\\
1&1&-1\\
1&-1&1
\end{pmatrix}$$
and its conjugates, via the Hadamard equivalence relation.
\end{proposition}

\begin{proof}
In order to get started, observe that Theorem 2.1 provides us with the following bound, which is of course not sharp, $\det H$ being an integer:
$$|\det H|\leq 3\sqrt{3}=5.1961..$$

Now observe that, $\det H$ being a sum of six $\pm1$ terms, it must be an even number. Thus, we obtain the estimate in the statement, namely:
$$|\det H|\leq4$$

Our claim now is that the following happens, with the nonzero situation appearing precisely for the matrix $Q_3$ in the statement, and its conjugates:
$$\det H\in\{-4,0,4\}$$

Indeed, let us try to find the matrices $H\in M_3(\pm1)$ having the property $\det H\neq0$. Up to equivalence, we can assume that the first row is $(1,1,1)$. Then, once again up to equivalence, we can assume that the second row is $(1,1,-1)$. And then, once again up to equivalence, we can assume that the third row is $(1,-1,1)$. Thus, we must have:
$$H=\begin{pmatrix}
1&1&1\\
1&1&-1\\
1&-1&1
\end{pmatrix}$$

The determinant of this matrix being $-4$, we have proved our claim, and the last assertion in the statement too, as a consequence of our study.
\end{proof}

In general, all this suggests the following definition:

\index{quasi-Hadamard matrix}
\index{maximizer of determinant}

\begin{definition}
A quasi-Hadamard matrix is a square binary matrix
$$H\in M_N(\pm1)$$
which maximizes the quantity $|\det H|$.
\end{definition}

We know from Theorem 2.1 that at $N\in 4\mathbb N$ such matrices are precisely the Hadamard matrices, provided that the Hadamard Conjecture holds at $N$. At values $N\notin 4\mathbb N$, what we have here are certain matrices which can be thought of as being ``generalized Hadamard matrices'', the simplest examples being the matrix $Q_3$ from Proposition 2.2, and its Hadamard conjugates. For more on all this, we refer to Park-Song \cite{pso}.

\bigskip

As a comment, however, Proposition 2.2 might look a bit dissapointing, because it is hard to imagine that the matrix $Q_3$ there, which is not a very interesting matrix, can really play the role of a ``generalized Hadamard matrix'' at $N=3$. We will come later with more interesting solutions to this latter problem, a first solution being as follows:
$$K_3=\frac{1}{\sqrt{3}}\begin{pmatrix}
-1&2&2\\
2&-1&2\\
2&2&-1
\end{pmatrix}$$

To be more precise, this matrix is of course not binary, but it is definitely an interesting matrix, that we will see to be sharing many properties with the Hadamard matrices. Also, we have as well another solution to the $N=3$ problem, which uses complex numbers, and more specifically the number $w=e^{2\pi i/3}$, which is as follows:
$$F_3=\begin{pmatrix}
1&1&1\\
1&w&w^2\\
1&w^2&w
\end{pmatrix}$$

\index{Fourier matrix}

As a conclusion to this study, looking at the maximizers $H\in M_N(\pm1)$ of the quantity $|\det H|$ is not exactly an ideal method, when looking for analogues of the Hadamard matrices at the forbidden size values $N\notin 4\mathbb N$, at least when $N$ is small. The situation changes, however, when looking at such questions at big values of $N\in\mathbb N$, where the determinant problematics for the binary matrices becomes very interesting, and quite technical. As a generic statement here, which is a bit informal, we have:

\index{determinant bound}

\begin{theorem}
We have, in the $N\to\infty$ limit, 
$$\max_{H\in M_N(\pm1)}|\det H|\simeq N^{N/2}$$
along with even finer estimates, modulo the Hadamard Conjectuere.
\end{theorem}

\begin{proof}
As mentioned, this is just an informal statement, standing here as a modest introduction to the subject, in the lack of something more precise, and elementary. There are basically two ways of dealing with such questions, namely:

\medskip

(1) A first idea, as mentioned, is that of using the existence of an Hadamard matrix $H_N\in M_N(\pm1)$, at values $N\in4\mathbb N$, modulo the Hadamard Conjecture of course, and then completing it into binary matrices $H_{N+k}\in M_{N+k}(\pm1)$, with $k=1,2,3$:
$$H_{N+k}=
\begin{pmatrix}
&&&*\\
&H_N&&*\\
&&&*\\
*&*&*&*
\end{pmatrix}$$

The determinant estimates for such matrices are however quite technical, and we refer here once again to Park-Song \cite{pso}, and related papers.

\medskip

(2) A second method is by using probability theory. The set of binary matrices $M_N(\pm1)$ is a probability space, when endowed with the counting measure rescaled by $1/2^{N^2}$, and the determinant can be regarded as a random variable on this space:
$$\det:M_N(\pm1)\to\mathbb Z$$

The point now is that the distribution of this variable can be computed, in the $N\to\infty$ limit, and as a consequence, we can investigate the maximizers of $|\det H|$. Once again, all this is quite technical, and we refer here to Tao-Vu \cite{tvu} and related papers.
\end{proof}

\index{Tao-Vu}

Summarizing, the Hadamard determinant bound provides us with an analytic method of locating the set $Y_N=M_N(\pm1)\cap\sqrt{N}O_N$ formed by the $N\times N$ Hadamard matrices inside $M_N(\pm1)$, and this leads to an interesting $N\to\infty$ theory.

\section*{2b. Norm maximizers}

From a ``dual'' point of view, the question of locating $Y_N$ inside $\sqrt{N}O_N$, once again via analytic methods, makes sense as well. The result here, from \cite{bcs}, is as follows:

\index{norm maximizer}

\begin{theorem}
Given a matrix $U\in O_N$ we have
$$||U||_1\leq N\sqrt{N}$$
with equality precisely when $H=\sqrt{N}U$ is Hadamard.
\end{theorem}

\begin{proof}
We have indeed the following estimate, for any $U\in O_N$, which uses the Cauchy-Schwarz inequality, and the trivial fact that we have $||U||_2=\sqrt{N}$:
\begin{eqnarray*}
||U||_1
&=&\sum_{ij}|U_{ij}|\\
&\leq&N\left(\sum_{ij}|U_{ij}|^2\right)^{1/2}\\
&=&N\sqrt{N}
\end{eqnarray*}

In addition, we know that the equality case holds when the variables are equal, and so when $|U_{ij}|=1/\sqrt{N}$, for any $i,j$. But this amounts in saying that $H=\sqrt{N}U$ must satisfy $H\in M_N(\pm1)$. Thus, this rescaled matrix $H$ must be Hadamard, as claimed.
\end{proof} 

We will need more general norms as well, so let record the following result:

\index{convex function}
\index{concave function}
\index{Jensen inequality}

\begin{proposition}
If $\psi:[0,\infty)\to\mathbb R$ is strictly concave/convex, the quantity
$$F(U)=\sum_{ij}\psi(U_{ij}^2)$$
over $U_N$ is maximized/minimized by the rescaled Hadamard matrices, $U=H/\sqrt{N}$.
\end{proposition}

\begin{proof}
We recall that the Jensen theorem states that for $\psi$ convex we have the following inequality, with equality, when $\psi$ is strictly convex, when $x_i$ are all equal:
$$\psi\left(\frac{x_1+\ldots+x_n}{n}\right)\leq\frac{\psi(x_1)+\ldots+\psi(x_n)}{n}$$

In our case, let us take $n=N^2$, and our variables to be as follows:
$$\left\{x_1,\ldots,x_n\right\}=\left\{U_{ij}^2\big|i,j=1,\ldots,N\right\}$$

We obtain that for any convex function $\psi$, the following holds:
$$\psi\left(\frac{1}{N}\right)\leq\frac{F(U)}{N^2}$$

Thus we have the following estimate, with $F$ being as in the statement:
$$F(U)\geq N^2\psi\left(\frac{1}{N}\right)$$

Now if $\psi$ is strictly convex, the equality case holds when the numbers $U_{ij}^2$ are all equal, so when $H=\sqrt{N}U$ is Hadamard. The proof for concave functions is similar.
\end{proof}

Of particular interest for us are the following consequences of Proposition 2.6: 

\index{maximizer of p-norm}
\index{minimizer of p-norm}

\begin{theorem}
The rescaled versions $U=H/\sqrt{N}$ of the Hadamard matrices $H\in M_N(\pm1)$ can be characterized as being:
\begin{enumerate}
\item The maximizers of the $p$-norm on $O_N$, at any $p\in[1,2)$.

\item The minimizers of the $p$-norm on $O_N$, at any $p\in(2,\infty]$.
\end{enumerate}
\end{theorem} 

\begin{proof}
Consider indeed the $p$-norm on $O_N$, which at $p\in[1,\infty)$ is given by:
$$||U||_p=\left(\sum_{ij}|U_{ij}|^p\right)^{1/p}$$

Since $\psi(x)=x^{p/2}$ is concave at $p\in[1,2)$, and convex at $p\in(2,\infty)$, Proposition 2.6 applies and gives the results at $p\in[1,\infty)$, the precise estimates being as follows:
$$||U||_p:
\begin{cases}
\leq N^{2/p-1/2}&{\rm if}\ p<2\\
=N^{1/2}&{\rm if}\ p=2\\
\geq N^{2/p-1/2}&{\rm if}\ p>2
\end{cases}$$

As for the case $p=\infty$, this follows either by letting $p\to\infty$ in the above estimates, or directly via Cauchy-Schwarz, a bit as in the proof of Theorem 2.5.
\end{proof}

As it was the case with the Hadamard determinant bound, all this suggests doing some further geometry and analysis, this time on the Lie group $O_N$, with a notion of ``almost Hadamard matrix'' at stake. Let us formulate indeed, in analogy with Definition 2.3:

\index{AHM}
\index{almost Hadamard matrix}
\index{optimal AHM}
\index{optimal almost Hadamard matrix}
\index{norm maximizer}

\begin{definition}
An optimal almost Hadamard matrix is a rescaled orthogonal matrix
$$H\in\sqrt{N}O_N$$
which maximizes the $1$-norm.
\end{definition}

Here the adjective ``optimal'' comes from the fact that, in contrast with what happens over $M_N(\pm1)$, in connection with the determinant bound, here over $\sqrt{N}O_N$ we have more flexibility, and we can talk if we want about the local maximizers of the 1-norm. These latter matrices are called ``almost Hadamard'', and we will investigate them in the next chapter. Also, we will talk there about more general $p$-norms as well.

\bigskip

We know from Theorem 2.6 that at $N\in 4\mathbb N$ the absolute almost Hadamard matrices are precisely the Hadamard matrices, provided that the Hadamard Conjecture holds at $N$. At values $N\notin 4\mathbb N$, what we have are certain matrices which can be thought of as being ``generalized Hadamard matrices'', and are waiting to be investigated. Let us begin with a preliminary study, at $N=3$. The result here, from \cite{bcs}, is as follows:

\begin{theorem}
For any matrix $U\in O_3$ we have the estimate
$$||U||_1\leq 5$$
and this is sharp, with the equality case being attained by the matrix
$$U=\frac{1}{3}\begin{pmatrix}
-1&2&2\\
2&-1&2\\
2&2&-1
\end{pmatrix}$$
and its conjugates, via the Hadamard equivalence relation.
\end{theorem}

\begin{proof}
By dividing by $\det U$, we can assume that we have $U\in SO_3$. We use the Euler-Rodrigues parametrization for the elements of $SO_3$, namely:
$$U=\begin{pmatrix}
x^2+y^2-z^2-t^2&2(yz-xt)&2(xz+yt)\\
2(xt+yz)&x^2+z^2-y^2-t^2&2(zt-xy)\\
2(yt-xz)&2(xy+zt)&x^2+t^2-y^2-z^2
\end{pmatrix}$$

\index{Euler-Rodrigues}

Here $(x,y,z,t)\in S^3$ come from the map $SU_2\to SO_3$. Now in order to obtain the estimate, we linearize. We must prove that for any numbers $x,y,z,t\in\mathbb R$ we have:
\begin{eqnarray*}
&&|x^2+y^2-z^2-t^2|+|x^2+z^2-y^2-t^2|+|x^2+t^2-y^2-z^2|\\
&&+2\left(|yz-xt|+|xz+yt|+|xt+yz|+|zt-xy|+|yt-xz|+|xy+zt|\right)\\
&&\leq 5(x^2+y^2+z^2+t^2)
\end{eqnarray*}

The problem being symmetric in $x,y,z,t$, and invariant under sign changes, we may assume that we have:
$$x\geq y\geq z\geq t\geq 0$$

Now if we look at the 9 absolute values in the above formula, in 7 of them the sign is known, and in the remaining 2 ones the sign is undetermined. More precisely, the inequality to be proved is as follows:
\begin{eqnarray*}
&&(x^2+y^2-z^2-t^2)+(x^2+z^2-y^2-t^2)+|x^2+t^2-y^2-z^2|\\
&&+2\left(|yz-xt|+(xz+yt)+(xt+yz)+(xy-zt)+(xz-yt)+(xy+zt)\right)\\
&&\leq 5(x^2+y^2+z^2+t^2)
\end{eqnarray*}

After simplification and rearrangement of the terms, this inequality reads:
\begin{eqnarray*}
&&|x^2+t^2-y^2-z^2|+2|xt-yz|\\
&\leq&3x^2+5y^2+5z^2+7t^2-4xy-4xz-2xt-2yz
\end{eqnarray*}

In principle we have now 4 cases to discuss, depending on the possible signs appearing at left. It is, however, easier to proceed simply by searching for the optimal case. First, by writing $y=\alpha+\varepsilon,z=\alpha-\varepsilon$ and by making $\varepsilon$ vary over the real line, we see that the optimal case is when $\varepsilon=0$, hence when $y=z$. The cases $y=z=0$ and $y=z=\infty$ being both clear, and not sharp, we can assume that we have:
$$y=z=1$$

Thus we must prove that for any numbers $x\geq 1\geq t\geq 0$ we have:
$$|x^2+t^2-2|+2|xt-1|\leq 3x^2+8+7t^2-8x-2xt$$

In the case $xt\geq 1$ we have $x^2+t^2\geq 2$, and the inequality becomes:
$$2xt+4x\leq x^2+3t^2+6$$

In the case $xt\leq 1,x^2+t^2\leq 2$ we get:
$$x^2+1+2t^2\geq 2x$$

In the remaining case $xt\leq 1,x^2+t^2\geq 2$ we get:
$$x^2+4+3t^2\geq 4x$$

But these inequalities are all true, and this finishes the proof of the estimate. Now regarding the maximum, we know that this is attained at $(xyzt)=(1110)$ or at $(xyzt)=(2110)$, plus permutations. The corresponding matrix is, modulo permutations:
$$V=\frac{1}{3}
\begin{pmatrix}
1&2&2\\
2&1&-2\\
-2&2&-1
\end{pmatrix}$$

But for this matrix we have indeed $||V||_1=5$, and we are done.
\end{proof}

In terms of Definition 2.8, the conclusion is as follows:

\begin{theorem}
The optimal almost Hadamard matrices at $N=3$ are
$$K_3=\frac{1}{\sqrt{3}}\begin{pmatrix}
-1&2&2\\
2&-1&2\\
2&2&-1
\end{pmatrix}$$
and its conjugates, via the Hadamard equivalence relation.
\end{theorem}

\begin{proof}
This is indeed a reformulation of Theorem 2.9, using Definition 2.8.
\end{proof}

The above result and the matrix $K_3$ appearing there are quite interesting, because they remind the Hadamard matrix $K_4$ studied in chapter 1, given by:
$$K_4=\begin{pmatrix}
-1&1&1&1\\
1&-1&1&1\\
1&1&-1&1\\
1&1&1&-1
\end{pmatrix}$$

To be more precise, all this suggests looking at the following remarkable family of matrices $K_N\in\sqrt{N}O_N$, having arbitrary size $N\in\mathbb N$:
$$K_N=\frac{1}{\sqrt{N}}\begin{pmatrix}
2-N&&2\\
&\ddots&\\
2&&2-N
\end{pmatrix}$$

These matrices are in general not optimal almost Hadamard, in the sense of Definition 2.8, for instance because at $N=2$ or at $N=8,12,16,\ldots$ they are obviously not Hadamard. We will see however in the next chapter that these matrices are ``almost Hadamard'', in the sense that they locally maximize the 1-norm on $\sqrt{N}O_N$.

\bigskip

To summarize, the computation of the maximizers of the 1-norm on $O_N$ is a difficult question, a bit like the computation of the maximizers of $|\det|$ on $M_N(\pm1)$ was, and looking instead at the local maximizers of the 1-norm on $O_N$ is the way to be followed, with some interesting examples and combinatorics at stake. We will be back to this.

\bigskip

Let us discuss now, as a continuation of all this, an analytic reformulation of the Hadamard Conjecture. Following \cite{bcs}, the starting statement here is:

\begin{proposition}
We have the following estimate,
$$\sup_{U\in O_N}||U||_1\leq N\sqrt{N}$$
with equality if and only if there exists an Hadamard matrix of order $N$.
\end{proposition}

\begin{proof}
This follows indeed from the inequality $||U||_1\leq N\sqrt{N}$, with equality in the rescaled Hadamard matrix case, $U=H/\sqrt{N}$, from Theorem 2.5.
\end{proof}

We begin our study with the following observation:

\begin{proposition}
If the Hadamard Conjecture holds, then
$$\sup_{U\in O_N}||U||_1\geq (N-4.5)\sqrt{N}$$
for any $N\in\mathbb N$.
\end{proposition}

\begin{proof}
If $N$ is a multiple of $4$ we can use an Hadamard matrix, and we are done. In general, we can write $N=M+k$ with $4|M$ and $0\leq k\leq 3$, and use an Hadamard matrix of order $N$, completed with an identity matrix of order $k$. This gives:
\begin{eqnarray*}
\sup_{U\in O_N}||U||_1
&\geq&
M\sqrt{M}+k\\
&\geq&(N-3)\sqrt{N-3}+3\\
&\geq&(N-4.5)\sqrt{N}+3
\end{eqnarray*}

Here the last inequality, which is something proved by taking squares, is valid for any $N\geq 5$. Thus, we are led to the conclusion in the statement.
\end{proof}

We would like to understand now which estimates on the quantity in Proposition 2.12 imply the Hadamard conjecture. We first have the following result:

\begin{proposition}
For any norm one vector $U\in\mathbb R^N$ we have the formula
$$||U||_1=\sqrt{N}\left(1-\frac{||U-H||^2}{2}\right)$$
where $H\in\mathbb R^N$ is the vector given by:
$$H_i=\frac{{\rm sgn}(U_i)}{\sqrt{N}}$$  
\end{proposition}

\begin{proof}
We indeed have the following computation:
\begin{eqnarray*}
||U-H||^2
&=&\sum_i\left(U_i-\frac{{\rm sgn}(U_i)}{\sqrt{N}}\right)^2\\
&=&\sum_iU_i^2-\frac{2|U_i|}{\sqrt{N}}+\frac{1}{N}\\
&=&||U||^2-\frac{2||U||_1}{\sqrt{N}}+1\\
&=&2-\frac{2||U||_1}{\sqrt{N}}
\end{eqnarray*}

But this gives the formula in the statement.
\end{proof}

Next, we have the following estimate, also from \cite{bcs}:

\begin{proposition}
Let $N$ be even, and let $U\in O_N$ be a matrix such that 
$$H=\frac{S}{\sqrt{N}}$$
is not Hadamard, where $S_{ij}={\rm sgn}(U_{ij})$. We have then the following estimate:
$$||U||_1\leq N\sqrt{N}-\frac{1}{N\sqrt{N}}$$
\end{proposition}

\begin{proof}
Since $H$ is not Hadamard, this matrix has two distinct rows $H_1,H_2$ which are not orthogonal. Since $N$ is even, we must have:
$$|<H_1,H_2>|\geq\frac{2}{N}$$

We obtain from this the following estimate:
\begin{eqnarray*}
||U_1-H_1||+||U_2-H_2||
&\geq&|<U_1-H_1,H_2>|+|<U_2-H_2,U_1>|\\
&\geq&|<U_1-H_1,H_2>+<U_2-H_2,U_1>|\\
&=&|<U_2,U_1>-<H_1,H_2>|\\
&=&|<H_1,H_2>|\\
&\geq&\frac{2}{N}
\end{eqnarray*}

Now by applying the estimate in Proposition 2.13 to $U_1,U_2$, we obtain:
\begin{eqnarray*}
||U_1||_1+||U_2||_1
&=&\sqrt{N}\left(2-\frac{||U_1-H_1||^2+||U_2-H_2||^2}{2}\right)\\
&\leq&\sqrt{N}\left(2-\left(\frac{||U_1-H_1||+||U_2-H_2||}{2}\right)^2\right)\\
&\leq&\sqrt{N}\left(2-\frac{1}{N^2}\right)\\
&=&2\sqrt{N}-\frac{1}{N\sqrt{N}}
\end{eqnarray*}

By adding to this inequality the 1-norms of the remaining $N-2$ rows, all bounded from above by $\sqrt{N}$, we obtain the result.  
\end{proof}

We can now answer the question raised above, as follows:

\index{HC}
\index{Hadamard conjecture}

\begin{theorem}
If $N$ is even and the following holds,
$$\sup_{U\in O_N}||U||_1\geq N\sqrt{N}-\frac{1}{N\sqrt{N}}$$
then the Hadamard Conjecture holds at $N$.
\end{theorem}

\begin{proof}
Indeed, if the Hadamard conjecture does not hold at $N$, then the assumption of Proposition 2.14 is satisfied for any $U\in O_N$, and this gives the result.
\end{proof}

As a related result now, also from \cite{bcs}, let us compute the average of the 1-norm on $O_N$. For this purpose, we will use the following well-known result:

\begin{proposition}
We have the following formulae,
$$\int_0^{\pi/2}\cos^pt\,dt=\int_0^{\pi/2}\sin^pt\,dt=\left(\frac{\pi}{2}\right)^{\varepsilon(p)}\frac{p!!}{(p+1)!!}$$
where $\varepsilon(p)=1$ if $p$ is even, and $\varepsilon(p)=0$ if $p$ is odd, and where
$$m!!=(m-1)(m-3)(m-5)\ldots$$
with the product ending at $2$ if $m$ is odd, and ending at $1$ if $m$ is even.
\end{proposition}

\begin{proof}
Let us first compute the integral on the left in the statement:
$$I_p=\int_0^{\pi/2}\cos^pt\,dt$$

We do this by partial integration. We have the following formula:
\begin{eqnarray*}
(\cos^pt\sin t)'
&=&p\cos^{p-1}t(-\sin t)\sin t+\cos^pt\cos t\\
&=&p\cos^{p+1}t-p\cos^{p-1}t+\cos^{p+1}t\\
&=&(p+1)\cos^{p+1}t-p\cos^{p-1}t
\end{eqnarray*}

By integrating between $0$ and $\pi/2$, we obtain the following formula:
$$(p+1)I_{p+1}=pI_{p-1}$$

But this gives the first formula in the statement. As for the second formula, regarding $\sin t$, this follows from the first formula, with the change of variables $t=\pi/2-s$.
\end{proof}

More generally, we have the following result, which is well-known as well:

\begin{proposition}
We have the following formula,
$$\int_0^{\pi/2}\cos^pt\sin^qt\,dt=\left(\frac{\pi}{2}\right)^{\varepsilon(p)\varepsilon(q)}\frac{p!!q!!}{(p+q+1)!!}$$
where $\varepsilon(p)=1$ if $p$ is even, and $\varepsilon(p)=0$ if $p$ is odd, as before.
\end{proposition}

\begin{proof}
Let $I_{pq}$ be the integral in the statement. Observe that we have:
\begin{eqnarray*}
(\cos^pt\sin^qt)'
&=&p\cos^{p-1}t(-\sin t)\sin^qt
+\cos^pt\cdot q\sin^{q-1}t\cos t\\
&=&-p\cos^{p-1}t\sin^{q+1}t+q\cos^{p+1}t\sin^{q-1}t
\end{eqnarray*}

By integrating between $0$ and $\pi/2$, we obtain, for $p,q>0$:
$$pI_{p-1,q+1}=qI_{p+1,q-1}$$

Thus, we can compute $I_{pq}$ by recurrence, and we obtain the above formula.
\end{proof}

Even more generally now, we have the following result, in $N$ dimensions:

\index{spherical integral}
\index{double factorial}

\begin{theorem}
For any exponents $k_1,\ldots,k_N\in\mathbb N$ we have
$$\int_{S^{N-1}}\left|x_1^{k_1}\ldots x_N^{k_N}\right|dx=\left(\frac{2}{\pi}\right)^{\Sigma(k_1,\ldots,k_N)}\frac{(N-1)!!k_1!!\ldots k_N!!}{(N+\Sigma k_i-1)!!}$$
with $\Sigma=[odds/2]$ if $N$ is odd and $\Sigma=[(odds+1)/2]$ if $N$ is even, where ``odds'' denotes the number of odd numbers in the sequence $k_1,\ldots,k_N$.
\end{theorem}

\begin{proof}
We use spherical coordinates, which are by definition as follows:
$$\begin{cases}
x_1\!\!\!&=\ r\cos t_1\\
x_2\!\!\!&=\ r\sin t_1\cos t_2\\
\vdots\\
x_{N-1}\!\!\!&=\ r\sin t_1\sin t_2\ldots\sin t_{N-2}\cos t_{N-1}\\
x_N\!\!\!&=\ r\sin t_1\sin t_2\ldots\sin t_{N-2}\sin t_{N-1}
\end{cases}$$

The corresponding Jacobian $J_N$ can be computed by developing the corresponding determinant over the last column, which gives the following formula:
\begin{eqnarray*}
J_N
&=&r\sin t_1\ldots\sin t_{N-2}\sin t_{N-1}\times \sin t_{N-1}J_{N-1}\\
&+&r\sin t_1\ldots \sin t_{N-2}\cos t_{N-1}\times\cos t_{N-1}J_{N-1}\\
&=&r\sin t_1\ldots\sin t_{N-2}(\sin^2 t_{N-1}+\cos^2 t_{N-1})J_{N-1}\\
&=&r\sin t_1\ldots\sin t_{N-2}J_{N-1}
\end{eqnarray*}

Thus, we obtain by recurrence the following formula:
$$J_N=r^{N-1}\sin^{N-2}t_1\sin^{N-3}t_2\,\ldots\,\sin^2t_{N-3}\sin t_{N-2}$$

With this in hand, the integral in the statement can be written in spherical coordinates, as follows, where $A$ is the area of the sphere, $J_N$ is the Jacobian, and the $2^N$ factor comes from the restriction to the $1/2^N$ part of the sphere where all coordinates are positive:
$$I=\frac{2^N}{A}\int_0^{\pi/2}\ldots\int_0^{\pi/2}x_1^{k_1}\ldots x_N^{k_N}J_N\,dt_1\ldots dt_{N-1}$$

The normalization constant in front of the integral is:
$$\frac{2^N}{A}=\left(\frac{2}{\pi}\right)^{[N/2]}(N-1)!!$$

As for the unnormalized integral, this is given by the following formula:
\begin{eqnarray*}
I'=\int_0^{\pi/2}\ldots\int_0^{\pi/2}
&&(\cos t_1)^{k_1}\\
&&(\sin t_1\cos t_2)^{k_2}\\
&&\vdots\\
&&(\sin t_1\sin t_2\ldots\sin t_{N-2}\cos t_{N-1})^{k_{N-1}}\\
&&(\sin t_1\sin t_2\ldots\sin t_{N-2}\sin t_{N-1})^{k_N}\\
&&\sin^{N-2}t_1\sin^{N-3}t_2\ldots\sin^2t_{N-3}\sin t_{N-2}\\
&&dt_1\ldots dt_{N-1}
\end{eqnarray*}

By rearranging the terms in the above product, we obtain:
\begin{eqnarray*}
I'
&=&\int_0^{\pi/2}\cos^{k_1}t_1\sin^{k_2+\ldots+k_N+N-2}t_1\,dt_1\\
&&\int_0^{\pi/2}\cos^{k_2}t_2\sin^{k_3+\ldots+k_N+N-3}t_2\,dt_2\\
&&\vdots\\
&&\int_0^{\pi/2}\cos^{k_{N-2}}t_{N-2}\sin^{k_{N-1}+k_N+1}t_{N-2}\,dt_{N-2}\\
&&\int_0^{\pi/2}\cos^{k_{N-1}}t_{N-1}\sin^{k_N}t_{N-1}\,dt_{N-1}
\end{eqnarray*}

Now by using the $N=2$ integration formula from Proposition 2.17, we obtain:
\begin{eqnarray*}
I'
&=&\frac{\pi}{2}\cdot\frac{k_1!!(k_2+\ldots+k_N+N-2)!!}{(k_1+\ldots+k_N+N-1)!!}\left(\frac{2}{\pi}\right)^{\delta(k_1,k_2+\ldots+k_N+N-2)}\\
&&\frac{\pi}{2}\cdot\frac{k_2!!(k_3+\ldots+k_N+N-3)!!}{(k_2+\ldots+k_N+N-2)!!}\left(\frac{2}{\pi}\right)^{\delta(k_2,k_3+\ldots+k_N+N-3)}\\
&&\vdots\\
&&\frac{\pi}{2}\cdot\frac{k_{N-2}!!(k_{N-1}+k_N+1)!!}{(k_{N-2}+k_{N-1}+k_N+2)!!}\left(\frac{2}{\pi}\right)^{\delta(k_{N-2},k_{N-1}+k_N+1)}\\
&&\frac{\pi}{2}\cdot\frac{k_{N-1}!!k_N!!}{(k_{N-1}+k_N+1)!!}\left(\frac{2}{\pi}\right)^{\delta(k_{N-1},k_N)}
\end{eqnarray*}

In order to compute this quantity, let us denote by $F$ the part involving the double factorials, and by $P$ the part involving the powers of $\pi/2$, so that we have:
$$I'=F\cdot P$$

Regarding $F$, there are many cancellations there, and we end up with:
$$F=\frac{k_1!!\ldots k_N!!}{(\Sigma k_i+N-1)!!}$$

As in what regards $P$,  the $\delta$ exponents on the right sum up to the following number:
$$\Delta(k_1,\ldots,k_N)=\sum_{i=1}^{N-1}\delta(k_i,k_{i+1}+\ldots+k_N+N-i-1)$$

In other words, with this notation, the above formula reads:
\begin{eqnarray*}
I'
&=&\left(\frac{\pi}{2}\right)^{N-1}\frac{k_1!!k_2!!\ldots k_N!!}{(k_1+\ldots+k_N+N-1)!!}\left(\frac{2}{\pi}\right)^{\Delta(k_1,\ldots,k_N)}\\
&=&\left(\frac{2}{\pi}\right)^{\Delta(k_1,\ldots,k_N)-N+1}\frac{k_1!!k_2!!\ldots k_N!!}{(k_1+\ldots+k_N+N-1)!!}\\
&=&\left(\frac{2}{\pi}\right)^{\Sigma(k_1,\ldots,k_N)-[N/2]}\frac{k_1!!k_2!!\ldots k_N!!}{(k_1+\ldots+k_N+N-1)!!}
\end{eqnarray*}

Here the formula relating $\Delta$ to $\Sigma$ follows from a number of simple observations, the first of which is the following one: due to obvious parity reasons, the sequence of $\delta$ numbers appearing in the definition of $\Delta$ cannot contain two consecutive zeroes. Thus, we have $I'$, and together with $I=(2^N/V)I'$, this gives the formula in the statement.
\end{proof}

As a technical observation, the exponent $\Sigma$ appearing in the statement of Theorem 2.18 can be written as well in the following more compact form:
$$\Sigma(k_1,\ldots,k_p)=\left[\frac{N+odds+1}{2}\right]-\left[\frac{N+1}{2}\right]$$

However, for concrete applications, the writing in Theorem 2.18 is more convenient. Now by using this result, we obtain the following estimate, from \cite{bcs}:

\index{average of 1-norm}

\begin{theorem}
We have the following estimate,
$$\int_{O_N}||U||_1\,dU\simeq\sqrt{\frac{2}{\pi}}\cdot N\sqrt{N}$$
valid in the $N\to\infty$ limit.
\end{theorem}

\begin{proof}
We use the well-known fact that the row slices of $O_N$ are all isomorphic to the sphere $S^{N-1}$, with the restriction of the Haar measure of $O_N$ corresponding in this way to the uniform measure on $S^{N-1}$. Together with a standard symmetry argument, this shows that the average of the 1-norm on $O_N$ is given by:
\begin{eqnarray*}
\int_{O_N}||U||_1\,dU
&=&\sum_{ij}\int_{O_N}|U_{ij}|\,dU\\
&=&N^2\int_{O_N}|U_{11}|\,dU\\
&=&N^2\int_{S^{N-1}}|x_1|\,dx
\end{eqnarray*}

We denote by $I$ the integral on the right. According to Theorem 2.18, we have:
\begin{eqnarray*}
I
&=&\left(\frac{2}{\pi}\right)^{\Sigma(1)}\frac{(N-1)!!}{N!!}\\
&=&\begin{cases}
\displaystyle{\frac{2}{\pi}\cdot\frac{2.4.6\ldots(N-2)}{3.5.7\ldots(N-1)}}& (N\ {\rm even})\\
\displaystyle{1\cdot\frac{3.5.7\ldots(N-2)}{2.4.6\ldots(N-1)}}& (N\ {\rm odd})
\end{cases}\\
&=&
\begin{cases}
\displaystyle{\frac{4^M}{\pi M}}\begin{pmatrix}2M\\ M\end{pmatrix}^{-1}& (N=2M)\\
4^{-M}\begin{pmatrix}2M\\ M\end{pmatrix}& (N=2M+1)
\end{cases}
\end{eqnarray*}

Now by using the Stirling formula, we get from this:
\begin{eqnarray*}
I
&\simeq&\begin{cases}
\displaystyle{\frac{4^M}{\pi M}\cdot\frac{\sqrt{\pi M}}{4^M}}& (N=2M)\\
\displaystyle{4^{-M}\cdot\frac{4^M}{\sqrt{\pi M}}}& (N=2M+1)
\end{cases}\\
&=&\begin{cases}
\displaystyle{\frac{1}{\sqrt{\pi M}}}& (N=2M)\\
\displaystyle{\frac{1}{\sqrt{\pi M}}}& (N=2M+1)
\end{cases}\\
&\simeq&\sqrt{\frac{2}{\pi N}}
\end{eqnarray*}

Thus, we are led to the conclusion in the statement.
\end{proof}

The above result gives in particular the following estimate, in the $N\to\infty$ limit:
$$\sup_{U\in O_N}||U||_1\,dU\simeq\sqrt{\frac{2}{\pi}}\cdot N\sqrt{N}$$

For better estimates, the problem is to compute the higher moments of the 1-norm:
$$I_k=\int_{O_N}||U||_1^k\,dU$$

Indeed, the supremum that we are interested in is given by the following formula:
$$\sup_{U\in O_N}||U||_1\,dU=\lim_{k\to\infty}I_k^{1/k}$$

However, the computation of the integrals $I_k$ is a difficult problem, and no concrete applications to the Hadamard Conjecture have been found so far. See \cite{bcs}.

\section*{2c. Bistochastic matrices}

Let us discuss now a third analytic topic. The motivation here comes from the fact that the bistochastic Hadamard matrices look better than their non-bistochastic counterparts. As an illustration, $F_2$ looks better in complex bistochastic form:
$$\begin{pmatrix}
1&1\\
1&-1
\end{pmatrix}\sim
\begin{pmatrix}
i&1\\
1&i
\end{pmatrix}$$

\index{Walsh matrix}

Also, the matrix $W_4$ looks better in its bistochastic form, which is the matrix $K_4$:
$$\begin{pmatrix}
1&1&1&1\\
1&-1&1&-1\\
1&1&-1&-1\\
1&-1&-1&1
\end{pmatrix}\sim\begin{pmatrix}
-1&1&1&1\\
1&-1&1&1\\
1&1&-1&1\\
1&1&1&-1
\end{pmatrix}$$

We have the following algebraic result on the subject, which shows in particular that we cannot put any Hadamard matrix in bistochastic form:

\index{bictochastic}
\index{row stochastic}
\index{column stochastic}

\begin{theorem}
For an Hadamard matrix $H\in M_N(\mathbb C)$, the following are equivalent:
\begin{enumerate}
\item $H$ is bistochastic, with sums $\lambda$.

\item $H$ is row-stochastic, with sums $\lambda$, and $\lambda^2=N$.
\end{enumerate}
In particular, is such a matrix exists, then $N\in 4\mathbb N$ must be a square.
\end{theorem}

\begin{proof}
Both the implications are elementary, as follows:

\medskip

$(1)\implies(2)$ If we denote by $H_1,\ldots,H_N\in(\pm1)^N$ the rows of $H$, we have indeed:
\begin{eqnarray*}
N&=&\sum_i<H_1,H_i>\\
&=&\sum_jH_{1j}\sum_iH_{ij}\\
&=&\sum_jH_{1j}\cdot\lambda\\
&=&\lambda^2
\end{eqnarray*}

$(2)\implies(1)$ Consider the all-one vector $\xi=(1)_i\in\mathbb R^N$. The fact that $H$ is row-stochastic with sums $\lambda$ reads:
\begin{eqnarray*}
\sum_jH_{ij}=\lambda,\forall i
&\iff&\sum_jH_{ij}\xi_j=\lambda\xi_i,\forall i\\
&\iff&H\xi=\lambda\xi
\end{eqnarray*}

Also, the fact that $H$ is column-stochastic with sums $\lambda$ reads:
\begin{eqnarray*}
\sum_iH_{ij}=\lambda,\forall j
&\iff&\sum_jH_{ij}\xi_i=\lambda\xi_j,\forall j\\
&\iff&H^t\xi=\lambda\xi
\end{eqnarray*}

We must prove that the first condition implies the second one, provided that the row sum $\lambda$ satisfies $\lambda^2=N$. But this follows from the following computation:
\begin{eqnarray*}
H\xi=\lambda\xi
&\implies&H^tH\xi=\lambda H^t\xi\\
&\implies&N\xi=\lambda H^t\xi\\
&\implies&H^t\xi=\lambda\xi
\end{eqnarray*}

Thus, we have proved both the implications, and we are done.
\end{proof}

In practice now, the even Walsh matrices, having size $N=4^n$, which is a square as required above, can be put in bistochastic form, as follows:
$$W_{4^n}\sim K_4^{\otimes n}$$

As for the odd Walsh matrices, having size $N=2\times 4^n$, these cannot be put in bistochastic form. However, we can do this over the complex numbers, with the equivalence being as follows at $N=2$, and then by tensoring with $K_4^{\otimes n}$ in general:
$$\begin{pmatrix}1&1\\ 1&-1\end{pmatrix}\sim
\begin{pmatrix}i&1\\1&i\end{pmatrix}$$

This is quite interesting, and in general now, it is known from Idel-Wolf \cite{iwo} that any complex Hadamard matrix can be put in bistochastic form, by a certain non-explicit method. Thus, we have here some theory to be developed. We will be back to this.

\bigskip

There is as well an analytic approach to these questions, based on:

\index{excess}

\begin{theorem}
For an Hadamard matrix $H\in M_N(\pm1)$, the excess, 
$$E(H)=\sum_{ij}H_{ij}$$
satisfies $|E(H)|\leq N\sqrt{N}$, with equality if and only if $H$ is bistochastic.
\end{theorem}

\begin{proof}
In terms of the all-one vector $\xi=(1)_i\in\mathbb R^N$, we have:
\begin{eqnarray*}
E(H)
&=&\sum_{ij}H_{ij}\\
&=&\sum_{ij}H_{ij}\xi_j\xi_i\\
&=&\sum_i(H\xi)_i\xi_i\\
&=&<H\xi,\xi>
\end{eqnarray*}

Now by using the Cauchy-Schwarz inequality, along with the fact that $U=H/\sqrt{N}$ is orthogonal, and hence of norm 1, we obtain, as claimed:
\begin{eqnarray*}
|E(H)|
&\leq&||H\xi||\cdot||\xi||\\
&\leq&||H||\cdot||\xi||^2\\
&=&N\sqrt{N}
\end{eqnarray*}

Regarding now the equality case, this requires the vectors $H\xi,\xi$ to be proportional, and so our matrix $H$ to be row-stochastic. But since $U=H/\sqrt{N}$ is orthogonal, we have: 
$$H\xi\sim\xi\iff H^t\xi\sim\xi$$

Thus our matrix $H$ must be bistochastic, as claimed.
\end{proof}

\section*{2d. The glow}

One interesting question, that we will discuss now, is that of computing the law of the excess over the equivalence class of $H$. Let us start with the following definition:

\index{glow}

\begin{definition}
The glow of $H\in M_N(\pm1)$ is the distribution of the excess,
$$E=\sum_{ij}H_{ij}$$
over the Hadamard equivalence class of $H$.
\end{definition}

Since the excess is invariant under permutations of rows and columns, we can restrict the attention to the matrices $\widetilde{H}\simeq H$ obtained by switching signs on rows and columns. More precisely, let $(a,b)\in\mathbb Z_2^N\times\mathbb Z_2^N$, and consider the following matrix:
$$\widetilde{H}_{ij}=a_ib_jH_{ij}$$

We can regard the sum of entries of $\widetilde{H}$ as a random variable, over the group $\mathbb Z_2^N\times\mathbb Z_2^N$, and we have the following equivalent description of the glow:

\index{Hadamard equivalence}

\begin{proposition}
Given a matrix $H\in M_N(\pm 1)$, if we define $\varphi:\mathbb Z_2^N\times\mathbb Z_2^N\to\mathbb Z$ as the excess of the corresponding Hadamard equivalent of $H$, 
$$\varphi(a,b)=\sum_{ij}a_ib_jH_{ij}$$
then the glow is the probability measure on $\mathbb Z$ given by $\mu(\{k\})=P(\varphi=k)$.
\end{proposition}

\begin{proof}
The function $\varphi$ in the statement can indeed be regarded as a random variable over the group $\mathbb Z_2^N\times\mathbb Z_2^N$, with this latter group being endowed with its uniform probability measure $P$. The distribution $\mu$ of this variable $\varphi$ is then given by:
$$\mu(\{k\})=\frac{1}{4^N}\#\left\{(a,b)\in \mathbb Z_2^N\times\mathbb Z_2^N\Big|\varphi(a,b)=k\right\}$$

By the above discussion, this distribution is exactly the glow.
\end{proof}

The terminology in Definition 2.22 comes from the following picture. Assume that we have a square city, with $N$ horizontal streets and $N$ vertical streets, and with street lights at each crossroads. When evening comes the lights are switched on at the positions $(i,j)$ where $H_{ij}=1$, and then, all night long, they are randomly switched on and off, with the help of $2N$ master switches, one at the end of each street:
$$\begin{matrix}
\to&&\diamondsuit&\diamondsuit&\diamondsuit&\diamondsuit\\
\to&&\diamondsuit&\times&\diamondsuit&\times\\
\to&&\diamondsuit&\diamondsuit&\times&\times\\
\to&&\diamondsuit&\times&\times&\diamondsuit\\
\\
&&\uparrow&\uparrow&\uparrow&\uparrow
\end{matrix}$$ 

With this picture in mind, $\mu$ describes indeed the glow of the city. At a more advanced level now, all this is related to the Gale-Berlekamp game, and this is where our main motivation for studying the glow comes from. We refer to Fishburn-Sloane \cite{fsl} and Roth-Viswanathan \cite{rvi} for details on the Gale-Berlekamp game.

\index{Gale-Berlekamp game}

\bigskip

In order to compute the glow, it is useful to have in mind the following picture:
$$\begin{matrix}
&&b_1&\ldots&b_N\\
&&\downarrow&&\downarrow\\
(a_1)&\to&H_{11}&\ldots&H_{1N}&\Rightarrow&S_1\\
\vdots&&\vdots&&\vdots&&\vdots\\
(a_N)&\to&H_{N1}&\ldots&H_{NN}&\Rightarrow&S_N
\end{matrix}$$

Here the columns of $H$ have been multiplied by the entries of the horizontal switching vector $b$, the resulting sums on rows are denoted $S_1,\ldots,S_N$, and the vertical switching vector $a$ still has to act on these sums, and produce the glow component at $b$. 

\bigskip

With this picture in mind, we first have the following result:

\index{Bernoulli law}

\begin{proposition}
The glow of a matrix $H\in M_N(\pm 1)$ is given by
$$\mu=\frac{1}{2^N}\sum_{b\in\mathbb Z_2^N}\beta_1(c_1)*\ldots*\beta_N(c_N)$$
where the measures on the right are convolution powers of Bernoulli laws,
$$\beta_r(c)=\left(\frac{\delta_r+\delta_{-r}}{2}\right)^{*c}$$
and where $c_r=\#\{r\in|S_1|,\ldots,|S_N|\}$, with $S=Hb$.
\end{proposition}

\begin{proof}
We use the interpretation of the glow explained above. So, consider the decomposition of the glow over $b$ components:
$$\mu=\frac{1}{2^N}\sum_{b\in\mathbb Z_2^N}\mu_b$$

With the notation $S=Hb$, as in the statement, the numbers $S_1,\ldots,S_N$ are the row sums of $\widetilde{H}_{ij}=H_{ij}a_ib_j$. Thus the glow components are given by:
$$\mu_b=law\left(\pm S_1\pm S_2\ldots\pm S_N\right)$$

By permuting now the sums on the right, we have the following formula:
$$\mu_b=law\big(\underbrace{\pm 0\ldots\pm 0}_{c_0}\ \underbrace{\pm 1\ldots\pm 1}_{c_1}\,\ldots\ldots\,\underbrace{\pm N\ldots\pm N}_{c_N}\big)$$

Now since the $\pm$ variables each follow a Bernoulli law, and these Bernoulli laws are independent, we obtain a convolution product as in the statement.
\end{proof}

We will need the following elementary fact:

\begin{proposition}
Let $H\in M_N(\pm1)$ be an Hadamard matrix of order $N\geq 4$. 
\begin{enumerate}
\item The sums of entries on rows $S_1,\ldots,S_N$ are even, and equal modulo $4$.

\item If the sums on the rows $S_1,\ldots,S_N$ are all $0$ modulo $4$, then the number of rows whose sum is $4$ modulo $8$ is odd for $N=4(8)$, and even for $N=0(8)$.
\end{enumerate}
\end{proposition}

\begin{proof}
This is something elementary, the proof being as follows:

\medskip

(1) Let us pick two rows of our matrix, and then permute the columns such that these two rows look as follows:
$$\begin{pmatrix}
1\ldots\ldots1&1\ldots\ldots1&-1\ldots-1&-1\ldots-1\\
\underbrace{1\ldots\ldots1}_a&\underbrace{-1\ldots-1}_b&\underbrace{1\ldots\ldots1}_c&\underbrace{-1\ldots-1}_d
\end{pmatrix}$$

We have $a+b+c+d=N$, and by orthogonality we obtain $a+d=b+c$. Thus $a+d=b+c=N/2$, and since $N/2$ is even we have $b=c(2)$, which gives the result.

\medskip

(2) In the case where $H$ is ``row-dephased'', in the sense that its first row consists of $1$ entries only, the row sums are $N,0,\ldots,0$, and so the result holds. In general now, by permuting the columns we can assume that our matrix looks as follows:
$$H=\begin{pmatrix}1\ldots\ldots1&-1\ldots-1\\ \underbrace{\vdots}_x&\underbrace{\ \vdots\ }_y\end{pmatrix}$$

We have $x+y=N=0(4)$, and since the first row sum $S_1=x-y$ is by assumption 0 modulo 4, we conclude that $x,y$ are even. In particular, since $y$ is even, the passage from $H$ to its row-dephased version $\widetilde{H}$ can be done via $y/2$ double sign switches. Now, in view of the above, it is enough to prove that the conclusion in the statement is stable under a double sign switch. So, let $H\in M_N(\pm1)$ be Hadamard, and let us perform to it a double sign switch, say on the first two columns. Depending on the values of the entries on these first two columns, the total sums on the rows change as follows:
\begin{eqnarray*}
\begin{pmatrix}+&+&\ldots&\ldots\end{pmatrix}&:&S\to S-4\\
\begin{pmatrix}+&-&\ldots&\ldots\end{pmatrix}&:&S\to S\\
\begin{pmatrix}-&+&\ldots&\ldots\end{pmatrix}&:&S\to S\\
\begin{pmatrix}-&-&\ldots&\ldots\end{pmatrix}&:&S\to S+4
\end{eqnarray*}

We can see that the changes modulo 8 of the row sum $S$ occur precisely in the first and in the fourth case. But, since the first two columns of our matrix $H\in M_N(\pm1)$ are orthogonal, the total number of these cases is even, and this finishes the proof.
\end{proof}

Observe that Proposition 2.24 and Proposition 2.25 (1) show that the glow of an Hadamard matrix of order $N\geq 4$ is supported by $4\mathbb Z$. With this in hand, we have:

\index{glow components}

\begin{theorem}
Let $H\in M_N(\pm1)$ be an Hadamard matrix of order $N\geq 4$, and denote by $\mu^{even},\mu^{odd}$ the mass one-rescaled restrictions of $\mu\in\mathcal P(4\mathbb Z)$ to $8\mathbb Z,8\mathbb Z+4$.
\begin{enumerate}
\item At $N=0(8)$ we have $\mu=\frac{3}{4}\mu^{even}+\frac{1}{4}\mu^{odd}$.

\item At $N=4(8)$ we have $\mu=\frac{1}{4}\mu^{even}+\frac{3}{4}\mu^{odd}$.
\end{enumerate}
\end{theorem}

\begin{proof}
We use the glow decomposition over $b$ components, from Proposition 2.24:
$$\mu=\frac{1}{2^N}\sum_{b\in\mathbb Z_2^N}\mu_b$$

The idea is that the decomposition formula in the statement will occur over averages of the following type, over truncated sign vectors $c\in\mathbb Z_2^{N-1}$:
$$\mu'_c=\frac{1}{2}(\mu_{+c}+\mu_{-c})$$

Indeed, we know from Proposition 2.25 (1) that modulo 4, the sums on rows are either $0,\ldots,0$ or $2,\ldots,2$. Now since these two cases are complementary when pairing switch vectors $(+c,-c)$, we can assume that we are in the case $0,\ldots,0$ modulo 4. Now by looking at this sequence modulo 8, and letting $x$ be the number of 4 components, so that the number of 0 components is $N-x$, we have:
$$\frac{1}{2}(\mu_{+c}+\mu_{-c})=\frac{1}{2}\left(law(\underbrace{\pm0\ldots\pm 0}_{N-x}\,\underbrace{\pm4\ldots\pm 4}_x)+law(\underbrace{\pm 2\ldots\pm 2}_N)\right)$$

Now by using Proposition 2.25 (2), the first summand splits $1-0$ or $0-1$ on $8\mathbb Z,8\mathbb Z+4$, depending on the class of $N$ modulo 8. As for the second summand, since $N$ is even this always splits $\frac{1}{2}-\frac{1}{2}$ on $8\mathbb Z,8\mathbb Z+4$. Thus, by making the average we obtain either a $\frac{3}{4}-\frac{1}{4}$ or a $\frac{1}{4}-\frac{3}{4}$ splitting on $8\mathbb Z,8\mathbb Z+4$, depending on the class of $N$ modulo 8, as claimed.
\end{proof}

Various computer simulations suggest that the above measures $\mu^{even},\mu^{odd}$ don't have further general properties, so that the basic algebraic theory stops here. However, analytically speaking now, we have an interesting result about the glow. We will need: 

\begin{proposition}
The moments of the normal law 
$$g_1=\frac{1}{\sqrt{2\pi}}e^{-x^2/2}dx$$
are the numbers $M_k=k!!$, with the convention $k!!=0$ when $k$ is odd. 
\end{proposition}

\begin{proof}
We have indeed the following computation:
\begin{eqnarray*}
M_k
&=&\frac{1}{\sqrt{2\pi}}\int_\mathbb Rx^ke^{-x^2/2}dx\\
&=&\frac{1}{\sqrt{2\pi}}\int_\mathbb R(x^{k-1})\left(-e^{-x^2/2}\right)'dx\\
&=&\frac{1}{\sqrt{2\pi}}\int_\mathbb R(k-1)x^{k-2}e^{-x^2/2}dx\\
&=&(k-1)\times\frac{1}{\sqrt{2\pi}}\int_\mathbb Rx^{k-2}e^{-x^2/2}dx\\
&=&(k-1)M_{k-2}
\end{eqnarray*}

On the other hand, we have $M_0=1$, $M_1=0$. Thus by recurrence, the even moments vanish, and the odd moments are given by the formula in the statement.
\end{proof}

We can now formulate our analytic result regarding the glow, as follows:

\index{normal variable}
\index{Gaussian variable}
\index{partition}

\begin{theorem}
The glow moments of $H\in M_N(\pm1)$ are given by:
$$\int_{\mathbb Z_2^N\times\mathbb Z_2^N}\left(\frac{E}{N}\right)^{2p}=(2p)!!+O(N^{-1})$$
In particular the normalized variable $F=E/N$ becomes Gaussian with $N\to\infty$.
\end{theorem}

\begin{proof}
Consider the variable in the statement, written as before, as a function of two vectors $a,b$, belonging to the group $\mathbb Z_2^N\times\mathbb Z_2^N$:
$$E=\sum_{ij}a_ib_jH_{ij}$$

Let $P_{even}(r)\subset P(r)$ be the set of partitions of $\{1,\ldots,r\}$ having all blocks of even size. The moments of $E$ are then given by:
\begin{eqnarray*}
\int_{\mathbb Z_2^N\times\mathbb Z_2^N}E^r
&=&\int_{\mathbb Z_2^N\times\mathbb Z_2^N}\sum_{ix}a_{i_1}\ldots a_{i_r}b_{x_1}\ldots b_{x_r}H_{i_1x_1}\ldots H_{i_rx_r}\\
&=&\sum_{ix}H_{i_1x_1}\ldots H_{i_rx_r}\int_{\mathbb Z_2^N}a_{i_1}\ldots a_{i_r}\int_{\mathbb Z_2^N}b_{x_1}\ldots b_{x_r}\\
&=&\sum_{\pi,\sigma\in P_{even}(r)}\sum_{\ker i=\pi,\ker x=\sigma}H_{i_1x_1}\ldots H_{i_rx_r}
\end{eqnarray*}

Thus the moments decompose over partitions $\pi\in P_{even}(r)$, with the contributions being obtained by integrating the following quantities:
$$C(\sigma)=\sum_{\ker x=\sigma}\sum_iH_{i_1x_1}\ldots H_{i_rx_r}\cdot a_{i_1}\ldots a_{i_r}$$

Now by M\"obius inversion, we obtain a formula as follows:
$$\int_{\mathbb Z_2^N\times\mathbb Z_2^N}E^r=\sum_{\pi\in P_{even}(r)}K(\pi)N^{|\pi|}I(\pi)$$

To be more precise, here the coefficients on the right are as follows, where $\mu$ is the M\"obius function of $P_{even}(r)$:
$$K(\pi)=\sum_{\sigma\in P_{even}(r)}\mu(\pi,\sigma)$$

As for the contributions on the right, with the convention that $H_1,\ldots,H_N\in\mathbb Z_2^N$ are the rows of our matrix $H$, these are as follows:
$$I(\pi)=\sum_i\prod_{b\in\pi}\frac{1}{N}\left\langle\prod_{r\in b}H_{i_r},1\right\rangle$$

With this formula in hand, the first assertion follows, because the biggest elements of the lattice $P_{even}(2p)$ are the $(2p)!!$ partitions consisting of $p$ copies of a $2$-block:
$$\int_{\mathbb Z_2^N\times\mathbb Z_2^N}\left(\frac{E}{N}\right)^{2p}=(2p)!!+O(N^{-1})$$

As for the second assertion, this follows from the moment formula, and from the fact that the glow of $H\in M_N(\pm1)$ is real, and symmetric with respect to $0$.
\end{proof}

All the above was of course a bit technical, using some familiarity with probability theory, and for an introduction to this, we refer for instance to Durrett \cite{dur}. We will be back to glow computations in chapter 11 below, in the complex setting.

\section*{2e. Exercises} 

We have seen a lot of calculus in the above, and most of our exercises will be about more calculus, precisely. To start with, however, we have:

\begin{exercise}
Briefly discuss how the theory of the determinant can be developed, as a signed volume.
\end{exercise}

This is something that we used in the above, in the proof of the Hadamard determinant bound. Make sure that everything is fine here, with your linear algebra knowledge.

\begin{exercise}
Prove that the following matrix belongs to $\sqrt{N}O_N$,
$$K_N=\frac{1}{\sqrt{N}}\begin{pmatrix}
2-N&&2\\
&\ddots&\\
2&&2-N
\end{pmatrix}$$
and is a critical point of the $1$-norm on $\sqrt{N}O_N$.
\end{exercise}

The first part is normally a standard linear algebra computation. As for the second part, this can only be something which can be done with Lagrange multipliers.

\begin{exercise}
Establish the following integration formula over the sphere,
$$\int_{S^{N-1}}x_1^{k_1}\ldots x_N^{k_N}\,dx=\frac{(N-1)!!k_1!!\ldots k_N!!}{(N+\Sigma k_i-1)!!}$$
by using spherical coordinates and Fubini.
\end{exercise}

Observe that this formula holds in the case where all the exponents $k_i$ are even, because here the quantity to be integrated equals its absolute value, and we have seen in the above how to integrate such absolute values. In general, the proof should be along the same lines as the proof for the formula with absolute values.

\chapter{Norm maximizers}

\section*{3a. Critical points}

We have seen in the previous chapter that the set $Y_N=M_N(\pm1)\cap\sqrt{N}O_N$ formed by the $N\times N$ Hadamard matrices can be located inside $\sqrt{N}O_N$ by using analytic techniques, and more precisely variations of the following result:

\index{norm maximizer}

\begin{theorem}
Given a matrix $H\in\sqrt{N}O_N$ we have:
\begin{enumerate}
\item $||H||_p\leq N^{2/p}$ for $p\in[1,2)$, with equality precisely when $H$ is Hadamard.

\item $||H||_p\geq N^{2/p}$ for $p\in(2,\infty]$, with equality precisely when $H$ is Hadamard.
\end{enumerate}
\end{theorem} 

\begin{proof}
This is something that we know from chapter 2, in rescaled reformulation. Consider indeed the $p$-norm on $\sqrt{N}O_N$, which at $p\in[1,\infty)$ is given by:
$$||H||_p=\left(\sum_{ij}|H_{ij}|^p\right)^{1/p}$$

We have then $||H||_2=N$, and by using this, together with the Jensen inequality for $\psi(x)=x^{p/2}$, or simply the H\"older inequality for the norms, we obtain the results. As for the case $p=\infty$, this follows with $p\to\infty$, or directly via Cauchy-Schwarz.
\end{proof}

Once again following the material in chapter 2, we have seen there that a nice result can be obtained along these lines at $N=3$ and $p=1$. To be more precise, the maximizers of the 1-norm on $\sqrt{3}O_3$ are the following matrix, and its Hadamard conjugates:
$$K_3=\frac{1}{\sqrt{3}}\begin{pmatrix}
-1&2&2\\
2&-1&2\\
2&2&-1
\end{pmatrix}$$

In general, however, computing the maximizers of the $1$-norm on $\sqrt{N}O_N$ remains a difficult question. So, based on the above, let us formulate the following definition:

\index{AHM}
\index{almost Hadamard matrix}
\index{optimal AHM}
\index{optimal almost Hadamard matrix}

\begin{definition}
A matrix $H\in\sqrt{N}O_N$ is called:
\begin{enumerate}
\item Almost Hadamard, if it locally maximizes the $1$-norm on $\sqrt{N}O_N$.

\item Optimal almost Hadamard, if it maximizes the $1$-norm on $\sqrt{N}O_N$.
\end{enumerate}
\end{definition}

More generally, we can talk about $p$-almost Hadamard matrices, exactly in the same way, at any $p\in[1,\infty]-\{2\}$, by using the results in Theorem 3.1. When a matrix $H\in\sqrt{N}O_N$ is almost Hadamard at any $p$, we call it ``absolute almost Hadamard''. We will see in what follows that, while the study of the optimal almost Hadamard matrices remains something quite difficult, in the general almost Hadamard setting there are many interesting things to be done, and some nice theory to be developed.

\bigskip

Needless to say, all this is motivated by the lack of Hadamard matrices at $N>2$, $N\notin4\mathbb N$. However, we will see that our theory is quite interesting even at values $N\in4\mathbb N$. Finally, let us mention that there is a long story with the almost Hadamard matrices, going back to the 2010 paper \cite{bcs}, then to the 2012 paper \cite{bnz}, and with the theory of such matrices having been further developed all over the 10s, in the series of papers \cite{bn1}, \cite{bn3}, \cite{bs1}, \cite{bs2}, \cite{moh}. We will try to explain here the basics of this theory.

\bigskip

In order to get started, let us study the local mazimizers of the 1-norm on $\sqrt{N}O_N$. It is technically convenient here to rescale by $1\sqrt{N}$, and work instead over the orthogonal group $O_N$, by using the avaliable tools here. Following \cite{bcs}, we first have:

\index{local maximizer}
\index{rotation trick}

\begin{theorem}
If $U\in O_N$ locally maximizes the $1$-norm, then 
$$U_{ij}\neq 0$$
must hold for any $i,j$.
\end{theorem}

\begin{proof}
Assume by contradiction that $U$ has a 0 entry. By permuting the rows we can assume that this 0 entry is in the first row, having under it a nonzero entry in the second row. We denote by $U_1,\ldots,U_N$ the rows of $U$. By permuting the columns we can assume that we have a block decomposition of the following type:
$$\begin{pmatrix}U_1\\ U_2\end{pmatrix}
=\begin{pmatrix}
0&0&Y&A&B\\
0&X&0&C&D
\end{pmatrix}$$

Here $X,Y,A,B,C,D$ are certain vectors with nonzero entries, with $A,B,C,D$ chosen such that each entry of $A$ has the same sign as the corresponding entry of $C$, and each entry of $B$ has sign opposite to the sign of the corresponding entry of $D$. Now for $t>0$ small consider the matrix $U^t$ obtained by rotating by an angle $t$ the first two rows of $U$. In row notation, this matrix is given by the following formula: 
$$U^t
=\begin{pmatrix}
\cos t&\sin t\\
-\sin t&\cos t\\
&&1\\
&&&\ddots\\
&&&&1\end{pmatrix}
\begin{pmatrix}
U_1\\ U_2\\ U_3\\ \vdots\\ U_N
\end{pmatrix}
=\begin{pmatrix}
\cos t\cdot U_1+\sin t\cdot U_2\\ -\sin t\cdot U_1+\cos t\cdot U_2\\ U_3\\ \vdots\\ U_N
\end{pmatrix}$$

We make the convention that the lower-case letters denote the 1-norms of the corresponding upper-case vectors. According to the above sign conventions, we have:
\begin{eqnarray*}
||U^t||_1
&=&||\cos t\cdot U_1+\sin t\cdot U_2||_1+||-\sin t\cdot U_1+\cos t\cdot U_2||_1+\sum_{i=3}^Nu_i\\
&=&(\cos t+\sin t)(x+y+b+c)+(\cos t-\sin t)(a+d)+\sum_{i=3}^Nu_i\\
&=&||U||_1+(\cos t+\sin t-1)(x+y+b+c)+(\cos t-\sin t-1)(a+d)
\end{eqnarray*}

By using $\sin t=t+O(t^2)$ and $\cos t=1+O(t^2)$ we obtain:
\begin{eqnarray*}
||U^t||_1
&=&||U||_1+t(x+y+b+c)-t(a+d)+O(t^2)\\
&=&||U||_1+t(x+y+b+c-a-d)+O(t^2)
\end{eqnarray*}

In order to conclude, we have to prove that $U$ cannot be a local maximizer of the $1$-norm. This will basically follow by comparing the norm of $U$ to the norm of $U^t$, with $t>0$ small or $t<0$ big. However, since in the above computation it was technically convenient to assume $t>0$, we actually have three cases:

\medskip

\underline{Case 1}: $b+c>a+d$. Here for $t>0$ small enough the above formula shows that we have $||U^t||_1>||U||_1$, and we are done.

\medskip

\underline{Case 2}: $b+c=a+d$. Here we use the fact that $X$ is not null, which gives $x>0$. Once again for $t>0$ small enough we have $||U^t||_1>||U||_1$, and we are done.

\medskip

\underline{Case 3}: $b+c<a+d$. In this case we can interchange the first two rows of $U$ and restart the whole procedure: we fall in Case 1, and we are done again.
\end{proof}

Let us study now the critical points. It is convenient here to talk about more general $p$-norms, or even more general functions of the quantities $|U_{ij}|$, because this will lead to some interesting combinatorics. Following \cite{bcs}, \cite{bn3}, we have the following result:

\index{critical point}

\begin{theorem}
Consider a differentiable function $\varphi:[0,\infty)\to\mathbb R$. An orthogonal matrix having nonzero entries, $U\in O_N^*$, is then a critical point of the function
$$F(U)=\sum_{ij}\varphi(|U_{ij}|)$$
precisely when the matrix $WU^t$ is symmetric, where:
$$W_{ij}={\rm sgn}(U_{ij})\varphi'(|U_{ij}|)$$
In particular, for $F(U)=||U||_1$ we need $SU^t$ to be symmetric, where $S_{ij}={\rm sgn}(U_{ij})$.
\end{theorem}

\begin{proof}
We regard $O_N$ as a real algebraic manifold, with coordinates $U_{ij}$. This manifold consists by definition of the zeroes of the following polynomials: 
$$A_{ij}=\sum_kU_{ik}U_{jk}-\delta_{ij}$$

Since $O_N$ is smooth, and so is a differential manifold in the usual sense, it follows from the general theory of Lagrange multipliers that a given matrix $U\in O_N$ is a critical point of $F$ precisely when the following condition is satisfied: 
$$dF\in span(dA_{ij})$$

Regarding the space $span(dA_{ij})$, this consists of the following quantities:
\begin{eqnarray*}
\sum_{ij}M_{ij}dA_{ij}
&=&\sum_{ijk}M_{ij}(U_{ik}dU_{jk}+U_{jk}dU_{ik})\\
&=&\sum_{jk}(M^tU)_{jk}dU_{jk}+\sum_{ik}(MU)_{ik}dU_{ik}\\
&=&\sum_{ij}(M^tU)_{ij}dU_{ij}+\sum_{ij}(MU)_{ij}dU_{ij}
\end{eqnarray*}

In order to compute $dF$, observe first that, with $S_{ij}=sgn(U_{ij})$, we have:
$$d|U_{ij}|
=d\sqrt{U_{ij}^2}
=\frac{U_{ij}dU_{ij}}{|U_{ij}|}
=S_{ij}dU_{ij}$$

Now let us set, as in the statement:
$$W_{ij}=sgn(U_{ij})\varphi'(|U_{ij}|)$$

In terms of these variables, we obtain:
$$dF
=\sum_{ij}d\left(\varphi(|U_{ij}|)\right)
=\sum_{ij}\varphi'(|U_{ij}|)d|U_{ij}|
=\sum_{ij}W_{ij}dU_{ij}$$

We conclude that $U\in O_N$ is a critical point of $F$ if and only if there exists a matrix $M\in M_N(\mathbb R)$ such that the following two conditions are satisfied:
$$W=M^tU\quad,\quad 
W=MU$$

Now observe that these two equations can be written as follows:
$$M^t=WU^t\quad,\quad 
M=WU^t$$

Thus, the matrix $WU^t$ must be symmetric, as claimed.
\end{proof}

In order to process the above result, we can use the following notion:

\index{color decomposition}
\index{balanced matrix}
\index{semi-balanced matrix}

\begin{definition}
Given $U\in O_N$, we consider its ``color decomposition'' 
$$U=\sum_{r>0}rU_r$$
with $U_r\in M_N(-1,0,1)$ containing the sign components at $r>0$, and we call $U$:
\begin{enumerate}
\item Semi-balanced, if $U_rU^t$ and $U^tU_r$, with $r>0$, are all symmetric.

\item Balanced, if $U_rU_s^t$ and $U_r^tU_s$, with $r,s>0$, are all symmetric.
\end{enumerate}
\end{definition}

These conditions are quite natural, because for an orthogonal matrix $U\in O_N$, the relations $UU^t=U^tU=1$ translate as follows, in terms of the color decomposition:
$$\sum_{r>0}rU_rU^t=\sum_{r>0}rU^tU_r=1$$
$$\sum_{r,s>0}rsU_rU_s^t=\sum_{r,s>0}rsU_r^tU_s=1$$

Thus, our balancing conditions express the fact that the various components of the above sums are all symmetric. Now back to our critical point questions, we have:

\begin{theorem}
For a matrix $U\in O_N^*$, the following are equivalent:
\begin{enumerate}
\item $U$ is a critical point of $F(U)=\sum_{ij}\varphi(|U_{ij}|)$, for any $\varphi:[0,\infty)\to\mathbb R$.

\item $U$ is a critical point of all the $p$-norms, with $p\in[1,\infty)$.

\item $U$ is semi-balanced, in the above sense.
\end{enumerate}
\end{theorem}

\begin{proof}
We use the critical point criterion found in Theorem 3.4. In terms of the color decomposition, the matrix constructed there is given by:
\begin{eqnarray*}
(WU^t)_{ij}
&=&\sum_k{\rm sgn}(U_{ik})\varphi'(|U_{ik}|)U_{jk}\\
&=&\sum_{r>0}\varphi'(r)\sum_{k,|U_{ik}|=r}{\rm sgn}(U_{ik})U_{jk}\\
&=&\sum_{r>0}\varphi'(r)\sum_k(U_r)_{ik}U_{jk}\\
&=&\sum_{r>0}\varphi'(r)(U_rU^t)_{ij}
\end{eqnarray*}

Thus we have the following formula:
$$WU^t=\sum_{r>0}\varphi'(r)U_rU^t$$

Now when the function $\varphi:[0,\infty)\to\mathbb R$ varies, either as an arbitrary differentiable function, or as a power function $\varphi(x)=x^p$ with $p\in[1,\infty)$, the individual components of this sum must be all self-adjoint, and this leads to the conclusion in the statement.
\end{proof}

In practice now, most of the known examples of semi-balanced matrices are actually balanced, so we will investigate instead this latter class of matrices. Following \cite{bn3}, we have the following collection of simple facts, regarding such matrices:

\begin{theorem}
The class of balanced matrices is as follows:
\begin{enumerate}
\item It contains the matrices $U=H/\sqrt{N}$, with $H\in M_N(\pm1)$ Hadamard.

\item It is stable under transposition.

\item It is stable under taking tensor products.

\item It is stable under Hadamard equivalence.

\item It contains the matrix $V_N=\frac{1}{N}(2\mathbb I_N-N1_N)$, where $\mathbb I_N$ is the all-$1$ matrix.
\end{enumerate}
\end{theorem}

\begin{proof}
All these results are elementary, the proof being as follows:

\medskip

(1) Here $U\in O_N$ follows from the Hadamard condition, and since there is only one color component, namely $U_{1/\sqrt{N}}=H$, the balancing condition is satisfied as well.

\medskip

(2) Assuming that $U=\sum_{r>0}rU_r$ is the color decomposition of a given matrix $U\in O_N$, the color decomposition of the transposed matrix $U^t$ is as follows:
$$U^t=\sum_{r>0}rU_r^t$$

It follows that if $U$ is balanced, so is the transposed matrix $U^t$. 

\medskip

(3) Assuming that $U=\sum_{r>0}rU_r$ and $V=\sum_{s>0}sV_s$ are the color decompositions of two given orthogonal matrices $U,V$, we have:
$$U\otimes V
=\sum_{r,s>0}rs\cdot U_r\otimes V_s
=\sum_{p>0}p\sum_{p=rs}U_r\otimes V_s$$

Thus the color components of $W=U\otimes V$ are the following matrices: 
$$W_p=\sum_{p=rs}U_r\otimes V_s$$

It follows that if $U,V$ are both balanced, then so is $W=U\otimes V$.

\medskip

(4) We recall that the Hadamard equivalence consists in permuting rows and columns, and switching signs on rows and columns. Since all these operations correspond to certain conjugations at the level of the matrices $U_rU_s^t,U_r^tU_s$, we obtain the result.

\medskip

(5) The matrix in the statement, which goes back to \cite{bnz}, is as follows:
$$V_N=\frac{1}{N}
\begin{pmatrix}
2-N&2&\ldots&2\\
2&2-N&\ldots&2\\
\ldots&\ldots&\ldots&\ldots\\
2&2&\ldots&2-N
\end{pmatrix}$$

Observe that this matrix is indeed orthogonal, its rows being of norm one, and pairwise orthogonal. The color components of this matrix being $V_{2/N-1}=1_N$ and $V_{2/N}=\mathbb I_N-1_N$, it follows that this matrix is balanced as well, as claimed.
\end{proof}

Let us look now more in detail at the matrix $V_N$ from the above statement, and at the matrices having similar properties. Following \cite{bnz}, let us start our study with:

\index{block design}
\index{BIBD}

\begin{definition}
An $(a,b,c)$ pattern is a matrix $M\in M_N(0,1)$, with $N=a+2b+c$, such that any two rows look as follows,
$$\begin{matrix}
0\ldots 0&0\ldots 0&1\ldots 1&1\ldots 1\\
\underbrace{0\ldots 0}_a&\underbrace{1\ldots 1}_b&\underbrace{0\ldots 0}_b&\underbrace{1\ldots 1}_c
\end{matrix}$$
up to a permutation of the columns.
\end{definition}

As explained in \cite{bnz}, there are many interesting examples of $(a,b,c)$ patterns, coming from the balanced incomplete block designs (BIBD), and all these examples can produce two-entry unitary matrices, by replacing the $0,1$ entries with suitable numbers $x,y$. For more on BIBD and design theory, we refer to Colbourn-Dinitz \cite{cdi} or Stinson \cite{sti}.

\bigskip

Now back to the matrix $V_N$ from Theorem 3.7 (5), observe that this matrix comes from a $(0,1,N-2)$ pattern. And also, independently of this, this matrix has the remarkable property of being at the same time circulant and self-adjoint. We have in fact:

\index{circulant matrix}

\begin{theorem}
The following matrices are balanced:
\begin{enumerate}
\item The orthogonal matrices coming from $(a,b,c)$ patterns.

\item The orthogonal matrices which are circulant and symmetric.
\end{enumerate}
\end{theorem}

\begin{proof}
These observations basically go back to \cite{bnz}, the proofs being as follows:

\medskip

(1) If we denote by $P,Q\in M_N(0,1)$ the matrices describing the positions of the $0,1$ entries inside the pattern, then we have the following formulae:
\begin{eqnarray*}
PP^t=P^tP&=&a\mathbb I_N+b1_N\\
QQ^t=Q^tQ&=&c\mathbb I_N+b1_N\\
PQ^t=P^tQ=QP^t=Q^tP&=&b\mathbb I_N-b1_N
\end{eqnarray*}

Since all these matrices are symmetric, $U$ is balanced, as claimed.

\medskip

(2) Assume that $U\in O_N$ is circulant, $U_{ij}=\gamma_{j-i}$, and in addition symmetric, which means $\gamma_i=\gamma_{-i}$. Consider the following sets, which must satisfy $D_r=-D_r$:
$$D_r=\{k:|\gamma_r|=k\}$$

In terms of these sets, we have the following formula:
\begin{eqnarray*}
(U_rU_s^t)_{ij}
&=&\sum_k(U_r)_{ik}(U_s)_{jk}\\
&=&\sum_k\delta_{|\gamma_{k-i}|,r}\,{\rm sgn}(\gamma_{k-i})\cdot\delta_{|\gamma_{k-j}|,s}\,{\rm sgn}(\gamma_{k-j})\\
&=&\sum_{k\in(D_r+i)\cap(D_s+j)}{\rm sgn}(\gamma_{k-i})\,{\rm sgn}(\gamma_{k-j})
\end{eqnarray*}

With $k=i+j-m$ we obtain, by using $D_r=-D_r$, and then $\gamma_i=\gamma_{-i}$:
\begin{eqnarray*}
(U_rU_s^t)_{ij}
&=&\sum_{m\in(-D_r+j)\cap(-D_s+i)}{\rm sgn}(\gamma_{j-m})\,{\rm sgn}(\gamma_{i-m})\\
&=&\sum_{m\in(D_r+i)\cap(D_r+j)}{\rm sgn}(\gamma_{j-m})\,{\rm sgn}(\gamma_{i-m})\\
&=&\sum_{m\in(D_r+i)\cap(D_r+j)}{\rm sgn}(\gamma_{m-j})\,{\rm sgn}(\gamma_{m-i})
\end{eqnarray*}

Now by interchanging $i\leftrightarrow j$, and with $m\to k$, this formula becomes:
$$(U_rU_s^t)_{ji}=\sum_{k\in(D_r+i)\cap(D_r+j)}{\rm sgn}(\gamma_{k-i})\,{\rm sgn}(\gamma_{k-j})$$

By comparing with the previous formula, we deduce that the matrix $U_rU_s^t$ is symmetric, as claimed. The proof for $U_r^tU_s$ is similar. 
\end{proof}

As a conclusion to all this, the study of the critical points of the various $p$-norms on $O_N$ has led us into the class of balanced matrices, which looks like an interesting class, which is waiting to be further investigated. We will be back to this.

\section*{3b. Second derivatives}

Let us get now into analytic questions. As in Theorem 3.4, it is convenient to do the computations in a general framework, with a function as follows:
$$F(U)=\sum_{ij}\psi(U_{ij}^2)$$

Consider the following function, depending on $t>0$ small:
$$f(t)
=F(Ue^{tA})
=\sum_{ij}\psi\left((Ue^{tA})_{ij}^2\right)$$

\index{Lie algebra}

Here $U\in O_N$ is an arbitrary orthogonal matrix, and $A\in M_N(\mathbb R)$ is assumed to be antisymmetric, $A^t=-A$, with this latter assumption needed for having $e^A\in O_N$. Let us first compute the derivative of $f$. Following \cite{bn3}, we have the following result:

\begin{proposition}
We have the following formula,
$$f'(t)=2\sum_{ij}\psi'((Ue^{tA})_{ij}^2)(UAe^{tA})_{ij}(Ue^{tA})_{ij}$$
valid for any $U\in O_N$, and any $A\in M_N(\mathbb R)$ antisymmetric.
\end{proposition}

\begin{proof}
The matrices $U,e^{tA}$ being both orthogonal, we have:
\begin{eqnarray*}
(Ue^{tA})_{ij}^2
&=&(Ue^{tA})_{ij}((Ue^{tA})^t)_{ji}\\
&=&(Ue^{tA})_{ij}(e^{tA^t}U^t)_{ji}\\
&=&(Ue^{tA})_{ij}(e^{-tA}U^t)_{ji}
\end{eqnarray*}

We can now differentiate our function $f$, and by using once again the orthogonality of the matrices $U,e^{tA}$, along with the formula $A^t=-A$, we obtain:
\begin{eqnarray*}
f'(t)
&=&\sum_{ij}\psi'((Ue^{tA})_{ij}^2)\left[(UAe^{tA})_{ij}(e^{-tA}U^t)_{ji}-(Ue^{tA})_{ij}(e^{-tA}AU^t)_{ji}\right]\\
&=&\sum_{ij}\psi'((Ue^{tA})_{ij}^2)\left[(UAe^{tA})_{ij}((e^{-tA}U^t)^t)_{ij}-(Ue^{tA})_{ij}((e^{-tA}AU^t)^t)_{ij}\right]\\
&=&\sum_{ij}\psi'((Ue^{tA})_{ij}^2)\left[(UAe^{tA})_{ij}(Ue^{tA})_{ij}+(Ue^{tA})_{ij}(UAe^{tA})_{ij}\right]
\end{eqnarray*}

But this gives the formula in the statement, and we are done.
\end{proof}

Before computing the second derivative, let us evaluate $f'(0)$. In terms of the color decomposition $U=\sum_{r>0}rU_r$ of our matrix, the result is:

\begin{proposition}
We have the following formula,
$$f'(0)=2\sum_{r>0}r\psi'(r^2)Tr(U_r^tUA)$$
where the matrices $U_r\in M_N(-1,0,1)$ are the color components of $U$.
\end{proposition}

\begin{proof}
We use the formula in Proposition 3.10. At $t=0$, we obtain:
$$f'(0)=2\sum_{ij}\psi'(U_{ij}^2)(UA)_{ij}U_{ij}$$

Consider now the color decomposition of $U$. We have the following formulae:
\begin{eqnarray*}
U_{ij}=\sum_{r>0}r(U_r)_{ij}
&\implies&U_{ij}^2=\sum_{r>0}r^2|(U_r)_{ij}|\\
&\implies&\psi'(U_{ij}^2)=\sum_{r>0}\psi'(r^2)|(U_r)_{ij}|
\end{eqnarray*}

Now by getting back to the above formula of $f'(0)$, we obtain:
$$f'(0)=2\sum_{r>0}\psi'(r^2)\sum_{ij}(UA)_{ij}U_{ij}|(U_r)_{ij}|$$

Our claim now is that we have the following formula:
$$U_{ij}|(U_r)_{ij}|=r(U_r)_{ij}$$

Indeed, in the case $|U_{ij}|\neq r$ this formula reads $U_{ij}\cdot 0=r\cdot 0$, which is true, and in the case $|U_{ij}|=r$ this formula reads $rS_{ij}\cdot 1=r\cdot S_{ij}$, which is once again true. Thus:
$$f'(0)=2\sum_{r>0}r\psi'(r^2)\sum_{ij}(UA)_{ij}(U_r)_{ij}$$

But this gives the formula in the statement, and we are done.
\end{proof}

Let us compute now the second derivative. The result here is as follows:

\begin{proposition}
We have the following formula,
\begin{eqnarray*}
f''(0)
&=&4\sum_{ij}\psi''(U_{ij}^2)\left[(UA)_{ij}U_{ij}\right]^2\\
&&+2\sum_{ij}\psi'(U_{ij}^2)\left[(UA^2)_{ij}U_{ij}\right]\\
&&+2\sum_{ij}\psi'(U_{ij}^2)(UA)_{ij}^2
\end{eqnarray*}
valid for any $U\in O_N$, and any $A\in M_N(\mathbb R)$ antisymmetric.
\end{proposition}

\begin{proof}
We use the formula in Proposition 3.10, namely:
$$f'(t)=2\sum_{ij}\psi'((Ue^{tA})_{ij}^2)(UAe^{tA})_{ij}(Ue^{tA})_{ij}$$

Since the term on the right, or rather its double, appears as the derivative of the quantity $(Ue^{tA})_{ij}^2$, when differentiating a second time, we obtain:
\begin{eqnarray*}
f''(t)
&=&4\sum_{ij}\psi''((Ue^{tA})_{ij}^2)\left[(UAe^{tA})_{ij}(Ue^{tA})_{ij}\right]^2\\
&&+2\sum_{ij}\psi'((Ue^{tA})_{ij}^2)\left[(UAe^{tA})_{ij}(Ue^{tA})_{ij}\right]'
\end{eqnarray*}

In order to compute now the missing derivative, observe that we have:
$$\left[(UAe^{tA})_{ij}(Ue^{tA})_{ij}\right]'
=(UA^2e^{tA})_{ij}(Ue^{tA})_{ij}+(UAe^{tA})_{ij}^2$$

Summing up, we have obtained the following formula:
\begin{eqnarray*}
f''(t)
&=&4\sum_{ij}\psi''((Ue^{tA})_{ij}^2)\left[(UAe^{tA})_{ij}(Ue^{tA})_{ij}\right]^2\\
&&+2\sum_{ij}\psi'((Ue^{tA})_{ij}^2)\left[(UA^2e^{tA})_{ij}(Ue^{tA})_{ij}\right]\\
&&+2\sum_{ij}\psi'((Ue^{tA})_{ij}^2)(UAe^{tA})_{ij}^2
\end{eqnarray*}

But at $t=0$ this gives the formula in the statement, and we are done.
\end{proof}

For the function $\psi(x)=\sqrt{x}$, corresponding to the functional $F(U)=||U||_1$, there are some simplifications, that we will work out now in detail. First, we have:

\begin{proposition}
For the function $F(U)=||U||_1$ we have the formula
$$f''(0)=Tr(S^tUA^2)$$
valid for any antisymmetric matrix $A$, where $S_{ij}={\rm sgn}(U_{ij})$.
\end{proposition}

\begin{proof}
We use the formula in Proposition 3.12, with the following data:
$$\psi(x)=\sqrt{x}\quad,\quad
\psi'(x)=\frac{1}{2\sqrt{x}}\quad,\quad
\psi''(x)=-\frac{1}{4x\sqrt{x}}$$

We therefore obtain the following formula:
\begin{eqnarray*}
f''(0)
&=&-\sum_{ij}\frac{\left[(UA)_{ij}U_{ij}\right]^2}{|U_{ij}|^3}
+\sum_{ij}\frac{(UA^2)_{ij}U_{ij}}{|U_{ij}|}
+\sum_{ij}\frac{(UA)_{ij}^2}{|U_{ij}|}\\
&=&-\sum_{ij}\frac{(UA)_{ij}^2}{|U_{ij}|}
+\sum_{ij}(UA^2)_{ij}S_{ij}
+\sum_{ij}\frac{(UA)_{ij}^2}{|U_{ij}|}\\
&=&\sum_{ij}(UA^2)_{ij}S_{ij}
\end{eqnarray*}

But this gives the formula in the statement, and we are done.
\end{proof}

We are therefore led to the following result, from \cite{bn3}, regarding the 1-norm:

\index{local maximizer}

\begin{theorem}
A matrix $U\in O_N$ locally maximizes the $1$-norm on $O_N$ precisely when the following conditions are satisfied:
\begin{enumerate}
\item The matrix $U$ has nonzero entries, $U\in O_N^*$.

\item The matrix $X=S^tU$ is symmetric, where $S_{ij}={\rm sgn}(U_{ij})$.

\item We have $Tr(XA^2)\leq0$, for any antisymmetric matrix $A\in M_N(\mathbb R)$.
\end{enumerate}
\end{theorem}

\begin{proof}
This follows the results that we have, with (1,2,3) coming respectively from Theorem 3.3, Theorem 3.4 and Proposition 3.13.
\end{proof}

In order to further improve the above result, we will need:

\index{smallest eigenvalues}

\begin{proposition}
For a symmetric matrix $X\in M_N(\mathbb R)$, the following are equivalent:
\begin{enumerate}
\item $Tr(XA^2)\leq0$, for any antisymmetric matrix $A$.

\item The sum of the two smallest eigenvalues of $X$ is positive. 
\end{enumerate}
\end{proposition}

\begin{proof}
Consider the following vector, which is antisymmetric:
$$a=\sum_{ij}A_{ij}e_i\otimes e_j$$

In terms of this vector, we have the following formula:
\begin{eqnarray*}
Tr(XA^2)
&=&<X,A^2>\\
&=&-<AX,A>\\
&=&-<a,(1\otimes X)a>
\end{eqnarray*}

Thus the condition (1) is equivalent to $P(1\otimes X)P$ being positive, with $P$ being the orthogonal projection on the antisymmetric subspace in $\mathbb R^N\otimes\mathbb R^N$. Now observe that for any two eigenvectors $x_i \perp x_j$ of $X$, with eigenvalues $\lambda_i, \lambda_j$, we have:
\begin{eqnarray*} 
P(1\otimes X)P(x_i\otimes x_j-x_j\otimes x_i)
&=&P(\lambda_j x_i\otimes x_j-\lambda_i x_j\otimes x_i)\\
&=&\frac{\lambda_i +\lambda_j}{2}(x_i\otimes x_j-x_j\otimes x_i)
\end{eqnarray*}

Thus, we are led to the conclusion in the statement.
\end{proof}

Following \cite{bn3}, we can now formulate a final result on the subject, which improves some previous findings from \cite{bcs}, and from \cite{bnz}, as follows:

\begin{theorem}
A matrix $U\in O_N$ locally maximizes the $1$-norm on $O_N$ precisely when it has nonzero entries, and when the following matrix, with $S_{ij}={\rm sgn}(U_{ij})$, 
$$X=S^tU$$
is symmetric, and the sum of its two smallest eigenvalues is positive.
\end{theorem}

\begin{proof}
This follows indeed from our main result so far, Theorem 3.14, by taking into account the positivity criterion from Proposition 3.15.
\end{proof}

In terms of the almost Hadamard matrices, as introduced in Definition 3.2, as rescaled versions of the above matrices, the above result reformulates as follows:

\index{AHM}
\index{almost Hadamard matrix}

\begin{theorem}
The almost Hadamard matrices are the matrices $H\in\sqrt{N}O_N$ having nonzero entries, and which are such that the following matrix, with $S_{ij}={\rm sgn}(H_{ij})$, 
$$X=S^tH$$
is symmetric, and the sum of its two smallest eigenvalues is positive.
\end{theorem}

\begin{proof}
This is a reformulation of Theorem 3.16, by rescaling everything by $\sqrt{N}$, as to reach to the objects axiomatized in Definition 3.2.
\end{proof}

Regarding now the examples of such matrices, which can be useful for various reasons, especially at values $N\notin4\mathbb N$, there are many of them, and we will discuss them gradually, in what follows. To start with, we have the following general result, from \cite{bcs}, \cite{bnz}:

\begin{theorem}
The class of almost Hadamard matrices has the following properties:
\begin{enumerate}
\item It contains all the Hadamard matrices.

\item It is stable under transposition.

\item It is stable under taking tensor products.

\item It is stable under Hadamard equivalence.

\item It contains the matrix $K_N=\frac{1}{\sqrt{N}}(2\mathbb I_N-N1_N)$.
\end{enumerate}
\end{theorem}

\begin{proof}
All the assertions are clear from what we have, as follows:

\medskip

(1) This follows either from Theorem 3.1, which shows that Hadamard implies almost Hadamard, without any need for further computations, or from the fact that if $H$ is Hadamard then $U=H/\sqrt{N}$ is orthogonal, and $SU^t=HU^t=\sqrt{N}1_N$ is positive.

\medskip

(2) This follows either from definitions, because the transposition operation preserves the local maximizers of the 1-norm, or from Theorem 3.17.

\medskip

(3) For a tensor product of almost Hadamard matrices $H=H'\otimes H''$ we have $U=U'\otimes U''$ and $S=S'\otimes S''$, so that $U$ is unitary and $SU^t$ is symmetric, with the sum of the two smallest eigenvalues being positive, as claimed. 

\medskip

(4) This follows either from definitions, because the Hadamard equivalence preserves the local maximizers of the 1-norm, or from Theorem 3.17.

\medskip

(5) We know from Theorem 3.7 that the matrix $U=K_N/\sqrt{N}$ is orthogonal. Also, we have $S=\mathbb I_N-21_N$, and so $SU^t$ is positive, because with $J_N=\mathbb I_N/N$ we have:
\begin{eqnarray*}
SU^t
&=&(NJ_N-21_N)(2J_N-1_N)\\
&=&(N-2)J_N+2(1_N-J_N)
\end{eqnarray*}

Thus, we are led to the conclusion in the statement.
\end{proof}

Observe the similarity between the above result and Theorem 3.7, which was about the balanced matrices. However, these two statements, even when properly rescaled, either both on $O_N$ or both on $\sqrt{N}O_N$, do not exactly cover the same class of matrices. Based on this analogy, however, we can look for explicit examples of almost Hadamard matrices by taking some inspiration from the main examples of balanced matrices, from Theorem 3.9. We will discuss this in the remainder of this chapter.

\section*{3c. Circulant matrices}

We have two classes of matrices to be investigated, generalizing the matrix $K_N$ from Theorem 3.18, namely the circulant matrices, and the 2-entry matrices. Following the work in \cite{bnz}, let us start with the circulant matrices. We let $F\in U_N$ be the normalized Fourier matrix, given by $F_{ij}=w^{ij}/\sqrt{N}$, where $w=e^{2\pi i/N}$. Also, we make the convention that associated to any vector $\alpha\in\mathbb C^N$ is the following diagonal matrix:
$$\alpha'=
\begin{pmatrix}
\alpha_0\\
&\ddots\\
&&\alpha_{N-1}
\end{pmatrix}$$

With these conventions, we have the following well-known result:

\index{circulant matrix}
\index{Fourier-diagonal}
\index{Fourier matrix}

\begin{proposition}
For a matrix $H\in M_N(\mathbb C)$, the following are equivalent:
\begin{enumerate}
\item $H$ is circulant, i.e. $H_{ij}=\gamma_{j-i}$, for a certain vector $\gamma\in\mathbb C^N$.

\item $H$ is Fourier-diagonal, i.e. $H=FDF^*$, with $D\in M_N(\mathbb C)$ diagonal.
\end{enumerate}
In addition, if so is the case, then with $D=\sqrt{N}\alpha'$ we have $\gamma=F\alpha$.
\end{proposition}

\begin{proof}
(1)$\implies$(2) The matrix $D=F^*HF$ is indeed diagonal, given by:
$$D_{ij}=\frac{1}{N}\sum_{kl}w^{jl-ik}\gamma_{l-k}=\delta_{ij}\sum_rw^{jr}\gamma_r$$ 

(2)$\implies$(1) The matrix $H=FDF^*$ is indeed circulant, given by:
$$H_{ij}=\sum_kF_{ik}D_{kk}\bar{F}_{jk}=\frac{1}{N}\sum_kw^{(i-j)k}D_{kk}$$

Finally, the last assertion is clear from the above formula of $H_{ij}$.
\end{proof}

Let us investigate now the circulant orthogonal matrices. We have:

\begin{proposition}
For a matrix $U\in M_N(\mathbb C)$, the following are equivalent:
\begin{enumerate}
\item $U$ is orthogonal and circulant.

\item $U=F\alpha'F^*$ with $\alpha\in\mathbb T^N$ satisfying $\bar{\alpha}_i=\alpha_{-i}$ for any $i$.
\end{enumerate}
\end{proposition}

\begin{proof}
We will use many times the fact that given a vector $\alpha\in\mathbb C^N$, the vector $\gamma=F\alpha$ is real if and only if the following happens, for any $i$:
$$\bar{\alpha}_i=\alpha_{-i}$$

This follows indeed from $\overline{F\alpha}=F\tilde{\alpha}$, with $\tilde{\alpha}_i=\bar{\alpha}_{-i}$.

\medskip

(1)$\implies$(2) Write $H_{ij}=\gamma_{j-i}$ with $\gamma\in\mathbb R^N$. By using Proposition 3.19 we obtain $H=FDF^*$ with $D=\sqrt{N}\alpha'$ and $\gamma=F\alpha$. Now since $U=F\alpha'F^*$ is unitary, so is $\alpha'$, so we must have $\alpha\in\mathbb T^N$. Finally, since $\gamma$ is real we have $\bar{\alpha}_i=\alpha_{-i}$, and we are done.

\medskip

(2)$\implies$(1) We know from Proposition 3.19 that $U$ is circulant. Also, from $\alpha\in\mathbb T^N$ we obtain that $\alpha'$ is unitary, and so must be $U$. Finally, since we have $\bar{\alpha}_i=\alpha_{-i}$, the vector $\gamma=F\alpha$ is real, and hence we have $U\in M_N(\mathbb R)$, which finishes the proof.
\end{proof}

Let us discuss now the almost Hadamard case. First, in the usual Hadamard case, the known examples and the corresponding $\alpha$-vectors are as follows:

\begin{proposition}
The known circulant Hadamard matrices, namely
$$\pm\begin{pmatrix}
-1\!\!&\!\!1\!\!&\!\!1\!\!&\!\!1\\
1\!\!&\!\!-1\!\!&\!\!1\!\!&\!\!1\\
1\!\!&\!\!1\!\!&\!\!-1\!\!&\!\!1\\
1\!\!&\!\!1\!\!&\!\!1\!\!&\!\!-1
\end{pmatrix}\qquad,\qquad
\pm\begin{pmatrix}
1\!\!&\!\!-1\!\!&\!\!1\!\!&\!\!1\\
1\!\!&\!\!1\!\!&\!\!-1\!\!&\!\!1\\
1\!\!&\!\!1\!\!&\!\!1\!\!&\!\!-1\\
-1\!\!&\!\!1\!\!&\!\!1\!\!&\!\!1
\end{pmatrix}$$
$$\pm\begin{pmatrix}
1\!\!&\!\!1\!\!&\!\!-1\!\!&\!\!1\\
1\!\!&\!\!1\!\!&\!\!1\!\!&\!\!-1\\
-1\!\!&\!\!1\!\!&\!\!1\!\!&\!\!1\\
1\!\!&\!\!-1\!\!&\!\!1\!\!&\!\!1
\end{pmatrix}\qquad,\qquad
\pm\begin{pmatrix}
1\!\!&\!\!1\!\!&\!\!1\!\!&\!\!-1\\
-1\!\!&\!\!1\!\!&\!\!1\!\!&\!\!1\\
1\!\!&\!\!-1\!\!&\!\!1\!\!&\!\!1\\
1\!\!&\!\!1\!\!&\!\!-1\!\!&\!\!1
\end{pmatrix}$$
come respectively from the following $\alpha$ vectors, via the above construction:
$$\pm(1,-1,-1,-1)\qquad,\qquad\pm(1,-i,1,i)$$
$$\pm(1,1,-1,1)\qquad\ \ \ \ ,\ \qquad\pm(1,i,1,-i)$$
\end{proposition}

\begin{proof}
At $N=4$ the conjugate of the Fourier matrix is given by:
$$F^*=\frac{1}{2}\begin{pmatrix}
1&1&1&1\\
1&-i&-1&i\\
1&-1&1&-1\\
1&i&-1&-i
\end{pmatrix}$$

Thus the vectors $\alpha=F^*\gamma$ are indeed those in the statement.
\end{proof}

Following \cite{bnz}, we have the following generalization of the above matrices:

\begin{proposition}
If $q^N=1$ then the vector 
$$\alpha=\pm(1,-q,-q^2,\ldots,-q^{N-1})$$
produces an almost Hadamard matrix, equivalent to $K_N=\frac{1}{\sqrt{N}}(2\mathbb I_N-N1_N)$. 
\end{proposition}

\begin{proof}
Observe first that these matrices generalize those in Proposition 3.21. Indeed, at $N=4$ the choices for $q$ are $1,i,-1,-i$, and this gives the above $\alpha$-vectors.

Assume that the $\pm$ sign in the statement is $+$. With $q=w^r$, we have:
\begin{eqnarray*}
\sqrt{N}\gamma_i
&=&\sum_{k=0}^{N-1}w^{ik}\alpha_k\\
&=&1-\sum_{k=1}^{N-1}w^{(i+r)k}\\
&=&2-\sum_{k=0}^{N-1}w^{(i+r)k}\\
&=&2-\delta_{i,-r}N
\end{eqnarray*}

In terms of the standard long cycle $(C_N)_{ij}=\delta_{i+1,j}$, we obtain:
$$H=\frac{1}{\sqrt{N}}(2\mathbb I_N-NC_N^{-r})$$

Thus $H$ is equivalent to $K_N$, and by Theorem 3.18, it is almost Hadamard.
\end{proof}

In general, the construction of circulant almost Hadamard matrices is quite a tricky problem. At the abstract level, we have the following result, from \cite{bnz}:

\begin{proposition}
A circulant matrix $H\in M_N(\mathbb R^*)$, written $H_{ij}=\gamma_{j-i}$, is almost Hadamard provided that the following conditions are satisfied:
\begin{enumerate}
\item The vector $\alpha=F^*\gamma$ satisfies $\alpha\in\mathbb T^N$.

\item With $\varepsilon={\rm sgn}(\gamma)$, $\rho_i=\sum_r\varepsilon_r\gamma_{i+r}$ and $\nu=F^*\rho$, we have $\nu>0$.
\end{enumerate}
In addition, if so is the case, then $\bar{\alpha}_i=\alpha_{-i}$, $\rho_i=\rho_{-i}$ and $\nu_i=\nu_{-i}$ for any $i$.
\end{proposition}

\begin{proof}
We know from Theorem 3.17 our matrix $H$ is almost Hadamard if the matrix $U=H/\sqrt{N}$ is orthogonal and $SU^t>0$, where $S_{ij}={\rm sgn}(U_{ij})$. By Proposition 3.19 the orthogonality of $U$ is equivalent to the condition (1). Regarding now the condition $SU^t>0$, this is equivalent to $S^tU>0$. But, with $k=i-r$, we have:
\begin{eqnarray*}
(S^tH)_{ij}
&=&\sum_kS_{ki}H_{kj}\\
&=&\sum_k\varepsilon_{i-k}\gamma_{j-k}\\
&=&\sum_r\varepsilon_r\gamma_{j-i+r}\\
&=&\rho_{j-i}
\end{eqnarray*}

Thus $S^tU$ is circulant, with $\rho/\sqrt{N}$ as first row. From Proposition 3.19 we get $S^tU=FLF^*$ with $L=\nu'$ and $\nu=F^*\rho$, so $S^tU>0$ iff $\nu>0$, which is the condition (2). Finally, the assertions about $\alpha,\nu$ follow from the fact that the vectors $F\alpha,F\nu$ are real. As for the assertion about $\rho$, this follows from the fact that $S^tU$ is symmetric.
\end{proof}

Here are now the main examples of such matrices, once again following \cite{bnz}:

\begin{theorem}
For $N$ odd the following matrix is almost Hadamard,
$$L_N=\frac{1}{\sqrt{N}}
\begin{pmatrix}
1&-\cos^{-1}\frac{\pi}{N}&\cos^{-1}\frac{2\pi}{N}&\ldots\ldots&\cos^{-1}\frac{(N-1)\pi}{N}\\
\cos^{-1}\frac{(N-1)\pi}{N}&1&-\cos^{-1}\frac{\pi}{N}&\ldots\ldots&-\cos^{-1}\frac{(N-2)\pi}{N}\\
\vdots&\vdots&\vdots&&\vdots\\
\vdots&\vdots&\vdots&&\vdots\\
-\cos^{-1}\frac{\pi}{N}&\cos^{-1}\frac{2\pi}{N}&-\cos^{-1}\frac{3\pi}{N}&\ldots\ldots&1
\end{pmatrix}$$
and comes from an $\alpha$-vector having all entries equal to $1$ or $-1$.
\end{theorem}

\begin{proof}
Write $N=2n+1$, and consider the following vector:
$$\alpha_i=\begin{cases}
(-1)^{n+i}&{\rm for }\ i=0,1,\ldots,n\\
(-1)^{n+i+1}&{\rm for}\ i=n+1,\ldots,2n
\end{cases}$$

Let us first prove that $(L_N)_{ij}=\gamma_{j-i}$, where $\gamma=F\alpha$. With $w=e^{2\pi i/N}$ we have:
\begin{eqnarray*}
\sqrt{N}\gamma_i
&=&\sum_{j=0}^{2n}w^{ij}\alpha_j\\
&=&\sum_{j=0}^n(-1)^{n+j}w^{ij}+\sum_{j=1}^n(-1)^{n+(N-j)+1}w^{i(N-j)}
\end{eqnarray*}

Now since $N$ is odd, and since $w^N=1$, we obtain:
\begin{eqnarray*}
\sqrt{N}\gamma_i
&=&\sum_{j=0}^n(-1)^{n+j}w^{ij}+\sum_{j=1}^n(-1)^{n-j}w^{-ij}\\
&=&\sum_{j=-n}^n(-1)^{n+j}w^{ij}
\end{eqnarray*}

By computing the sum on the right, with $\xi=e^{\pi i/N}$ we get, as claimed:
\begin{eqnarray*}
\sqrt{N}\gamma_i
&=&\frac{2w^{-ni}}{1+w^i}\\
&=&\frac{2\xi^{-2ni}}{1+\xi^{2i}}\\
&=&\frac{2\xi^{-Ni}}{\xi^{-i}+\xi^i}\\
&=&(-1)^i\cos^{-1}\frac{i\pi}{N}
\end{eqnarray*}

In order to prove now that $L_N$ is almost Hadamard, we use Proposition 3.23. Since the sign vector is simply $\varepsilon=(-1)^n\alpha$, the vector $\rho_i=\sum_r\varepsilon_r\gamma_{i+r}$ is given by:
\begin{eqnarray*}
\sqrt{N}\rho_i
&=&(-1)^n\sum_{r=0}^{2n}\alpha_r\sum_{j=-n}^n(-1)^{n+j}w^{(i+r)j}\\
&=&\sum_{j=-n}^n(-1)^jw^{ij}\sum_{r=0}^{2n}\alpha_rw^{rj}
\end{eqnarray*}

Now since the last sum on the right is $(\sqrt{N}F\alpha)_j=\sqrt{N}\gamma_j$, we obtain:
\begin{eqnarray*}
\rho_i
&=&\sum_{j=-n}^n(-1)^jw^{ij}\gamma_j\\
&=&\frac{1}{\sqrt{N}}\sum_{j=-n}^n(-1)^jw^{ij}\sum_{k=-n}^n(-1)^{n+k}w^{jk}
\end{eqnarray*}

Thus we have the following formula:
$$\rho_i=\frac{(-1)^n}{\sqrt{N}}\sum_{j=-n}^n\sum_{k=-n}^n(-1)^{j+k}w^{(i+k)j}$$

Let us compute now the vector $\nu=F^*\rho$. We have:
\begin{eqnarray*}
\nu_l
&=&\frac{1}{\sqrt{N}}\sum_{i=0}^{2n}w^{-il}\rho_i\\
&=&\frac{(-1)^n}{N}\sum_{j=-n}^n\sum_{k=-n}^n(-1)^{j+k}w^{jk}\sum_{i=0}^{2n}w^{i(j-l)}
\end{eqnarray*}

The sum on the right is $N\delta_{jl}$, with both $j,l$ taken modulo $N$, so it is equal to $N\delta_{jL}$, where $L=l$ for $l\leq n$, and $L=l-N$ for $l>n$. We obtain:
\begin{eqnarray*}
\nu_l
&=&(-1)^n\sum_{k=-n}^n(-1)^{L+k}w^{Lk}\\
&=&(-1)^{n+L}\sum_{k=-n}^n(-w^L)^k
\end{eqnarray*}

With $\xi=e^{\pi i/N}$ as before, this gives the following formula:
\begin{eqnarray*}
\nu_l
&=&(-1)^{n+L}\frac{2(-w^L)^{-n}}{1+w^L}\\
&=&(-1)^L\frac{2w^{-nL}}{1+w^L}
\end{eqnarray*}

In terms of the variable $\xi=e^{\pi i/N}$, we obtain the following formula:
\begin{eqnarray*}
\nu_l
&=&(-1)^L\frac{2\xi^{-2nL}}{1+\xi^{2L}}\\
&=&(-1)^L\frac{2\xi^{-NL}}{\xi^{-L}+\xi^L}\\
&=&\cos^{-1}\frac{L\pi}{N}
\end{eqnarray*}

Now since $L\in[-n,n]$, all the entries of $\nu$ are positive, and we are done.
\end{proof}

At the level of examples now, at $N=3$ we obtain the matrix $L_3=-K_3$:
$$L_3=\frac{1}{\sqrt{3}}\begin{pmatrix}
1&-2&-2\\
-2&1&-2\\
-2&-2&1
\end{pmatrix}$$

At $N=5$ we obtain the following matrix, with $x=-\cos^{-1}\frac{\pi}{5}$, $y=\cos^{-1}\frac{2\pi}{5}$:
$$L_5=\frac{1}{\sqrt{5}}\begin{pmatrix}
1&x&y&y&x\\
x&1&x&y&y\\
y&x&1&x&y\\
y&y&x&1&x\\
x&y&y&x&1
\end{pmatrix}$$

For further examples of matrices of this type, and for a discussion of their 1-norms, which happen quite often to be optimal, or almost, we refer to \cite{bnz}.

\index{optimal AHM}
\index{optimal almost Hadamard matrix}

\section*{3d. Block designs}

Let us study now the almost Hadamard matrices having two entries, $H\in M_N(x,y)$, with $x,y\in\mathbb R$. These are related to design theory, so let us start with:

\begin{definition}
A filled $(a,b,c)$ pattern is a matrix $M\in M_N(x,y)$, with $N=a+2b+c$, such that any two rows look as follows, up to a permutation of columns:
$$\begin{matrix}
x\ldots x&x\ldots x&y\ldots y&y\ldots y\\
\underbrace{x\ldots x}_a&\underbrace{y\ldots y}_b&\underbrace{x\ldots x}_b&\underbrace{y\ldots y}_c
\end{matrix}$$
When the entries $x,y$ are the numbers $0,1$, we say that we have an $(a,b,c)$ pattern.
\end{definition}

There are many interesting examples of patterns coming from block designs, that we can use in order to construct almost Hadamard matrices. Let us begin with:

\index{block design}

\begin{definition}
A $(v,k,\lambda)$ symmetric balanced incomplete block design is a collection $B$ of subsets of a set $X$, called blocks, with the following properties:
\begin{enumerate}
\item $|X|=|B|=v$.

\item Each block contains exactly $k$ points from $X$.

\item Each pair of distinct points is contained in exactly $\lambda$ blocks of $B$.
\end{enumerate} 
\end{definition}

This is a standard definition in design theory, and for more we refer to Colbourn-Dinitz \cite{cdi} and Stinson \cite{sti}. In relation with our linear algebra questions, we will be interested in the incidence matrix of such a block design, which is the $v\times v$ matrix given by:
$$M_{bx}=\begin{cases}
1&\text{if }x\in b\\
0&\text{if }x\notin b
\end{cases}$$

The connection between designs and patterns comes from:

\begin{proposition}
If $N=a+2b+c$ then the adjacency matrix of any $(N,a+b,a)$ symmetric balanced incomplete block design is an $(a,b,c)$ pattern.
\end{proposition}

\begin{proof}
Let us replace the $0-1$ values in the adjacency matrix $M$ by abstract $x-y$ values. Then each row of $M$ contains $a+b$ copies of $x$ and $b+c$ copies of $y$, and since every pair of distinct blocks intersect in exactly $a$ points, we see that every pair of rows has exactly $a$ variables $x$ in matching positions, so that $M$ is an $(a,b,c)$ pattern.
\end{proof}

\index{Fano plane}

As a first example for all this, consider the Fano plane, which is the simplest instance of ``discrete geometry'', consisting of 7 points and 7 lines, as follows:
$$\xymatrix@R=9pt@C=10pt{
&&&&\bullet\ar@{-}[ddddd]\\
&&&&&&&\\
&&&&&&&\\
&&&&\ar@{-}@/^/[drr]\\
&&\bullet\ar@{-}[uuuurr]\ar@{-}@/^/[urr]\ar@{-}@/_/[dd]&&&&\bullet\ar@{-}[uuuull]&&\\
&&&&\bullet\ar@{-}[urr]\ar@{-}[ull]&&&&\\
&&\ar@{-}@/_/[drr]&&&&\ar@{-}@/^/[dll]\ar@{-}@/_/[uu]&&&\\
\bullet\ar@{-}[uuurr]\ar@{-}[rrrr]\ar@{-}[uurrrr]&&&&\bullet\ar@{-}[rrrr]\ar@{-}[uu]&&&&\bullet\ar@{-}[uuull]\ar@{-}[uullll]
}$$

Here the circle in the middle is by definition a line, and with this convention, the basic axioms of elementary geometry are satisfied, in the sense that any two points determine a line, and any two lines determine a point. Which is something really beautiful.

\bigskip

Now observe that the sets $X,B$ of points and lines of the Fano plane form a $(7,3,1)$ block design, corresponding to the following filled $(1,2,2)$ pattern:
$$I_7=\begin{pmatrix}
x&x&y&y&y&x&y\\
y&x&x&y&y&y&x\\
x&y&x&x&y&y&y\\
y&x&y&x&x&y&y\\
y&y&x&y&x&x&y\\
y&y&y&x&y&x&x\\
x&y&y&y&x&y&x
\end{pmatrix}$$

In order to construct now more general examples, along the same lines, observe that the Fano plane is the projective plane over the finite field $\mathbb F_2=\{0,1\}$. The same method works with $\mathbb F_2$ replaced by an arbitrary finite field $\mathbb F_q$, and we have:

\begin{proposition}
Assume that $q=p^k$ is a prime power. Then the point-line incidence matrix of the projective plane over $\mathbb F_q$ is a $(1,q,q^2-q)$ pattern.
\end{proposition}

\begin{proof}
The sets $X,B$ of points and lines of the projective plane over $\mathbb F_q$ are indeed known to form a $(q^2+q+1,q+1,1)$ block design, and this gives the result.
\end{proof}

\index{Paley biplane}

There are many other interesting examples of block designs giving rise to patterns, via Proposition 3.27. For instance the Paley biplane, which is a famous object in combinatorics, is a $(11,5,2)$ block design, giving rise to a $(2,3,3)$ pattern. See \cite{bnz}.

\bigskip

Let us discuss now the problem of associating real values to the symbols $x,y$ in an $(a,b,c)$ pattern such that the resulting matrix $U(x,y)$ is orthogonal. We have:

\begin{proposition}
Given $a,b,c\in\mathbb N$, there exists an orthogonal matrix having pattern $(a,b,c)$ iff $b^2\geq ac$. In this case the solutions are $U(x,y)$ and $-U(x,y)$, where
$$x=-\frac{t}{\sqrt{b}(t+1)}\quad,\quad
y=\frac{1}{\sqrt{b}(t+1)}$$
with $t=(b\pm\sqrt{b^2-ac})/a$ being one of the solutions of $at^2-2bt+c=0$.
\end{proposition}

\begin{proof}
Consider a filled $(a,b,c)$ pattern $U\in M_N(x,y)$, as in Definition 3.25. In order for this matrix $U$ to be orthogonal, the following conditions must be satisfied:
$$ax^2+2bxy+cy^2=0$$
$$(a+b)x^2+(b+c)y^2=1$$

The first condition, coming from the orthogonality of rows, tells us that $t=-x/y$ must be the variable in the statement. As for the second condition, this becomes:
\begin{eqnarray*}
y^2
&=&\frac{1}{(a+b)t^2+(b+c)}\\
&=&\frac{1}{(at^2+c)+(bt^2+b)}\\
&=&\frac{1}{2bt+bt^2+b}\\
&=&\frac{1}{b(t+1)^2}
\end{eqnarray*}

This gives the above formula of $y$, and hence the formula of $x=-ty$ as well.
\end{proof}

Next in line, following \cite{bn3}, \cite{bnz}, we have the following result:

\begin{proposition}
Let $U=U(x,y)$ be orthogonal, corresponding to an $(a,b,c)$ pattern. Then $H=\sqrt{N}U$ is almost Hadamard if:
$$(N(a-b)+2b)|x|+(N(c-b)+2b)|y|\geq 0$$
\end{proposition}

\begin{proof}
Let $S_{ij}={\rm sgn}(U_{ij})$. Since any row of $U$ consists of $a+b$ copies of $x$ and $b+c$ copies of $y$, we have:
$$(SU^t)_{ii}=\sum_k{\rm sgn}(U_{ik})U_{ik}=(a+b)|x|+(b+c)|y|$$

Regarding now $(SU^t)_{ij}$ with $i\neq j$, we can assume in the computation that the $i$-th and $j$-th row of $U$ are exactly those pictured in Definition 3.25. Thus:
\begin{eqnarray*}
(SU^t)_{ij}
&=&\sum_k{\rm sgn}(U_{ik})U_{jk}\\
&=&a\,{\rm sgn}(x)x+b\,{\rm sgn}(x)y+b\,{\rm sgn}(y)x+c\,{\rm sgn}(y)y\\
&=&a|x|-b|y|-b|x|+c|y|\\
&=&(a-b)|x|+(c-b)|y|
\end{eqnarray*}

We obtain the following formula for the matrix $SU^t$ itself, with $J_N=\mathbb I_N/N$:
\begin{eqnarray*}
SU^t
&=&2b(|x|+|y|)1_N+((a-b)|x|+(c-b)|y|)NJ_N\\
&=&2b(|x|+|y|)(1_N-J_N)+((N(a-b)+2b)|x|+(N(c-b)+2b)|y|))J_N
\end{eqnarray*}

Now since the matrices $1_N-J_N,J_N$ are orthogonal projections, we have $SU^t>0$ if and only if the coefficients of these matrices in the above expression are both positive. Since the coefficient of $1_N-J_N$ is clearly positive, the condition left is:
$$(N(a-b)+2b)|x|+(N(c-b)+2b)|y|\geq 0$$

So, we have obtained the condition in the statement, and we are done.
\end{proof}

Once again following \cite{bn3}, \cite{bnz}, we have the following result:

\begin{theorem}
Assume that $a,b,c\in\mathbb N$ satisfy $c\geq a$ and $b(b-1)=ac$, and consider the $(a,b,c)$ pattern $U=U(x,y)$, where:
$$x=\frac{a+(1-a-b)\sqrt{b}}{Na}\quad,\quad 
y=\frac{b+(a+b)\sqrt{b}}{Nb}$$
Then $H=\sqrt{N}U$ is an almost Hadamard matrix. 
\end{theorem}

\begin{proof}
We have $b^2-ac=b$, so Proposition 3.30 applies, and shows that with $t=(b-\sqrt{b})/a$ we have an orthogonal matrix $U=U(x,y)$, where:
$$x=-\frac{t}{\sqrt{b}(t+1)}\quad,\quad
y=\frac{1}{\sqrt{b}(t+1)}$$

But this gives the formulae of $x,y$ in the statement. Now, observe that we have:
\begin{eqnarray*}
N(a-b)+2b
&=&(a+2b+c)(a-b)+2b\\
&=&a^2+ab-2b^2+ac-bc+2b\\
&=&a^2+ab-ac-bc\\
&=&(a-c)(a+b)
\end{eqnarray*}

Similarly, we have the following formula:
$$N(c-b)+2b=(c-a)(c+b)$$

Thus the quantity in Proposition 3.30 is $Ky$, with:
\begin{eqnarray*}
K
&=&(a-c)(a+b)t+(c-a)(c+b)\\
&=&(c-a)(c+b-(a+b)t)\\
&=&\frac{c-a}{a}(ac+ab-(a+b)(b-\sqrt{b}))\\
&=&\frac{c-a}{a}((ac-b^2)+(a+b)\sqrt{b})\\
&=&\frac{c-a}{a}((a+b)\sqrt{b}-b)
\end{eqnarray*}

Since this quantity is positive, Proposition 3.30 applies and gives the result.
\end{proof}

As a main application, we have the following result, also from \cite{bn3}, \cite{bnz}:

\index{projective plane}
\index{finite field}

\begin{theorem}
Assume that $q=p^k$ is a prime power. Then the matrix $I_N\in M_N(x,y)$, where $N=q^2+q+1$ and
$$x=\frac{1-q\sqrt{q}}{\sqrt{N}}\quad,\quad
y=\frac{q+(q+1)\sqrt{q}}{q\sqrt{N}}$$
having $(1,q,q^2-q)$ pattern coming from the point-line incidence of the projective plane over $\mathbb F_q$ is an almost Hadamard matrix.
\end{theorem}

\begin{proof}
Indeed, the conditions $c\geq a$ and $b(b-1)=ac$ in Theorem 3.31 are satisfied, and the variables constructed there are $x'=x/\sqrt{N}$ and $y'=y/\sqrt{N}$.
\end{proof}

We refer to \cite{bn3}, \cite{bnz} for more on such matrices, including examples and norm numerics, in relation with the optimization question for the 1-norm. In what concerns us, we will be back to this in chapter 12 below, with a similar discussion in the complex case.

\section*{3e. Exercises}

There are many interesting questions in relation with the above, and especially with the circulant matrices, and the block designs. Let us start with:

\begin{exercise}
Work out the formula of the basic circulant almost Hadamard matrix
$$L_N=\frac{1}{\sqrt{N}}
\begin{pmatrix}
1&-\cos^{-1}\frac{\pi}{N}&\cos^{-1}\frac{2\pi}{N}&\ldots\ldots&\cos^{-1}\frac{(N-1)\pi}{N}\\
\cos^{-1}\frac{(N-1)\pi}{N}&1&-\cos^{-1}\frac{\pi}{N}&\ldots\ldots&-\cos^{-1}\frac{(N-2)\pi}{N}\\
\vdots&\vdots&\vdots&&\vdots\\
\vdots&\vdots&\vdots&&\vdots\\
-\cos^{-1}\frac{\pi}{N}&\cos^{-1}\frac{2\pi}{N}&-\cos^{-1}\frac{3\pi}{N}&\ldots\ldots&1
\end{pmatrix}$$
at $N=3,5,7,9,11$, and compute its $1$-norm.
\end{exercise}

The interest in these computations comes from the fact that $L_N$ is believed to be optimal in many cases, although there is no known proof for this.

\begin{exercise}
Compute the almost Hadamard matrix associated to the Fano plane,
$$I_7=\begin{pmatrix}
x&x&y&y&y&x&y\\
y&x&x&y&y&y&x\\
x&y&x&x&y&y&y\\
y&x&y&x&x&y&y\\
y&y&x&y&x&x&y\\
y&y&y&x&y&x&x\\
x&y&y&y&x&y&x
\end{pmatrix}$$
and its $1$-norm. Then do the same with the Paley biplane.
\end{exercise}

Here the picture of the Paley biplane can be found of course with an internet search. As a bonus exercise, try to find out if these almost Hadamard matrices are optimal.

\begin{exercise}
Draw the projective planes over $\mathbb F_q$ with $q=p^k$ small, and compute the associated almost Hadamard matrices, and their $1$-norm.
\end{exercise}

Here we have chosen not to give a precise bound for $q$. The more, the better.

\chapter{Partial matrices}

\section*{4a. Partial matrices}

In this chapter we discuss a number of more specialized questions in the real case, regarding the square or rectangular submatrices of the Hadamard matrices $H\in M_N(\pm1)$, and some related classes of square or rectangular real matrices. There are many things to be done here, going in various directions, and our plan will be as follows:

\bigskip

(1) We will first review the material from chapter 1 regarding the partial Hadamard matrices, with some further algebraic results, and with a few analytic things added too, inspired from the theory developed in the square matrix case in chapters 2-3.

\bigskip

(2) Then, we will get into the question of counting the partial Hadamard matrices $H\in M_{M\times N}(\pm1)$, at small values of $M$, and with $N\to\infty$. This is a question  having no square counterpart, and following de Launey-Levin \cite{dle}, interesting things can be said.

\bigskip

(3) Finally, we will go back to the square matrix case, and present some results from \cite{bs2} regarding the square submatrices of the usual Hadamard matrices $H\in M_N(\pm1)$, making the connection with the almost Hadamard matrices from chapter 3.

\bigskip

All in all, many things to be done. Let us mention right away that the most important thing in all this is (2), with the counting result of de Launey and Levin in \cite{dle} being something truly remarkable, and providing a viable alternative to the whole HC problematics, developed by countless people since the papers of Sylvester \cite{syl} and Hadamard \cite{had}. 

\bigskip

Getting started now, let us begin by reviewing what we know about the partial Hadamard matrices, from chapter 1. The definition of these matrices is as follows:

\index{PHM}
\index{partial Hadamard matrix}

\begin{definition}
A partial Hadamard matrix (PHM) is a rectangular matrix 
$$H\in M_{M\times N}(\pm1)$$
whose rows are pairwise orthogonal, with respect to the scalar product of $\mathbb R^N$. 
\end{definition}

The motivating examples are the usual Hadamard matrices $H\in M_N(\pm1)$, and their various $M\times N$ submatrices, with $M\leq N$. However, there are as well examples which are not of this form, and the PHM are interesting combinatorial objects, on their own. 

\bigskip

Following the study from the square matrix case, we first have:

\begin{proposition}
The set $Y_{M,N}$ of the $M\times N$ partial Hadamard matrices is
$$Y_{M,N}=M_{M\times N}(\pm1)\cap\sqrt{N}O_{M,N}$$
where $O_{M,N}$ is the following space of rectangular matrices:
$$O_{M,N}=\left\{U\in M_{M\times N}(\mathbb R)\Big|UU^t=1_M\right\}$$
At $M=N$, we recover in this way the previous formula $Y_N=M_N(\pm1)\cap\sqrt{N}O_N$.
\end{proposition}

\begin{proof}
This follows exactly as in the square matrix case. Indeed, given a rectangular matrix $U\in M_{M\times N}(\mathbb R)$ having rows $R_1,\ldots,R_M\in\mathbb R^N$, we have:
$$(UU^t)_{ij}=\sum_kU_{ik}U_{jk}=<R_i,R_j>$$

Thus, the condition $UU^t=1_M$ expresses the fact that the vectors $R_1,\ldots,R_M$  are pairwise orthogonal, and of norm 1, and this gives the formula in the statement.
\end{proof}

As a remark here, at $M=1$ we have of course $Y_{1,N}=M_{1\times N}(\pm1)$, and this because of an automatic inclusion $M_{1\times N}(\pm1)\subset\sqrt{N}O_{1,N}$. Indeed, given $H\in M_{1\times N}(\pm1)$, the matrix $U=H/\sqrt{N}$ satisfies $UU^t=\frac{1}{N}\cdot N=1$, and so we have, as claimed:
$$U\in\sqrt{N}O_{1,N}$$

In general, the space $O_{M,N}$ appearing above can be thought of as being a joint generalization of the unit sphere $S^{N-1}$, which appears at $M=1$, and of the orthogonal group $O_N$, which appears in the square case, $M=N$. Based on this analogy, the space $O_{M,N}$ has several useful interpretations, which can be summarized as follows:

\index{homogeneous space}

\begin{proposition}
The space $O_{M,N}$ has the following properties:
\begin{enumerate}
\item Its elements are the transposes of the isometries $g:\mathbb R^M\to\mathbb R^N$.

\item It is the space of vectors $R_1,\ldots,R_M\in S^{N-1}$ which are pairwise orthogonal.

\item It is also an homogeneous space, given by $O_{M,N}\simeq O_N/O_{N-M}$.

\item It is also the space determined by the first $M$ rows of coordinates on $O_N$.
\end{enumerate}
\end{proposition}

\begin{proof}
All this is standard algebra and geometry, the idea being as follows:

\medskip

(1) Each matrix $U\in M_{M\times N}(\mathbb R)$ determines a linear map $f:\mathbb R^N\to\mathbb R^M$, given by $f(x)=Ux$, whose transpose is the linear map $g:\mathbb R^M\to\mathbb R^N$ given by $g(x)=U^tx$. Now observe that for any two vectors $x,y\in\mathbb R^M$ we have:
$$<g(x),g(y)>=<U^tx,U^ty>=<x,UU^ty>$$

Thus the condition $UU^t=1$ is equivalent to the following condition:
$$<g(x),g(y)>=<x,y>$$

But this latter condition tells us that $g$ must be an isometry, as desired.

\medskip

(2) This follows from the fact, that we know from the proof of Proposition 4.2, that the condition $UU^t=1_M$ tells us that the row vectors $R_1,\ldots,R_M\in\mathbb R^N$ of our matrix $U\in M_{M\times N}(\mathbb R)$ must be pairwise orthogonal, and of norm 1.

\medskip

(3) Since the condition $UU^t=1$ defining $O_{M,N}$ implies $(UA^t)(UA^t)^t=1$, for any orthogonal matrix $A\in O_N$, we have an action, as follows:
$$O_N\curvearrowright O_{M,N}\quad,\quad 
A\to[U\to UA^t]$$

Let us compute now the stabilizer of the following particular element:
$$U=\begin{pmatrix}
1&&0&0&\ldots&0\\
&\ddots\\
0&&1&0&\ldots&0
\end{pmatrix}$$

Given an orthogonal matrix $A\in O_N$, we have the following formula:
$$UA^t=\begin{pmatrix}
A_{11}&\ldots&A_{N1}\\
\vdots&&\vdots\\
A_{1M}&\ldots&A_{NM}
\end{pmatrix}$$

Thus $U=UA^t$ means that the matrix $A^t\in O_N$ must be of the following form:
$$A^t=\begin{pmatrix}
1_M&0\\
*&*
\end{pmatrix}$$

Now since $A^t$ is orthogonal, it must be of the following form, with $B\in O_{N-M}$:
$$A^t=\begin{pmatrix}
1_M&0\\
0&B^t
\end{pmatrix}$$

Thus the stabilizer is $O_{N-M}$, and we obtain $O_{M,N}\simeq O_N/O_{N-M}$.

\medskip

(4) This follows from some basic functional analysis. Consider indeed the algebra $C(O_N)$ of continuous functions $f:O_N\to\mathbb C$. By Stone-Weierstrass, this algebra is generated by the coordinate functions $u_{ij}:O_N\to\mathbb C$, which are given by:
$$u_{ij}(U)=U_{ij}$$

Consider now the following closed subalgebra of the algebra $C(O_N)$:
$$A=\left<u_{ij}\Big|i=1,\ldots,M,j=1,\ldots,N\right>$$

We have then $A\simeq C(O_{M,N})$, coming from the homogeneous space result in (3).
\end{proof}

Let us discuss now, as a continuation of the study from the real case, some basic analytic aspects. In what regards the 1-norm bound, we have the following result:

\index{norm maximizer}

\begin{theorem}
Given a matrix $U\in O_{M,N}$ we have
$$||U||_1\leq M\sqrt{N}$$
with equality precisely when $H=\sqrt{N}U$ is partial Hadamard.
\end{theorem}

\begin{proof}
We have indeed the following estimate, valid for any $U\in O_{M,N}$:
\begin{eqnarray*}
||U||_1
&=&\sum_{ij}|U_{ij}|\\
&\leq&\sqrt{MN}\left(\sum_{ij}|U_{ij}|^2\right)^{1/2}\\
&=&M\sqrt{N}
\end{eqnarray*}

In this estimate the equality case holds when $|U_{ij}|=1/\sqrt{N}$ for any $i,j$. But this amounts in saying that the rescaled matrix $H=\sqrt{N}U$ must satisfy $H\in M_{M\times N}(\pm1)$, and so that this rescaled matrix must be partial Hadamard, as claimed.
\end{proof}

Observe that in terms of the rescaled matrix $H\in\sqrt{N}O_{M,N}$, the inequality found above reformulates as $||H||_1\leq MN$, with equality precisely when $H$ is partial Hadamard. Thus, in analogy with the square matrix case, we can formulate:

\index{almost PHM}

\begin{definition}
A matrix $H\in\sqrt{N}O_{M,N}$ is called:
\begin{enumerate}
\item Almost PHM, when it locally maximizes the $1$-norm on $\sqrt{N}O_{M,N}$.

\item Optimal almost PHM, when it maximizes the $1$-norm on $\sqrt{N}O_{M,N}$.
\end{enumerate} 
\end{definition} 

Some similar estimates hold for the $p$-norms, with $p\neq2$. The whole subject, while being potentially quite interesting, is for the moment largely unexplored. So, let us turn instead to algebra. Still following the study from the square case, let us formulate:

\index{dephased PHM}
\index{equivalent PHM}
\index{standard form}

\begin{definition}
Two PHM are called equivalent when we can pass from one to the other by permuting the rows or columns, or multiplying rows or columns by $-1$. Also:
\begin{enumerate}
\item We say that a PHM is in dephased form when its first row and its first column consist of $1$ entries.

\item We say that a PHM is in standard form when it is dephased, with the $1$ entries moved to the left as much as possible, by proceeding from top to bottom. 
\end{enumerate}
\end{definition}

Unlike in the square case, where the standard form is generally not used, putting a rectangular matrix in standard form is something quite useful, in practice. As an illustration here, here is a result that we already know, from chapter 1, regarding the partial Hadamard matrices in standard form, at small values of $M$:

\begin{proposition}
The standard form of dephased PHM at $M=2,3,4$ is
$$H=\begin{pmatrix}+&+\\\underbrace{+}_{N/2}&\underbrace{-}_{N/2}\end{pmatrix}$$
$$H=\begin{pmatrix}+&+&+&+\\+&+&-&-\\\underbrace{+}_{N/4}&\underbrace{-}_{N/4}&\underbrace{+}_{N/4}&\underbrace{-}_{N/4}\end{pmatrix}$$
$$H=\begin{pmatrix}
+&+&+&+&+&+&+&+\\
+&+&+&+&-&-&-&-\\
+&+&-&-&+&+&-&-\\
\underbrace{+}_a&\underbrace{-}_b&\underbrace{+}_b&\underbrace{-}_a&\underbrace{+}_b&\underbrace{-}_a&\underbrace{+}_a&\underbrace{-}_b
\end{pmatrix}$$
where the numbers $a,b\in\mathbb N$ satisfy $a+b=N/4$.
\end{proposition}

\begin{proof}
This is something that we know from chapter 1, the idea being that the $M=2$ result is obvious, that the $M=3$ result follows from the orthogonality conditions between the rows, and that the $M=4$ result follows from the $M=3$ result.
\end{proof}

At $M=5$ and higher the situation is more complicated, and we will be back to this. For the moment, let us stay with $M=4$. We can fine-tune our result, as follows:

\begin{theorem}
The $4\times N$ partial Hadamard matrices are of the form
$$H=(\underbrace{W_4\ \ldots\ W_4}_a\ \underbrace{K_4\ \ldots\ K_4}_b)$$
with $a+b=N/4$. Moreover, we can assume $a\geq b$.
\end{theorem}

\begin{proof}
Let $H\in M_{4\times N}(\pm1)$ be as in Proposition 4.7. The matrix formed by the $a$ type columns, one from each block, is equivalent to $W_4$, via a permutation of columns:
$$\begin{pmatrix}
+&+&+&+\\
+&+&-&-\\
+&-&+&-\\
+&-&-&+
\end{pmatrix}\sim W_4$$

Also, the matrix formed by the $b$ type columns, one from each block, is equivalent to $K_4$, via a first column sign switch, plus a certain permutation of the columns:
$$\begin{pmatrix}
+&+&+&+\\
+&+&-&-\\
+&-&+&-\\
-&+&+&-
\end{pmatrix}\sim K_4$$

Thus, just by performing operations on the columns, we obtain, as desired:
$$H\sim(\underbrace{W_4\ \ldots\ W_4}_a\ \underbrace{K_4\ \ldots\ K_4}_b)$$

In order to prove now the last assertion, we must prove that we have:
$$(\underbrace{W_4\ \ldots\ W_4}_a\ \underbrace{K_4\ \ldots\ K_4}_b)\sim(\underbrace{K_4\ \ldots\ K_4}_a\ \underbrace{W_4\ \ldots\ W_4}_b)$$

But this can be seen by performing a sign switch on the last row, and then permuting the columns. Equivalently, we can start with the original matrix, in standard form, and perform a sign switch on the last row. The matrix becomes:
$$H\sim\begin{pmatrix}
+&+&+&+&+&+&+&+\\
+&+&+&+&-&-&-&-\\
+&+&-&-&+&+&-&-\\
\underbrace{-}_a&\underbrace{+}_b&\underbrace{-}_b&\underbrace{+}_a&\underbrace{-}_b&\underbrace{+}_a&\underbrace{-}_a&\underbrace{+}_b
\end{pmatrix}$$

Now by putting this matrix in standard form, we obtain:
$$H=\begin{pmatrix}
+&+&+&+&+&+&+&+\\
+&+&+&+&-&-&-&-\\
+&+&-&-&+&+&-&-\\
\underbrace{+}_b&\underbrace{-}_a&\underbrace{+}_a&\underbrace{-}_b&\underbrace{+}_a&\underbrace{-}_b&\underbrace{+}_b&\underbrace{-}_a
\end{pmatrix}$$

Thus $a,b$ got interchanged, and this gives the result.
\end{proof}

At $M=5$ now, as already mentioned above, the combinatorics becomes quite complicated, and we will see in a moment that there are $5\times N$ partial Hadamard matrices which do not complete into Hadamard matrices. We first have the following result:

\begin{proposition}
The $5\times N$ partial Hadamard matrices are of the form
$$H=\begin{pmatrix}
W_4&\ldots&W_4&&K_4&\ldots&K_4\\
v_1&\ldots&v_a&&x_1&\ldots&x_b
\end{pmatrix}$$
with $a\geq b$, $a+b=N/4$ and with $v_i,x_j\in(\pm1)^4$ satisfying
$$W_4\begin{pmatrix}r_1\\r_2\\r_3\\r_4\end{pmatrix}
=-K_4\begin{pmatrix}s_1\\s_2\\s_3\\s_4\end{pmatrix}$$
where $r_t=\sum_i(v_i)_t$ and $s_t=\sum_j(v_j)_t$.
\end{proposition}

\begin{proof}
This is something that we already worked out at $N=8$, in chapter 1, in both of the cases that can appear, namely $a=2,b=0$ and $a=1,b=1$. The proof in general is similar, via some routine computations, with the equations in the statement coming by processing the orthogonality conditions between the 5th row and the first 4 rows.
\end{proof}

As a first observation, the equations in the above statement can be written in the following more convenient form:
$$K_4^{-1}W_4\begin{pmatrix}r_1\\r_2\\r_3\\r_4\end{pmatrix}=-\begin{pmatrix}s_1\\s_2\\s_3\\s_4\end{pmatrix}$$

Now observe that the matrix of this system is as follows:
$$K_4^{-1}W_4
=\frac{1}{2}\begin{pmatrix}
-&+&+&+\\
-&-&+&-\\
-&+&-&-\\
-&-&-&+
\end{pmatrix}$$

Thus, the system can be written as follows:
$$\begin{pmatrix}
-&+&+&+\\
-&-&+&-\\
-&+&-&-\\
-&-&-&+
\end{pmatrix}
\begin{pmatrix}r_1\\r_2\\r_3\\r_4\end{pmatrix}
=-2\begin{pmatrix}s_1\\s_2\\s_3\\s_4\end{pmatrix}$$

Thus, we are led into parity and positivity questions, regarding the vectors $r_t=\sum_i(v_i)_t$ and $s_t=\sum_j(v_j)_t$. It is possible to further go along these lines, but the structure of the $5\times N$ partial Hadamard matrices remains something quite complicated. As an explicit consequence of our study, however, we have the following result:

\begin{theorem}
Consider an arbitrary $4\times N$ partial Hadamard matrix, written as
$$H=(\underbrace{W_4\ \ldots\ W_4}_a\ \underbrace{K_4\ \ldots\ K_4}_b)$$
with $a\geq b$, $a+b=N/4$, up to equivalence. In order for this matrix to complete into a $5\times N$ partial Hadamard matrix, the following condition must be satisfied:
$$ab=0\implies N=0(8)$$
In particular, the following $4\times N$ partial Hadamard matrix,
$$Z=(W_4\ W_4\ W_4)$$
does not complete into a $5\times N$ partial Hadamard matrix.
\end{theorem}

\begin{proof}
This follows from Proposition 4.9, because with the notations there, the condition $b=0$ implies that the system there is simply:
$$W_4\begin{pmatrix}r_1\\r_2\\r_3\\r_4\end{pmatrix}=0$$

Since $W_4$ is invertible, the solution is $r=0$. Now observe that, by definition of the numbers $r_i$, we have $r_i=a(2)$ for any $i$. Thus, we must have $a=0(2)$, and since we have $a=N/4$, this gives $N=0(8)$, as desired. The proof in the case $a=0$ is similar.
\end{proof}

In general, the full classification of all the possible $5\times 8$ completions of a given  $4\times N$ partial Hadamard matrix is something quite difficult, and we have already seen this at $N=8$, where a careful study is needed, the result being as follows:

\begin{theorem}
The two $4\times 8$ partial Hadamard matrices, namely
$$A=(W_4\ W_4)\quad,\quad 
B=(W_4\ K_4)$$
both complete into $5\times 8$ partial Hadamard matrices, with the solutions being those coming from the lower rows of the following matrices, which are Hadamard:
$$\begin{pmatrix}W_4&W_4\\ W_4&-W_4\end{pmatrix}
\quad,\quad
\begin{pmatrix}W_4&W_4\\ K_4&-K_4\end{pmatrix}\quad,\quad 
\begin{pmatrix}W_4&K_4\\ W_4&-K_4\end{pmatrix}
\quad,\quad
\begin{pmatrix}W_4&K_4\\ K_4&-W_4\end{pmatrix}$$
This gives as well the higher completions, $M\times 8$ with $M=6,7,8$.
\end{theorem}

\begin{proof}
This is something that we already know, from chapter 1.
\end{proof}

At $N=12$ now, we have only one matrix to be studied, which is as follows, and with at least 8 solutions to the completion problem, coming from the Paley matrix $P_{12}$:
$$P=(W_4\ W_4\ K_4)$$

Generally speaking, all this leads to quite complicated algebra and combinatorics. We refer to Hall \cite{hal}, Ito \cite{ito} and Verheiden \cite{ver} for more on the combinatorics of the PHM. Finally, let us end this discussion with an elementary result, from \cite{bsk}:

\begin{theorem}
For a partial Hadamard matrix $H\in M_{(N-1)\times N}(\pm1)$, with rows $R_1,\ldots,R_{N-1}$ and columns $C_1,\ldots,C_N$, the following are equivalent:
\begin{enumerate}
\item $H$ is completable into a $N\times N$ Hadamard matrix.

\item $|\det H^{(j)}|$ is independent from $j$, where $H^{(j)}$ is obtained from $H$ by removing $C_j$.

\item $|\det H^{(j)}|=N^{N/2-1}$ for any $i$, where $H^{(j)}$ is as above.
\end{enumerate}
Moreover, if these conditions hold, the completion is obtained by setting
$$H_{Nj} = (-1)^{j}N^{1-N/2}\,\det H^{(j)}$$
with $H^{(j)}$ being as above, obtained from $H$ by removing the column $C_j$.
\end{theorem}

\begin{proof}
This follows from some basic linear algebra, the idea being as follows:

\medskip

$(1)\iff(2)$. Consider the following vector, having integer entries:
$$Z_j=(-1)^j\,\det H^{(j)}$$

Our claim is that we have the following equality of vector spaces:
$$span(R_1,\ldots,R_{N-1})^\perp=\{\lambda Z|\lambda \in \mathbb R\}$$

Indeed, if we denote by $H_i$ the square matrix obtained from $H$ by adding a first row equal to $R_i$, then we have the following computation, which proves our claim:
\begin{eqnarray*}
<R_i,Z>
&=&\sum_jH_{ij}Z_j\\
&=&\sum_j(-1)^jH_{ij}\det H^{(j)}\\
&=&\det H_i\\
&=&0
\end{eqnarray*}

But this gives $(1)\iff(2)$, since the existence of a completion is equivalent to the fact that $span(R_1,\ldots,R_{N-1})^\perp$ contains a vector with all entries having absolute value $1$. 

\medskip

$(1)\implies(3)$. Write $c=|\det H^{(j)}|$ and let $M\in M_N(\pm1)$ be the Hadamard matrix completing $H$. The proof of $(1)\iff(2)$ above shows that the last row of $M$ must be the vector $c^{-1} Z$. Also, since the matrix $M\in M_N(\pm1)$ is Hadamard, we have:
$$|\det M|=N^{N/2}$$

Thus, it remains to compute this determinant by expansion with respect to the last row, and the computation here gives:
\begin{eqnarray*}
\det M
&=&\sum_{j=1}^N c^{-1} (-1)^{N+j} (-1)^j\det H^{(j)}\cdot\det H^{(j)}\\
&=&(-1)^N c N
\end{eqnarray*}

But this means that we have $c= N^{N/2-1}$, which proves the implication $(1)\implies(3)$, and also proves the last assertion of our theorem.

\medskip

$(3)\implies(2)$. This is something obvious, and so we are done.
\end{proof}

We will be back to the algebraic properties of the PHM on several occasions in this book, but directly in the complex matrix case, or sometimes in the general root of unity case, where more things can be said. In relation with the real case, of particular interest will be the material in chapter 15 below, where, following \cite{bsk}, we will associate a quantum semigroup of partial permutations of $\{1,\ldots,N\}$ to each such matrix, real or complex. 

\section*{4b. Counting results}

Let us try now to count the partial Hadamard matrices $H\in M_{M\times N}(\pm1)$. This is an easy task at $M=2,3,4$, where the answer is as follows:

\index{multinomial coefficient}

\begin{proposition}
The number of PHM at $M=2,3,4$ is
\begin{eqnarray*}
\#PHM_{2\times N}&=&2^N\binom{N}{N/2}\\
\#PHM_{3\times N}&=&2^N\binom{N}{N/4,N/4,N/4,N/4}\\
\#PHM_{4\times N}&=&2^N\sum_{a+b=N/4}\binom{N}{a,b,b,a,b,a,a,b}
\end{eqnarray*}
with the quantities on the right being multinomial coefficients.
\end{proposition}

\begin{proof}
We use the structure results for the PHM in standard form at $M\leq4$ found above, which are as follows, with the numbers $a,b\in\mathbb N$ satisfing $a+b=N/4$:
$$H=\begin{pmatrix}+&+\\\underbrace{+}_{N/2}&\underbrace{-}_{N/2}\end{pmatrix}$$
$$H=\begin{pmatrix}+&+&+&+\\+&+&-&-\\\underbrace{+}_{N/4}&\underbrace{-}_{N/4}&\underbrace{+}_{N/4}&\underbrace{-}_{N/4}\end{pmatrix}$$
$$H=\begin{pmatrix}
+&+&+&+&+&+&+&+\\
+&+&+&+&-&-&-&-\\
+&+&-&-&+&+&-&-\\
\underbrace{+}_a&\underbrace{-}_b&\underbrace{+}_b&\underbrace{-}_a&\underbrace{+}_b&\underbrace{-}_a&\underbrace{+}_a&\underbrace{-}_b
\end{pmatrix}$$

But this gives the formulae in the statement, with the multinomial coefficients counting the matrices having the first row consisting of 1 entries only, obtained by permuting the columns of the above solutions, and with the $2^N$ factors coming from this.
\end{proof}

In order to convert the above result into $N\to\infty$ estimates, we will need the following technical result regarding the multinomial coefficients, from Richmond-Shallit \cite{rsh}:

\begin{theorem}
We have the estimate
$$\sum_{n_1+\ldots+n_k=N}\binom{N}{n_1,\ldots,n_k}^p
\simeq k^{pN}\sqrt{\frac{k^{k(p-1)}}{p^{k-1}(2\pi N)^{(k-1)(p-1)}}}$$
in the $N\to\infty$ limit.
\end{theorem}

\index{Richmond-Shallit}

\begin{proof}
This is proved by Richmond and Shallit in \cite{rsh} at $p=2$, and the proof in the general case, $p\in\mathbb N$, is similar, the idea being as follows:

\medskip

(1) In order to do some analysis, we agree to use the convention $x!=\Gamma(x+1)$ for $x>0$ real. Since the multinomial coefficient in the statement attains its maximum when the numbers $n_i$ are all equal, it is natural to make a change of variables, as follows:
$$n_i=\frac{N}{k}+\sqrt{N}x_i$$

Observe that, since we have $n_1+\ldots+n_k=N$, the numbers $x_1,\ldots,x_k$ satisfy:
$$x_1+\ldots+x_k=0$$

(2) Let us first estimate, in terms of the numbers $x_1,\ldots,x_k$, the multinomial coefficient in the statement. By using the Taylor formula $\log(1+y)\simeq y-y^2/2$, we obtain:
\begin{eqnarray*}
\log n_i
&=&\log\left(\frac{N}{k}\left(1+\frac{kx_i}{\sqrt{N}}\right)\right)\\
&\simeq&\log\frac{N}{k}+\frac{kx_i}{\sqrt{N}}-\frac{k^2x_i^2}{2N}
\end{eqnarray*}

By multiplying by $n_i$, this gives the following estimate:
\begin{eqnarray*}
n_i\log n_i
&\simeq&\left(\frac{N}{k}+\sqrt{N}x_i\right)\log\frac{N}{k}+\left(\frac{N}{k}+\sqrt{N}x_i\right)\left(\frac{kx_i}{\sqrt{N}}-\frac{k^2x_i^2}{2N}\right)\\
&\simeq&\left(\frac{N}{k}+\sqrt{N}x_i\right)\log\frac{N}{k}+\sqrt{N}x_i+\frac{kx_i^2}{2}
\end{eqnarray*}

Now by further substracting $n_i$, we obtain the following estimate:
$$n_i\log n_i-n_i
\simeq\left(\frac{N}{k}+\sqrt{N}x_i\right)\log\frac{N}{k}+\frac{kx_i^2}{2}-\frac{N}{k}$$

(3) We are now ready to estimate the multinomial coefficient in the statement. By summing over $i$, and using $x_1+\ldots+x_k=0$, the formula found above gives:
$$\sum_in_i\log n_i-n_i\simeq N\log\frac{N}{k}-N+\frac{k}{2}\sum_ix_i^2$$

By using the Stirling formula $n!\simeq e^{n\log n-n}\sqrt{2\pi n}$, we obtain from this:
\begin{eqnarray*}
n_1!\ldots n_k!
&\simeq&\exp\left(\sum_in_i\log n_i-n_i\right)\sqrt{2\pi n_1}\ldots\sqrt{2\pi n_k}\\
&\simeq&\exp\left(N\log\frac{N}{k}-N+\frac{k}{2}\sum_ix_i^2\right)\left(\frac{2\pi N}{k}\right)^{k/2}\\
&=&\left(\frac{N}{ke}\right)^N\exp\left(\frac{k}{2}\sum_ix_i^2\right)\left(\frac{2\pi N}{k}\right)^{k/2}
\end{eqnarray*}

Thus, the multinomial coefficient in the statement is:
\begin{eqnarray*}
\binom{N}{n_1,\ldots,n_k}
&\simeq&\left(\frac{N}{e}\right)^N\sqrt{2\pi N}\left(\frac{ke}{N}\right)^N\exp\left(-\frac{k}{2}\sum_ix_i^2\right)\left(\frac{2\pi N}{k}\right)^{-k/2}\\
&=&k^N\exp\left(-\frac{k}{2}\sum_ix_i^2\right)(2\pi N)^{(1-k)/2}k^{k/2}
\end{eqnarray*}

(4) Raising now to the power $p$ gives the following formula:
$$\binom{N}{n_1,\ldots,n_k}^p\simeq k^{pN}\exp\left(-\frac{kp}{2}\sum_ix_i^2\right)(2\pi N)^{(1-k)p/2}k^{kp/2}$$

Getting now to what we want to do, the point is that, by using the above estimate for the summands, we can estimate their sum by a multiple integral, as follows:
\begin{eqnarray*}
&&\sum_{n_1+\ldots+n_k=N}\binom{N}{n_1,\ldots,n_k}^p\\
&\simeq&k^{pN}(2\pi N)^{\frac{(1-k)p}{2}}k^{\frac{kp}{2}}\int_0^N\ldots\int_0^N\exp\left(-\frac{kp}{2}\sum_{i=1}^kx_i^2\right)dn_1\ldots dn_{k-1}\\
&=&k^{pN}(2\pi N)^{\frac{(1-k)p}{2}}k^{\frac{kp}{2}}N^{\frac{k-1}{2}}\\
&&\times\int_0^N\ldots\int_0^N\exp\left(-\frac{kp}{2}\sum_{i=1}^{k-1}x_i^2-\frac{kp}{2}\left(\sum_{i=1}^{k-1}x_i\right)^2\right)dx_1\ldots dx_{k-1}
\end{eqnarray*}

(5) We are almost there. By doing the calculus, as explained in \cite{rsh}, this gives:
\begin{eqnarray*}
\sum_{n_1+\ldots+n_k=N}\binom{N}{n_1,\ldots,n_k}^p
&\simeq&k^{pN}(2\pi N)^{\frac{(1-k)p}{2}}k^{\frac{kp}{2}}N^{\frac{k-1}{2}}\times\pi^{\frac{k-1}{2}}k^{-\frac{1}{2}}\left(\frac{kp}{2}\right)^{\frac{1-k}{2}}\\
&=&k^{pN}(2\pi N)^{\frac{(1-k)p}{2}}k^{\frac{kp}{2}-\frac{1}{2}+\frac{1-k}{2}}\left(\frac{p}{2\pi N}\right)^{\frac{1-k}{2}}\\
&=&k^{pN}(2\pi N)^{\frac{(1-k)(p-1)}{2}}k^{\frac{kp-k}{2}}p^{\frac{1-k}{2}}\\
&=&k^{pN}\sqrt{\frac{k^{k(p-1)}}{p^{k-1}(2\pi N)^{(k-1)(p-1)}}}
\end{eqnarray*}

Thus we have obtained the formula in the statement, and we are done.
\end{proof}

The above formula is something very useful, that we will heavily use in what follows. Getting back now to the PHM, we have the following result:

\begin{theorem}
The probability for a random $H\in M_{M\times N}(\pm1)$ to be a PHM is
$$P_2\simeq\frac{2}{\sqrt{2\pi N}}\quad,\quad 
P_3\simeq\frac{16}{\sqrt{(2\pi N)^3}}\quad,\quad 
P_4\simeq\frac{512}{(2\pi N)^3}$$
in the $N\in4\mathbb N$, $N\to\infty$ limit.
\end{theorem}

\begin{proof}
Since there are exactly $2^{MN}$ sign matrices of size $N\times M$, the probability $P_M$ for a random $H\in M_{M\times N}(\pm1)$ to be a PHM is given by:
$$P_M=\frac{1}{2^{MN}}\#PHM_{M\times N}$$

With this formula in hand, the result follows from Proposition 4.13, by using the estimates for sums of multinomial coefficients from Theorem 4.14. 
\end{proof}

\section*{4c. Asymptotic count}

In their remarkable paper \cite{dle}, de Launey and Levin were able to count the PHM, in the asymptotic limit $N\in 4\mathbb N$, $N\to\infty$. Their method is based on:

\begin{proposition}
The probability for a random $H\in M_{M\times N}(\pm1)$ to be partial Hadamard equals the probability for a length $N$ random walk with increments drawn from
$$E=\left\{(e_i\bar{e}_j)_{i<j}\Big|e\in\mathbb Z_2^M\right\}$$
regarded as a subset of $\mathbb Z_2^{\binom{M}{2}}$ to return at the origin.
\end{proposition}

\begin{proof}
Indeed, with $T(e)=(e_i\bar{e}_j)_{i<j}$, a matrix $X=[e_1,\ldots,e_N]\in M_{M\times N}(\mathbb Z_2)$ is partial Hadamard precisely when $T(e_1)+\ldots+T(e_N)=0$. But this gives the result.
\end{proof}

As explained in \cite{dle}, the above probability can be indeed computed, and we have:

\begin{theorem}
The probability for a random $H\in M_{M\times N}(\pm1)$ to be PHM is
$$P_M\simeq\frac{2^{(M-1)^2}}{\sqrt{(2\pi N)^{\binom{M}{2}}}}$$
in the $N\in 4\mathbb N$, $N\to\infty$ limit.
\end{theorem}

\index{de Launey-Levin}
\index{asymptotic count}
\index{real PHM}
\index{PHM}

\begin{proof}
According to Proposition 4.16, we have:
\begin{eqnarray*}
P_M
&=&\frac{1}{q^{(M-1)N}}\#\left\{\xi_1,\ldots,\xi_N\in E\Big|\sum_i\xi_i=0\right\}\\
&=&\frac{1}{q^{(M-1)N}}\sum_{\xi_1,\ldots,\xi_N\in E}\delta_{\Sigma\xi_i,0}
\end{eqnarray*}

By using the Fourier inversion formula we have, with $D=\binom{M}{2}$:
$$\delta_{\Sigma\xi_i,0}=\frac{1}{(2\pi)^D}\int_{[-\pi,\pi]^D}e^{i<\lambda,\Sigma\xi_i>}d\lambda$$

After many non-trivial computations, this leads to the result. See \cite{dle}.
\end{proof}

All this is quite interesting, because it provides a viable alternative to the HC problematics. To be more precise, after long decades of work on the HC, the conclusion that emerges is that this is probably an analytic question, at least in the $N>>0$ regime, with the thing to be done being that of conjecturing something of type $C_N\simeq f(N)$ about the asymptotics of the number $C_N$ of the $N\times N$ Hadamard matrices, with $f(N)$ being some kind of known function, and then proving this conjecture, with $C_N>0$ coming as consequence. But, no one knows what the conjecture of type $C_N\simeq f(N)$ should be.

\bigskip

In contrast to this, the work of de Launey and Levin \cite{dle} explained above puts us on a clear track, in order to deal with such questions. Indeed, when enlarging the attention to the partial Hadamard matrices $H\in M_{M\times N}(\pm1)$, we do have their counting result, at any $M\in\mathbb N$, and in the $N\to\infty$ limit, as a non-trivial and rock-solid starting point, and the problem is that of slowly fine-tuning their methods, as to get towards asymptotic counting results in the square matrix case, $M=N$. But this is a quite tough mix of probability and combinatorics, and no one managed so far to go beyond \cite{dle}.

\section*{4d. Square submatrices}

Following now \cite{bs2}, and some previous work of Koukouvinos, Mitrouli, Seberry  \cite{kms} and Sz\"oll\H{o}si \cite{sz2}, let us discuss now another topic, namely the square submatrices of the usual, square Hadamard matrices. We will see that all this is related, in a quite subtle way, to the notion of almost Hadamard matrix (AHM), discussed in chapter 3. Let us start with some basic linear algebra. We will need the following standard result:

\index{polar decomposition}

\begin{theorem}
Any matrix $D\in M_N(\mathbb R)$ can be written as
$$D=UT$$
with positive semidefinite $T=\sqrt{D^tD}$, and with orthogonal $U\in O_N$. Moreover:
\begin{enumerate}
\item If $D$ is invertible, then $U$ is uniquely determined, and we write:
$$U=Pol(D)$$

\item If $D=V\Delta W^t$ with $V,W$ being orthogonal and $\Delta$ being diagonal is the singular value decomposition of $D$, then $Pol(D)=VW^t$.
\end{enumerate}
\end{theorem}

\begin{proof}
All this is very standard, and can be found in any linear algebra book, one method for instance being that of deducing (2), and then the whole result, from the singular value decomposition theorem for the matrices $D\in M_N(\mathbb R)$.
\end{proof}

We start analyzing the square submatrices of the Hadamard matrices. By permuting rows and columns, we can always reduce the problem to the following situation:

\index{submatrix}

\begin{definition}
$D\in M_d(\pm1)$ is called a submatrix of $H\in M_N(\pm1)$ if we have
$$H=\begin{pmatrix}A&B\\C&D\end{pmatrix}$$
up to a permutation of the rows and columns of $H$. In this case we set:
$$r=size(A)=N-d$$
\end{definition}

Observe that any $D\in M_2(\pm1)$ having distinct columns appears as a submatrix of $W_4$, and that any $D\in M_2(\pm1)$ appears as a submatrix of $W_8$. In fact, we have:

\begin{proposition}
Let $D\in M_d(\pm1)$ be an arbitrary sign matrix.
\begin{enumerate}
\item If $D$ has distinct columns, then $D$ is as submatrix of $W_N$, with $N=2^d$.

\item In general, $D$ appears as submatrix of $W_M$, with $M=2^{d+[\log_2d]}$.
\end{enumerate}
\end{proposition}

\begin{proof}
This is something elementary, as follows:

\medskip

(1) Set $N=2^d$. If we use length $d$ bit strings $x,y\in\{0,1\}^d$ as indices, then:
$$(W_N)_{xy}=(-1)^{\sum x_iy_i}$$

Let $\widetilde{W}_N\in M_{d\times N}(\pm1)$ be the submatrix of $W_N$ having as row indices the strings of the following type:
$$x_i=(\underbrace{0\ldots 0}_i\,1\,\underbrace{0\ldots0}_{N-i-1})$$

Then for $i\in\{1,\ldots,d\}$ and $y\in\{0,1\}^d$, we have:
\vskip-3mm
$$(\widetilde{W}_N)_{iy}=(-1)^{y_i}$$

Thus the columns of $\widetilde{W}_N$ are the $N$ elements of $\{\pm 1\}^d$, which gives the result. 

\medskip

(2) Set $R=2^{[\log_2d]}\geq d$. Since the first row of $W_R$ contains only ones, $W_R\otimes W_N$ contains as a submatrix $R$ copies of $\widetilde{W}_N$, in which $D$ can be embedded, as desired.
\end{proof}

Let us go back now to Definition 4.19, and try to relate the matrices $A,D$ appearing there. The following result, due to Sz\"{o}ll\H{o}si \cite{sz2}, is a first one in this direction:

\index{singular values}
\index{minors}

\begin{theorem}
Assuming that a square matrix
$$U=\begin{pmatrix}A&B\\C&D\end{pmatrix}$$
is unitary, with $A\in M_r(\mathbb C)$, $D\in M_d(\mathbb C)$, then:
\begin{enumerate}
\item The singular values of $A,D$ are identical, up to $|r-d|$ values of $1$. 

\item $\det A=\det U\cdot\overline{\det D}$, so in particular, $|\det A|=|\det D|$.
\end{enumerate}
\end{theorem}

\begin{proof}
Here is a simplified proof. From the unitarity of $U$ we get:
\begin{align*}
A^*A+C^*C&=I_r\\
CC^*+DD^*&=I_d\\
AC^*+BD^*&=0_{r\times d}
\end{align*}

(1) This follows from the first two equations, and from the well-known fact that the matrices $CC^*,C^*C$ have the same eigenvalues, up to $|r-d|$ values of $0$.

\medskip

(2) By using the above unitarity equations, we have:
$$\begin{pmatrix}A&0\\C&I\end{pmatrix}
=\begin{pmatrix}A&B\\C&D\end{pmatrix}
\begin{pmatrix}I&C^*\\0&D^*\end{pmatrix}$$

The result follows then by taking determinants.
\end{proof}

We want to find a formula for the polar decomposition of $D$. Let us introduce:

\begin{definition}
Associated to any $A\in M_r(\pm1)$ are the matrices
\begin{eqnarray*}
X_A&=&(\sqrt{N}I_r+\sqrt{A^tA})^{-1}Pol(A)^t\\
Y_A&=&(\sqrt{N}I_r+\sqrt{AA^t})^{-1}
\end{eqnarray*}
depending on a parameter $N$.
\end{definition}

Observe that, in terms of the polar decomposition $A=VP$, we have:
\begin{eqnarray*}
X_A&=&(\sqrt{N}+P)^{-1}V^t\\
Y_A&=&V(\sqrt{N}+P)^{-1}V^t
\end{eqnarray*}

The idea now will be that, under the assumptions of Theorem 4.21, the polar parts of the matrices $A,D$ appearing there should be related by a simple formula, with the passage $Pol(A)\to Pol(D)$ involving the above matrices $X_A,Y_A$. 

\bigskip

In what follows we will focus on the case where $U\in U_N$ is replaced by $U=\sqrt{N}H$ with $H\in M_N(\pm1)$ Hadamard. In the non-singular case, following \cite{bs2}, we have:

\begin{proposition}
Assuming that a square matrix 
$$H=\begin{pmatrix}A&B\\C & D\end{pmatrix}\in M_N(\pm1)$$
is Hadamard, with $A\in M_r(\pm1)$ invertible, $D\in M_d(\pm1)$, and $||A||< \sqrt N$, the polar decomposition $D=UT$ is given by the formulae
$$U=\frac{1}{\sqrt{N}}(D-E)\quad,\quad 
T=\sqrt{N}I_d-S$$ 
where $E=CX_AB$ and $S=B^tY_AB$.
\end{proposition}

\begin{proof}
Since $H$ is Hadamard, we can use the formulae coming from:
$$\begin{pmatrix}A&B\\C&D\end{pmatrix}
\begin{pmatrix}A^t&C^t\\B^t&D^t\end{pmatrix}
=\begin{pmatrix}A^t&C^t\\B^t&D^t\end{pmatrix}
\begin{pmatrix}A&B\\C&D\end{pmatrix}
=\begin{pmatrix}N&0\\0&N\end{pmatrix}$$

We start from the singular value decomposition of $A$: 
$$A=Vdiag(s_i)X^t$$

Here $V,X \in O_r$ and $s_i\in(0,||A||]$. From $AA^t+BB^t=NI_r$ we get:
$$BB^t=Vdiag(N-s_i^2)V^t$$

Thus, the singular value decomposition of $B$ is as follows, with $Y\in O_d$:
$$B=V\begin{pmatrix}diag(\sqrt{N-s_i^2})&0_{r\times(d-r)}\end{pmatrix}Y^t$$

Similarly, from $A^tA+C^tC=I_r$ we deduce the singular value decomposition for $C$, the result being that there exists an orthogonal matrix $\widetilde{Z}\in O_d$ such that: 
$$C=-\widetilde Z\begin{pmatrix}diag(\sqrt{N-s_i^2})\\0_{(d-r)\times r}\end{pmatrix}X^t$$

From $B^tB+D^tD=NI_d$ we obtain: 
$$D^tD = Y (diag(s_i^2)\oplus NI_{(d-r)})Y^t$$

Thus the polar decomposition of $D$ reads:
$$D=UY(diag(s_i)\oplus\sqrt NI_{(d-r)})Y^t$$

Let $Z=UY$. By using the orthogonality relation $CA^t+DB^t=0_{d \times r}$, we obtain:
$$\widetilde Z \begin{pmatrix}diag(s_i\sqrt{N-s_i^2})\\0_{(d-r)\times r}\end{pmatrix}
=Z\begin{pmatrix}diag(s_i\sqrt{N-s_i^2})\\0_{(d-r)\times r}\end{pmatrix}$$

From the assumptions of our theorem, we have the following inequality:
$$s_i\sqrt{N-s_i^2}>0$$

Thus  $Z^t\widetilde Z = I_r \oplus Q$, for some orthogonal matrix $Q \in O_d$. Plugging $\widetilde Z = Z(I_r \oplus Q)$ in the singular value decomposition formula for $C$, we obtain:
\begin{eqnarray*}
C
&=&-Z(I_r \oplus Q)\begin{pmatrix}diag(\sqrt{N-s_i^2})\\0_{(d-r)\times r}\end{pmatrix}X^t\\
&=&-Z\begin{pmatrix}diag(\sqrt{N-s_i^2})\\0_{(d-r)\times r}\end{pmatrix}X^t
\end{eqnarray*}

To summarize, we have found $V,X \in O_r$ and $Y,Z \in O_d$ such that:
\begin{align*}
A &=Vdiag(s_i)X^t\\
B &= V\begin{pmatrix}diag(\sqrt{N-s_i^2})&0_{r\times(d-r)}
\end{pmatrix}Y^t\\
C &=-Z\begin{pmatrix}diag(\sqrt{N-s_i^2})\\0_{(d-r)\times r}\end{pmatrix}X^t\\
D &= Z (diag(s_i)\oplus \sqrt N I_{(d-r)}) Y^t
\end{align*}

Now with $U,T,E,S$ defined as in the statement, we obtain:
\begin{eqnarray*}
U&=&ZY^t\\
E&=&Z(diag(\sqrt{N}-s_i)\oplus0_{d-r})Y^t\\
\sqrt{A^tA}&=&Xdiag(s_i)X^t\\
(\sqrt{N}I_r+\sqrt{A^tA})^{-1}&=&Xdiag(1/(\sqrt{N}+s_i))X^t\\
X_A&=&Xdiag(1/(\sqrt{N}+s_i))V^t\\
CX_AB&=&Z(diag(\sqrt{N}-s_i)\oplus0_{d-r})Y^t
\end{eqnarray*}

Thus we have $E=CX_AB$, as claimed. Also, we have:
\begin{eqnarray*}
T&=&Y(diag(s_i)\oplus\sqrt{N}I_{d-r})Y^t\\
S&=&Y(diag(\sqrt{N}-s_i)\oplus0_{d-r})Y^t\\
\sqrt{AA^t}&=&Vdiag(s_i)V^t\\
Y_A&=&Vdiag(1/(\sqrt{N}+s_i))V^t\\
B^tY_AB&=&Y(diag(\sqrt{N}-s_i)\oplus0_{d-r})Y^t
\end{eqnarray*}

Thus, we have as well $S=B^tY_AB$, as claimed, and we are done.
\end{proof}

Observe that, in the above statement, in the case where the size of the upper left block satisfies $r<\sqrt N$, the condition $||A||< \sqrt N$ is automatically satisfied. Our claim now is that all this is related to the notion of almost Hadamard matrix, from chapter 3. To be more precise, still following \cite{bs2}, let us introduce the following notion:

\index{AHP}
\index{AHM sign pattern}

\begin{definition}
A sign matrix $S\in M_N(\pm1)$ is called an almost Hadamard sign pattern (AHP) if it appears as
$$S_{ij}=sgn(H_{ij})$$
for a certain almost Hadamard matrix $H\in M_N(\mathbb R)$.
\end{definition}

Observe that, due to the theory in chapter 3, if a sign matrix $S$ is an AHP, then there exists a unique almost Hadamard matrix $H$ such that $S_{ij}=sgn(H_{ij})$, namely:
$$H =\sqrt{N}Pol(S)$$

Getting back to Proposition 4.23, let us try to find out when $D$ is AHP. For this purpose, we must estimate the quantity $||E||_\infty=\max_{ij}|E_{ij}|$, and we have here:

\begin{proposition}
Assuming that a matrix
$$H=\begin{pmatrix}A&B\\C & D\end{pmatrix}\in M_N(\pm1)$$
is an Hadamard matrix, with $A\in M_r(\pm1)$, $D\in M_d(\pm1)$ and $r\leq d$, then 
$$Pol(D)=\frac{1}{\sqrt{N}}(D-E)$$
with $E$ satisfying:
\begin{enumerate}
\item $||E||_\infty\leq\frac{r\sqrt{r}}{\sqrt{r}+\sqrt{N}}$ when $A$ is Hadamard.

\item $||E||_\infty\leq\frac{r^2c\sqrt{N}}{N-r^2}$ if $r^2<N$, with $c=||Pol(A)-\frac{A}{\sqrt{N}}||_\infty$.

\item $||E||_\infty\leq\frac{r^2(1+\sqrt{N})}{N-r^2}$ if $r^2<N$.
\end{enumerate}
\end{proposition}

\begin{proof}
We use the basic fact that for two rectangular matrices which are multipliable, $X\in M_{p\times r}(\mathbb C)$ and $Y\in M_{r\times q}(\mathbb C)$, we have the following estimate:
$$||XY||_\infty\leq r||X||_\infty||Y||_\infty$$

Thus, according to Proposition 4.23, we have:
\begin{eqnarray*}
||E||_\infty
&=&||CX_AB||_\infty\\
&\leq&r^2||C||_\infty||X_A||_\infty||B||_\infty\\
&=&r^2||X_A||_\infty
\end{eqnarray*}

(1) If $A$ is Hadamard, $AA^t =rI_r$, $Pol(A)=A/\sqrt{r}$ and thus: 
\begin{eqnarray*}
X_A
&=&(\sqrt{N}I_r+\sqrt{r}I_r)^{-1}\frac{A^t}{\sqrt{r}}\\
&=&\frac{A^t}{r+\sqrt{rN}}
\end{eqnarray*}

We therefore obtain from this:
$$||X_A||_\infty=\frac{1}{r+\sqrt{rN}}$$

But this gives the result.

\medskip

(2) According to the definition of $X_A$, we have:
\begin{eqnarray*}
X_A
&=&(\sqrt{N}I_r+\sqrt{A^tA})^{-1}Pol(A)^t\\
&=&(NI_r-A^tA)^{-1}(\sqrt{N}I_r-\sqrt{A^tA})Pol(A)^t\\
&=&(NI_r-A^tA)^{-1}(\sqrt{N}Pol(A)-A)^t
\end{eqnarray*}

We therefore obtain the following estimate:
\begin{eqnarray*}
||X_A||_\infty
&\leq&r||(NI_r-A^tA)^{-1}||_\infty||\sqrt{N}Pol(A)-A||_\infty\\
&=&\frac{rc}{\sqrt{N}}\Big|\Big|\left(I_r-\frac{A^tA}{N}\right)^{-1}\Big|\Big|_\infty
\end{eqnarray*}

Now by using $||A^tA||_\infty\leq r$, we obtain:
\begin{eqnarray*}
\Big|\Big|\left(I_r-\frac{A^tA}{N}\right)^{-1}\Big|\Big|_\infty
&\leq&\sum_{k=0}^\infty\frac{||(A^tA)^k||_\infty}{N^k}\\
&\leq&\sum_{k=0}^\infty\frac{r^{2k-1}}{N^k}\\
&=&\frac{1}{r}\cdot\frac{1}{1-r^2/N}\\
&=&\frac{N}{rN-r^3}
\end{eqnarray*}

Thus we have the following estimate:
$$||X_A||_\infty
\leq\frac{rc}{\sqrt{N}}\cdot\frac{N}{rN-r^3}
=\frac{c\sqrt{N}}{N-r^2}$$

But this gives the result.

\medskip

(3) This follows from (2), because:
$$c
\leq||Pol(A)||_\infty+||A/\sqrt{N}||_\infty
\leq1+\frac{1}{\sqrt{N}}$$

The proof is now complete.
\end{proof}

Following \cite{bs2}, we can now state and prove a main result, as follows:

\index{submatrix}
\index{minors}
\index{AHP}

\begin{theorem}
Assume that a matrix
$$H=\begin{pmatrix}A&B\\C & D\end{pmatrix}$$
is Hadamard, with $A\in M_r(\pm1),H\in M_N(\pm1)$.
\begin{enumerate}
\item If $A$ is Hadamard, and $N>r(r-1)^2$, then $D$ is AHP. 

\item If $N>\frac{r^2}{4}(x+\sqrt{x^2+4})^2$, where $x=r||Pol(A)-\frac{A}{\sqrt{N}}||_\infty$, then $D$ is AHP.

\item If $N>\frac{r^2}{4}(r+\sqrt{r^2+8})^2$, then $D$ is AHP.
\end{enumerate}
\end{theorem}

\begin{proof}
This follows from the various estimates that we have, as follows:

\medskip

(1) This follows from Proposition 4.25 (1), because:
\begin{eqnarray*}
\frac{r\sqrt{r}}{\sqrt{r}+\sqrt{N}}<1
&\iff&r<1+\sqrt{N/r}\\
&\iff&r(r-1)^2<N
\end{eqnarray*}

(2) This follows from Proposition 4.25 (2), because:
\begin{eqnarray*}
\frac{r^2c\sqrt{N}}{N-r^2}<1
&\iff&N-r^2c\sqrt{N}>r^2\\
&\iff&(2\sqrt{N}-r^2c)^2>r^4c^2+4r^2
\end{eqnarray*}

Indeed, this is equivalent to:
\begin{eqnarray*}
2\sqrt{N}
&>&r^2c+r\sqrt{r^2c^2+4}\\
&=&r(x+\sqrt{x^2+4})
\end{eqnarray*}

Here the value of $x$ is as follows:
$$x
=rc
=r\left|\left|Pol(A)-\frac{A}{\sqrt{N}}\right|\right|_\infty$$

(3) This follows from Proposition 4.25 (3), because:
\begin{eqnarray*}
\frac{r^2(1+\sqrt N)}{N-r^2}<1
&\iff&N-r^2\sqrt{N}>2r^2\\
&\iff&(2\sqrt{N}-r^2)^2>r^4+8r^2
\end{eqnarray*}

Indeed, this is equivalent to:
$$2\sqrt{N}>r^2+r\sqrt{r^2+8}$$

But this gives the result.
\end{proof}

As a technical comment, for $A\in M_r(\pm1)$ Hadamard, Proposition 4.25 (2) gives:
$$||E||_\infty
\leq\frac{r^2\sqrt{N}}{N-r^2}\left(\frac{1}{\sqrt{r}}-\frac{1}{\sqrt{N}}\right)
=\frac{r\sqrt{r}N-r^2}{N-r^2}$$

Thus $||E||_\infty<1$ for $N>r^3$, which is slightly weaker than Theorem 4.26 (1). 

\bigskip

In view of the results above, it is convenient to make the following convention:

\begin{definition}
We denote by $\{x\}_{m \times n} \in M_{m \times n}(\mathbb R)$ the all-$x$ matrix, and by
$$\begin{Bmatrix}x_{11}&\ldots&x_{1l}\\ \ldots&\ldots&\ldots\\ x_{k1}&\ldots&x_{kl}\end{Bmatrix}_{(m_1,\ldots,m_k) \times (n_1, \ldots, n_l)}$$
the matrix having all-$x_{ij}$ rectangular blocks $X_{ij}=\{x_{ij}\}_{m_i \times n_j} \in M_{m_i \times n_j}(\mathbb R)$, of prescribed size. In the case of square diagonal blocks, we simply write $\{x\}_n= \{x\}_{n \times n}$ and: 
$$\begin{Bmatrix}x_{11}&\ldots&x_{1k}\\ \ldots&\ldots&\ldots\\ x_{kk}&\ldots&x_{kk}\end{Bmatrix}_{n_1, \ldots n_k} = \begin{Bmatrix}x_{11}&\ldots&x_{1k}\\ \ldots&\ldots&\ldots\\ x_{k1}&\ldots&x_{kk}\end{Bmatrix}_{(n_1,\ldots,n_k) \times (n_1, \ldots, n_k)}$$
\end{definition}

Modulo equivalence, the $\pm1$ matrices of size $r=1,2$ are as follows:
$$\begin{pmatrix}+\end{pmatrix}_{(1)}\quad,\quad
\begin{pmatrix}+&+\\+&-\end{pmatrix}_{(2)}\quad,\quad
\begin{pmatrix}+&+\\+&+\end{pmatrix}_{(2')}$$

In the cases $(1)$ and $(2)$ above, where the matrix $A$ is invertible, the spectral properties of their complementary matrices are as follows:

\begin{theorem}
For the $N\times N$ Hadamard matrices of type
$$\begin{pmatrix}+&+\\ +&D\end{pmatrix}_{(1)}\quad,\quad
\begin{pmatrix}
+&+&+&+\\
+&-&+&-\\
+&+&D_{00}&D_{01}\\
+&-&D_{10}&D_{11}
\end{pmatrix}_{(2)}$$
the polar decomposition $D=UT$ with 
$$U=\frac{1}{\sqrt{N}}(D-E)\quad,\quad T=\sqrt{N}I-S$$
is given by the following formulae:
$$E_{(1)}=\begin{Bmatrix}\frac{1}{1+\sqrt{N}}\end{Bmatrix}_{N-1}\quad,\quad
E_{(2)}=\frac{2}{2+\sqrt{2N}}\begin{Bmatrix}1&1\\1&-1\end{Bmatrix}_{N/2-1,N/2-1}$$
$$S_{(1)}=\begin{Bmatrix}\frac{1}{1+\sqrt{N}}\end{Bmatrix}_{N-1}\quad,\quad
S_{(2)}=\frac{2}{\sqrt{2}+\sqrt{N}}\begin{Bmatrix}1&0\\0&1\end{Bmatrix}_{N/2-1,N/2-1}$$
In particular, all the matrices $D$ above are AHP.
\end{theorem}

\begin{proof}
For $A\in M_r(\pm1)$ Hadamard, the quantities in Definition 4.22 are:
$$X_A=\frac{A^t}{r+\sqrt{rN}}$$
$$Y_A=\frac{I_r}{\sqrt{r}+\sqrt{N}}$$

These formulae follow indeed from the following equalities:
$$AA^t=A^tA=rI_r$$
$$Pol(A)=A/\sqrt{r}$$

(1) Using the notation introduced in Definition 4.27, we have here:
$$B_{(1)} = \{1\}_{1 \times N-1}$$
$$C_{(1)} = B_{(1)}^t$$

Since the matrix $A_{(1)}=[+]$ is Hadamard we have:
$$X_{A_{(1)}}=Y_{A_{(1)}}=\frac{1}{1+\sqrt{N}}$$

We therefore obtain that: 
\begin{eqnarray*}
E_{(1)}
&=&\frac{1}{1+\sqrt{N}}\{1\}_{N-1 \times 1}[1]\{1\}_{1 \times N-1}\\
&=&\frac{1}{1+\sqrt{N}}\{1\}_{N-1}
\end{eqnarray*}

Similarly, we obtain that:
\begin{eqnarray*}
S_{(1)}
&=&\frac{1}{1+\sqrt{N}}\{1\}_{N-1 \times 1}\{1\}_{1 \times N-1}\\
&=&\frac{1}{1+\sqrt{N}}\{1\}_{N-1}
\end{eqnarray*}

(2) Using the orthogonality of the first two rows of $H_{(2)}$, we find that the matrices $D_{00}$ and $D_{11}$ have size $N/2-1$. Since since the matrix $A_{(2)}=[^+_+{\ }^+_-]$ is Hadamard we have:
$$X_{A_{(2)}}=\frac{A}{2+\sqrt{2N}}$$
$$Y_{A_{(2)}}=\frac{I_2}{\sqrt{2}+\sqrt{N}}$$

But this gives the following formula:
\begin{eqnarray*}
&&E_{(2)}\\
&=&\frac{1}{2+\sqrt{2N}}\begin{Bmatrix}1&1\\1&-1\end{Bmatrix}_{(N/2-1,N/2-1) \times (1,1)}\begin{pmatrix}1&1\\1&-1\end{pmatrix}\begin{Bmatrix}1&1\\1&-1\end{Bmatrix}_{(1,1) \times (N/2-1,N/2-1)}\\
&=&\frac{2}{2+\sqrt{2N}}\begin{Bmatrix}1&1\\1&-1\end{Bmatrix}_{N/2-1,N/2-1}
\end{eqnarray*}

Similarly, we obtain the following formula:
\begin{eqnarray*}
&&S_{(2)}\\
&=&\frac{1}{\sqrt{2}+\sqrt{N}}\begin{Bmatrix}1&1\\1&-1\end{Bmatrix}_{(N/2-1,N/2-1) \times (1,1)}\begin{Bmatrix}1&1\\1&-1\end{Bmatrix}_{(1,1) \times (N/2-1,N/2-1)}\\
&=&\frac{2}{\sqrt{2}+\sqrt{N}}\begin{Bmatrix}1&0\\0&1\end{Bmatrix}_{N/2-1,N/2-1}
\end{eqnarray*}

Thus, we have obtained the formulae in the statement.
\end{proof}

We refer to \cite{bs2} for more on all the above.

\section*{4e. Exercises} 

Here is a first exercise, in connection with the PHM:

\begin{exercise}
Find the almost PHM in the cases $M=1,2$. 
\end{exercise}

To start with, there is some differential geometry to be done here, in analogy with the differential geometry computations done in chapter 3.

\begin{exercise}
Work out the asymptotic count for the $5\times N$ PHM.
\end{exercise}

To be more precise, the problem here is that of completing the $M=5$ work that we started above, and recovering from this the de Launey-Levin formula, at $M=5$.

\begin{exercise}
Write down the axioms and basic theory of the AHP.
\end{exercise}

To be more precise, we know from chapter 3 the axioms and basic theory of the AHM, and the problem is that of converting that material in AHP terms.

\part{Complex matrices}

\ \vskip50mm

\begin{center}
{\em Beulah Land, I'm longing for you

And some day on thee I'll stand

There my home shall be eternal

Beulah Land, sweet Beulah Land}
\end{center}

\chapter{Complex matrices}

\section*{5a. Basic theory}

We have seen that the Hadamard matrices $H\in M_N(\pm1)$ are very interesting objects. In what follows, we will be interested in their complex versions:

\index{Hadamard matrix}
\index{complex Hadamard matrix}

\begin{definition}
A complex Hadamard matrix is a square matrix whose entries belong to the unit circle in the complex plane, 
$$H\in M_N(\mathbb T)$$
and whose rows are pairwise orthogonal, with respect to the scalar product of $\mathbb C^N$.
\end{definition}

Here, and in what follows, the scalar product is the usual one on $\mathbb C^N$, taken to be linear in the first variable and antilinear in the second one:
$$<x,y>=\sum_ix_i\bar{y}_i$$

As basic examples of complex Hamadard matrices, we have the real Hadamard matrices, $H\in M_N(\pm1)$, which have sizes $N\in\{2\}\cup4\mathbb N$. Here is now a new, motivating example, with $w=e^{2\pi i/3}$, which appears at the forbidden size value $N=3$:
$$F_3=\begin{pmatrix}1&1&1\\ 1&w&w^2\\ 1&w^2&w\end{pmatrix}$$

And here is another example, which appears at $N=4$, and whose combinatorics is different from the one of the unique $4\times4$ real Hadamard matrix, $W_4\sim K_4$:
$$F_4=\begin{pmatrix}
1&1&1&1\\
1&i&-1&-i\\
1&-1&1&-1\\
1&-i&-1&i
\end{pmatrix}$$

We will see that there are many other examples, and in particular that there are such matrices at any $N\in\mathbb N$, which in addition can be chosen to be circulant. Thus, the HC and CHC problematics will dissapear in the general complex setting. And we will also see that many other questions about the real Hadamard matrices $H\in M_N(\pm1)$ become far more clear, and sometimes even solvable, when passing to the complex case.

\bigskip

Before anything, however, let us recommend some reading. Although the field of complex numbers $\mathbb C$ is something very familiar in mathematics, and there are plenty of good reasons for sometimes using it, instead of the field of real numbers $\mathbb R$, in what concerns the matrices, things are more tricky. Why, after all, looking at $M_N(\mathbb C)$?

\bigskip

The answer to this question comes from physics, and more specifically from quantum mechanics. Remember Newton, Leibnitz and others who started talking about functions, derivatives, integrals, and all sorts of other things, that we learn now in 1st year at the university, motivated by classical mechanics? Well, pretty much the same happened with Heisenberg, Schr\"odinger, Dirac and others, who all of the sudden started to talk about complex matrices, motivated by quantum mechanics. And with these complex matrices being now part of the mathematical landscape too, starting with the 3rd year or so.

\bigskip

So, quantum mechanics. This is, and we repeat, something that you need to know a bit, in order to love the complex matrices, and appreciate the remainder of this book. Standard places for learning it are the books of Feynman \cite{fey}, Griffiths \cite{gri}, Weinberg \cite{wei}. There are some delightful good old books as well, if you prefer, such as Dirac \cite{dir}, von Neumann \cite{von}, Weyl \cite{wey}. And for more fancy stuff, if you're really into action, teaching you how to win a war by totally paralyzing the enemy, with a powerful quantum computer, go with Bengtsson-\.Zyczkowski \cite{bzy}, Nielsen-Chuang \cite{nch}, Watrous \cite{wat}.

\bigskip

Getting back now to the complex Hadamard matrices, although these originate in a 1962 paper by Butson \cite{but}, motivated by pure mathematics, their study only really took off in the 90s, under the influence of people like Haagerup \cite{ha1}, Jones \cite{jo3}, Popa \cite{pop}, all mathematicians interested in quantum mechanics. Later on physicists joined too, of course. And so again, conclusion to this, to be kept in mind: quantum mechanics.

\bigskip

In what follows we will take Definition 5.1 as it is, as a nice and natural mathematical definition, which is fully motivated, mathematically speaking, by the few remarks made afterwards. Let us start our study of the complex Hadamard matrices by extending some basic results from the real case, from chapter 1. First, we have: 

\index{Hadamard matrix manifold}
\index{unitary group}

\begin{proposition}
The set formed by the $N\times N$ complex Hadamard matrices is the real algebraic manifold
$$X_N=M_N(\mathbb T)\cap\sqrt{N}U_N$$
where $U_N$ is the unitary group, the intersection being taken inside $M_N(\mathbb C)$.
\end{proposition}

\begin{proof}
Let $H\in M_N(\mathbb T)$. Then $H$ is Hadamard if and only if its rescaling $U=H/\sqrt{N}$ belongs to the unitary group $U_N$, and so when $H\in X_N$, as claimed.
\end{proof}

We should mention that the above manifold $X_N$, while appearing by definition as an intersection of smooth manifolds, is very far from being smooth. We will be back to this, later on. As a basic consequence now of the above result, we have:

\begin{proposition}
Let $H\in M_N(\mathbb C)$ be an Hadamard matrix.
\begin{enumerate}
\item The columns of $H$ must be pairwise orthogonal.

\item The matrices $H^t,\bar{H},H^*\in M_N(\mathbb C)$ are Hadamard as well. 
\end{enumerate}
\end{proposition}

\begin{proof}
We use the well-known fact that if a matrix is unitary, $U\in U_N$, then so is its complex conjugate $\bar{U}=(\bar{U}_{ij})$, the inversion formulae being as follows:
$$U^*=U^{-1}\quad,\quad
U^t=\bar{U}^{-1}$$

Thus the unitary group $U_N$ is stable under the following operations:
$$U\to U^t\quad,\quad
U\to\bar{U}\quad,\quad 
U\to U^*$$

It follows that the algebraic manifold $X_N$ constructed in Proposition 5.2 is stable as well under these operations. But this gives all the assertions.
\end{proof}

Let us introduce now the following equivalence notion for the complex Hadamard matrices, taking into account some basic operations which can be performed:

\index{equivalence}
\index{Hadamard equivalence}
\index{dephased matrix}

\begin{definition}
Two complex Hadamard matrices are called equivalent, and we write $H\sim K$, when it is possible to pass from $H$ to $K$ via the following operations:
\begin{enumerate}
\item Permuting the rows, or permuting the columns.

\item Multiplying the rows or columns by numbers in $\mathbb T$.
\end{enumerate}
Also, we say that $H$ is dephased when its first row and column consist of $1$ entries.
\end{definition}

The same remarks as in the real case apply. First of all, we have not taken into account the results in Proposition 5.3 when formulating the above definition, because the operations $H\to H^t,\bar{H},H^*$ are far more subtle than those in (1,2) above.

\bigskip

Regarding the equivalence, there is a certain group $G$ acting there, made of two copies of $S_N$, one for the rows and one for the columns, and of two copies of $\mathbb T^N$, once again one for the rows, and one for the columns. It is possible to be a bit more explicit here, with a formula for $G$ and so on, but we will not need this, in what follows next.

\bigskip

Observe that, up to the above equivalence relation, any complex Hadamard matrix $H\in M_N(\mathbb T)$ can be put in dephased form. Moreover, the dephasing operation is unique, if we allow only the operations (2) in Definition 5.4, namely row and column multiplications by numbers in $\mathbb T$. In what follows, ``dephasing the matrix'' will have precisely this meaning, namely dephasing by using the operations (2) in Definition 5.4.

\bigskip

Regarding analytic aspects, once again in analogy with the study from the real case, we can locate the complex Hadamard matrices inside $M_N(\mathbb T)$, as follows:

\index{determinant bound}

\begin{theorem}
Given a matrix $H\in M_N(\mathbb T)$, we have
$$|\det(H)|\leq N^{N/2}$$
with equality precisely when $H$ is Hadamard.
\end{theorem}

\begin{proof}
By using the basic properties of the determinant, as in the real case, we have indeed the following estimate, valid for any vectors $H_1,\ldots,H_N\in\mathbb T^N$:
\begin{eqnarray*}
|\det(H_1,\ldots,H_N)|
&\leq&||H_1||\times\ldots\times||H_N||\\
&=&(\sqrt{N})^N
\end{eqnarray*}

Moreover, again as in the real case, the equality situation appears precisely when our vectors $H_1,\ldots,H_N\in\mathbb T^N$ are pairwise orthogonal, and this gives the result.
\end{proof}

From a ``dual'' point of view, the question of locating $X_N$ inside $\sqrt{N}U_N$, once again via analytic methods, makes sense as well, and we have here the following result:

\index{norm maximizer}

\begin{theorem}
Given a matrix $U\in U_N$ we have
$$||U||_1\leq N\sqrt{N}$$
with equality precisely when $H=\sqrt{N}U$ is Hadamard.
\end{theorem}

\begin{proof}
We have indeed the following estimate, valid for any $U\in U_N$:
\begin{eqnarray*}
||U||_1
&=&\sum_{ij}|U_{ij}|\\
&\leq&N\left(\sum_{ij}|U_{ij}|^2\right)^{1/2}\\
&=&N\sqrt{N}
\end{eqnarray*}

The equality case holds when $|U_{ij}|=\sqrt{N}$, for any $i,j$. But this amounts in saying that the rescaled matrix $H=\sqrt{N}U$ must satisfy $H\in M_N(\mathbb T)$, as desired.
\end{proof} 

The above Cauchy-Schwarz estimate can be improved with a H\"older estimate, the conclusion being that the rescaled Hadamard matrices maximize the $p$-norm on $U_N$ at any $p\in[1,2)$, and minimize it at any $p\in(2,\infty]$. We will be back to this.

\section*{5b. Fourier matrices}

At the level of the examples now, we have the following basic construction:

\index{Fourier matrix}

\begin{theorem}
The Fourier matrix, $F_N=(w^{ij})$ with $w=e^{2\pi i/N}$, which in standard matrix form, with indices $i,j=0,1,\ldots,N-1$, is as follows,
$$F_N=\begin{pmatrix}
1&1&1&\ldots&1\\
1&w&w^2&\ldots&w^{N-1}\\
1&w^2&w^4&\ldots&w^{2(N-1)}\\
\vdots&\vdots&\vdots&&\vdots\\
1&w^{N-1}&w^{2(N-1)}&\ldots&w^{(N-1)^2}
\end{pmatrix}$$
is a complex Hadamard matrix, in dephased form.
\end{theorem}

\begin{proof}
By using the standard fact that the averages of complex numbers correspond to barycenters, we conclude that the scalar products between the rows of $F_N$ are:
\begin{eqnarray*}
<R_a,R_b>
&=&\sum_jw^{aj}w^{-bj}\\
&=&\sum_jw^{(a-b)j}\\
&=&N\delta_{ab}
\end{eqnarray*}

Thus $F_N$ is indeed a complex Hadamard matrix. As for the fact that $F_N$ is dephased, this follows from our convention $i,j=0,1,\ldots,N-1$, which is there for this.
\end{proof}

\index{HC}
\index{CHC}

As an obvious consequence of the above result, there is no analogue of the HC in the complex case. We will see later on, in chapter 9 below, that the Fourier matrix $F_N$ can be put in circulant form, so there is no analogue of the CHC either, in this setting. As a first classification result now, in the complex case, we have:

\begin{proposition}
The Fourier matrices $F_2,F_3$, which are given by
$$F_2=\begin{pmatrix}1&1\\ 1&-1\end{pmatrix}\quad,\quad
F_3=\begin{pmatrix}1&1&1\\ 1&w&w^2\\ 1&w^2&w\end{pmatrix}$$
with $w=e^{2\pi i/3}$ are the only Hadamard matrices at $N=2,3$, up to equivalence.
\end{proposition}

\begin{proof}
The proof at $N=2$ is similar to the proof from the real case, from chapter 1. Indeed, given $H\in M_N(\mathbb T)$ Hadamard, we can dephase it, as follows:
$$\begin{pmatrix}a&b\\c&d\end{pmatrix}
\to\begin{pmatrix}1&1\\\bar{a}c&\bar{b}d\end{pmatrix}
\to\begin{pmatrix}1&1\\1&a\bar{b}\bar{c}d\end{pmatrix}$$

Thus, we obtain by dephasing the matrix $F_2$. Regarding now the case $N=3$, consider an Hadamard matrix $H\in M_3(\mathbb T)$, assumed to be in dephased form:
$$H=\begin{pmatrix}1&1&1\\ 1&x&y\\ 1&z&t\end{pmatrix}$$

The orthogonality conditions between the rows of this matrix read:
$$(1\perp2)\quad:\quad x+y=-1$$
$$(1\perp3)\quad:\quad z+t=-1$$
$$\ \ \ \,(2\perp3)\quad:\quad x\bar{z}+y\bar{t}=-1$$

In order to process this, consider an arbitrary equation of the following type:
$$p+q=-1\quad,\quad p,q\in\mathbb T$$

This equation tells us that the triangle having vertices at $1,p,q$ must be equilateral, and so that we must have $\{p,q\}=\{w,w^2\}$, with $w=e^{2\pi i/3}$. By using this fact, for the first two equations, we conclude that we must have:
$$\{x,y\}=\{w,w^2\}\quad,\quad 
\{z,t\}=\{w,w^2\}$$

As for the third equation, this gives $x\neq z$. Thus, $H$ is either the Fourier matrix $F_3$, or the matrix obtained from $F_3$ by permuting the last two columns, and we are done.
\end{proof}

In order to deal now with the case $N=4$, we already know, from our study in the real case, that we will need tensor products. So, let us formulate:

\index{tensor product}
\index{double indices}
\index{lexicographic order}

\begin{definition}
The tensor product of complex Hadamard matrices is given, in double indices, by $(H\otimes K)_{ia,jb}=H_{ij}K_{ab}$. In other words, we have the formula
$$H\otimes K=
\begin{pmatrix}
H_{11}K&\ldots&H_{1M}K\\ 
\vdots&&\vdots\\ 
H_{M1}K&\ldots&H_{MM}K
\end{pmatrix}$$
by using the lexicographic order on the double indices.
\end{definition}

Here the fact that $H\otimes K$ is indeed Hadamard comes from the fact that its rows $R_{ia}$ are pairwise orthogonal, as shown by the following computation:
\begin{eqnarray*}
<R_{ia},R_{kc}>
&=&\sum_{jb}H_{ij}K_{ab}\cdot \bar{H}_{kj}\bar{K}_{cb}\\
&=&\sum_jH_{ij}\bar{H}_{kj}\sum_bK_{ab}\bar{K}_{cb}\\
&=&M\delta_{ik}\cdot N\delta_{ac}\\
&=&MN\delta_{ia,kc}
\end{eqnarray*}

In order to advance now, our first task will be that of tensoring the Fourier matrices. We have here the following statement, refining and generalizing Theorem 5.7:

\index{abelian group}
\index{dual group}
\index{Fourier coupling}
\index{Fourier matrix}
\index{generalized Fourier matrix}

\begin{theorem}
Given a finite abelian group $G$, with dual group $\widehat{G}=\{\chi:G\to\mathbb T\}$, consider the Fourier coupling $\mathcal F_G:G\times\widehat{G}\to\mathbb T$, given by $(i,\chi)\to\chi(i)$.
\begin{enumerate}
\item Via the standard isomorphism $G\simeq\widehat{G}$, this Fourier coupling can be regarded as a square matrix, $F_G\in M_G(\mathbb T)$, which is a complex Hadamard matrix.

\item In the case of the cyclic group $G=\mathbb Z_N$ we obtain in this way, via the standard identification $\mathbb Z_N=\{1,\ldots,N\}$, the Fourier matrix $F_N$.

\item In general, when using a decomposition $G=\mathbb Z_{N_1}\times\ldots\times\mathbb Z_{N_k}$, the corresponding Fourier matrix is given by $F_G=F_{N_1}\otimes\ldots\otimes F_{N_k}$.
\end{enumerate}
\end{theorem}

\begin{proof}
This follows indeed from some basic facts from group theory:

\medskip

(1) With the identification $G\simeq\widehat{G}$ made our matrix is given by $(F_G)_{i\chi}=\chi(i)$, and the scalar products between the rows are then, as desired:
\begin{eqnarray*}
<R_i,R_j>
&=&\sum_\chi\chi(i)\overline{\chi(j)}\\
&=&\sum_\chi\chi(i-j)\\
&=&|G|\cdot\delta_{ij}
\end{eqnarray*}

(2) This follows from the well-known and elementary fact that, via the identifications $\mathbb Z_N=\widehat{\mathbb Z_N}=\{1,\ldots,N\}$, the Fourier coupling here is as follows, with $w=e^{2\pi i/N}$:
$$(i,j)\to w^{ij}$$

(3) We use here the following well-known formula, for the duals of products: 
$$\widehat{H\times K}=\widehat{H}\times\widehat{K}$$

At the level of the corresponding Fourier couplings, we obtain from this:
$$F_{H\times K}=F_H\otimes F_K$$

Now by decomposing $G$ into cyclic groups, as in the statement, and by using (2) for the cyclic components, we obtain the formula in the statement.
\end{proof}

As a first application of the above result, we have:

\index{Walsh matrix}
\index{Klein group}

\begin{proposition}
The Walsh matrix, $W_N$ with $N=2^n$, which is given by
$$W_N=\begin{pmatrix}1&1\\1&-1\end{pmatrix}^{\otimes n}$$
is the Fourier matrix of the finite abelian group $K_N=\mathbb Z_2^n$.
\end{proposition}

\begin{proof}
We know that the first Walsh matrix is a Fourier matrix:
$$W_2=F_2=F_{K_2}$$

Now by taking tensor powers we obtain from this that we have, for any $N=2^n$:
$$W_N
=W_2^{\otimes n}
=F_{K_2}^{\otimes n}
=F_{K_2^n}
=F_{K_N}$$

Thus, we are led to the conclusion in the statement.
\end{proof}

By getting back to classification, we will need the following result of Di\c t\u a \cite{dit}:

\index{Di\c t\u a deformation}
\index{deformed tensor product}
\index{affine deformation}

\begin{theorem}
If $H\in M_M(\mathbb T)$ and $K\in M_N(\mathbb T)$ are Hadamard, then so are the following two matrices, for any choice of a parameter matrix $Q\in M_{M\times N}(\mathbb T)$:
\begin{enumerate}
\item $H\otimes_QK\in M_{MN}(\mathbb T)$, given by $(H\otimes_QK)_{ia,jb}=Q_{ib}H_{ij}K_{ab}$.

\item $H\!\!{\ }_Q\!\otimes K\in M_{MN}(\mathbb T)$, given by $(H\!\!{\ }_Q\!\otimes K)_{ia,jb}=Q_{ja}H_{ij}K_{ab}$.
\end{enumerate}
These are called right and left Di\c t\u a deformations of $H\otimes K$, with parameter $Q$.
\end{theorem}

\begin{proof}
These results follow from the same computations as in the usual tensor product case, the idea being that the $Q$ parameters will cancel:

\medskip

(1) The rows $R_{ia}$ of the matrix $H\otimes_QK$ are indeed pairwise orthogonal, because:
\begin{eqnarray*}
<R_{ia},R_{kc}>
&=&\sum_{jb}Q_{ib}H_{ij}K_{ab}\cdot\bar{Q}_{kb}\bar{H}_{kj}\bar{K}_{cb}\\
&=&M\delta_{ik}\sum_bK_{ab}\bar{K}_{cb}\\
&=&M\delta_{ik}\cdot N\delta_{ac}\\
&=&MN\delta_{ik,ac}
\end{eqnarray*}

(2) The rows $L_{ia}$ of the matrix $H\!\!{\ }_Q\!\otimes K$ are orthogonal as well, because:
\begin{eqnarray*}
<L_{ia},L_{kc}>
&=&\sum_{jb}Q_{ja}H_{ij}K_{ab}\cdot\bar{Q}_{jc}\bar{H}_{kj}\bar{K}_{cb}\\
&=&N\delta_{ac}\sum_jH_{ij}\bar{H}_{kj}\\
&=&N\delta_{ac}\cdot M\delta_{ik}\\
&=&MN\delta_{ik,ac}
\end{eqnarray*}

Thus, both the matrices in the statement are Hadamard, as claimed.
\end{proof}

As a first observation, when the parameter matrix is the all-one matrix $\mathbb I\in M_{M\times N}(\mathbb T)$, we obtain in this way the usual tensor product of our matrices:
$$H\otimes_{\mathbb I}K
=H\!\!{\ }_{\mathbb I}\!\otimes K
=H\otimes K$$

As a non-trivial example now, let us compute the right deformations of the Walsh matrix $W_4=F_2\otimes F_2$, with arbitrary parameter matrix $Q=(^p_r{\ }^q_s)$. We have:
\begin{eqnarray*}
F_2\otimes_QF_2
&=&\begin{pmatrix}
1&1\\
1&-1
\end{pmatrix}
\otimes_{\begin{pmatrix}
p&q\\
r&s
\end{pmatrix}}
\begin{pmatrix}
1&1\\
1&-1
\end{pmatrix}\\
&=&\begin{pmatrix}
p&q&p&q\\
p&-q&p&-q\\
r&s&-r&-s\\ 
r&-s&-r&s
\end{pmatrix}
\end{eqnarray*}

This follows indeed by carefully working out what happens, by using the lexicographic order on the double indices, as explained in chapter 1. To be more precise, the usual tensor product $W_4=F_2\otimes F_2$ appears as follows:
$$W_4=
\begin{pmatrix}
ia\backslash jb&&00&01&10&11\\
\\
00&&1&1&1&1\\
01&&1&-1&1&-1\\
10&&1&1&-1&-1\\
11&&1&-1&-1&1
\end{pmatrix}$$

The corresponding values of the parameters $Q_{ib}$ to be inserted are as follows:
$$(Q_{ib})=\begin{pmatrix}
ia\backslash jb&&00&01&10&11\\
\\
00&&Q_{00}&Q_{01}&Q_{00}&Q_{01}\\
01&&Q_{00}&Q_{01}&Q_{00}&Q_{01}\\
10&&Q_{10}&Q_{11}&Q_{10}&Q_{11}\\
11&&Q_{10}&Q_{11}&Q_{10}&Q_{11}
\end{pmatrix}$$

With the notation $Q=(^p_r{\ }^q_s)$, this latter matrix becomes:
$$(Q_{ib})=\begin{pmatrix}
ia\backslash jb&&00&01&10&11\\
\\
00&&p&q&p&q\\
01&&p&q&p&q\\
10&&r&s&r&s\\
11&&r&s&r&s
\end{pmatrix}$$

Now by pointwise multiplying this latter matrix with the matrix $W_4$ given above, we obtain the announced formula for the deformed tensor product $F_2\otimes_QF_2$.

\bigskip

As for the left deformations of $W_4=F_2\otimes F_2$, once again with arbitrary parameter matrix $Q=(^p_r{\ }^q_s)$, these are given by a similar formula, as follows:
\begin{eqnarray*}
F_2\!\!{\ }_Q\!\otimes F_2
&=&\begin{pmatrix}
1&1\\
1&-1
\end{pmatrix}
\!{\ }_{\begin{pmatrix}
p&q\\
r&s
\end{pmatrix}}\!\otimes
\begin{pmatrix}
1&1\\
1&-1
\end{pmatrix}\\
&=&\begin{pmatrix}
p&p&r&r\\
q&-q&s&-s\\
p&p&-r&-r\\ 
q&-q&-s&s
\end{pmatrix}
\end{eqnarray*}

Observe that this latter matrix is transpose to $F_2\otimes_QF_2$. However, this is something accidental, coming from the fact that $F_2$, and so $W_4$ as well, are self-transpose.

\bigskip

With the above constructions in hand, we have the following result:

\begin{theorem}
The only complex Hadamard matrices at $N=4$ are, up to the standard equivalence relation, the matrices 
$$F_4^s=\begin{pmatrix}
1&1&1&1\\
1&-1&1&-1\\
1&s&-1&-s\\ 
1&-s&-1&s
\end{pmatrix}$$
with $s\in\mathbb T$, which appear as right Di\c t\u a deformations of $W_4=F_2\otimes F_2$. Moreover,
$$F_4^s\sim F_4^{-s}\sim F_4^{\bar{s}}\sim F_4^{-\bar{s}}$$
so we can assume, up to equivalence, that we have $s=e^{it}$ with $t\in[0,\pi/2]$.
\end{theorem}

\begin{proof}
There are several things to be done here, the idea being as follows:

\medskip

(1) First of all, the matrix $F_4^s$ is indeed Hadamard, appearing from the construction in Theorem 5.12, assuming that the parameter matrix $Q\in M_2(\mathbb T)$ is dephased:
$$Q=\begin{pmatrix}1&1\\1&s\end{pmatrix}$$

Observe also that, conversely, any right Di\c t\u a deformation of $W_4=F_2\otimes F_2$ is of this form. Indeed, if we consider such a deformation, with general parameter matrix $Q=(^p_r{\ }^q_s)$ as above, by dephasing we obtain an equivalence with $F_4^{s'}$, where $s'=ps/qr$:
\begin{eqnarray*}
\begin{pmatrix}
p&q&p&q\\
p&-q&p&-q\\
r&s&-r&-s\\ 
r&-s&-r&s
\end{pmatrix}
&\to&
\begin{pmatrix}
1&1&1&1\\
1&-1&1&-1\\
r/p&s/q&-r/p&-s/q\\ 
r/p&-s/q&-r/p&s/q
\end{pmatrix}\\
&\to&
\begin{pmatrix}
1&1&1&1\\
1&-1&1&-1\\
1&ps/qr&-1&-ps/qr\\ 
1&-ps/qr&-1&ps/qr
\end{pmatrix}
\end{eqnarray*}

(2) Summarizing, in what regards the first assertion, we must prove that any complex Hadamard matrix $H\in M_4(\mathbb T)$ is equivalent to one of the matrices $F_4^s$. But this follows by using the same arguments as in the proof from the real case, from chapter 1, at $N=4$, and from the proof of Proposition 5.8. Indeed, let us first dephase our matrix:
$$H=\begin{pmatrix}1&1&1&1\\ 1&a&b&c\\ 1&d&e&f\\ 1&g&h&i\end{pmatrix}$$

We use now the fact, coming from plane geometry, that the solutions $x,y,z,t\in\mathbb T$ of the equation $x+y+z+t=0$ are as follows, with $p,q\in\mathbb T$:
$$\{x,y,z,t\}=\{p,q,-p,-q\}$$ 

In our case, we have $1+a+d+g=0$, and so up to a permutation of the last 3 rows, our matrix must look at follows, for a certain $s\in\mathbb T$:
$$H=\begin{pmatrix}1&1&1&1\\ 1&-1&b&c\\ 1&s&e&f\\ 1&-s&h&i\end{pmatrix}$$

(3) In the case $s=\pm1$ we can permute the middle two columns, then repeat the same reasoning, and we end up with the matrix in the statement.

\medskip

(4) In the case $s\neq\pm1$ we have $1+s+e+f=0$, and so $-1\in\{e,f\}$. Up to a permutation of the last columns, we can assume $e=-1$, and our matrix becomes:
$$H=\begin{pmatrix}1&1&1&1\\ 1&-1&b&c\\ 1&s&-1&-s\\ 1&-s&h&i\end{pmatrix}$$

Similarly, from $1-s+h+i=0$ we deduce that $-1\in\{h,i\}$. In the case $h=-1$ our matrix must look as follows, and we are led to the matrix in the statement:
$$H=\begin{pmatrix}1&1&1&1\\ 1&-1&b&c\\ 1&s&-1&-s\\ 1&-s&-1&i\end{pmatrix}$$

As for the remaining case $i=-1$, here our matrix must look as follows:
$$H=\begin{pmatrix}1&1&1&1\\ 1&-1&b&c\\ 1&s&-1&-s\\ 1&-s&h&-1\end{pmatrix}$$

We obtain from the last column $c=s$, then from the second row $b=-s$, then from the third column $h=s$, and so our matrix must be as follows:
$$H=\begin{pmatrix}1&1&1&1\\ 1&-1&-s&s\\ 1&s&-1&-s\\ 1&-s&s&-1\end{pmatrix}$$

But, in order for the second and third row to be orthogonal, we must have $s\in\mathbb R$, and so $s=\pm1$, which contradicts our above assumption $s\neq\pm1$.

\medskip

(5) Thus, we are done with the proof of the main assertion. Regarding now the second assertion, observe first that by permuting the last two rows we have:
$$F_4^s=\begin{pmatrix}
1&1&1&1\\
1&-1&1&-1\\
1&s&-1&-s\\ 
1&-s&-1&s
\end{pmatrix}
\sim\begin{pmatrix}
1&1&1&1\\
1&-1&1&-1\\
1&-s&-1&s\\
1&s&-1&-s
\end{pmatrix}
=F_4^{-s}$$

Also, by starting with $F_4^{\bar{s}}$ and multiplying the last three rows by $-1,s,-s$, then intechanging the first two columns, and the last two columns, we have:
$$F_4^{\bar{s}}
=\begin{pmatrix}
1&1&1&1\\
1&-1&1&-1\\
1&\bar{s}&-1&-\bar{s}\\ 
1&-\bar{s}&-1&\bar{s}
\end{pmatrix}
\sim\begin{pmatrix}
1&1&1&1\\
-1&1&-1&1\\
s&1&-s&-1\\ 
-s&1&s&-1
\end{pmatrix}
\sim\begin{pmatrix}
1&1&1&1\\
1&-1&1&-1\\
1&s&-1&-s\\ 
1&-s&-1&s
\end{pmatrix}
=F_4^s$$

Thus, we are led to the final conclusion in the statement too.
\end{proof}

As a comment here, Theorem 5.13 does not close the discussion at $N=4$, because we would still like to prove that the matrices $F_4^s$ are non-equivalent, up to identifying $\{s,-s,\bar{s},-\bar{s}\}$. However, this is something undobale with bare hands, so we must trick. To be more precise, we would like to have an invariant which distinguishes the matrices $F_s$, and a natural candidate here is the ``complex glow'', which should be by definition the law over the equivalence class of the following quantity, called excess:
$$E(H)=\sum_{ij}H_{ij}$$

We will discuss this later, in chapters 10-11, but as an advertisement for the material there, let us mention that the quantities to look at are the moments $\int|E|^{2p}$, which are Laurent polynomials in $s\in\mathbb T$, and with $p=3$ doing the job. More on this later.

\section*{5c. Haagerup theorem}

At $N=5$ now, the situation is considerably more complicated, with $F_5$ being the only matrix. The key technical result here, due to Haagerup \cite{ha1}, is as follows:

\index{Haagerup lemma}

\begin{proposition}
Given an Hadamard matrix $H\in M_5(\mathbb T)$, chosen dephased,
$$H=\begin{pmatrix}
1&1&1&1&1\\
1&a&x&*&*\\
1&y&b&*&*\\
1&*&*&*&*\\
1&*&*&*&*
\end{pmatrix}$$
the numbers $a,b,x,y$ must satisfy the equation $(x-y)(x-ab)(y-ab)=0$.
\end{proposition}

\begin{proof}
This is something quite surprising, and tricky, the proof in \cite{ha1} being as follows. Let us look at the upper 3-row truncation of $H$, which is of the following form:
$$H'=\begin{pmatrix}
1&1&1&1&1\\
1&a&x&p&q\\
1&y&b&r&s
\end{pmatrix}$$

By using the orthogonality of the rows, we have:
$$(1+a+x)(1+\bar{b}+\bar{y})(1+\bar{a}y+b\bar{x})
=-(p+q)(r+s)(\bar{p}r+\bar{q}s)$$

On the other hand, by using $p,q,r,s\in\mathbb T$, we have:
\begin{eqnarray*}
(p+q)(r+s)(\bar{p}r+\bar{q}s)
&=&(r+p\bar{q}s+\bar{p}qr+s)(\bar{r}+\bar{s})\\
&=&1+p\bar{q}\bar{r}s+\bar{p}q+\bar{r}s+r\bar{s}+p\bar{q}+\bar{p}qr\bar{s}+1\\
&=&2Re(1+p\bar{q}+r\bar{s}+p\bar{q}r\bar{s})\\
&=&2Re[(1+p\bar{q})(1+r\bar{s})]
\end{eqnarray*}

We conclude that we have the following formula, involving $a,b,x,y$ only:
$$(1+a+x)(1+\bar{b}+\bar{y})(1+\bar{a}y+b\bar{x})\in\mathbb R$$

Now this is a product of type $(1+\alpha)(1+\beta)(1+\gamma)$, with the first summand being 1, and with the last summand, namely $\alpha\beta\gamma$, being real as well, as shown by the above general $p,q,r,s\in\mathbb T$ computation. Thus, when expanding, and we are left with:
\begin{eqnarray*}
&&(a+x)+(\bar{b}+\bar{y})+(\bar{a}y+b\bar{x})+(a+x)(\bar{b}+\bar{y})\\
&+&(a+x)(\bar{a}y+b\bar{x})+(\bar{b}+\bar{y})(\bar{a}y+b\bar{x})\in\mathbb R
\end{eqnarray*}

By expanding all the products, our formula looks as follows:
\begin{eqnarray*}
&&a+x+\bar{b}+\bar{y}+\bar{a}y+b\bar{x}+a\bar{b}+a\bar{y}+\bar{b}x+x\bar{y}\\
&+&1+ab\bar{x}+\bar{a}xy+b+\bar{a}\bar{b}y+\bar{x}+\bar{a}+b\bar{x}\bar{y}\in\mathbb R
\end{eqnarray*}

By removing from this all terms of type $z+\bar{z}$, we are left with:
$$a\bar{b}+x\bar{y}+ab\bar{x}+\bar{a}\bar{b}y+\bar{a}xy+b\bar{x}\bar{y}\in\mathbb R$$

Now by getting back to our Hadamard matrix, all this remains true when transposing it, which amounts in interchanging $x\leftrightarrow y$. Thus, we have as well:
$$a\bar{b}+\bar{x}y+ab\bar{y}+\bar{a}\bar{b}x+\bar{a}xy+b\bar{x}\bar{y}\in\mathbb R$$

By substracting now the two equations that we have, we obtain:
$$x\bar{y}-\bar{x}y+ab(\bar{x}-\bar{y})+\bar{a}\bar{b}(y-x)\in\mathbb R$$

Now observe that this number, say $Z$, is purely imaginary, because $\bar{Z}=-Z$. Thus our equation reads $Z=0$. On the other hand, we have the following formula:
\begin{eqnarray*}
abxyZ
&=&abx^2-aby^2+a^2b^2(y-x)+xy(y-x)\\
&=&(y-x)(a^2b^2+xy-ab(x+y))\\
&=&(y-x)(ab-x)(ab-y)
\end{eqnarray*}

Thus, our equation $Z=0$ corresponds to the formula in the statement.
\end{proof}

We are led in this way to the following theorem, also from Haagerup \cite{ha1}:

\index{Haagerup theorem}

\begin{theorem}
The only Hadamard matrix at $N=5$ is the Fourier matrix,
$$F_5=\begin{pmatrix}
1&1&1&1&1\\
1&w&w^2&w^3&w^4\\
1&w^2&w^4&w&w^3\\
1&w^3&w&w^4&w^2\\
1&w^4&w^3&w^2&w
\end{pmatrix}$$
with $w=e^{2\pi i/5}$, up to the standard equivalence relation for such matrices. 
\end{theorem}

\begin{proof}
Assume that have an Hadamard matrix $H\in M_5(\mathbb T)$, chosen dephased, and written as in Proposition 5.14, with emphasis on the upper left $2\times2$ subcorner: 
$$H=\begin{pmatrix}
1&1&1&1&1\\
1&a&x&*&*\\
1&y&b&*&*\\
1&*&*&*&*\\
1&*&*&*&*
\end{pmatrix}$$

(1) We know from Proposition 5.14, applied to $H$ itself, and to its transpose $H^t$ as well, that the entries $a,b,x,y$ must satisfy the following equations:
$$(a-b)(a-xy)(b-xy)=0$$
$$(x-y)(x-ab)(y-ab)=0$$

Our first claim is that, by doing some combinatorics, we can actually obtain from this $a=b$ and $x=y$, up to the equivalence relation for the Hadamard matrices: 
$$H\sim\begin{pmatrix}
1&1&1&1&1\\
1&a&x&*&*\\
1&x&a&*&*\\
1&*&*&*&*\\
1&*&*&*&*
\end{pmatrix}$$

Indeed, the above two equations lead to 9 possible cases, the first of which is, as desired, $a=b$ and $x=y$. As for the remaining 8 cases, here again things are determined by 2 parameters, and in practice, we can always permute the first 3 rows and 3 columns, and then dephase our matrix, as for our matrix to take the above special form.

\medskip

(2) With this result in hand, the combinatorics of the scalar products between the first 3 rows, and between the first 3 columns as well, becomes something which is quite simple to investigate. By doing a routine study here, and then completing it with a study of the lower right $2\times2$ corner as well, we are led to 2 possible cases, as follows:
$$H\sim\begin{pmatrix}
1&1&1&1&1\\
1&a&b&c&d\\
1&b&a&d&c\\
1&c&d&a&b\\
1&d&c&b&a
\end{pmatrix}\quad,\quad 
H\sim\begin{pmatrix}
1&1&1&1&1\\
1&a&b&c&d\\
1&b&a&d&c\\
1&c&d&b&a\\
1&d&c&a&b
\end{pmatrix}$$

(3) Our claim now is that the first case is in fact not possible. Indeed, we must have:
\begin{eqnarray*}
a+b+c+d&=&-1\\
2Re(a\bar{b})+2Re(c\bar{d})&=&-1\\
2Re(a\bar{c})+2Re(b\bar{d})&=&-1\\
2Re(a\bar{d})+2Re(b\bar{c})&=&-1
\end{eqnarray*}

Now since $|Re(x)|\leq1$ for any $x\in\mathbb T$, we deduce from the second equation that:
$$Re(a\bar{b})\leq 1/2$$

In other words, the arc length between $a,b$ satisfies:
$$\theta(a,b)\geq\pi/3$$

The same argument applies to $c,d$, and to the other pairs of numbers in the last 2 equations. Now since our equations are invariant under permutations of $a,b,c,d$, we can assume that $a,b,c,d$ are ordered in this way on the unit circle, and by the above, separated by $\geq\pi/3$ arc lengths. But this tells us that we have the following inequalities:
$$\theta(a,c)\geq 2\pi/3\quad,\quad 
\theta(b,d)\geq 2\pi/3$$

These two inequalities give the following estimates:
$$Re(a\bar{c})\leq-1/2\quad,\quad 
Re(b\bar{d})\leq-1/2$$

But these estimates contradict the third equation. Thus, our claim is proved.

\medskip

(4) Summarizing, we have proved so far that our matrix must be as follows:
$$H\sim\begin{pmatrix}
1&1&1&1&1\\
1&a&b&c&d\\
1&b&a&d&c\\
1&c&d&b&a\\
1&d&c&a&b
\end{pmatrix}$$

We are now in position of finishing. The orthogonality equations are as follows:
\begin{eqnarray*}
a+b+c+d&=&-1\\
2Re(a\bar{b})+2Re(c\bar{d})&=&-1\\
a\bar{c}+c\bar{b}+b\bar{d}+d\bar{a}&=&-1
\end{eqnarray*}

The third equation can be written in the following equivalent form:
\begin{eqnarray*}
Re[(a+b)(\bar{c}+\bar{d})]&=&-1\\
Im[(a-b)(\bar{c}-\bar{d})]&=&0
\end{eqnarray*}

By using now $a,b,c,d\in\mathbb T$, we obtain from this:
$$\frac{a+b}{a-b}\in i\mathbb R\quad,\quad 
\frac{c+d}{c-d}\in i\mathbb R$$

Thus we can find $s,t\in\mathbb R$ such that:
$$a+b=is(a-b)\quad,\quad 
c+d=it(c-d)$$

By plugging in these values, our system of equations simplifies, as follows:
\begin{eqnarray*}
(a+b)+(c+d)&=&-1\\
|a+b|^2+|c+d|^2&=&3\\
(a+b)(\bar{c}+\bar{d})&=&-1
\end{eqnarray*}

Now observe that the last equation implies in particular that we have:
$$|a+b|^2\cdot|c+d|^2=1$$

Thus $|a+b|^2,|c+d|^2$ must be roots of the following polynomial:
$$X^2-3X+1=0$$

But this gives the following equality of sets:
$$\Big\{|a+b|\,,\,|c+d|\Big\}=\left\{\frac{\sqrt{5}+1}{2}\,,\,\frac{\sqrt{5}-1}{2}\right\}$$

This is good news, because we are now into 5-th roots of unity. To be more precise, we have 2 cases to be considered, the first one being as follows, with $z\in\mathbb T$:
$$a+b=\frac{\sqrt{5}+1}{2}\,z\quad,\quad 
c+d=-\frac{\sqrt{5}-1}{2}\,z$$

From $a+b+c+d=-1$ we obtain $z=-1$, and by using this we obtain $b=\bar{a}$, $d=\bar{c}$. Thus we have the following formulae:
$$Re(a)=\cos(2\pi/5)\quad,\quad 
Re(c)=\cos(\pi/5)$$

We conclude that we have $H\sim F_5$, as claimed. As for the second case, with $a,b$ and $c,d$ interchanged, this leads to $H\sim F_5$ as well.
\end{proof}

\section*{5d. Further matrices}

At $N=6$ now, the situation becomes very complicated, with lots of ``exotic'' solutions, and with the structure of the Hadamard manifold $X_6$ being not understood yet, despite years of efforts. In fact, $X_6$ looks as complicated as the real algebraic manifolds can get. The simplest examples of Hadamard matrices at $N=6$ are as follows:

\begin{theorem}
We have the following basic Hadamard matrices, at $N=6$:
\begin{enumerate}
\item The Fourier matrix $F_6$.

\item The Di\c t\u a deformations of $F_2\otimes F_3$ and of $F_3\otimes F_2$.

\item The Haagerup matrix $H_6^q$.

\item The Tao matrix $T_6$.
\end{enumerate}
\end{theorem}

\begin{proof}
All this is elementary, the idea, and formulae of the matrices, being as follows:

\medskip

(1) This is something that we know well.

\medskip

(2) Consider indeed the dephased Di\c t\u a deformations of $F_2\otimes F_3$ and $F_3\otimes F_2$:
$$F_6^{(rs)}=F_2
\otimes_{\begin{pmatrix}
1&1&1\\
1&r&s
\end{pmatrix}}
F_3\qquad,\qquad 
F_6^{(^r_s)}=F_3
\otimes_{\begin{pmatrix}
1&1\\
1&r\\
1&s
\end{pmatrix}}F_2$$

Here $r,s$ are two parameters on the unit circle, $r,s\in\mathbb T$. In matrix form:
$$F_6^{(rs)}=\begin{pmatrix}
1&1&1&&1&1&1\\
1&w&w^2&&1&w&w^2\\
1&w^2&w&&1&w^2&w\\ 
\\
1&r&s&&-1&-r&-s\\
1&wr&w^2s&&-1&-wr&-w^2s\\
1&w^2r&ws&&-1&-w^2r&-ws
\end{pmatrix}$$

As for the other deformation, this is given by:
$$F_6^{(^r_s)}
=\begin{pmatrix}
1&1&&1&1&&1&1\\
1&-1&&1&-1&&1&-1\\
\\
1&r&&w&wr&&w^2&w^2r\\ 
1&-r&&w&-wr&&w^2&-w^2r\\
\\
1&s&&w^2&w^2s&&w&ws\\
1&-s&&w^2&-w^2s&&w&-ws
\end{pmatrix}$$

(3) The matrix here, from Haagerup's paper \cite{ha1}, is as follows, with $q\in\mathbb T$:
$$H_6^q=\begin{pmatrix}
1&1&1&1&1&1\\
1&-1&i&i&-i&-i\\ 
1&i&-1&-i&q&-q\\ 
1&i&-i&-1&-q&q\\
1&-i&\bar{q}&-\bar{q}&i&-1\\ 
1&-i&-\bar{q}&\bar{q}&-1&i
\end{pmatrix}$$

\index{Haagerup matrix}
\index{Tao matrix}

(4) The matrix here, from Tao's paper \cite{tao}, is as follows, with $w=e^{2\pi i/3}$:
$$T_6=\begin{pmatrix}
1&1&1&1&1&1\\ 
1&1&w&w&w^2&w^2\\ 
1&w&1&w^2&w^2&w\\
1&w&w^2&1&w&w^2\\ 
1&w^2&w^2&w&1&w\\ 
1&w^2&w&w^2&w&1
\end{pmatrix}$$

Observe that both $H_6^q$ and $T_6$ are indeed complex Hadamard matrices.
\end{proof}

\index{regular matrix}
\index{sums of roots}
\index{vanishing sum of roots}

The matrices in Theorem 5.16 are ``regular'', in the sense that the scalar products between rows appear in the simplest possible way, namely from vanishing sums of roots of unity, possibly rotated by a scalar. We will be back to this in chapter 6 below, with a result stating that these matrices are the only regular ones, at $N=6$.

\bigskip

In the non-regular case now, there are many known constructions at $N=6$. Here is one such construction, found by Bj\"orck and Fr\"oberg in \cite{bjf}:

\index{circulant matrix}
\index{Bj\"orck-Fr\"oberg matrix}

\begin{proposition}
The following is a complex Hadamard matrix,
$$BF_6=\begin{pmatrix}
1&ia&-a&-i&-\bar{a}&i\bar{a}\\
i\bar{a}&1&ia&-a&-i&-\bar{a}\\
-\bar{a}&i\bar{a}&1&ia&-a&-i\\
-i&-\bar{a}&i\bar{a}&1&ia&-a\\
-a&-i&-\bar{a}&i\bar{a}&1&ia\\
ia&-a&-i&-\bar{a}&i\bar{a}&1
\end{pmatrix}$$
where $a\in\mathbb T$ is one of the roots of $a^2+(\sqrt{3}-1)a+1=0$.
\end{proposition}

\begin{proof}
The matrix in the statement is circulant, in the sense that the rows appear by cyclically permuting the first row. Thus, we only have to check that the first row is orthogonal to the other 5 rows. But this follows from $a^2+(\sqrt{3}-1)a+1=0$.
\end{proof}

The obvious question here is how Bj\"orck and Fr\"oberg were able to construct the above matrix. This was done via some general theory for the circulant Hadamard matrices, and some computer simulations. We will discuss this in chapter 9 below.

\bigskip

Further study in the $N=6$ case leads to fairly complicated things, and we have here, as an illustrating example, the following result of Beauchamp-Nicoara \cite{ben}:

\index{self-adjoint matrix}
\index{Beauchamp-Nicoara matrix}

\begin{theorem}
The self-adjoint $6\times6$ Hadamard matrices are, up to equivalence
$$BN_6^q=
\begin{pmatrix}
1&1&1&1&1&1\\
1&-1&\bar{x}&-y&-\bar{x}&y\\
1&x&-1&t&-t&-x\\
1&-\bar{y}&\bar{t}&-1&\bar{y}&-\bar{t}\\
1&-x&-\bar{t}&y&1&\bar{z}\\
1&\bar{y}&-\bar{x}&-t&z&1
\end{pmatrix}$$
with $x,y,z,t\in\mathbb T$ depending on a parameter $q\in\mathbb T$, in a complicated way.
\end{theorem}

\begin{proof}
The study here can be done via a lot of work, the equations being:
\begin{eqnarray*}
x&=&\frac{1+2q+q^2-\sqrt{2}\sqrt{1+2q+2q^3+q^4}}{1+2q-q^2}\\
y&=&q\\
z&=&\frac{1+2q-q^2}{q(-1+2q+q^2)}\\
t&=&\frac{1+2q+q^2-\sqrt{2}\sqrt{1+2q+2q^3+q^4}}{-1+2q+q^2}
\end{eqnarray*}

All this is quite technical, and we refer here to \cite{ben}.
\end{proof}

There are many other examples at $N=6$, and no classification known. For a recent discussion on this subject, we refer to the survey paper of Tadej-\.Zyczkowski \cite{tz1}. 

\bigskip

Let us discuss now the case $N=7$. We will restrict the attention to case where the combinatorics comes from roots of unity. We use the following result of Sz\"oll\H{o}si \cite{sz2}:

\index{Sz\"oll\H{o}si construction}
\index{block design}

\begin{theorem}
If $H\in M_N(\pm 1)$ with $N\geq 8$ is dephased symmetric Hadamard, and
$$w=\frac{(1\pm i\sqrt{N-5})^2}{N-4}$$
then the following procedure yields a complex Hadamard matrix $M\in M_{N-1}(\mathbb T)$:
\begin{enumerate}
\item Erase the first row and column of $H$.

\item Replace all diagonal $1$ entries with $-w$.

\item Replace all off-diagonal $-1$ entries with $w$.
\end{enumerate} 
\end{theorem}

\begin{proof}
We know from chapter 1 that the scalar product between any two rows of $H$, normalized as there, appears as follows:
\begin{eqnarray*}
P
&=&\frac{N}{4}\cdot1\cdot1+\frac{N}{4}\cdot1\cdot(-1)+\frac{N}{4}\cdot(-1)\cdot1+\frac{N}{4}\cdot(-1)\cdot(-1)\\
&=&0
\end{eqnarray*}

Let us peform now the above operations (1,2,3), in reverse order. When replacing $-1\to w$, all across the matrix, the above scalar product becomes:
\begin{eqnarray*}
P'
&=&\frac{N}{4}\cdot1\cdot1+\frac{N}{4}\cdot1\cdot\bar{w}+\frac{N}{4}\cdot w\cdot1+\frac{N}{4}\cdot(-1)\cdot(-1)\\
&=&\frac{N}{2}(1+Re(w))
\end{eqnarray*}

Now when adjusting the diagonal via $w\to-1$ back, and $1\to-w$, this amounts in adding the quantity $-2(1+Re(w))$ to our product. Thus, our product becomes:
\begin{eqnarray*}
P''
&=&\left(\frac{N}{2}-2\right)(1+Re(w))\\
&=&\frac{N-4}{2}\left(1+\frac{6-N}{N-4}\right)\\
&=&1
\end{eqnarray*}

Finally, erasing the first row and column amounts in substracting 1 from our scalar product. Thus, our scalar product becomes $P'''=1-1=0$, and we are done.
\end{proof}

Observe that the number $w$ in the above statement is a root of unity precisely at $N=8$, where the only matrix satisfying the conditions in the statement is the Walsh matrix $W_8$. So, let us apply, as in \cite{sz2}, the above construction to this matrix, namely:
$$W_8=\begin{pmatrix}
1&1&1&1&1&1&1&1\\
1&-1&1&-1&1&-1&1&-1\\
1&1&-1&-1&1&1&-1&-1\\
1&-1&-1&1&1&-1&-1&1\\
1&1&1&1&-1&-1&-1&-1\\
1&-1&1&-1&-1&1&-1&1\\
1&-1&-1&-1&-1&-1&1&1\\
1&-1&-1&1&-1&1&1&-1
\end{pmatrix}$$

We obtain in this way the following matrix:
$$W_8'=\begin{pmatrix}
*&*&*&*&*&*&*&*\\
*&-1&1&w&1&w&1&w\\
*&1&-1&w&1&1&w&w\\
*&w&w&-w&1&w&w&1\\
*&1&1&1&-1&w&w&w\\
*&w&1&w&w&-w&w&1\\
*&1&w&w&w&w&-w&1\\
*&w&w&1&w&1&1&-1
\end{pmatrix}$$

The Hadamard matrix obtained in this way, by deleting the $*$ entries, is the Petrescu matrix $P_7$, found in \cite{pet}. Thus, we have the following result: 

\index{Petrescu matrix}

\begin{theorem}
$P_7$ is the unique matrix formed by roots of unity that can be obtained by the Sz\"oll\H{o}si construction. It appears at $N=8$, from $H=W_8$. Its formula is
$$(P_7)_{ijk,abc}=
\begin{cases}
-w&{\rm if}\ (ijk)=(abc),\ ia+jb+kc=0(2)\\
w&{\rm if}\ (ijk)\neq(abc),\ ia+jb+kc\neq 0(2)\\
(-1)^{ia+jb+kc}&{\rm otherwise}
\end{cases}$$
where $w=e^{2\pi i/3}$, and with the indices belonging to the set $\{0,1\}^3-\{(0,0,0)\}$.
\end{theorem}

\begin{proof}
We know that the Sz\"oll\H{o}si construction maps $W_8\to P_7$. Since the formula of the second Fourier matrix is $(F_2)_{ij}=(-1)^{ij}$, the formula of the Walsh matrix $W_8$ is:
$$(W_8)_{ijk,abc}=(-1)^{ia+jb+kc}$$

But this gives the formula in the statement.
\end{proof}

Now observe that we are in the quite special situation $H=F_2\otimes K$, with $K$ being dephased and symmetric. Thus, we can search for a one-parameter affine deformation $K(q)$ which is dephased and symmetric, and then build the following matrix:
$$H(q)=\begin{pmatrix}K(q)&K\\ K&-K(\bar{q})\end{pmatrix}$$

In our case, such a deformation $K(q)=W_4(q)$ can be obtained by putting the $q$ parameters in the $2\times 2$ middle block. Now by performing the Sz\"oll\H{o}si construction, with the parameters $q,\bar{q}$ left untouched, we obtain the parametric Petrescu matrix \cite{pet}:

\index{Petrescu matrix}

\begin{theorem}
The following is a complex Hadamard matrix,
$$P_7^q
=\begin{pmatrix}
-q&q&w&1&w&1&w\\
q&-q&w&1&1&w&w\\
w&w&-w&1&w&w&1\\
1&1&1&-1&w&w&w\\
w&1&w&w&-\bar{q}w&\bar{q}w&1\\
1&w&w&w&\bar{q}w&-\bar{q}w&1\\
w&w&1&w&1&1&-1
\end{pmatrix}$$
where $w=e^{2\pi i/3}$, and $q\in\mathbb T$.
\end{theorem}

\begin{proof}
This follows from the above considerations, or from a direct verification of the orthogonality of the rows, which uses either $1-1=0$, or $1+w+w^2=0$.
\end{proof}

Observe that the above matrix $P_7^q$ has the property of being ``regular'', in the sense that the scalar products between rows appear from vanishing sums of roots of unity, possibly rotated by a scalar. We will be back to this in the next chapter, with the conjectural statement that $F_7,P_7^q$ are the only regular Hadamard matrices at $N=7$.

\section*{5e. Exercises} 

In connection with the Fourier matrices, we first have:

\begin{exercise}
Prove the following formula, with $w=e^{2\pi i/N}$,
$$\frac{1}{N}\sum_kw^{jk}=\delta_{0j}$$
where all the indices, and the Kronecker symbol too, are taken modulo $N$.
\end{exercise}

This is something that we have used in the above, in order to prove that $F_N$ is indeed Hadamard, and the argument there, which was quick and correct, was that the above average is the barycenter of the regular polygon formed by the numbers $w^{jk}$ in the complex plane, which is 0 generically, and is 1 if the polygon is degenerate. The problem now is that of finding another proof of this fact, by using abstract mathematics only. 

\begin{exercise}
Compute the determinant of the Fourier matrix $F_N$.
\end{exercise}

This certainly looks like something that can be done, by using standard linear algebra tricks. The problem is that of finding the trick which applies.

\begin{exercise}
Diagonalize the Fourier matrix $F_N$.
\end{exercise}

There is actually a lot of work here, and the answer is not trivial. In case you do not find the answer, a study at $N=2,3,4,5,6$ will do too.

\begin{exercise}
Prove that the deformed Fourier matrices $F_4^s$ are not equivalent to each other, up to identifying $\{s,-s,\bar{s},-\bar{s}\}$.
\end{exercise}

A natural idea here would be to look for an invariant $\varphi$ of the complex Hadamard matrices, or rather of the equivalence classes of such matrices, which gives $\varphi(F_4^s)=s$, but this is not obvious. In the lack of a good idea here, the best is to assume $F_4^s\sim F_4^t$, do computations, and look for a contradiction. And in case all this leads you nowhere, do not worry, we will come back later to this problem, with a clever invariant.

\begin{exercise}
Find a simple formula for the Tao matrix $T_6$.
\end{exercise}

To be more precise, the problem here is that of finding a simple formula for $(T_6)_{ij}$, as function of $i,j$. This is actually quite difficult. We will be back to this.

\begin{exercise}
Prove that the Beauchamp-Nicoara matrix $BN_6^q$ is indeed Hadamard.
\end{exercise}

There are some computations to be done here, which do not look very difficult. In case you are done with them quickly, you can try then proving the converse, namely that any self-adjoint Hadamard matrix at $N=6$ is equivalent to a matrix of type $BN_6^q$. 

\chapter{Roots of unity}

\section*{6a. Basic obstructions}

Many interesting examples of complex Hadamard matrices $H\in M_N(\mathbb T)$, including the real ones $H\in M_N(\pm1)$, have as entries roots of unity, of finite order. We discuss here this case, and more generally the ``regular'' case, where the combinatorics of the scalar products between the rows comes from vanishing sums of roots of unity. Let us begin with the following definition, going back to the work of Butson \cite{but}:

\index{Butson matrix}
\index{level}
\index{roots of unity}

\begin{definition}
An Hadamard matrix is called of Butson type if its entries are roots of unity of finite order. The Butson class $H_N(l)$ consists of the Hadamard matrices
$$H\in M_N(\mathbb Z_l)$$
where $\mathbb Z_l$ is the group of the $l$-th roots of unity. The level of a Butson matrix $H\in M_N(\mathbb T)$ is the smallest integer $l\in\mathbb N$ such that $H\in H_N(l)$.
\end{definition}

As basic examples, we have the real Hadamard matrices, which form the Butson class $H_N(2)$. The Fourier matrices are Butson matrices as well, because we have $F_N\in H_N(N)$, and more generally $F_G\in H_N(l)$, with $N=|G|$, and with $l\in\mathbb N$ being the smallest common order of the elements of $G$. There are many other examples, as for instance most of those at $N=6$ discussed in chapter 5, at 1 values of the various parameters $q,r,s$ there.

\bigskip

Generally speaking, the main question regarding the Butson matrices is that of understanding when $H_N(l)\neq 0$, via a theorem providing obstructions, and then a result or conjecture stating that these obstructions are the only ones. Let us begin with:

\index{Sylvester obstruction}

\begin{proposition}[Sylvester obstruction]
The following holds, 
$$H_N(2)\neq\emptyset\implies N\in\{2\}\cup 4\mathbb N$$
due to the orthogonality of the first $3$ rows.
\end{proposition}

\begin{proof}
This is something that we know from chapter 1, with the obstruction, going back to Sylvester's paper \cite{syl}, being explained there.
\end{proof}

The above obstruction is fully satisfactory, because according to the HC, its converse should hold. Thus, we are fully done with the case $l=2$. Our purpose now will be that of finding analogous statements at $l\geq3$, theorem plus conjecture. At very small values of $l$ this is certainly possible, and in what regards the needed obstructions, we can get away with the following simple fact, from Butson \cite{but} and Winterhof \cite{win}:

\index{vanishing sum of roots}
\index{sum of roots}

\begin{proposition}
For a prime power $l=p^a$, the vanishing sums of $l$-th roots of unity
$$\lambda_1+\ldots+\lambda_N=0\quad,\quad\lambda_i\in\mathbb Z_l$$
appear as formal sums of rotated full sums of $p$-th roots of unity. 
\end{proposition}

\begin{proof}
This is something elementary, coming from basic number theory. Consider indeed the full sum of $p$-th roots of unity, taken in a formal sense:
$$S=\sum_{k=1}^p(e^{2\pi i/p})^k$$

Let also $w=e^{2\pi i/l}$, and for $r\in\{1,2,\ldots ,l/p\}$ let us denote by $S_p^r=w^r\cdot S$ the above formal sum of roots of unity, rotated by $w^r$:
$$S_p^r=\sum_{k=1}^pw^r(e^{2\pi i/p})^k$$

We must show that any vanishing sum of $l$-th roots of unity appears as a sum of such quantities $S_p^r$. For this purpose, consider the following map, which assigns to the abstract elements of the group ring $\mathbb Z[\mathbb Z_l]$ their precise numeric values, inside $\mathbb Z(w)\subset\mathbb C$:
$$\Phi:\mathbb Z[\mathbb Z_l]\to\mathbb Z(w)$$

Our claim is that the elements $\{S_p^r\}$ form a basis of the vector space $\ker\Phi$. In order to prove this claim, observe first that we have:
$$S_p^r\in\ker\Phi$$

Also, the elements $S_p^r$ are linearly independent, because the support of $S_p^r$ contains a unique element of the subset $\{1,2,\ldots ,p^{a-1}\}\subset\mathbb Z_l$, namely the element $r\in\mathbb Z_l$, so all the coefficients of a vanishing linear combination of such sums $S_p^r$ must vanish. Thus, we are left with proving that $\ker\Phi$ is spanned by the elements $\{S_p^r\}$. For this purpose, let us recall the well-known fact that the minimal polynomial of $w$ is as follows:
$$\frac{X^{p^{a}}-1}{X^{{p^{a-1}}}-1}=1+X^{p^{a-1}}+X^{2p^{a-1}}+\ldots+X^{(p-1)p^{a-1}}$$ 

We conclude that the dimension of $\ker\Phi$ is given by:
$$\dim(\ker\Phi)
=p^a-(p^a-p^{a-1})
=p^{a-1}$$

Now since this is exactly the number of the sums $S_p^r$, this finishes the proof of our claim. Thus, any vanishing sum of $l$-th roots of unity must be of the form $\sum\pm S_p^r$, and the above support considerations show the coefficients must be positive, as desired.
\end{proof}

We can now formulate a result in the spirit of Proposition 6.2, as follows:

\index{Butson obstruction}

\begin{proposition}[Butson obstruction]
The following holds, 
$$H_N(p^a)\neq\emptyset\implies N\in p\mathbb N$$
due to the orthogonality of the first $2$ rows.
\end{proposition}

\begin{proof}
This follows indeed from Proposition 6.3, because the scalar product between the first 2 rows of our matrix is a vanishing sum of $l$-th roots of unity.
\end{proof}

WIth these obstructions in hand, we can discuss the case $l\leq5$, as follows:

\begin{theorem}
We have the following results,
\begin{enumerate}
\item $H_N(2)\neq\emptyset\implies N\in\{2\}\cup 4\mathbb N$,

\item $H_N(3)\neq\emptyset\implies N\in3\mathbb N$,

\item $H_N(4)\neq\emptyset\implies N\in2\mathbb N$,

\item $H_N(5)\neq\emptyset\implies N\in5\mathbb N$,
\end{enumerate}
with in cases $(1,3)$, a conjecture stating that the converse should hold as well.
\end{theorem}

\begin{proof}
In this statement (1) is the Sylvester obstruction, and (2,3,4) are particular cases of the Butson obstruction. As for the last assertion, which is of course something rather informal, but which is important for our purposes, the situation is as follows:

\medskip

(1) At $l=2$, as already mentioned, we have the Hadamard Conjecture, which comes with solid evidence, as explained in chapter 1 above.

\medskip

(2) At $l=4$ we have an old conjecture, dealing with complex Hadamard matrices over $\{\pm1,\pm i\}$, going back to the work of Turyn in \cite{tur}, and called Turyn Conjecture.
\end{proof}

\index{Turyn Conjecture}

At $l=3$ things are complicated, due to the following result of de Launey \cite{del}:

\index{de Launey obstruction}

\begin{proposition}[de Launey obstruction]
The following holds, 
$$H_N(l)\neq\emptyset\implies\exists\,d\in\mathbb Z[e^{2\pi i/l}],\,|d|^2=N^N$$
due to the orthogonality of all $N$ rows. In particular, we have
$$5|N\implies H_N(6)=\emptyset$$
so in particular $H_{15}(3)=\emptyset$, showing that the Butson obstruction is too weak at $l=3$.
\end{proposition}

\begin{proof}
The obstruction follows from the unitarity condition $HH^*=N$ for the complex Hadamard matrices, by applying the determinant, which gives:
$$|{\rm det}(H)|^2=N^N$$

Regarding the second assertion, let $w=e^{2\pi i/3}$, and assume that $d=a+bw+cw^2$ with $a,b,c\in\mathbb Z$ satisfies $|d|^2=0(5)$. We have the following computation:
\begin{eqnarray*}
|d|^2
&=&(a+bw+cw^2)(a+bw^2+cw)\\
&=&a^2+b^2+c^2-ab-bc-ac\\
&=&\frac{1}{2}[(a-b)^2+(b-c)^2+(c-a)^2]
\end{eqnarray*}

Thus our condition $|d|^2=0(5)$ leads to the following system, modulo 5:
$$x+y+z=0$$
$$x^2+y^2+z^2=0$$

But this system has no solutions. Indeed, let us look at $x^2+y^2+z^2=0$:

\medskip

(1) If this equality appears as $0+0+0=0$ we can divide $x,y,z$ by $5$ and redo the computation.

\medskip

(2) Otherwise, this equality can only appear as $0+1+(-1)=0$. 

\medskip

Thus, modulo permutations, we must have $x=0,y=\pm1,z=\pm2$, which contradicts $x+y+z=0$. Finally, the last assertion follows from $H_{15}(3)\subset H_{15}(6)=\emptyset$.
\end{proof}

At $l=5$ now, things are a bit unclear, with the converse of Theorem 6.5 (4) being something viable, at the conjectural level, at least to our knowledge. At $l=6$, however, the situation becomes again complicated, as follows:

\index{Haagerup obstruction}

\begin{proposition}[Haagerup obstruction]
The following holds, due to Haagerup's $N=5$ classification result, involving the orthogonality of all $5$ rows of the matrix:
$$H_5(l)\neq\emptyset\implies 5|l$$
In particular we have $H_5(6)=\emptyset$, which follows by the way from the de Launey obstruction as well, in contrast with the fact that we generally have $H_N(6)\neq\emptyset$.
\end{proposition}

\begin{proof}
In this statement the obstruction $H_5(l)=\emptyset\implies 5|l$ comes indeed from Haagerup's classification result in \cite{ha1}, explained in chapter 5. As for the last assertion, this is something informal, the situation at small values of $N$ being as follows:

\medskip

-- At $N=2,3,4$ we have the matrices $F_2,F_3,W_4$.

\medskip

-- At $N=6,7,8,9$ we have the matrices $F_6,P_7^1,W_8,F_3\otimes F_3$.

\medskip

-- At $N=10$ we have the following matrix, found in \cite{bbs} by using a computer, and written in logarithmic form, with $k$ standing for $e^{k\pi i/3}$:
$$X^6_{10}=
\left(\begin{array}{cccccccccccccc}
0&0&0&0&0&0&0&0&0&0\\
0&4&1&5&3&1&3&3&5&1\\
0&1&2&3&5&5&1&3&5&3\\
0&5&3&2&1&5&3&5&3&1\\
0&3&5&1&4&1&1&5&3&3\\
0&3&3&3&3&3&0&0&0&0\\
0&1&1&5&3&4&3&0&2&4\\
0&1&5&3&5&2&4&3&2&0\\
0&5&3&5&1&2&0&2&3&4\\
0&3&5&1&1&4&4&2&0&3
\end{array}\right)$$

We refer to \cite{bbs} for more details on this topic.
\end{proof}

All this is not good news. Indeed, there is no hope of conjecturally solving our $H_N(l)\neq\emptyset$ problem in general, because this would have to take into account, and in a simple and conceptual way, both the subtle arithmetic consequences of the de Launey obstruction, and the Haagerup classification result at $N=5$, and this does not seem feasible.

\section*{6b. Sums of roots}

Let us discuss now a generalization of the Butson obstruction from Proposition 6.4, which has been our main source of obstructions, so far. Let us start with:

\index{cycle}
\index{sum of cycles}

\begin{definition}
A cycle is a full sum of roots of unity, possibly rotated by a scalar, 
$$C=q\sum_{k=1}^lw^k\quad,\quad w=e^{2\pi i/l}\quad,\quad q\in\mathbb T$$
and taken in a formal sense. A sum of cycles is a formal sum of cycles.
\end{definition}

The actual sum of a cycle, or of a sum of cycles, is of course 0. This is why the word ``formal'' is there, for reminding us that we are working with formal sums. As an example, here is a sum of cycles, with $w=e^{2\pi i/6}$, and with $|q|=1$:
$$1+w^2+w^4+qw+qw^4=0$$

We know from Proposition 6.3 above that any vanishing sum of $l$-th roots of unity must be a sum of cycles, at least when $l=p^a$ is a prime power. However, this is not the case in general, the simplest counterexample being as follows, with $w=e^{2\pi i/30}$:
$$w^5+w^6+w^{12}+w^{18}+w^{24}+w^{25}=0$$

Indeed, this sum is obviously not a sum a cycles. However, this sum vanishes indeed, as shown by the following computation:
\begin{eqnarray*}
w^5+w^6+w^{12}+w^{18}+w^{24}+w^{25}
&=&w^5+w^{15}+w^{25}\\
&+&w^0+w^6+w^{12}+w^{18}+w^{24}\\
&-&w^0-w^{15}\\
&=&0+0-0\\
&=&0
\end{eqnarray*}

The following deep result on the subject is due to Lam and Leung \cite{lle}:

\index{Lam-Leung theorem}
\index{vanishing sum of roots}
\index{sum of roots of unity}

\begin{theorem}
Let $l=p_1^{a_1}\ldots p_k^{a_k}$, and assume that $\lambda_i\in\mathbb Z_l$ satisfy:
$$\lambda_1+\ldots+\lambda_N=0$$
\begin{enumerate}
\item $\sum\lambda_i$ is a sum of cycles, with $\mathbb Z$ coefficients.

\item If $k\leq 2$ then $\sum\lambda_i$ is a sum of cycles, with $\mathbb N$ coefficients.

\item If $k\geq 3$ then $\sum\lambda_i$ might not decompose as a sum of cycles.

\item $\sum\lambda_i$ has the same length as a sum of cycles: $N\in p_1\mathbb N+\ldots+p_k\mathbb N$.
\end{enumerate}
\end{theorem}

\begin{proof}
This is something that we will not really need in what follows, but that we included here, in view of its importance. The idea of the proof is as follows:

\medskip

(1) This is a well-known result, which follows from basic number theory, by using arguments in the spirit of those in the proof of Proposition 6.3.

\medskip

(2) This is something that we already know at $k=1$, from Proposition 6.3. At $k=2$ the proof is more technical, along the same lines. See \cite{lle}.

\medskip

(3) The smallest possible $l$ potentially producing a counterexample is $l=2\cdot3\cdot 5=30$, and we have here indeed the sum given above, with $w=e^{2\pi i/30}$.

\medskip

(4) This is a deep result, due to Lam and Leung, relying on advanced number theory knowledge. We refer to their paper \cite{lle} for the proof.
\end{proof}

As a side comment here, with such results we are now into rather advanced number theory. We warmly recommend at this point the reading of the paper of Lam-Leung \cite{lle}, not that we will really need this in what follows, but for getting a taste of the subject. As a consequence now of the above result, we have the following generalization of the Butson obstruction, which is something final and optimal on this subject:

\index{Lam-Leung obstruction}

\begin{theorem}[Lam-Leung obstruction]
Assuming the we have
$$l=p_1^{a_1}\ldots p_k^{a_k}$$
the following must hold, due to the orthogonality of the first $2$ rows:
$$H_N(l)\neq\emptyset\implies N\in p_1\mathbb N+\ldots+p_k\mathbb N$$
In the case $k\geq2$, the latter condition is automatically satisfied at $N>>0$.
\end{theorem}

\begin{proof}
Here the first assertion, which generalizes the $l=p^a$ obstruction from Proposition 6.4 above, comes from Theorem 6.9 (4), applied to the vanishing sum of $l$-th roots of unity coming from the scalar product between the first 2 rows. As for the second assertion, this is something well-known, coming from basic number theory.
\end{proof}

Summarizing, our study so far of the condition $H_N(l)\neq\emptyset$ has led us into an optimal obstruction coming from the first 2 rows, namely the Lam-Leung one, then an obstruction coming from the first 3 rows, namely the Sylvester one, and then two subtle obstructions coming from all $N$ rows, namely the de Launey one, and the Haagerup one. As an overall conclusion, by contemplating all these obstructions, nothing good in relation with our problem $H_N(l)\neq\emptyset$ is going on at small $N$. So, as a natural and more modest objective, we should perhaps try instead to solve this problem at $N>>0$. 

\bigskip

The point indeed is that everything simplifies at $N>>0$, with some of the above obstructions dissapearing, and with some other known obstructions, not to be discussed here, dissapearing as well. We are therefore led to the following statement:

\index{Asymptotic Butson Conjecture}

\begin{conjecture}[Asymptotic Butson Conjecture (ABC)]
The following equivalences should hold, in an asymptotic sense, at $N>>0$,
\begin{enumerate}
\item $H_N(2)\neq\emptyset\iff 4|N$,

\item $H_N(p^a)\neq\emptyset\iff p|N$, for $p^a\geq3$ prime power,

\item $H_N(l)\neq\emptyset\iff\emptyset$, for $l\in\mathbb N$ not a prime power,
\end{enumerate}
modulo the de Launey obstruction, $|d|^2=N^N$ for some $d\in\mathbb Z[e^{2\pi i/l}]$.
\end{conjecture}

In short, our belief is that when imposing the condition $N>>0$, only the Sylvester, Butson and de Launey obstructions survive. This is of course something quite nice, but in what regards a possible proof, this looks difficult. Indeed, our above conjecture generalizes the HC in the $N>>0$ regime, which is so far something beyond reach. One idea, however, in dealing with such questions, coming from the de Launey-Levin result from \cite{dle}, is that of looking at the partial Butson matrices, at $N>>0$. Observe in particular that restricting the attention to the rectangular case, and this not even in  the $N>>0$ regime, would make dissapear the de Launey obstruction from the ABC, which uses the orthogonality of all $N$ rows. We will discuss this later. For a number of related considerations, we refer as well to de Launey \cite{del} and de Launey-Gordon \cite{dgo}.

\section*{6c. Regularity}

Getting away now from all the above arithmetic difficulties, let us discuss, following \cite{bbs}, the classification of the regular complex Hadamard matrices of small order. The definition here, which already appeared in the above, is as follows:

\index{regular matrix}

\begin{definition}
A complex Hadamard matrix $H\in M_N(\mathbb T)$ is called regular if the scalar products between rows decompose as sums of cycles.
\end{definition}

We should mention that there is some terminology clash here, with the word ``regular'' being sometimes used in order to designate the bistochastic matrices. In this book we use the above notion of regularity, and we call bistochastic the bistochastic matrices.

\bigskip

Our purpose in what follows will be that of showing that the notion of regularity can lead to full classification results at $N\leq6$, and perhaps at $N=7$ too, and all this while covering most of the interesting complex Hadamard matrices that we met, so far. As a first observation, supporting this last claim, we have the following result:

\begin{proposition}
The following complex Hadamard matrices are regular:
\begin{enumerate}
\item The matrices at $N\leq5$, namely $F_2,F_3,F_4^s,F_5$.

\item The main examples at $N=6$, namely $F_6^{(rs)},F_6^{(^r_s)},H_6^q,T_6$.

\item The main examples at $N=7$, namely $F_7,P_7^q$.
\end{enumerate}
\end{proposition}

\begin{proof}
The Fourier matrices $F_N$ are all regular, with the scalar products between rows appearing as certain sums of full sums of $l$-th roots of unity, with $l|N$. As for the other matrices appearing in the statement, with the convention that ``cycle structure'' means the lengths of the cycles in the regularity property, the situation is as follows:

\medskip

(1) $F_4^s$ has cycle structure $2+2$, and this because the verification of the Hadamard condition is always based on the formula $1+(-1)=0$, rotated by scalars.

\medskip

(2) $F_6^{(rs)},F_6^{(^r_s)}$ have mixed cycle structure $2+2+2/3+3$, in the sense that both cases appear, $H_6^q$ has cycle structure $2+2+2$, and $T_6$ has cycle structure $3+3$.

\medskip

(3) $P_7^q$ has cycle structure $3+2+2$, its Hadamard property coming from $1+w+w^2=0$, with $w=e^{2\pi i/3}$, and from $1+(-1)=0$, applied twice, rotated by scalars.
\end{proof}
 
Let us discuss now the classification of regular matrices. We first have:

\begin{theorem}
The regular Hadamard matrices at $N\leq 5$ are 
$$F_2,F_3,F_4^s,F_5$$
up to the equivalence relation for the complex Hadamard matrices. 
\end{theorem}

\begin{proof}
This is something that we already know, coming from the classification results from chapter 5, and from Proposition 6.13 (1). However, and here comes our point, proving this result does not need in fact all this, the situation being as follows:

\medskip

(1) At $N=2$ the cycle structure can be only 2, and we obtain $F_2$.

\medskip

(2) At $N=3$ the cycle structure can be only 3, and we obtain $F_3$.

\medskip

(3) At $N=4$ the cycle structure can be only $2+2$, and we obtain $F_4^s$.

\medskip

(4) At $N=5$ some elementary combinatorics shows that the cycle structure $3+2$ is excluded. Thus we are left with the cycle structure $5$, and we obtain $F_5$.
\end{proof}

Let us discuss now the classification at $N=6$. The result here, from \cite{bbs}, states that the matrices $F_6^{(rs)},F_6^{(^r_s)},H_6^q,T_6$ are the only solutions. The proof is quite long and technical, but we will present here its main ideas. Let us start with:

\begin{proposition}
The regular Hadamard matrices at $N=6$ fall into $3$ classes:
\begin{enumerate}
\item Cycle structure $3+3$, with $T_6$ being an example.

\item Cycle structure $2+2+2$, with $H_6^q$ being an example.

\item Mixed cycle structure $3+3/2+2+2$, with $F_6^{(rs)},F_6^{(^r_s)}$ being examples.
\end{enumerate}
\end{proposition}

\begin{proof}
This is a bit of an empty statement, with the above (1,2,3) possibilities being the only ones, and with the various examples coming from Proposition 6.13 (2).
\end{proof}

In order to do the classification, we must prove that the examples in (1,2,3) are the only ones. Let us start with the Tao matrix. The result here is as follows:

\index{Tao matrix}

\begin{proposition}
The Tao matrix, namely
$$T_6=\begin{pmatrix}
1&1&1&1&1&1\\ 
1&1&w&w&w^2&w^2\\ 
1&w&1&w^2&w^2&w\\
1&w&w^2&1&w&w^2\\ 
1&w^2&w^2&w&1&w\\ 
1&w^2&w&w^2&w&1
\end{pmatrix}$$
with $w=e^{2\pi i/3}$ is the only one with cycle structure $3+3$.
\end{proposition}

\begin{proof}
The proof of this fact, from \cite{bbs}, is quite long and technical, the idea being that of studying first the $3\times 6$ case, then the $4\times6$ case, and finally the $6\times6$ case:

\medskip

(1) Consider first a partial Hadamard matrix $A\in M_{3\times 6}(\mathbb T)$, with the scalar products between rows assumed to be all of type $3+3$. By doing some elementary combinatorics, explained in \cite{bbs}, we can see that, modulo equivalence, either all entries of $A$ belong to $\mathbb Z_3=\{1,w,w^2\}$, or $A$ has the following special form, for certain parameters $r,s\in\mathbb T$:
$$A=\begin{pmatrix}
1&1&1&1&1&1\\
1&w&w^2&r&wr&w^2r\\
1&w^2&w&s&w^2s&ws
\end{pmatrix}$$

(2) With this result in hand, we can now investigate the $4\times6$ case.  Assume indeed that we have a partial Hadamard matrix $B\in M_{4\times 6}(\mathbb T)$, with the scalar products between rows assumed to be all of type $3+3$. By looking at the 4 submatrices $A^{(1)},A^{(2)},A^{(3)},A^{(4)}$ obtained from $B$ by deleting one row, and applying the above $3\times 6$ result, we see that all the possible parameters dissapear. Thus, our matrix must be of the following type:
$$B\in M_{4\times 6}(\mathbb Z_3)$$

(3) With this, we can now go for the general case. Indeed, an Hadamard matrix $M\in M_6(\mathbb T)$ having cycle structure $3+3$ must be of the form $M\in M_6(\mathbb Z_3)$. But the study of such matrices is elementary, with $T_6$ as the only solution. See \cite{bbs}.
\end{proof}

Regarding now the Haagerup matrix, the result is similar, as follows:

\index{Haagerup matrix}

\begin{proposition}
The Haagerup matrix, namely
$$H_6^q=\begin{pmatrix}
1&1&1&1&1&1\\
1&-1&i&i&-i&-i\\ 
1&i&-1&-i&q&-q\\ 
1&i&-i&-1&-q&q\\
1&-i&\bar{q}&-\bar{q}&i&-1\\ 
1&-i&-\bar{q}&\bar{q}&-1&i
\end{pmatrix}$$
with $q\in\mathbb T$ is the only one with cycle structure $2+2+2$.
\end{proposition}

\begin{proof}
The proof here, from \cite{bbs}, uses the same idea as in the proof of Proposition 6.16, namely a detailed combinatorial study, by increasing the number of rows. First of all, the study of the $3\times 6$ partial Hadamard matrices with cycle structure $2+2+2$ leads, up to equivalence, to the following 4 solutions, with $q\in\mathbb T$ being a parameter:
$$A_1=\begin{pmatrix}
1&1&1&1&1&1\\
1&-i&1&i&-1&-1\\
1&-1&i&-i&q&-q
\end{pmatrix}$$
$$A_2=\begin{pmatrix}
1&1&1&1&1&1\\
1&1&-1&i&-1&-i\\
1&-1&q&-q&iq&-iq
\end{pmatrix}$$
$$A_3=\begin{pmatrix}
1&1&1&1&1&1\\
1&-1&i&-i&q&-q\\
1&-i&i&-1&-q&q
\end{pmatrix}$$
$$A_4=\begin{pmatrix}
1&1&1&1&1&1\\
1&-i&-1&i&q&-q\\
1&-1&-q&-iq&iq&q
\end{pmatrix}$$

With this result in hand, we can go directly for the $6\times6$ case. Indeed, a careful examination of the $3\times6$ submatrices, and of the way that different parameters can overlap vertically, shows that our matrix must have a $3\times 3$ block decomposition as follows:
$$M=\begin{pmatrix}
A&B&C\\
D&xE&yF\\
G&zH&tI
\end{pmatrix}$$

Here $A,\ldots,I$ are $2\times 2$ matrices over $\{\pm 1,\pm i\}$, and $x,y,z,t$ are in $\{1,q\}$. A more careful examination shows that the solution must be of the following form:
$$M=\begin{pmatrix}
A&B&C\\
D&E&qF\\
G&qH&qI
\end{pmatrix}$$

More precisely, the matrix must be as follows:
$$M=\begin{pmatrix}
1&1&1&1&1&1\\
1&1&-i&i&-1&-1\\ 
1&i&-1&-i&-q&q\\ 
1&-i&i&-1&-iq&iq\\
1&-1&q&-iq&iq&-q\\ 
1&-1&-q&iq&q&-iq
\end{pmatrix}$$

But this matrix is equivalent to $H_6^q$, and we are done. See \cite{bbs}.
\end{proof}

Regarding now the mixed case, where both $2+2+2$ and $3+3$ situations can appear, this is a bit more complicated. We can associate to any mixed Hadamard matrix $M\in M_6(\mathbb C)$ its ``row graph'', having the 6 rows as vertices, and with each edge being called ``binary'' or ``ternary'', depending on whether the corresponding scalar product is of type $2+2+2$ or $3+3$. With this convention, we have the following result:

\index{row graph}

\begin{proposition}
The row graph of a mixed matrix $M\in M_6(\mathbb C)$ can be:
\begin{enumerate}
\item Either the bipartite graph having $3$ binary edges.

\item Or the bipartite graph having $2$ ternary triangles.
\end{enumerate}
\end{proposition}

\begin{proof}
Let $X$ be the row graph in the statement. By doing some combinatorics, of rather elementary type, we are led to the following conclusions about $X$:

\medskip

-- $X$ has no binary triangle.

\smallskip

-- $X$ has no ternary square.

\smallskip

-- $X$ has at least one ternary triangle.

\medskip

With these results in hand, we see that there are only two types of squares in our graph $X$, namely those having 1 binary edge and 5 ternary edges, and those consisting of a ternary triangle, connected to the 4-th point with 3 binary edges. By looking at pentagons, then hexagons that can be built with these squares, we see that the above two types of squares cannot appear at the same time, at that at the level of hexagons, we have the two solutions in the statement. For details regarding all this, we refer to \cite{bbs}.
\end{proof}

We can now complete our classification results at $N=6$ with:

\begin{proposition}
The deformed Fourier matrices, namely
$$F_6^{(rs)}=\begin{pmatrix}
1&1&1&&1&1&1\\
1&w&w^2&&1&w&w^2\\
1&w^2&w&&1&w^2&w\\ 
\\
1&r&s&&-1&-r&-s\\
1&wr&w^2s&&-1&-wr&-w^2s\\
1&w^2r&ws&&-1&-w^2r&-ws
\end{pmatrix}$$
$$F_6^{(^r_s)}
=\begin{pmatrix}
1&1&&1&1&&1&1\\
1&-1&&1&-1&&1&-1\\
\\
1&r&&w&wr&&w^2&w^2r\\ 
1&-r&&w&-wr&&w^2&-w^2r\\
\\
1&s&&w^2&w^2s&&w&ws\\
1&-s&&w^2&-w^2s&&w&-ws
\end{pmatrix}$$
with $r,s\in\mathbb T$ are the only ones with mixed cycle structure.
\end{proposition}

\begin{proof}
According to Proposition 6.18, we have two cases:

\medskip

(1) Assume first that the row graph is the bipartite one with 3 binary edges. By permuting the rows, the upper $4\times6$ submatrix of our matrix must be as follows:
$$B=\begin{pmatrix}
1&1&1&1&1&1\\
1&w&w^2&r&wr&w^2r\\
1&w^2&w&s&w^2s&ws\\ 
1&1&1&t&t&t
\end{pmatrix}$$

Now since the scalar product between the first and the fourth row is binary, we must have $t=-1$, so the solution is:
$$B=\begin{pmatrix}
1&1&1&1&1&1\\
1&w&w^2&r&wr&w^2r\\
1&w^2&w&s&w^2s&ws\\ 
1&1&1&-1&-1&-1
\end{pmatrix}$$

We can use the same argument for finding the fifth and sixth row, by arranging the matrix formed by the first three rows such as the second, respectively third row consist only of 1's. This will make appear some parameters of the form $w,w^2,r,s$ in the extra row, and we obtain in this way a matrix which is equivalent to $F_6^{(rs)}$. See \cite{bbs}.

\medskip

(2) Assume now that the row graph is the bipartite one with 2 ternary triangles. By permuting the rows, the upper $4\times6$ submatrix of our matrix must be as follows:
$$B=\begin{pmatrix}
1&1&1&1&1&1\\
1&1&w&w&w^2&w^2\\ 
1&1&w^2&w^2&w&w\\
1&-1&r&-r&s&-s
\end{pmatrix}$$

We can use the same argument for finding the fifth and sixth row, and we conclude that the matrix is of the following type:
$$M=\begin{pmatrix}
1&1&1&1&1&1\\
1&1&w&w&w^2&w^2\\ 
1&1&w^2&w^2&w&w\\
1&-1&r&-r&s&-s\\
1&-1&a&-a&b&-b\\
1&-1&c&-c&d&-d
\end{pmatrix}$$

Now since the last three rows must form a ternary triangle, we conclude that the matrix must be of the following form:
$$M=\begin{pmatrix}
1&1&1&1&1&1\\
1&1&w&w&w^2&w^2\\ 
1&1&w^2&w^2&w&w\\
1&-1&r&-r&s&-s\\
1&-1&wr&-wr&w^2s&-w^2s\\
1&-1&w^2r&-w^2r&ws&-ws
\end{pmatrix}$$

But this matrix is equivalent to $F_6^{(^r_s)}$, and we are done. See \cite{bbs}.
\end{proof}

All this was quite technical, but good news, we are done with our study. Indeed, summing up all the above, we have proved the following theorem, from \cite{bbs}:

\index{regular matrix}

\begin{theorem}
The regular complex Hadamard matrices at $N=6$ are:
\begin{enumerate}
\item The deformations $F_6^{(rs)},F_6^{(^r_s)}$ of the Fourier matrix $F_6$.

\item The Haagerup matrix $H_6^q$.

\item The Tao matrix $T_6$.
\end{enumerate}
\end{theorem}

\begin{proof}
This follows indeed from the trichotomy from Proposition 6.15, and from the results in Proposition 6.16, Proposition 6.17 and Proposition 6.19. See \cite{bbs}.
\end{proof}

All this is quite nice, bringing some fresh air into the classification question for the complex Hadamard matrices at $N=6$, which is stuck, as explained in chapter 5. As a continuation of this, our belief is that the $N=7$ classification is doable as well. Here we have 3 possible cycle structures, namely $3+2+2$, $5+2$, $7$, and our first job is that of understanding what cycle structures are indeed possible, in practice. 

\bigskip

In order to deal with this latter question, we use the same idea as at $N=6$, namely looking at $3\times N$ submatrices. Let us start with the following definition:

\begin{definition}
Given numbers $p_i,q_i,r_i$ with $\sum p_i=\sum q_i=\sum r_i$, we write
$$(p_1+\ldots+p_k)\perp_{(r_1+\ldots+r_s)}(q_1+\ldots+q_l)$$
if there exist sums of cycles $P,Q$ having cycle structure $\sum p_i,\sum q_i$, such that the scalar product $R=<P,Q>$ vanishes, and has $\sum r_i$ as cycle structure. Otherwise, we write:
$$(p_1+\ldots+p_k)\not\perp_{(r_1+\ldots+r_s)}(q_1+\ldots+q_l)$$
If there are no numbers $r_i$ such that $\sum p_i\perp_{\sum r_i}\sum q_i$ holds, we write $\sum p_i\not\perp\sum q_i$.
\end{definition}

In other words, we write $\sum p_i\perp_{\sum r_i}\sum q_i$ if there exist complex numbers $a_k,b_k,c_k\in\mathbb T$ such that $\sum a_k,\sum b_k,\sum c_k$ have cycle structure $\sum p_i,\sum q_i,\sum r_i$ respectively, and such that $a_k\bar{b}_k=c_k$ for any $k$, and we use as well the related notations $\sum p_i\not\perp_{\sum r_i}\sum q_i$ and $\sum p_i\not\perp\sum q_i$, taken in an obvious sense. Now with these notions in hand, we have:

\begin{proposition}
Assume that $p,q\geq 3$ are primes.
\begin{enumerate}
\item If $p=q+2$ we have $p\not\perp(q+2)$.

\item If $p=q+2$ then $(p+2)\not\perp_{(q+2+2)}(p+2)$.

\item We have $(p+2)\not\perp_{(p+2)}(p+2)$.

\item If $p=q+4$ then $(q+2+2)\not\perp_p(q+2+2)$.

\item If $p=q+2$ then $(q+2+2)\not\perp_{(p+2)}(q+2+2)$.
\end{enumerate}
\end{proposition}

\begin{proof}
All this follows from some basic number theory, the idea being as follows:

\medskip

(1) By multiplying by scalars and permuting columns, we can assume that the $2\times p$ matrix formed by our sums is as follows, with $w=e^{2\pi i/q}$ and $\xi=e^{2\pi i/p}$:
$$\begin{pmatrix}
1&w&\ldots&w^{q-1}&a&-a\\
1&\xi^{r_1}&\ldots&\xi^{r_{q-1}}&\xi^s&\xi^t
\end{pmatrix}$$

Now since the scalar product between rows vanishes, we obtain:
$$a=\frac{1+w\xi^{-r_1}+\ldots+w^{q-1}\xi^{-r_{q-1}}}{\xi^{-t}-\xi^{-s}}$$

On the other hand we have $|a|=1$, and the equation $a=\bar{a}^{-1}$ reads:
$$\frac{1+w\xi^{-r_1}+\ldots+w^{q-1}\xi^{-r_{q-1}}}{\xi^{-t}-\xi^{-s}}
=\frac{\xi^t-\xi^s}{1+w^{-1}\xi^{r_1}+\ldots+w^{-q+1}\xi^{r_{q-1}}}$$

Now by developing, we obtain a formula of type $(q-2)+S=0$, where $S$ is a certain sum of $q^2-q+2$ roots of unity of order $pq$. Now since $pq$ has $k=2$ prime factors, Theorem 6.9 (2) applies, and shows that $(q-2)+S$ must be a sum of cycles. But since $q\geq 3$, some of the terms of $S$ must be roots of unity of order $p$, or of order $q$, and this shows that $\xi$ is a power of $w$ or vice versa, which is a contradiction, as desired.

\medskip

(2-5) Here the study goes along the same lines, with the needed technical ingredient being the well-known Galois theory fact that a number $\lambda\in\mathbb Z[w]$, where $w=e^{2\pi i/n}$, satisfies $|\lambda|=1$ precisely when it is of the form $\lambda=\pm w^k$, for some $k\in\mathbb N$.
\end{proof}

Getting back now to our $N=7$ questions, we have the following result:

\begin{proposition}
We have the following obstructions:
\begin{enumerate}
\item $7\not\perp(5+2)$.

\item $(5+2)\not\perp(5+2)$.

\item $7\not\perp(3+2+2)$.

\item $(5+2)\not\perp(3+2+2)$.
\end{enumerate}
\end{proposition}

\begin{proof}
This follows from Proposition 6.22, as follows:

\medskip

(1) This follows from Proposition 6.22 (1), at $p=7$.

\medskip

(2) We have indeed $(5+2)\not\perp_7(5+2)$ from Proposition 6.22 (1), $(5+2)\not\perp_{(5+2)}(5+2)$ from Proposition 6.22 (2), and $(5+2)\not\perp_{(3+2+2)}(5+2)$ from Proposition 6.22 (3).

\medskip

(3) First, $7\not\perp_7(3+2+2)$ is clear. Also, we have $7\not\perp_{(5+2)}(3+2+2)$ from Proposition 6.22 (1) and $7\not\perp_{(3+2+2)}(3+2+2)$ from Proposition 6.22 (4), and this gives the result.

\medskip

(4) We have $(5+2)\not\perp_7(3+2+2)$ from Proposition 6.22, $(5+2)\not\perp_{(5+2)}(3+2+2)$ from Proposition 6.22 (3), and $(5+2)\not\perp_{(3+2+2)}(3+2+2)$ from Proposition 6.22 (5).
\end{proof}

In the context of the regular complex Hadamard matrices $H\in M_7(\mathbb T)$, the above result shows that the cycle structure $5+2$ is excluded, and that the cases $3+2+2$ and $7$ cannot interact. Thus we have a dichotomy, and our conjecture is as follows:

\begin{conjecture}
The regular complex Hadamard matrices at $N=7$ are:
\begin{enumerate}
\item The Fourier matrix $F_7$.

\item The Petrescu matrix $P_7^q$.
\end{enumerate}
\end{conjecture}

\index{Petrescu matrix}

Regarding (1), one can show indeed that $F_7$ is the only matrix having cycle structure 7, with this being related to more general results of Hiranandani-Schlenker \cite{hsc}. As for (2), the problem is that of proving that $P_7^q$ is the only matrix having cycle structure $3+2+2$. The computations here are unfortunately far more involved than those at $N=6$, briefly presented above, and finishing the classification work here is not an easy question.

\bigskip

Besides the classification questions, there are as well a number of theoretical questions in relation with the notion of regularity, that we believe to be very interesting. We have for instance the following conjecture, going back to \cite{bbs}:

\index{regularity conjecture}

\begin{conjecture}[Regularity Conjecture]
The following hold:
\begin{enumerate}
\item Any Butson matrix $H\in M_N(\mathbb C)$ is regular.

\item Any regular matrix $H\in M_N(\mathbb C)$ is an affine deformation of a Butson matrix.
\end{enumerate}
\end{conjecture}

In order to comment on the first conjecture, let us recall from Theorem 6.9 that in the case where the level of the Butson matrix has at most 2 prime factors, $l=p^a$ or $l=p^aq^b$, any vanishing sum of roots of unity, and in particular the various scalar products between rows, decompose as a sum of cycles. Thus, in this case, the conjecture holds. 

\bigskip

The problem appears when the level $l$ has at least 3 prime factors, for instance when $l=30$. Here we have ``exotic'' vanishing sums of roots of unity, such as the following one, with $w=e^{2\pi i/30}$, discussed after Definition 6.8:
$$w^5+w^6+w^{12}+w^{18}+w^{24}+w^{25}=0$$

\index{exotic sum of roots}

To be more precise, our above conjecture (1) says that such an exotic vanishing sum of roots of unity cannot be used in order to construct a complex Hadamard matrix, as part of the arithmetics leading to the vanishing of the various scalar products between rows. This looks like a quite difficult question, coming however with substantial computer evidence. We have no idea on how to approach it, abstractly. See \cite{bbs}.

\bigskip

As for the second conjecture, (2) above, this simply comes from the known examples of regular Hadamard matrices, which all appear from certain Butson matrices, by inserting parameters, in an affine way. We will further discuss the notion of affine deformation, with some general results on the subject, in chapters 7-8 below.

\section*{6d. Partial matrices}

As already mentioned after Conjecture 6.11, one way of getting away from the above algebraic difficulties is by doing $N>>0$ analysis for the partial Hadamard matrices, with counting results in the spirit of those of de Launey-Levin \cite{dle}. Let us start with:

\index{partial Butson matrix}
\index{PBM}

\begin{definition}
A partial Butson matrix (PBM) is a matrix 
$$H\in M_{M\times N}(\mathbb Z_q)$$ 
having its rows pairwise orthogonal, where $\mathbb Z_q\subset\mathbb C^\times$ is the group of $q$-roots of unity.
\end{definition}

Two PBM are called equivalent if one can pass from one to the other by permuting the rows and columns, or by multiplying the rows and columns by numbers in $\mathbb Z_q$. Up to this equivalence, we can assume that $H$ is dephased, in the sense that its first row consists of $1$ entries only. We can also put $H$ in ``standard form'', as follows:

\index{standard form}

\begin{definition}
We say that that a partial Butson matrix $H\in M_{M\times N}(\mathbb Z_q)$ is in standard form if the low powers of 
$$w=e^{2\pi i/q}$$
are moved to the left as much as possible, by proceeding from top to bottom.
\end{definition}

Let us first try to understand the case $M=2$. Here a dephased partial Butson matrix $H\in M_{2\times N}(\mathbb Z_q)$ must look as follows, with $\lambda_i\in\mathbb Z_q$ satisfying $\lambda_1+\ldots+\lambda_N=0$:
$$H=\begin{pmatrix}1&\ldots&1\\ \lambda_1&\ldots&\lambda_N\end{pmatrix}$$

With $q=p_1^{k_1}\ldots p_s^{k_s}$, we must have, according to Lam and Leung \cite{lle}:
$$N\in p_1\mathbb N+\ldots+p_s\mathbb N$$

Observe however that at $s\geq 2$ this obstruction dissapears at $N\geq p_1p_2$. With this discussion made, let us get now into the prime power case. We have:

\begin{proposition}
When $q=p^k$ is a prime power, the standard form of the dephased partial Butson matrices at $M=2$ is
$$H=\begin{pmatrix}
1&1&\ldots&1&\ldots&\ldots&1&1&\ldots&1\\
\underbrace{1}_{a_1}&\underbrace{w}_{a_2}&\ldots&\underbrace{w^{q/p-1}}_{a_{q/p}}&\ldots&\ldots&\underbrace{w^{q-q/p}}_{a_1}&\underbrace{w^{q-q/p+1}}_{a_2}&\ldots&\underbrace{w^{q-1}}_{a_{q/p}}
\end{pmatrix}$$
where $w=e^{2\pi i/q}$ and where $a_1,\ldots,a_{q/p}\in\mathbb N$ are multiplicities, summing up to $N/p$.
\end{proposition}

\begin{proof}
Indeed, it is well-known that for $q=p^k$ the solutions of $\lambda_1+\ldots+\lambda_N=0$ with $\lambda_i\in\mathbb Z_q$ are, up to permutations of the terms, exactly those in the statement.
\end{proof}

Now with Proposition 6.28 in hand, we can prove:

\begin{theorem}
When $q=p^k$ is a prime power, the probability for a randomly chosen $M\in M_{2\times N}(\mathbb Z_q)$, with $N\in p\mathbb N$, $N\to\infty$, to be partial Butson is:
$$P_2\simeq\sqrt{\frac{p^{2-\frac{q}{p}}q^{q-\frac{q}{p}}}{(2\pi N)^{q-\frac{q}{p}}}}$$
\end{theorem}

\begin{proof}
First, the probability $P_M$ for a random $M\in M_{M\times N}(\mathbb Z_q)$ to be PBM is:
$$P_M=\frac{1}{q^{MN}}\#PBM_{M\times N}$$

Thus, according to Proposition 6.28, we have the following formula:
\begin{eqnarray*}
P_2
&=&\frac{1}{q^N}\sum_{a_1+\ldots +a_{q/p}=N/p}\binom{N}{\underbrace{a_1\ldots a_1}_p\ldots\ldots\underbrace{a_{q/p}\ldots a_{q/p}}_p}\\
&=&\frac{1}{q^N}\binom{N}{\underbrace{N/p\ldots N/p}_p}\sum_{a_1+\ldots +a_{q/p}=N/p}\binom{N/p}{a_1\ldots a_{q/p}}^p\\
&=&\frac{1}{p^N}\binom{N}{\underbrace{N/p\ldots N/p}_p}\times\frac{1}{(q/p)^N}\sum_{a_1+\ldots +a_{q/p}=N/p}\binom{N/p}{a_1\ldots a_{q/p}}^p
\end{eqnarray*}

Now by using the Stirling formula for the left term, and the basic multinomial sum estimate from chapter 4 with $s=q/p$ and $n=N/p$ for the right term, we obtain:
\begin{eqnarray*}
P_2
&=&\sqrt{\frac{p^p}{(2\pi N)^{p-1}}}\times\sqrt{\frac{(q/p)^{\frac{q}{p}(p-1)}}{p^{\frac{q}{p}-1}(2\pi N/p)^{(\frac{q}{p}-1)(p-1)}}}\\
&=&\sqrt{\frac{p^{p-\frac{q}{p}(p-1)-\frac{q}{p}+1+(\frac{q}{p}-1)(p-1)}q^{\frac{q}{p}(p-1)}}{(2\pi N)^{p-1+(\frac{q}{p}-1)(p-1)}}}\\
&=&\sqrt{\frac{p^{2-\frac{q}{p}}q^{q-\frac{q}{p}}}{(2\pi N)^{q-\frac{q}{p}}}}
\end{eqnarray*}

Thus we have obtained the formula in the statement, and we are done.
\end{proof}

Let us discuss now the case where $M=2$, and $q=p_1^{k_1}p_2^{k_2}$ has two prime factors. We first examine the simplest such case, namely $q=p_1p_2$, with $p_1,p_2$ primes:

\begin{proposition}
When $q=p_1p_2$ is a product of distinct primes, the standard form of the dephased partial Butson matrices at $M=2$ is
$$H=\begin{pmatrix}
1&1&\ldots&1&\ldots&\ldots&1&1&\ldots&1\\
\underbrace{1}_{A_{11}}&\underbrace{w}_{A_{12}}&\ldots&\underbrace{w^{p_2-1}}_{A_{1p_2}}&\ldots&\ldots&\underbrace{w^{q-p_2}}_{A_{p_11}}&\underbrace{w^{q-p_2+1}}_{A_{p_12}}&\ldots&\underbrace{w^{q-1}}_{A_{p_1p_2}}
\end{pmatrix}$$
where $w=e^{2\pi i/q}$, and $A\in M_{p_1\times p_2}(\mathbb N)$ is of the form $A_{ij}=B_i+C_j$, with $B_i,C_j\in\mathbb N$.
\end{proposition}

\begin{proof}
We use the fact that for $q=p_1p_2$ any vanishing sum of $q$-roots of unity decomposes as a sum of cycles. Now if we denote by $B_i,C_j\in\mathbb N$ the multiplicities of the various $p_2$-cycles and $p_1$-cycles, then we must have $A_{ij}=B_i+C_j$, as claimed.
\end{proof}

Regarding now the matrices of type $A_{ij}=B_i+C_j$, when taking them over integers, $B_i,C_j\in\mathbb Z$, these form a vector space of dimension $d=p_1+p_2-1$. Given $A\in M_{p_1\times p_2}(\mathbb Z)$, the ``test'' for deciding if we have $A_{ij}=B_i+C_j$ or not is:
$$A_{ij}+A_{kl}=A_{il}+A_{jk}$$

The problem comes of course from the assumption $B_i,C_j\geq0$, which is quite a subtle one. In what follows we restrict the attention to the case $p_1=2$. Here we have:

\index{random walk}

\begin{theorem}
For $q=2p$ with $p\geq 3$ prime, $P_2$ equals the probability for a random walk on $\mathbb Z^p$ to end up on the diagonal, i.e. at a position of type $(t,\ldots,t)$, with $t\in\mathbb Z$.
\end{theorem}

\begin{proof}
According to Proposition 6.30, we must understand the structure of the matrices $A\in M_{2\times p}(\mathbb N)$ which decompose as follows, with $B_i,C_j\geq0$:
$$A_{ij}=B_i+C_j$$

But this is an easy task, because depending on the value of $A_{11}$ compared to the value of $A_{21}$ we have 3 types of solutions, as follows:
$$\begin{pmatrix}
a_1&\ldots&a_p\\
a_1&\ldots&a_p
\end{pmatrix}\quad,\quad
\begin{pmatrix}
a_1&\ldots&a_p\\
a_1+t&\ldots&a_p+t
\end{pmatrix}\quad,\quad 
\begin{pmatrix}
a_1+t&\ldots&a_p+t\\
a_1&\ldots&a_p
\end{pmatrix}$$

Here $a_i\geq0$ and $t\geq1$. Now since cases 2,3 contribute in the same way, we obtain:
\begin{eqnarray*}
P_2
&=&\frac{1}{(2p)^N}\sum_{2\Sigma a_i=N}\binom{N}{a_1,a_1,\ldots,a_p,a_p}\\
&+&\frac{2}{(2p)^N}\sum_{t\geq1}\sum_{2\Sigma a_i+pt=N}\binom{N}{a_1,a_1+t,\ldots,a_p,a_p+t}
\end{eqnarray*}

We can write this formula in a more compact way, as follows:
$$P_2=\frac{1}{(2p)^N}\sum_{t\in\mathbb Z}\sum_{2\Sigma a_i+p|t|=N}\binom{N}{a_1,a_1+|t|,\ldots,a_p,a_p+|t|}$$

Now since the sum on the right, when rescaled by $\frac{1}{(2p)^N}$, is exactly the probability for a random walk on $\mathbb Z^p$ to end up at $(t,\ldots,t)$, this gives the result.
\end{proof}

According to the above result we have $P_2=\sum_{t\in\mathbb Z}P_2^{(t)}$, where $P_2^{(t)}$ with $t\in\mathbb Z$ is the probability for a random walk on $\mathbb Z^p$ to end up at $(t,\ldots,t)$. By using the basic binomial sum estimate of Richmond-Shallit \cite{rsh}, explained in chapter 4, we obtain:
\begin{eqnarray*}
P_2^{(0)}
&=&\frac{1}{(2p)^N}\binom{N}{N/2}\sum_{a_1+\ldots+a_p=N/2}\binom{N/2}{a_1,\ldots,a_p}^2\\
&\simeq&\sqrt{\frac{2}{\pi N}}\times\sqrt{\frac{p^p}{2^{p-1}(\pi N)^{p-1}}}\\
&=&2\sqrt{\left(\frac{p}{2\pi N}\right)^p}
\end{eqnarray*}

Regarding now the probability $P_2^{(t)}$ of ending up at $(t,\ldots,t)$, in principle for small $t$ this can be estimated by using a modification of the method in \cite{rsh}. However, it is not clear how to compute the full diagonal return probability in Theorem 6.31.

\bigskip

Let us discuss now the exponents $q=3p$. The same method as in the proof of Theorem 6.31 works, with the ``generic'' solution for $A$ being as follows:
$$A=\begin{pmatrix}
a_1&\ldots&a_p\\
a_1+t&\ldots&a_p+t\\
a_1+s+t&\ldots&a_p+s+t\\
\end{pmatrix}$$

More precisely, this type of solution, with $s,t\geq1$, must be counted 6 times, then its $s=0,t\geq1$ and $s\geq1,t=0$ particular cases must be counted 3 times each, and finally the $s=t=0$ case must be counted once. Observe that the $s=t=0$ contribution is:
\begin{eqnarray*}
P_3^{(0,0)}
&=&\frac{1}{(3p)^N}\binom{N}{N/3,N/3,N/3}\sum_{a_1+\ldots+a_p=N/3}\binom{N/3}{a_1,\ldots,a_p}^3\\
&\simeq&\sqrt{\frac{27}{(2\pi N)^2}}\times\sqrt{\frac{p^{2p}}{3^{p-1}(2\pi N/3)^{2(p-1)}}}\\
&=&3\sqrt{3^p}\left(\frac{p}{2\pi N}\right)^p
\end{eqnarray*} 

Finally, regarding arbitrary exponents with two prime factors, we have:

\begin{proposition}
When $q=p_1^{k_1}p_2^{k_2}$ has exactly two prime factors, the dephased partial Butson matrices at $M=2$ are indexed by the solutions of
$$A_{ij,xy}=B_{ijy}+C_{jxy}$$
with $B_{ijy},C_{jxy}\in\mathbb N$, with $i\in\mathbb Z_{p_1}$, $j\in\mathbb Z_{p_1^{k_1-1}}$, $x\in\mathbb Z_{p_2}$, $y\in\mathbb Z_{p_2^{k_2-1}}$.
\end{proposition}

\begin{proof}
We follow the method in the proof of Proposition 6.30. First, according to Lam-Leung \cite{lle}, for $q=p_1^{k_1}p_2^{k_2}$ any vanishing sum of $q$-roots of unity decomposes as a sum of cycles. Let us first work out a simple particular case, namely $q=4p$. Here the multiplicity matrices $A\in M_{4\times p}(\mathbb N)$ appear as follows:
$$A=\begin{pmatrix}B_1&\ldots&B_1\\ B_2&\ldots&B_2\\ B_3&\ldots&B_3\\ B_4&\ldots&B_4\end{pmatrix}+
\begin{pmatrix}C_1&\ldots&C_p\\ D_1&\ldots&D_p\\ C_1&\ldots&C_p\\ D_1&\ldots&D_p\end{pmatrix}$$

Thus, if we use double binary indices for the elements of $\{1,2,3,4\}$, the condition is:
$$A_{ij,x}=B_{ij}+C_{jx}$$

The same method works for any exponent of type $q=p_1^{k_1}p_2^{k_2}$, the formula being:
$$A_{i_1\ldots i_{k_1},x_1\ldots x_{k_2}}=B_{i_1\ldots i_{k_1},x_2\ldots x_{k_2}}+C_{i_2\ldots i_{k_1},x_1\ldots x_{k_2}}$$

But this gives the formula in the statement, and we are done.
\end{proof}

At $M=3$ now, we first restrict attention to the case where $q=p$ is prime. In this case, the general result in Proposition 6.32 becomes simply:
$$H=\begin{pmatrix}
1&1&\ldots&1\\
\underbrace{1}_a&\underbrace{w}_a&\ldots&\underbrace{w^{p-1}}_a
\end{pmatrix}$$

\index{tristochastic matrix}

We call a matrix $A\in M_p(\mathbb N)$ ``tristochastic'' if the sums on its rows, columns and diagonals are all equal. Here, and in what follows, we call ``diagonals'' the main diagonal, and its $p-1$ translates to the right, obtained by using modulo $p$ indices. With this convention, here is now the result at $M=3$:

\begin{proposition}
For $p$ prime, the standard form of dephased PBM at $M=3$ is
$$H=\begin{pmatrix}
1&1&\ldots&1&\ldots&\ldots&1&1&\ldots&1\\
1&1&\ldots&1&\ldots&\ldots&w^{p-1}&w^{p-1}&\ldots&w^{p-1}\\
\underbrace{1}_{A_{11}}&\underbrace{w}_{A_{12}}&\ldots&\underbrace{w^{p-1}}_{A_{1p}}&\ldots&\ldots&\underbrace{1}_{A_{p1}}&\underbrace{w}_{A_{p2}}&\ldots&\underbrace{w^{p-1}}_{A_{pp}}
\end{pmatrix}$$
where $w=e^{2\pi i/p}$ and where $A\in M_p(\mathbb N)$ is tristochastic, with sums $N/p$.
\end{proposition}

\begin{proof}
Consider a dephased matrix $H\in M_{3\times N}(\mathbb Z_p)$, written in standard form as in the statement. Then the orthogonality conditions between the rows are as follows:

\medskip

$1\perp 2$ means $A_{11}+\ldots+A_{1p}=A_{21}+\ldots+A_{2p}=\ldots\ldots=A_{p1}+\ldots+A_{pp}$.

\medskip

$1\perp 3$ means $A_{11}+\ldots+A_{p1}=A_{12}+\ldots+A_{p2}=\ldots\ldots=A_{1p}+\ldots+A_{pp}$.

\medskip

$2\perp 3$ means $A_{11}+\ldots+A_{pp}=A_{12}+\ldots+A_{p1}=\ldots\ldots=A_{1p}+\ldots+A_{p,p-1}$.

\medskip

Thus $A$ must have constant sums on rows, columns and diagonals, as claimed.
\end{proof}

It is quite unobvious on how to deal with the tristochastic matrices with bare hands. For the moment, let us just record a few elementary results:

\begin{proposition}
For $p=2,3$, the standard form of the dephased PBM at $M=3$ is respectively as follows, with $w=e^{2\pi i/3}$ and $a+b+c=N/3$ at $p=3$:
$$H=\begin{pmatrix}+&+&+&+\\+&+&-&-\\\underbrace{+}_{N/4}&\underbrace{-}_{N/4}&\underbrace{+}_{N/4}&\underbrace{-}_{N/4}\end{pmatrix}$$
$$H=\begin{pmatrix}
1&1&1&1&1&1&1&1&1\\
1&1&1&w&w&w&w^2&w^2&w^2\\
\underbrace{1}_a&\underbrace{w}_b&\underbrace{w^2}_c&\underbrace{1}_b&\underbrace{w}_c&\underbrace{w^2}_a&\underbrace{1}_c&\underbrace{w}_a&\underbrace{w^2}_b
\end{pmatrix}$$
Also, for $p\geq 3$ prime and $N\in p\mathbb N$, there is at least one Butson matrix $H\in M_{3\times N}(\mathbb Z_p)$.
\end{proposition}

\begin{proof}
The idea is that the $p=2$ assertion follows from Proposition 6.33, and from the fact that the $2\times 2$ tristochastic matrices are as follows:
$$A=\begin{pmatrix}a&a\\a&a\end{pmatrix}$$

As for the $p=3$ assertion, once again the idea is that this follows from Proposition 6.33, and from the fact that the $3\times 3$ tristochastic matrices are as follows:
$$A=\begin{pmatrix}a&b&c\\ b&c&a\\ c&a&b\end{pmatrix}$$

Indeed, the $p=2$ assertion is clear. Regarding now the $p=3$ assertion, consider an arbitary $3\times 3$ bistochastic matrix, written as follows:
$$A=\begin{pmatrix}a&b&n-a-b\\ d&c&n-c-d\\ n-a-d&n-b-c&*\end{pmatrix}$$

Here $*=a+b+c+d-n$, but we won't use this value, because one of the 3 diagonal equations is redundant anyway. With these notations in hand, the conditions are:
$$b+(n-c-d)+(n-a-d)=n$$
$$(n-a-b)+d+(n-b-c)=n$$

Now since substracting these equations gives $b=d$, we obtain the result. Regarding now the last assertion, consider the following $p\times p$ permutation matrix:
$$A=\begin{pmatrix}
1&&&&\\ 
&&&&1\\
&&&1\\
&&\ldots\\
&1
\end{pmatrix}$$

Since this matrix is tristochastic, for any $p\geq 3$ odd, this gives the result.
\end{proof}

Regarding now the asymptotic count, we have here:

\begin{theorem}
For $p=2,3$, the probability for a randomly chosen 
$$M\in M_{3\times N}(\mathbb Z_p)$$
with $N\in p\mathbb N$, $N\to\infty$, to be partial Butson is respectively given by
$$P_3^{(2)}\simeq\begin{cases}
\frac{16}{\sqrt{(2\pi N)^3}}&{\rm if}\ N\in4\mathbb N\\
0&{\rm if}\ N\notin 4\mathbb N\end{cases}$$
at $p=2$, and
$$P_3^{(3)}\simeq\frac{243\sqrt{3}}{(2\pi N)^3}$$
at $p=3$. In addition, we have $P_3^{(p)}>0$ for any $N\in p\mathbb N$, for any $p\geq 3$ prime.
\end{theorem}

\begin{proof}
According to Proposition 6.34, and then to the Stirling formula, we have:
$$P_3^{(2)}
=\frac{1}{4^N}\binom{N}{N/4,N/4,N/4,N/4}\\
\simeq\frac{16}{\sqrt{(2\pi N)^3}}$$

Similarly, by using the basic estimate with $s=p=3$, $n=N/3$, we have:
\begin{eqnarray*}
P_3^{(3)}
&=&\frac{1}{9^N}\sum_{a+b+c=N/3}\binom{N}{a,b,c,b,c,a,c,a,b}\\
&=&\frac{1}{3^N}\binom{N}{N/3,N/3,N/3}\times\frac{1}{3^N}\sum_{a+b+c=N/3}\binom{N/3}{a,b,c}^3\\
&\simeq&\frac{3\sqrt{3}}{2\pi N}\cdot\sqrt{\frac{81}{(2\pi N/3)^4}}\\
&=&\frac{243\sqrt{3}}{(2\pi N)^3}
\end{eqnarray*}

Finally, the last assertion is clear from the last assertion in Proposition 6.33.
\end{proof}

It is possible to establish a few more results in this direction, making interesting connections with probability. However, the main question regarding the partial Butson matrices remains that of adapting the asymptotic counting methods of de Launey-Levin \cite{dle} to the root of unity case. As a preliminary observation here, we have:

\index{random walk}
\index{de Launey-Levin}

\begin{proposition}
The probability $P_M$ for a random $H\in M_{M\times N}(\mathbb Z_q)$ to be partial Butson equals the probability for a length $N$ random walk with increments drawn from
$$E=\left\{(e_i\bar{e}_j)_{i<j}\Big|e\in\mathbb Z_q^M\right\}$$
regarded as a subset $\mathbb Z_q^{\binom{M}{2}}$, to return at the origin.
\end{proposition}

\begin{proof}
Indeed, with $T(e)=(e_i\bar{e}_j)_{i<j}$, a matrix $X=[e_1,\ldots,e_N]\in M_{M\times N}(\mathbb Z_q)$ is partial Butson if and only if the following condition is satisfied:
$$T(e_1)+\ldots+T(e_N)=0$$

But this leads to the conclusion in the statement.
\end{proof}

Observe now that, according to the above result, we have:
\begin{eqnarray*}
P_M
&=&\frac{1}{q^{(M-1)N}}\#\left\{\xi_1,\ldots,\xi_N\in E\Big|\sum_i\xi_i=0\right\}\\
&=&\frac{1}{q^{(M-1)N}}\sum_{\xi_1,\ldots,\xi_N\in E}\delta_{\sum\xi_i,0}
\end{eqnarray*}

The problem is to continue the computation in the proof of the inversion formula. More precisely, the next step at $q=2$, which is the key one, is as follows:
$$\delta_{\sum\xi_i,0}=\frac{1}{(2\pi)^D}\int_{[-\pi,\pi]^D}e^{i<\lambda,\sum\xi_i>}d\lambda$$

Here $D=\binom{M}{2}$. The problem is that this formula works when $\sum\xi_i$ is real, as is the case in the context of \cite{dle}, but not when $\sum\xi_i$ is complex, as is the case in Proposition 6.36. As before with other open questions, this is a good question for you, reader.

\section*{6e. Exercises} 

There are many interesting things to be done in connection with the roots of unity, and the corresponding Hadamard matrices, and here is a first exercise on this:

\begin{exercise}
Find the minimal polynomial of an arbitrary root of unity $w\in\mathbb T$.
\end{exercise}

This is standard algebra, that we used in the proof of the Butson obstruction.

\begin{exercise}
Develop the theory of the conjecture $H_{3n}(3)\neq\emptyset$, in analogy with the theory of the Hadamard conjecture, namely $H_{4n}(2)\neq\emptyset$.
\end{exercise}

This is of course a bit loosely formulated, the problem being that of finding some good results here, including evidence at small values of $n\in\mathbb N$, and so on.

\begin{exercise}
Prove that for any $l\in\mathbb N$, any vanishing sum of $l$-roots of unity appears as a sum of cycles, with $\mathbb Z$ coefficients.
\end{exercise}

This is something that we already discussed in the above, but very briefly, with the indication that this should follow from basic number theory, via arguments which are similar to those from the proof of the Butson obstruction.

\begin{exercise}
Prove that for $l=p^aq^b$, any vanishing sum of $l$-roots of unity appears as a sum of cycles.
\end{exercise}

To be more precise here, we already know that the conclusion in the statement holds in the case $l=p^a$. The problem is that of adapting that proof, from the case $l=p^a$, to the case $l=p^aq^b$. This is not exactly easy, but with some work, can be done.

\begin{exercise}
Read the proof of the Lam-Leung theorem, stating that the lenght of a vanishing sum of roots of unity should equal the length of a sum of cycles, and write down a brief account of that proof, explaining the main ideas there. 
\end{exercise}

Obviously, this is something quite time-consuming. However, this is worth the effort, the paper of Lam-Leung being an excellent introduction to advanced algebra.

\begin{exercise}
Work out all the details for the dichotomy in Proposition 6.18.
\end{exercise}

To be more precise here, Proposition 6.18 above comes with 1/2 page of proof, which is quite brief, and the problem is that of adding 1 page or so of details.

\begin{exercise}
Prove that the $7\times7$ regular matrices can only have
$$3+2+2\quad,\quad 5+2\quad,\quad 7$$
as cycle structure, then prove that the case $5+2$ is actually excluded.
\end{exercise}

Here the first assertion is something trivial, and the problem is that of finding the good number theoretic argument for excluding the case $5+2$.

\begin{exercise}
In the context of the previous exercise, prove that the cases
$$3+2+2\quad,\quad 7$$
do not interact, in the sense that a regular $7\times7$ Hadamard matrix has either all scalar products between the rows of type $3+2+2$, or of type $7$.
\end{exercise}

As before, with the previous exercise, the problem is that of finding the good number theoretic argument which applies, and gives the result.

\begin{exercise}
Prove that the Fourier matrix $F_7$ is the only $7\times7$ complex Hadamard matrix having cycle structure $7$.
\end{exercise}

This exercise is independent from the previous exercises, and is of different nature too, the problem here being not number theoretical, but rather purely combinatorial.

\chapter{Geometry, defect}

\section*{7a. Affine deformations}

We have seen so far that some theory for the complex Hadamard matrices $H\in M_N(\mathbb T)$ can be developed with some inspiration for the real case, $H\in M_N(\pm1)$, by looking at the Butson matrix case, $H\in M_N(\mathbb Z_l)$ with $l<\infty$. However, all this root of unity business ultimately leads into questions of HC flavor, and to put it squarely, wrong way. In this chapter we take a radically different approach to the study of the complex Hadamard matrices. Let us recall that the complex Hadamard manifold appears as:
$$X_N=M_N(\mathbb T)\cap\sqrt{N}U_N$$

This intersection is far from being smooth, and given a point $H\in X_N$, the problem is that of understanding the structure of $X_N$ around $H$, which is often singular. And this is what we will do, real algebraic geometry, for studying $X_N$ and its singularities. For this purpose, let us begin with some notations. We denote by $X_p$ an unspecified neighborhood of a point in a manifold, $p\in X$. Also, for $q\in\mathbb T_1$, meaning that $q\in\mathbb T$ is close to $1$, we define $q^r$ with $r\in\mathbb R$ by $(e^{it})^r=e^{itr}$. With these conventions, we have:

\begin{proposition}
For $H\in X_N$ and $A\in M_N(\mathbb R)$, the following are equivalent:
\begin{enumerate}
\item The following is an Hadamard matrix, for any $q\in\mathbb T_1$:
$$H_{ij}^q=H_{ij}q^{A_{ij}}$$

\item The following equations hold, for any $i\neq j$ and any $q\in\mathbb T_1$:
$$\sum_kH_{ik}\bar{H}_{jk}q^{A_{ik}-A_{jk}}=0$$

\item The following equations hold, for any $i\neq j$ and any $\varphi:\mathbb R\to\mathbb C$:
$$\sum_kH_{ik}\bar{H}_{jk}\varphi(A_{ik}-A_{jk})=0$$

\item For any $i\neq j$ and any $r\in\mathbb R$, with $E_{ij}^r=\{k|A_{ik}-A_{jk}=r\}$, we have:
$$\sum_{k\in E_{ij}^r}H_{ik}\bar{H}_{jk}=0$$
\end{enumerate}
If these conditions are satisfied, we call the matrix $H^q$ an affine deformation of $H$.
\end{proposition}

\begin{proof}
These equivalences are all elementary, and can be proved as follows:

\medskip

$(1)\iff(2)$ Indeed, the scalar products between the rows of $H^q$ are:
\begin{eqnarray*}
<H^q_i,H^q_j>
&=&\sum_kH_{ik}q^{A_{ik}}\bar{H}_{jk}\bar{q}^{A_{jk}}\\
&=&\sum_kH_{ik}\bar{H}_{jk}q^{A_{ik}-A_{jk}}
\end{eqnarray*}

$(2)\implies(4)$ This follows from the following formula, and from the fact that the power functions $\{q^r|r\in\mathbb R\}$ over the unit circle $\mathbb T$ are linearly independent:
$$\sum_kH_{ik}\bar{H}_{jk}q^{A_{ik}-A_{jk}}=\sum_{r\in\mathbb R}q^r\sum_{k\in E_{ij}^r}H_{ik}\bar{H}_{jk}$$

$(4)\implies(3)$ This follows from the following formula:
$$\sum_kH_{ik}\bar{H}_{jk}\varphi(A_{ik}-A_{jk})=\sum_{r\in\mathbb R}\varphi(r)\sum_{k\in E_{ij}^r}H_{ik}\bar{H}_{jk}$$

$(3)\implies(2)$ This simply follows by taking $\varphi(r)=q^r$.
\end{proof}

In order to understand the above deformations, which are ``affine'' in a certain sense, as suggested at the end of the statement, it is convenient to enlarge the attention to all types of deformations. We keep using the neighborhood notation $X_p$ introduced above, and we consider functions of type $f:X_p\to Y_q$, which by definition satisfy $f(p)=q$. We have the following definition, further clarifying the terminology in Proposition 7.1:

\index{affine deformation}
\index{trivial deformation}

\begin{definition}
Let $H\in M_N(\mathbb C)$ be a complex Hadamard matrix.
\begin{enumerate}
\item A deformation of $H$ is a smooth function $f:\mathbb T_1\to (X_N)_H$.

\item The deformation is called ``affine'' if $f_{ij}(q)=H_{ij}q^{A_{ij}}$, with $A\in M_N(\mathbb R)$.

\item We call ``trivial'' the deformations of type $f_{ij}(q)=H_{ij}q^{a_i+b_j}$, with $a,b\in\mathbb R^N$.
\end{enumerate}
\end{definition}

Here the adjective ``affine'', which is used in the same context as in Proposition 7.1, comes from the formula $f_{ij}(e^{it})=H_{ij}e^{iA_{ij}t}$, because the function $t\to A_{ij}t$ which produces the exponent is indeed affine. As for the adjective ``trivial'', this comes from the fact that the affine deformations of type $f(q)=(H_{ij}q^{a_i+b_j})_{ij}$ are obtained from $H$ by multiplying the rows and columns by certain numbers in $\mathbb T$, so are automatically Hadamard.

\bigskip

The basic example of an affine deformation comes from the Di\c t\u a deformations $H\otimes_QK$, by taking all parameters $q_{ij}\in\mathbb T$ to be powers of $q\in\mathbb T$. As an example, here are the exponent matrices coming from the left and right Di\c t\u a deformations of $F_2\otimes F_2$:
$$A_l=
\begin{pmatrix}
a&a&b&b\\
c&c&d&d\\ 
a&a&b&b\\
c&c&d&d
\end{pmatrix}\qquad,\qquad
A_r=
\begin{pmatrix}
a&b&a&b\\
a&b&a&b\\ 
c&d&c&d\\
c&d&c&d
\end{pmatrix}$$

There are of course many other examples, which are less trivial, as for instance the Haagerup matrix, that we met in chapters 5-6, which is as follows:
$$H_6^q=\begin{pmatrix}
1&1&1&1&1&1\\
1&-1&i&i&-i&-i\\ 
1&i&-1&-i&q&-q\\ 
1&i&-i&-1&-q&q\\
1&-i&\bar{q}&-\bar{q}&i&-1\\ 
1&-i&-\bar{q}&\bar{q}&-1&i
\end{pmatrix}$$

\index{Haagerup matrix}

Observe that this is indeed an affine deformation of $H_6=H_6^1$, in the sense of Definition 7.2 (2), the corresponding matrix of exponents being as follows:
$$A=\begin{pmatrix}
0&0&0&0&0&0\\
0&0&0&0&0&0\\ 
0&0&0&0&1&1\\ 
0&0&0&0&1&1\\
0&0&-1&-1&0&0\\ 
0&0&-1&-1&0&0
\end{pmatrix}$$

We will see that there are many other interesting examples of affine deformations, and that some general theory for such deformations can be developed. In order to investigate now the above types of deformations, we will use the corresponding tangent vectors. So, let us recall that the complex Hadamard matrix manifold $X_N$ is given by:
$$X_N=M_N(\mathbb T)\cap\sqrt{N}U_N$$

This observation leads to the following definition, where in the first part we denote by $T_pX$ the tangent space to a point in a smooth manifold, $p\in X$:

\index{tangent space}
\index{enveloping tangent space}
\index{tangent cone}
\index{enveloping tangent cone}
\index{trivial tangent cone}
\index{affine deformation}

\begin{definition}
Associated to a point $H\in X_N$ are the following objects:
\begin{enumerate}
\item The enveloping tangent space: $\widetilde{T}_HX_N=T_HM_N(\mathbb T)\cap T_H\sqrt{N}U_N$.

\item The tangent cone $T_HX_N$: the set of tangent vectors to the deformations of $H$.

\item The affine tangent cone $T_H^\circ X_N$: same as above, using affine deformations only.

\item The trivial tangent cone $T_H^\times X_N$: as above, using trivial deformations only.
\end{enumerate}
\end{definition}

Observe that $\widetilde{T}_HX_N,T_H^\times X_N$ are real linear spaces, and that $T_HX_N,T_H^\circ X_N$ are two-sided cones, in the sense that they satisfy the following condition:
$$\lambda\in\mathbb R,A\in T\implies\lambda A\in T$$

Observe also that we have inclusions of cones, as follows:
$$T_H^\times X_N\subset T_H^\circ X_N\subset T_HX_N\subset\widetilde{T}_HX_N$$

In more algebraic terms now, these various tangent cones are best described by the corresponding matrices, and we have here the following result:

\begin{theorem}
The cones $T_H^\times X_N\subset T_H^\circ X_N\subset T_HX_N\subset\widetilde{T}_HX_N$ are as follows:
\begin{enumerate}
\item $\widetilde{T}_HX_N$ can be identified with the linear space formed by the matrices $A\in M_N(\mathbb R)$ satisfying the following condition: 
$$\sum_kH_{ik}\bar{H}_{jk}(A_{ik}-A_{jk})=0$$.

\item $T_HX_N$ consists of those matrices $A\in M_N(\mathbb R)$ appearing as $A_{ij}=g_{ij}'(0)$, where $g:M_N(\mathbb R)_0\to M_N(\mathbb R)_0$ satisfies: 
$$\sum_kH_{ik}\bar{H}_{jk}e^{i(g_{ik}(t)-g_{jk}(t))}=0$$

\item $T^\circ_HX_N$ is formed by the matrices $A\in M_N(\mathbb R)$ satisfying the following condition, for any $i\neq j$ and any $q\in\mathbb T$:
$$\sum_kH_{ik}\bar{H}_{jk}q^{A_{ik}-A_{jk}}=0$$

\item $T^\times_HX_N$ is formed by the matrices $A\in M_N(\mathbb R)$ which are of the form $A_{ij}=a_i+b_j$, for certain vectors $a,b\in\mathbb R^N$.
\end{enumerate}
\end{theorem}

\begin{proof}
All these assertions can be deduced by using basic differential geometry:

\medskip

(1) This result is well-known, the idea being as follows. First, $M_N(\mathbb T)$ is defined by the algebraic relations $|H_{ij}|^2=1$, and with $H_{ij}=X_{ij}+iY_{ij}$ we have:
$$d|H_{ij}|^2
=d(X_{ij}^2+Y_{ij}^2)
=2(X_{ij}\dot{X}_{ij}+Y_{ij}\dot{Y}_{ij})$$

Consider now an arbitrary vector $\xi\in T_HM_N(\mathbb C)$, written as follows:
$$\xi=\sum_{ij}\alpha_{ij}\dot{X}_{ij}+\beta_{ij}\dot{Y}_{ij}$$

This vector belongs then to the tangent space $T_HM_N(\mathbb T)$ if and only if we have:
$$<\xi,d|H_{ij}|^2>=0$$

We therefore obtain the following formula, for the tangent cone:
$$T_HM_N(\mathbb T)=\left\{\sum_{ij}A_{ij}(Y_{ij}\dot{X}_{ij}-X_{ij}\dot{Y}_{ij})\Big|A_{ij}\in\mathbb R\right\}$$

We also know that the rescaled unitary group $\sqrt{N}U_N$ is defined by the following algebraic relations, where $H_1,\ldots,H_N$ are the rows of $H$:
$$<H_i,H_j>=N\delta_{ij}$$
 
The relations $<H_i,H_i>=N$ being automatic for the matrices $H\in M_N(\mathbb T)$, if for $i\neq j$ we let $L_{ij}=<H_i,H_j>$, then we have:
$$\widetilde{T}_HC_N=\left\{\xi\in T_HM_N(\mathbb T)\Big|<\xi,\dot{L}_{ij}>=0,\,\forall i\neq j\right\}$$

On the other hand, differentiating the formula of $L_{ij}$ gives:
$$\dot{L}_{ij}=\sum_k(X_{ik}+iY_{ik})(\dot{X}_{jk}-i\dot{Y}_{jk})+(X_{jk}-iY_{jk})(\dot{X}_{ik}+i\dot{Y}_{ik})$$

Now if we pick $\xi\in T_HM_N(\mathbb T)$, written as above in terms of $A\in M_N(\mathbb R)$, we obtain:
$$<\xi,\dot{L}_{ij}>=i\sum_k\bar{H}_{ik}H_{jk}(A_{ik}-A_{jk})$$

Thus we have reached to the description of $\widetilde{T}_HX_N$ in the statement. 

\medskip

(2) We pick an arbitrary deformation, written as $f_{ij}(e^{it})=H_{ij}e^{ig_{ij}(t)}$. Observe first that the Hadamard condition corresponds to the equations in the statement, namely:
$$\sum_kH_{ik}\bar{H}_{jk}e^{i(g_{ik}(t)-g_{jk}(t))}=0$$

Observe also that by differentiating this formula at $t=0$, we obtain:
$$\sum_kH_{ik}\bar{H}_{jk}(g_{ik}'(0)-g_{jk}'(0))=0$$

Thus the matrix $A_{ij}=g_{ij}'(0)$ belongs indeed to $\widetilde{T}_HX_N$, so we obtain in this way a certain map, as follows:
$$T_HX_N\to\widetilde{T}_HX_N$$

In order to check that this map is indeed the correct one, we have to verify that, for any $i,j$, the tangent vector to our deformation is given by:
$$\xi_{ij}=g_{ij}'(0)(Y_{ij}\dot{X}_{ij}-X_{ij}\dot{Y}_{ij})$$

But this latter verification is just a one-variable problem. So, by dropping all $i,j$ indices, which is the same as assuming $N=1$, we have to check that for any point $H\in\mathbb T$, written $H=X+iY$, the tangent vector to the deformation $f(e^{it})=He^{ig(t)}$ is:
$$\xi=g'(0)(Y\dot{X}-X\dot{Y})$$

But this is clear, because the unit tangent vector at $H\in\mathbb T$ is $\eta=-i(Y\dot{X}-X\dot{Y})$, and  its coefficient coming from the deformation is:
$$(e^{ig(t)})'_{|t=0}=-ig'(0)$$

(3) Observe first that by taking the derivative at $q=1$ of the condition (2) in Proposition 7.1, of just by using the condition (3) there with the function $\varphi(r)=r$, we get:
$$\sum_kH_{ik}\bar{H}_{jk}\varphi(A_{ik}-A_{jk})=0$$

Thus we have a map $T_H^\circ X_N\to\widetilde{T}_HX_N$, and the fact that is map is indeed the correct one comes for instance from the computation in (2), with $g_{ij}(t)=A_{ij}t$.

\medskip

(4) Observe first that the Hadamard matrix condition is satisfied, because:
$$\sum_kH_{ik}\bar{H}_{jk}q^{A_{ik}-A_{jk}}
=q^{a_i-a_j}\sum_kH_{ik}\bar{H}_{jk}
=\delta_{ij}$$

As for the fact that $T_H^\times X_N$ is indeed the space in the statement, this is clear.
\end{proof}

Let $Z_N\subset X_N$ be the real algebraic manifold formed by all the dephased $N\times N$ complex Hadamard matrices. Observe that we have a quotient map $X_N\to Z_N$, obtained by dephasing. With this notation, we have the following refinement of (4) above:

\begin{proposition}
We have a direct sum decomposition of cones
$$T_H^\circ X_N=T_H^\times X_N\oplus T_H^\circ Z_N$$
where at right we have the affine tangent cone to the dephased manifold $X_N\to Z_N$.
\end{proposition}

\begin{proof}
If we denote by $M_N^\circ(\mathbb R)$ the set of matrices having $0$ outside the first row and column, we have a direct sum decomposition, as follows:
$$\widetilde{T}_H^\circ X_N=M_N^\circ(\mathbb R)\oplus\widetilde{T}_H^\circ Z_N$$

Now by looking at the affine cones, and using Theorem 7.4, this gives the result.
\end{proof}

Summarizing, we have so far a number of theoretical results about the tangent cones $T_HX_N$ that we are interested in, and their versions coming from the trivial and affine deformations, and from the intersection formula $X_N=M_N(\mathbb T)\cap\sqrt{N}U_N$ as well. In practice now, passed a few special cases where all these cones collapse to the trivial cone $T_N^\times X_N$, which by Proposition 7.5 means that the image of $H\in X_N$ must be isolated in the dephased manifold $X_N\to Z_N$, things are quite difficult to compute. However, as a concrete numerical invariant arising from all this, which can be effectively computed in many cases of interest, we have, following Tadej-\.Zyczkowski \cite{tz2}:

\index{defect}
\index{undephased defect}

\begin{definition}
The real dimension $d(H)$ of the enveloping tangent space 
$$\widetilde{T}_HX_N=T_HM_N(\mathbb T)\cap T_H\sqrt{N}U_N$$
is called undephased defect of a complex Hadamard matrix $H\in X_N$. 
\end{definition}

In view of Proposition 7.5, it is sometimes convenient to replace $d(H)$ by the following related quantity, also introduced in \cite{tz2}, and called dephased defect of $H$: 
$$d'(H)=d(H)-2N+1$$
 
In what follows we will rather use the quantity $d(H)$ defined before, which behaves better with respect to a number of operations, and simply call it ``defect'' of $H$. We already know, from Theorem 7.4, what is the precise geometric meaning of the defect, and how to compute it. Let us record again these results, that we will use many times in what follows, in a slightly different form, closer to the spirit of \cite{tz2}:

\index{defect equations}

\begin{theorem}
The defect $d(H)$ is the real dimension of the linear space
$$\widetilde{T}_HX_N=\left\{A\in M_N(\mathbb R)\Big|\sum_kH_{ik}\bar{H}_{jk}(A_{ik}-A_{jk})=0,\forall i,j\right\}$$
and the elements of this space are those making $H^q_{ij}=H_{ij}q^{A_{ij}}$ Hadamard at order $1$.
\end{theorem}

\begin{proof}
Here the first assertion is something that we already know, from Theorem 7.4 (1), and the second assertion follows either from Theorem 7.4 and its proof, or directly from the definition of the enveloping tangent space $\widetilde{T}_HX_N$, as used in Definition 7.6.
\end{proof}

Still following \cite{tz2}, here are a few basic properties of the defect:

\index{isolated matrix}

\begin{proposition}
Let $H\in X_N$ be a complex Hadamard matrix.
\begin{enumerate}
\item If $H\simeq\widetilde{H}$ then $d(H)=d(\widetilde{H})$.

\item We have $2N-1\leq d(H)\leq N^2$.

\item If $d(H)=2N-1$, the image of $H$ in the dephased manifold $X_N\to Z_N$ is isolated.
\end{enumerate}
\end{proposition}

\begin{proof}
All these results are elementary, the proof being as follows:

\medskip

(1) If we let $K_{ij}=a_ib_jH_{ij}$ with $|a_i|=|b_j|=1$ be a trivial deformation of our matrix $H$, the equations for the enveloping tangent space for $K$ are:
$$\sum_ka_ib_kH_{ik}\bar{a}_j\bar{b}_k\bar{H}_{jk}(A_{ik}-A_{jk})=0$$

By simplifying we obtain the equations for $H$, so $d(H)$ is invariant under trivial deformations. Since $d(H)$ is invariant as well by permuting rows or columns, we are done.

\medskip

(2) Consider the inclusions $T_H^\times X_N\subset T_HX_N\subset\widetilde{T}_HX_N$. Since $\dim(T_H^\times X_N)=2N-1$, the inequality at left holds indeed. As for the inequality at right, this is clear.

\medskip

(3) If $d(H)=2N-1$ then $T_HX_N=T_H^\times X_N$, so any deformation of $H$ is trivial. Thus the image of $H$ in the quotient manifold $X_N\to Z_N$ is indeed isolated, as stated.
\end{proof}

\section*{7b. Defect computations}

As an illustration for the above notions, let us discuss now the computation of the defect for the most basic examples of complex Hadamard matrices that we know, namely the real ones, and the Fourier ones. In order to deal with the real case, it is convenient to modify the general formula from Theorem 7.7, via a change of variables, as follows:

\begin{proposition}
We have a linear space isomorphism as follows,
$$\widetilde{T}_HX_N\simeq\left\{E\in M_N(\mathbb C)\Big|E=E^*,(EH)_{ij}\bar{H}_{ij}\in\mathbb R,\forall i,j\right\}$$
the correspondences $A\to E$ and $E\to A$ being given by the formulae
$$E_{ij}=\sum_kH_{ik}\bar{H}_{jk}A_{ik}\quad,\quad 
A_{ij}=(EH)_{ij}\bar{H}_{ij}$$
with $A\in\widetilde{T}_HX_N$ being the usual components, from Theorem 7.7.
\end{proposition}

\begin{proof}
Given a matrix $A\in M_N(\mathbb C)$, if we set $R_{ij}=A_{ij}H_{ij}$ and $E=RH^*$, the correspondence $A\to R\to E$ is then bijective onto $M_N(\mathbb C)$, and we have:
$$E_{ij}=\sum_kH_{ik}\bar{H}_{jk}A_{ik}$$

In terms of these new variables, the equations in Theorem 7.7 become:
$$E_{ij}=\bar{E}_{ji}$$

Thus, when taking into account these conditions, we are simply left with the conditions $A_{ij}\in\mathbb R$. But these correspond to the conditions $(EH)_{ij}\bar{H}_{ij}\in\mathbb R$, as claimed.
\end{proof}

With the above result in hand, we can now compute the defect of the real Hadamard matrices. The result here, from Sz\"oll\H{o}si \cite{sz1}, is as follows:

\index{real Hadamard matrix}

\begin{theorem}
For any real Hadamard matrix $H\in M_N(\pm1)$ we have
$$\widetilde{T}_HX_N\simeq M_N(\mathbb R)^{symm}$$
and so the corresponding defect is $d(H)=N(N+1)/2$.
\end{theorem}

\begin{proof}
We use Proposition 7.9. Since $H$ is now real the condition $(EH)_{ij}\bar{H}_{ij}\in\mathbb R$ there simply tells us that $E$ must be real, and this gives the result.
\end{proof}

As another computation now, let us discuss the case $N=4$. Here we know from chapter 5 that the only complex Hadamard matrices are, up to equivalence, the Di\c t\u a deformations of $F_4$. To be more precise, we have the following result:

\index{Di\c t\u a deformation}

\begin{proposition}
The complex Hadamard matrices at $N=4$ are, up to equivalence, the following matrices, appearing as Di\c t\u a deformations of $F_4$:
$$F_{2,2}^q
=\begin{pmatrix}
1&1\\
1&-1
\end{pmatrix}
\otimes_{\begin{pmatrix}
1&1\\
1&q
\end{pmatrix}}
\begin{pmatrix}
1&1\\
1&-1
\end{pmatrix}
=\begin{pmatrix}
1&1&1&1\\
1&-1&q&-q\\
1&1&-1&-1\\ 
1&-1&-q&q
\end{pmatrix}$$
At $q\in\{1,i,-1,-i\}$ we obtain tensor products of Fourier matrices, as follows:
\begin{enumerate}
\item At $q=1$ we have $F_{2,2}^q=F_2\otimes F_2$.

\item At $q=-1$ we have $F_{2,2}^q\simeq F_2\otimes F_2$.

\item At $q=\pm i$ we have $F_{2,2}^q\simeq F_4$.
\end{enumerate}
\end{proposition}

\begin{proof}
The first assertion is something that we already know, from chapter 5. Regarding now the $q=1,i,-1,-i$ specializations, the situation here is as follows:

\medskip

(1) This is clear from definitions.

\medskip

(2) This follows from (1), by permuting the third and the fourth columns:
$$F_{2,2}^{-1}
=\begin{pmatrix}
1&1&1&1\\
1&-1&-1&1\\
1&1&-1&-1\\ 
1&-1&1&-1
\end{pmatrix}
\sim\begin{pmatrix}
1&1&1&1\\
1&-1&1&-1\\
1&1&-1&-1\\ 
1&-1&-1&1
\end{pmatrix}
=F_{2,2}^1$$

(3) This follows from the following computation:
$$F_{2,2}^{\pm i}
=\begin{pmatrix}
1&1&1&1\\
1&-1&\pm i&\mp i\\
1&1&-1&-1\\ 
1&-1&\mp i&\pm i
\end{pmatrix}
\sim\begin{pmatrix}
1&1&1&1\\
1&i&-1&-i\\
1&-1&1&-1\\ 
1&-i&-1&i
\end{pmatrix}
=F_4$$

Here we have interchanged the second column with the third one in the case $q=i$, and we have used a cyclic permutation of the last 3 columns in the case $q=-i$. 
\end{proof}

Let us compute now the defect of the above matrices. We will work out everything in detail, as an illustration for how the equations in Theorem 7.7 work. The result is:

\begin{theorem}
The defect of the $4\times4$ complex Hadamard matrices is given by
$$d(F_{2,2}^q)=
\begin{cases}
10&(q=\pm1)\\
8&(q\neq\pm1)
\end{cases}$$
with $F_{2,2}^q$, depending on $q\in\mathbb T$, being the matrix in Proposition 7.11.
\end{theorem}

\begin{proof}
Our starting point are the equations in Theorem 7.7, namely:
$$\sum_hH_{ik}\bar{H}_{jk}(A_{ik}-A_{jk})=0$$

Since the $i>j$ equations are equivalent to the $i<j$ ones, and the $i=j$ equations are trivial, we just have to write down the equations corresponding to indices $i<j$. And, with $ij=01,02,03,12,13,23$, these equations are:
\begin{eqnarray*}
(A_{00}-A_{10})-(A_{01}-A_{11})+\bar{q}(A_{02}-A_{12})-\bar{q}(A_{03}-A_{13})&=&0\\
(A_{00}-A_{20})+(A_{01}-A_{21})-(A_{02}-A_{22})-(A_{03}-A_{23})&=&0\\
(A_{00}-A_{30})-(A_{01}-A_{31})-\bar{q}(A_{02}-A_{32})+\bar{q}(A_{03}-A_{33})&=&0\\
(A_{10}-A_{20})-(A_{11}-A_{21})-q(A_{12}-A_{22})+q(A_{13}-A_{23})&=&0\\
(A_{10}-A_{30})+(A_{11}-A_{31})-(A_{12}-A_{32})-(A_{13}-A_{33})&=&0\\
(A_{20}-A_{30})-(A_{21}-A_{31})+\bar{q}(A_{22}-A_{32})-\bar{q}(A_{23}-A_{33})&=&0
\end{eqnarray*}

Assume first $q\neq\pm 1$. Then $q$ is not real, and appears in 4 of the above equations. But these 4 equations can be written in the following way:
\begin{eqnarray*}
(A_{00}-A_{01})-(A_{10}-A_{11})+\bar{q}((A_{02}-A_{03})-(A_{12}-A_{13}))&=&0\\
(A_{00}-A_{01})-(A_{30}-A_{31})-\bar{q}((A_{02}-A_{03})-(A_{32}-A_{33}))&=&0\\
(A_{10}-A_{11})-(A_{20}-A_{21})-q((A_{12}-A_{13})-(A_{22}-A_{23}))&=&0\\
(A_{20}-A_{21})-(A_{30}-A_{31})+\bar{q}((A_{22}-A_{23})-(A_{32}-A_{33}))&=&0
\end{eqnarray*}

Now since the unknowns are real, and $q$ is not, we conclude that the terms between braces in the left part must be all equal, and that the same must happen at right:
\begin{eqnarray*}
A_{00}-A_{01}&=&A_{10}-A_{11}=A_{20}-A_{21}=A_{30}-A_{31}\\
A_{02}-A_{03}&=&A_{12}-A_{13}=A_{22}-A_{23}=A_{32}-A_{33}
\end{eqnarray*}

Thus, the equations involving $q$ tell us that $A$ must be of the following form:
$$A=\begin{pmatrix}
a&a+x&e+y&e\\
b&b+x&f+y&f\\
c&c+x&g+y&g\\
d&d+x&h+y&h\end{pmatrix}$$

Let us plug now these values in the remaining 2 equations. We obtain:
\begin{eqnarray*}
a-c+a+x-c-x-e-y+g+y-e+g&=&0\\
b-d+b+x-d-x-f-y+h+y-f+h&=&0
\end{eqnarray*}

Thus we must have $a+g=c+e$ and $b+h=d+f$, which are independent conditions. We conclude that the dimension of the space of solutions is $10-2=8$, as claimed.

Assume now $q=\pm 1$. For simplicity we set $q=1$, and we compute the dephased defect. The dephased equations, obtained by setting $A_{i0}=A_{0j}=0$ in our system, are:
\begin{eqnarray*}
A_{11}-A_{12}+A_{13}&=&0\\
-A_{21}+A_{22}+A_{23}&=&0\\
A_{31}+A_{32}-A_{33}&=&0\\
-A_{11}+A_{21}-A_{12}+A_{22}+A_{13}-A_{23}&=&0\\
A_{11}-A_{31}-A_{12}+A_{32}-A_{13}+A_{33}&=&0\\
-A_{21}+A_{31}+A_{22}-A_{32}-A_{23}+A_{33}&=&0
\end{eqnarray*}

The first three equations tell us that our matrix must be of the following form:
$$A=\begin{pmatrix}a&a+b&b\\ c+d&c&d\\ e&f&e+f\end{pmatrix}$$

Now by plugging these values in the last three equations, these become:
\begin{eqnarray*}
-a+c+d-a-b+c+b-d&=&0\\
a-e-a-b+f-b+e+f&=&0\\
-c-d+e+c-f-d+e+f&=&0
\end{eqnarray*}

Thus we must have $a=c$, $b=f$, $d=e$, and since these conditions are independent, the dephased defect is 3, and so the undephased defect is $3+7=10$, as claimed.
\end{proof}

In general, the defect computation for the Di\c t\u a deformations, of even for the usual tensor products, is a difficult question. We will be back to this in chapter 8 below.

\section*{7c. Fourier matrices}

Let us discuss now a fundamental question, namely the computation of the defect of the Fourier matrix $F_G$. The main idea here goes back to a 1989 preprint of Karabegov \cite{kar}, with some supplementary contributions from Nicoara \cite{nic}, in 2006, and then the main formula, in the cyclic group case, was obtained by Tadej-\.Zyczkowski in \cite{tz2}, and the corresponding deformations of $F_G$ were studied by Nicoara-White in \cite{nwh}. As a first result on this subject, we have, following Tadej-\.Zyczkowski \cite{tz2}:

\begin{theorem}
For a Fourier matrix $F=F_G$, the matrices $A\in\widetilde{T}_FX_N$
with $N=|G|$, are those of the form $A=PF^*$, with $P\in M_N(\mathbb C)$ satisfying
$$P_{ij}=P_{i+j,j}=\bar{P}_{i,-j}$$
where the indices $i,j$ are by definition taken in the group $G$.
\end{theorem}

\begin{proof}
We use the system of equations in Theorem 7.7, namely:
$$\sum_kF_{ik}\bar{F}_{jk}(A_{ik}-A_{jk})=0$$

By decomposing our finite abelian group as $G=\mathbb Z_{N_1}\times\ldots\times\mathbb Z_{N_r}$ we can assume:
$$F=F_{N_1}\otimes\ldots\otimes F_{N_r}$$

Thus with $w_k=e^{2\pi i/k}$ we have the following formula:
$$F_{i_1\ldots i_r,j_1\ldots j_r}=(w_{N_1})^{i_1j_1}\ldots (w_{N_r})^{i_rj_r}$$

With $N=N_1\ldots N_r$ and $w=e^{2\pi i/N}$, we obtain the following formula:
$$F_{i_1\ldots i_r,j_1\ldots j_r}=w^{\left(\frac{i_1j_1}{N_1}+\ldots+\frac{i_rj_r}{N_r}\right)N}$$

Thus the matrix of our system of equations is given by:
$$F_{i_1\ldots i_r,k_1\ldots k_r}\bar{F}_{j_1\ldots j_r,k_1\ldots k_r}=w^{\left(\frac{(i_1-j_1)k_1}{N_1}+\ldots+\frac{(i_r-j_r)k_r}{N_r}\right)N}$$

Now by plugging in a multi-indexed matrix $A$, our system becomes:
$$\sum_{k_1\ldots k_r}w^{\left(\frac{(i_1-j_1)k_1}{N_1}+\ldots+\frac{(i_r-j_r)k_r}{N_r}\right)N}(A_{i_1\ldots i_r,k_1\ldots k_r}-A_{j_1\ldots j_r,k_1\ldots k_r})=0$$

Now observe that in the above formula we have in fact two matrix multiplications, so our system can be simply written as:
$$(AF)_{i_1\ldots i_r,i_1-j_1\ldots i_r-j_r}-(AF)_{j_1\ldots j_r,i_1-j_1\ldots i_r-j_r}=0$$

Now recall that our indices have a ``cyclic'' meaning, so they belong in fact to the group $G$. So, with $P=AF$, and by using multi-indices, our system is simply:
$$P_{i,i-j}=P_{j,i-j}$$

With $i=I+J,j=I$ we obtain the condition $P_{I+J,J}=P_{IJ}$ in the statement. In addition, $A=PF^*$ must be a real matrix. But, if we set $\tilde{P}_{ij}=\bar{P}_{i,-j}$, we have:
\begin{eqnarray*}
\overline{(PF^*)}_{i_1\ldots i_r,j_1\ldots j_r}
&=&\sum_{k_1\ldots k_r}\bar{P}_{i_1\ldots i_r,k_1\ldots k_r}F_{j_1\ldots j_r,k_1\ldots k_r}\\
&=&\sum_{k_1\ldots k_r}\tilde{P}_{i_1\ldots i_r,-k_1\ldots -k_r}(F^*)_{-k_1\ldots -k_r,j_1\ldots j_r}\\
&=&(\tilde{P}F^*)_{i_1\ldots i_r,j_1\ldots j_r}
\end{eqnarray*}

Thus we have $\overline{PF^*}=\tilde{P}F^*$, so the fact that the matrix $PF^*$ is real, which means by definition that we have $\overline{PF^*}=PF^*$, can be reformulated as $\tilde{P}F^*=PF^*$, and hence as $\tilde{P}=P$. So, we obtain the conditions $P_{ij}=\bar{P}_{i,-j}$ in the statement.
\end{proof}

We can now compute the defect, and we are led to the following formula:

\index{defect}
\index{number of 1 entries}

\begin{theorem}
The defect of a Fourier matrix $F_G$ is given by
$$d(F_G)=\sum_{g\in G}\frac{|G|}{ord(g)}$$
and equals as well the number of $1$ entries of the matrix $F_G$.
\end{theorem}

\begin{proof}
According to the formula $A=PF^*$ from Theorem 7.13, the defect $d(F_G)$ is the dimension of the real vector space formed by the matrices $P\in M_N(\mathbb C)$ satisfying:
$$P_{ij}=P_{i+j,j}=\bar{P}_{i,-j}$$

Here, and in what follows, the various indices $i,j,\ldots$ will be taken in $G$. Now the point is that, in terms of the columns of our matrix $P$, the above conditions are:

\medskip

(1) The entries of the $j$-th column of $P$, say $C$, must satisfy $C_i=C_{i+j}$.

\medskip

(2) The $(-j)$-th column of $P$ must be conjugate to the $j$-th column of $P$.

\medskip

Thus, in order to count the above matrices $P$, we can basically fill the columns one by one, by taking into account the above conditions. In order to do so, consider the subgroup $G_2=\{j\in G|2j=0\}$, and then write $G$ as a disjoint union, as follows:
$$G=G_2\sqcup X\sqcup(-X)$$ 

With this notation, the algorithm is as follows. First, for any $j\in G_2$ we must fill the $j$-th column of $P$ with real numbers, according to the periodicity rule:
$$C_i=C_{i+j}$$

Then, for any $j\in X$ we must fill the $j$-th column of $P$ with complex numbers, according to the same periodicity rule $C_i=C_{i+j}$. And finally, once this is done, for any $j\in X$ we just have to set the $(-j)$-th column of $P$ to be the conjugate of the $j$-th column.

\medskip

So, let us compute the number of choices for filling these columns. Our claim is that, when uniformly distributing the choices for the $j$-th and $(-j)$-th columns, for $j\notin G_2$, there are exactly $[G:<j>]$ choices for the $j$-th column, for any $j$. Indeed:

\medskip

(1) For the $j$-th column with $j\in G_2$ we must simply pick $N$ real numbers subject to the condition $C_i=C_{i+j}$ for any $i$, so we have indeed $[G:<j>]$ such choices.

\medskip

(2) For filling the $j$-th and $(-j)$-th column, with $j\notin G_2$, we must pick $N$ complex numbers subject to the condition $C_i=C_{i+j}$ for any $i$. Now since there are $[G:<j>]$ choices for these numbers, so a total of $2[G:<j>]$ choices for their real and imaginary parts, on average over $j,-j$ we have $[G:<j>]$ choices, and we are done again.

\medskip

Summarizing, the dimension of the vector space formed by the matrices $P$, which is equal to the number of choices for the real and imaginary parts of the entries of $P$, is:
$$d(F_G)=\sum_{j\in G}[G:<j>]$$

But this is exactly the number in the statement. Regarding now the second assertion, according to the definition of $F_G$, the number of $1$ entries of $F_G$ is given by:
\begin{eqnarray*}
\#(1\in F_G)
&=&\#\left\{(g,\chi)\in G\times\widehat{G}\Big|\chi(g)=1\right\}\\
&=&\sum_{g\in G}\#\left\{\chi\in\widehat{G}\Big|\chi(g)=1\right\}\\
&=&\sum_{g\in G}\frac{|G|}{ord(g)}
\end{eqnarray*}

Thus, the second assertion follows from the first one.
\end{proof}

Let us finish now the work, and explicitely compute the defect of $F_G$. It is convenient to consider the following quantity, which behaves better:
$$\delta(G)=\sum_{g\in G}\frac{1}{ord(g)}$$

As a first example, consider a cyclic group $G=\mathbb Z_N$, with $N=p^a$ power of a prime. The count here is very simple, over sets of elements having a given order:
\begin{eqnarray*}
\delta(\mathbb Z_{p^a})
&=&1+(p-1)p^{-1}+(p^2-p)p^{-2}+\ldots+(p^a-p^{a-1})p^{-1}\\
&=&1+a-\frac{a}{p}
\end{eqnarray*}

In order to extend this kind of count to the general abelian case, we use two ingredients. First is the following result, which splits the computation over isotypic components:

\begin{proposition}
For any finite groups $G,H$ we have:
$$\delta(G\times H)\geq\delta(G)\delta(H)$$
In addition, if $(|G|,|H|)=1$, we have equality. 
\end{proposition}

\begin{proof}
Indeed, we have the following estimate, coming from definitions:
\begin{eqnarray*}
\delta(G\times H)
&=&\sum_{gh}\frac{1}{ord(g,h)}\\
&=&\sum_{gh}\frac{1}{[ord(g),ord(h)]}\\
&\geq&\sum_{gh}\frac{1}{ord(g)\cdot ord(h)}\\
&=&\delta(G)\delta(H)
\end{eqnarray*}

Regarding the last assertion, in the case $(|G|,|H|)=1$, the least common multiple appearing on the right becomes a product:
$$[ord(g),ord(h)]=ord(g)\cdot ord(h)$$

Thus, we have equality in this case, as desired. 
\end{proof}

We deduce from this that we have the following result:

\begin{proposition}
For a finite abelian group $G$ we have
$$\delta(G)=\prod_p\delta(G_p)$$
where $G_p$ with $G=\times_pG_p$ are the isotypic components of $G$.
\end{proposition}

\begin{proof}
This is clear from Proposition 7.15, the order of $G_p$ being a power of $p$.
\end{proof}

As an illustration for the above results, we can recover in this way the following key defect computation, from Tadej-\.Zyczkowski \cite{tz2}:

\index{Tadej-\.Zyczkowski formula}

\begin{theorem}
The defect of a usual Fourier matrix $F_N$ is given by
$$d(F_N)=N\prod_{i=1}^s\left(1+a_i-\frac{a_i}{p_i}\right)$$
where $N=p_1^{a_1}\ldots p_s^{a_s}$ is the decomposition of $N$ into prime factors.
\end{theorem}

\begin{proof}
The underlying group here is the cyclic group $G=\mathbb Z_N$, whose isotypic components are the following cyclic groups:
$$G_{p_i}=\mathbb Z_{p_i^{a_i}}$$

By applying now Proposition 7.16, and by using the computation for cyclic $p$-groups performed before Proposition 7.15, we obtain:
$$d(F_N)=N\prod_{i=1}^s\left(1+p_i^{-1}(p_i-1)a_i\right)$$

But this is exactly the formula in the statement.
\end{proof}

Now back to the general case, where we have an arbitrary Fourier matrix $F_G$, we will need, as a second ingredient for our computation, the following result:

\begin{proposition}
For the $p$-groups, the quantities
$$c_k=\#\left\{g\in G\Big|ord(g)\leq p^k\right\}$$
are multiplicative, in the sense that $c_k(G\times H)=c_k(G)c_k(H)$.
\end{proposition}

\begin{proof}
Indeed, for a product of $p$-groups we have:
\begin{eqnarray*}
c_k(G\times H)
&=&\#\left\{(g,h)\Big|ord(g,h)\leq p^k\right\}\\
&=&\#\left\{(g,h)\Big|ord(g)\leq p^k,ord(h)\leq p^k\right\}\\
&=&\#\left\{g\Big|ord(g)\leq p^k\right\}\#\left\{h\Big|ord(h)\leq p^k\right\}
\end{eqnarray*}

We recognize at right $c_k(G)c_k(H)$, and we are done.
\end{proof}

Let us compute now $\delta$ in the general isotypic case. We have here:

\begin{proposition}
For $G=\mathbb Z_{p^{a_1}}\times\ldots\times\mathbb Z_{p^{a_r}}$ with $a_1\leq a_2\leq\ldots\leq a_r$ we have
$$\delta(G)=1+\sum_{k=1}^rp^{(r-k)a_{k-1}+(a_1+\ldots+a_{k-1})-1}(p^{r-k+1}-1)[a_k-a_{k-1}]_{p^{r-k}}$$
with the convention $a_0=0$, and with the notation $[a]_q=1+q+q^2+\ldots+q^{a-1}$.
\end{proposition}

\begin{proof}
First, in terms of the numbers $c_k$, we have the following formula:
$$\delta(G)=1+\sum_{k\geq 1}\frac{c_k-c_{k-1}}{p^k}$$

In the case of a cyclic group $G=\mathbb Z_{p^a}$ we have $c_k=p^{\min(k,a)}$. Thus, in the general isotypic case $G=\mathbb Z_{p^{a_1}}\times\ldots\times\mathbb Z_{p^{a_r}}$ we have the following formula:
\begin{eqnarray*}
c_k
&=&p^{\min(k,a_1)}\ldots p^{\min(k,a_r)}\\
&=&p^{\min(k,a_1)+\ldots+\min(k,a_r)}
\end{eqnarray*}

Now observe that the exponent on the right is a piecewise linear function of $k$. More precisely, by assuming $a_1\leq a_2\leq\ldots\leq a_r$ as in the statement, the exponent is linear on each of the intervals $[0,a_1],[a_1,a_2],\ldots,[a_{r-1},a_r]$. So, the quantity $\delta(G)$ to be computed will be 1 plus the sum of $2r$ geometric progressions, 2 for each interval.

\medskip

In practice now, the numbers $c_k$ are as follows:
$$c_0=1,c_1=p^r,c_2=p^{2r},\ldots,c_{a_1}=p^{ra_1},$$
$$c_{a_1+1}=p^{a_1+(r-1)(a_1+1)},c_{a_1+2}=p^{a_1+(r-1)(a_1+2)},\ldots,c_{a_2}=p^{a_1+(r-1)a_2},$$
$$c_{a_2+1}=p^{a_1+a_2+(r-2)(a_2+1)},c_{a_2+2}=p^{a_1+a_2+(r-2)(a_2+2)},\ldots,c_{a_3}=p^{a_1+a_2+(r-2)a_3},$$
$$\vdots$$
$$c_{a_{r-1}+1}=p^{a_1+\ldots+a_{r-1}+(a_{r-1}+1)},c_{a_{r-1}+2}=p^{a_1+\ldots+a_{r-1}+(a_{r-1}+2)},\ldots,c_{a_r}=p^{a_1+\ldots+a_r}$$

Now by separating the positive and negative terms in the above formula of $\delta(G)$, we have indeed $2r$ geometric progressions to be summed, as follows:
\begin{eqnarray*}
\delta(G)
&=&1+(p^{r-1}+p^{2r-2}+p^{3r-3}+\ldots+p^{a_1r-a_1})\\
&&-(p^{-1}+p^{r-2}+p^{2r-3}+\ldots+p^{(a_1-1)r-a_1})\\
&&+(p^{(r-1)(a_1+1)-1}+p^{(r-1)(a_1+2)-2}+\ldots+p^{a_1+(r-2)a_2})\\
&&-(p^{a_1r-a_1-1}+p^{(r-1)(a_1+1)-2}+\ldots+p^{a_1+(r-1)(a_2-1)-a_2})\\
&&\vdots\\
&&+(p^{a_1+\ldots+a_{r-1}}+p^{a_1+\ldots+a_{r-1}}+\ldots+p^{a_1+\ldots+a_{r-1}})\\
&&-(p^{a_1+\ldots+a_{r-1}-1}+p^{a_1+\ldots+a_{r-1}-1}+\ldots+p^{a_1+\ldots+a_{r-1}-1})
\end{eqnarray*}

Now by performing all the sums, we obtain the following formula:
\begin{eqnarray*}
\delta(G)
&=&1+p^{-1}(p^r-1)\frac{p^{(r-1)a_1}-1}{p^{r-1}-1}\\
&&+p^{(r-2)a_1+(a_1-1)}(p^{r-1}-1)\frac{p^{(r-2)(a_2-a_1)}-1}{p^{r-2}-1}\\
&&+p^{(r-3)a_2+(a_1+a_2-1)}(p^{r-2}-1)\frac{p^{(r-3)(a_3-a_2)}-1}{p^{r-3}-1}\\
&&\vdots\\
&&+p^{a_1+\ldots+a_{r-1}-1}(p-1)(a_r-a_{r-1})
\end{eqnarray*}

By looking now at the general term, we get the formula in the statement.
\end{proof}

Let us go back now to the general defect formula in Theorem 7.14. By putting it together with the various results above, we obtain:

\begin{theorem}
For a finite abelian group $G$, decomposed as $G=\times_pG_p$, we have
$$d(F_G)=|G|\prod_p\left( 1+\sum_{k=1}^rp^{(r-k)a_{k-1}+(a_1+\ldots+a_{k-1})-1}(p^{r-k+1}-1)[a_k-a_{k-1}]_{p^{r-k}}\right)$$
where $a_0=0$ and $a_1\leq a_2\leq\ldots\leq a_r$ are such that $G_p=\mathbb Z_{p^{a_1}}\times\ldots\times\mathbb Z_{p^{a_r}}$.
\end{theorem}

\begin{proof}
Indeed, we know from Theorem 7.14 that we have:
$$d(F_G)=|G|\delta(G)$$

The result follows then from Proposition 7.16 and Proposition 7.19.
\end{proof}

As a first illustration, we can recover in this way the formula in Theorem 7.17. Indeed, assuming that $N=p_1^{a_1}\ldots p_s^{a_s}$ is the decomposition of $N$ into prime factors, we have:
\begin{eqnarray*}
d(F_N)
&=&N\prod_{i=1}^s\left(1+p_i^{-1}(p_i-1)a_i\right)\\
&=&N\prod_{i=1}^s\left(1+a_i-\frac{a_i}{p_i}\right)
\end{eqnarray*}

As a second illustration, for the group $G=\mathbb Z_{p^{a_1}}\times\mathbb Z_{p^{a_2}}$ with $a_1\leq a_2$ we obtain:
\begin{eqnarray*}
d(F_G)
&=&p^{a_1+a_2}(1+p^{-1}(p^2-1)[a_1]_p+p^{a_1-1}(p-1)(a_2-a_1))\\
&=&p^{a_1+a_2-1}(p+(p^2-1)\frac{p^{a_1}-1}{p-1}+p^{a_1}(p-1)(a_2-a_1))\\
&=&p^{a_1+a_2-1}(p+(p+1)(p^{a_1}-1)+p^{a_1}(p-1)(a_2-a_1))
\end{eqnarray*}

Finally, let us mention that for general non-abelian groups, there does not seem to be any reasonable algebraic formula for the quantity $\delta(G)$. As an example, consider the dihedral group $D_N$, consisting of $N$ symmetries and $N$ rotations. We have:
$$\delta(D_N)=\frac{N}{2}+\delta(\mathbb Z_N)$$

Now remember the formula for $\mathbb Z_N$ established above, namely:
$$\delta(\mathbb Z_N)=\prod_i(1+p_i^{-1}(p_i-1)a_i)$$

It is quite clear that the $N/2$ factor can not be incorporated in any nice way, and so, as indicated above, the quantity $\delta(G)$ remains something quite complicated.

\section*{7d. Explicit deformation}

Let us discuss now, following the paper of Nicoara and White \cite{nwh}, the key fact that for the Fourier matrices the defect is ``attained'', in the sense that the deformations at order 0 are true deformations, at order $\infty$. This is something quite surprising, and non-trivial. Let us begin with some generalities. We first recall that we have:

\index{Lie algebra}

\begin{proposition}
The unitary matrices $U\in U_N$ around $1$ are of the form
$$U=e^A$$
with $A$ being an antihermitian matrix, $A=-A^*$, around $0$.
\end{proposition}

\begin{proof}
This is something well-known. Indeed, assuming that a matrix $A$ is antihermitian, $A=-A^*$, the matrix $U=e^A$ follows to be unitary:
\begin{eqnarray*}
UU^*
&=&e^A(e^A)^*\\
&=&e^Ae^{A^*}\\
&=&e^Ae^{-A}\\
&=&1
\end{eqnarray*}

As for the converse, this follows either by using a dimension argument, which shows that the space of antihermitian matrices is the correct one, or by diagonalizing $U$. 
\end{proof}

Now back to the Hadamard matrices, we will need to rewrite a part of the basic theory of the defect, using deformations of type $t\to U_tH$. First, we have:

\begin{theorem}
Assume that $H\in M_N(\mathbb C)$ is Hadamard, let $A\in M_N(\mathbb C)$ be antihermitian, and consider the matrix $UH$, where $U=e^{tA}$, with $t\in\mathbb R$.
\begin{enumerate}
\item $UH$ is Hadamard when, for any $p,q$:
$$|\sum_{rs}H_{rq}\bar{H}_{sq}(e^{tA})_{pr}(e^{-tA})_{sp}|=1$$

\item $UH$ is Hadamard at order $0$ when, for any $p,q$:
$$|(AH)_{pq}|=1$$
\end{enumerate}
\end{theorem}

\begin{proof}
We already know that $UH$ is unitary, so we must find the conditions which guarantee that we have $UH\in M_N(\mathbb T)$, in general, and then at order 0.

\medskip

(1) We have the following computation, valid for any unitary $U$:
\begin{eqnarray*}
|(UH)_{pq}|^2
&=&(UH)_{pq}\overline{(UH)_{pq}}\\
&=&(UH)_{pq}(H^*U^*)_{qp}\\
&=&\sum_{rs}U_{pr}H_{rq}(H^*)_{qs}(U^*)_{sp}\\
&=&\sum_{rs}H_{rq}\bar{H}_{sq}U_{pr}\bar{U}_{ps}
\end{eqnarray*}

Now with $U=e^{tA}$ as in the statement, we obtain:
$$|(e^{tA}H)_{pq}|^2=\sum_{rs}H_{rq}\bar{H}_{sq}(e^{tA})_{pr}(e^{-tA})_{sp}$$

Thus, we are led to the conclusion in the statement.

\medskip

(2) The derivative of the function computed above, taken at $0$, is as follows:
\begin{eqnarray*}
\frac{\partial |(e^{tA}H)_{pq}|^2}{\partial t}_{|t=0}
&=&\sum_{rs}H_{rq}\bar{H}_{sq}(e^{tA}A)_{pr}(-e^{tA}A)_{sp}{\,}_{|t=0}\\
&=&\sum_{rs}H_{rq}\bar{H}_{sq}A_{pr}(-A)_{sp}\\
&=&\sum_rA_{pr}H_{rq}\sum_s(H^*)_{qs}(A^*)_{sp}\\
&=&(AH)_{pq}(H^*A^*)_{qp}\\
&=&|(AH)_{pq}|^2
\end{eqnarray*}

Thus, we are led to the conclusion in the statement.
\end{proof}

In the Fourier matrix case we can go beyond this, and we have:

\begin{proposition}
Given a Fourier matrix $F_G\in M_G(\mathbb C)$, and an antihermitian matrix $A\in M_G(\mathbb C)$, the matrix $H=UF_G$, where $U=e^{tA}$ with $t\in\mathbb R$, is Hadamard when
$$\left|\sum_s\sum_m\frac{t^m}{m!}\sum_{k+l=m}\binom{m}{l}\sum_sA^k_{p,s+n}(-A)^l_{sp}\right|=\delta_{n0}$$
for any $p$, with the indices being $k,l,m\in\mathbb N$, and $n,p,s\in G$.
\end{proposition}

\begin{proof}
According to the formula in the proof of Theorem 7.22 (1), we have:
\begin{eqnarray*}
|(UF_G)_{pq}|^2
&=&\sum_{rs}(F_G)_{rq}(\overline{F_G})_{sq}(e^{tA})_{pr}(e^{-tA})_{sp}\\
&=&\sum_{rs}<r,q><-s,q>(e^{tA})_{pr}(e^{-tA})_{sp}\\
&=&\sum_{rs}<r-s,q>(e^{tA})_{pr}(e^{-tA})_{sp}
\end{eqnarray*}

By setting $n=r-s$, can write this formula in the following way:
\begin{eqnarray*}
|(UF_G)_{pq}|^2
&=&\sum_{ns}<n,q>(e^{tA})_{p,s+n}(e^{-tA})_{sp}\\
&=&\sum_n<n,q>\sum_s(e^{tA})_{p,s+n}(e^{-tA})_{sp}
\end{eqnarray*}

Since this quantity must be 1 for any $q$, we must have:
$$\sum_s(e^{tA})_{p,s+n}(e^{-tA})_{sp}=\delta_{n0}$$

On the other hand, we have the following computation:
\begin{eqnarray*}
\sum_s(e^{tA})_{p,s+n}(e^{-tA})_{sp}
&=&\sum_s\sum_{kl}\frac{(tA)^k_{p,s+n}}{k!}\,\cdot\,\frac{(-tA)^l_{sp}}{l!}\\
&=&\sum_s\sum_{kl}\frac{1}{k!l!}\sum_s(tA)^k_{p,s+n}(-tA)^l_{sp}\\
&=&\sum_s\sum_{kl}\frac{t^{k+l}}{k!l!}\sum_sA^k_{p,s+n}(-A)^l_{sp}\\
&=&\sum_s\sum_mt^m\sum_{k+l=m}\frac{1}{k!l!}\sum_sA^k_{p,s+n}(-A)^l_{sp}\\
&=&\sum_s\sum_m\frac{t^m}{m!}\sum_{k+l=m}\binom{m}{l}\sum_sA^k_{p,s+n}(-A)^l_{sp}
\end{eqnarray*}

Thus, we are led to the conclusion in the statement.
\end{proof}

Following Nicoara-White \cite{nwh}, let us construct now the deformations of $F_G$. The result here, which came a long time after the original defect paper of Tadej-\.Zyczkowski \cite{tz2}, and even more time after the early computations of Karabegov \cite{kar}, appearing somewhat as a total surprise, puzzling all known experts at that time, is as follows:

\index{Nicoara-White theorem}

\begin{theorem}
Let $G$ be a finite abelian group, and for any $g,h\in G$, let us set:
$$B_{pq}=\begin{cases}
1&{\rm if}\ \exists k\in\mathbb N,p=h^kg,q=h^{k+1}g\\
0&{\rm otherwise}
\end{cases}$$
When $(g,h)\in G^2$ range in suitable cosets, the unitary matrices
$$e^{it(B+B^t)}F_G\quad,\quad e^{t(B-B^t)}F_G$$
are both Hadamard, and make the defect of $F_G$ to be attained.
\end{theorem}

\begin{proof}
The proof of this result, from \cite{nwh}, is quite long and technical, based on the Fourier computation from Proposition 7.23, the idea being as follows:

\medskip

(1) First of all, an elementary algebraic study shows that when $(g,h)\in G^2$ range in some suitable cosets, coming from the proof of Theorem 7.14, the various matrices $B=B^{gh}$ constructed above are distinct, the matrices $A=i(B+B^t)$ and $A'=B-B^t$ are linearly independent, and the number of such matrices equals the defect of $F_G$.

\medskip

(2) It is also standard to check that each $B=(B_{pq})$ is a partial isometry, and that $B^k,B^{*k}$ are given by simple formulae. With this ingredients in hand, the Hadamard property follows from the Fourier computation from the proof of Proposition 7.23. Indeed, we can compute the exponentials there, and eventually use the binomial formula.

\medskip

(3) Finally, the matrices in the statement can be shown to be non-equivalent, and this is something more technical, for which we refer to \cite{nwh}. With this last ingredient in hand, a comparison with Theorem 7.14 shows that the defect of $F_G$ is indeed attained, in the sense that all order 0 deformations are actually true deformations. See \cite{nwh}.
\end{proof}

Finally, let us mention that the paper of Nicoara-White \cite{nwh} was written in terms of subfactor-theoretic commuting squares, which is a quite technical operator algebra notion, and with a larger class of commuting squares being actually under investigation. 

\bigskip

We will discuss a bit the relation between Hadamard matrices and commuting squares in chapter 14 below, but in what regards the Nicoara-White theorem, which is the main known theorem regarding the geometry of the complex Hadamard matrices, this definitely remains something to be learned, from their paper \cite{nwh} and their follow-up papers, which are quite technical, and that we would like however to warmly recommend here.

\section*{7e. Exercises} 

Before anything, in connection with the material from the present chapter, we recommend some general geometry reading, with this meaning learning some basic differential and algebraic geometry, if needed. Here is now a first exercise, in connection with the general geometric aspects of the complex Hamadard matrices:

\begin{exercise}
Prove that the Hadamard matrix manifold
$$X_N=M_N(\mathbb T)\cap\sqrt{N}U_N$$
is in general not smooth, and nor it is a complex algebraic manifold.
\end{exercise}

In order to deal with such questions, the best is to try at small values of $N\in\mathbb N$, by using the various classification results from chapter 5. To be more precise, the values $N=2,3$ will certainly not work, so $N=4$ is the case to look at.

\begin{exercise}
Prove that the dephased Hadamard matrix manifold
$$Z_N=\left\{H\in X_N\Big|H_{1j}=H_{i1}=1\right\}$$
is in general not smooth, and not a complex algebraic manifold either.
\end{exercise}

As with the previous exercise, trying $N\in\mathbb N$ small is the way to go, and again, $N=4$ is the precise case to look at, by using the classification results from chapter 5.

\begin{exercise}
Prove that the set $E_N$ formed by the $N\times N$ complex Hadamard matrices modulo the equivalence relation is given by
$$E_N=Z_N\Big/(S_{N-1}\times S_{N-1})$$
and compute this set at $N=2,3,4,5$.
\end{exercise}

As before, in order to solve this problem, the best idea is that of using the various classification results from chapter 5.

\begin{exercise}
Work out the formula of the dephased defect of the Fourier matrix $F_N$, and then of the generalized Fourier matrix $F_G$.
\end{exercise}

As a comment here, if the final formulae do not look very good, this is normal. This exercise is precisely there for showing that the undephased defect is the good quantity to look at, and so that what we did in the above is indeed the thing to do.

\begin{exercise}
Find an alternative proof for the formula
$$d(H)=\frac{N(N+1)}{2}$$
for the real Hadamard matrices, $H\in M_N(\pm1)$.
\end{exercise}

To be more precise here, the above formula was fully proved in the above, by using the general defect equations from the complex case, and then a number of tricks. The problem is that of finding a purely combinatorial proof of this.

\begin{exercise}
Find the defect of the following matrix,
$$K_4=\begin{pmatrix}
-1&1&1&1\\
1&-1&1&1\\
1&1&-1&1\\
1&1&1&-1
\end{pmatrix}$$
via the simplest possible proof.
\end{exercise}

There are many things that can be tried here, such as solving the previous exercise first, and then trying to see if there are simplifications in the case $H=K_4$, or using the general computations that we did for $F_{2,2}^q$, at a suitable value of $q\in\mathbb T$.

\begin{exercise}
Prove that the Tao matrix,
$$T_6=\begin{pmatrix}
1&1&1&1&1&1\\ 
1&1&w&w&w^2&w^2\\ 
1&w&1&w^2&w^2&w\\
1&w&w^2&1&w&w^2\\ 
1&w^2&w^2&w&1&w\\ 
1&w^2&w&w^2&w&1
\end{pmatrix}$$
with $w=e^{2\pi i/3}$, is isolated in the dephased Hadamard matrix manifold.
\end{exercise}

To be more precise, the problem here is that of computing the defect of this matrix $T_6$. Normally this can be done with the defect equations that we have, and some time invested into this problem, or a computer. Alternatively, one can try to find the affine deformations of $T_6$, by using combinatorics and ad-hoc techniques.

\begin{exercise}
Is the defect always equal to the number of $1$ entries?
\end{exercise}

It is of course hard to believe that it is so, and the problem is that of finding the simplest counterexample to this, knowing that the Fourier matrices won't work.

\begin{exercise}
Prove that given two Hadamard matrices $H,K$, we have:
$$d(H\otimes K)\geq d(H)d(K)$$
Is this actually always an equality, or not?
\end{exercise}

Here the first part does not look very difficult, and for the second part we just need a counterexample, based on the various defect computations performed so far.

\begin{exercise}
Develop a defect theory for the partial Hadamard matrices
$$H\in M_{M\times N}(\mathbb T)$$
notably by finding the defect equations, in this setting.
\end{exercise}

This is actually something that we will discuss later in this book, but with no complete proof for the defect equations. Thus, this is a good exercise to be solved now.

\chapter{Special matrices}

\section*{8a. Deformed products}

We have seen in the previous chapter that the defect theory of Tadej-\.Zyczkowski \cite{tz2} can be successfully applied to the real Hadamard matrices, and to the generalized Fourier matrices. Following Avan et al. \cite{aff}, McNulty-Weigert \cite{mwe}, Tadej-\.Zyczkowski \cite{tz1}, \cite{tz2}, and \cite{bop} and other papers, we discuss here a number of more specialized questions, once again in relation with deformations and the defect, regarding the following matrices:

\medskip

-- The tensor products. The main problem here, which quite surprisingly is non-trivial, and even open, is that of computing the defect of the tensor products.

\smallskip

-- The Di\c t\u a deformations of such tensor products. Here the problem is more complicated than for the tensor products, but a few things, however, can be said.

\smallskip

-- The Butson and the regular matrices. Here we have already met, in chapter 6, a conjecture about regular matrices and deformation, so again, things to be done.

\smallskip

-- The master Hadamard matrices. These are some interesting complex Hadamard matrices, introduced by Avan et al. in \cite{aff}, generalizing the Fourier matrices.

\smallskip

-- The McNulty-Weigert matrices. These are again interesting complex Hadamard matrices, introduced by McNulty-Weigert in \cite{mwe}, which are quite often isolated.

\smallskip

-- The partial Hadamard matrices. Here there are, again, many things to be done, following \cite{bop}, inspired by the theory from the square matrix case.

\medskip

Let us begin with the tensor products. As already mentioned, this is a very interesting topic, which is far from being trivial, and to start with, we have the following result, coming straight from the general defect equations, found in chapter 7:

\index{tensor product}

\begin{proposition}
For a tensor product $L=H\otimes K$ we have 
$$d(L)\geq d(H)d(K)$$
coming from an inclusion of linear spaces, as follows:
$$\widetilde{T}_HX_M\otimes\widetilde{T}_KX_N\subset\widetilde{T}_LX_{MN}$$
Moreover, the above inequality is not an equality, in general.
\end{proposition}

\begin{proof}
We have several things to be proved, the idea being as follows:

\medskip

(1) Let us first prove that we have the inclusion of linear spaces in the statement. For this purpose, we use the defect equations found in chapter 7, namely:
$$\sum_kL_{ik}\bar{L}_{jk}(A_{ik}-A_{jk})=0$$

For a tensor product $A=B\otimes C$, we have the following formula:
\begin{eqnarray*}
\sum_{kc}(H\otimes K)_{ia,kc}\overline{(H\otimes K)}_{jb,kc}A_{ia,kc}
&=&\sum_{kc}H_{ik}K_{ac}\cdot\bar{H}_{jk}\bar{K}_{bc}\cdot B_{ik}C_{ac}\\
&=&\sum_kH_{ik}\bar{H}_{jk}B_{ik}\sum_cK_{ac}\bar{K}_{bc}C_{ac}
\end{eqnarray*}

On the other hand, we have as well the following formula:
\begin{eqnarray*}
\sum_{kc}(H\otimes K)_{ia,kc}\overline{(H\otimes K)}_{jb,kc}A_{jb,kc}
&=&\sum_{kc}H_{ik}K_{ac}\cdot\bar{H}_{jk}bar{K}_{bc}\cdot B_{jk}C_{bc}\\
&=&\sum_kH_{ik}\bar{H}_{jk}B_{jk}\sum_cK_{ac}\bar{K}_{bc}C_{bc}
\end{eqnarray*}

Now by assuming $B\in\widetilde{T}_HX_M$ and $C\in\widetilde{T}_KX_N$, the two quantities on the right in the above formulae are equal. Thus we have indeed $A\in\widetilde{T}_LX_{MN}$, as desired.

\medskip

(2) The defect inequality $d(L)\geq d(H)d(K)$ follows from (1).

\medskip

(3) Regarding now the equality case, this does not happen, even in very simple cases. For instance if we consider two Fourier matrices $F_2$, we know from chapter 7 that:
$$d(F_2\otimes F_2)=10>9=d(F_2)^2$$

There are of course many other counterexamples that can be constructed.
\end{proof}

Generally speaking, it is quite hard to go beyond the above result. In fact, besides the isotypic decomposition results from chapter 7, valid for the Fourier matrices, there does not seem to be anything conceptual on this subject. We will be back to this, however, in Theorem 8.3 below, with a slight advance on all this.

\bigskip

In what regards now the computation of the defect for the Di\c t\u a deformations, which generalize the usual tensor products, this is an even more difficult question. Our only result here will concern the case where the deformation matrix is generic:

\index{generic deformation}

\begin{definition}
A rectangular matrix $Q\in M_{M\times N}(\mathbb T)$ is called ``dephased and elsewhere generic'' if the entries on its first row and column are all equal to $1$, and the remaining $(M-1)(N-1)$ entries are algebrically independent over $\mathbb Q$.
\end{definition}

Here the last condition takes of course into account the fact that the entries of $Q$ themselves have modulus 1, the independence assumption being modulo this fact. With this convention made, we have the following result:

\begin{theorem}
Assume that $H\in X_M,K\in X_N$ are dephased, of Butson type, and that $Q\in M_{M\times N}(\mathbb T)$ is dephased and elsewhere generic. We have then
$$A=(A_{ia,kc})\in \widetilde{T}_{H\otimes_QK}X_{MN}$$
when the following equations are satisfied,
$$A_{ac}^{ij}=A_{bc}^{ij}\quad,\quad 
A_{ac}^{ij}=\overline{A_{ac}^{ji}}\quad,\quad 
(A_{xy}^{ii})_{xy}\in\widetilde{T}_KX_N$$
for any $a,b,c$ and $i\neq j$, where:
$$A_{ac}^{ij}=\sum_kH_{ik}\bar{H}_{jk}A_{ia,kc}$$
\end{theorem}

\begin{proof}
Consider the standard system of equations for the enveloping tangent space in the statement, coming from the results in chapter 7, namely:
$$\sum_{kc}(H\otimes_QK)_{ia,kc}\overline{(H\otimes_QK)}_{jb,kc}(A_{ia,kc}-A_{jb,kc})=0$$

We have the following formula, for our matrix:
$$(H\otimes_QK)_{ia,jb}=q_{ib}H_{ij}K_{ab}$$ 

Thus, our system of equations is as follows:
$$\sum_cq_{ic}\bar{q}_{jc}K_{ac}\bar{K}_{bc}\sum_kH_{ik}\bar{H}_{jk}(A_{ia,kc}-A_{jb,kc})=0$$

Consider now the variables in the statement, namely:
$$A_{ac}^{ij}=\sum_kH_{ik}\bar{H}_{jk}A_{ia,kc}$$

The conjugates of these variables are given by:
$$\overline{A_{ac}^{ij}}
=\sum_k\bar{H}_{ik}H_{jk}A_{ia,kc}
=\sum_kH_{jk}\bar{H}_{ik}A_{ia,kc}$$

Thus, in terms of these variables, our system becomes simply:
$$\sum_cq_{ic}\bar{q}_{jc}K_{ac}\bar{K}_{bc}(A_{ac}^{ij}-\overline{A_{bc}^{ji}})=0$$

More precisely, the above equations must hold for any $i,j,a,b$. By distinguishing now two cases, depending on whether $i,j$ are equal or not, the situation is as follows:

\medskip

(1) Case $i\neq j$. In this case, let us look at the row vector of parameters, namely:
$$(q_{ic}\bar{q}_{jc})_c=(1,q_{i1}\bar{q}_{j1},\ldots,q_{iM}\bar{q}_{jM})$$

Since the matrix $Q$ was assumed to be dephased and elsewhere generic, and because of our assumption $i\neq j$, the entries of the above vector are linearly independent over $\bar{\mathbb Q}$. But, since by linear algebra we can restrict the attention to the computation of the solutions over $\bar{\mathbb Q}$, the $i\neq j$ part of our system simply becomes:
$$A_{ac}^{ij}=\overline{A_{bc}^{ji}}\quad,\quad\forall a,b,c,\forall i\neq j$$

Now by making now $a,b,c$ vary, we are led to the following equations:
$$A_{ac}^{ij}=A_{bc}^{ij},\quad A_{ac}^{ij}=\overline{A_{ac}^{ji}},\quad\forall a,b,c,i\neq j$$

(2) Case $i=j$. In this case the $q$ parameters cancel, and our equations become:
$$\sum_cK_{ac}\bar{K}_{bc}(A_{ac}^{ii}-\overline{A_{bc}^{ii}})=0,\quad\forall a,b,c,i$$

Now observe that we have the following formula:
$$A_{ac}^{ii}=\sum_kA_{ia,kc}$$

Thus, our equations simply become:
$$\sum_cK_{ac}\bar{K}_{bc}(A_{ac}^{ii}-A_{bc}^{ii})=0,\quad\forall a,b,c,i$$

But these are precisely the equations for the space $\widetilde{T}_KX_N$, and we are done.
\end{proof}

Let us go back now to usual tensor products, and look at the affine cones. In view of the inclusion from Proposition 8.1, the problem is that of finding the biggest subcone of $T_{H\otimes K}^\circ X_{MN}$, obtained by gluing $T_H^\circ X_M,T_K^\circ X_N$. Our answer here, taking into account the two ``semi-trivial'' cones coming from left and right Di\c t\u a deformations, is as follows:

\index{tangent cone gluing}

\begin{theorem}
The cones $T_H^\circ X_M=\{B\}$ and $T_K^\circ X_N=\{C\}$ glue via the formulae
$$A_{ia,jb}=\lambda B_{ij}+\psi_jC_{ab}+X_{ia}+Y_{jb}+F_{aj}$$
$$A_{ia,jb}=\phi_bB_{ij}+\mu C_{ab}+X_{ia}+Y_{jb}+E_{ib}$$
producing in this way two subcones of the affine cone $T_{H\otimes K}^\circ X_{MN}=\{A\}$.
\end{theorem}

\begin{proof}
The idea will be that $X_{ia},Y_{jb}$ are the trivial parameters, and that $E_{ib},F_{aj}$ are the Di\c t\u a parameters. Given a matrix $A=(A_{ia,jb})$, consider the following quantity:
$$P=\sum_{kc}H_{ik}\bar{H}_{jk}K_{ac}\bar{K}_{bc}q^{A_{ia,kc}-A_{jb,kc}}$$

Let us prove now the first statement, namely that for any choice of matrices $B\in T_H^\circ X_M,C\in T_H^\circ X_N$ and of parameters $\lambda,\psi_j,X_{ia},Y_{jb},F_{aj}$, the first matrix $A=(A_{ia,jb})$ constructed in the statement belongs indeed to $T_{H\otimes K}^\circ X_{MN}$. We have: 
$$A_{ia,kc}=\lambda B_{ik}+\psi_kC_{ac}+X_{ia}+Y_{kc}+F_{ak}$$
$$A_{jb,kc}=\lambda B_{jk}+\psi_kC_{bc}+X_{jb}+Y_{kc}+F_{bk}$$

Now by substracting these equations, we obtain:
$$A_{ia,kc}-A_{jb,kc}
=\lambda(B_{ik}-B_{jk})+\psi_k(C_{ac}-C_{bc})+(X_{ia}-X_{jb})+(F_{ak}-F_{bk})$$

It follows that the above quantity $P$ is given by:
\begin{eqnarray*}
P
&=&\sum_{kc}H_{ik}\bar{H}_{jk}K_{ac}\bar{K}_{bc}q^{\lambda(B_{ik}-B_{jk})+\psi_k(C_{ac}-C_{bc})+(X_{ia}-X_{jb})+(F_{ak}-F_{bk})}\\
&=&q^{X_{ia}-X_{jb}}\sum_kH_{ik}\bar{H}_{jk}q^{F_{ak}-F_{bk}}q^{\lambda(B_{ik}-B_{jk})}\sum_cK_{ac}\bar{K}_{bc}(q^{\psi_k})^{C_{ac}-C_{bc}}\\
&=&\delta_{ab}q^{X_{ia}-X_{ja}}\sum_kH_{ik}\bar{H}_{jk}(q^\lambda)^{B_{ik}-B_{jk}}\\
&=&\delta_{ab}\delta_{ij}
\end{eqnarray*}

We conclude that we have, as claimed:
$$A\in T_{H\otimes K}^\circ X_{MN}$$

In the second case now, the proof is similar. First, we have:
$$A_{ia,kc}=\phi_cB_{ik}+\mu C_{ac}+X_{ia}+Y_{kc}+E_{ic}$$
$$A_{jb,kc}=\phi_cB_{jk}+\mu C_{bc}+X_{jb}+Y_{kc}+E_{jc}$$

Thus by substracting, we obtain:
$$A_{ia,kc}-A_{jb,kc}
=\phi_c(B_{ik}-B_{jk})+\mu(C_{ac}-C_{bc})+(X_{ia}-X_{jb})+(E_{ic}-E_{jc})$$

It follows that the above quantity $P$ is given by:
\begin{eqnarray*}
P
&=&\sum_{kc}H_{ik}\bar{H}_{jk}K_{ac}\bar{K}_{bc}q^{\phi_c(B_{ik}-B_{jk})+\mu(C_{ac}-C_{bc})+(X_{ia}-X_{jb})+(E_{ic}-E_{jc})}\\
&=&q^{X_{ia}-X_{jb}}\sum_cK_{ac}\bar{K}_{bc}q^{E_{ic}-E_{jc}}q^{\mu(C_{ac}-C_{bc})}\sum_kH_{ik}\bar{H}_{jk}(q^{\phi_c})^{B_{ik}-B_{jk}}\\
&=&\delta_{ij}q^{X_{ia}-X_{ib}}\sum_cK_{ac}\bar{K}_{bc}(q^\mu)^{C_{ac}-C_{bc}}\\
&=&\delta_{ij}\delta_{ab}
\end{eqnarray*}

Thus, we are led to the conclusion in the statement.
\end{proof}

We believe Theorem 8.4 to be ``optimal'', in the sense that nothing more can be said about the affine tangent spaces of type $T_{H\otimes K}^\circ X_{MN}$, in the general case, besides what has been said there. However, this is something rather conjectural. As a continuation now of all this, bringing us into some concrete, interesting mathematics, let us discuss some rationality questions, in relation with the following definition:

\index{rational defect}

\begin{definition}
The rational defect of $H\in X_N$ is the following number:
$$d_\mathbb Q(H)=\dim_\mathbb Q(\widetilde{T}_HC_N\cap M_N(\mathbb Q))$$
The vector space on the right is called rational enveloping tangent space at $H$.
\end{definition}

As a first observation, this notion can be extended to all the tangent cones at $H$, and by using an arbitrary field $\mathbb K\subset\mathbb C$ instead of $\mathbb Q$. Indeed, we can set:
$$T_H^*X_N(\mathbb K)=T_H^*X_N\cap M_N(\mathbb K)$$

However, in what follows we will be interested only in the objects constructed in Definition 8.5. It follows from definitions that $d_\mathbb Q(H)\leq d(H)$, and we have:

\index{Rationality Conjecture}

\begin{conjecture}[Rationality]
For the Butson matrices we have:
$$d_\mathbb Q(H)=d(H)$$
That is, for such matrices, the defect equals the rational defect. 
\end{conjecture}

More generally, we believe that the above equality should hold in the regular matrix case. However, since the regular matrix case is not known to fully cover the Butson matrix case, as explained in chapter 6, we prefer to state our conjecture as above. As a first piece of evidence now, we have the following elementary result:

\begin{theorem}
The rationality conjecture holds for $H\in H_N(l)$ with $l=2,3,4,6$.
\end{theorem}

\begin{proof}
Let us recall that the equations for the enveloping tangent space are:
$$\sum_kH_{ik}\bar{H}_{jk}(A_{ik}-A_{jk})=0$$

With these equations in hand, the proof goes as follows:

\medskip

\underline{Case $l=2$}. Here the above equations are all real, and have $\pm1$ coefficients, so in particular, have rational coefficients. 

\medskip

\underline{Case $l=3$}. Here we can use the fact that, with $w=e^{2\pi i/3}$, the real solutions of $x+wy+w^2z=0$ are those satisfying $x=y=z$. We conclude that once again our system, after some manipulations, is equivalent to a real system having rational coefficients. 

\medskip

\underline{Case $l=4$}. Here the coefficients are $1,i,-1,-i$, so by taking the real and imaginary parts, we reach once again to a system with rational coefficients.

\medskip

\underline{Case $l=6$}. Here the study is similar to the study at $l=3$.

\medskip

Thus, in all cases under investigation, $l=2,3,4,6$, we have a real system with rational coefficients, and the result follows from standard linear algebra.
\end{proof}

Observe that the above method cannot work at $l=5$, where the equation $a+wb+w^2c+w^3d+w^4e=0$ with $w=e^{2\pi i/5}$ and $a,b,c,d,e\in\mathbb R$ can have exotic solutions. Let us prove now that Conjecture 8.6 is verified for the Fourier matrices. We say that a matrix $L^{rs}$ over the group $\mathbb Z_{p^r}\times\mathbb Z_{p^s}$ is dephased if its nonzero entries belong to:
$$X_{rs}=(\mathbb Z_{p^r}-\mathbb Z_{p^{r-1}})\times(\mathbb Z_{p^s}-\mathbb Z_{p^{s-1}})$$

Here, and in what follows, we use the convention $\mathbb Z_{p^{-1}}=\emptyset$. We have:

\begin{proposition}
For $F=F_{p^a}$, the elements $A\in\widetilde{T}_FC_N$ are the solutions of
$$A_{ij}=\sum_{r+s\leq a}L^{rs}_{p^{a-r}i,p^{a-s}j}$$
where the $L$ variables are free, and form dephased matrices $L^{rs}$.
\end{proposition}

\begin{proof}
The number of $L$ variables is given by:
\begin{eqnarray*}
d
&=&\sum_{r+s\leq a}|\mathbb Z_{p^r}-\mathbb Z_{p^{r-1}}|\cdot|\mathbb Z_{p^s}-\mathbb Z_{p^{s-1}}|\\
&=&\sum_{r\leq a}p^{a-r}|\mathbb Z_{p^r}-\mathbb Z_{p^{r-1}}|\\
&=&p^a+\sum_{r=1}^ap^{a-r}(p^r-p^{r-1})\\
&=&p^a+a(p-1)p^{a-1}\\
&=&(p+ap-a)p^{a-1}
\end{eqnarray*}

Thus the number of $L$ variables equals the defect $d(F)$, so it is indeed the good one. As for the proof now, in the general case, this is quite similar to the one at $a=1,2$. More precisely, consider the map $L\to A$. This map is linear, and in view of the above calculation, it is enough to prove that this map is injective, and has the correct target:

\medskip

(1) For the injectivity part, recall that at $a=2$ the formula in the statement reads:
$$A_{ij}=L^{00}_{00}+L^{01}_{0,pj}+L^{10}_{pi,0}+L^{02}_{0j}+L^{20}_{i0}+L^{11}_{pi,pj}$$

Now assume $A=0$. Then with $i=j=0$ we get $L^{00}_{00}=0$. Using this, with $i=0$ and $pj=0,j\neq 0$ we get $L^{00}_{00}+L^{02}_{0j}=0$, and so $L^{02}_{0j}=0$. So, with $i=0$ and $pj\neq 0$ we therefore obtain $L^{00}_{00}+L^{02}_{0j}+L^{01}_{0,pj}=0$, and so $L^{01}_{0,pj}=0$. Now the same method gives as well succesively $L^{20}_{i0}=0$ and $L^{10}_{pi,0}=0$, so we are left with $A_{ij}=L^{11}_{pi,pj}$, so we must have $L^{11}_{pi,pj}=0$ as well, and we are done. This method works of course for any $a\in\mathbb N$.

\medskip

(2) Regarding now the ``target'' part, we must prove $A\in\widetilde{T}_FC_N$. The equations are:
$$\sum_kw^{(i-j)k}\left(\sum_{r+s\leq a}L^{rs}_{p^{a-r}i,p^{a-s}k}-L^{rs}_{p^{a-r}j,p^{a-s}k}\right)=0$$

So, for any indices $i,j$ and any $r+s\leq a$, we must prove that we have:
$$\sum_kw^{(i-j)k}\left(L^{rs}_{p^{a-r}i,p^{a-s}k}-L^{rs}_{p^{a-r}j,p^{a-s}k}\right)=0$$

In order to do this, consider the following quantity:
$$X_{il}=\frac{1}{p^a}\sum_kw^{lk}L^{rs}_{p^{a-r}i,p^{a-s}k}$$

We must prove $X_{i,i-j}=X_{j,i-j}$. But, with $k=m+p^sn$, we have:
\begin{eqnarray*}
X_{il}
&=&\frac{1}{p^a}\sum_nw^{lp^sn}\sum_mw^{lm}L^{rs}_{p^{a-r}i,p^{a-s}m}\\
&=&\delta_{l0}\sum_mw^{lm}L^{rs}_{p^{a-r}i,p^{a-s}m}
\end{eqnarray*}

Thus we have $l\neq 0\implies X_{il}=0$, and so $X_{i,i-j}=X_{j,i-j}$ and we are done.
\end{proof}

By using the above result, we obtain:

\index{isotypic Fourier matrix}

\begin{proposition}
For an isotypic Fourier matrix, $H=F_N$ with $N=p^a$, we have
$$T_H^\circ C_N=T_HC_N=\widetilde{T}_HC_N=\left\{A\in M_N(\mathbb R)\Big|A_{ij}=\sum_{r+s\leq a}L^{rs}_{p^{a-r}i,p^{a-s}j}\right\}$$
where the $L$ variables are free, and form dephased matrices $L^{rs}$. 
\end{proposition}

\begin{proof}
We just have to show that the defect of $F_N$ is exhausted by affine deformations. With $k=m+p^sn$, as in the proof of Proposition 8.8, we have:
\begin{eqnarray*}
\sum_kH_{ik}\bar{H}_{jk}q^{A_{ik}-A_{jk}}
&=&\sum_kw^{(i-j)k}\prod_{r+s\leq a}q^{L^{rs}_{p^{a-r}i,p^{a-s}k}-L^{rs}_{p^{a-r}j,p^{a-s}k}}\\
&=&\sum_nw^{(i-j)p^sn}\sum_mw^{(i-j)m}\prod_{r+s\leq a}q^{L^{rs}_{p^{a-r}i,p^{a-s}m}-L^{rs}_{p^{a-r}j,p^{a-s}m}}\\
&=&\delta_{ij}p^a\sum_mw^{(i-j)m}\prod_{r+s\leq a}q^{L^{rs}_{p^{a-r}i,p^{a-s}m}-L^{rs}_{p^{a-r}j,p^{a-s}m}}
\end{eqnarray*}

Now since this quantity vanishes for $i\neq j$, this gives the result.
\end{proof}

Observe that the above result shows that Conjecture 8.6 holds for the isotypic Fourier matrices. We will see in what follows that the same happens for any Fourier matrix. In order now to discuss the general case, $H=F_N$, we will need:

\begin{proposition}
If $G=H\times K$ is such that $(|H|,|K|)=1$, the canonical inclusion
$$\widetilde{T}_{F_H}C_{|H|}\otimes\widetilde{T}_{F_K}C_{|K|}\subset\widetilde{T}_{F_G}C_{|G|}$$
constructed in Proposition 8.1 is an isomorphism.
\end{proposition}

\begin{proof}
We have $F_G=F_{H\times K}$, and the defect of this matrix is given by:
\begin{eqnarray*}
d(F_{H\times K})
&=&\sum_{(h,k)\in H\times K}\frac{|H\times K|}{ord(h,k)}\\
&=&\sum_{(h,k)\in H\times K}\frac{|H\times K|}{ord(h)ord(k)}\\
&=&d(F_H)d(F_K)
\end{eqnarray*}

Thus the inclusion in the statement must be indeed an isomorphism.
\end{proof}

With the above result in hand, the idea now will be simply to ``glue'' the various isotypic formulae coming from Proposition 8.9. Indeed, let us recall from there that in the isotypic case, $N=p^a$, the parameter set for the enveloping tangent space is:
$$X(p^a)=\bigsqcup_{r+s\leq a}(\mathbb Z_{p^r}-\mathbb Z_{p^{r-1}})\times (\mathbb Z_{p^s}-\mathbb Z_{p^{s-1}})$$

Now since the defect is multiplicative over isotypic components, the parameter set in the general case, $N=p_1^{a_1}\ldots p_k^{a_k}$, will be simply given by:
$$X(p_1^{a_1}\ldots p_k^{a_k})=X(p_1^{a_1})\times\ldots\times X(p_k^{a_k})$$

We can obtain from this an even simpler description of the parameter set, just by expanding the product, and gluing the group components. Indeed, let us start with:

\begin{definition}
Given a finite abelian group $G=\mathbb Z_{p_1^{r_1}}\times\ldots\times\mathbb Z_{p_k^{r_k}}$ we set:
$$G^\circ=(\mathbb Z_{p_1^{r_1}}-\mathbb Z_{p_1^{r_1-1}})\times\ldots\times(\mathbb Z_{p_k^{r_k}}-\mathbb Z_{p_k^{r_k-1}})$$
A matrix $L\in M_{G\times H}(\mathbb R)$ will be called dephased if $L_{ij}=0$ for any $(i,j)\not\in G^\circ\times H^\circ$.
\end{definition}

Observe now that, with the above notation $G^\circ$, the parameter set discussed above is given by the following simple formula:
$$X(N)=\bigsqcup_{G\times H\subset\mathbb Z_N}G^\circ\times H^\circ$$

In addition, we can see that the collection of dephased matrices $L\in M_{G\times H}(\mathbb R)$ , over all possible configurations $G\times H\subset\mathbb Z_N$, takes its parameters precisely in $X(N)$. In order now to formulate our main result, we will need one more definition, as follows:

\begin{definition}
Given $N=p_1^{a_1}\ldots p_k^{a_k}$ and a subgroup $G\subset \mathbb Z_N$, we set 
$$\varphi_G(i_1,\ldots,i_k)=(p_1^{a_1-r_1}i_1,\ldots p_k^{a_k-r_k}i_k)$$
where the exponents $r_i\leq a_i$ are given by $G=\mathbb Z_{p_1^{r_1}}\times\ldots\times\mathbb Z_{p_k^{r_k}}$.
\end{definition}

Observe that in the case $k=1$ this function is precisely the one appearing in Proposition 8.9. In fact, we have the following generalization of Proposition 8.9:

\begin{theorem}
For $H=F_N$ the vectors $A\in\widetilde{T}_HC_N$ appear as plain sums of type
$$A_{ij}=\sum_{G\times H\subset\mathbb Z_N}L^{GH}_{\varphi_G(i)\varphi_H(j)}$$
where the $L$ variables form dephased matrices $L^{GH}\in M_{G\times H}(\mathbb R)$. 
\end{theorem}

\begin{proof}
According to the above discussion, we just have to glue the various isotypic formulae coming from Proposition 8.9. The gluing formula reads:
\begin{eqnarray*}
A_{i_1\ldots i_k,j_1\ldots j_k}
&=&A_{i_1j_1}\ldots A_{i_kj_k}\\
&=&\left(\sum_{r_1+s_1\leq a_1}L^{r_1s_1p_1}_{p_1^{a_1-r_1}i_1,p_1^{a_1-s_1}j_1}\ldots\sum_{r_k+s_k\leq a_k}L^{r_ks_kp_k}_{p_k^{a_k-r_k}i_k,p_k^{a_k-s_k}j_k}\right)\\
&=&\sum_{r_1+s_1\leq a_1}\ldots \sum_{r_k+s_k\leq a_k}L^{r_1s_1p_1}_{p_1^{a_1-r_1}i_1,p_1^{a_1-s_1}j_1}\ldots L^{r_ks_kp_k}_{p_k^{a_k-r_k}i_k,p_k^{a_k-s_k}j_k}
\end{eqnarray*}

Now, let us introduce the following variables:
$$L^{r_1\ldots r_k,s_1\ldots s_k}_{i_1\ldots i_k,j_1\ldots j_k}=L^{r_1s_1}_{i_1j_1}\ldots L^{r_ks_k}_{i_kj_k}$$

In terms of these new variables, the gluing formula reads:
$$A_{i_1\ldots i_k,j_1\ldots j_k}=\sum_{r_1+s_1\leq a_1}\ldots \sum_{r_k+s_k\leq a_k}L^{r_1\ldots r_k,s_1\ldots s_k}_{p_1^{a_1-r_1}i_1,\ldots p_k^{a_k-r_k}i_k,p_1^{a_1-r_1}j_1\ldots p_k^{a_k-r_k}j_k}$$

Together with the fact that the new $L$ variables form dephased matrices, in the sense of Definition 8.11, this gives the result.
\end{proof}

As a main consequence, we have the following result:

\index{Rationality Conjecture}

\begin{theorem}
The rationality conjecture holds for the Fourier matrices.
\end{theorem}

\begin{proof}
Indeed, the formula in Theorem 8.13 shows that for $H=F_N$ the rational defect, as constructed in Definition 8.5, counts the same variables as the usual defect.
\end{proof}

\section*{8b. Master matrices}

Let us discuss now some defect computations for an interesting class of Hadamard matrices, namely the ``master'' ones, introduced by Avan et al. in \cite{aff}:

\index{master Hadamard matrix}
\index{master function}

\begin{definition}
A master Hadamard matrix is an Hadamard matrix of the form 
$$H_{ij}=\lambda_i^{n_j}$$
with $\lambda_i\in\mathbb T,n_j\in\mathbb R$. The associated ``master function'' is $f(z)=\sum_jz^{n_j}$.
\end{definition}

Observe that with $\lambda_i=e^{im_i}$ we have $H_{ij}=e^{im_in_j}$. The basic example of such a matrix is the Fourier matrix $F_N$, having master function as follows:
$$f(z)=\frac{z^N-1}{z-1}$$

Observe that, in terms of $f$, the Hadamard condition on $H$ is simply:
$$f\left(\frac{\lambda_i}{\lambda_j}\right)=N\delta_{ij}$$

These matrices were introduced in \cite{aff}, the motivating remark there being the fact that the following operator defines a representation of the Temperley-Lieb algebra \cite{tli}:
$$R=\sum_{ij}e_{ij}\otimes\Lambda^{n_i-n_j}$$

\index{Temperley-Lieb algebra}

At the level of examples, the first observation, from \cite{aff}, is that the standard $4\times 4$ complex Hadamard matrices are, with 2 exceptions, master Hadamard matrices:

\begin{proposition}
The following complex Hadamard matrix, with $|q|=1$,
$$F_{2,2}^q=\begin{pmatrix}
1&1&1&1\\
1&-1&1&-1\\
1&q&-1&-q\\
1&-q&-1&q
\end{pmatrix}$$
is a master Hadamard matrix, for any $q\neq\pm1$.
\end{proposition}

\begin{proof}
We use the exponentiation convention $(e^{it})^r=e^{itr}$, for $t\in[0,2\pi)$ and $r\in\mathbb R$. Since we have $q^2\neq1$, we can find $k\in\mathbb R$ such that:
$$q^{2k}=-1$$

In terms of this parameter $k\in\mathbb R$, our matrix becomes:
$$F_{2,2}^q=\begin{pmatrix}
1^0&1^1&1^{2k}&1^{2k+1}\\
(-1)^0&(-1)^1&(-1)^{2k}&(-1)^{2k+1}\\
q^0&q^1&q^{2k}&q^{2k+1}\\
(-q)^0&(-q)^1&(-q)^{2k}&(-q)^{2k+1}\\
\end{pmatrix}$$

Now let us pick $\lambda\neq1$ and write, by using our exponentiation convention above:
$$1=\lambda^x\quad,\quad -1=\lambda^y$$
$$q=\lambda^z\quad,\quad -q=\lambda^t$$

But this gives the formula in the statement.
\end{proof}

Observe that the above result shows that any Hamadard matrix at $N\leq5$ is master Hadamard. We have the following generalization of it, once again from \cite{aff}:

\index{deformed Fourier matrix}

\begin{theorem}
The deformed Fourier matrices $F_M\otimes_QF_N$ are master Hadamard, for any parameter matrix $Q\in M_{M\times N}(\mathbb T)$ of the form
$$Q_{ib}=q^{i(Np_b+b)}$$
where $q=e^{2\pi i/MNk}$ with $k\in\mathbb N$, and $p_0,\ldots,p_{N-1}\in\mathbb R$.
\end{theorem}

\begin{proof}
The main construction in \cite{aff}, in connection with deformations, that we will follow here, is in terms of master functions as follows:
$$f(z)=f_M(z^{Nk})f_N(z)$$

Here $k\in\mathbb N$, and the functions on the right are by definition as follows:
$$f_M(z)=\sum_iz^{Mr_i+i}\quad,\quad 
f_N(z)=\sum_az^{Np_a+a}$$

We use the eigenvalues $\lambda_{ia}=q^iw^a$, where $w=e^{2\pi i/N}$, and where $q^{Nk}=\nu$, where $\nu^M=1$. We have $f(z)=f_M(z^{Nk})f_N(z)$, so the exponents are:
$$n_{jb}=Nk(Mr_j+j)+Np_b+b$$

Thus the associated master Hadamard matrix is given by:
\begin{eqnarray*}
H_{ia,jb}
&=&(q^iw^a)^{Nk(Mr_j+j)+Np_b+b}\\
&=&\nu^{ij}q^{i(Np_b+b)}w^{a(Np_b+b)}\\
&=&\nu^{ij}w^{ab}q^{i(Np_b+b)}
\end{eqnarray*}

Now let us recall that we have the following formula, for the tensor product:
$$(F_M\otimes F_N)_{ia,jb}=\nu^{ij}w^{ab}$$

Thus we have as claimed $H=F_M\otimes_QF_N$, with $Q$ being as follows:
$$Q_{ib}=q^{i(Np_b+b)}$$

Finally, observe that $Q$ itself is a ``master matrix'' in our sense, because the indices split. Thus, we are led to the conclusions in the statement.
\end{proof}

In view of the above examples, and of the lack of other known examples of master Hadamard matrices, the following conjecture was made in \cite{aff}:

\index{Master Hadamard Conjecture}

\begin{conjecture}[Master Hadamard Conjecture]
The master Hadamard matrices appear as Di\c t\u a deformations of $F_N$.
\end{conjecture}

There is a relation here with the notions of defect and isolation, that we would like to discuss now. First, we have the following defect computation:

\index{defect}

\begin{theorem}
The defect of a master Hadamard matrix is given by
$$d(H)=\dim_\mathbb R\left\{B\in M_N(\mathbb C)\Big|\bar{B}=\frac{1}{N}BL, (BR)_{i,ij}=(BR)_{j,ij}\ \forall i,j\right\}$$
where the matrices on the right are given by
$$L_{ij}=f\left(\frac{1}{\lambda_i\lambda_j}\right)\quad,\quad
R_{i,jk}=f\left(\frac{\lambda_j}{\lambda_i\lambda_k}\right)$$
with $f$ being the master function.
\end{theorem}

\begin{proof}
The first order deformation equations from chapter 7 are as follows:
$$\sum_kH_{ik}\bar{H}_{jk}(A_{ik}-A_{jk})=0$$

In our case, with $H_{ij}=\lambda_i^{n_j}$ we have the following formula:
$$H_{ij}\bar{H}_{jk}=\left(\frac{\lambda_i}{\lambda_j}\right)^{n_k}$$

Thus, the defect is given by the following formula:
$$d(H)=\dim_\mathbb R\left\{A\in M_N(\mathbb R)\Big|\sum_kA_{ik}\left(\frac{\lambda_i}{\lambda_j}\right)^{n_k}=\sum_kA_{jk}\left(\frac{\lambda_i}{\lambda_j}\right)^{n_k}\ \forall i,j\right\}$$

Now, pick $A\in M_N(\mathbb C)$ and set $B=AH^t$. We have the following formula:
$$A=\frac{1}{N}B\bar{H}$$

By using this formula, we have the following computation:
\begin{eqnarray*}
A\in M_N(\mathbb R)
&\iff&B\bar{H}=\bar{B}H\\
&\iff&\bar{B}=\frac{1}{N}B\bar{H}H^*
\end{eqnarray*}

On the other hand, the matrix on the right is given by:
$$(\bar{H}H^*)_{ij}
=\sum_k\bar{H}_{ik}\bar{H}_{jk}
=\sum_k(\lambda_i\lambda_j)^{-n_k}
=L_{ij}$$

Thus $A\in M_N(\mathbb R)$ if and only the condition $\bar{B}=\frac{1}{N}BL$ in the statement is satisfied. Regarding now the second condition on $A$, observe that with $A=\frac{1}{N}B\bar{H}$ we have:
\begin{eqnarray*}
\sum_kA_{ik}\left(\frac{\lambda_i}{\lambda_j}\right)^{n_k}
&=&\frac{1}{N}\sum_{ks}B_{is}\left(\frac{\lambda_i}{\lambda_j\lambda_s}\right)^{n_k}\\
&=&\frac{1}{N}\sum_sB_{is}R_{s,ij}\\
&=&\frac{1}{N}(BR)_{i,ij}
\end{eqnarray*}

Thus the second condition on $A$ simply reads:
$$(BR)_{i,ij}=(BR)_{j,ij}$$

But this leads to the conclusions in the statement.
\end{proof}

We refer to \cite{aff} and related papers for more on the master Hadamard matrices. In what follows we will not discuss them further, but we will be back to a related topic, namely Temperley-Lieb algebra representations coming from the complex Hadamard matrices, in chapters 13-16 below, when talking about quantum permutation groups.

\section*{8c. Isolated matrices}

Let us discuss now yet another interesting construction of complex Hadamard matrices, due to McNulty and Weigert \cite{mwe}. The matrices constructed there generalize the Tao matrix $T_6$, and usually have the interesting feature of being isolated. The construction in \cite{mwe} uses the theory of MUB, as developed in \cite{bbe}, \cite{deb}, but we will follow here a more direct approach, from \cite{bop}. The starting observation from \cite{mwe} is as follows:

\index{MUB}
\index{McNulty-Weigert matrix}

\begin{theorem}
Assuming that $K\in M_N(\mathbb C)$ is Hadamard, so is the matrix
$$H_{ia,jb}=\frac{1}{\sqrt{Q}}K_{ij}(L_i^*R_j)_{ab}$$
provided that $\{L_1,\ldots,L_N\}\subset\sqrt{Q}U_Q$ and $\{R_1,\ldots,R_N\}\subset\sqrt{Q}U_Q$ are such that
$$\frac{1}{\sqrt{Q}}L_i^*R_j\in\sqrt{Q}U_Q$$
with $i,j=1,\ldots,N$, are complex Hadamard.
\end{theorem}

\begin{proof}
The check of the unitarity of the matrix in the statement can be done as follows, by using our various assumptions on the various matrices involved:
\begin{eqnarray*}
<H_{ia},H_{kc}>
&=&\frac{1}{Q}\sum_{jb}K_{ij}(L_i^*R_j)_{ab}\bar{K}_{kj}\overline{(L_k^*R_j)}_{cb}\\
&=&\sum_jK_{ij}\bar{K}_{kj}(L_i^*L_k)_{ac}\\
&=&N\delta_{ik}(L_i^*L_k)_{ac}\\
&=&NQ\delta_{ik}\delta_{ac}
\end{eqnarray*}

The entries of our matrix being in addition on the unit circle, we are done.
\end{proof}

The above construction is of course something quite abstract, but as a very concrete input for it, we can use the following well-known Fourier analysis construction:

\begin{proposition}
For $q\geq3$ prime, the matrices 
$$\{F_q,DF_q,\ldots,D^{q-1}F_q\}$$
where $F_q$ is the Fourier matrix, and where
$$D=diag\left(1,1,w,w^3,w^6,w^{10},\ldots,w^{\frac{q^2-1}{8}},\ldots,w^{10},w^6,w^3,w\right)$$
with $w=e^{2\pi i/q}$, are such that $\frac{1}{\sqrt{q}}E_i^*E_j$ is complex Hadamard, for any $i\neq j$.
\end{proposition}

\begin{proof}
With by definition $0,1,\ldots,q-1$ as indices for our matrices, as usual in a Fourier analysis context, the formula of the above matrix $D$ is:
$$D_c
=w^{0+1+\ldots+(c-1)}
=w^{\frac{c(c-1)}{2}}$$

Since we have $\frac{1}{\sqrt{q}}E_i^*E_j\in\sqrt{q}U_q$, we just need to check that these matrices have entries belonging to $\mathbb T$, for any $i\neq j$. With $k=j-i$, these entries are given by:
\begin{eqnarray*}
\frac{1}{\sqrt{q}}(E_i^*E_j)_{ab}
&=&\frac{1}{\sqrt{q}}(F_q^*D^kF_q)_{ab}\\
&=&\frac{1}{\sqrt{q}}\sum_cw^{c(b-a)}D_c^k
\end{eqnarray*}

Now observe that with $s=b-a$, we have the following formula:
\begin{eqnarray*}
\left|\sum_cw^{cs}D_c^k\right|^2
&=&\sum_{cd}w^{cs-ds}w^{\frac{c(c-1)}{2}\cdot k-\frac{d(d-1)}{2}\cdot k}\\
&=&\sum_{cd}w^{(c-d)\left(\frac{c+d-1}{2}\cdot k+s\right)}\\
&=&\sum_{de}w^{e\left(\frac{2d+e-1}{2}\cdot k+s\right)}\\
&=&\sum_e\left(w^{\frac{e(e-1)}{2}\cdot k+es}\sum_dw^{edk}\right)\\
&=&\sum_ew^{\frac{e(e-1)}{2}\cdot k+es}\cdot q\delta_{e0}\\
&=&q
\end{eqnarray*}

Thus the entries are on the unit circle, and we are done.
\end{proof}

\index{Legendre symbol}

We recall that the Legendre symbol is defined as follows:
$$\left(\frac{s}{q}\right)=\begin{cases}
0&{\rm if}\ s=0\\
1&{\rm if}\ \exists\,\alpha,s=\alpha^2\\
-1&{\rm if}\not\!\exists\,\alpha,s=\alpha^2
\end{cases}$$

With this convention, we have the following result, following \cite{mwe}:

\begin{proposition}
The following matrices, 
$$G_k=\frac{1}{\sqrt{q}}F_q^*D^kF_q$$
with the matrix $D$ being as above,
$$D=diag\left(w^{\frac{c(c-1)}{2}}\right)$$
and with $k\neq0$ are circulant, their first row vectors $V^k$ being given by
$$V^k_i=\delta_q\left(\frac{k/2}{q}\right)w^{\frac{q^2-1}{8}\cdot k}\cdot w^{-\frac{\frac{i}{k}(\frac{i}{k}-1)}{2}}$$
where $\delta_q=1$ if $q=1(4)$ and $\delta_q=i$ if $q=3(4)$, and with all inverses being taken in $\mathbb Z_q$.
\end{proposition}

\begin{proof}
This is a standard exercice on quadratic Gauss sums. First of all, the matrices $G_k$ in the statement are indeed circulant, their first vectors being given by:
$$V^k_i=\frac{1}{\sqrt{q}}\sum_cw^{\frac{c(c-1)}{2}\cdot k+ic}$$

Let us first compute the square of this quantity. We have:
$$(V_i^k)^2=\frac{1}{q}\sum_{cd}w^{\left[\frac{c(c-1)}{2}+\frac{d(d-1)}{2}\right]k+i(c+d)}$$

The point now is that the sum $S$ on the right, which has $q^2$ terms, decomposes as follows, where $x$ is a certain exponent, depending on $q,i,k$:
$$S=\begin{cases}
(q-1)(1+w+\ldots+w^{q-1})+qw^x&{\rm if}\ q=1(4)\\
(q+1)(1+w+\ldots+w^{q-1})-qw^x&{\rm if}\ q=3(4)
\end{cases}$$

We conclude that we have a formula as follows, where $\delta_q\in\{1,i\}$ is as in the statement, so that $\delta_q^2\in\{1,-1\}$ is given by $\delta_q^2=1$ if $q=1(4)$ and $\delta_q^2=-1$ if $q=3(4)$:
$$(V_i^k)^2=\delta_q^2\,w^x$$

In order to compute now the exponent $x$, we must go back to the above calculation of the sum $S$.  We succesively have:

\smallskip

-- First of all, at $k=1,i=0$ we have $x=\frac{q^2-1}{4}$.

\smallskip

-- By translation we obtain $x=\frac{q^2-1}{4}-i(i-1)$, at $k=1$ and any $i$.

\smallskip

-- By replacing $w\to w^k$ we obtain $x=\frac{q^2-1}{4}\cdot k-\frac{i}{k}(\frac{i}{k}-1)$, at any $k\neq0$ and any $i$. 

\smallskip

Summarizing, we have computed the square of the quantity that we are interested in, the formula being as follows, with $\delta_q$ being as in the statement:
$$(V^k_i)^2=\delta_q^2\cdot w^{\frac{q^2-1}{4}\cdot k}\cdot w^{-\frac{i}{k}(\frac{i}{k}-1)}$$

By extracting now the square root, we obtain a formula as follows:
$$V^k_i=\pm\delta_q\cdot w^{\frac{q^2-1}{8}\cdot k}\cdot w^{-\frac{\frac{i}{k}(\frac{i}{k}-1)}{2}}$$

The computation of the missing sign is non-trivial, but by using the theory of quadratic Gauss sums, and more specifically a result of Gauss, computing precisely this kind of sign, we conclude that we have indeed a Legendre symbol, $\pm=\left(\frac{k/2}{q}\right)$, as claimed. 
\end{proof}

\index{quadratic Gauss sum}

Let us combine now all the above results. We obtain the following statement:

\begin{theorem}
Let $q\geq3$ be prime, consider two subsets 
$$S,T\subset\{0,1,\ldots,q-1\}$$
satisfying the conditions $|S|=|T|$ and $S\cap T=\emptyset$, and write:
$$S=\{s_1,\ldots,s_N\}\quad,\quad 
T=\{t_1,\ldots,t_N\}$$
Then, with the matrix $V$ being as above, the matrix
$$H_{ia,jb}=K_{ij}V^{t_j-s_i}_{b-a}$$
is complex Hadamard, provided that the matrix $K\in M_N(\mathbb C)$ is complex Hadamard.
\end{theorem}

\begin{proof}
This follows indeed by using the general construction in Theorem 8.20, with input coming from Proposition 8.21 and Proposition 8.22.
\end{proof}

\index{McNulty-Weigert matrix}
\index{Tao matrix}
\index{isolated matrix}

As explained by McNulty-Weigert in \cite{mwe}, the above construction covers many interesting examples of Hadamard matrices, previously known from Tadej-\.Zyczkowski \cite{tz1}, \cite{tz2} to be isolated, such as the Tao matrix, which is as follows, with $w=e^{2\pi i/3}$:
$$T_6=\begin{pmatrix}
1&1&1&1&1&1\\ 
1&1&w&w&w^2&w^2\\ 
1&w&1&w^2&w^2&w\\
1&w&w^2&1&w&w^2\\ 
1&w^2&w^2&w&1&w\\ 
1&w^2&w&w^2&w&1
\end{pmatrix}$$

In general, in order to find isolated matrices, the idea from \cite{mwe} is that of starting with an isolated matrix, and then use suitable sets $S,T$. The defect computations are, however, quite difficult. As a concrete statement, however, we have the following conjecture:

\begin{conjecture}
The complex Hadamard matrix constructed in Theorem 8.23 is isolated, provided that:
\begin{enumerate}
\item $K$ is an isolated Fourier matrix, of prime order.

\item $S,T$ consist of consecutive odd numbers, and consecutive even numbers.
\end{enumerate}
\end{conjecture}

This statement is supported by the isolation result for $T_6$, and by several computer simulations from \cite{mwe}. For further details on all this, we refer to \cite{bop}, \cite{mwe}.

\section*{8d. Partial matrices}

As a final topic now, we would like to discuss an extension of a part of our results, from here and from chapter 7, to the case of the partial Hadamard matrices (PHM). The extension, from \cite{bop}, is quite straightforward, but there are a number of subtleties appearing. First of all, we can talk about deformations of PHM, as follows:

\index{PHM}
\index{deformed PHM}

\begin{definition}
Let $H\in X_{M,N}$ be a partial complex Hadamard matrix.
\begin{enumerate}
\item A deformation of $H$ is a smooth function $f:\mathbb T_1\to (X_{M,N})_H$.

\item The deformation is called ``affine'' if $f_{ij}(q)=H_{ij}q^{A_{ij}}$, with $A\in M_{M\times N}(\mathbb R)$.

\item We call ``trivial'' the deformations $f_{ij}(q)=H_{ij}q^{a_i+b_j}$, with $a\in\mathbb R^M,b\in\mathbb R^N$.
\end{enumerate}
\end{definition}

Observe that we have the following equality, where $U_{M,N}\subset M_{M\times N}(\mathbb C)$ is the set of matrices having all rows of norm 1, and pairwise orthogonal:
$$X_{M,N}=M_{M\times N}(\mathbb T)\cap\sqrt{N}U_{M,N}$$

As in the square matrix case, this leads to the following definition:

\index{tangent cone}
\index{affine tangent cone}
\index{enveloping tangent space}

\begin{definition}
Associated to a point $H\in X_{M,N}$ are the enveloping tangent space
$$\widetilde{T}_HX_{M,N}=T_HM_{M\times N}(\mathbb T)\cap T_H\sqrt{N}U_{M,N}$$
as well as the following subcones of this enveloping tangent space:
\begin{enumerate}
\item The tangent cone $T_HX_{M,N}$: the set of tangent vectors to the deformations of $H$.

\item The affine tangent cone $T_H^\circ X_{M,N}$: same as above, using affine deformations only.

\item The trivial tangent cone $T_H^\times X_{M,N}$: as above, using trivial deformations only.
\end{enumerate}
\end{definition}

Observe that $\widetilde{T}_HX_{M,N},T_HX_{M,N}$ are real vector spaces, and that $T_HX_{M,N},T_H^\circ X_{M,N}$ are two-sided cones, in the sense that they satisfy the following condition:
$$\lambda\in\mathbb R,A\in T\implies\lambda A\in T$$

Also, we have inclusions of cones as follows:
$$T_H^\times X_{M,N}\subset T_H^\circ X_{M,N}\subset T_HX_{M,N}\subset\widetilde{T}_HX_{M,N}$$

As in the square matrix case, we can formulate the following definition:

\index{defect}
\index{PHM defect}

\begin{definition}
The defect of a matrix $H\in X_{M,N}$ is the dimension
$$d(H)=\dim(\widetilde{T}_HX_{M,N})$$
of the real vector space $\widetilde{T}_HX_{M,N}$ constructed above.
\end{definition}

The basic remarks and comments regarding the defect from the square matrix case extend then to this setting. In particular, we have the following basic result:

\begin{theorem}
The enveloping tangent space at $H\in X_{M,N}$ is given by
$$\widetilde{T}_HX_{M,N}\simeq\left\{A\in M_{M\times N}(\mathbb R)\Big|\sum_kH_{ik}\bar{H}_{jk}(A_{ik}-A_{jk})=0,\forall i,j\right\}$$
and the defect of $H$ is the dimension of this real vector space.
\end{theorem}

\begin{proof}
In the square matrix case this was done in chapter 7, and the extension of the computations there to the rectangular case is straightforward. First, the manifold $M_{M\times N}(\mathbb T)$ is defined by the following algebraic relations:
$$|H_{ij}|^2=1$$

In terms of real and imaginary parts, $H_{ij}=X_{ij}+iY_{ij}$, we have:
$$d|H_{ij}|^2
=d(X_{ij}^2+Y_{ij}^2)
=2(X_{ij}\dot{X}_{ij}+Y_{ij}\dot{Y}_{ij})$$

Consider now an arbitrary vector $\xi\in T_HM_{M\times N}(\mathbb C)$, written as follows:
$$\xi=\sum_{ij}\alpha_{ij}\dot{X}_{ij}+\beta_{ij}\dot{Y}_{ij}$$

This vector belongs then to $T_HM_{M\times N}(\mathbb T)$ if and only if we have:
$$<\xi,d|H_{ij}|^2>=0$$

We therefore obtain the following formula, for the tangent cone:
$$T_HM_{M\times N}(\mathbb T)=\left\{\sum_{ij}A_{ij}(Y_{ij}\dot{X}_{ij}-X_{ij}\dot{Y}_{ij})\Big|A_{ij}\in\mathbb R\right\}$$

We also know that the manifold $\sqrt{N}U_{M,N}$ is defined by the following algebraic relations, where $H_1,\ldots,H_N$ are the rows of $H$:
$$<H_i,H_j>=N\delta_{ij}$$
 
The relations $<H_i,H_i>=N$ being automatic for the matrices $H\in M_{M\times N}(\mathbb T)$, if for $i\neq j$ we let $L_{ij}=<H_i,H_j>$, then we have:
$$\widetilde{T}_HC_N=\left\{\xi\in T_HM_N(\mathbb T)\Big|<\xi,\dot{L}_{ij}>=0,\,\forall i\neq j\right\}$$

On the other hand, differentiating the formula of $L_{ij}$ gives:
$$\dot{L}_{ij}=\sum_k(X_{ik}+iY_{ik})(\dot{X}_{jk}-i\dot{Y}_{jk})+(X_{jk}-iY_{jk})(\dot{X}_{ik}+i\dot{Y}_{ik})$$

Now if we pick a vector $\xi\in T_HM_{M\times N}(\mathbb T)$, written as above in terms of $A\in M_{M\times N}(\mathbb R)$, we obtain the following formula:
$$<\xi,\dot{L}_{ij}>=i\sum_k\bar{H}_{ik}H_{jk}(A_{ik}-A_{jk})$$

Thus we have reached to the description of $\widetilde{T}_HX_{M,N}$ in the statement. 
\end{proof}

Summarizing, the extension of the basic defect theory, from the square matrix case to the rectangular matrix case, appears to be quite straightforward. By using the above defect equations, most of the general comments and remarks from chapter 7 regarding the square matrix case extend to the rectangular matrix case. See \cite{bop}. At the level of non-trivial results now, we first have the following statement:

\begin{theorem}
Let $H\in X_{M,N}$, and pick a square matrix
$$K\in\sqrt{N}U_N$$
extending $H$. We have then the following formula,
$$\widetilde{T}_HX_{M,N}\simeq\left\{E=(X\ Y)\in M_{M\times N}(\mathbb C)\Big|X=X^*,(EK)_{ij}\bar{H}_{ij}\in\mathbb R,\forall i,j\right\}$$
with the correspondence $A\to E$ being constructed as follows:
$$E_{ij}=\sum_kH_{ik}\bar{K}_{jk}A_{ik}\quad,\quad 
A_{ij}=(EK)_{ij}\bar{H}_{ij}$$
\end{theorem}

\begin{proof}
Let us set indeed $R_{ij}=A_{ij}H_{ij}$ and $E=RK^*$. The correspondence $A\to R\to E$ is then bijective, and we have the following formula:
$$E_{ij}=\sum_kH_{ik}\bar{K}_{jk}A_{ik}$$

With these changes, the system of equations in Theorem 8.28 becomes $E_{ij}=\bar{E}_{ji}$ for any $i,j$ with $j\leq M$. But this shows that we must have $E=(X\ Y)$ with $X=X^*$, and the condition $A_{ij}\in\mathbb R$ corresponds to the condition $(EK)_{ij}\bar{H}_{ij}\in\mathbb R$, as claimed.
\end{proof}

As an illustration, in the real case we obtain the following result:

\index{real PHM}

\begin{theorem}
For an Hadamard matrix $H\in M_{M\times N}(\pm1)$ we have
$$\widetilde{T}_HX_{M,N}\simeq M_M(\mathbb R)^{symm}\oplus M_{M\times(N-M)}(\mathbb R)$$
and so the defect is given by 
$$d(H)=\frac{N(N+1)}{2}+M(N-M)$$
independently of the precise value of $H$.
\end{theorem}

\begin{proof}
We use Theorem 8.29. Since $H$ is now real we can pick $K\in\sqrt{N}U_N$ extending it to be real too, and with nonzero entries, so the last condition appearing there, namely $(EK)_{ij}\bar{H}_{ij}\in\mathbb R$, simply tells us that $E$ must be real. Thus we have:
$$\widetilde{T}_HX_{M,N}\simeq\left\{E=(X\ Y)\in M_{M\times N}(\mathbb R)\Big|X=X^*\right\}$$

But this is the formula in the statement, and we are done.
\end{proof}

A matrix $H\in X_{M,N}$ cannot be isolated, simply because the space of its Hadamard equivalents provides a copy $\mathbb T^{MN}\subset X_{M,N}$, passing through $H$. However, if we restrict the attention to the matrices which are dephased, the notion of isolation makes sense:

\begin{proposition}
The defect $d(H)=\dim(\widetilde{T}_HX_{M,N})$ satisfies
$$d(H)\geq M+N-1$$
and if $d(H)=M+N-1$ then $H$ is isolated inside the dephased quotient $X_{M,N}\to Z_{M,N}$.
\end{proposition}

\begin{proof}
Once again, the known results in the square case extend:

\medskip

(1) We have indeed $\dim(T_H^\times X_{M,N})=M+N-1$, and since the tangent vectors to these trivial deformations belong to $\widetilde{T}_HX_{M,N}$, this gives the first assertion.

\medskip

(2) Since $d(H)=M+N-1$, the inclusions $T_H^\times X_{M,N}\subset T_HX_{M,N}\subset\widetilde{T}_HX_{M,N}$ must be equalities, and from $T_HX_{M,N}=T_H^\times X_{M,N}$ we obtain the result.
\end{proof}

Finally, still at the theoretical level, we have the following conjecture:

\index{isolated PHM}

\begin{conjecture}
An isolated partial Hadamard matrix $H\in Z_{M,N}$ must have minimal defect, namely $d(H)=M+N-1$.
\end{conjecture}

In other words, the conjecture is that if $H\in Z_{M,N}$ has only trivial first order deformations, then it has only trivial deformations at any order, including at $\infty$. In the square matrix case this statement comes with solid evidence, all known examples of complex Hadamard matrices $H\in X_N$ having non-minimal defect being known to admit one-parameter deformations. For more on this subject, see \cite{bop}, \cite{tz1}, \cite{tz2}.

\bigskip

Let us discuss now some examples of isolated partial Hadamard matrices, and provide some evidence for Conjecture 8.32. We are interested in the following matrices:

\index{truncated Fourier matrix}

\begin{definition}
The truncated Fourier matrix $F_{S,G}$, with $G$ being a finite abelian group, and with $S\subset G$ being a subset, is constructed as follows:
\begin{enumerate}
\item Given $N\in\mathbb N$, we set $F_N=(w^{ij})_{ij}$, where $w=e^{2\pi i/N}$.

\item Assuming $G=\mathbb Z_{N_1}\times\ldots\times\mathbb Z_{N_s}$, we set $F_G=F_{N_1}\otimes\ldots\otimes F_{N_s}$.

\item We let $F_{S,G}$ be the submatrix of $F_G$ having $S\subset G$ as row index set. 
\end{enumerate}
\end{definition}

Observe that $F_N$ is the Fourier matrix of the cyclic group $\mathbb Z_N$. More generally, $F_G$ is the Fourier matrix of the finite abelian group $G$. Observe also that $F_{G,G}=F_G$. We can compute the defect of $F_{S,G}$ by using Theorem 8.28, and we obtain:

\begin{theorem}
For a truncated Fourier matrix $F=F_{S,G}$ we have the formula
$$\widetilde{T}_FX_{M,N}=\left\{A\in M_{M\times N}(\mathbb R)\Big|P=AF^t\ {\rm satisfies}\ P_{ij}=P_{i+j,j}=\bar{P}_{i,-j},\forall i,j\right\}$$
where $M=|S|,N=|G|$, and with all the indices being regarded as group elements. 
\end{theorem}

\begin{proof}
We use Theorem 8.28. The defect equations there are as follows:
$$\sum_kF_{ik}\bar{F}_{jk}(A_{ik}-A_{jk})=0$$

For $F=F_{S,G}$ we have the following formula:
$$F_{ik}\bar{F}_{jk}=(F^t)_{k,i-j}$$

We therefore obtain the following formula:
$$\widetilde{T}_FX_{M,N}=\left\{A\in M_{M\times N}(\mathbb R)\Big|(AF^t)_{i,i-j}=(AF^t)_{j,i-j},\forall i,j\right\}$$

Now observe that for an arbitrary matrix $P\in M_M(\mathbb C)$, we have:
\begin{eqnarray*}
P_{i,i-j}=P_{j,i-j},\forall i,j
&\iff&P_{i+j,i}=P_{ji},\forall i,j\\
&\iff&P_{i+j,j}=P_{ij},\forall i,j
\end{eqnarray*}

We therefore conclude that we have the following equality:
$$\widetilde{T}_FX_{M,N}=\left\{A\in M_{M\times N}(\mathbb R)\Big|
P=AF^t\ {\rm satisfies}\ P_{ij}=P_{i+j,j},\forall i,j\right\}$$

Now observe that with $A\in M_{M\times N}(\mathbb R)$ and $P=AF^t\in M_M(\mathbb C)$ as above, we have:
\begin{eqnarray*}
\bar{P}_{ij}
&=&\sum_kA_{ik}(F^*)_{kj}\\
&=&\sum_kA_{ik}(F^t)_{k,-j}\\
&=&P_{i,-j}
\end{eqnarray*}

Thus, we obtain the formula in the statement, and we are done. 
\end{proof}

Let us try to find some explicit examples of isolated matrices, of truncated Fourier type. For this purpose, we can use the following improved version of Theorem 8.34:

\begin{theorem}
The defect of $F=F_{S,G}$ is the number 
$$d(F)=\dim(K)+\dim(I)$$
where $K,I$ are the following linear spaces,
\begin{eqnarray*}
K&=&\left\{A\in M_{M\times N}(\mathbb R)\Big|AF^t=0\right\}\\
I&=&\left\{P\in L_M\Big|\exists A\in M_{M\times N}(\mathbb R),P=AF^t\right\}
\end{eqnarray*}
with $L_M$ being the following linear space,
$$L_M=\left\{P\in M_M(\mathbb C)\Big|P_{ij}=P_{i+j,j}=\bar{P}_{i,-j},\forall i,j\right\}$$
with all the indices belonging by definition to the group $G$.
\end{theorem}

\begin{proof}
We use the general formula in Theorem 8.34. With the notations there, and with the linear space $L_M$ being as above, we have a linear map as follows:
$$\Phi:\widetilde{T}_FX_{M,N}\to L_M\quad,\quad 
\Phi(A)=AF^t$$

By using this map, we obtain the following equality:
$$\dim(\widetilde{T}_FX_{M,N})=\dim(\ker\Phi)+\dim({\rm Im}\,\Phi)$$

Now since the spaces on the right are precisely those in the statement, we have:
$$\ker\Phi=K\quad,\quad 
{\rm Im}\, \Phi=I$$

Thus by applying Theorem 8.34 we obtain the result.
\end{proof}

In order to look now for isolated matrices, the first remark is that since a deformation of $F_G$ will produce a deformation of $F_{S,G}$ too, we must  restrict the attention to the case where $G=\mathbb Z_p$, with $p$ prime. And here, we have the following conjecture:

\index{truncated Fourier matrix}

\begin{conjecture}
There exists a constant $\varepsilon>0$ such that $F_{S,p}$ is isolated, for any $p$ prime, once $S\subset\mathbb Z_p$ satisfies $|S|\geq(1-\varepsilon)p$.
\end{conjecture}

In principle this conjecture can be approached by using the formula in Theorem 8.35, and we have for instance evidence towards the fact that $F_{p-1,p}$ should be always isolated, that $F_{p-2,p}$ should be isolated too, provided that $p$ is big enough, and so on. However, finding a number $\varepsilon>0$ as above looks like a quite difficult question. See \cite{bop}.

\section*{8e. Exercises} 

There has been a lot of material in this chapter, regarding many types of Hadamard matrices. As a first exercise, in connection with the tensor products, we have:

\begin{exercise}
Write down a list of examples where we have equality case,
$$d(H\otimes K)=d(H)d(K)$$
in the general inequality $d(H\otimes K)\geq d(H)d(K)$ established above.
\end{exercise}

To be more precise, there is some work to be done in the Fourier matrix case, and passed that, the problem is to see which other of our defect computations can help.

\begin{exercise}
Prove that the Tao matrix, namely
$$T_6=\begin{pmatrix}
1&1&1&1&1&1\\ 
1&1&w&w&w^2&w^2\\ 
1&w&1&w^2&w^2&w\\
1&w&w^2&1&w&w^2\\ 
1&w^2&w^2&w&1&w\\ 
1&w^2&w&w^2&w&1
\end{pmatrix}$$
with $w=e^{2\pi i/3}$, is indeed a McNulty-Weigert matrix.
\end{exercise}

Observe in particular that a solution to this exercise would provide a solution to one of our previous exercises, which was probably difficult, asking for an explicit formula for $T_6$, with the matrix entries $(T_6)_{ij}$ expressed as explicit functions of the indices $i,j$.

\begin{exercise}
Compute the defect of the truncated Fourier matrices, at small values of the truncation parameter.
\end{exercise}

The problem here is that of applying the various results established above.

\part{Analytic aspects}

\ \vskip50mm

\begin{center}
{\em Look what they've done to my song, ma

It was the only thing I could do half right

And it's turning out all wrong, ma

Look what they've done to my song}
\end{center}

\chapter{Circulant matrices}

\section*{9a. Cyclic roots}

After some 200 pages of analysis, time to do some analysis. In this third part of the present book we discuss a number of more specialized analytic topics, in relation with the following questions, regarding the complex Hadamard matrices:

\bigskip

-- Circulant Hadamard matrices. We will discuss here Bj\"orck's cyclic root formalism \cite{bjo}, the Butson matrix analogues of the CHC, the Haagerup counting result in \cite{ha2}, and, following \cite{bs1}, an analytic approach to the CHC, using the 4-norm.

\bigskip

-- Bistochastic Hadamard matrices. These matrices, covering all the circulant ones, and very interesting objects, due to a result of Idel-Wolf \cite{iwo}, stating that any unitary matrix, and so any complex Hadamard matrix, can be put in bistochastic form.

\bigskip

-- The glow of Hadamard matrices. This is another interesting theme, related on one hand to the glow computations from the real case, that we did in chapter 1, motivated by the Gale-Berlekamp game, and on the other hand, by the Idel-Wolf theorem.

\bigskip

-- Almost Hadamard matrices. The study here, from \cite{bn3}, initially paralleling the study from the real case, from chapter 3, leads to an unexpected and potentially far-reaching conjecture, stating that ``being complex Hadamard is a local property''.

\bigskip

All in all, many things to be discussed, and we should mention too that all this will be rather research-grade material, quite recent, and with more conjectures than theorems, and with all this waiting for some enthusiastic young people. Like you.

\bigskip

Getting started now, we will first discuss an important class of complex Hadamard matrices, namely the circulant ones. There has been a lot of work here, starting with the Circulant Hadamard Conjecture (CHC) in the real case, and with many results in the complex case as well. We will present here the main techniques in dealing with such matrices. It is convenient to introduce the circulant matrices as follows:

\index{circulant matrix}

\begin{definition}
A complex matrix $H\in M_N(\mathbb C)$ is called circulant when we have
$$H_{ij}=\gamma_{j-i}$$
for some $\gamma\in\mathbb C^N$, with  the matrix indices $i,j\in \{0,1,\ldots,N-1\}$ taken modulo $N$. 
\end{definition}

Here the index convention is quite standard, as for the Fourier matrices $F_N$, and with this coming from some Fourier analysis considerations, that we will get into later on. In practice now, the fact that a matrix is circulant means that it has the following pattern, with the entries in the first row ``circulating'' downwards and to the right:
$$H=\begin{pmatrix}
a&b&c&d\\
d&a&b&c\\
c&d&a&b\\
b&c&d&a
\end{pmatrix}$$

As a basic example of a circulant Hadamard matrix, in the real case, we have the matrix $K_4$. The circulant Hadamard conjecture states that this matrix is, up to equivalence, the only circulant Hadamard matrix $H\in M_N(\pm1)$, regardless of the value of $N\in\mathbb N$:

\index{CHC}
\index{Circulant Hadamard Conjecture}

\begin{conjecture}[Circulant Hadamard Conjecture (CHC)]
The only circulant real Hadamard matrices $H\in M_N(\pm1)$ are the matrix
$$K_4=\begin{pmatrix}-1&1&1&1\\1&-1&1&1\\1&1&-1&1\\1&1&1&-1\end{pmatrix}$$
and its Hadamard conjugates, and this regardless of the value of $N\in\mathbb N$.
\end{conjecture}

As explained in chapter 1, this conjecture is something of different nature from the Hadamard Conjecture (HC). Indeed, while the HC might look like something simple, at the first glance, working a bit on it quickly reveals that this is certainly something quite complicated, or even worse, that this might be one of these ``black holes'' in the mathematical landscape, including too the Riemann Hypothesis, the Jacobian Conjecture, the Collatz Problem and so on, all questions having little to do with modern mathematics as we know it, since Newton and others, and better to be avoided.

\bigskip

Regarding the CHC, however, it is quite unclear where the difficulty comes from. Indeed, if we denote by $S\subset\{1,\ldots,N\}$ the set of positions of the $-1$ entries of the first row vector $\gamma\in(\pm1)^N$, the Hadamard matrix condition reads, for any $k\neq0$:
$$|S\cap(S+k)|=|S|-N/4$$

Thus, the CHC simply states that at $N\neq 4$, such a set $S$ cannot exist. Let us record here this latter statement, originally due to Ryser \cite{rys}:

\index{Ryser Conjecture}

\begin{conjecture}[Ryser Conjecture]
Given an integer $N>4$, there is no set 
$$S\subset\{1,\ldots,N\}$$
satisfying $|S\cap(S+k)|=|S|-N/4$ for any $k\neq 0$, taken modulo $N$.
\end{conjecture}

And prove this if you can. This question is 60 years old, and many competent people have looked at it, with basically 0 serious advances. So, most likely, what we have here is the same type of annoying question as the HC, Riemann, Collatz and so on.

\bigskip

Erd\H{o}s famously said about Collatz that ``mathematics is not ready for such things''. But, will it ever be ready? Probably not. Never. It is always good to remember here that modern mathematics as we know it was developed by Newton and others, with inspiration from classical mechanics. And so, want it or not, mathematics as we know it ``is'' classical mechanics. And this might explain why the HC, CHC, Riemann, Collatz and so on are so inaccessible, these are probably simply questions which are orthogonal to classical mechanics, and so are orthogonal to mathematics as we know it too.

\bigskip

You might say then, why not trying mathematics inspired from some other physics, like quantum mechanics. Well, the problem is that quantum mechanics, or at least quantum mechanics as we know it, is in fact not that far from classical mechanics. Same types of beasts, like functions, derivatives, integrals and so on, all good old stuff going back to Newton, doing most of the mathematics that we know, in the quantum world.

\bigskip

But then you would say why not sending to trash all modern mathematics, and developing some new, original mathematics, especially tailored for problems like the HC, CHC, Riemann, Collatz and so on. Well, people have tried, for instance with design theory for the HC, CHC, and this does not work either. And why? No one really knows the answer here, but this is probably because there is no physics that you can rely upon, and intuition in general, for making that original mathematics of yours strong and reliable.

\bigskip

Looks like we are in a kind of vicious circle, with all these questions. Math needs physics, and so, want it or not, the physics surrounding us ultimately dictates what's doable and what's not, mathematically speaking. As a conjecture, in some alien world where the physics is different, the HC, CHC, Riemann, Collatz and so on might be all trivial. But that little green men who know how to solve all these questions might, on the other hand, have things like partial integration as longstanding, open problems.

\bigskip

And let us end this discussion with a famous quote by Dirac, ``shut up and compute''. This is what he used to say to students asking too many questions about quantum mechanics. Computation is our only tool, so let's compute some more. After all, there is still a chance that the HC, CHC might be related to mechanics. And so, be doable.

\bigskip

Back to work now, we will in fact not start with computations for the CHC, which looks quite scary. Our first purpose will be that of showing that the CHC dissapears in the complex case, where we have examples at any $N\in\mathbb N$. As a first result, we have:

\begin{proposition}
The following are circulant and symmetric Hadamard matrices,
$$F_2'=\begin{pmatrix}i&1\\1&i\end{pmatrix}\qquad,\qquad
F_3'=\begin{pmatrix}w&1&1\\1&w&1\\1&1&w\end{pmatrix}$$
$$F_4''=\begin{pmatrix}-1&\nu&1&\nu\\\nu&-1&\nu&1\\1&\nu&-1&\nu\\ \nu&1&\nu&-1\end{pmatrix}$$
where $w=e^{2\pi i/3},\nu=e^{\pi i/4}$, equivalent to the Fourier matrices $F_2,F_3,F_4$.
\end{proposition}

\begin{proof}
The orthogonality between rows being clear, we have here complex Hadamard matrices. The fact that we have an equivalence $F_2\sim F_2'$ follows from:
$$\begin{pmatrix}1&1\\1&-1\end{pmatrix}
\sim\begin{pmatrix}i&i\\1&-1\end{pmatrix}
\sim\begin{pmatrix}i&1\\1&i\end{pmatrix}$$

At $N=3$ now, the equivalence $F_3\sim F_3'$ can be constructed as follows:
$$\begin{pmatrix}1&1&1\\1&w&w^2\\1&w^2&w\end{pmatrix}
\sim\begin{pmatrix}1&1&w\\1&w&1\\w&1&1\end{pmatrix}
\sim\begin{pmatrix}w&1&1\\1&w&1\\1&1&w\end{pmatrix}$$

As for the case $N=4$, here the equivalence $F_4\sim F_4''$ can be constructed as follows, where we use the logarithmic notation $[k]_s=e^{2\pi ki/s}$, with respect to $s=8$:
$$\begin{bmatrix}0&0&0&0\\0&2&4&6\\0&4&0&4\\0&6&4&2\end{bmatrix}_8
\sim\begin{bmatrix}0&1&4&1\\1&4&1&0\\4&1&0&1\\1&0&1&4\end{bmatrix}_8
\sim\begin{bmatrix}4&1&0&1\\1&4&1&0\\0&1&4&1\\1&0&1&4\end{bmatrix}_8
$$

Thus, the Fourier matrices $F_2,F_3,F_4$ can be put indeed in circulant form.
\end{proof}

We will explain later the reasons for denoting the above matrix $F_4''$, instead of $F_4'$, the idea being that $F_4'$, not introduced yet, will be a matrix belonging to a certain series. Getting back now to the real circulant matrix $K_4$, this is equivalent to the Fourier matrix $F_G=F_2\otimes F_2$ of the Klein group $G=\mathbb Z_2\times\mathbb Z_2$, as shown by:
$$\begin{pmatrix}-1&1&1&1\\1&-1&1&1\\1&1&-1&1\\1&1&1&-1\end{pmatrix}
\sim\begin{pmatrix}
1&1&1&-1\\
1&-1&1&1\\
1&1&-1&1\\
-1&1&1&1
\end{pmatrix}
\sim\begin{pmatrix}
1&1&1&1\\
1&-1&1&-1\\
1&1&-1&-1\\
1&-1&-1&1
\end{pmatrix}$$

In fact, we have the following construction of circulant and symmetric Hadamard matrices at $N=4$, which involves an extra parameter $q\in\mathbb T$:

\begin{proposition}
The following circulant and symmetric matrix is Hadamard,
$$K_4^q=\begin{pmatrix}-1&q&1&q\\q&-1&q&1\\1&q&-1&q\\q&1&q&-1\end{pmatrix}$$
for any $q\in\mathbb T$. At $q=1,e^{\pi i/4}$ recover respectively the matrices $K_4,F_4''$.
\end{proposition}

\begin{proof}
The rows of the above matrix are pairwise orthogonal for any $q\in\mathbb C$, and so at $q\in\mathbb T$ we obtain an Hadamard matrix. As for the last assertion, this is clear.
\end{proof}

As a first conclusion, coming from the above considerations, we have:

\index{circulant form}
\index{circulant and symmetric form}

\begin{theorem}
The complex Hadamard matrices of order $N=2,3,4,5$, namely 
$$F_2,F_3,F_4^s,F_5$$
can be put, up to equivalence, in circulant and symmetric form.
\end{theorem}

\begin{proof}
As explained in chapter 5, the complex Hadamard matrices at $N=2,3,4,5$ are, up to equivalence, those in the statement, with the classification being something elementary at $N=2,3,4$, and with the $N=5$ result being due to Haagerup \cite{ha1}. 

\medskip

(1) At $N=2,3$ the problem is solved by Proposition 9.4. 

\medskip

(2) At $N=4$ now, our claim is that, with $s=q^{-2}$, we have:
$$K_4^q\sim F_4^s$$

Indeed, by multiplying the rows and columns of $K_4^q$ by suitable scalars, we have:
$$K_4^q
=\begin{pmatrix}-1&q&1&q\\q&-1&q&1\\1&q&-1&q\\q&1&q&-1\end{pmatrix}\sim\begin{pmatrix}
1&-q&-1&-q\\
1&-\bar{q}&1&\bar{q}\\
1&q&-1&q\\
1&\bar{q}&1&-\bar{q}\end{pmatrix}
\sim\begin{pmatrix}
1&1&1&1\\
1&s&-1&-s\\
1&-1&1&-1\\
1&-s&-1&s\end{pmatrix}$$

On the other hand, by permuting the second and third rows of $F_4^s$, we obtain:
$$F_4^s
=\begin{pmatrix}
1&1&1&1\\
1&-1&1&-1\\
1&s&-1&-s\\ 
1&-s&-1&s
\end{pmatrix}
\sim
\begin{pmatrix}
1&1&1&1\\
1&s&-1&-s\\
1&-1&1&-1\\
1&-s&-1&s\end{pmatrix}$$

Thus these matrices are equivalent, and the result follows from Proposition 9.5. 

\medskip

(3) At $N=5$ now, the matrix that we are looking for is as follows, with $w=e^{2\pi i/5}$:
$$F_5'=\begin{pmatrix}
w^2&1&w^4&w^4&1\\ 
1&w^2&1&w^4&w^4\\ 
w^4&1&w^2&1&w^4\\ 
w^4&w^4&1&w^2&1\\ 
1&w^4&w^4&1&w^2
\end{pmatrix}$$

It is indeed clear that this matrix is circulant, symmetric, and complex Hadamard, and the fact that we have $F_5\sim F_5'$ follows either directly, or by using Haagerup \cite{ha1}.
\end{proof}

Summarizing, many interesting examples of complex Hadamard matrices are circulant. This is in stark contrast with the real case, where the CHC, discussed above, states that the only circulant real matrices should be those appearing at $N=4$. Let us prove now, as a generalization of all this, that any Fourier matrix $F_N$ can be put in circulant and symmetric form. We use Bj\"orck's cyclic root formalism \cite{bjo}, which is as follows:

\index{cyclic root}
\index{Bj\"orck cyclic root}

\begin{theorem}
Assume that a matrix $H\in M_N(\mathbb T)$ is circulant, $H_{ij}=\gamma_{j-i}$. Then $H$ is a complex Hadamard matrix precisely when the vector 
$$z=(z_0,z_1,\ldots,z_{N-1})$$
given by $z_i=\gamma_i/\gamma_{i-1}$ satisfies the following equations:
\begin{eqnarray*}
z_0+z_1+\ldots+z_{N-1}&=&0\\
z_0z_1+z_1z_2+\ldots+z_{N-1}z_0&=&0\\
\vdots\\
z_0z_1\ldots z_{N-2}+\ldots+z_{N-1}z_0\ldots z_{N-3}&=&0\\
z_0z_1\ldots z_{N-1}&=&1
\end{eqnarray*}
If so is the case, we say that $z=(z_0,\ldots,z_{N-1})$ is a cyclic $N$-root.
\end{theorem}

\begin{proof}
Assume that a matrix of type $H\in M_N(\mathbb T)$ is circulant, $H_{ij}=\gamma_{j-i}$, and set $z_i=\gamma_i/\gamma_{i-1}$, as in the statement. Observe that we have:
$$z_0z_1\ldots z_{N-1}=1$$

Up to a multiplication by a scalar $w\in\mathbb T$, our matrix is then as follows:
$$H=\begin{pmatrix}
z_0&z_0z_1&z_0z_1z_2&\ldots\ldots&z_0\ldots z_{N-1}\\
z_0\ldots z_{N-1}&z_0&z_0z_1&\ldots\ldots&z_0\ldots z_{N-2}\\
z_0\ldots z_{N-2}&z_0\ldots z_{N-1}&z_0&\ldots\ldots&z_0\ldots z_{N-3}\\
\vdots&\vdots&\vdots&&\vdots\\
z_0z_1&z_0z_1z_2&z_0z_1z_2z_3&\ldots\ldots&z_0
\end{pmatrix}$$

Since this matrix is circulant, it is Hadamard precisely when the first row $R_0$ is orthogonal to the other rows $R_1,\ldots,R_{N-1}$. And the equations here are as follows:

\medskip

$(R_0\perp R_1)$. Here the orthogonality condition is as follows:
$$\overline{z_1\ldots z_{N-1}}+z_1+z_2+\ldots+z_{N-1}=0$$

Now by using $z_0z_1\ldots z_{N-1}=1$, this is the 1st equation for cyclic roots, namely:
$$z_0+z_1+z_2+\ldots+z_{N-1}=0$$

$(R_0\perp R_2)$. Here the orthogonality condition is as follows:
$$\overline{z_1\ldots z_{N-2}}+\overline{z_2\ldots z_{N-1}}+z_1z_2+\ldots+z_{N-2}z_{N-1}=0$$

By using again $z_0z_1\ldots z_{N-1}=1$, this is the 2nd equation for cyclic roots, namely:
$$z_{N-1}z_0+z_0z_1+z_1z_2+\ldots+z_{N-2}z_{N-1}=0$$

\ \ \ \ $\vdots$

\medskip

$(R_0\perp R_{N-1})$. Here the orthogonality condition is as follows:
$$\overline{z}_1+\overline{z}_2+\overline{z}_3+\ldots+z_1\ldots z_{N-1}=0$$

And again by using $z_0z_1\ldots z_{N-1}=1$, this is the last equation for cyclic roots, namely:
$$z_2\ldots z_{N-1}z_0+z_3\ldots z_{N-1}z_0z_1+z_4\ldots z_{N-1}z_0z_1z_2+\ldots+z_1\ldots z_{N-1}=0$$

Thus, we are led to the conclusion in the statement.
\end{proof}

The above manipulation might look like something very simple, but in practice this considerably simplifies things, and leads to non-trivial results. Technically speaking now, observe that, up to a multiplication by a scalar $w\in\mathbb T$, the first row vector $\gamma=(\gamma_0,\ldots,\gamma_{N-1})$ of the matrix $H\in M_N(\mathbb T)$ constructed in Theorem 9.7 is as follows:
$$\gamma=(z_0,z_0z_1,z_0z_1z_2,\ldots\ldots,z_0z_1\ldots z_{N-1})$$

We will use this observation several times, in what follows. Now back to the Fourier matrices, we have the following result:

\index{Fourier matrix}
\index{circulant form}
\index{circulant and symmetric form}

\begin{theorem}
Given $N\in\mathbb N$, construct the following complex numbers:
$$\nu=e^{\pi i/N}\quad,\quad 
q=\nu^{N-1}\quad,\quad 
w=\nu^2$$
We have then a cyclic $N$-root as follows, in the above sense,
$$(q,qw,qw^2,\ldots,qw^{N-1})$$
and the corresponding complex Hadamard matrix $F_N'$ is circulant and symmetric, and equivalent to the Fourier matrix $F_N$.
\end{theorem}

\begin{proof}
Given two numbers $q,w\in\mathbb T$, let us find out when $(q,qw,qw^2,\ldots,qw^{N-1})$ is a cyclic root. We have two conditions to be verified, as follows:

\medskip

(1) In order for the $=0$ equations in Theorem 9.7 to be satisfied, the value of $q$ is irrelevant, and $w$ must be a primitive $N$-root of unity. 

\medskip

(2) As for the $=1$ equation in Theorem 9.7, this states in our case that we must have: 
$$q^Nw^{\frac{N(N-1)}{2}}=1$$

We conclude from this that we must have:
$$q^N=(-1)^{N-1}$$

Thus, with the values of $q,w\in\mathbb T$ in the statement, we have indeed a cyclic $N$-root. Now construct $H_{ij}=\gamma_{j-i}$ as in Theorem 9.7. We have:
\begin{eqnarray*}
\gamma_k=\gamma_{-k}
&\iff&q^{k+1}w^{\frac{k(k+1)}{2}}=q^{-k+1}w^{\frac{k(k-1)}{2}}\\
&\iff&q^{2k}w^k=1\\
&\iff&q^2=w^{-1}
\end{eqnarray*}

But this latter condition holds indeed, because we have:
$$q^2
=\nu^{2N-2}
=\nu^{-2}
=w^{-1}$$

We conclude that our circulant matrix $H$ is symmetric as well, as claimed. It remains to construct an equivalence $H\sim F_N$. In order to do this, observe that, due to our conventions $q=\nu^{N-1},w=\nu^2$, the first row vector of $H$ is given by:
\begin{eqnarray*}
\gamma_k
&=&q^{k+1}w^{\frac{k(k+1)}{2}}\\
&=&\nu^{(N-1)(k+1)}\nu^{k(k+1)}\\
&=&\nu^{(N+k-1)(k+1)}
\end{eqnarray*}

Thus, the entries of $H$ are given by the following formula:
\begin{eqnarray*}
H_{-i,j}
&=&H_{0,i+j}\\
&=&\nu^{(N+i+j-1)(i+j+1)}\\
&=&\nu^{i^2+j^2+2ij+Ni+Nj+N-1}\\
&=&\nu^{N-1}\cdot\nu^{i^2+Ni}\cdot\nu^{j^2+Nj}\cdot\nu^{2ij}
\end{eqnarray*}

With this formula in hand, we can now finish the proof. Indeed, this shows that the matrix $H=(H_{ij})$ is equivalent to the following matrix:
$$H'=(H_{-i,j})$$

Now regarding this latter matrix $H'$, observe that in the above formula, the factors $\nu^{N-1}$, $\nu^{i^2+Ni}$, $\nu^{j^2+Nj}$ correspond respectively to a global multiplication by a scalar, and to row and column multiplications by scalars. Thus this matrix $H'$ is equivalent to the matrix $H''$ obtained from it by deleting these factors. But this latter matrix is:
$$H''_{ij}=\nu^{2ij}\quad,\quad \nu=e^{\pi i/N}$$

Since this is precisely the Fourier matrix $F_N$, we are done.
\end{proof}

As an illustration, let us work out the cases $N=2,3,4,5$. We have here:

\begin{proposition}
The matrices $F_N'$ are as follows:
\begin{enumerate}
\item At $N=2,3$ we obtain the old matrices $F_2',F_3'$.

\item At $N=4$ we obtain the following matrix, with $\nu=e^{\pi i/4}$:
$$F_4'=\begin{pmatrix}
\nu^3&1&\nu^7&1\\
1&\nu^3&1&\nu^7\\
\nu^7&1&\nu^3&1\\
1&\nu^7&1&\nu^3
\end{pmatrix}$$

\item At $N=5$ we obtain the old matrix $F_5'$.
\end{enumerate}
\end{proposition}

\begin{proof}
With notations from Theorem 9.8, the proof goes as follows:

\medskip

(1) At $N=2$ we have $\nu=i,q=i,w=-1$, so the cyclic root is $(i,-i)$. The first row vector is $(i,1)$, and we obtain indeed the old matrix $F_2'$. 

\medskip 

At $N=3$ we have $\nu=e^{\pi i/3}$ and $q=w=\nu^2=e^{2\pi i/3}$, the cyclic root is $(w,w^2,1)$. The first row vector is $(w,1,1)$, and we obtain indeed the old matrix $F_3'$.

\medskip 

(2) At $N=4$ we have $\nu=e^{\pi i/4}$ and $q=\nu^3,w=\nu^2$, the cyclic root is $(\nu^3,\nu^5,\nu^7,\nu)$. The first row vector is $(\nu^3,1,\nu^7,1)$, and we obtain the matrix in the statement.

\medskip 

(3) At $N=5$ we have $\nu=e^{\pi i/5}$ and $q=\nu^4=w^2$, with $w=\nu^2=e^{2\pi i/5}$, and the cyclic root is therefore $(w^2,w^3,w^4,1,w)$. The first row vector is $(w^2,1,w^4,w^4,1)$, and we obtain in this way the old matrix $F_5'$, as claimed.
\end{proof}

Regarding the above matrix $F_4'$, observe that this is equivalent to the matrix $F_4''$ from Proposition 9.4, with the equivalence $F_4'\sim F_4''$ being obtained by multiplying everything by $\nu=e^{\pi i/4}$.  While both these matrices are circulant and symmetric, and of course equivalent to $F_4$, one of them, namely $F_4'$, is ``better'' than the other, because the corresponding cyclic root comes from a progression. This is the reason for our notations $F_4',F_4''$.

\bigskip

Let us discuss now the case of the generalized Fourier matrices $F_G$. In this context, the assumption of being circulant is somewhat unnatural, because this comes from a $\mathbb Z_N$ symmetry, and the underlying group is no longer $\mathbb Z_N$. It is possible to fix this issue by talking about $G$-patterned Hadamard matrices, with $G$ being a finite abelian group, but for our purposes here, the best is to formulate the result in a weaker form, as follows:

\index{generalized Fourier matrix}
\index{bistochastic form}

\begin{theorem}
The generalized Fourier matrices $F_G$, associated to the finite abelian groups $G$, can be put in symmetric and bistochastic form.
\end{theorem} 

\begin{proof}
We know from Theorem 9.8 that any usual Fourier matrix $F_N$ can be put in circulant and symmetric form. Since circulant implies bistochastic, in the sense that the sums on all rows and all columns must be equal, the result holds for $F_N$. In general now, if we decompose our group as $G=\mathbb Z_{N_1}\times\ldots\times\mathbb Z_{N_k}$, we have:
$$F_G=F_{N_1}\otimes\ldots\otimes F_{N_k}$$

Now since the property of being circulant is stable under taking tensor products, and so is the property of being bistochastic, we therefore obtain the result.
\end{proof}

We have as well the following alternative generalization of Theorem 9.8, coming from Backelin's work in \cite{bac}, and remaining in the circulant and symmetric setting:

\index{Backelin construction}

\begin{theorem}
Let $M|N$, and set $w=e^{2\pi i/N}$. We have a cyclic root as follows,
$$(\,\underbrace{q_1,\ldots,q_M}_M\,,\,\underbrace{q_1w,\ldots,q_Mw}_M\,,\ldots\ldots,\,\underbrace{q_1w^{N-1},\ldots,q_Mw^{N-1}}_M\,)$$
provided that $q_1,\ldots,q_M\in\mathbb T$ satisfy the following condition:
$$(q_1\ldots q_M)^N=(-1)^{M(N-1)}$$
Moreover, assuming that the following conditions are satisfied,
$$q_1q_2=1\quad,\quad q_3q_M=q_4q_{M-1}=\ldots=w$$
which imply $(q_1\ldots q_M)^N=(-1)^{M(N-1)}$, the Hadamard matrix is symmetric.
\end{theorem}

\begin{proof}
We have several things to be proved, the idea being as follows:

\medskip

(1) Let us first check the $=0$ equations for a cyclic root. Given arbitrary numbers $q_1,\ldots,q_M\in\mathbb T$, if we denote by $(z_i)$ the vector in the statement, we have:
\begin{eqnarray*}
\sum_iz_{i+1}\ldots z_{i+K}
&=&\begin{pmatrix}q_1\ldots q_K+q_2\ldots q_{K+1}+\ldots\ldots+q_{M-K+1}\ldots q_M\\
+q_{M-K+2}\ldots q_Mq_1w+\ldots\ldots+q_Mq_1\ldots q_{K-1}w^{K-1}\end{pmatrix}\\
&&\times(1+w^K+w^{2K}+\ldots+w^{(N-1)K})
\end{eqnarray*}

Now since the sum on the right vanishes, the $=0$ conditions are satisfied. 

\medskip

(2) Regarding now the $=1$ condition, the total product of the numbers $z_i$ is given by:
\begin{eqnarray*}
\prod_iz_i
&=&(q_1\ldots q_M)^N(1\cdot w\cdot w^2\ldots w^{N-1})^M\\
&=&(q_1\ldots q_M)^Nw^{\frac{MN(N-1)}{2}}
\end{eqnarray*}

By using $w=e^{2\pi i/N}$ we obtain that the coefficient on the right is:
\begin{eqnarray*}
w^{\frac{MN(N-1)}{2}}
&=&e^{\frac{2\pi i}{N}\cdot\frac{MN(N-1)}{2}}\\
&=&e^{\pi iM(N-1)}\\
&=&(-1)^{M(N-1)}
\end{eqnarray*}

Thus, if $(q_1\ldots q_M)^N=(-1)^{M(N-1)}$, we obtain a cyclic root, as stated. For further details on all this, we refer to the papers of Backelin \cite{bac} and Faug\`ere \cite{fau}.

\medskip

(3) The corresponding first row vector can be written as follows:
$$V=\left(\underbrace{q_1,q_1q_2,\ldots,q_1\ldots q_M}_M\,,\ldots\ldots\ldots,\underbrace{\frac{w^{M-1}}{q_2\ldots q_M},\ldots,\frac{w^2}{q_{M-1}q_M},\frac{w}{q_M},1}_M\right)$$

Thus, the corresponding circulant complex Hadamard matrix is as follows:
$$H=\begin{pmatrix}
q_1&q_1q_2&q_1q_2q_3&q_1q_2q_3q_4&q_1q_2q_3q_4q_5&\ldots\\
1&q_1&q_1q_2&q_1q_2q_3&q_1q_2q_3q_4&\ldots\\
\frac{w}{q_M}&1&q_1&q_1q_2&q_1q_2q_3&\ldots\\
\frac{w^2}{q_{M-1}q_M}&\frac{w}{q_M}&1&q_1&q_1q_2&\ldots\\
\frac{w^3}{q_{M-2}q_{M-1}q_M}&\frac{w^2}{q_{M-1}q_M}&\frac{w}{q_M}&1&q_1&\ldots\\
\vdots&\vdots&\vdots&\vdots&\vdots&\ddots
\end{pmatrix}$$

We are therefore led to the symmetry conditions in the statement, and we are done.
\end{proof}

Observe that the story is not over here, because Theorem 9.11 still remains to be unified with Theorem 9.10. There are many interesting questions here.

\section*{9b. Butson matrices}

Still in relation with the CHC, the problem of investigating the existence of circulant Butson matrices of a given level appears. Following Turyn \cite{tur}, we first have:

\index{size of circulant matrix}

\begin{proposition}
The size of a circulant real Hadamard matrix 
$$H\in M_N(\pm 1)$$
must be of the form $N=4n^2$, with $n\in\mathbb N$.
\end{proposition}

\begin{proof}
Let $a,b\in\mathbb N$ with $a+b=N$ be the number of $1,-1$ entries in the first row of $H$. If we denote by $H_0,\ldots,H_{N-1}$ the rows of $H$, by summing over columns we get:
\begin{eqnarray*}
\sum_{i=0}^{N-1}<H_0,H_i>
&=&a(a-b)+b(b-a)\\
&=&(a-b)^2
\end{eqnarray*}

On the other hand, by orthogonality of the rows, the quantity on the left is:
$$<H_0,H_0>=N$$

Thus $N=(a-b)^2$ is a square, and since $N\in 2\mathbb N$, this gives $N=4n^2$, with $n\in\mathbb N$.
\end{proof}

Also found by Turyn in \cite{tur} is the fact that the above number $n\in\mathbb N$ must be odd, and not a prime power. In the general Butson matrix setting now, we have:

\index{circulant Butson matrix}
\index{Turyn obstruction}

\begin{proposition}
Assume that $H\in H_N(l)$ is circulant, let $w=e^{2\pi {\rm i}/l}$. If 
$$a_0,\ldots,a_{l-1}\in\mathbb N$$
with $\sum a_i=N$ are the number of $1,w,\ldots,w^{l-1}$ entries in the first row of $H$, then:
$$\sum_{ik}w^ka_ia_{i+k}=N$$
This condition, with $\sum a_i=N$, will be called ``Turyn obstruction'' on $(N,l)$.
\end{proposition}

\begin{proof}
Indeed, by summing over the columns of $H$, we obtain:
\begin{eqnarray*}
\sum_i<H_0,H_i>
&=&\sum_{ij}<w^i,w^j>a_ia_j\\
&=&\sum_{ij}w^{i-j}a_ia_j
\end{eqnarray*}

Now since the left term is $<H_0,H_0>=N$, this gives the result.
\end{proof}

We can deduce from this a number of concrete obstructions, as follows:

\begin{theorem}
When $l$ is prime, the Turyn obstruction is 
$$\sum_i(a_i-a_{i+k})^2=2N$$
for any $k\neq 0$. Also, for small values of $l$, the Turyn obstruction is as follows:
\begin{enumerate}
\item At $l=2$ the condition is:
$$(a_0-a_1)^2=N$$

\item At $l=3$ the condition is:
$$(a_0-a_1)^2+(a_1-a_2)^2+(a_2-a_3)^2=2N$$

\item At $l=4$ the condition is:
$$(a_0-a_2)^2+(a_1-a_3)^2=N$$

\item At $l=5$ the condition is:
$$\sum_i(a_i-a_{i+1})^2=\sum_i(a_i-a_{i+2})^2=2N$$
\end{enumerate}
\end{theorem}

\begin{proof}
We use the fact, from chapter 6, that when $l$ is prime, the vanishing sums of $l$-roots of unity are exactly the sums of the following type, with $c\in\mathbb N$:
$$S=c+cw+\ldots+cw^{l-1}$$

We conclude that the Turyn obstruction is equivalent to the following system of equations, one for each $k\neq 0$:
$$\sum_ia_i^2-\sum_ia_ia_{i+k}=N$$

Now by forming squares, this gives the equations in the statement. Regarding now the $l=2,3,4,5$ assertions, these follow from the first assertion when $l$ is prime, $l=2,3,5$. Also, at $l=4$ we have $w=i$, so the Turyn obstruction reads:
$$(a_0^2+a_1^2+a_2^2+a_3^2)+i\sum a_ia_{i+1}-2(a_0a_2+a_1a_3)-i\sum a_ia_{i+1}=N$$

Thus the imaginary terms cancel, and we obtain the formula in the statement.
\end{proof}

The above results are of course just some basic observations on the subject, and the massive amount of work on the CHC has a number of interesting Butson matrix extensions. For some more advanced theory on all this, we refer to \cite{bs1}, \cite{ckh}.

\section*{9c. Haagerup count}

Let us go back now to the pure complex case, and discuss Fourier analytic aspects. From a traditional linear algebra viewpoint, the circulant matrices are best understood as being the matrices which are Fourier-diagonal, and we will exploit this here. Let us fix $N\in\mathbb N$, and denote by $F=(w^{ij})/\sqrt{N}$ with $w=e^{2\pi i/N}$ the rescaled Fourier matrix, with indices $i,j=0,1,\ldots,N-1$, which is unitary, given by the following formula:
$$F=\frac{1}{\sqrt{N}}\begin{pmatrix}
1&1&1&\ldots&1\\
1&w&w^2&\ldots&w^{N-1}\\
1&w^2&w^4&\ldots&w^{2(N-1)}\\
\vdots&\vdots&\vdots&&\vdots\\
1&w^{N-1}&w^{2(N-1)}&\ldots&w^{(N-1)^2}
\end{pmatrix}$$

Also, given a vector $q\in\mathbb C^N$, once again with cyclic indices, $i=0,1,\ldots,N-1$, we denote by $Q\in M_N(\mathbb C)$ the diagonal matrix having $q$ as vector of diagonal entries:
$$Q=\begin{pmatrix}
q_0\\
&\ddots\\
&&q_{N-1}
\end{pmatrix}$$

With these conventions, we have the following well-known result, that we have already used in this book, but that we reproduce here for convenience:

\index{circulant matrix}
\index{Fourier-diagonal}

\begin{theorem}
For a complex matrix $H\in M_N(\mathbb C)$, the following are equivalent:
\begin{enumerate}
\item $H$ is circulant, $H_{ij}=\xi_{j-i}$ for some $\xi\in\mathbb C^N$.

\item $H$ is Fourier-diagonal, $H=FQF^*$ with $Q$ diagonal.
\end{enumerate}
In addition, the first row vector of $FQF^*$ is given by $\xi=Fq/\sqrt{N}$.
\end{theorem}

\begin{proof}
If $H_{ij}=\xi_{j-i}$ is circulant then $Q=F^*HF$ is diagonal, given by:
\begin{eqnarray*}
Q_{ij}
&=&\frac{1}{N}\sum_{kl}w^{jl-ik}\xi_{l-k}\\
&=&\delta_{ij}\sum_rw^{jr}\xi_r
\end{eqnarray*}

Also, if $Q=diag(q)$ is diagonal then $H=FQF^*$ is circulant, given by:
\begin{eqnarray*}
H_{ij}
&=&\sum_kF_{ik}Q_{kk}\bar{F}_{jk}\\
&=&\frac{1}{N}\sum_kw^{(i-j)k}q_k
\end{eqnarray*}

Thus, we have proved the equivalence between the conditions in the statement. Finally, regarding $\xi=Fq/\sqrt{N}$, this follows from the last formula established above.
\end{proof}

The above result is useful in connection with any question regarding the circular matrices, and in relation with the orthogonal and unitary cases, we have:

\index{circulant orthogonal matrix}
\index{circulant unitary matrix}

\begin{proposition}
The various sets of circulant matrices are as follows:
\begin{enumerate}
\item The set of all circulant matrices is:
$$M_N(\mathbb C)^{circ}=\left\{FQF^*\Big|q\in\mathbb C^N\right\}$$

\item The set of all circulant unitary matrices is:
$$U_N^{circ}=\left\{FQF^*\Big|q\in\mathbb T^N\right\}$$

\item The set of all circulant orthogonal matrices is:
$$O_N^{circ}=\left\{FQF^*\Big|q\in\mathbb T^N,\bar{q}_i=q_{-i},\forall i\right\}$$
\end{enumerate}
In addition, the first row vector of $FQF^*$ is given by $\xi=Fq/\sqrt{N}$.
\end{proposition}

\begin{proof}
All this follows from Theorem 9.15, as follows:

\medskip

(1) This assertion, along with the last one, is Theorem 9.15 itself.

\medskip

(2) This is clear from (1), because the eigenvalues must be on the unit circle $\mathbb T$.

\medskip

(3) In order to prove this result, observe first that for a vector $q\in\mathbb C^N$ we have the following formula, with $\tilde{q}_i=\bar{q}_{-i}$:
$$\overline{Fq}=F\tilde{q}$$

We conclude from this that the vector $\xi=Fq$ is real if and only if $\bar{q}_i=q_{-i}$ for any $i$. Together with (2), this gives the result.
\end{proof}

Observe that in Proposition 9.16 (3), the equations for the parameter space for $O_N^{circ}$ are as follows, going until $[N/2]+1$:
$$q_0=\bar{q}_0\quad,\quad 
\bar{q}_1=q_{n-1}\quad,\quad 
\bar{q}_2=q_{n-2}\quad,\quad
\ldots$$
 
Thus, with the convention $\mathbb Z_\infty=\mathbb T$, we have the following formula:
$$O_N^{circ}\simeq
\begin{cases}
\mathbb Z_2\times\mathbb Z_\infty^{{(N-1)}/2}&(N\ {\rm odd})\\
\mathbb Z_2^2\times\mathbb Z_\infty^{(N-2)/2}&(N\ {\rm even})
\end{cases}$$

In terms of circulant Hadamard matrices, we have the following statement:

\begin{theorem}
The sets of complex and real circulant Hadamard matrices are:
$$X_N^{circ}=\left\{\sqrt{N}FQF^*\Big|q\in\mathbb T^N\right\}\cap M_N(\mathbb T)$$
$$Y_N^{circ}=\left\{\sqrt{N}FQF^*\Big|q\in\mathbb T^N,\bar{q}_i=q_{-i}\right\}\cap M_N(\pm1)$$
In addition, the sets of $q$ parameters are invariant under cyclic permutations, and also under mutiplying by numbers in $\mathbb T$, respectively under multiplying by $-1$. 
\end{theorem}

\begin{proof}
All the assertions are indeed clear from Proposition 9.16, by intersecting the sets computed there with $M_N(\mathbb T)$.
\end{proof}

The above statement is of course something quite theoretical in the real case, where the CHC states that we should have $Y_N^{circ}=\emptyset$, at any $N\neq 4$. However, in the complex case all this is useful, and complementary to Bj\"orck's cyclic root formalism. Indeed, let us discuss now a number of geometric and analytic aspects, in the complex matrix case. First, we have the following deep counting result, due to Haagerup \cite{ha2}:

\index{Haagerup counting theorem}
\index{Chebotarev theorem}
\index{minors}
\index{circulant Hadamard matrix}

\begin{theorem}
When $N$ is prime, the number of circulant $N\times N$ complex Hadamard matrices, counted with certain multiplicities, is exactly: 
$$N_{circ}=\binom{2N-2}{N-1}$$
\end{theorem}

\begin{proof}
This is something advanced, using a variety of techiques from Fourier analysis, number theory, complex analysis and algebraic geometry. The idea is as follows:

\medskip

(1) As explained in \cite{ha2}, when $N$ is prime, Bj\"orck's cyclic root formalism, explained above, can be further manipulated, by using discrete Fourier transforms, and we are eventually led to a simpler system of equations. 

\medskip

(2) This simplified system can be shown then to have a finite number of solutions, the key ingredient here being a well-known theorem of Chebotarev, which states that when $N$ is prime, all the minors of the Fourier matrix $F_N$ are nonzero. 

\medskip

(3) With this finiteness result in hand, the precise count can be done as well, by using various techniques from classical algebraic geometry, and we are led to the formula in the statement. For the details here, we refer to Haagerup's paper \cite{ha2}.
\end{proof}

When $N$ is not prime, the situation is considerably more complicated, with some values leading to finitely many solutions, and with other values leading to an infinite number of solutions, and with many other new phenomena appearing. We refer here to the papers of  Bj\"orck \cite{bjo}, Bj\"orck-Fr\"oberg \cite{bjf}, Bj\"orck-Haagerup \cite{bha} and Haagerup \cite{ha2}.

\section*{9d. Analytic aspects}

Let us discuss now an alternative take on these questions, based on the $p$-norm considerations from chapter 3. As explained in \cite{bs1}, the most adapted exponent for the circulant case is $p=4$. So, as a starting point, let us formulate:

\index{norm minimizer}
\index{four-norm}

\begin{proposition}
Given a matrix $U\in U_N$ we have
$$||U||_4\geq 1$$
with equality precisely when $H=U/\sqrt{N}$ is Hadamard.
\end{proposition}

\begin{proof}
This is something that we already know, from chapter 3, as a particular case of our results there regarding $p$-norms, obtained by using the Jensen inequality. However, this follows as well directly from the Cauchy-Schwarz inequality, as follows:
\begin{eqnarray*}
||U||_4^4
&=&\sum_{ij}|U_{ij}|^4\\
&\geq&\frac{1}{N^2}\left(\sum_{ij}|U_{ij}|^2\right)^2\\
&=&1
\end{eqnarray*}

Thus we have $||U||_4\geq 1$, with equality if and only if $H=\sqrt{N}U$ is Hadamard.
\end{proof}

In the circulant case now, and in Fourier formulation, the estimate is as follows:

\begin{theorem}
Given a vector $q\in\mathbb T^N$, written $q=(q_0,\ldots,q_{N-1})$ consider the following quantity, with all the indices being taken modulo $N$:
$$\Phi=\sum_{i+k=j+l}\frac{q_iq_k}{q_jq_l}$$
Then this quantity $\Phi$ is real, and we have the estimate
$$\Phi\geq N^2$$
with the equality case happening precisely when $\sqrt{N}q$ is the eigenvalue vector of a circulant Hadamard matrix $H\in M_N(\mathbb C)$. 
\end{theorem}

\begin{proof}
By conjugating the formula of $\Phi$ we see that this quantity is indeed real, as stated. In fact, $\Phi$ appears by definition as a sum of $N^3$ terms, consisting of $N(2N-1)$ values of $1$ and of $N(N-1)^2$ other complex numbers of modulus 1, coming in pairs $(a,\bar{a})$. Regarding now the second assertion, by using the various identifications in Theorem 9.15 and Proposition 9.16, and the formula $\xi=Fq/\sqrt{N}$ there, we have:
\begin{eqnarray*}
||U||_4^4
&=&N\sum_s|\xi_s|^4\\
&=&\frac{1}{N^3}\sum_s\left|\sum_iw^{si}q_i\right|^4\\
&=&\frac{1}{N^3}\sum_s\sum_iw^{si}q_i\sum_jw^{-sj}\bar{q}_j\sum_kw^{sk}q_k\sum_lw^{-sl}\bar{q}_l\\
&=&\frac{1}{N^3}\sum_s\sum_{ijkl}w^{(i-j+k-l)s}\frac{q_iq_k}{q_jq_l}\\
&=&\frac{1}{N^2}\sum_{i+k=j+l}\frac{q_iq_k}{q_jq_l}
\end{eqnarray*}

Thus Proposition 9.19 gives the following estimate:
$$\Phi
=N^2||U||_4^4
\geq N^2$$

Moreover, we have equality precisely in the Hadamard matrix case, as claimed.
\end{proof}

The above result is something quite subtle, and even surprising, at the level of the consequences. We have as well the following more direct explanation for it:

\begin{proposition}
With the above notations, we have the formula
$$\Phi=N^2+\sum_{i\neq j}(|\nu_i|^2-|\nu_j|^2)^2$$
where $\nu=(\nu_0,\ldots,\nu_{N-1})$ is the vector given by $\nu=Fq$.
\end{proposition}

\begin{proof}
This follows by replacing in the above proof the Cauchy-Schwarz estimate by the corresponding sum of squares. More precisely, we know from the above proof that:
$$\Phi=N^3\sum_i|\xi_i|^4$$

On the other hand the matrix $U_{ij}=\xi_{j-i}$ being unitary, we have:
$$\sum_i|\xi_i|^2=1$$

We therefore have the following computation:
\begin{eqnarray*}
1
&=&\sum_i|\xi_i|^4+\sum_{i\neq j}|\xi_i|^2\cdot|\xi_j|^2\\
&=&N\sum_i|\xi_i|^4-\left((N-1)\sum_i|\xi_i|^4-\sum_{i\neq j}|\xi_i|^2\cdot|\xi_j|^2\right)\\
&=&\frac{1}{N^2}\Phi-\sum_{i\neq j}(|\xi_i|^2-|\xi_j|^2)^2
\end{eqnarray*}

Now by multiplying by $N^2$, this gives the formula in the statement.
\end{proof}

Let us explore now the minimization problem for $\Phi$, by using various combinatorial and analytic methods. As an illustration for the difficulties in dealing with this problem, let us work out the case where $N$ is small. At $N=1$ our inequality $\Phi\geq N^2$ is simply:
$$\Phi=1\geq 1$$

At $N=2$ our inequality is also clearly true, as follows:
$$\Phi=6+\left(\frac{q_0}{q_1}\right)^2+\left(\frac{q_1}{q_0}\right)^2\geq 4$$

At $N=3$ now, the inequality is something more subtle:
$$\Phi=15+4Re\left(\frac{q_0^3+q_1^3+q_2^3}{q_0q_1q_2}\right)\geq 9$$

Observe that in terms of $a=q_0^2/(q_1q_2)$, $b=q_1^2/(q_0q_2)$, $c=q_2^2/(q_0q_1)$, which satisfy $|a|=|b|=|c|=1$ and $abc=1$, our function is:
$$\Phi=15+4Re(a+b+c)$$

Thus at $N=3$ our inequality still has a quite tractable form, namely:
$$|a|=|b|=|c|=1,abc=1\implies Re(a+b+c)\geq-\frac{3}{2}$$

At $N=4$ however, the formula of $\Phi$ is as follows:
\begin{eqnarray*}\Phi
&=&28+4\left(\frac{q_0q_1}{q_2q_3}+\frac{q_2q_3}{q_0q_1}+\frac{q_0q_3}{q_1q_2}+\frac{q_1q_2}{q_0q_3}\right)+\left(\frac{q_0^2}{q_2^2}+\frac{q_2^2}{q_0^2}+\frac{q_1^2}{q_3^2}+\frac{q_3^2}{q_1^2}\right)\\
&&+2\left(\frac{q_0q_2}{q_1^2}+\frac{q_1^2}{q_0q_2}+\frac{q_0q_2}{q_3^2}+\frac{q_3^2}{q_0q_2}+\frac{q_1q_3}{q_0^2}+\frac{q_0^2}{q_1q_3}+\frac{q_1q_3}{q_2^2}+\frac{q_2^2}{q_1q_3}\right)
\end{eqnarray*}

It is not clear how to obtain a simple, direct proof of $\Phi\geq 16$, based on this formula. This is actually a quite challenging calculus problem, and we will be back to it, most likely on the occasion of our next exercise session.

\bigskip

As an application of the above considerations, in the real Hadamard matrix case, we have the following analytic reformulation of the CHC, from \cite{bs1}:

\index{CHC}
\index{Circulant Hadamard Conjecture}

\begin{theorem}
For a vector $q\in\mathbb T^N$ satisfying $\bar{q}_i=q_{-i}$ the following quantity is real,
$$\Phi=\sum_{i+j+k+l=0}q_iq_jq_kq_l$$
and satisfies the following inequality:
$$\Phi\geq N^2$$
The CHC states that we cannot have equality at $N>4$.
\end{theorem}

\begin{proof}
This follows indeed from Theorem 9.20, via the identifications from Proposition 9.16, the parameter space in the real case being $\left\{q\in\mathbb T^N|\bar{q}_i=q_{-i}\right\}$.
\end{proof}

Following \cite{bs1}, let us further discuss all this. We first have the following result:

\begin{theorem}
Let us decompose the above function as 
$$\Phi=\Phi_0+\ldots+\Phi_{N-1}$$
with each $\Phi_i$ being given by the same formula as $\Phi$, namely 
$$\Phi=\sum_{i+k=j+l}\frac{q_iq_k}{q_jq_k}$$
but keeping the index $i$ fixed. Then:
\begin{enumerate}
\item The critical points of $\Phi$ are those where $\Phi_i\in\mathbb R$, for any $i$.

\item In the Hadamard case we have $\Phi_i=N$, for any $i$.
\end{enumerate}
\end{theorem}

\begin{proof}
This follows by doing some elementary computations, as follows:

\medskip

(1) The first observation is that the non-constant terms in the definition of $\Phi$ involving the variable $q_i$ are the terms of the sum $K_i+\bar{K}_i$, where:
$$K_i=\sum_{2i=j+l}\frac{q_i^2}{q_jq_l}+2\sum_{k\neq i}\sum_{i+k=j+l}\frac{q_iq_k}{q_jq_l}$$

Thus if we fix $i$ and we write $q_i=e^{i\alpha_i}$, we obtain:
\begin{eqnarray*}
\frac{d\Phi}{d\alpha_i}
&=&4Re\left(\sum_k\sum_{i+k=j+l}i\cdot\frac{q_iq_k}{q_jq_l}\right)\\
&=&4Im\left(\sum_{i+k=j+l}\frac{q_iq_k}{q_jq_l}\right)\\
&=&4Im(\Phi_i)
\end{eqnarray*}

Now since the derivative must vanish for any $i$, this gives the result.

\medskip

(2) We first perform the end of the Fourier computation in the proof of Theorem 9.20 above backwards, by keeping the index $i$ fixed. We obtain: 
\begin{eqnarray*}
\Phi_i
&=&\sum_{i+k=j+l}\frac{q_iq_k}{q_jq_l}\\
&=&\frac{1}{N}\sum_s\sum_{ijkl}w^{(i-j+k-l)s}\frac{q_iq_k}{q_jq_l}\\
&=&\frac{1}{N}\sum_sw^{si}q_i\sum_jw^{-sj}\bar{q}_j\sum_kw^{sk}q_k\sum_lw^{-sl}\bar{q}_l\\
&=&N^2\sum_sw^{si}q_i\bar{\xi}_s\xi_s\bar{\xi}_s
\end{eqnarray*}

Here we have used the formula $\xi=Fq/\sqrt{N}$. Now by assuming that we are in the Hadamard case, we have $|\xi_s|=1/\sqrt{N}$ for any $s$, and so we obtain:
\begin{eqnarray*}
\Phi_i
&=&N\sum_s w^{si}q_i\bar{\xi}_s\\
&=&N\sqrt{N}q_i\overline{(F^*\xi)}_i\\
&=&Nq_i\bar{q}_i\\
&=&N
\end{eqnarray*}

Thus, we have obtained the conclusion in the statement.
\end{proof}

Let us discuss now a probabilistic approach to all this. Given a compact manifold $X$ endowed with a probability measure, and a bounded function $\Theta:X\to[0,\infty)$, the maximum of this function can be recaptured via following well-known formula:
$$\max\Theta=\lim_{p\to\infty}\left(\int_X\Theta(x)^p\,dx\right)^{1/p}$$

In our case, we are rather interested in computing a minimum, and we have:

\begin{proposition}
We have the formula
$$\min\Phi=N^3-\lim_{p\to\infty}\left(\int_{\mathbb T^N}(N^3-\Phi)^p\,dq\right)^{1/p}$$
where the torus $\mathbb T^N$ is endowed with its usual probability measure.
\end{proposition}

\begin{proof}
This follows from the above formula, with $\Theta=N^3-\Phi$. Observe that $\Theta$ is indeed positive, because $\Phi$ is a sum of $N^3$ complex numbers of modulus 1. 
\end{proof}

Let us restrict now the attention to the problem of computing the moments of $\Phi$, which is more or less the same as computing those of $N^3-\Phi$. We have here:

\begin{proposition}
The moments of $\Phi$ are given by
$$\int_{\mathbb T^N}\Phi^p\,dq=\#\left\{ \begin{pmatrix}i_1k_1\ldots i_pk_p\\ j_1l_1\ldots j_pl_p\end{pmatrix}\Big|i_s+k_s=j_s+l_s,[i_1k_1\ldots i_pk_p]=[j_1l_1\ldots j_pl_p]\right\}$$
where the sets between brackets are by definition sets with repetition. 
\end{proposition}

\begin{proof}
This is indeed clear from the formula of $\Phi$. See \cite{bs1}.
\end{proof}

Regarding now the real case, an analogue of Proposition 9.25 holds, but the combinatorics does not get any simpler. One idea in dealing with this problem is by considering the ``enveloping sum'', obtained from $\Phi$ by dropping the condition $i+k=j+l$:
$$\tilde{\Phi}=\sum_{ijkl}\frac{q_iq_k}{q_jq_l}$$

The point is that the moments of $\Phi$ appear as ``sub-quantities'' of the moments of $\tilde{\Phi}$, so perhaps the question to start with is to understand very well the moments of $\tilde{\Phi}$. And this latter problem sounds like a quite familiar one, because:
$$\tilde{\Phi}=\left|\sum_iq_i\right|^4$$

We will be back to this later. For the moment, let us do some combinatorics:

\index{partition}

\begin{proposition}
We have the moment formula
$$\int_{\mathbb T^N}\tilde{\Phi}^p\,dq=\sum_{\pi\in P(2p)}\binom{2p}{\pi}\frac{N!}{(N-|\pi|)!}$$
where the coefficients on the right are given by
$$\binom{2p}{\pi}=\binom{2p}{b_1,\ldots,b_{|\pi|}}$$
with $b_1,\ldots,b_{|\pi|}$ being the lengths of the blocks of $\pi$.
\end{proposition}

\begin{proof}
Indeed, by using the same method as for $\Phi$, we obtain:
$$\int_{\mathbb T^N}\tilde{\Phi}(q)^p\,dq=\#\left\{ \begin{pmatrix}i_1k_1\ldots i_pk_p\\ j_1l_1\ldots j_pl_p\end{pmatrix}\Big|[i_1k_1\ldots i_pk_p]=[j_1l_1\ldots j_pl_p]\right\}$$

The sets with repetitions on the right are best counted by introducing the corresponding partitions $\pi=\ker\begin{pmatrix}i_1k_1\ldots i_pk_p\end{pmatrix}$, and this gives the formula in the statement.
\end{proof}

In order to discuss now the real case, we have to slightly generalize the above result, by computing all the half-moments of $\widetilde{\Phi}$. The result here is best formulated as:

\begin{proposition}
We have the moment formula
$$\int_{\mathbb T^N}\left|\sum_iq_i\right|^{2p}\,dq=\sum_kC_{pk}\frac{N!}{(N-k)!}$$
with the coefficients being given by
$$C_{pk}=\sum_{\pi\in P(p),|\pi|=k}\binom{p}{b_1,\ldots,b_{|\pi|}}$$
where $b_1,\ldots,b_{|\pi|}$ are the lengths of the blocks of $\pi$.
\end{proposition}

\begin{proof}
This follows indeed exactly as Proposition 9.26, by replacing the exponent $p$ by the exponent $p/2$, and by splitting the resulting sum as in the statement.
\end{proof}

Finally, here is a random walk formulation of the problem:

\index{random walk}
\index{piecewise balanced}

\begin{theorem}
The moments of $\Phi$ have the following interpretation:
\begin{enumerate}
\item First, the moments of the enveloping sum $\int\widetilde{\Phi}^p$ count the loops of length $4p$ on the standard lattice $\mathbb Z^N\subset\mathbb R^N$, based at the origin.

\item $\int\Phi^p$ counts those loops which are ``piecewise balanced'', in the sense that each of the $p$ consecutive $4$-paths forming the loop satisfy $i+k=j+l$ modulo $N$.
\end{enumerate}
\end{theorem}

\begin{proof}
The first assertion follows from the formula in the proof of Proposition 9.26, and the second assertion follows from the formula in Proposition 9.25. 
\end{proof}

There are many interesting questions here. We refer to \cite{bs1} for more on all this.

\section*{9e. Exercises} 

In relation with the Butson matrices, we have the following exercise:

\begin{exercise}
Work out the Turyn obstruction for the circulant Butson matrices at the exponent values $l=6,7,8$.
\end{exercise}

To be more precise, we have seen in the above how to deal with such questions at the exponent values $l=2,3,4,5$, and the problem now is that continuing that work.

\begin{exercise}
Work out formulae or estimates for the number of circulant complex $N\times N$ complex Hadamard matrices, at small values of $N\in\mathbb N$, not prime.
\end{exercise}

This is something quite tricky, normally requiring some computer programming.

\begin{exercise}
Find a proof for the estimate $\Phi\geq16$ at $N=4$.
\end{exercise}

This question was already mentioned in the above, with the comment that there is no obvious proof. The problem is that of finding a reasonably elementary proof.

\chapter{Bistochastic form}

\section*{10a. Basic theory}

In this chapter and in the next one we discuss some further analytic aspects of the complex Hadamard matrices, which this time are brand new or almost, going back to the mid 10s and onwards, and are very exciting too. The general idea is that any Hadamard matrix, real or complex, can be put in bistochastic form over the complex numbers $\mathbb C$, and with this bistochastic form looking much better than the original form.

\bigskip

Thus, we have here a potentially far-reaching idea, consisting in reformulating everything that we know, including our favorite questions from the real case, the HC and CHC, in complex bistochastic form. But, and here comes the second point, putting an Hadamard matrix in bistochastic form is something non-trivial, in general done by a non-explicit result of Idel-Wolf \cite{iwo}, based on some non-trivial symplectic geometry results of Biran-Entov-Polterovich \cite{bep} and Cho \cite{cho}, motivated by a deep conjecture of Arnold.

\bigskip

And isn't this exciting. We have been commenting in the last chapter on open questions in mathematics, our point being that the closer you get to classical mechanics, the better that is, for the fate of your open problem. And since in classical mechanics all roads lead to Arnold, we are probably on the right track here. Perhaps for the first time, since the beginning of this book. That is, plenty of reasons to be optimistic.

\bigskip

All this is however very new, and our presentation here will be quite modest. Lots of further work are needed, and it is a pity that nothing much is going on here, so far. Young reader, if I have an excellent question to recommend to you, this is the one, continuation of what will be said here. Get to know and love classical mechanics, which is the mother of everything, in mathematics and physics, than read some books of Arnold, starting with \cite{arn}, which are a must-read anyway, no matter what mathematics or physics you want to do, and then start solving some Hadamard matrix questions, using this technology.

\bigskip

In order to get started now, we have already talked about bistochastic Hadamard matrices, in the real case, on several occasions, in chapters 1-4. Our first purpose will be that of carefully reviewing and extending that material, in the complex Hadamard matrix case. Let us start our discussion with the following definition:

\index{bistochastic matrix}

\begin{definition}
A complex Hadamard matrix $H\in M_N(\mathbb C)$ is called bistochastic when the sums on all rows and all columns are equal,
$$\sum_iH_{ij}=\sum_jH_{ij}=\lambda$$
for a certain number $\lambda\in\mathbb C$. We denote by
$$X_N^{bis}=\left\{H\in X_N\Big|\,H={\rm bistochastic}\right\}$$
the real algebraic manifold formed by such matrices.
\end{definition}

The bistochastic Hadamard matrices are quite interesting objects, and include for instance all the circulant Hadamard matrices, that we discussed in chapter 9. Indeed, assuming that $H_{ij}=\xi_{j-i}$ is circulant, all rows and columns sum up to $\lambda=\sum_i\xi_i$:
$$\sum_i\xi_{j-i}=\sum_j\xi_{j-i}=\sum_i\xi_i$$

We will be back to this, in a moment. Let us begin, however, with some considerations regarding the real case. Our point here is that the real Hadamard matrices often ``look better'' in complex bistochastic form, and that there is some potentially interesting mathematics behind all this. As a first and trivial remark, the first Walsh matrix $W_2=F_2$ looks better in complex bistochastic form, modulo the standard equivalence relation:
$$W_2=\begin{pmatrix}1&1\\1&-1\end{pmatrix}
\sim\begin{pmatrix}i&i\\1&-1\end{pmatrix}
\sim\begin{pmatrix}i&1\\1&i\end{pmatrix}$$

To be more precise, the matrix on the right, while having the slight disadvantage of being complex instead of real, is something very nice, circulant and symmetric. Regarding the second Walsh matrix $W_4=W_2\otimes W_2$, this looks as well better in bistochastic form, because it becomes in this way equivalent to $K_4$, the most beautiful matrix ever:
\begin{eqnarray*}
W_4&=&\begin{pmatrix}1&1&1&1\\ 1&-1&1&-1\\ 1&1&-1&-1\\ 1&-1&-1&1\end{pmatrix}
\sim\begin{pmatrix}
1&1&1&1\\
1&1&-1&-1\\
1&-1&1&-1\\
1&-1&-1&1
\end{pmatrix}\\
&\sim&\begin{pmatrix}
1&-1&-1&-1\\
1&-1&1&1\\
1&1&-1&1\\
1&1&1&-1
\end{pmatrix}
\sim\begin{pmatrix}
-1&1&1&1\\
1&-1&1&1\\
1&1&-1&1\\
1&1&1&-1
\end{pmatrix}
\end{eqnarray*}

As before, the matrix on the right looks better than the one on the left, because it is circulant and symmetric. And all this is quite interesting, philosophically speaking. Indeed, we have here a new idea, in connection with the various questions explained in chapters 1-4, namely that of studying the real Hadamard matrices $H\in M_N(\pm1)$ by putting them in complex bistochastic form, $H'\in M_N(\mathbb T)$, and then studying these latter matrices. Let us record here, as a partial conclusion, the following simple fact:

\index{Walsh matrix}

\begin{theorem}
All Walsh matrices can be put in bistocastic form, as follows:
\begin{enumerate}
\item The matrices $W_N$ with $N=4^n$ admit a real bistochastic form, namely:
$$W_N\sim\begin{pmatrix}
-1&1&1&1\\
1&-1&1&1\\
1&1&-1&1\\
1&1&1&-1
\end{pmatrix}^{\otimes n}$$

\item The matrices $W_N$ with $N=2\times4^n$ admit a complex bistochastic form, namely:
$$W_N\sim\begin{pmatrix}i&1\\1&i\end{pmatrix}\otimes\begin{pmatrix}
-1&1&1&1\\
1&-1&1&1\\
1&1&-1&1\\
1&1&1&-1
\end{pmatrix}^{\otimes n}$$
\end{enumerate}
\end{theorem}

\begin{proof}
This follows indeed from the above discussion.
\end{proof}

Let us review now the material in chapter 9. According to the results there, and to the above-mentioned fact that circulant implies bistochastic, we have:

\index{Backelin matrix}
\index{circulant and symmetric form}

\begin{theorem}
The class of bistochastic Hadamard matrices is stable under permuting rows and columns, and under taking tensor products. As examples, we have:
\begin{enumerate}
\item The circulant and symmetric forms $F_N'$ of the Fourier matrices $F_N$.

\item The bistochastic and symmetric forms $F_G'$ of the Fourier matrices $F_G$.

\item The circulant and symmetric Backelin matrices, having size $MN$ with $M|N$.
\end{enumerate}
\end{theorem}

\begin{proof}
In this statement the claim regarding permutations of rows and columns is clear. Assuming now that $H,K$ are bistochastic, with sums $\lambda,\mu$, we have:
\begin{eqnarray*}
\sum_{ia}(H\otimes K)_{ia,jb}
&=&\sum_{ia}H_{ij}K_{ab}\\
&=&\sum_iH_{ij}\sum_aK_{ab}\\
&=&\lambda\mu
\end{eqnarray*}

We have as well the following computation:
\begin{eqnarray*}
\sum_{jb}(H\otimes K)_{ia,jb}
&=&\sum_{jb}H_{ij}K_{ab}\\
&=&\sum_jH_{ij}\sum_bK_{ab}\\
&=&\lambda\mu
\end{eqnarray*}

Thus, the matrix $H\otimes K$ is bistochastic as well. As for the assertions (1,2,3), we already know all this, coming from our study from from chapter 9.
\end{proof}

In the above list of examples, those in (2), which are not necessarily circulant, are the key ones. Indeed, while many interesting complex Hadamard matrices, such as the usual Fourier ones $F_N$, can be put in circulant form, this is something quite exceptional, which does not work any longer when looking for instance at the general Fourier matrices $F_G$. To be more precise, consider a finite abelian group, written as follows:
$$G=\mathbb Z_{N_1}\times\ldots\times\mathbb Z_{N_k}$$

We can then consider the following matrix, with $F_N'$ standing as before for the circulant and symmetric form of the Fourier matrix $F_N$, which is equivalent to $F_G$:
$$F_G'=F_{N_1}'\otimes\ldots\otimes F_{N_k}'$$

Now since the tensor product of circulant matrices is bistochastic, but not necessarily circulant, we can only say that this matrix $F_G'$ is bistochastic, as stated in (2) above.

\bigskip

As a conclusion to all this, the bistochastic complex Hadamard matrices are interesting objects, covering all the generalized Fourier matrices, up to equivalence, and which are definitely worth some study. So, let us develop now some general theory, for such matrices. As a first result, regarding the unitary bistochastic matrices in general, we have:

\begin{proposition}
The real and complex bistochastic groups, which are the sets
$$B_N\subset O_N\quad,\quad 
C_N\subset U_N$$
consisting of matrices which are bistochastic, are isomorphic to $O_{N-1}$, $U_{N-1}$.
\end{proposition}

\begin{proof}
Let us pick a unitary matrix $F\in U_N$ satisfying the following condition, where $e_0,\ldots,e_{N-1}$ is the standard basis of $\mathbb C^N$, and where $\xi$ is the all-one vector:
$$Fe_0=\frac{1}{\sqrt{N}}\xi$$ 

Observe that such matrices $F\in U_N$ exist indeed, the basic example being the normalized Fourier matrix $F_N/\sqrt{N}$. We have then, by using the above property of $F$:
\begin{eqnarray*}
u\xi=\xi
&\iff&uFe_0=Fe_0\\
&\iff&F^*uFe_0=e_0\\
&\iff&F^*uF=diag(1,w)
\end{eqnarray*}

Thus we have isomorphisms as in the statement, given by:
$$w_{ij}\to(F^*uF)_{ij}$$

But this gives both the assertions.
\end{proof}

Now back to the Hadamard matrices, we first have the following elementary result:

\begin{proposition}
For a complex Hadamard matrix $H\in M_N(\mathbb C)$, the following conditions are equivalent:
\begin{enumerate}
\item $H$ is bistochastic, with sums $\lambda$.

\item $H$ is row-stochastic, with sums $\lambda$, and $|\lambda|^2=N$.
\end{enumerate}
\end{proposition}

\begin{proof}
Both the implications are elementary, as follows:

\medskip

$(1)\implies(2)$ If we denote by $H_1,\ldots,H_N\in\mathbb T^N$ the rows of $H$, we have indeed:
\begin{eqnarray*}
N
&=&\sum_i<H_1,H_i>\\
&=&\sum_i\sum_jH_{1j}\bar{H}_{ij}\\
&=&\sum_jH_{1j}\sum_i\bar{H}_{ij}\\
&=&\sum_jH_{1j}\cdot\bar{\lambda}\\
&=&|\lambda|^2
\end{eqnarray*}

$(2)\implies(1)$ Consider the all-one vector $\xi=(1)_i\in\mathbb C^N$. The fact that $H$ is row-stochastic with sums $\lambda$ reads:
\begin{eqnarray*}
\sum_jH_{ij}=\lambda,\forall i
&\iff&\sum_jH_{ij}\xi_j=\lambda\xi_i,\forall i\\
&\iff&H\xi=\lambda\xi
\end{eqnarray*}

Also, the fact that $H$ is column-stochastic with sums $\lambda$ reads:
\begin{eqnarray*}
\sum_iH_{ij}=\lambda,\forall j
&\iff&\sum_jH_{ij}\xi_i=\lambda\xi_j,\forall j\\
&\iff&H^t\xi=\lambda\xi
\end{eqnarray*}

We must prove that the first condition implies the second one, provided that the row sum $\lambda$ satisfies $|\lambda|^2=N$. But this follows from the following computation:
\begin{eqnarray*}
H\xi=\lambda\xi
&\implies&H^*H\xi=\lambda H^*\xi\\
&\implies&N^2\xi=\lambda H^*\xi\\
&\implies&N^2\xi=\bar{\lambda}H^t\xi\\
&\implies&H^t\xi=\lambda\xi
\end{eqnarray*}

Thus, we have proved both the implications, and we are done.
\end{proof}

Here is another basic result, that we will need as well in what follows:

\begin{proposition}
For a complex Hadamard matrix $H\in M_N(\mathbb C)$, and a number $\lambda\in\mathbb C$ satisfying $|\lambda|^2=N$, the following are equivalent:
\begin{enumerate}
\item We have $H\sim H'$, with $H'$ being bistochastic, with sums $\lambda$.

\item $K_{ij}=a_ib_jH_{ij}$ is bistochastic with sums $\lambda$, for some $a,b\in\mathbb T^N$. 

\item The equation $Hb=\lambda\bar{a}$ has solutions $a,b\in\mathbb T^N$. 
\end{enumerate}
\end{proposition}

\begin{proof}
Once again, this is an elementary result, the proof being as follows:

\medskip

$(1)\iff(2)$ Since the permutations of the rows and columns preserve the bistochasticity condition, the equivalence $H\sim H'$ that we are looking for can be assumed to come only from multiplying the rows and columns by numbers in $\mathbb T$. Thus, we are looking for scalars $a_i,b_j\in\mathbb T$ such that the following matrix is bistochastic with sums $\lambda$:
$$K_{ij}=a_ib_jH_{ij}$$

Thus, we are led to the conclusion that (1) and (2) are equivalent, as claimed.

\medskip

$(2)\iff(3)$ The row sums of the matrix $K_{ij}=a_ib_jH_{ij}$ are given by:
$$\sum_jK_{ij}
=\sum_ja_ib_jH_{ij}
=a_i(Hb)_i$$

Thus $K$ is row-stochastic with sums $\lambda$ precisely when $Hb=\lambda\bar{a}$, and by using the equivalence in Proposition 10.5, we obtain the result.
\end{proof}

Finally, here is an extension of the excess inequality from chapter 2:

\index{excess}

\begin{theorem}
For a complex Hadamard matrix $H\in M_N(\mathbb C)$, the excess, 
$$E(H)=\sum_{ij}H_{ij}$$
satisfies $|E(H)|\leq N\sqrt{N}$, with equality precisely when $H$ is bistochastic.
\end{theorem}

\begin{proof}
In terms of the all-one vector $\xi=(1)_i\in\mathbb C^N$, we have:
\begin{eqnarray*}
E(H)
&=&\sum_{ij}H_{ij}\\
&=&\sum_{ij}H_{ij}\xi_j\bar{\xi}_i\\
&=&\sum_i(H\xi)_i\bar{\xi}_i\\
&=&<H\xi,\xi>
\end{eqnarray*}

Now by using the Cauchy-Schwarz inequality, along with the fact that $U=H/\sqrt{N}$ is unitary, and hence of norm 1, we obtain, as claimed:
\begin{eqnarray*}
|E(H)|
&\leq&||H\xi||\cdot||\xi||\\
&\leq&||H||\cdot||\xi||^2\\
&=&N\sqrt{N}
\end{eqnarray*}

Regarding now the equality case, this requires the vectors $H\xi,\xi$ to be proportional, and so our matrix $H$ to be row-stochastic. Now, let us assume:
$$H\xi=\lambda\xi$$

We have then $|\lambda|^2=N$, and by Proposition 10.5 we obtain the result.
\end{proof}

The above result was just an introduction to what can be said about the excess, and we refer to Kharaghani-Seberry \cite{kse} for more on all this. In what concerns us, we will be back to the excess in chapter 11 below, with some probabilistic computations.

\bigskip

Let us go back now to the fundamental question of putting an arbitrary Hadamard matrix in bistochastic form. As already explained in the above, we are interested in solving this question in general, and in particular in the real case, with potential complex reformulations of the HC and CHC, and other real Hadamard questions, at stake. What we know so far on this subject can be summarized as follows:

\begin{proposition}
An Hadamard matrix $H\in M_N(\mathbb C)$ can be put in bistochastic form when one of the following conditions is satisfied:
\begin{enumerate}
\item The equations $|Ha|_i=\sqrt{N}$, with $i=1,\ldots,N$, have solutions $a\in\mathbb T^N$.

\item The quantity $|E|$ attains its maximum $N\sqrt{N}$ over the equivalence class of $H$.
\end{enumerate}
\end{proposition}

\begin{proof}
This follows indeed from Proposition 10.5 and Proposition 10.6, which altogether gives the equivalence between the two conditions in the statement.
\end{proof}

Thus, we have two approaches to the problem, one algebraic, and one analytic. 

\section*{10b. Idel-Wolf theorem}

Let us first discuss the algebraic approach, coming from Proposition 10.8 (1). What we have there is a certain system of $N$ equations, having as unknowns $N$ real variables, namely the phases of $a_1,\ldots,a_N$. This system is highly non-linear, but can be solved, however, via a certain non-explicit method, as explained by Idel-Wolf in \cite{iwo}. In order to discuss this material, which is quite advanced, let us begin with some preliminaries. The complex projective space appears by definition as follows:
$$P^{N-1}_\mathbb C=(\mathbb C^N-\{0\})\big/<x=\lambda y>$$

Inside this projective space, we have the Clifford torus, constructed as follows:
$$\mathbb T^{N-1}=\left\{(z_1,\ldots,z_N)\in P^{N-1}_\mathbb C\Big||z_1|=\ldots=|z_N|\right\}$$

\index{complex projective space}
\index{Clifford torus}

With these conventions, we have the following result, from \cite{iwo}:

\begin{proposition}
For a unitary matrix $U\in U_N$, the following are equivalent:
\begin{enumerate}
\item There exist $L,R\in U_N$ diagonal such that the following matrix is bistochastic:
$$U'=LUR$$ 

\item The standard torus $\mathbb T^N\subset\mathbb C^N$ satisfies:
$$\mathbb T^N\cap U\mathbb T^N\neq\emptyset$$

\item The Clifford torus $\mathbb T^{N-1}\subset P^{N-1}_\mathbb C$ satisfies:
$$\mathbb T^{N-1}\cap U\mathbb T^{N-1}\neq\emptyset$$
\end{enumerate}
\end{proposition}

\begin{proof}
These equivalences are all elementary, as follows:

\medskip

$(1)\implies(2)$ Assuming that $U'=LUR$ is bistochastic, which in terms of the all-1 vector $\xi$ means $U'\xi=\xi$, if we set $f=R\xi\in\mathbb T^N$ we have:
$$Uf
=\bar{L}U'\bar{R}f
=\bar{L}U'\xi
=\bar{L}\xi\in\mathbb T^N$$

Thus we have $Uf\in\mathbb T^N\cap U\mathbb T^N$, which gives the conclusion.

\medskip

$(2)\implies(1)$ Given $g\in\mathbb T^N\cap U\mathbb T^N$, we can define $R,L$ as follows:
$$R=\begin{pmatrix}g_1\\&\ddots\\&&g_N\end{pmatrix}\quad,\quad 
\bar{L}=\begin{pmatrix}(Ug)_1\\&\ddots\\&&(Ug)_N\end{pmatrix}$$

With these values for $L,R$, we have then the following formulae:
$$R\xi=g\quad,\quad 
\bar{L}\xi=Ug$$

Thus the matrix $U'=LUR$ is bistochastic, because:
$$U'\xi
=LUR\xi
=LUg
=\xi$$

$(2)\implies(3)$ This is clear, because $\mathbb T^{N-1}\subset P^{N-1}_\mathbb C$ appears as the projective image of $\mathbb T^N\subset\mathbb C^N$, and so $\mathbb T^{N-1}\cap U\mathbb T^{N-1}$ appears as the projective image of $\mathbb T^N\cap U\mathbb T^N$.

\medskip

$(3)\implies(2)$ We have indeed the following equivalence:
$$\mathbb T^{N-1}\cap U\mathbb T^{N-1}\neq\emptyset
\iff\exists\lambda\neq 0,\lambda\mathbb T^N\cap U\mathbb T^N\neq\emptyset$$

But $U\in U_N$ implies $|\lambda|=1$, and this gives the result.
\end{proof}

The point now is that the condition (3) above is something familiar in symplectic geometry, and known to hold for any $U\in U_N$. Thus, following \cite{iwo}, we have:

\index{Idel-Wolf theorem}
\index{Lagrangian submanifold}
\index{Hamiltonian isotopy}
\index{symplectic manifold}
\index{bistochastic form}
\index{Sinkhorn normal form}
\index{unitary matrix}

\begin{theorem}
Any unitary matrix $U\in U_N$ can be put in bistochastic form,
$$U'=LUR$$
with $L,R\in U_N$ being both diagonal, via a certain non-explicit method.
\end{theorem}

\begin{proof}
As already mentioned, the condition $\mathbb T^{N-1}\cap U\mathbb T^{N-1}\neq\emptyset$ in Proposition 10.9 (3) is something quite natural in symplectic geometry. To be more precise:

\medskip

-- The Clifford torus $\mathbb T^{N-1}\subset P^{N-1}_\mathbb C$ is a Lagrangian submanifold.

\medskip

-- The map $\mathbb T^{N-1}\to U\mathbb T^{N-1}$ is a Hamiltonian isotopy.

\medskip

-- A non-trivial result of Biran-Entov-Polterovich \cite{bep} and Cho \cite{cho} states that $\mathbb T^{N-1}$ cannot be displaced from itself via a Hamiltonian isotopy. 

\medskip

Thus, the results in \cite{bep}, \cite{cho} tells us that  $\mathbb T^{N-1}\cap U\mathbb T^{N-1}\neq\emptyset$ holds indeed, for any $U\in U_N$. We therefore obtain the result, via Proposition 10.9. See Idel-Wolf \cite{iwo}.
\end{proof}

In relation now with our Hadamard matrix questions, we have:

\begin{theorem}
Any complex Hadamard matrix can be put in bistochastic form, up to the standard equivalence relations for such matrices.
\end{theorem}

\begin{proof}
This follows indeed from Theorem 10.10, because if $H=\sqrt{N}U$ is Hadamard then so is $H'=\sqrt{N}U'$, and with the remark that, in what regards the equivalence relation, we just need the multiplication of the rows and columns by scalars in $\mathbb T$.
\end{proof}

There are many further things that can be said here. As explained in \cite{iwo}, the various technical results from \cite{bep}, \cite{cho} show that in the generic, ``transverse'' situation, there are at least $2^{N-1}$ ways of putting a unitary matrix $U\in U_N$ in bistochastic form, and this modulo the obvious transformation $U\to zU$, with $|z|=1$. 

\bigskip

Thus, the question of explicitely putting the Hadamard matrices $H\in M_N(\mathbb C)$ in bistochastic form remains open, and open as well is the question of finding a simpler proof for the fact that this can be done indeed, without using \cite{bep}, \cite{cho}.

\section*{10c. Complex glow}

Regarding the above questions, a possible approach comes from the excess result from Theorem 10.7. Indeed, we know from there that the excess $E(H)=\sum_{ij}H_{ij}$ satisfies the following inequality, with equality precisely when $H$ is bistochastic:
$$|E(H)|\leq N\sqrt{N}$$

Thus, in order to put a complex Hadamard matrix $H\in M_N(\mathbb T)$ in bistochastic form, it is enough to show that the law of $|E|$ over the equivalence class of $H$ has $N\sqrt{N}$ as upper support bound. In order to comment on this, let us first formulate:

\index{excess}
\index{complex glow}

\begin{definition}
The glow of $H\in M_N(\mathbb C)$ is the measure $\mu\in\mathcal P(\mathbb C)$ given by:
$$\int_\mathbb C\varphi(x)d\mu(x)=\int_{\mathbb T^N\times\mathbb T^N}\varphi\left(\sum_{ij}a_ib_jH_{ij}\right)d(a,b)$$
That is, the glow is the law of the following quantity, called excess 
$$E=\sum_{ij}H_{ij}$$
computed over the Hadamard equivalence class of $H$.
\end{definition}

Here $H$ can be any complex matrix, but the equivalence relation is the one for the complex Hadamard matrices. To be more precise, let us call two complex matrices $H,K\in M_N(\mathbb C)$ Hadamard equivalent if one can pass from one to the other by permuting rows and columns, or by multiplying the rows and columns by numbers in $\mathbb T$. Now since permuting rows and columns does not change the quantity $E=\sum_{ij}H_{ij}$, we can restrict attention from the full equivalence group $G=(S_N\rtimes\mathbb T^N)\times(S_N\rtimes\mathbb T^N)$ to the smaller group $G'=\mathbb T^N\times\mathbb T^N$, and we obtain in this way the measure $\mu$ in Definition 10.12.

\bigskip

As in the real case, the terminology comes from a picture of the following type, with the stars $*$ representing the entries of our matrix, and with the switches being supposed now to be continuous, randomly changing the phases of the concerned entries:
\medskip
$$\begin{matrix}
\to&&*&*&*&*\\
\to&&*&*&*&*\\
\to&&*&*&*&*\\
\to&&*&*&*&*\\
\\
&&\uparrow&\uparrow&\uparrow&\uparrow
\end{matrix}$$ 
\smallskip

In short, what we have here is a complex generalization of the Gale-Berlekamp game \cite{fsl}, \cite{rvi}, and this is where a main motivation for studying the glow comes from.

\index{Gale-Berlekamp game}
\index{switching lights}

\bigskip

As a first remark, simplifying our study, exactly as in the real case, we are in fact interested in computing a real measure, due to  the following simple fact: 

\index{multiplicative convolution}
\index{polar decomposition}
\index{rotational invariance}

\begin{proposition}
With $E=\sum_{ij}H_{ij}$, the laws $\mu,\mu^+$ of the variables
$$E,|E|$$
over the torus $\mathbb T^N\times\mathbb T^N$ are related by the formula
$$\mu=\varepsilon\times\mu^+$$
where $\times$ is the multiplicative convolution, and $\varepsilon$ is the uniform measure on $\mathbb T$.
\end{proposition}

\begin{proof}
By definition of the excess $E$, as being the total sum of the entries of the matrix, we have the following equality, valid for any $\lambda\in\mathbb T$:
$$E(\lambda H)=\lambda E(H)$$

We conclude from this that $\mu=law(E)$ is invariant under the action of $\mathbb T$. Thus $\mu$ must decompose as follows, with $\mu^+$ being a certain probability measure on $[0,\infty)$:
$$\mu=\varepsilon\times\mu^+$$

But, according to our definitions, this measure $\mu^+$ is precisely the measure in the statement, that of variable $|E|$, and this gives the result.
\end{proof}

In particular, we can see from the above result that the glow is invariant under rotations. With this observation made, we can formulate the following result:

\index{glow support}

\begin{theorem}
The glow of any Hadamard matrix $H\in M_N(\mathbb C)$, or more generally of any $H\in\sqrt{N}U_N$, satisfies the following conditions, where $\mathbb D$ is the unit disk, 
$$N\sqrt{N}\,\mathbb T\subset supp(\mu)\subset N\sqrt{N}\,\mathbb D$$
with the inclusion on the right coming from Cauchy-Schwarz, and with the inclusion on the left corresponding to the fact that $H$ can be put in bistochastic form.
\end{theorem}

\begin{proof}
We have two inclusions to be proved, the idea being as follows:

\medskip

(1) The inclusion on the right comes indeed from Cauchy-Schwarz, as explained in the proof of Theorem 10.7, with the remark that the computation there only uses the fact that the rescaled matrix $U=H/\sqrt{N}$ is unitary.

\medskip

(2) Regarding now the inclusion on the left, we know from Theorem 10.10 that $H$ can be put in bistochastic form. According to Proposition 10.8, this tells us that we have:
$$N\sqrt{N}\,\mathbb T\cap supp(\mu)\neq\emptyset$$

Now by using the rotational invariance of the glow, and hence of its support, coming from Proposition 10.13, we obtain from this:
$$N\sqrt{N}\,\mathbb T\subset supp(\mu)$$

Thus, we are led to the conclusions in the statement.
\end{proof}

The challenging question now is that of proving the above result, which comes from heavy symplectic geometry, by using standard probabilistic techniques. Indeed, as explained in chapter 9, in the context of the questions investigated there, the support of a real measure can be recaptured from the moments, by computing a limit. Thus, knowing the moments of the glow well enough would solve the problem. 

\bigskip

Regarding these moments, the general formula is as follows:

\index{excess moments}
\index{glow moments}

\begin{proposition}
For $H\in M_N(\mathbb T)$ the even moments of $|E|$ are given by
$$\int_{\mathbb T^N\times\mathbb T^N}|E|^{2p}=\sum_{[i]=[k],[j]=[l]}\frac{H_{i_1j_1}\ldots H_{i_pj_p}}{H_{k_1l_1}\ldots H_{k_pl_p}}$$
where the sets between brackets are by definition sets with repetition.
\end{proposition}

\begin{proof}
We have indeed the following computation:
\begin{eqnarray*}
\int_{\mathbb T^N\times\mathbb T^N}|E|^{2p}
&=&\int_{\mathbb T^N\times\mathbb T^N}\Big|\sum_{ij}H_{ij}a_ib_j\Big|^{2p}\\
&=&\int_{\mathbb T^N\times\mathbb T^N}\left(\sum_{ijkl}\frac{H_{ij}}{H_{kl}}\cdot\frac{a_ib_j}{a_kb_l}\right)^p\\
&=&\sum_{ijkl}\frac{H_{i_1j_1}\ldots H_{i_pj_p}}{H_{k_1l_1}\ldots H_{k_pl_p}}\int_{\mathbb T^N}\frac{a_{i_1}\ldots a_{i_p}}{a_{k_1}\ldots a_{k_p}}\int_{\mathbb T^N}\frac{b_{j_1}\ldots b_{j_p}}{b_{l_1}\ldots b_{l_p}}
\end{eqnarray*}

Now since the integrals at right equal respectively the Kronecker symbols $\delta_{[i],[k]}$ and $\delta_{[j],[l]}$, we are led to the formula in the statement.
\end{proof}

With this formula in hand, the main result, regarding the fact that the complex Hadamard matrices can be put in bistochastic form, reformulates as follows:

\begin{theorem}
For a complex Hadamard matrix $H\in M_N(\mathbb T)$ we have
$$\lim_{p\to\infty}\left(\sum_{[i]=[k],[j]=[l]}\frac{H_{i_1j_1}\ldots H_{i_pj_p}}{H_{k_1l_1}\ldots H_{k_pl_p}}\right)^{1/p}=N^3$$
coming from the fact that $H$ can be put in bistochastic form.
\end{theorem}

\begin{proof}
This follows from the well-known fact that the maximum of a bounded function $\Theta:X\to[0,\infty)$ can be recaptured via following formula:
$$\max(\Theta)=\lim_{p\to\infty}\left(\int_X\Theta(x)^p\,dx\right)^{1/p}$$

We can use this estimate for the following function, over $X=\mathbb T^N\times\mathbb T^N$:
$$\Theta=|E|^2$$

We conclude that the limit in the statement is the square of the upper bound of the glow. But, according to Theorem 10.14, this upper bound is known to be $\leq N^3$ by Cauchy-Schwarz, and the equality holds by the results in \cite{iwo}.
\end{proof}

To conclude now, the challenging question is that of finding a direct proof for Theorem 10.16. All this would provide an alternative aproach to the results in \cite{iwo}, which would be of course still not explicit, but which would use at least some more familiar tools. We will discuss such questions in chapter 11 below, with the remark however that the problems at $N\in\mathbb N$ fixed being quite difficult, we will do a $N\to\infty$ study only.

\section*{10d. Fourier matrices}

Getting away now from these difficult questions, we have nothing concrete so far, besides the list of examples from Theorem 10.3, coming from the circulant matrix considerations in chapter 9. So, our purpose will be that of extending that list. A first natural question is that of looking at the Butson matrix case. To start with, we have:

\index{bistochastic Butson matrix}

\begin{proposition}
Assuming that the Butson class $H_N(l)$ contains a bistochastic matrix, the equations
\begin{eqnarray*}
a_0+a_1+\ldots+a_{l-1}&=&N\\
|a_0+a_1w+\ldots+a_{l-1}w^{l-1}|^2&=&N
\end{eqnarray*}
must have solutions, over the positive integers.
\end{proposition}

\begin{proof}
This is a reformulation of the following equality, from Proposition 10.5, regarding the row sums of a bistochastic Hadamard matrix:
$$|\lambda|^2=N$$

Indeed, if we set $w=e^{2\pi i/l}$, and we denote by $a_i\in\mathbb N$ the number of $w^i$ entries appearing in the first row of our matrix, then the row sum of the matrix is given by:
$$\lambda=a_0+a_1w+\ldots+a_{l-1}w^{l-1}$$

Thus, we obtain the system of equations in the statement.
\end{proof}

The point now is that, in practice, we are led precisely to the Turyn obstructions from chapter 9. At small values of $l$, the obstructions are as follows:

\index{Turyn obstruction}

\begin{theorem}
Assuming that $H_N(l)$ contains a bistochastic matrix, the following equations must have solutions, over the integers:
\begin{enumerate}
\item $l=2$: $4n^2=N$.

\item $l=3$: $x^2+y^2+z^2=2N$, with $x+y+z=0$.

\item $l=4$: $a^2+b^2=N$.
\end{enumerate}
\end{theorem}

\begin{proof}
This follows indeed from the results that we have:

\medskip

(1) This is something well-known, which follows from Proposition 10.17.

\medskip

(2) This is best viewed by using Proposition 10.17, and the following formula, that we already know, from chapter 5 above:
$$\left|a+bw+cw^2\right|^2=\frac{1}{2}[(a-b)^2+(b-c)^2+(c-a)^2]$$

At the level of the concrete obstructions, we must have for instance $5\!\!\not|N$. Indeed, this follows as in the proof of the de Launey obstruction for $H_N(3)$ with $5|N$.

\medskip

(3) This follows again from Proposition 10.17, and from $|a+ib|^2=a^2+b^2$.
\end{proof}

As a conclusion, nothing much interesting is going on in the Butson matrix case, with various arithmetic obstructions, that we partly already met, appearing here. In order to reach, however, to a number of positive results, beyond those in Theorem 10.3, we can investigate various special classes of matrices, such as the Di\c t\u a products. In order to formulate our results, we will use the following notion:

\index{almost bistochastic form}

\begin{definition}
We say that a complex Hadamard matrix $H\in M_N(\mathbb C)$ is in ``almost bistochastic form'' when all the row sums belong to $\sqrt{N}\cdot\mathbb T$.
\end{definition}

Observe that, assuming that this condition holds, the matrix $H$ can be put in bistochastic form, just by multiplying its rows by suitable numbers from $\mathbb T$. We will be particularly interested here in the special situation where the affine deformations $H^q\in M_N(\mathbb C)$ of a given complex Hadamard matrix $H\in M_N(\mathbb C)$ can be put in almost bistochastic form, independently of the value of the parameter $q$. For the simplest deformations, namely those of $F_2\otimes F_2$, this is indeed the case, as shown by the following result:

\begin{proposition}
The deformations of $F_2\otimes F_2$, with parameter matrix $Q=(^p_r{\ }^q_s)$,
$$F_2\otimes_QF_2=
\begin{pmatrix}
p&q&p&q\\
p&-q&p&-q\\
r&s&-r&-s\\ 
r&-s&-r&s
\end{pmatrix}$$
can be put in almost bistochastic form, independently of the value of $Q$.
\end{proposition}

\begin{proof}
By multiplying the columns of the matrix in the statement with $1,1,-1,1$ respectively, we obtain the following matrix:
$$F_2\otimes''_QF_2=
\begin{pmatrix}
p&q&-p&q\\
p&-q&-p&-q\\
r&s&r&-s\\ 
r&-s&r&s
\end{pmatrix}$$

The row sums of this matrix are as follows:
$$2q,-2q,2r,2r\in2\mathbb T$$

Thus, by multiplying by suitable scalars, namely the complex conjugates of these numbers, we can put our matrix in bistochastic form, as desired.
\end{proof}

We will see later that $F_2\otimes''_QF_2$ is equivalent to a certain matrix $F_2\otimes'F_2$, which is part of a series $F_N\otimes'F_N$. Now back to the general case, we have:

\begin{theorem}
A deformed tensor product $H\otimes_QK$ can be put in bistochastic form when there exist numbers $x^i_a\in\mathbb T$ such that with
$$G_{ib}=\frac{(K^*x^i)_b}{Q_{ib}}$$
we have $|(H^*G)_{ib}|=\sqrt{MN}$, for any $i,b$.
\end{theorem}

\begin{proof}
According to our tensor product conventions, the deformed tensor product $L=H\otimes_QK$ is given by the following formula:
$$L_{ia,jb}=Q_{ib}H_{ij}K_{ab}$$

By multiplying the columns by scalars $R_{jb}\in\mathbb T$, this matrix becomes:
$$L'_{ia,jb}=R_{jb}Q_{ib}H_{ij}K_{ab}$$

The row sums of this matrix are given by:
\begin{eqnarray*}
S_{ia}'
&=&\sum_{jb}R_{jb}Q_{ib}H_{ij}K_{ab}\\
&=&\sum_bK_{ab}Q_{ib}\sum_jH_{ij}R_{jb}\\
&=&\sum_bK_{ab}Q_{ib}(HR)_{ib}
\end{eqnarray*}

Consider now the following variables:
$$C^i_b=Q_{ib}(HR)_{ib}$$

In terms of these variables, the rows sums are given by:
$$S_{ia}'
=\sum_bK_{ab}C^i_b
=(KC^i)_a$$

Thus $H\otimes_QK$ can be put in bistochastic form when we can find scalars $R_{jb}\in\mathbb T$ and $x^i_a\in\mathbb T$ such that, with  $C^i_b=Q_{ib}(HR)_{ib}$, the following condition is satisfied:
$$(KC^i)_a=\sqrt{MN}x^i_a\quad,\quad\forall i,a$$

But this condition is equivalent to the following condition:
$$KC^i=\sqrt{MN}x^i\quad,\quad\forall i$$

Now by multiplying to the left by $K^*$, we are led to the following condition:
$$\sqrt{N}C^i=\sqrt{M}K^*x^i\quad,\quad\forall i$$

Now by recalling that $C^i_b=Q_{ib}(HR)_{ib}$, this condition is equivalent to:
$$\sqrt{N}Q_{ib}(HR)_{ib}=\sqrt{M}(K^*x^i)_b\quad,\quad\forall i,b$$

Consider now the variables in the statement, namely:
$$G_{ib}=\frac{(K^*x^i)_b}{Q_{ib}}$$

In terms of these variables, the above condition reads:
$$\sqrt{N}(HR)_{ib}=\sqrt{M}G_{ib}\quad,\quad\forall i,b$$

But this condition is equivalent to:
$$\sqrt{N}HR=\sqrt{M}G$$

Now by multiplying to the left by $H^*$, we are led to the following condition:
$$\sqrt{MN}R=H^*G$$

Thus, we have obtained the condition in the statement.
\end{proof}

As an illustration for the above result, assume that $H,K$ can be put in bistochastic form, by using vectors $y\in\mathbb T^M,z\in\mathbb T^N$, and let us set:
$$x^i_a=y_iz_a$$

Then with the choice $Q=1$ for our parameter matrix, we have:
\begin{eqnarray*}
G_{ib}
&=&(K^*x^i)_b\\
&=&[K^*(y_iz)]_b\\
&=&y_i(K^*z)_b
\end{eqnarray*}

We therefore obtain the following formula:
\begin{eqnarray*}
(H^*G)_{ib}
&=&\sum_j(H^*)_{ij}G_{jb}\\
&=&\sum_j(H^*)_{ij}y_j(K^*z)_b\\
&=&(H^*y)_i(K^*z)_b
\end{eqnarray*}

Thus the usual tensor product $H\otimes K$ can be put in bistochastic form as well, which is of course something that we already know, from the above. Now back to the general case, that of the arbitrary Di\c t\u a deformations in Theorem 10.21, the point is that in the particular case $H=F_M$ the equations simplify, and we have the following result:

\begin{proposition}
A deformed tensor product $F_M\otimes_QK$ can be put in bistochastic form when there exist numbers $x^i_a\in\mathbb T$ such that with
$$G_{ib}=\frac{(K^*x^i)_b}{Q_{ib}}$$
we have the following formulae, with $l$ being taken modulo $M$:
$$\sum_jG_{jb}\bar{G}_{j+l,b}=MN\delta_{l,0}\quad,\quad\forall l,b$$
Moreover, the $M\times N$ matrix $|G_{jb}|^2$ is row-stochastic with sums $N^2$, and the $l=0$ equations state that this matrix must be column-stochastic, with sums $MN$.
\end{proposition}

\begin{proof}
With notations from Theorem 10.21, and with $w=e^{2\pi i/M}$, we have:
$$(H^*G)_{ib}=\sum_jw^{-ij}G_{jb}$$

The absolute value of this number can be computed as follows:
\begin{eqnarray*}
|(H^*G)_{ib}|^2
&=&\sum_{jk}w^{i(k-j)}G_{jb}\bar{G}_{kb}\\
&=&\sum_{jl}w^{il}G_{jb}\bar{G}_{j+l,b}\\
&=&\sum_lw^{il}\sum_jG_{jb}\bar{G}_{j+l,b}
\end{eqnarray*}

If we denote by $v^b_l$ the sum on the right, we obtain:
$$|(H^*G)_{ib}|^2
=\sum_lw^{il}v^b_l
=(F_Mv^b)_i$$

Now if we denote by $\xi$ the all-one vector in $\mathbb C^M$, the condition $|(H^*G)_{ib}|=\sqrt{MN}$ for any $i,b$ found in Theorem 10.21 reformulates as follows:
$$F^Mv^b=MN\xi\quad,\quad\forall b$$

By multiplying to the left by $F_M^*/M$, this condition is equivalent to:
$$v^b
=NF_M^*\xi
=\begin{pmatrix}MN\\0\\ \vdots\\0\end{pmatrix}$$

Let us examine the first equation, $v^b_0=MN$. By definition of $v^b_l$, we have:
$$v^b_0
=\sum_jG_{jb}\bar{G}_{jb}
=\sum_j|G_{jb}|^2$$

Now recall from Theorem 10.21 that we have, for certain numbers $x^j_b\in\mathbb T$:
$$G_{jb}=\frac{(K^*x^j)_b}{Q_{jb}}$$

Since we have $Q_{jb}\in\mathbb T$ and $K^*/\sqrt{N}\in U_N$, we obtain:
\begin{eqnarray*}
\sum_b|G_{jb}|^2
&=&\sum_b|(K^*x^j)_b|^2\\
&=&||K^*x^j||_2^2\\
&=&N||x^j||_2^2\\
&=&N^2
\end{eqnarray*}

Thus the $M\times N$ matrix $|G_{jb}|^2$ is row-stochastic, with sums $N^2$, and our equations $v^b_0=MN$ for any $b$ state that this matrix must be column-stochastic, with sums $MN$. Regarding now the other equations that we found, namely $v^b_l=0$ for $l\neq0$, by definition of $v^b_l$ and of the variables $G_{jb}$, these state that we must have:
$$\sum_jG_{jb}\bar{G}_{j+l,b}=0\quad,\quad\forall l\neq0,\forall b$$

Thus, we are led to the conditions in the statement.
\end{proof}

As an illustration for this result, let us go back to the $Q=1$ situation, explained after Theorem 10.21. By using the formula $G_{ib}=y_i(K^*z)_b$ there, we have:
\begin{eqnarray*}
\sum_jG_{jb}\bar{G}_{j+l,b}
&=&\sum_jy_j(K^*z)_b\,\overline{y}_{j+l}\overline{(K^*z)_b}\\
&=&|(K^*z)_b|^2\sum_j\frac{y_j}{y_{j+l}}\\
&=&M\cdot N\delta_{l,0}
\end{eqnarray*}

Thus, if $K$ can be put in bistochastic form, then so can be put $F_M\otimes K$. As a second illustration now, let us go back to the matrices $F_2\otimes'_QF_2$ from the proof of Proposition 10.20. For these matrices, the vector of the row sums is as follows: 
$$S=(2q,-2q,2r,2r)$$

Thus, with the above notations, we have the following formula:
$$x=(q,-q,r,r)$$

We therefore obtain the following formulae for the upper entries of $G$:
$$G_{0b}
=\frac{\left[\begin{pmatrix}1&1\\1&-1\end{pmatrix}\begin{pmatrix}q\\-q\end{pmatrix}\right]_b}{Q_{0b}}
=\frac{\begin{pmatrix}0\\2q\end{pmatrix}_b}{Q_{0b}}$$

As for the lower entries of $G$, these are as follows:
$$G_{1b}
=\frac{\left[\begin{pmatrix}1&1\\1&-1\end{pmatrix}\begin{pmatrix}r\\r\end{pmatrix}\right]_b}{Q_{1b}}
=\frac{\begin{pmatrix}2r\\0\end{pmatrix}_b}{Q_{1b}}$$

Thus, in this case the matrix $G$ is as follows, independently of $Q$:
$$G=\begin{pmatrix}0&2\\2&0\end{pmatrix}$$

In particular, we see that the conditions in Proposition 10.22 are satisfied. Now back to the general case, as a main application of our results so far, we have:

\index{deformed Fourier matrix}

\begin{theorem}
The Di\c t\u a deformations of tensor squares of Fourier matrices,
$$F_N\otimes_QF_N$$
can be put in almost bistochastic form, independently of the value of $Q\in M_N(\mathbb T)$.
\end{theorem}

\begin{proof}
We use Proposition 10.22, with $M=N$, and with $K=F_N$. Let $w=e^{2\pi i/N}$, and consider the vectors $x^i\in\mathbb T^N$ given by: 
$$x^i=(w^{(i-1)a})_a$$

Since $K^*K=N1_N$, and $x^i$ are the column vectors of $K$, shifted by 1, we have:
$$K^*x^0=\begin{pmatrix}0\\0\\ \vdots\\0\\N\end{pmatrix}\quad,\quad
K^*x^1=\begin{pmatrix}N\\0\\ \vdots\\0\\0\end{pmatrix}\quad,\ \ldots\ ,\quad
K^*x^{N-1}=\begin{pmatrix}0\\0\\ \vdots\\N\\0\end{pmatrix}$$

We conclude that we have the following formula:
$$(K^*x^i)_b=N\delta_{i-1,b}$$

Thus the matrix $G$ is given by the following formula:
$$G_{ib}=\frac{N\delta_{i-1,b}}{Q_{ib}}$$

With this formula in hand, the sums in Proposition 10.22 are given by:
$$\sum_jG_{jb}\bar{G}_{j+l,b}=\sum_j\frac{N\delta_{j-1,b}}{Q_{jb}}\cdot\frac{N\delta_{j+l-1,b}}{Q_{j+l,b}}$$

In the case $l\neq0$ we clearly get $0$, because the products of Kronecker symbols are $0$. In the case $l=0$ the denominators are $|Q_{jb}|^2=1$, and we obtain:
$$\sum_jG_{jb}\bar{G}_{jb}
=N^2\sum_j\delta_{j-1,b}
=N^2$$

Thus, the conditions in Proposition 10.21 are satisfied, and we obtain the result.
\end{proof}

In relation with the various questions raised above, regarding the Di\c t\u a deformations of the Fourier matrices, this is best result that we have, so far. Here is an equivalent formulation of the above result, which is quite useful, in practice:

\begin{theorem}
The matrix $F_N\otimes'_QF_N$, with $Q\in M_N(\mathbb T)$, defined by
$$(F_N\otimes'_QF_N)_{ia,jb}=\frac{w^{ij+ab}}{w^{bj+j}}\cdot\frac{Q_{ib}}{Q_{b+1,b}}$$
where $w=e^{2\pi i/N}$ is almost bistochastic, and equivalent to $F_N\otimes_QF_N$.
\end{theorem}

\begin{proof}
Our claim is that this is the matrix constructed in the proof of Theorem 10.23. Indeed, let us first go back to the proof of Theorem 10.21. In the case $M=N$ and $H=K=F_N$, the Di\c t\u a deformation $L=H\otimes_QK$ studied there is given by:
$$L_{ia,jb}
=Q_{ib}H_{ij}K_{ab}
=w^{ij+ab}Q_{ib}$$

As explained in the proof of Theorem 10.23, if the conditions in the statement there are satisfied, then the matrix $L_{ia,jb}'=R_{jb}L_{ia,jb}$ is almost bistochastic, where:
$$\sqrt{MN}\cdot R=H^*G$$

In our case now, $M=N$ and $H=K=F_N$, we know from the proof of Proposition 10.22 that the choice of $G$ which makes work Theorem 10.23 is as follows:
$$G_{ib}=\frac{N\delta_{i-1,b}}{Q_{ib}}$$

With this formula in hand, we can compute the matrix $R$, as follows:
\begin{eqnarray*}
R_{jb}
&=&\frac{1}{N}(H^*G)_{jb}\\
&=&\frac{1}{N}\sum_iw^{-ij}G_{ib}\\
&=&\sum_iw^{ij}\cdot\frac{\delta_{i-1,b}}{Q_{ib}}\\
&=&\frac{w^{-(b+1)j}}{Q_{b+1,b}}
\end{eqnarray*}

Thus, the modified version of $F_N\otimes_QF_N$ which is almost bistochastic is given by:
\begin{eqnarray*}
L_{ia,jb}'
&=&R_{jb}L_{ia,jb}\\
&=&\frac{w^{-(b+1)j}}{Q_{b+1,b}}\cdot w^{ij+ab}Q_{ib}\\
&=&\frac{w^{ij+ab}}{w^{bj+j}}\cdot\frac{Q_{ib}}{Q_{b+1,b}}
\end{eqnarray*}

Thus we have obtained the formula in the statement, and we are done.
\end{proof}

As an illustration, let us work out the case $N=2$. Here the root of unity is $w=-1$. Let us denote the deformation matrix as follows:
$$Q=\begin{pmatrix}p&q\\ r&s\end{pmatrix}$$

With the notations $u=p/r$, $v=s/q$, we obtain the following matrix:
\begin{eqnarray*}
F_2\otimes_QF_2
&=&\begin{pmatrix}
\frac{p}{r}&\frac{q}{q}&-\frac{p}{r}&\frac{q}{q}\\
\frac{p}{r}&-\frac{q}{q}&-\frac{p}{r}&-\frac{q}{q}\\
\frac{r}{r}&\frac{s}{q}&\frac{r}{r}&-\frac{s}{q}\\
\frac{r}{r}&-\frac{s}{q}&\frac{r}{r}&\frac{s}{q}
\end{pmatrix}\\
&=&\begin{pmatrix}
u&1&-u&1\\
u&-1&-u&-1\\
1&v&1&-v\\
1&-v&1&v
\end{pmatrix}
\end{eqnarray*}

In general, the question of putting the Di\c t\u a deformations of the tensor products in explicit bistochastic form remains open. Open as well is the question of putting the arbitrary affine deformations of the Fourier matrices in explicit bistochastic form.

\bigskip

We would like to end this chapter by discussing a related interesting question, which can serve as a very good motivation for all this, namely the question on whether the real Hadamard matrices, $H\in M_N(\pm1)$, can be put or not in bistochastic form, in an explicit way. This is certainly true for the Walsh matrices, as explained before, but for the other basic examples, such as the Paley or the Williamson matrices, no results seem to be known so far. Having such a theory would be potentially very interesting, with a complex reformulation of the HC and of the other real Hadamard questions at stake. 

\bigskip

We already know that we are done with the case $N\leq8$. The next problem regards the Paley matrix at $N=12$, which is the unique real Hadamard matrix there:
$$P_{12}\sim P_{12}^1\sim P_{12}^2$$

This matrix is as follows, with the $\pm$ signs standing for $\pm1$ entries:
$$P_{12}=\left(
\begin{array}{ccccccccccccccc}
+&+&+&+&&-&+&+&+&&-&+&+&+\\
+&-&+&-&&+&-&+&+&&+&-&+&+\\
+&+&-&-&&+&+&-&+&&+&+&-&+\\
+&-&-&+&&+&+&+&-&&+&+&+&-\\
\\
-&+&-&-&&+&-&+&+&&-&+&+&-\\
-&+&+&+&&+&+&-&+&&+&-&+&-\\
-&-&-&+&&+&+&+&+&&-&-&-&+\\
+&+&+&-&&+&+&+&-&&-&-&-&-\\
\\
-&+&+&+&&+&-&+&-&&+&+&-&+\\
+&+&-&+&&+&-&-&-&&-&-&+&+\\
+&-&+&+&&+&-&-&+&&-&+&-&-\\
-&-&+&-&&+&+&-&-&&-&+&+&+
\end{array}\right)$$

\index{Paley matrix}

This matrix cannot be put of course in real bistochastic form, its size being not of the form $N=4n^2$. Nor can it be put in bistochastic form over $\{\pm1,\pm i\}$, because the Turyn obstruction for matrices over $\{\pm1,\pm i\}$  is $N=a^2+b^2$, and we have:
$$12\neq a^2+b^2$$

However, the question of putting $P_{12}$ in bistochastic form over the 3-roots of unity makes sense, because the Turyn obstruction here is:
$$x+y+z=0\quad,\quad 
x^2+y^2+z^2=2N$$

And, we do have solutions to these equations at $N=12$, as follows:
$$4^2+(-2)^2+(-2)^2=24$$

Another question is whether $P_{12}$ can be put in bistochastic form over the 8-roots of unity. In order to comment on this, let us first work out the Turyn obstruction, for the bistochastic matrices having as entries the 8-roots of unity. The result is as follows:

\begin{proposition}
The Turyn obstruction for the bistochastic matrices having as entries the $8$-roots of unity is
$$x^2+y^2+z^2+t^2=N\quad,\quad 
xy+yz+zt=xt$$
which must hold for certain numbers $x,y,z,t\in\mathbb Z$.
\end{proposition}

\begin{proof}
The 8-roots of unity are as follows, with $w=e^{\pi i/4}$:
$$1,w,i,iw,-1,-w,-i,-iw$$

Thus, we are led to an equation as follows, with $x,y,z,t\in\mathbb Z$:
$$\left|x+wy+iz+iwt\right|^2=N$$

We have the following computation:
\begin{eqnarray*}
\left|x+wy+iz+iwt\right|^2
&=&(x+wy+iz+iwt)(x-iwy-iz-wt)\\
&=&x^2+y^2+z^2+t^2+w(1-i)(xy+yz+zt-xt)\\
&=&x^2+y^2+z^2+t^2-\sqrt{2}(xy+yz+zt-xt)
\end{eqnarray*}

Thus, we are led to the conclusion in the statement.
\end{proof}

In relation with the above, the point now is that the equations in Proposition 10.25 do have solutions at $N=12$, namely:
$$x=0,y=2,z=-2,t=\pm2$$

Summarizing, the Paley matrix $P_{12}$ cannot be put in bistochastic form over the 4-roots, but the question makes sense over the 3-roots, and over the 8-roots. However, the computations here are not exactly trivial, and the answer is not known.

\bigskip

There are many interesting questions here, and as already mentioned above, the interest in this subject comes from the fact that all this can potentially lead to a complex reformulation of the HC and of the other real Hadamard matrix questions. 

\index{HC}
\index{CHC}

\section*{10e. Exercises}

The material in the present chapter has often gone into research matters, and our exercises here will be of the same type, more difficult than usual. First, we have:

\begin{exercise}
Learn more about the real and complex bistochastic groups $B_N,C_N$, and write down a brief account of what you learned.
\end{exercise}

To be more precise, we have already seen in the above that $B_N,C_N$ are isomorphic respectively to $O_{N-1},U_{N-1}$, via a Fourier transform type operation. However, there are many other interesting things which can be said about $B_N,C_N$, which can be potentially useful in connection with our Hadamard matrix problems, and it is up to you here to check the literature, and learn what can be potentially good to know.

\begin{exercise}
Check the symplectic geometry literature, and write down a concise proof for the Idel-Wolf theorem, based on that, by explaining the main ideas involved.
\end{exercise} 

An even better question would be of course that of writing down a concise proof for the Idel-Wolf theorem, in the rescaled complex Hadamard matrix case, that we are interested in here. We do not know if this is really possible, in the sense that if the Hadamard matrix assumption can really bring some simplifications. Bonus question.

\begin{exercise}
Find the best bound for the support of the glow of the complex Hadamard matrices, by using the moment method, and combinatorics.
\end{exercise}

As with the previous exercise, this is rather a research question.

\begin{exercise}
Study the deformations of the Fourier matrix $F_6$, with the aim of putting them in bistochastic form, and write down what you found.
\end{exercise}

To be more precise here, we know from the above that the deformations of the tensor products of type $F_N\otimes F_N$ can be put in bistochastic form, and in order to get beyond this, the case of the matrices $F_N\otimes F_M$ with $M\neq N$, which numerically starts with the case of the matrix $F_6=F_2\otimes F_3=F_3\otimes F_2$, is the one to be investigated first.

\begin{exercise}
Study the Paley matrix $P_{12}$, with the aim of putting it in bistochastic form, over the complex numbers, and write down what you found.
\end{exercise}

And this is all we have. Only research exercises for this chapter. Sorry for this, and enjoy. Working on difficult exercises can be more fun than working on easy ones, and in any case, any type of work always leads to ``things'', that can be written down.

\chapter{Glow computations}

\section*{11a. Basic results}

We discuss here the computation of the glow of the complex Hadamard matrices, as a continuation of the material from chapter 2, where we discussed the basics of the glow in the real case, and as a continuation as well of the material from chapter 10. As a first motivation for all this, we have the Gale-Berlekamp game \cite{fsl}, \cite{rvi}. Another motivation comes from the questions regarding the bistochastic matrices, in relation with the Ideal-Wolf theorem \cite{iwo}, explained in chapter 10. Finally, we have the question of connecting the defect, and other invariants of the Hadamard matrices, to the glow. 

\bigskip

Let us begin by reviewing the few theoretical things that we know about the glow, from chapter 10. The main results there can be summarized as follows:

\index{excess}
\index{glow}
\index{glow support}

\begin{theorem}
The glow of $H\in M_N(\mathbb C)$, which is the law $\mu\in\mathcal P(\mathbb C)$ of the excess
$$E=\sum_{ij}H_{ij}$$
over the Hadamard equivalence class of $H$, has the following properties:
\begin{enumerate}

\item $\mu=\varepsilon\times\mu^+$, where $\mu^+=law(|E|)$.

\item $\mu$ is invariant under rotations.

\item $H\in\sqrt{N}U_N$ implies $supp(\mu)\subset N\sqrt{N}\,\mathbb D$. 

\item $H\in\sqrt{N}U_N$ implies as well $N\sqrt{N}\,\mathbb T\subset supp(\mu)$.
\end{enumerate}
\end{theorem}

\begin{proof}
We already know all this from chapter 10, the idea being as follows:

\medskip

(1) This follows indeed by using $H\to zH$ with $|z|=1$.

\medskip

(2) This follows from (1), the convolution with $\varepsilon$ bringing the invariance.

\medskip

(3) This follows indeed from Cauchy-Schwarz.

\medskip

(4) This is something highly non-trivial, coming from \cite{iwo}.
\end{proof}

In what follows we will be mainly interested in the Hadamard matrix case, but since the computations here are quite difficult, let us begin our study with other matrices. It is convenient to normalize our matrices, as to make them a bit similar to the complex Hadamard ones. To be more precise, consider the $2$-norm on the vector space of the complex $N\times N$ matrices, which is given by the following formula:
$$||H||_2=\sqrt{\sum_{ij}|H_{ij}|^2}$$ 

We will assume in what follows, by multiplying our matrix $H\in M_N(\mathbb C)$ by a suitable scalar, that this norm takes the same value as for the Hadamard matrices, namely:
$$||H||_2=N$$

We know from chapter 2 that in the real case, the real glow is asymptotically Gaussian. In the complex matrix case, we will reach to the conclusion that the glow is asymptotically complex Gaussian, with the complex Gaussian distribution being as follows:

\index{complex normal variable}
\index{complex Gaussian variable}

\begin{proposition}
The complex Gaussian distribution $\mathcal C$ is the law of the variable
$$z=\frac{1}{\sqrt{2}}(x+iy)$$
with $x,y$ being independent standard Gaussian variables. We have
$$\mathbb E(|z|^{2p})=p!$$
and this moment formula, along with rotational invariance, determines $\mathcal C$.
\end{proposition}

\begin{proof}
This is standard probability theory, with the main result, namely the moment formula in the statement, coming from some routine computations. For more on all this, we refer to any standard probability book, such as Durrett \cite{dur}.
\end{proof}

Finally, we use in what follows the symbol $\sim$ to denote an equality of distributions. With these conventions, we have the following result, to start with:

\index{flat matrix}

\begin{proposition}
We have the following computations:
\begin{enumerate}
\item For the rescaled identity $\widetilde{I}_N=\sqrt{N}I_N$ we have 
$$E\sim\sqrt{N}(q_1+\ldots +q_N)$$
with $q\in\mathbb T^N$ random. With $N\to\infty$ we have $E/N\sim\mathcal C$.

\item For the flat matrix $J_N=(1)_{ij}$ we have 
$$E\sim(a_1+\ldots+a_N)(b_1+\ldots+b_N)$$
with $(a,b)\in\mathbb T^N\times\mathbb T^N$ random. With $N\to\infty$ we have $E/N\sim\mathcal C\times\mathcal C$.
\end{enumerate}
\end{proposition}

\begin{proof}
We use Theorem 11.1, and the moment method:

\medskip

(1) Here we have $E=\sqrt{N}\sum_{i}a_ib_i$, with $a,b\in\mathbb T^N$ random. With $q_i=a_ib_i$ this gives the first assertion. Let us estimate now the moments of $|E|^2$. We have:
\begin{eqnarray*}
\int_{\mathbb T^N\times\mathbb T^N}|E|^{2p}
&=&N^p\int_{\mathbb T^N}|q_1+\ldots+q_N|^{2p}dq\\
&=&N^p\int_{\mathbb T^N}\sum_{ij}\frac{q_{i_1}\ldots q_{i_p}}{q_{j_1}\ldots q_{j_p}}\,dq\\
&=&N^p\#\left\{(i,j)\in\{1,\ldots,N\}^p\times\{1,\ldots,N\}^p\Big|[i_1,\ldots,i_p]=[j_1,\ldots,j_p]\right\}\\
&\simeq&N^p\cdot p!N(N-1)\ldots(N-p+1)\\
&\simeq&N^p\cdot p!N^p\\
&=&p!N^{2p}
\end{eqnarray*}

Here, and in what follows, the sets between brackets are by defintion sets with repetition, and the middle estimate comes from the fact that, with $N\to\infty$, only the multi-indices $i=(i_1,\ldots,i_p)$ having distinct entries contribute. But this gives the result.

\medskip

(2) Here we have the following formula, which gives the first assertion:
$$E
=\sum_{ij}a_ib_j
=\sum_ia_i\sum_jb_j$$

Now since $a,b\in\mathbb T^N$ are independent, so are the quantities $\sum_ia_i,\sum_jb_j$, so we have:
$$\int_{\mathbb T^N\times\mathbb T^N}|E|^{2p}
=\left(\int_{\mathbb T^N}|q_1+\ldots+q_N|^{2p}dq\right)^2
\simeq(p!N^{p})^2$$

Here we have used the estimate in the proof of (1), and this gives the result.
\end{proof}

As a conclusion, the glow is intimately related to the basic hypertoral law, namely the law of the variable $q_1+\ldots+q_N$, with $q\in\mathbb T^N$ being random. Observe that at $N=1$ this hypertoral law is the Dirac mass $\delta_1$, and that at $N=2$ we obtain the following law:
\begin{eqnarray*}
law|1+q|
&=&law\sqrt{(1+e^{it})(1+e^{-it})}\\
&=&law\sqrt{2+2\cos t}\\
&=&law\left(2\cos\frac{t}{2}\right)
\end{eqnarray*}

In general, the law of $\sum q_i$ is known to be related to the P\'olya random walk \cite{pol}. Also, as explained for instance in chapter 9, the moments of this law are:
$$\int_{\mathbb T^N}|q_1+\ldots+q_N|^{2p}dq=\sum_{\pi\in P(p)}\binom{p}{\pi}\frac{N!}{(N-|\pi|)!}$$

\index{P\'olya random walk}

As a second conclusion, even under the normalization $||H||_2=N$, the glow can behave quite differently in the $N\to\infty$ limit. So, let us restrict now the attention to the complex Hadamard matrices. At $N=2$ we only have $F_2$ to be invesigated, the result being:

\begin{proposition}
For the Fourier matrix $F_2$ we have
$$|E|^2=4+2Re(\alpha-\beta)$$
for certain variables $\alpha,\beta\in\mathbb T$ which are uniform, and independent.
\end{proposition}

\begin{proof}
The matrix that we interested in, namely the Fourier matrix $F_2$ altered by a vertical switching vector $(a,b)$ and an horizontal switching vector $(c,d)$, is:
$$\widetilde{F}_2=\begin{pmatrix}ac&ad\\bc&-bd\end{pmatrix}$$

With this notation, we have the following formula:
\begin{eqnarray*}
|E|^2
&=&|ac+ad+bc-bd|^2\\
&=&4+\frac{ad}{bc}+\frac{bc}{ad}-\frac{bd}{ac}-\frac{ac}{bd}
\end{eqnarray*}

For proving that the variables $\alpha=\frac{ad}{bc}$ and $\beta=\frac{bd}{ac}$ are independent, we can use the moment method, as follows:
\begin{eqnarray*}
\int_{\mathbb T^4}\left(\frac{ad}{bc}\right)^p\left(\frac{bd}{ac}\right)^q
&=&\int_{\mathbb T}a^{p-q}\int_{\mathbb T}b^{q-p}\int_{\mathbb T}c^{-p-q}\int_{\mathbb T}d^{p+q}\\
&=&\delta_{pq}\delta_{pq}\delta_{p,-q}\delta_{p,-q}\\
&=&\delta_{p,q,0}
\end{eqnarray*}

Thus $\alpha,\beta$ are indeed independent, and we are done.
\end{proof}

It is possible of course to derive from this some more concrete formulae, but let us look instead at the case $N=3$. Here the matrix that we are interested in is:
$$\widetilde{F}_3=\begin{pmatrix}ad&ae&af\\ bd&wbe&w^2bf\\ cd&w^2ce&wcf\end{pmatrix}$$

Thus, we would like to compute the law of the following quantity:
$$|E|=|ad+ae+af+bd+wbe+w^2bf+cd+w^2ce+wcf|$$

The problem is that when trying to compute $|E|^2$, the terms won't cancel much. More precisely, we have a formula of the following type:
$$|E|^2=9+C_0+C_1w+C_2w^2$$

Here the quantities $C_0,C_1,C_2$ are as follows:
\begin{eqnarray*}
C_0&=&\frac{ae}{bd}+\frac{ae}{cd}+\frac{af}{bd}+\frac{af}{cd}+\frac{bd}{ae}+\frac{bd}{af}
+\frac{be}{cf}+\frac{bf}{ce}+\frac{cd}{ae}+\frac{cd}{af}+\frac{ce}{bf}+\frac{cf}{be}\\
C_1&=&\frac{ad}{bf}+\frac{ad}{ce}+\frac{ae}{bf}+\frac{af}{ce}+\frac{bd}{ce}+\frac{be}{ad}
+\frac{be}{af}+\frac{be}{cd}+\frac{cd}{bf}+\frac{cf}{ad}+\frac{cf}{ae}+\frac{cf}{bd}\\
C_2&=&\frac{ad}{be}+\frac{ad}{cf}+\frac{ae}{cf}+\frac{af}{be}+\frac{bd}{cf}+\frac{bf}{ad}
+\frac{bf}{ae}+\frac{bf}{cd}+\frac{cd}{be}+\frac{ce}{ad}+\frac{ce}{af}+\frac{ce}{bd}
\end{eqnarray*}

In short, all this obviously leads nowhere, and the exact study stops at $F_2$. In general now, one idea is that of using Bernoulli-type variables coming from the row sums, a bit as we did in chapter 2 in the real case, the result here being as follows:

\begin{theorem}
The glow of $H\in M_N(\mathbb C)$ is given by the formula
$$law(E)=\int_{a\in\mathbb T^N}B((Ha)_1,\ldots,(Ha)_N)$$
where the quantities on the right are
$$B(c_1,\ldots,c_N)=law\left(\sum_i\lambda_ic_i\right)$$
with $\lambda\in\mathbb T^N$ being random.
\end{theorem}

\begin{proof}
This is clear indeed from the following formula:
$$E=<a,Hb>$$

To be more precise, when the vector $a\in\mathbb T^N$ is assumed to be fixed, this variable $E$ follows the law $B((Ha)_1,\ldots,(Ha)_N)$ in the statement.
\end{proof}

Observe that, in what regards the laws appearing in Theorem 11.5, we can write a formula for them of the following type, with $\times$ being a multiplicative convolution:
$$B(c_1,\ldots,c_N)=\varepsilon\times\beta(|c_1|,\ldots,|c_N|)$$

To be more precise, such a formula holds indeed, with the measure $\beta(r_1,\ldots,r_N)\in\mathcal P(\mathbb R_+)$ with $r_1,\ldots,r_N\geq 0$ being given by the following formula:
$$\beta(r_1,\ldots,r_N)=law\left|\sum_i\lambda_ir_i\right|$$

Regarding now the explicit computation of $\beta$, observe we have:
$$\beta(r_1,\ldots,r_N)=law\sqrt{\sum_{ij}\frac{\lambda_i}{\lambda_j}\cdot r_ir_j}$$

Consider now the following variable, which is easily seen, for instance by using the moment method, to be uniform over the projective torus $\mathbb T^{N-1}=\mathbb T^N/\mathbb T$:
$$(\mu_1,\mu_2,\ldots,\mu_N)=\left(\frac{\lambda_1}{\lambda_2},\frac{\lambda_2}{\lambda_3},\ldots,\frac{\lambda_N}{\lambda_1}\right)$$

Now since we have $\lambda_i/\lambda_j=\mu_i\mu_{i+1}\ldots\mu_j$, with the convention $\mu_i\ldots\mu_j=\overline{\mu_j\ldots\mu_i}$ for $i>j$, this gives the following formula, with $\mu\in\mathbb T^{N-1}$ random:
$$\beta(r_1,\ldots,r_N)=law\sqrt{\sum_{ij}\mu_i\mu_{i+1}\ldots\mu_j\cdot r_ir_j}$$

It is possible to further study the laws $\beta$ by using this formula. However, in practice, it is more convenient to use the complex measures $B$ from Theorem 11.5.

\bigskip

Let us end these preliminaries with a discussion of the ``arithmetic'' version of the problem, which makes the link with the Gale-Berlekamp game \cite{fsl}, \cite{rvi} and with the work in the real case, from chapter 2. We have the following unifying formalism:

\index{arithmetic glow}

\begin{definition}
Given $H\in M_N(\mathbb C)$ and $s\in\mathbb N\cup\{\infty\}$, we define a measure 
$$\mu_s\in\mathcal P(\mathbb C)$$
by the following formula, valid for any continuous function $\varphi$,
$$\int_\mathbb C\varphi(x)d\mu_s(x)=\int_{\mathbb Z^N_s\times\mathbb Z^N_s}\varphi\left(\sum_{ij}a_ib_jH_{ij}\right)d(a,b)$$
where $\mathbb Z_s\subset\mathbb T$ is the group of the $s$-roots of unity, with the convention $\mathbb Z_\infty=\mathbb T$.
\end{definition}

Observe that at $s=\infty$ we obtain the measure in Theorem 11.1. Also, at $s=2$ and for a usual Hadamard matrix, $H\in M_N(\pm1)$, we obtain the measure from chapter 2. Observe also that for $H\in M_N(\pm1)$, knowing $\mu_2$ is the same as knowing the statistics of the number of one entries, $|1\in H|$. This follows indeed from the following formula:
\begin{eqnarray*}
\sum_{ij}H_{ij}
&=&|1\in H|-|-1\in H|\\
&=&2|1\in H|-N^2
\end{eqnarray*}

More generally, at $s=p$ prime, we have the following result:

\index{number of 1 entries}

\begin{theorem}
When $s$ is prime and $H\in M_N(\mathbb Z_s)$, the statistics of the number of one entries, $|1\in H|$, can be recovered from that of the total sum, $E=\sum_{ij}H_{ij}$.
\end{theorem}

\begin{proof}
The problem here is of vectorial nature, so given $V\in\mathbb Z_s^n$, we would like to compare the quantities $|1\in V|$ and $\sum V_i$. Let us write, up to permutations:
$$V=(\underbrace{1\ldots1}_{a_0}\,\,\underbrace{w\ldots w}_{a_1}\,\ldots\ldots\,\underbrace{w^{s-1}\ldots w^{s-1}}_{a_{s-1}})$$

We have then $|1\in V|=a_0$, as well as: 
$$\sum V_i=a_0+a_1w+\ldots+a_{s-1}w^{s-1}$$

We also know that $a_0+a_1+\ldots+a_{s-1}=n$. Now when $s$ is prime, the only ambiguity in recovering $a_0$ from $a_0+a_1w+\ldots+a_{s-1}w^{s-1}$ can come from: 
$$1+w+\ldots+w^{s-1}=0$$

But since the sum of the numbers $a_i$ is fixed,  $a_0+a_1+\ldots+a_{s-1}=n$, this ambiguity dissapears, and this gives the result.
\end{proof}

\section*{11b. Glow moments}

Let us investigate now the glow of the complex Hadamard matrices, by using the moment method. We use the moment formula from chapter 10, namely:

\index{excess moments}
\index{glow moments}

\begin{proposition}
For $H\in M_N(\mathbb T)$ the even moments of $|E|$ are given by
$$\int_{\mathbb T^N\times\mathbb T^N}|E|^{2p}=\sum_{[i]=[k],[j]=[l]}\frac{H_{i_1j_1}\ldots H_{i_pj_p}}{H_{k_1l_1}\ldots H_{k_pl_p}}$$
where the sets between brackets are by definition sets with repetition.
\end{proposition}

\begin{proof}
As explained in chapter 10, with $E=\sum_{ij}H_{ij}a_ib_j$ we obtain:
\begin{eqnarray*}
\int_{\mathbb T^N\times\mathbb T^N}|E|^{2p}
&=&\int_{\mathbb T^N\times\mathbb T^N}\left(\sum_{ijkl}\frac{H_{ij}}{H_{kl}}\cdot\frac{a_ib_j}{a_kb_l}\right)^p\\
&=&\sum_{ijkl}\frac{H_{i_1j_1}\ldots H_{i_pj_p}}{H_{k_1l_1}\ldots H_{k_pl_p}}\int_{\mathbb T^N}\frac{a_{i_1}\ldots a_{i_p}}{a_{k_1}\ldots a_{k_p}}\int_{\mathbb T^N}\frac{b_{j_1}\ldots b_{j_p}}{b_{l_1}\ldots b_{l_p}}
\end{eqnarray*}

The integrals on the right being $\delta_{[i],[k]}$ and $\delta_{[j],[l]}$, we obtain the result.
\end{proof}

As a first application, let us investigate the tensor products. We have:

\index{tensor product}

\begin{proposition}
The even moments of the variable $|E|$ for a tensor product 
$$L=H\otimes K$$
are given by the following formula,
$$\int_{\mathbb T^{NM}\times\mathbb T^{NM}}|E|^{2p}=\sum_{[ia]=[kc],[jb]=[ld]}\frac{H_{i_1j_1}\ldots H_{i_pj_p}}{H_{k_1l_1}\ldots H_{k_pl_p}}\cdot\frac{K_{a_1b_1}\ldots K_{a_pb_p}}{K_{c_1d_1}\ldots K_{c_pd_p}}$$
where the sets between brackets are as usual sets with repetition.
\end{proposition}

\begin{proof}
With $L=H\otimes K$, the formula in Proposition 11.8 reads:
$$\int_{\mathbb T^{NM}\times\mathbb T^{NM}}|E|^{2p}=\sum_{[ia]=[kc],[jb]=[ld]}\frac{L_{i_1a_1,j_1b_1}\ldots L_{i_pa_p,j_pb_p}}{L_{k_1c_1,l_1d_1}\ldots L_{k_pc_p,l_pd_p}}$$

But this gives the formula in the statement, and we are done.
\end{proof}

Thus, we cannot reconstruct the glow of $H\otimes K$ from that of $H,K$, because the indices ``get mixed''. We have as well a result regarding the deformations, as follows:

\begin{proposition}
The even moments of $|E|$ for a deformed tensor product 
$$L=H\otimes_Q K$$
are given by the following formula,
$$\int_{\mathbb T^{NM}\times\mathbb T^{NM}}|E|^{2p}
=\sum_{[ia]=[kc],[jb]=[ld]}\frac{Q_{i_1b_1}\ldots Q_{i_pb_p}}{Q_{k_1d_1}\ldots Q_{k_pb_p}}\cdot\frac{H_{i_1j_1}\ldots H_{i_pj_p}}{H_{k_1l_1}\ldots H_{k_pl_p}}\cdot\frac{K_{a_1b_1}\ldots K_{a_pb_p}}{K_{c_1d_1}\ldots K_{c_pd_p}}$$
where the sets between brackets are as usual sets with repetition.
\end{proposition}

\begin{proof}
As before, we use the formula in Proposition 11.8. We have:
$$L_{ia,jb}=Q_{ib}H_{ij}K_{ab}$$

Thus, we obtain the following formula for the moments:
\begin{eqnarray*}
\int_{\mathbb T^{NM}\times\mathbb T^{NM}}|E|^{2p}
&=&\sum_{[ia]=[kc],[jb]=[ld]}\frac{L_{i_1a_1,j_1b_1}\ldots L_{i_pa_p,j_pb_p}}{L_{k_1c_1,l_1d_1}\ldots L_{k_pc_p,l_pd_p}}\\
&=&\sum_{[ia]=[kc],[jb]=[ld]}\frac{Q_{i_1b_1}\ldots Q_{i_pb_p}}{Q_{k_1d_1}\ldots Q_{k_pb_p}}\cdot\frac{H_{i_1j_1}\ldots H_{i_pj_p}}{H_{k_1l_1}\ldots H_{k_pl_p}}\cdot\frac{K_{a_1b_1}\ldots K_{a_pb_p}}{K_{c_1d_1}\ldots K_{c_pd_p}}
\end{eqnarray*}

Thus, we are led to the conclusion in the statement.
\end{proof}

The above formulae might look quite complicated, but they have some practical use. Let us go back indeed to a question that we had open since chapter 5, namely classifying the $4\times4$ complex Hadamard matrices, up to equivalence. We can now formulate:

\begin{theorem}
The complex Hadamard matrices at $N=4$ are, up to equivalence, the following matrices, with $s=e^{it}$ with $t\in[0,\pi/2]$,
$$F_4^s=\begin{pmatrix}
1&1&1&1\\
1&-1&1&-1\\
1&s&-1&-s\\ 
1&-s&-1&s
\end{pmatrix}$$
and these matrices are distinguished by the third moment of the glow. Alternatively, these matrices are distinguished by the third order term of the glow.
\end{theorem}

\begin{proof}
We know from chapter 5 that the $4\times4$ complex Hadamard matrices are those in the statement, with $s\in\mathbb T$, and we also know from there that we have equivalences as follows, which in practice means that we can assume $s=e^{it}$ with $t\in[0,\pi/2]$:
$$F_4^s\sim F_4^{-s}\sim F_4^{\bar{s}}\sim F_4^{-\bar{s}}$$

It remains to prove that these matrices, namely $F_4^s$ with $s=e^{it}$ with $t\in[0,\pi/2]$, are not equivalent. For this purpose, let us look at the moments of the glow:

\medskip

(1) Regarding the first moment, this is not something useful, because we have the following formula, coming from Proposition 11.8, valid for any $H\in M_N(\mathbb T)$:
$$\int_{\mathbb T^N\times\mathbb T^N}|E|^2=\sum_{i=k,j=l}\frac{H_{ij}}{H_{kl}}=\sum_{ij}\frac{H_{ij}}{H_{ij}}=N^2$$

(2) Regarding the second moment, this is something not useful either, because once again by using Proposition 11.8, we obtain a formula as follows, for any $H\in M_N(\mathbb T)$:
\begin{eqnarray*}
\int_{\mathbb T^N\times\mathbb T^N}|E|^4
&=&\sum_{[ik]=[mp],[jl]=[nq]}\frac{H_{ij}H_{kl}}{H_{mn}H_{pq}}\\
&=&P(N)+2\sum_{i\neq k,j\neq l}\frac{H_{ij}H_{kl}}{H_{il}H_{kj}}\\
&=&Q(N)+2\sum_{ijkl}\frac{H_{ij}H_{kl}}{H_{il}H_{kj}}\\
&=&Q(N)+2\sum_{ik}|<H_i,H_k>|^2\\
&=&Q(N)+2N^3
\end{eqnarray*}

To be more precise, here $P$ is a certain polynomial, not depending on $H$, collecting the contributions from the ``trivial'' solutions of $[ik]=[mp]$, $[jl]=[nq]$, and then $Q$ is another polynomial, again not depending on $H$, obtained from $P$ via a summing trick.

\medskip

(3) However, when getting to the third moment, or higher, things become interesting. Indeed, the equivalences $F_4^s\sim F_4^{-s}\sim F_4^{\bar{s}}\sim F_4^{-\bar{s}}$ tell us that the $p$-th moment of $|E|^2$ is a degree $p$ even, symmetric Laurent polynomial in $s\in\mathbb T$, and a direct computation at $p=3$, based on the formula in Proposition 11.10, shows that the parameter $s\in\mathbb T$ can be recaptured, up to identifying $\{s,-s,\bar{s},-\bar{s}\}$, from the knowledge of this polynomial.

\medskip

(4) Alternatively, we can say that the parameter $s\in\mathbb T$ can be recaptured, again up to identifying $\{s,-s,\bar{s},-\bar{s}\}$, from the knowledge of the third order term of the glow, with this meaning by definition the $N^{-2}$ factor in the $N^{-1}$ expansion of the law of $|E|/N$.
\end{proof}

Summarizing, some interesting things going on here, which will actually need some time to be fully understood. So, let us develop now some systematic moment machinery for the glow, along the above lines. Let $P(p)$ be the set of partitions of $\{1,\ldots,p\}$, with its standard order relation $\leq$, which is such that, for any $\pi\in P(p)$:
$$\sqcap\hskip-1.6mm\sqcap\ldots\leq\pi\leq|\ |\ldots|\ |$$ 

We denote by $\mu(\pi,\sigma)$ the associated M\"obius function, given by:
$$\mu(\pi,\sigma)=\begin{cases}
1&{\rm if}\ \pi=\sigma\\
-\sum_{\pi\leq\tau<\sigma}\mu(\pi,\tau)&{\rm if}\ \pi<\sigma\\
0&{\rm if}\ \pi\not\leq\sigma
\end{cases}$$

\index{partition}
\index{lattice of partitions}
\index{M\"obius function}
\index{M\"obius inversion}
\index{order of partitions}

To be more precise, the M\"obius function is defined by recurrence, by using this formula. The main interest in the M\"obius function comes from the M\"obius inversion formula, which states that the following happens, at the level of the functions on $P(p)$:
$$f(\sigma)=\sum_{\pi\leq\sigma}g(\pi)
\quad\implies\quad g(\sigma)=\sum_{\pi\leq\sigma}\mu(\pi,\sigma)f(\pi)$$

For $\pi\in P(p)$ we use the following notation, where $b_1,\ldots,b_{|\pi|}$ are the block lenghts:
$$\binom{p}{\pi}=\binom{p}{b_1\ldots b_{|\pi|}}=\frac{p!}{b_1!\ldots b_{|\pi|}!}$$

Finally, we use the following notation, where $H_1,\ldots,H_N\in\mathbb T^N$ are the rows of $H$:
$$H_\pi(i)=\bigotimes_{\beta\in\pi}\prod_{r\in\beta}H_{i_r}$$

With these notations, we have the following result:

\begin{theorem}
The glow moments of a matrix $H\in M_N(\mathbb T)$ are given by
$$\int_{\mathbb T^N\times\mathbb T^N}|E|^{2p}=\sum_{\pi\in P(p)}K(\pi)N^{|\pi|}I(\pi)$$
where the coefficients are given by
$$K(\pi)=\sum_{\sigma\in P(p)}\mu(\pi,\sigma)\binom{p}{\sigma}$$
and where the contributions are given by
$$I(\pi)=\frac{1}{N^{|\pi|}}\sum_{[i]=[j]}<H_\pi(i),H_\pi(j)>$$
by using the above notations and conventions.
\end{theorem}

\begin{proof}
We know from Proposition 11.8 that the moments are given by:
$$\int_{\mathbb T^N\times\mathbb T^N}|E|^{2p}
=\sum_{[i]=[j],[x]=[y]}\frac{H_{i_1x_1}\ldots H_{i_px_p}}{H_{j_1y_1}\ldots H_{j_py_p}}$$

With $\sigma=\ker x,\rho=\ker y$, we deduce that the moments of $|E|^2$ decompose over partitions, according to a formula as follows:
$$\int_{\mathbb T^N\times\mathbb T^N}|E|^{2p}=\int_{\mathbb T^N}\sum_{\sigma,\rho\in P(p)}C(\sigma,\rho)$$

To be more precise, the contributions are as follows:
$$C(\sigma,\rho)=\sum_{\ker x=\sigma,\ker y=\rho}\delta_{[x],[y]}\sum_{ij}
\frac{H_{i_1x_1}\ldots H_{i_px_p}}{H_{j_1y_1}\ldots H_{j_py_p}}
\cdot\frac{a_{i_1}\ldots a_{i_p}}{a_{j_1}\ldots a_{j_p}}$$

We have $C(\sigma,\rho)=0$ unless $\sigma\sim\rho$, in the sense that $\sigma,\rho$ must have the same block structure. The point now is that the sums of type $\sum_{\ker x=\sigma}$ can be computed by using the M\"obius inversion formula. We obtain a formula as follows:
$$C(\sigma,\rho)=\delta_{\sigma\sim\rho}\sum_{\pi\leq\sigma}\mu(\pi,\sigma)\prod_{\beta\in\pi}C_{|\beta|}(a)$$

Here the functions on the right are by definition given by:
\begin{eqnarray*}
C_r(a)
&=&\sum_x\sum_{ij}\frac{H_{i_1x}\ldots H_{i_rx}}{H_{j_1x}\ldots H_{j_rx}}\cdot\frac{a_{i_1}\ldots a_{i_r}}{a_{j_1}\ldots a_{j_r}}\\
&=&\sum_{ij}<H_{i_1}\ldots H_{i_r},H_{j_1}\ldots H_{j_r}>\cdot\frac{a_{i_1}\ldots a_{i_r}}{a_{j_1}\ldots a_{j_r}}
\end{eqnarray*}

Now since there are $\binom{p}{\sigma}$ partitions having the same block structure as $\sigma$, we obtain:
\begin{eqnarray*}
&&\int_{\mathbb T^N\times\mathbb T^N}|\Omega|^{2p}\\
&=&\int_{\mathbb T^N}\sum_{\pi\in P(p)}\left(\sum_{\sigma\sim\rho}\sum_{\mu\leq\sigma}\mu(\pi,\sigma)\right)\prod_{\beta\in\pi}C_{|\beta|}(a)\\
&=&\sum_{\pi\in P(p)}\left(\sum_{\sigma\in P(p)}\mu(\pi,\sigma)\binom{p}{\sigma}\right)\int_{\mathbb T^N}\prod_{\beta\in\pi}C_{|\beta|}(a)
\end{eqnarray*}

But this gives the formula in the statement, and we are done.
\end{proof}

Let us discuss now the asymptotic behavior of the glow. For this purpose, we first study the coefficients $K(\pi)$ in Theorem 11.12. We have here the following result:

\begin{proposition}
The coeffients appearing in the above, namely
$$K(\pi)=\sum_{\pi\leq\sigma}\mu(\pi,\sigma)\binom{p}{\sigma}$$
have the following properties:
\begin{enumerate}
\item The function $\widetilde{K}(\pi)=\frac{K(\pi)}{p!}$ is multiplicative, in the sense that: 
$$\widetilde{K}(\pi\pi')=\widetilde{K}(\pi)\widetilde{K}(\pi')$$

\item On the one-block partitions, we have:
$$K(\sqcap\!\!\sqcap\ldots\sqcap)=\sum_{\sigma\in P(p)}(-1)^{|\sigma|-1}(|\sigma|-1)!\binom{p}{\sigma}$$

\item We have as well the following fomula,
$$K(\sqcap\!\!\sqcap\ldots\sqcap)=\sum_{r=1}^p(-1)^{r-1}(r-1)!C_{pr}$$
where the coefficients on the right are given by:
$$C_{pr}=\sum_{p=a_1+\ldots+a_r}\binom{p}{a_1,\ldots,a_r}^2$$
\end{enumerate}
\end{proposition}

\begin{proof}
This follows from some standard computations, as follows:

\medskip

(1) We can use here the following formula, which is a well-known property of the M\"obius function, which can be proved by recurrence:
$$\mu(\pi\pi',\sigma\sigma')=\mu(\pi,\sigma)\mu(\pi',\sigma')$$

Now if $b_1,\ldots,b_s$ and $c_1,\ldots,c_t$ are the block lengths of $\sigma,\sigma'$, we obtain, as claimed:
\begin{eqnarray*}
\widetilde{K}(\pi\pi')
&=&\sum_{\pi\pi'\leq\sigma\sigma'}\mu(\pi\pi',\sigma\sigma')\cdot\frac{1}{b_1!\ldots b_s!}\cdot\frac{1}{c_1!\ldots c_t!}\\
&=&\sum_{\pi\leq\sigma,\pi'\leq\sigma'}\mu(\pi,\sigma)\mu(\pi',\sigma')\cdot\frac{1}{b_1!\ldots b_s!}\cdot\frac{1}{c_1!\ldots c_t!}\\
&=&\widetilde{K}(\pi)\widetilde{K}(\pi')
\end{eqnarray*}

(2) We can use here the following formula, which once again is well-known, and can be proved by recurrence on $|\sigma|$:
$$\mu(\sqcap\!\!\sqcap\ldots\sqcap,\sigma)=(-1)^{|\sigma|-1}(|\sigma|-1)!$$

We therefore obtain, as claimed:
$$K(\sqcap\!\!\sqcap\ldots\sqcap)
=\sum_{\sigma\in P(p)}\mu(\sqcap\!\!\sqcap\ldots\sqcap,\sigma)\binom{p}{\sigma}
=\sum_{\sigma\in P(p)}(-1)^{|\sigma|-1}(|\sigma|-1)!\binom{p}{\sigma}$$

(3) By using the formula in (2), and summing over $r=|\sigma|$, we obtain:
$$K(\sqcap\!\!\sqcap\ldots\sqcap)
=\sum_{r=1}^p(-1)^{r-1}(r-1)!\sum_{|\sigma|=r}\binom{p}{\sigma}$$

Now if we denote by $a_1,\ldots,a_r$ with $a_i\geq1$ the block lengths of $\sigma$, then: 
$$\binom{p}{\sigma}=\binom{p}{a_1,\ldots,a_r}$$

On the other hand, given $a_1,\ldots,a_r\geq1$ with $a_1+\ldots+a_r=p$, the number of partitions $\sigma$ having these numbers as block lengths is:
$$N_{a_1,\ldots,a_r}=\binom{p}{a_1,\ldots,a_r}$$

Thus, we are led to the conclusion in the statement.
\end{proof}

Now let us take a closer look at the integrals $I(\pi)$ from Theorem 11.12, namely:
$$I(\pi)=\frac{1}{N^{|\pi|}}\sum_{[i]=[j]}<H_\pi(i),H_\pi(j)>$$

We have here the following result:

\begin{proposition}
Consider the one-block partition $\sqcap\!\!\sqcap\ldots\sqcap\in P(p)$.
\begin{enumerate}
\item $I(\sqcap\!\!\sqcap\ldots\sqcap)=\#\{i,j\in\{1,\ldots,N\}^p|[i]=[j]\}$.

\item $I(\sqcap\!\!\sqcap\ldots\sqcap)=\int_{\mathbb T^N}|\sum_ia_i|^{2p}da$.

\item $I(\sqcap\!\!\sqcap\ldots\sqcap)=\sum_{\sigma\in P(p)}\binom{p}{\sigma}\frac{N!}{(N-|\sigma|)!}$.

\item $I(\sqcap\!\!\sqcap\ldots\sqcap)=\sum_{r=1}^{p-1}C_{pr}\frac{N!}{(N-r)!}$, where $C_{pr}=\sum_{p=b_1+\ldots+b_r}\binom{p}{b_1,\ldots,b_r}^2$.
\end{enumerate}
\end{proposition}

\begin{proof}
Once again, these formulae follow from some standard combinatorics:

\medskip

(1) This follows indeed from the following computation:
$$I(\sqcap\!\!\sqcap\ldots\sqcap)
=\sum_{[i]=[j]}\frac{1}{N}<H_{i_1}\ldots H_{i_r},H_{j_1}\ldots H_{j_r}>
=\sum_{[i]=[j]}1$$

(2) This follows from the following computation:
$$\int_{\mathbb T^N}\left|\sum_ia_i\right|^{2p}
=\int_{\mathbb T^N}\sum_{ij}\frac{a_{i_1}\ldots a_{i_p}}{a_{j_1}\ldots a_{j_p}}da
=\#\left\{i,j\Big|[i]=[j]\right\}$$

(3) If we let $\sigma=\ker i$ in the above formula of $I(\sqcap\!\!\sqcap\ldots\sqcap)$, we obtain:
$$I(\sqcap\!\!\sqcap\ldots\sqcap)=\sum_{\sigma\in P(p)}\#\left\{i,j\Big|\ker i=\sigma,[i]=[j]\right\}$$

Now since there are $\frac{N!}{(N-|\sigma|)!}$ choices for the multi-index $i$, and then $\binom{p}{\sigma}$ choices for the multi-index $j$, this gives the result.

\medskip

(4) If we set $r=|\sigma|$, the formula in (3) becomes:
$$I(\sqcap\!\!\sqcap\ldots\sqcap)=\sum_{r=1}^{p-1}\frac{N!}{(N-r)!}\sum_{\sigma\in P(p),|\sigma|=r}\binom{p}{\sigma}$$

Now since there are exactly $\binom{p}{b_1,\ldots,b_r}$ permutations $\sigma\in P(p)$ having $b_1,\ldots,b_r$ as block lengths, the sum on the right is given by:
$$\sum_{\sigma\in P(p),|\sigma|=r}\binom{p}{\sigma}=\sum_{p=b_1+\ldots+b_r}\binom{p}{b_1,\ldots,b_r}^2$$

Thus, we are led to the conclusion in the statement.
\end{proof}

In general, the integrals $I(\pi)$ can be estimated as follows:

\begin{proposition}
Let $H\in M_N(\mathbb T)$, having its rows pairwise orthogonal.
\begin{enumerate}
\item $I(|\,|\,\ldots|)=N^p$.

\item $I(|\,|\,\ldots|\ \pi)=N^aI(\pi)$, for any $\pi\in P(p-a)$.

\item $|I(\pi)|\lesssim p!N^p$, for any $\pi\in P(p)$. 
\end{enumerate}
\end{proposition}

\begin{proof}
This is something elementary, as follows:

\medskip

(1) Since the rows of $H$ are pairwise orthogonal, we have:
\begin{eqnarray*}
I(|\,|\ldots|)
&=&\sum_{[i]=[j]}\prod_{r=1}^p\delta_{i_r,j_r}\\
&=&\sum_{[i]=[j]}\delta_{ij}\\
&=&\sum_i1\\
&=&N^p
\end{eqnarray*}

(2) This follows by the same computation as the above one for (1).

\medskip

(3) We have indeed the following estimate:
\begin{eqnarray*}
|I(\pi)|
&\leq&\sum_{[i]=[j]}\prod_{\beta\in\pi}1\\
&=&\sum_{[i]=[j]}1\\
&=&\#\left\{i,j\in\{1,\ldots,N\}\Big|[i]=[j]\right\}\\
&\simeq&p!N^p
\end{eqnarray*}

Thus we have obtained the formula in the statement, and we are done.
\end{proof}

We have now all needed ingredients for a universality result:

\index{glow universality}
\index{complex normal variable}
\index{complex Gaussian variable}

\begin{theorem}
The glow of a complex Hadamard matrix $H\in M_N(\mathbb T)$ is given by:
$$\frac{1}{p!}\int_{\mathbb T^N\times\mathbb T^N}\left|\frac{E}{N}\right|^{2p}=1-\binom{p}{2}N^{-1}+O(N^{-2})$$
In particular, $E/N$ becomes complex Gaussian in the $N\to\infty$ limit.
\end{theorem}

\begin{proof}
We use the moment formula in Theorem 11.12, namely:
$$\int_{\mathbb T^N\times\mathbb T^N}|E|^{2p}=\sum_{\pi\in P(p)}K(\pi)N^{|\pi|}I(\pi)$$

By using Proposition 11.15 (3), we conclude that only the $p$-block and $(p-1)$-block partitions contribute at order 2, so:
\begin{eqnarray*}
\int_{\mathbb T^N\times\mathbb T^N}|E|^{2p}
&=&K(|\,|\ldots|)N^pI(|\,|\ldots|)\\
&+&\binom{p}{2}K(\sqcap|\ldots|)N^{p-1}I(\sqcap|\ldots|)\\
&+&O(N^{2p-2})
\end{eqnarray*}

Now by dividing by $N^{2p}$ and then by using the various formulae in Proposition 11.13, Proposition 11.14 and Proposition 11.15, we obtain, as claimed:
$$\int_{\mathbb T^N\times\mathbb T^N}\left|\frac{E}{N}\right|^{2p}
=p!
-\binom{p}{2}\frac{p!}{2}\cdot\frac{2N-1}{N^2}
+O(N^{-2})$$

Finally, since the law of $E$ is invariant under centered rotations in the complex plane, this moment formula gives as well the last assertion.
\end{proof}

Summarizing, the complex glow of the complex Hadamard matrices appears to have similar properties to the real glow of the real Hadamard matrices.

\section*{11c. Fourier matrices}

\index{Fourier matrix}

Let us study now the glow of the Fourier matrices, $F=F_G$. We use the following standard formulae, which all come from definitions:
$$F_{ix}F_{iy}=F_{i,x+y}\quad,\quad 
\overline{F}_{ix}=F_{i,-x}\quad,\quad 
\sum_xF_{ix}=N\delta_{i0}$$

We first have the following result:

\begin{proposition}
For a Fourier matrix $F_G$ we have
$$I(\pi)=\#\left\{i,j\Big|[i]=[j],\sum_{r\in\beta}i_r=\sum_{r\in\beta}j_r,\forall\beta\in\pi\right\}$$
with all the indices, and with the sums at right, taken inside $G$.
\end{proposition}

\begin{proof}
The basic components of the integrals $I(\pi)$ are given by:
\begin{eqnarray*}
\frac{1}{N}\left\langle\prod_{r\in\beta}F_{i_r},\prod_{r\in\beta}F_{j_r}\right\rangle
&=&\frac{1}{N}\left\langle F_{\sum_{r\in\beta}i_r},F_{\sum_{r\in\beta}i_r}\right\rangle\\
&=&\delta_{\sum_{r\in\beta}i_r,\sum_{r\in\beta}j_r}
\end{eqnarray*}

But this gives the formula in the statement, and we are done.
\end{proof}

We have the following interpretation of the above integrals:

\begin{proposition}
For any partition $\pi$ we have the formula
$$I(\pi)=\int_{\mathbb T^N}\prod_{b\in\pi}\left(\frac{1}{N^2}\sum_{ij}|H_{ij}|^{2|\beta|}\right)da$$
where $H=FAF^*$, with $F=F_G$ and $A=diag(a_0,\ldots,a_{N-1})$.
\end{proposition}

\begin{proof}
We have the following computation:
\begin{eqnarray*}
H=F^*AF
&\implies&|H_{xy}|^2=\sum_{ij}\frac{F_{iy}F_{jx}}{F_{ix}F_{jy}}\cdot\frac{a_i}{a_j}\\
&\implies&|H_{xy}|^{2p}=\sum_{ij}\frac{F_{j_1x}\ldots F_{j_px}}{F_{i_1x}\ldots F_{i_px}}\cdot\frac{F_{i_1y}\ldots F_{i_py}}{F_{j_1y}\ldots F_{j_py}}\cdot\frac{a_{i_1}\ldots a_{i_p}}{a_{j_1}\ldots a_{j_p}}\\
&\implies&\sum_{xy}|H_{xy}|^{2p}=\sum_{ij}\left|<H_{i_1}\ldots H_{i_p},H_{j_1}\ldots H_{j_p}>\right|^2\cdot\frac{a_{i_1}\ldots a_{i_p}}{a_{j_1}\ldots a_{j_p}}
\end{eqnarray*}

But this gives the formula in the statement, and we are done.
\end{proof}

We must estimate now the quantities $I(\pi)$. We first have the following result:

\begin{proposition}
For $F_G$ we have the estimate
$$I(\pi)=b_1!\ldots b_{|\pi|}!N^p+O(N^{p-1})$$
where the numbers $b_1,\ldots,b_{|\pi|}$ with $b_1+\ldots+b_{|\pi|}=p$ are the block lengths of $\pi$.
\end{proposition}

\begin{proof}
With $\sigma=\ker i$ we obtain:
$$I(\pi)=\sum_{\sigma\in P(p)}\#\left\{i,j\Big|\ker i=\sigma,[i]=[j],\sum_{r\in\beta}i_r=\sum_{r\in\beta}j_r,\forall\beta\in\pi\right\}$$

The number of choices for $i$ satisfying $\ker i=\sigma$ is:
$$\frac{N!}{(N-|\sigma|)!}\simeq N^{|\sigma|}$$

Then, the number of choices for $j$ satisfying $[i]=[j]$ is:
$$\binom{p}{\sigma}=O(1)$$

We conclude that the main contribution comes from the following partition:
$$\sigma=|\,|\ldots|$$

Thus, we have the following formula:
$$I(\pi)=\#\left\{i,j\Big|\ker i=|\,|\ldots|,[i]=[j],\sum_{r\in\beta}i_r=\sum_{r\in\beta}j_r,\forall\beta\in\pi\right\}+O(N^{p-1})$$

Now $\ker i=|\,|\ldots|$ tells us that $i$ must have distinct entries, and there are $\frac{N!}{(N-p)!}\simeq N^p$ choices for such multi-indices $i$. Regarding the indices $j$, the main contribution comes from those obtained from $i$ by permuting the entries over the blocks of $\pi$, and there are $b_1!\ldots b_{|\pi|}!$ choices here. Thus, we are led to the conclusion in the statement.
\end{proof}

At the second order now, the estimate is as follows:

\begin{proposition}
For $F_G$ we have the formula
$$\frac{I(\pi)}{b_1!\ldots b_s!N^p}=1+\left(\sum_{i<j}\sum_{c\geq2}\binom{b_i}{c}\binom{b_j}{c}-\frac{1}{2}\sum_i\binom{b_i}{2}\right)N^{-1}+O(N^{-2})$$
where $b_1,\ldots,b_s$ being the block lengths of $\pi\in P(p)$.
\end{proposition}

\begin{proof}
Let us define the ``non-arithmetic'' part of $I(\pi)$ as follows: 
$$I^\circ(\pi)=\#\left\{i,j\Big|[i_r|r\in\beta]=[j_r|r\in\beta],\forall\beta\in\pi\right\}$$

We then have the following formula:
$$I^\circ(\pi)=\prod_{\beta\in\pi}\left\{i,j\in I^{|\beta|}\Big|[i]=[j]\right\}=\prod_{\beta\in\pi}I(\beta)$$

Also, Proposition 11.19 shows that we have the following estimate:
$$I(\pi)=I^\circ(\pi)+O(N^{p-1})$$

Our claim now is that we have the following formula:
$$\frac{I(\pi)-I^\circ(\pi)}{b_1!\ldots b_s!N^p}=\sum_{i<j}\sum_{c\geq2}\binom{b_i}{c}\binom{b_j}{c}N^{-1}+O(N^{-2})$$

Indeed, according to Proposition 11.19, we have a formula of the following type:
$$I(\pi)=I^\circ(\pi)+I^1(\pi)+O(N^{p-2})$$

More precisely, this formula holds indeed, with $I^1(\pi)$ coming from $i_1,\ldots,i_p$ distinct, $[i]=[j]$, and with one constraint of type:
$$\sum_{r\in\beta}i_r=\sum_{j\in\beta}j_r\quad,\quad 
[i_r|r\in\beta]\neq[j_r|r\in\beta]$$ 

Now observe that for a two-block partition $\pi=(a,b)$ this constraint is implemented, up to permutations which leave invariant the blocks of $\pi$, as follows:
$$\begin{matrix}
i_1\ldots i_c&k_1\ldots k_{a-c}&&j_1\ldots j_c&l_1\ldots l_{a-c}\\
\underbrace{j_1\ldots j_c}_c&\underbrace{k_1\ldots k_{a-c}}_{a-c}&&\underbrace{i_1\ldots i_c}_c&\underbrace{l_1\ldots l_{a-c}}_{b-c}
\end{matrix}$$

Let us compute now $I^1(a,b)$. We cannot have $c=0,1$, and once $c\geq2$ is given, we have $\binom{a}{c},\binom{b}{c}$ choices for the positions of the $i,j$ variables in the upper row, then $N^{p-1}+O(N^{p-2})$ choices for the variables in the upper row, and then finally we have $a!b!$ permutations which can produce the lower row. We therefore obtain:
$$I^1(a,b)=a!b!\sum_{c\geq2}\binom{a}{c}\binom{b}{c}N^{p-1}+O(N^{p-2})$$

In the general case now, a similar discussion applies. Indeed, the constraint of type $\sum_{r\in\beta}i_r=\sum_{r\in\beta}j_r$ with $[i_r|r\in\beta]\neq[j_r|r\in\beta]$ cannot affect $\leq1$ blocks, because we are not in the non-arithmetic case, and cannot affect either $\geq3$ blocks, because affecting $\geq3$ blocks would require $\geq2$ constraints. Thus this condition affects exactly $2$ blocks, and if we let $i<j$ be the indices in $\{1,\ldots,s\}$ corresponding to these 2 blocks, we obtain:
$$I^1(\pi)=b_1!\ldots b_s!\sum_{i<j}\sum_{c\geq2}\binom{b_i}{c}\binom{b_j}{c}N^{p-1}+O(N^{p-2})$$

But this proves the above claim. Let us estimate now $I(\sqcap\!\!\sqcap\ldots\sqcap)$. We have:
\begin{eqnarray*}
&&I(\sqcap\!\!\sqcap\ldots\sqcap)\\
&=&p!\frac{N!}{(N-p)!}+\binom{p}{2}\frac{p!}{2}\cdot\frac{N!}{(N-p+1)!}+O(N^{p-2})\\
&=&p!N^r\left(1-\binom{p}{2}N^{-1}+O(N^{-2})\right)+\binom{p}{2}\frac{p!}{2}N^{p-1}+O(N^{p-2})\\
&=&p!N^p\left(1-\frac{1}{2}\binom{p}{2}N^{-1}+O(N^{-2})\right)
\end{eqnarray*}

Now recall that we have:
$$I^\circ(\pi)=\prod_{\beta\in\pi}I(\beta)$$

We therefore obtain:
$$I^\circ(\pi)=b_1!\ldots b_s!N^p\left(1-\frac{1}{2}\sum_i\binom{b_i}{2}N^{-1}+O(N^{-2})\right)$$

By plugging this quantity into the above estimate, we obtain the result.
\end{proof}

In order to estimate glow, we will need the explicit formula of $I(\sqcap\sqcap)$:

\begin{proposition}
For $F_G$ with $G=\mathbb Z_{N_1}\times\ldots\times\mathbb Z_{N_k}$ we have the formula
$$I(\sqcap\sqcap)=N(4N^3-11N+2^e+7)$$
where $e\in\{0,1,\ldots,k\}$ is the number of even numbers among $N_1,\ldots,N_k$.
\end{proposition}

\begin{proof}
The conditions defining the quantities $I(\pi)$ are as follows:
$$\sum_{r\in\beta}i_r=\sum_{r\in\beta}j_r$$

We use the fact that, when dealing with these conditions, one can always erase some of the variables $i_r,j_r$, as to reduce to the ``purely arithmetic'' case, namely:
$$\{i_r|r\in\beta\}\cap\{j_r|r\in\beta\}=\emptyset$$

We deduce from this that we have:
$$I(\sqcap\sqcap)=I^\circ(\sqcap\sqcap)+I^{ari}(\sqcap\sqcap)$$

Let us compute now $I^{ari}(\sqcap\sqcap)$. There are 3 contributions to this quantity, namely:

\medskip

(1) \underline{Case $(^{iijj}_{jjii})$}, with $i\neq j$, $2i=2j$. Since $2(i_1,\ldots,i_k)=2(j_1,\ldots,j_k)$ corresponds to the collection of conditions $2i_r=2j_r$, inside $\mathbb Z_{N_r}$, which each have 1 or 2 solutions, depending on whether $N_r$ is odd or even, the contribution here is:
\begin{eqnarray*}
I^{ari}_1(\sqcap\sqcap)
&=&\#\{i\neq j|2i=2j\}\\
&=&\#\{i,j|2i=2j\}-\#\{i,j|i=j\}\\
&=&2^eN-N\\
&=&(2^e-1)N
\end{eqnarray*}

(2) \underline{Case $(^{iijk}_{jkii})$}, with $i,j,k$ distinct, $2i=j+k$. The contribution here is:
\begin{eqnarray*}
I^{ari}_2(\sqcap\sqcap)
&=&4\#\{i,j,k\ {\rm distinct}|2i=j+k\}\\
&=&4\#\{i\neq j|2i-j\neq i,j\}\\
&=&4\#\{i\neq j|2i\neq 2j\}\\
&=&4(\#\{i,j|i\neq j\}-\#\{i\neq j|2i=2j\})\\
&=&4(N(N-1)-(2^e-1)N)\\
&=&4N(N-2^e)
\end{eqnarray*}

(3) \underline{Case $(^{ijkl}_{klij})$}, with $i,j,k,l$ distinct, $i+j=k+l$. The contribution here is:
\begin{eqnarray*}
I^{ari}_3(\sqcap\sqcap)
&=&4\#\{i,j,k,l\ {\rm distinct}|i+j=k+l\}\\
&=&4\#\{i,j,k\ {\rm distinct}|i+j-k\neq i,j,k\}\\
&=&4\#\{i,j,k\ {\rm distinct}|i+j-k\neq k\}\\
&=&4\#\{i,j,k\ {\rm distinct}|i\neq 2k-j\}
\end{eqnarray*}

We can split this quantity over two cases, $2j\neq 2k$ and $2j=2k$, and we obtain:
\begin{eqnarray*}
I^{ari}_3(\sqcap\sqcap)
&=&4(\#\{i,j,k\ {\rm distinct}|2j\neq 2k,i\neq 2k-j\}\\
&&+\#\{i,j,k\ {\rm distinct}|2j=2k,i\neq 2k-j\})
\end{eqnarray*}

The point now is that in the first case, $2j\neq 2k$, the numbers $j,k,2k-j$ are distinct, while in the second case, $2j=2k$, we simply have $2k-j=j$. Thus, we obtain:
\begin{eqnarray*}
I^{ari}_3(\sqcap\sqcap)
&=&4\left(\sum_{j\neq k,2j\neq 2k}\#\{i|i\neq j,k,2k-j\}+\sum_{j\neq k,2j=2k}\#\{i|i\neq j,k\}\right)\\
&=&4(N(N-2^e)(N-3)+N(2^e-1)(N-2))\\
&=&4N(N(N-3)-2^e(N-3)+2^e(N-2)-(N-2))\\
&=&4N(N^2-4N+2^e+2)
\end{eqnarray*}

We can now compute the arithmetic part. This is given by:
\begin{eqnarray*}
I^{ari}(\sqcap\sqcap)
&=&(2^e-1)N+4N(N-2^e)+4N(N^2-4N+2^e+2)\\
&=&N(2^e-1+4(N-2^e)+4(N^2-4N+2^e+2))\\
&=&N(4N^2-12N+2^e+7)
\end{eqnarray*}

Thus the integral to be computed is given by:
\begin{eqnarray*}
I(\sqcap\sqcap)
&=&N^2(2N-1)^2+N(4N^2-12N+2^e+7)\\
&=&N(4N^3-4N^2+N+4N^2-12N+2^e+7)\\
&=&N(4N^3-11N+2^e+7)
\end{eqnarray*}

Thus we have reached to the formula in the statement, and we are done.
\end{proof}

\section*{11d. Universality}

We have the following asymptotic result:

\index{glow universality}
\index{complex normal variable}
\index{complex Gaussian variable}

\begin{theorem}
The glow of $F_G$, with $|G|=N$, is given by 
$$\frac{1}{p!}\int_{\mathbb T^N\times\mathbb T^N}\left|\frac{E}{N}\right|^{2p}=1-K_1N^{-1}+K_2N^{-2}-K_3N^{-3}+O(N^{-4})$$ 
with the coefficients being as follows:
$$K_1=\binom{p}{2}\quad,\quad 
K_2=\binom{p}{2}\frac{3p^2+p-8}{12}\quad,\quad 
K_3=\binom{p}{3}\frac{p^3+4p^2+p-18}{8}$$
Thus, the rescaled complex glow is asymptotically complex Gaussian,
$$\frac{E}{N}\sim\mathcal C$$
and we have in fact universality at least up to order $3$.
\end{theorem}

\begin{proof}
We use the following quantities:
$$\widetilde{K}(\pi)=\frac{K(\pi)}{p!}\quad,\quad 
\widetilde{I}(\pi)=\frac{I(\pi)}{N^p}$$

These are subject to the following formulae:
$$\widetilde{K}(\pi|\ldots|)=\widetilde{K}(\pi)\quad,\quad 
\widetilde{I}(\pi|\ldots|)=\widetilde{I}(\pi)$$

Consider as well the following quantities:
$$J(\sigma)=\binom{p}{\sigma}\widetilde{K}(\sigma)\widetilde{I}(\sigma)$$

In terms of these quantities, we have:
\begin{eqnarray*}
\frac{1}{p!}\int_{\mathbb T^N\times\mathbb T^N}|E|^{2p}
&=&J(\emptyset)\\
&+&N^{-1}J(\sqcap)\\
&+&N^{-2}\left(J(\sqcap\!\sqcap)+J(\sqcap\sqcap)\right)\\
&+&N^{-3}\left(J(\sqcap\!\!\sqcap\!\!\sqcap)+J(\sqcap\!\!\sqcap\sqcap)+J(\sqcap\sqcap\sqcap)\right)\\
&+&O(N^{-4})
\end{eqnarray*}

We have the following formulae:
$$\widetilde{K}_0=1$$
$$\widetilde{K}_1=1$$
$$\widetilde{K}_2=\frac{1}{2}-1=-\frac{1}{2}$$ 
$$\widetilde{K}_3=\frac{1}{6}-\frac{3}{2}+2=\frac{2}{3}$$
$$\widetilde{K}_4=\frac{1}{24}-\frac{4}{6}-\frac{3}{4}+\frac{12}{2}-6=-\frac{11}{8}$$

Regarding now the numbers $C_{pr}$ in Proposition 11.19, these are given by:
$$C_{p1}=1$$
$$C_{p2}=\frac{1}{2}\binom{2p}{p}-1$$
$$\vdots$$
$$C_{p,p-1}=\frac{p!}{2}\binom{p}{2}$$
$$C_{pp}=p!$$

We deduce that we have the following formulae:
$$I(|)=N$$
$$I(\sqcap)=N(2N-1)$$
$$I(\sqcap\!\sqcap)=N(6N^2-9N+4)$$
$$I(\sqcap\!\!\sqcap\!\!\sqcap)=N(24N^3-72N^2+82N-33)$$

By using Proposition 11.20 and Proposition 11.21, we obtain the following formula:
\begin{eqnarray*}
\frac{1}{p!}\int_{\mathbb T^N\times\mathbb T^N}|E|^{2p}
&=&1-\frac{1}{2}\binom{p}{2}(2N^{-1}-N^{-2})+\frac{2}{3}\binom{p}{3}(6N^{-2}-9N^{-3})\\
&+&3\binom{p}{4}N^{-2}-33\binom{p}{4}N^{-3}-40\binom{p}{5}N^{-3}\\
&-&15\binom{p}{6}N^{-3}+O(N^{-4})
\end{eqnarray*}

But this gives the formulae of $K_1,K_2,K_3$ in the statement, and we are done.
\end{proof}

It is possible to compute the next term as well, the result being as follows:

\index{glow universality}
\index{complex normal variable}
\index{complex Gaussian variable}

\begin{theorem}
Let $G=\mathbb Z_{N_1}\times\ldots\times\mathbb Z_{N_k}$ be a finite abelian group, and set:
$$N=N_1\ldots N_k$$
Then the glow of the associated Fourier matrix $F_G$ is given by 
$$\frac{1}{p!}\int_{\mathbb T^N\times\mathbb T^N}\left|\frac{E}{N}\right|^{2p}=1-K_1N^{-1}+K_2N^{-2}-K_3N^{-3}+K_4N^{-4}+O(N^{-5})$$ where the quantities $K_1,K_2,K_3,K_4$ are given by
\begin{eqnarray*}
K_1&=&\binom{p}{2}\\
K_2&=&\binom{p}{2}\frac{3p^2+p-8}{12}\\
K_3&=&\binom{p}{3}\frac{p^3+4p^2+p-18}{8}\\
K_4&=&\frac{8}{3}\binom{p}{3}+\frac{3}{4}\left(121+\frac{2^e}{N}\right)\binom{p}{4}+416\binom{p}{5}+\frac{2915}{2}\binom{p}{6}+40\binom{p}{7}+105\binom{p}{8}
\end{eqnarray*}
where $e\in\{0,1,\ldots,k\}$ is the number of even numbers among $N_1,\ldots,N_k$.
\end{theorem}

\begin{proof}
This is something that we already know, up to order 3, and the next coefficient $K_4$ can be computed in a similar way, based on results that we already have. 
\end{proof}

The passage to Theorem 11.23 is quite interesting, because it shows that the glow of the Fourier matrices $F_G$ is not polynomial in $N=|G|$. When restricting the attention to the usual Fourier matrices $F_N$, the glow up to order 4 is polynomial both in $N$ odd, and in $N$ even, but it is not clear what happens at higher order. An interesting question here is that of computing the complex glow of the Walsh matrices. Indeed, for the Walsh matrices the integrals $I(\pi)$, and hence the glow itself, might be polynomial in $N$.

\index{Walsh matrix}

\section*{11e. Exercises} 

There had been a lot of advanced combinatorics and probability in this chapter, and our exercises here will be the most about this, advanced combinatorics and probability. Let us start however with a very standard exercise, as follows:

\begin{exercise}
Establish the M\"obius inversion formula, namely
$$f(\sigma)=\sum_{\pi\leq\sigma}g(\pi)
\quad\implies\quad g(\sigma)=\sum_{\pi\leq\sigma}\mu(\pi,\sigma)f(\pi)$$
for the functions on $P(p)$.
\end{exercise}

The idea here is that the formula on the left should normally allow the computation of $g$ in terms of $f$, by some kind of recurrence. And the point is that when working out the coefficients, we are normally led to the recurrence formula for the M\"obius function. As a bonus exercise, try to find as well some basic applications of this.

\begin{exercise}
Prove that the inverse of the adjacency matrix of $P(k)$, given by
$$A_k(\pi,\sigma)=\begin{cases}
1&{\rm if}\ \pi\leq\sigma\\
0&{\rm if}\ \pi\not\leq\sigma
\end{cases}$$
is the M\"obius matrix of $P$, given by $M_k(\pi,\sigma)=\mu(\pi,\sigma)$.
\end{exercise}

This exercise is indeed equivalent to the first exercice, and with this equivalence being an instructive preliminary exercise. As for the proof, the idea here is that the matrix $A_k$ is upper triangular, with respect to a suitably chosen order on the partitions, that you will have to find, and so when inverting, we are led into the above recurrence for $\mu$.

\begin{exercise}
Prove that given independent normal variables $x,y$, by setting
$$z=\frac{1}{\sqrt{2}}(x+iy)$$
the even moments of the variable $|z|$ are given by the formula $\mathbb E(|z|^{2p})=p!$.
\end{exercise}

This is something well-known, that we have been heavily using in the above. As for the proof of this fact, this depends on your knowledge of calculus.

\begin{exercise}
Establish the following formulae,
$$F_{ix}F_{iy}=F_{i,x+y}\quad,\quad 
\overline{F}_{ix}=F_{i,-x}\quad,\quad 
\sum_xF_{ix}=N\delta_{i0}$$
valid for any generalized Fourier matrix, $F=F_G$.
\end{exercise}

As before with the previous exercise, this is something well-known, that we have been heavily using in the above. As for the proof, this should not be difficult.

\begin{exercise}
Compute the glow of the Walsh matrices 
$$W_N=F_2^{\otimes n}$$ 
with $N=2^n$, and check if this glow is polynomial or not in $N$.
\end{exercise}

There are some interesting computations here, and as before with previous research-level exercises, doing them at least partly, or even very partly, can be source of joy.

\chapter{Local estimates}

\section*{12a. Norm maximizers}

We discuss here some further analytic questions, regarding the complex Hadamard matrices, following \cite{bn3}, in analogy with the considerations from chapter 3. We will be interested in the complex analogue of the notion of almost Hadamard matrix. This looks more as a routine topic, and for a long time it was believed that there is no hurry in developing all this, since complex Hadamard matrices exist anyway at any $N\in\mathbb N$, and so there is no really need for almost Hadamard matrices, in the complex setting.

\bigskip

However, some work on this subject was eventually done in \cite{bn3}, and surprise, it turned out that, at least conjecturally, there are no almost Hadamard matrices, in the complex sense. Which is very good news, because this shows, again conjecturally, that for a matrix $H\in\sqrt{N}U_N$, the property of being complex Hadamard is ``local''. Which itself is a surprising and potentially far-reaching statement, suggesting reformulating all the Hadamard matrix problematics, including the HC and CHC, in local terms.

\bigskip

We will explain this exciting material in this chapter. To start with, we have the following basic estimate, that we already know, from chapter 11:

\index{convex function}
\index{concave function}
\index{Jensen inequality}

\begin{proposition}
Given $\psi:[0,\infty)\to\mathbb R$, the following function over $U_N$,
$$F(U)=\sum_{ij}\psi(|U_{ij}|^2)$$
satisfies the following inequality, when $\psi$ is convex/concave,
$$F(U)\geq N^2\psi\left(\frac{1}{N}\right)\quad\big/\quad F(U)\leq N^2\psi\left(\frac{1}{N}\right)$$
and assuming that $\psi$ is strictly convex/concave, the equality case appears precisely for the rescaled Hadamard matrices, $U=H/\sqrt{N}$ with $H\in M_N(\mathbb T)$ Hadamard.
\end{proposition}

\begin{proof}
This follows indeed from the Jensen inequality applied to the function in the statement, exactly as in the real case, as explained in chapter 2.
\end{proof}

Of particular interest for us are the power functions $\psi(x)=x^{p/2}$, which are concave at $p\in[1,2)$, and convex at $p\in(2,\infty)$. These lead to the following statement:

\index{norm maximizer}
\index{norm minimizer}

\begin{theorem}
Let $U\in U_N$, and set $H=\sqrt{N}U$.
\begin{enumerate}
\item For $p\in[1,2)$ we have $||U||_p\leq N^{2/p-1/2}$,

\item For $p\in(2,\infty]$ we have $||U||_p\geq N^{2/p-1/2}$.
\end{enumerate}
In both cases, the equality situation happens precisely when $H$ is Hadamard.
\end{theorem}

\begin{proof}
Consider indeed the $p$-norm on $U_N$, which at $p\in[1,\infty)$ is given by:
$$||U||_p=\left(\sum_{ij}|U_{ij}|^p\right)^{1/p}$$

By the above discussion, involving the functions $\psi(x)=x^{p/2}$, Proposition 12.1 applies and gives the results at $p\in[1,\infty)$, the precise estimates being as follows:
$$||U||_p:
\begin{cases}
\leq N^{2/p-1/2}&{\rm if}\ p<2\\
=N^{1/2}&{\rm if}\ p=2\\
\geq N^{2/p-1/2}&{\rm if}\ p>2
\end{cases}$$

As for the case $p=\infty$, this follows with $p\to\infty$, or directly via Cauchy-Schwarz.
\end{proof}

For future reference, let us record as well the particular cases $p=1,4,\infty$ of the above result, that we already met before, and which are of particular interest:

\begin{theorem}
For any matrix $U\in U_N$ we have the estimates
$$||U||_1\leq N\sqrt{N}\quad,\quad
||U||_4\geq 1\quad,\quad 
||U||_\infty\geq\frac{1}{\sqrt{N}}$$
which in terms of the rescaled matrix $H=\sqrt{N}U$ read
$$||H||_1\leq N^2\quad,\quad
||H||_4\geq\sqrt{N}\quad,\quad 
||H||_\infty\geq1$$
and in each case, the equality case holds when $H$ is Hadamard.
\end{theorem}

\begin{proof}
These results follow from Theorem 12.2 at $p=1,4,\infty$, with the remark that for each of these particular exponents, we do not really need the H\"older inequality, with a basic application of the Cauchy-Schwarz inequality doing the job.
\end{proof}

The above results suggest the following definition:

\index{AHM}
\index{almost Hadamard matrix}
\index{p-almost Hadamard matrix}
\index{absolute almost Hadamard matrix}
\index{p-AHM}
\index{absolute AHM}

\begin{definition}
Given $U\in U_N$, the matrix $H=\sqrt{N}U$ is called:
\begin{enumerate}
\item Almost Hadamard, if $U$ locally maximizes the $1$-norm on $U_N$.

\item $p$-almost Hadamard, with $p<2$, if $U$ locally maximizes the $p$-norm on $U_N$.

\item $p$-almost Hadamard, with $p>2$, if $U$ locally minimizes the $p$-norm on $U_N$.

\item Absolute almost Hadamard, if it is $p$-almost Hadamard at any $p\neq2$.
\end{enumerate}
We have as well real versions of these notions, with $U_N$ replaced by $O_N$.
\end{definition}

All this might seem a bit complicated, but this is the best way of presenting things. We are mainly interested in (1), but as explained in chapter 9, the exponent $p=4$ from (3) is interesting as well, and once we have (3) we must formulate (2) as well, and finally (4) is a useful thing too, because the absolute case is sometimes easier to study. As for the ``doubling'' of all these notions, via the last sentence, this is necessary too, because given a function $F:U_N\to\mathbb R$, an element $U\in O_N$ can be a local extremum of the restriction $F_{|O_N}:O_N\to\mathbb R$,  but not of the function $F$ itself, and we will see examples of this.

\bigskip

Let us first study the critical points. Things are quite tricky here, and complete results are available so far only at $p=1$. Following \cite{bn3}, we first have the following result:

\index{rotation trick}
\index{local maximizer}
\index{norm maximizer}

\begin{theorem}
If $U\in U_N$ locally maximizes the $1$-norm, then 
$$U_{ij}\neq 0$$
must hold for any $i,j$.
\end{theorem}

\begin{proof}
We use the same method as in the real case, namely a rotation trick. Let us denote by $U_1,\ldots,U_N$ the rows of $U$, and let us perform a rotation of $U_1,U_2$:
$$\begin{bmatrix}U^t_1\\ U^t_2\end{bmatrix}
=\begin{bmatrix}\cos t\cdot U_1-\sin t\cdot U_2\\ \sin t\cdot U_1+\cos t\cdot U_2\end{bmatrix}$$

In order to compute the 1-norm, let us permute the columns of $U$, in such a way that the first two rows look as follows, with $X,Y,A,B$ having nonzero entries:
$$\begin{bmatrix}U_1\\ U_2\end{bmatrix}
=\begin{bmatrix}0&0&Y&A\\0&X&0&B\end{bmatrix}$$

The rotated matrix will look then as follows:
$$\begin{bmatrix}U_1^t\\ U_2^t\end{bmatrix}
=\begin{bmatrix}
0&-\sin t\cdot X&\cos t\cdot Y&\cos t\cdot A-\sin t\cdot B\\
0&\cos t\cdot X&\sin t\cdot y&\sin t\cdot A+\cos t\cdot B\end{bmatrix}$$

Our claim is that $X,Y$ must be empty. Indeed, if $A$ and $B$ are not empty, let us fix a column index $k$ for both $A,B$, and set $\alpha=A_k$, $\beta=B_k$. We have then:
\begin{eqnarray*}
|(U_1^t)_k|+|(U_2^t)_k|
&=&|\cos t\cdot\alpha-\sin t\cdot\beta|+|\sin t\cdot\alpha+\cos t\cdot\beta|\\
&=&\sqrt{\cos^2t\cdot|\alpha|^2+\sin^2t\cdot|\beta|^2-\sin t\cos t(\alpha\bar{\beta}+\beta\bar{\alpha})}\\
&+&\sqrt{\sin^2t\cdot|\alpha|^2+\cos^2t\cdot|\beta|^2+\sin t\cos t(\alpha\bar{\beta}+\beta\bar{\alpha})}
\end{eqnarray*}

Since $\alpha,\beta\neq 0$, the above function is differentiable at $t=0$, and we obtain:
\begin{eqnarray*}
\frac{d\left(|(U_1^t)_k|+|(U_2^t)_k|\right)}{dt}
&=&\frac{\sin 2t(|\beta|^2-|\alpha|^2)-\cos 2t(\alpha\bar{\beta}+\beta\bar{\alpha})}{2\sqrt{\cos^2t\cdot|\alpha|^2+\sin^2t\cdot|\beta|^2-\sin t\cos t(\alpha\bar{\beta}+\beta\bar{\alpha})}}\\
&+&\frac{\sin 2t(|\alpha|^2-|\beta|^2)+\cos 2t(\alpha\bar{\beta}+\beta\bar{\alpha})}{2\sqrt{\sin^2t\cdot|\alpha|^2+\cos^2t\cdot|\beta|^2+\sin t\cos t(\alpha\bar{\beta}+\beta\bar{\alpha})}}
\end{eqnarray*}

Thus at $t=0$, we obtain the following formula:
$$\frac{d\left(|(U_1^t)_k|+|(U_2^t)_k|\right)}{dt}(0)=\frac{\alpha\bar{\beta}+\beta\bar{\alpha}}{2}\left(\frac{1}{|\beta|}-\frac{1}{|\alpha|}\right)$$

Now since our matrix $U$ locally maximizes the 1-norm, both directional derivatives of $||U^t||_1$ must be negative in the limit $t\to 0$. On the other hand, if we denote by $C$ the contribution coming from the right, which might be zero in the case where $A$ and $B$ are empty, i.e. the sum over $k$ of the above quantities, we have:
\begin{eqnarray*}
\frac{d||U^t||_1}{dt}_{\big|t=0^+}
&=&\frac{d}{dt}_{\big|t=0^+}(|\cos t|+|\sin t|)(||X||_1+||Y||_1)+C\\
&=&(-\sin t+\cos t)_{\big|t=0}(||X||_1+||Y||_1)+C\\
&=&||X||_1+||Y||_1+C
\end{eqnarray*}

As for the derivative at left, this is given by the following formula:
\begin{eqnarray*}
\frac{d||U^t||_1}{dt}_{\big|t=0^-}
&=&\frac{d}{dt}_{\big|t=0^-}(|\cos t|+|\sin t|)(||X||_1+||Y||_1)+C\\
&=&(-\sin t-\cos t)_{\big|t=0}(||X||_1+||Y||_1)+C\\
&=&-||X||_1-||Y||_1+C
\end{eqnarray*}

We therefore obtain the following inequalities, where $C$ is as above:
\begin{eqnarray*}
||X||_1+||Y||_1+C &\leq& 0\\
-||X||_1-||Y||_1+C&\leq& 0
\end{eqnarray*}

Consider now the matrix obtained from $U$ by interchanging $U_1,U_2$. Since this matrix must be as well a local maximizer of the 1-norm, and since the above formula shows that $C$ changes its sign when interchanging $U_1,U_2$, we obtain:
\begin{eqnarray*}
||X||_1+||Y||_1-C &\leq& 0\\
-||X||_1-||Y||_1-C&\leq& 0
\end{eqnarray*}

The four inequalities that we have give altogether the following conclusion:
$$||X||_1+||Y||_1=C=0$$

Now from $||X||_1+||Y||_1=0$ we obtain that both the vectors $X,Y$ must be empty, as claimed. As a conclusion, up to a permutation of the columns, the first two rows of our matrix $U$ must be of the following form, with $A,B$ having only nonzero entries:
$$\begin{bmatrix}U_1\\ U_2\end{bmatrix}
=\begin{bmatrix}0&A\\0&B\end{bmatrix}$$

By permuting the rows of $U$, the same must hold for any two rows $U_i,U_j$. Now since $U$ cannot have a zero column, we conclude that $U$ cannot have zero entries, as claimed.
\end{proof}

Let us compute now the critical points. Following \cite{bn3}, we have:

\index{critical point}
\index{Lagrange multipliers}

\begin{theorem}
Let $\varphi:[0,\infty)\to\mathbb R$ be a differentiable function. A unitary matrix with nonzero entries $U\in U_N^*$ is a critical point of the quantity
$$F(U)=\sum_{ij}\varphi(|U_{ij}|)$$
precisely when $WU^*$ is self-adjoint, where $W_{ij}={\rm sgn}(U_{ij})\varphi'(|U_{ij}|)$.
\end{theorem}

\begin{proof}
Again, this follows like in the real case, by performing modifications where needed. We regard $U_N$ as a real algebraic manifold, with coordinates $U_{ij},\bar{U}_{ij}$. This manifold consists by definition of the zeroes of the following polynomials: $$A_{ij}=\sum_kU_{ik}\bar{U}_{jk}-\delta_{ij}$$

A given matrix $U\in U_N$ is then a critical point of $F$ precisely when $dF\in span(dA_{ij})$. Regarding the space $span(dA_{ij})$, this consists of the following quantities:
\begin{eqnarray*}
\sum_{ij}M_{ij}dA_{ij}
&=&\sum_{ijk}M_{ij}(U_{ik}d\bar{U}_{jk}+\bar{U}_{jk}dU_{ik})\\
&=&\sum_{jk}(M^tU)_{jk}d\bar{U}_{jk}+\sum_{ik}(M\bar{U})_{ik}dU_{ik}\\
&=&\sum_{ij}(M^tU)_{ij}d\bar{U}_{ij}+\sum_{ij}(M\bar{U})_{ij}dU_{ij}
\end{eqnarray*}

In order to compute $dF$, observe first that, with $S_{ij}=sgn(U_{ij})$, we have:
\begin{eqnarray*}
d|U_{ij}|
&=&d\sqrt{U_{ij}\bar{U}_{ij}}\\
&=&\frac{U_{ij}d\bar{U}_{ij}+\bar{U}_{ij}dU_{ij}}{2|U_{ij}|}\\
&=&\frac{1}{2}(S_{ij}d\bar{U}_{ij}+\bar{S}_{ij}dU_{ij})
\end{eqnarray*}

In terms of $W_{ij}=sgn(U_{ij})\varphi'(|U_{ij}|)$, as in the statement, we obtain:
\begin{eqnarray*}
dF
&=&\sum_{ij}d\left(\varphi(|U_{ij}|)\right)\\
&=&\sum_{ij}\varphi'(|U_{ij}|)d|U_{ij}|\\
&=&\frac{1}{2}\sum_{ij}W_{ij}d\bar{U}_{ij}+\bar{W}_{ij}dU_{ij}
\end{eqnarray*}

We conclude that $U\in U_N$ is a critical point of $F$ if and only if there exists a matrix $M\in M_N(\mathbb C)$ such that the following two conditions are satisfied:
$$W=2M^tU\quad,\quad 
\bar{W}=2M\bar{U}$$

But this means $WU^*=UW^*$, and so that $WU^*$ must be self-adjoint, as claimed.
\end{proof}

\section*{12b. Balanced matrices}

In order to process the above result, we proceed exactly as in chapter 3, by adding some complex conjugates where needed. We can use the following notion:

\index{color decomposition}
\index{balanced matrix}
\index{semi-balanced matrix}

\begin{definition}
Given $U\in U_N$, we consider its ``color decomposition'' 
$$U=\sum_{r>0}rU_r$$
with $U_r\in M_N(\mathbb T\cup\{0\})$ containing the phase components at $r>0$, and we call $U$:
\begin{enumerate}
\item Semi-balanced, if $U_rU^*$ and $U^*U_r$, with $r>0$, are all self-adjoint.

\item Balanced, if $U_rU_s^*$ and $U_r^*U_s$, with $r,s>0$, are all self-adjoint.
\end{enumerate}
\end{definition}

These conditions are quite natural, because for a unitary matrix $U\in U_N$, the relations $UU^*=U^*U=1$ translate as follows, in terms of the color decomposition:
$$\sum_{r>0}rU_rU^*=\sum_{r>0}rU^*U_r=1$$
$$\sum_{r,s>0}rsU_rU_s^*=\sum_{r,s>0}rsU_r^*U_s=1$$

Thus, our balancing conditions express the fact that the various components of the above sums all self-adjoint. Now back to our critical point questions, we have:

\begin{theorem}
For a matrix $U\in U_N^*$, the following are equivalent:
\begin{enumerate}
\item $U$ is a critical point of $F(U)=\sum_{ij}\varphi(|U_{ij}|)$, for any $\varphi:[0,\infty)\to\mathbb R$.

\item $U$ is a critical point of all the $p$-norms, with $p\in[1,\infty)$.

\item $U$ is semi-balanced, in the above sense.
\end{enumerate}
\end{theorem}

\begin{proof}
We use Theorem 12.6. The matrix constructed there is given by:
\begin{eqnarray*}
(WU^*)_{ij}
&=&\sum_k{\rm sgn}(U_{ik})\varphi'(|U_{ik}|)\bar{U}_{jk}\\
&=&\sum_{r>0}\varphi'(r)\sum_{k,|U_{ik}|=r}{\rm sgn}(U_{ik})\bar{U}_{jk}\\
&=&\sum_{r>0}\varphi'(r)\sum_k(U_r)_{ik}\bar{U}_{jk}\\
&=&\sum_{r>0}\varphi'(r)(U_rU^*)_{ij}
\end{eqnarray*}

We conclude that we have the following formula for this matrix:
$$WU^*=\sum_{r>0}\varphi'(r)U_rU^*$$

Now when $\varphi:[0,\infty)\to\mathbb R$ varies, as a differentiable function, or as a power function $\varphi(x)=x^p$ with $p\in[1,\infty)$, the individual components must be self-adjoint, as desired.
\end{proof}

In practice now, most of the known examples of semi-balanced matrices are actually balanced. We have the following collection of simple facts, regarding such matrices:

\begin{proposition}
The class of balanced matrices is as follows:
\begin{enumerate}
\item It contains the matrices $U=H/\sqrt{N}$, with $H\in M_N(\mathbb C)$ Hadamard.

\item It is stable under transposition, complex conjugation, and taking adjoints.

\item It is stable under taking tensor products.

\item It is stable under the Hadamard equivalence relation.

\item It contains the matrix $V_N=\frac{1}{N}(2\mathbb I_N-N1_N)$, where $\mathbb I_N$ is the all-$1$ matrix.
\end{enumerate}
\end{proposition}

\begin{proof}
All these results are elementary, the proof being as follows:

\medskip

(1) Here $U\in U_N$ follows from the Hadamard condition, and since there is only one color component, namely $U_{1/\sqrt{N}}=H$, the balancing condition is satisfied as well.

\medskip

(2) Assuming that $U=\sum_{r>0}rU_r$ is a color decomposition of a given matrix $U\in U_N$, the following are color decompositions too, and this gives the assertions:
$$U^t=\sum_{r>0}rU_r^t\quad,\quad 
\bar{U}=\sum_{r>0}r\bar{U}_r\quad,\quad 
U^*=\sum_{r>0}rU_r^*$$

(3) Assuming that $U=\sum_{r>0}rU_r$ and $V=\sum_{s>0}sV_s$ are the color decompositions of two given unitary matrices $U,V$, we have the following formula:
\begin{eqnarray*}
U\otimes V
&=&\sum_{r,s>0}rs\cdot U_r\otimes V_s\\
&=&\sum_{p>0}p\sum_{p=rs}U_r\otimes V_s
\end{eqnarray*}

Thus the color components of $W=U\otimes V$ are the following matrices: 
$$W_p=\sum_{p=rs}U_r\otimes V_s$$

It follows that if $U,V$ are both balanced, then so is $W=U\otimes V$.

\medskip

(4) We recall that the Hadamard equivalence consists in permuting rows and columns, and switching signs on rows and columns. Since all these operations correspond to certain conjugations at the level of the matrices $U_rU_s^*,U_r^*U_s$, we obtain the result.

\medskip

(5) The matrix in the statement, which goes back to \cite{bnz}, is as follows:
$$V_N=\frac{1}{N}
\begin{pmatrix}
2-N&2&\ldots&2\\
2&2-N&\ldots&2\\
\ldots&\ldots&\ldots&\ldots\\
2&2&\ldots&2-N
\end{pmatrix}$$

Observe that this matrix is indeed unitary, its rows being of norm one, and pairwise orthogonal. The color components of this matrix are:
$$V_{2/N-1}=1_N\quad,\quad 
V_{2/N}=\mathbb I_N-1_N$$

It follows that this matrix is balanced as well, as claimed.
\end{proof}

Let us look now more in detail at $V_N$, and at the matrices having similar properties. Following \cite{bnz}, let us call $(a,b,c)$ pattern any matrix $M\in M_N(0,1)$, with $N=a+2b+c$, such that any two rows look as follows, up to a permutation of the columns:
$$\begin{matrix}
0\ldots 0&0\ldots 0&1\ldots 1&1\ldots 1\\
\underbrace{0\ldots 0}_a&\underbrace{1\ldots 1}_b&\underbrace{0\ldots 0}_b&\underbrace{1\ldots 1}_c
\end{matrix}$$

As explained in \cite{bnz}, there are many interesting examples of $(a,b,c)$ patterns, coming from the balanced incomplete block designs (BIBD), and all these examples can produce two-entry unitary matrices, by replacing the $0,1$  entries with suitable numbers $x,y$. Now back to the matrix $V_N$ from Proposition 12.9 (5), observe that this matrix comes from a $(0,1,N-2)$ pattern, in the above sense. And also, independently of this, this matrix has the remarkable property of being at the same time circulant and self-adjoint. We have in fact the following result, generalizing Proposition 12.9 (5):

\index{block design}

\begin{theorem}
The following matrices are balanced:
\begin{enumerate}
\item The orthogonal matrices coming from $(a,b,c)$ patterns.

\item The unitary matrices which are circulant and self-adjoint.
\end{enumerate}
\end{theorem}

\begin{proof}
These observations basically go back to \cite{bnz}, the proofs being as follows:

\medskip

(1) If we denote by $P,Q\in M_N(0,1)$ the matrices describing the positions of the $0,1$ entries inside the pattern, then we have the following formulae:
\begin{eqnarray*}
PP^t=P^tP&=&a\mathbb I_N+b1_N\\
QQ^t=Q^tQ&=&c\mathbb I_N+b1_N\\
PQ^t=P^tQ=QP^t=Q^tP&=&b\mathbb I_N-b1_N
\end{eqnarray*}

Since all these matrices are symmetric, $U$ is balanced, as claimed.

\medskip

(2) Assume that $U\in U_N$ is circulant, $U_{ij}=\gamma_{j-i}$, and in addition self-adjoint, which means $\bar{\gamma}_i=\gamma_{-i}$. Consider the following sets, which must satisfy $D_r=-D_r$:
$$D_r=\{k:|\gamma_r|=k\}$$

In terms of these sets, we have the following formula:
\begin{eqnarray*}
(U_rU_s^*)_{ij}
&=&\sum_k(U_r)_{ik}(\bar{U}_s)_{jk}\\
&=&\sum_k\delta_{|\gamma_{k-i}|,r}\,{\rm sgn}(\gamma_{k-i})\cdot\delta_{|\gamma_{k-j}|,s}\,{\rm sgn}(\bar{\gamma}_{k-j})\\
&=&\sum_{k\in(D_r+i)\cap(D_s+j)}{\rm sgn}(\gamma_{k-i})\,{\rm sgn}(\bar{\gamma}_{k-j})
\end{eqnarray*}

With $k=i+j-m$ we obtain, by using $D_r=-D_r$, and then $\bar{\gamma}_i=\gamma_{-i}$:
\begin{eqnarray*}
(U_rU_s^*)_{ij}
&=&\sum_{m\in(-D_r+j)\cap(-D_s+i)}{\rm sgn}(\gamma_{j-m})\,{\rm sgn}(\bar{\gamma}_{i-m})\\
&=&\sum_{m\in(D_r+i)\cap(D_r+j)}{\rm sgn}(\gamma_{j-m})\,{\rm sgn}(\bar{\gamma}_{i-m})\\
&=&\sum_{m\in(D_r+i)\cap(D_r+j)}{\rm sgn}(\bar{\gamma}_{m-j})\,{\rm sgn}(\gamma_{m-i})
\end{eqnarray*}

Now by interchanging $i\leftrightarrow j$, and with $m\to k$, this formula becomes:
$$(U_rU_s^*)_{ji}=\sum_{k\in(D_r+i)\cap(D_r+j)}{\rm sgn}(\bar{\gamma}_{k-i})\,{\rm sgn}(\gamma_{k-j})$$

We recognize here the complex conjugate of $(U_rU_s^*)_{ij}$, as previously computed above, and we therefore deduce that $U_rU_s^*$ is self-adjoint. The proof for $U_r^*U_s$ is similar. 
\end{proof}

\section*{12c. Hessian computations}

Let us compute now derivatives. As in Theorem 12.6, it is convenient to do the computations in a more general framework, where we have a function as follows:
$$F(U)=\sum_{ij}\psi(|U_{ij}|^2)$$

In order to study the local extrema of these quantities, consider the following function:
$$f(t)
=F(Ue^{tA})
=\sum_{ij}\psi(|(Ue^{tA})_{ij}|^2)$$

Here $U\in U_N$ is a unitary matrix, and $A\in M_N(\mathbb C)$ is assumed to be anti-hermitian, $A^*=-A$, as for having $e^A\in U_N$. Let us first compute the derivative of $f$. We have:

\begin{proposition}
We have the following formula,
$$f'(t)=2\sum_{ij}\psi'(|(Ue^{tA})_{ij}|^2)Re\left[(UAe^{tA})_{ij}\overline{(Ue^{tA})_{ij}}\right]$$
valid for any $U\in U_N$, and any $A\in M_N(\mathbb C)$ anti-hermitian.
\end{proposition}

\begin{proof}
The matrices $U,e^{tA}$ being both unitary, we have:
\begin{eqnarray*}
|(Ue^{tA})_{ij}|^2
&=&(Ue^{tA})_{ij}\overline{(Ue^{tA})_{ij}}\\
&=&(Ue^{tA})_{ij}((Ue^{tA})^*)_{ji}\\
&=&(Ue^{tA})_{ij}(e^{tA^*}U^*)_{ji}\\
&=&(Ue^{tA})_{ij}(e^{-tA}U^*)_{ji}
\end{eqnarray*}

We can now differentiate our function $f$, and by using once again the unitarity of the matrices $U,e^{tA}$, along with the formula $A^*=-A$, we obtain:
\begin{eqnarray*}
f'(t)
&=&\sum_{ij}\psi'(|(Ue^{tA})_{ij}|^2)\left[(UAe^{tA})_{ij}(e^{-tA}U^*)_{ji}-(Ue^{tA})_{ij}(e^{-tA}AU^*)_{ji}\right]\\
&=&\sum_{ij}\psi'(|(Ue^{tA})_{ij}|^2)\left[(UAe^{tA})_{ij}\overline{((e^{-tA}U^*)^*)_{ij}}-(Ue^{tA})_{ij}\overline{((e^{-tA}AU^*)^*)_{ij}}\right]\\
&=&\sum_{ij}\psi'(|(Ue^{tA})_{ij}|^2)\left[(UAe^{tA})_{ij}\overline{(Ue^{tA})_{ij}}+(Ue^{tA})_{ij}\overline{(UAe^{tA})_{ij}}\right]
\end{eqnarray*}

But this gives the formula in the statement, and we are done.
\end{proof}

Before computing the second derivative, let us evaluate $f'(0)$. We have:

\begin{proposition}
We have the following formula,
$$f'(0)=2\sum_{r>0}r\psi'(r^2)Re\left[Tr(U_r^*UA)\right]$$
where the matrices $U_r\in M_N(\mathbb T\cup\{0\})$ are the color components of $U$.
\end{proposition}

\begin{proof}
We use the formula in Proposition 12.11. At $t=0$, we obtain:
$$f'(0)=2\sum_{ij}\psi'(|U_{ij}|^2)Re\left[(UA)_{ij}\overline{U}_{ij}\right]$$

Consider now the color decomposition of $U$. We have the following formulae:
\begin{eqnarray*}
U_{ij}=\sum_{r>0}r(U_r)_{ij}
&\implies&|U_{ij}|^2=\sum_{r>0}r^2|(U_r)_{ij}|\\
&\implies&\psi'(|U_{ij}|^2)=\sum_{r>0}\psi'(r^2)|(U_r)_{ij}|
\end{eqnarray*}

Now by getting back to the above formula of $f'(0)$, we obtain:
$$f'(0)=2\sum_{r>0}\psi'(r^2)\sum_{ij}Re\left[(UA)_{ij}\overline{U}_{ij}|(U_r)_{ij}|\right]$$

Our claim now is that we have the following formula:
$$\overline{U}_{ij}|(U_r)_{ij}|=r\overline{(U_r)}_{ij}$$

Indeed, in the case $|U_{ij}|\neq r$ this formula reads $\overline{U}_{ij}\cdot 0=r\cdot 0$, which is true, and in the case $|U_{ij}|=r$ this formula reads $r\bar{S}_{ij}\cdot 1=r\cdot\bar{S}_{ij}$, which is once again true. Thus:
$$f'(0)=2\sum_{r>0}r\psi'(r^2)\sum_{ij}Re\left[(UA)_{ij}\overline{(U_r)}_{ij}\right]$$

But this gives the formula in the statement, and we are done.
\end{proof}

Let us compute now the second derivative. The result here is as follows:

\begin{proposition}
We have the following formula,
\begin{eqnarray*}
f''(0)
&=&4\sum_{ij}\psi''(|U_{ij}|^2)Re\left[(UA)_{ij}\overline{U}_{ij}\right]^2\\
&&+2\sum_{ij}\psi'(|U_{ij}|^2)Re\left[(UA^2)_{ij}\overline{U}_{ij}\right]\\
&&+2\sum_{ij}\psi'(|U_{ij}|^2)|(UA)_{ij}|^2
\end{eqnarray*}
valid for any $U\in U_N$, and any $A\in M_N(\mathbb C)$ anti-hermitian.
\end{proposition}

\begin{proof}
We use the formula in Proposition 12.11, namely:
$$f'(t)=2\sum_{ij}\psi'(|(Ue^{tA})_{ij}|^2)Re\left[(UAe^{tA})_{ij}\overline{(Ue^{tA})_{ij}}\right]$$

Since the real part on the right, or rather its double, appears as the derivative of the quantity $|(Ue^{tA})_{ij}|^2$, when differentiating a second time, we obtain:
\begin{eqnarray*}
f''(t)
&=&4\sum_{ij}\psi''(|(Ue^{tA})_{ij}|^2)Re\left[(UAe^{tA})_{ij}\overline{(Ue^{tA})_{ij}}\right]^2\\
&&+2\sum_{ij}\psi'(|(Ue^{tA})_{ij}|^2)Re\left[(UAe^{tA})_{ij}\overline{(Ue^{tA})_{ij}}\right]'
\end{eqnarray*}

In order to compute now the missing derivative, observe that we have:
\begin{eqnarray*}
\left[(UAe^{tA})_{ij}\overline{(Ue^{tA})_{ij}}\right]'
&=&(UA^2e^{tA})_{ij}\overline{(Ue^{tA})_{ij}}+(UAe^{tA})_{ij}\overline{(UAe^{tA})_{ij}}\\
&=&(UA^2e^{tA})_{ij}\overline{(Ue^{tA})_{ij}}+|(UAe^{tA})_{ij}|^2
\end{eqnarray*}

Summing up, we have obtained the following formula:
\begin{eqnarray*}
f''(t)
&=&4\sum_{ij}\psi''(|(Ue^{tA})_{ij}|^2)Re\left[(UAe^{tA})_{ij}\overline{(Ue^{tA})_{ij}}\right]^2\\
&&+2\sum_{ij}\psi'(|(Ue^{tA})_{ij}|^2)Re\left[(UA^2e^{tA})_{ij}\overline{(Ue^{tA})_{ij}}\right]\\
&&+2\sum_{ij}\psi'(|(Ue^{tA})_{ij}|^2)|(UAe^{tA})_{ij}|^2
\end{eqnarray*}

But at $t=0$ this gives the formula in the statement, and we are done.
\end{proof}

By using the function $\psi(x)=\sqrt{x}$, corresponding to $F(U)=||U||_1$, we obtain:

\begin{proposition}
Let $U \in U_N^*$. For the function $F(U)=||U||_1$ we have the formula
$$f''(0)=Re\left[Tr(S^*UA^2)\right]+\sum_{ij}\frac{Im\left[(UA)_{ij}\overline{S}_{ij}\right]^2}{|U_{ij}|}$$
valid for any anti-hermitian matrix $A$, where $U_{ij}=S_{ij}|U_{ij}|$.
\end{proposition}

\begin{proof}
We use the formula in Proposition 12.13, with the following data:
$$\psi(x)=\sqrt{x}\quad,\quad 
\psi'(x)=\frac{1}{2\sqrt{x}}\quad,\quad 
\psi''(x)=-\frac{1}{4x\sqrt{x}}$$

We obtain the following formula:
\begin{eqnarray*}
f''(0)
&=&-\sum_{ij}\frac{Re\left[(UA)_{ij}\overline{U}_{ij}\right]^2}{|U_{ij}|^3}
+\sum_{ij}\frac{Re\left[(UA^2)_{ij}\overline{U}_{ij}\right]}{|U_{ij}|}
+\sum_{ij}\frac{|(UA)_{ij}|^2}{|U_{ij}|}\\
&=&-\sum_{ij}\frac{Re\left[(UA)_{ij}\overline{S}_{ij}\right]^2}{|U_{ij}|}
+\sum_{ij}Re\left[(UA^2)_{ij}\overline{S}_{ij}\right]
+\sum_{ij}\frac{|(UA)_{ij}|^2}{|U_{ij}|}\\
&=&Re\left[Tr(S^*UA^2)\right]+\sum_{ij}\frac{|(UA)_{ij}|^2-Re\left[(UA)_{ij}\overline{S}_{ij}\right]^2}{|U_{ij}|}
\end{eqnarray*}

But this gives the formula in the statement, and we are done.
\end{proof}

We are therefore led to the following result, regarding the 1-norm:

\begin{theorem}
A matrix $U\in U_N^*$ locally maximizes the one-norm on $U_N$ precisely when $S^*U$ is self-adjoint, where $S_{ij}={\rm sgn}(U_{ij})$, and when
$$Tr(S^*UA^2)+\sum_{ij}\frac{Im\left[(UA)_{ij}\overline{S}_{ij}\right]^2}{|U_{ij}|}\leq0$$
holds, for any anti-hermitian matrix $A\in M_N(\mathbb C)$.
\end{theorem}

\begin{proof}
According to Theorem 12.6 and Proposition 12.14, the local maximizer condition requires $X=S^*U$ to be self-adjoint, and the following inequality to be satisfied:
$$Re\left[Tr(S^*UA^2)\right]+\sum_{ij}\frac{Im\left[(UA)_{ij}\overline{S}_{ij}\right]^2}{|U_{ij}|}\leq0$$

Now observe that since both $X$ and $A^2$ are self-adjoint, we have:
\begin{eqnarray*}
Re\left[Tr(XA^2)\right]
&=&\frac{1}{2}\left[Tr(XA^2)+Tr(A^2X)\right]\\
&=&Tr(XA^2)
\end{eqnarray*}

Thus we can remove the real part, and we obtain the inequality in the statement.
\end{proof}

As a comment here, the above computations can be of course interpreted by using more advanced differential geometric language. The unitary group $U_N$ is a Lie group, and its tangent space at $U\in U_N$ is isomorphic to the corresponding Lie algebra, which consists of the anti-hermitian matrices $A\in M_N(\mathbb C)$. With this picture in hand, our formulae for $f'(0)$ translate into the fact that the gradient of the 1-norm is given by:
$$\nabla||U||_1=\frac{1}{2}(S-US^*U)$$

Regarding now the second derivative, $f''(0)$, our computations here provide us with a formula for the Hessian of the 1-norm. Indeed, with the change of variables $A=iB$ on the tangent space, the Hessian $H$ of the 1-norm is given by the following formula, where $\Phi(U,iA)$ is the quantity appearing in Theorem 12.15:
$$<B,H(B)>=-\Phi(U,B)$$

Getting back to Theorem 12.15 as stated, the story is of course not over there. In order to further improve this result, we will need the following standard fact:

\begin{proposition}
For a self-adjoint $X\in M_N(\mathbb C)$, the following are equivalent:
\begin{enumerate}
\item $Tr(XA^2)\leq0$, for any anti-hermitian matrix $A\in M_N(\mathbb C)$.

\item $Tr(XB^2)\geq0$, for any hermitian matrix $B\in M_N(\mathbb C)$.

\item $Tr(XC)\geq0$, for any positive matrix $C\in M_N(\mathbb C)$.

\item $X\geq0$.
\end{enumerate}
\end{proposition}

\begin{proof}
These equivalences are well-known, the proof being as follows:

\medskip

$(1)\implies(2)$ follows by taking $B=iA$. 

\medskip

$(2)\implies(3)$ follows by taking $C=B^2$. 

\medskip

$(3)\implies(4)$ follows by diagonalizing $X$, and then taking $C$ to be diagonal.

\medskip

$(4)\implies(1)$ is clear as well, because with $Y=\sqrt{X}$ we have:
\begin{eqnarray*}
Tr(XA^2)
&=&Tr(Y^2A^2)\\
&=&Tr(YA^2Y)\\
&=&-Tr((YA)(YA)^*)\\
&\leq&0
\end{eqnarray*}

Thus, the above four conditions are indeed equivalent.
\end{proof}

Following \cite{bn3}, we can now formulate a final result on the subject, as follows:

\index{AHM}
\index{complex AHM}
\index{local maximizer}

\begin{theorem}
Given $U\in U_N$, set $S_{ij}={\rm sgn}(U_{ij})$, and let:
$$X=S^*U$$
Then $U$ locally maximizes the $1$-norm on $U_N$ precisely when $X\geq0$, and when
$$\Phi(U,B)=Tr(XB^2)-\sum_{ij}\frac{Re\left[(UB)_{ij}\overline{S}_{ij}\right]^2}{|U_{ij}|}$$ 
is positive, for any hermitian matrix $B\in M_N(\mathbb C)$.
\end{theorem}

\begin{proof}
This follows from Theorem 12.15, by setting $A=iB$, and by using Proposition 12.16, which shows that we must have indeed $X\geq0$. 
\end{proof}

\section*{12d. The conjecture}

In relation with the above, quite surprisingly, the basic real almost Hadamard matrix $K_N$ is not an almost Hadamard matrix in the complex sense. That is, while $K_N/\sqrt{N}$ locally maximizes the 1-norm on $O_N$, it does not do so over $U_N$. Moreover, as we will see in a moment, the same happens for the other basic real almost Hadamard matrices discussed in chapter 3, such as the circulant ones, and the 2-entry ones studied there. Thus, the situation in the complex case is drastically different from the one in the real case, and we are led in this way to the following remarkable statement:

\begin{conjecture}[Almost Hadamard conjecture (AHC)]
Any local maximizer of the $1$-norm on $U_N$, 
$$||U||_1=\sum_{ij}|U_{ij}|$$
must be a global maximizer, i.e. must be a rescaled Hadamard matrix.
\end{conjecture}

\index{AHC}
\index{Almost Hadamard Conjecture}
\index{gradient descent method}

In other words, our conjecture is that, in the complex setting, almost Hadamard implies Hadamard. This would be something very useful, because we would have here a new approach to the complex Hadamard matrices, which is analytic and local. Which new approach, importantly, could potentially shed some new light on all the Hadamard matrix problems, be them real or complex, including the HC and CHC.

\bigskip

In order to explain all this, and the evidence that we have for the above conjecture, let us study more in detail the quantity $\Phi(U,B)$ from Theorem 12.17, namely:
$$\Phi(U,B)=Tr(XB^2)-\sum_{ij}\frac{Re\left[(UB)_{ij}\overline{S}_{ij}\right]^2}{|U_{ij}|}$$ 

As a first observation here, we have the following result:

\begin{proposition}
With $S_{ij}=sgn(U_{ij})$ and $X=S^*U$ as above, we have
$$\Phi(U,B)=\Phi(U,B+D)$$
for any $D\in M_N(\mathbb R)$ diagonal.
\end{proposition}

\begin{proof}
The matrices $X,B,D$ being all self-adjoint, we have:
$$(XBD)^*=DBX$$

Thus when computing $\Phi(U,B+D)$, the trace term decomposes as follows:
\begin{eqnarray*}
Tr(X(B+D)^2)
&=&Tr(XB^2)+Tr(XBD)+Tr(XDB)+Tr(XD^2)\\
&=&Tr(XB^2)+Tr(XBD)+Tr(DBX)+Tr(XD^2)\\
&=&Tr(XB^2)+2Re[Tr(XBD)]+Tr(XD^2)
\end{eqnarray*}

Regarding now the second term, in order to compute it, observe that with the notation $D=diag(\lambda_1,\ldots,\lambda_N)$, with $\lambda_i\in\mathbb R$, we have the following formula:
$$(UD)_{ij}\overline{S}_{ij}
=U_{ij}\lambda_j\overline{S}_{ij}
=\lambda_j|U_{ij}|$$

Thus the second term decomposes as follows:
\begin{eqnarray*}
&&\sum_{ij}\frac{Re\left[(UB+UD)_{ij}\overline{S}_{ij}\right]^2}{|U_{ij}|}\\
&=&\sum_{ij}\frac{Re\left[(UB)_{ij}\overline{S}_{ij}+\lambda_j|U_{ij}|\right]^2}{|U_{ij}|}\\
&=&\sum_{ij}\frac{\left[Re\left[(UB)_{ij}\overline{S}_{ij}\right]+\lambda_j|U_{ij}|\right]^2}{|U_{ij}|}\\
&=&\sum_{ij}\frac{Re\left[(UB)_{ij}\overline{S}_{ij}\right]^2}{|U_{ij}|}+2\sum_{ij}\lambda_jRe\left[(UB)_{ij}\overline{S}_{ij}\right]+\sum_{ij}\lambda_j^2|U_{ij}|
\end{eqnarray*}

Now observe that the middle term in this expression is given by:
\begin{eqnarray*}
2\sum_{ij}\lambda_jRe\left[(UB)_{ij}\overline{S}_{ij}\right]
&=&2Re\left[\sum_{ij}\lambda_j(UB)_{ij}\overline{S}_{ij}\right]\\
&=&2Re\left[\sum_{ij}(S^*)_{ji}(UB)_{ij}D_{jj}\right]\\
&=&2Re[Tr(XBD)]
\end{eqnarray*}

As for the term on the right in the above expression, this is given by:
\begin{eqnarray*}
\sum_{ij}\lambda_j^2|U_{ij}|
&=&\sum_{ij}\lambda_j^2\overline{S}_{ij}U_{ij}\\
&=&\sum_{ij}\overline{S}_{ij}(UD^2)_{ij}\\
&=&Tr(XD^2)
\end{eqnarray*}

Thus when doing the substraction we obtain $\Phi(U,B+D)=\Phi(U,B)$, as claimed.
\end{proof}

Observe that with $B=0$ we obtain $\Phi(U,D)=0$, for any $D\in M_N(\mathbb R)$ diagonal, so the inequality is Theorem 12.17 is an equality, when $B$ is diagonal. Getting now to the real thing, we have the following result, providing the first piece of evidence for the AHC:

\begin{theorem}
Consider the matrix $U=(2\mathbb I_N-N1_N)/N$. Assuming that a matrix $B\in M_N(\mathbb R)$ is symmetric and satisfies $UB=\lambda B$, we have
$$\Phi(U,B)=\lambda\cdot\frac{N-4}{2}\left[Tr(B^2)+\frac{\lambda N}{N-2}\sum_iB_{ii}^2\right]$$
and in particular, $K_N=\sqrt{N}U$ is not complex AHM at $N\neq4$, because:
\begin{enumerate}
\item For $B=\mathbb I_N$ we have the formula
$$\Phi(U,B)=\frac{N^2(N-1)(N-4)}{2(N-2)}$$
and this quantity is negative at $N=3$.

\item For $B\in M_N(\mathbb R)$ nonzero, symmetric, with $B\,\mathbb I_N=0$, $diag(B)=0$ we have 
$$\Phi(U,B)=\left(2-\frac{N}{2}\right)Tr(B^2)$$
and this quantity is negative at $N\geq5$.
\end{enumerate}
\end{theorem} 

\begin{proof}
With $U\in O_N$, $B\in M_N(\mathbb R)$, the formula in Theorem 12.17 reads:
$$\Phi(U,B)=Tr(S^tUB^2)-\sum_{ij}\frac{(UB)_{ij}^2}{|U_{ij}|}$$

Asusming now $U=\frac{1}{N}(2\mathbb I_N-N1_N)$ and $UB=\lambda B$, this formula becomes:
$$\Phi(U,B)=\lambda\left[Tr(S^tB^2)-\lambda N\sum_{ij}\frac{B_{ij}^2}{|2-N\delta_{ij}|}\right]$$

Now observe that in our case, we have the following formula:
$$\mathbb I_NB=\frac{N}{2}(U+1_N)B=\frac{(\lambda+1)N}{2}B$$

Thus the trace term is given by the following formula:
\begin{eqnarray*}
Tr(S^tB^2)
&=&Tr\left[(\mathbb I_N-21_N)B^2\right]\\
&=&\left(\frac{(\lambda+1)N}{2}-2\right)Tr(B^2)
\end{eqnarray*}

Regarding now the sum on the right, this can be computed as follows:
\begin{eqnarray*}
\sum_{ij}\frac{B_{ij}^2}{|2-N\delta_{ij}|}
&=&\sum_{ij}B_{ij}^2\left(\frac{1}{2}+\left(\frac{1}{N-2}-\frac{1}{2}\right)\delta_{ij}\right)\\
&=&\sum_{ij}B_{ij}^2\left(\frac{1}{2}-\frac{N-4}{2(N-2)}\delta_{ij}\right)\\
&=&\frac{1}{2}Tr(B^2)-\frac{N-4}{2(N-2)}\sum_iB_{ii}^2
\end{eqnarray*}

We obtain the following formula, which gives the one in the statement:
$$\Phi(U,B)=\lambda\left[\left(\frac{(\lambda+1)N}{2}-2-\frac{\lambda N}{2}\right)Tr(B^2)+\frac{\lambda N(N-4)}{2(N-2)}\sum_iB_{ii}^2\right]$$

We can now prove our various results, as follows:

\medskip

(1) For $B=\mathbb I_N$ we have $\lambda=1$, and we obtain, as claimed:
$$\Phi(U,B)
=\frac{N-4}{2}\left[N^2+\frac{N^2}{N-2}\right]
=\frac{N^2(N-4)(N-1)}{2(N-2)}$$

(2) For $B\in M_N(\mathbb R)$ nonzero, symmetric, and satisfying $B\,\mathbb I_N=0$ and $diag(B)=0$, we have $\lambda=-1$, and we obtain, as claimed:
$$\Phi(U,B)=\left(2-\frac{N}{2}\right)Tr(B^2)$$

It remains to prove that matrices $B$ as in the statement exist, at any $N\geq5$. As a first remark, such matrices cannot exist at $N=2,3$. At $N=4$, however, we have solutions, which are as follows, with $x+y+z=0$, not all zero:
$$B=\begin{pmatrix}
0&x&y&z\\
x&0&z&y\\
y&z&0&x\\
z&y&x&0
\end{pmatrix}$$

At $N\geq5$ now, we can simply use this matrix, completed with $0$ entries, and we are led to the conclusion in the statement.
\end{proof}

Let us go back now to the inequality in Theorem 12.17. When $U$ is a rescaled complex Hadamard matrix we have of course equality, and in addition, the following happens:

\begin{proposition}
For a rescaled complex Hadamard matrix, a stronger version of the inequality in Theorem 12.17, namely $\Phi(U,B)\geq0$ with
$$\Phi(U,B)=Tr(XB^2)-\sum_{ij}\frac{Re\left[(UB)_{ij}\overline{S}_{ij}\right]^2}{|U_{ij}|}$$ 
holds, with the real part replaced by the absolute value.
\end{proposition}

\begin{proof}
Indeed, for a rescaled Hadamard matrix $U=H/\sqrt{N}$ we have:
$$S=H=\sqrt{N}U$$

Thus $X=\sqrt{N}1_N$. We therefore obtain the following estimate: 
\begin{eqnarray*}
\Phi(U,B) 
&=&\sqrt{N}\left[Tr(B^2)-\sum_{ij}Re\left[(UB)_{ij}\overline{S}_{ij}\right]^2\right]\\
&\geq&\sqrt{N}\left[Tr(B^2)-\sum_{ij}|(UB)_{ij}\overline{S}_{ij}|^2\right]\\
&=&\sqrt{N}\left[Tr(B^2)-\sum_{ij}|(UB)_{ij}|^2\right]\\
&=&\sqrt{N}\left[Tr(B^2)-Tr(UB^2U^*)\right]\\
&=&0
\end{eqnarray*}

But this proves our claim, and we are done.
\end{proof}

In relation with the Tadej-\.Zyczkowski notion of defect \cite{tz2}, we have:

\index{defect}

\begin{theorem}
For a rescaled complex Hadamard matrix, the space
$$E_U=\left\{B\in M_N(\mathbb C)\Big|B=B^*,\Phi(U,B)=0\right\}$$
is isomorphic, via $B\to[(UB)_{ij}\overline{U}_{ij}]_{ij}$, to the following space:
$$D_U=\left\{A \in M_N(\mathbb R)\Big|\sum_k\bar{U}_{ki}U_{kj}(A_{ki}-A_{kj})=0,\forall i,j\right\}$$
In particular the two ``defects'' $\dim_\mathbb RE_U$ and $\dim_\mathbb RD_U$ coincide.
\end{theorem}

\begin{proof}
Since a self-adjoint matrix $B\in M_N(\mathbb C)$ belongs to $E_U$ precisely when the only inequality in the proof of Proposition 12.21 is saturated, we have:
$$E_U=\left\{B\in M_N(\mathbb C)\Big|B=B^*,Im\left[(UB)_{ij}\overline{U}_{ij}\right]=0,\forall i,j\right\}$$

The condition on the right tells us that the matrix $A=(UB)_{ij}\bar{U}_{ij}$ must be real. Now since the construction $B\to A$ is injective, we obtain an isomorphism, as follows:
$$E_U\simeq\left\{A\in M_N(\mathbb R)\Big|A_{ij}=(UB)_{ij}\bar{U}_{ij}\implies B=B^*\right\}$$

Our claim is that the space on the right is $D_U$. Indeed, let us pick $A\in M_N(\mathbb R)$. The condition $A_{ij}=(UB)_{ij}\bar{U}_{ij}$ is then equivalent to:
$$(UB)_{ij}=NU_{ij}A_{ij}$$

Thus in terms of the matrix $C_{ij}=U_{ij}A_{ij}$ we have $(UB)_{ij}=NC_{ij}$, and so: 
$$UB=NC$$

Thus we have $B=NU^*C$, and we can now perform the study of the self-adjointness condition $B=B^*$, as follows:
\begin{eqnarray*}
B=B^*
&\iff&U^*C=C^*U\\
&\iff&\sum_k\bar{U}_{ki}C_{kj}=\sum_k\bar{C}_{ki}U_{kj},\forall i,j\\
&\iff&\sum_k\bar{U}_{ki}U_{kj}A_{kj}=\sum_k\bar{U}_{ki}A_{ki}U_{kj},\forall i,j
\end{eqnarray*}

Thus we have reached to the condition defining $D_U$, and we are done. 
\end{proof}

Regarding now the known verifications of the AHC, these basically concern the natural ``candidates'' coming from Theorem 12.10, as well as some straightforward complex generalizations of these candidates. All this is quite technical, and generally speaking, we refer here to \cite{bn3}. As a first illustration, however, which is of theoretical importance, in the circulant orthogonal case, we have the following result, from \cite{bn3}:

\begin{proposition}
If $U\in O_N$ is circulant, $U_{ij}=\gamma_{j-i}$, we have
$$\Phi(U,\mathbb I_N)=Nu(Ns-uw)$$
where $u,s,w$ are the row sums of $U,S$ and $W_{ij}=\frac{1}{|U_{ij}|}$. Thus $\Phi(U,\mathbb I_N)<0$ when
$$\mathbb E({\rm sgn}(\gamma_i))<\mathbb E(\gamma_i)\,\mathbb E\left(\frac{1}{|\gamma_i|}\right)$$ 
where the symbol $\mathbb E$ stands as usual for ``average''.
\end{proposition}

\begin{proof}
We have $U\mathbb I_N=u\mathbb I_N$, which gives the following formula:
$$Tr(S^tU\mathbb I_N^2)=NuTr(S^t\mathbb I_N)=N^2us$$

Similarly, once again from $U\mathbb I_N=u\mathbb I_N$, we obtain the following formula:
$$\sum_{ij}\frac{(U\mathbb I_N)_{ij}^2}{|U_{ij}|}=\sum_{ij}\frac{u^2}{|U_{ij}|}=Nu^2w$$

By substracting, we obtain the formula in the statement, which gives the result. 
\end{proof}

Here is another exclusion criterion, also from \cite{bn3}, which is useful as well:

\begin{proposition}
If $U\in U_N$ is circulant, $U_{ij}=\gamma_{i-j}$, and self-adjoint, we have
$$\Phi(U,U)=N\left(-\frac{1}{|\gamma_0|}+\sum_i|\gamma_i|\right)$$
and so $\Phi(U,U)<0$ when $\sum_i|\gamma_i|<1/|\gamma_0|$.
\end{proposition}

\begin{proof}
Since $U$ is circulant and self-adjoint, we have $U=Fdiag(\beta)F^*$, for some vector $\beta\in\{\pm 1\}^N$. The first term in the expression of $\Phi(U,U)$ reads:
$$Tr[S^*U\cdot U^2]=\sum_{ij}|U_{ij}|=N\sum_i|\gamma_i|$$

For the second term in the formula of $\Phi$, we have the following formula:
$$S_{ii}={\rm sgn}(\gamma_0)={\rm sgn}\left(\sum_i\beta_i\right)\in\{\pm1\}$$

We therefore obtain the following formula:
$$\sum_{ij}\frac{Re[(U^2)_{ij}\bar S_{ij}]^2}{|U_{ij}|}=\sum_i\frac{1}{|\gamma_0|}=\frac{N}{|\gamma_0|}$$

But this leads to the conclusion in the statement.
\end{proof}

Still following \cite{bn3}, here is now a more advanced result, also in the circulant self-adjoint case, making this time use of a random derivative method:

\begin{theorem}
If $U\in U_N$ is circulant, $U_{ij}=\gamma_{j-i}$, and self-adjoint, we have
$$\mathbb E(\Phi(U,B))=N\sum_i|\gamma_i|-\frac{1}{2}\left(\frac{1}{|\gamma_0|}+\frac{1-e}{|\gamma_{N/2}|}+\sum_i\frac{1}{|\gamma_i|}\right)$$
where $e=0,1$ is the parity of $N$ and $\mathbb E$ denotes the expectation with respect to the uniform measure on the set of circulant self-adjoint unitary matrices $B$.
\end{theorem}

\begin{proof}
Since $B$ is circulant, this matrix is Fourier-diagonal. That is, we can diagonalize it with the help of the normalized Fourier matrix $F=F_N/\sqrt{N}$, as follows:
$$B=Fdiag(\alpha_i) F^*$$

The requirement that $B$ is unitary and self-adjoint amounts then to $\alpha_i=\pm 1$. The expectation is taken in the probability space where the random variables $\alpha_i$ are i.i.d., with symmetric Bernoulli distributions $(\delta_{-1}+\delta_1)/2$. In particular, we have:
$$\mathbb E[\alpha_i\alpha_j]=\delta_{ij}$$

By using $B^2=1_N$, the first term in the expression of $\Phi(U,B)$ reads:
$$Tr(S^*UB^2)
=Tr(S^*U)
=\sum_{ij}|U_{ij}|
=N\sum_i|\gamma_i|$$

For the second term in the formula of $\Phi$, observe first that we have:
$$Re[(UB)_{ij}\bar{S}_{ij}]^2=\frac{1}{4}\left[(UB)_{ij}^2\bar{S}_{ij}^2+ \overline{(UB)}_{ij}^2S_{ij}^2+2(UB)_{ij}\overline{(UB)}_{ij}\right]$$

We have the following computation, by using the formula $\mathbb E[\alpha_i\alpha_j]=\delta_{ij}$:
\begin{eqnarray*}
\mathbb E\left[(UB)_{ij}^2\right]
&=&\mathbb E\left[(F diag(q) diag(\alpha)F^*)^2_{ij}\right]\\
&=&\frac{1}{N^2}\sum_{kl}w^{(k+l)(i-j)}q_kq_l\mathbb E(\alpha_k\alpha_l)\\
&=&\frac{1}{N^2}\sum_{kl}w^{(k+l)(i-j)}q_kq_l\delta_{kl}\\
&=&\frac{1}{N^2}\sum_kw^{2k(i-j)} 
\end{eqnarray*}

We therefore obtain the following formula, for the above quantity:
$$\mathbb E\left[(UB)_{ij}^2\right]=\begin{cases}
\frac{1}{N}&\quad\ {\rm if}\ 2(i-j)=0\ ({\rm mod}\ N)\\
0  &\quad {\rm otherwise}
\end{cases}$$

Similarly, we have the following formula, for the last term:
\begin{eqnarray*}
\mathbb E\left[(UB)_{ij}\overline{(UB)}_{ij}\right]
&=&\frac{1}{N^2}\sum_{kl}w^{(k-l)(i-j)}q_k\bar q_l\mathbb E(\alpha_k\alpha_l)\\
&=&\frac{1}{N^2}\sum_k|q_k|^2\\
&=&\frac{1}{N}
\end{eqnarray*}

Since in both the cases $i=j$ and $i=j+N/2$, when $N$ is even, we have $S_{ij}\in\{\pm 1\}$, the above two formulae are all that we need, and we obtain the following formula:
$$\mathbb E\left[Re[(UB)_{ij}\bar S_{ij}]^2\right]=\frac{1}{4}\left[\frac{2\delta_{ij}}{N}+\frac{2(1-e)\delta_{i,j+N/2}}{N}+\frac{2}{N}\right]$$

Now by summing over $i,j$, and taking into account as well the first term in the expression of $\Phi(U,B)$, computed above, we obtain the formula in the statement. 
\end{proof}

In the orthogonal case now, we have a similar result, also from \cite{bn3}, as follows:

\index{circulant symmetric matrix}

\begin{theorem}
If $U\in O_N$ is circulant, $U_{ij}=\gamma_{j-i}$, and symmetric, we have
$$\mathbb E(\Phi(U,B))=N\sum_i |\gamma_i| - \left( \frac{1}{|\gamma_0|} +\frac{1-e}{|\gamma_{N/2}|} + \frac{N-2+e}{N}\sum_i \frac{1}{|\gamma_i|}\right)$$
where $e=0,1$ is the parity of $N$ and $\mathbb E$ denotes the expectation with respect to the uniform measure on the set of circulant symmetric orthogonal matrices $B$.
\end{theorem}

\begin{proof}
As before, in the proof of Theorem 12.25, the expectation is taken with respect to the distribution of the eigenvalues $\alpha_0,\ldots,\alpha_{N-1}=\pm 1$ of the matrix $B$, which are now, in the present real case, subject to the following extra condition:
$$\alpha_i=\alpha_{-i}$$

The first term in the expression of $\Phi(U,B)$ is then equal to $N \sum_i |\gamma_i|$. For the second term in $\Phi$, we need the following covariance term, in the present real case: 
$$\mathbb E(\alpha_k\alpha_l)=\begin{cases}
1 &\quad \text{ if } k \pm l = 0\\
0 &\quad \text{ otherwise}
\end{cases}$$

Since all quantities are real in this case, we have the following formula:
\begin{eqnarray*}
\mathbb E\left[(UB)_{ij}^2\right]
&=&\frac{1}{N^2}\sum_{kl}w^{(k+l)(i-j)}q_kq_l \mathbb E(\alpha_k\alpha_l)\\
&=&\frac{1}{N^2}\sum_{kl}w^{(k+l)(i-j)}q_kq_l(\delta_{kl}+\delta_{k,-l}-\delta_{2k,2l,0})\\
&=&\frac{1}{N^2}\left[\sum_kw^{2k(i-j)}q_k^2+\sum_kq_kq_{-k}-q_0^2-(1-e)q_{N/2}^2\right]  \\
&=&\frac{1}{N^2}\left[N\delta_{2i,2j}+N-2+e\right] 
\end{eqnarray*}

We have then the following formula:
$$\sum_{ij} N^{-1} |U_{ij}|^{-1} \delta_{2i,2j} = \sum_k |\gamma_{k}|^{-1} \delta_{2k,0} = \frac{1}{|\gamma_0|} + \frac{1-e}{|\gamma_{N/2}|}$$

On the other hand, we have as well the following formula:
$$\sum_{ij} N^{-2}(N-2+e) |U_{ij}|^{-1}  = \frac{N-2+e}{N} \sum_{i} \frac{1}{|\gamma_i|}$$

Now by putting everything together, gives the formula in the statement. 
\end{proof}

As an illustration for the above methods, we can now go back to the matrices in Theorem 12.20, and find a better proof for the fact that these matrices are not complex AHM. Indeed, we have the following result, which basically solves the problem:

\begin{proposition}
With $U=\frac{1}{N}(2\mathbb I_N-N1_N)$ we have the formula
$$\mathbb E(\Phi(U,B))=\frac{4-N}{2}\left(N-4-\frac{2+e}{N-2}\right)$$
where $e=0,1$ is the parity of $N$, and where $B$ varies over the space of orthogonal circulant symmetric matrices. This quantity equals $-2,0,0,-\frac{3}{2},-\frac{18}{5},\ldots$ at $N=3,4,5,6,7\ldots$
\end{proposition}

\begin{proof}
This follows indeed from the general formula in Theorem 12.26.
\end{proof}

We can therefore recover Theorem 12.20, modulo a bit of extra work still needed at $N=5$. Regarding the case $N=5$, here the above expectation vanishes, but by using Proposition 12.23 or Proposition 12.24, we conclude that the vanishing of the expectation must come from both positive and negative contributions, and we are done. 

\bigskip

In fact, the above results can be used for excluding all the explicit examples of circulant AHM found in \cite{bnz}. All these verifications suggest the following conjecture:

\index{random derivative}

\begin{conjecture}
For any $U\in O_N$ which is circulant and symmetric we have
$$\mathbb E(\Phi(U,B))\leq0$$
where $B$ varies over the space of orthogonal circulant symmetric matrices. In addition, a similar result should hold in the unitary, circulant and self-adjoint case.
\end{conjecture}

This looks like a subtle Fourier analysis question. In fact, the main idea that emerges from the computations in \cite{bn3}, including the block design ones, is that of using a random derivative, pointing towards a suitable homogeneous space coset. However, no one really knows how to do that. And so we will have it as an exercise for you, reader.

\section*{12e. Exercises} 

The material in the present chapter has been quite research-oriented, and our exercises here will be of the same type, rather difficult. First, we have:

\begin{exercise}
Establish the rotation trick, stating that we must have
$$U_{ij}\neq0$$
for the local maxima/minima of the $p$-norms on $U_N$, at values $p\neq1$.
\end{exercise}

The cases $p<2$ and $p>2$ are of quite different nature, at least when using a straightforward approach to the problem, in the spirit of the one that we used in the above, at $p=1$. The first problem is that of deciding which case is the one to go with.

\begin{exercise}
Establish the Hessian formula for the second derivative of the $1$-norm by using advanced differential geometry techniques.
\end{exercise}

To be more precise here, the formula for the second derivative that we obtained in the above was based on some straightforward computations, which are quite long. The problem is that of replacing these computations by something more conceptual, based on advanced knowledge of differential geometry, or of calculus in several variables.

\begin{exercise}
Verify the AHC for the various examples of almost Hadamard matrices, in the real sense, from chapter 3, coming from block designs.
\end{exercise}

There are many things that can be done here, and as a bottom line, your computations should generalize those that we have for $K_N$, explained in the above.

\begin{exercise}
Reformulate the verifications of the AHC for circulant matrices presented in the above in a more conceptual way, by using a random derivative method, pointing towards a suitable homogeneous space coset.
\end{exercise}

To be more precise here, the homogeneous space coset in question should appear by applying a discrete Fourier transform to the circulant matrices.

\part{Quantum algebra}

\ \vskip50mm

\begin{center}
{\em Many things about tomorrow

I don't seem to understand

But I know who holds tomorrow

And I know who holds my hand}
\end{center}

\chapter{Quantum groups}

\section*{13a. Operator algebras}

Welcome to this fourth and last part of the present book. We discuss here yet another idea in order to deal with the Hadamard matrices, be them real or complex, this time in relation with quantum groups. What we will be doing here will be deeply related to all sorts of advanced algebraic considerations regarding the Hadamard matrices, from chapters 1-12 above, and also to a quite good deal of deep considerations from operator algebras, following Haagerup \cite{ha1}, Jones \cite{jo3}, Popa \cite{pop} and others. So, we will be here working at a foundational level in mathematical physics. In fact, all the potential applications of the complex Hadamard matrices to questions in physics, be them from general quantum mechanics, quantum information, statistical mechanics, and many more, are expected to come via the link with the quantum groups.

\bigskip

The idea is extremely simple, namely that associated to any complex Hadamard matrix $H\in M_N(\mathbb C)$ is a certain quantum permutation group $G\subset S_N^+$, which describes the ``symmetries'' of the matrix. As a basic illustration, for a Fourier matrix $H=F_G$ we obtain the group $G$ itself, acting on itself, $G\subset S_G$. In general, however, we obtain non-classical quantum groups, whose computation is a key problem.

\bigskip

In order to discuss this, we will need many preliminaries, namely operator theory, operator algebras and quantum spaces, compact quantum groups, quantum permutation groups, and finally matrix models for such quantum groups, which produce the above correspondence. Before getting started, some references. For functional analysis, operator theory and operator algebras you have Lax \cite{lax}, and also Connes \cite{con}, if you want to learn more. For quantum groups you have the papers of Woronowicz \cite{wo1}, \cite{wo2} or my book \cite{ba1}, but we will explain the needed material here. For tools for dealing with such quantum groups, these will often come from Jones \cite{jo1}, \cite{jo2}, \cite{jo3} and Voiculescu \cite{vdn}.

\bigskip

Also, importantly, there are no quantum groups or quantum mechanics without quantum mechanics. In order to appreciate what will follow, get to learn some, standard places being Feynman \cite{fey}, Griffiths \cite{gri}, Weinberg \cite{wei}. In case you would rather enjoy a rigorous text written by a mathematician, you can go with my book \cite{ba2}. Although that is not an inch more clever, or even rigorous, than what physicists are doing.

\bigskip

Getting started now, we first have the following standard result:

\index{linear operator}
\index{bounded operator}
\index{operator norm}
\index{Banach algebra}
\index{adjoint operator}

\begin{theorem}
Given a complex Hilbert space $H$, the linear operators $T:H\to H$ which are bounded, in the sense that the quantity
$$||T||=\sup_{||x||\leq1}||Tx||$$
is finite, form a complex algebra with unit, denoted $B(H)$, having the following properties:
\begin{enumerate}
\item $B(H)$ is complete with respect to $||.||$, so we have a Banach algebra. 

\item $B(H)$ has an involution $T\to T^*$, given by $<Tx,y>=<x,T^*y>$.
\end{enumerate}
In addition, the norm and involution are related by the formula $||TT^*||=||T||^2$.
\end{theorem}

\begin{proof}
The fact that we have indeed an algebra follows from:
$$||S+T||\leq||S||+||T||\quad,\quad 
||\lambda T||=|\lambda|\cdot||T||\quad,\quad 
||ST||\leq||S||\cdot||T||$$

Regarding now (1), if $\{T_n\}\subset B(H)$ is Cauchy then $\{T_nx\}$ is Cauchy for any $x\in H$, so we can define the limit $T=\lim_{n\to\infty}T_n$ simply by setting:
$$Tx=\lim_{n\to\infty}T_nx$$

As for (2), here the existence of $T^*$ comes from the fact that $\varphi(x)=<Tx,y>$ being a linear map $H\to\mathbb C$, we must have, for a certain vector $T^*y\in H$:
$$\varphi(x)=<x,T^*y>$$

Moreover, since this vector is unique, $T^*$ is unique too, and we have as well:
$$(S+T)^*=S^*+T^*\quad,\quad 
(\lambda T)^*=\bar{\lambda}T^*\quad,\quad 
(ST)^*=T^*S^*\quad,\quad
(T^*)^*=T$$

Observe also that we have indeed $T^*\in B(H)$, because:
\begin{eqnarray*}
||T||
&=&\sup_{||x||=1}\sup_{||y||=1}<Tx,y>\\
&=&\sup_{||y||=1}\sup_{||x||=1}<x,T^*y>\\
&=&||T^*||
\end{eqnarray*}

Regarding now the last assertion, observe first that we have:
$$||TT^*||
\leq||T||\cdot||T^*||
=||T||^2$$

On the other hand, we have as well the following estimate:
\begin{eqnarray*}
||T||^2
&=&\sup_{||x||=1}|<Tx,Tx>|\\
&=&\sup_{||x||=1}|<x,T^*Tx>|\\
&\leq&||T^*T||
\end{eqnarray*}

By replacing $T\to T^*$ we obtain from this $||T||^2\leq||TT^*||$, and we are done.
\end{proof}

We will be interested in the algebras of operators, rather than in the operators themselves. The basic axioms here, inspired from Theorem 13.1, are as follows:

\index{operator algebra}

\begin{definition}
A $C^*$-algebra is a complex algebra with unit $A$, having:
\begin{enumerate}
\item A norm $a\to||a||$, making it a Banach algebra (the Cauchy sequences converge).

\item An involution $a\to a^*$, which satisfies $||aa^*||=||a||^2$, for any $a\in A$.
\end{enumerate}
\end{definition}

According to Theorem 13.1, the operator algebra $B(H)$ itself is a $C^*$-algebra. More generally, we have as examples all the closed $*$-subalgebras $A\subset B(H)$. We will see later on (the ``GNS theorem'') that any $C^*$-algebra appears in fact in this way. However, even before knowing that, in view of the examples that we have, we can think of the elements $a\in A$ of an arbitrary $C^*$-algebra as being some kind of ``generalized beounded operators'', on some Hilbert space which is not necessarily present. By using this idea, one can emulate spectral theory in this setting, and we have the following result:

\index{spectrum}
\index{spectral radius}
\index{normal element}

\begin{theorem}
Given $a\in A$, define its spectrum as being the set
$$\sigma(a)=\left\{\lambda\in\mathbb C\Big|a-\lambda\not\in A^{-1}\right\}$$
and its spectral radius $\rho(a)$ as the radius of the smallest centered disk containing $\sigma(a)$.
\begin{enumerate}
\item The spectrum of a norm one element is in the unit disk.

\item The spectrum of a unitary element $(a^*=a^{-1}$) is on the unit circle. 

\item The spectrum of a self-adjoint element ($a=a^*$) consists of real numbers. 

\item The spectral radius of a normal element ($aa^*=a^*a$) is equal to its norm.
\end{enumerate}
\end{theorem}

\begin{proof}
Our first claim is that for any polynomial $f\in\mathbb C[X]$, and more generally for any rational function $f\in\mathbb C(X)$ having poles outside $\sigma(a)$, we have:
$$\sigma(f(a))=f(\sigma(a))$$

This indeed something well-known for the usual matrices. In the general case, assume first that we have a polynomial, $f\in\mathbb C[X]$. If we pick an arbitrary number $\lambda\in\mathbb C$, and write $f(X)-\lambda=c(X-r_1)\ldots(X-r_k)$, we have then, as desired:
\begin{eqnarray*}
\lambda\notin\sigma(f(a))
&\iff&f(a)-\lambda\in A^{-1}\\
&\iff&c(a-r_1)\ldots(a-r_k)\in A^{-1}\\
&\iff&a-r_1,\ldots,a-r_k\in A^{-1}\\
&\iff&r_1,\ldots,r_k\notin\sigma(a)\\
&\iff&\lambda\notin f(\sigma(a))
\end{eqnarray*}

Assume now that we are in the general case, $f\in\mathbb C(X)$. We pick $\lambda\in\mathbb C$, we write $f=P/Q$, and we consider the following polynomial:
$$F=P-\lambda Q$$

By using the above finding, for this polynomial $F$, we obtain, as desired:
\begin{eqnarray*}
\lambda\in\sigma(f(a))
&\iff&F(a)\notin A^{-1}\\
&\iff&0\in\sigma(F(a))\\
&\iff&0\in F(\sigma(a))\\
&\iff&\exists\mu\in\sigma(a),F(\mu)=0\\
&\iff&\lambda\in f(\sigma(a))
\end{eqnarray*}

Regarding now the assertions in the statement, these basically follow from this:

\medskip

(1) This comes from the following formula, valid when $||a||<1$:
$$\frac{1}{1-a}=1+a+a^2+\ldots$$

(2) Assuming $a^*=a^{-1}$, we have the following norm computations:
$$||a||=\sqrt{||aa^*||}=\sqrt{1}=1$$
$$||a^{-1}||=||a^*||=||a||=1$$

If we denote by $D$ the unit disk, we obtain from this, by using (1):
$$||a||=1\implies\sigma(a)\subset D$$
$$||a^{-1}||=1\implies\sigma(a^{-1})\subset D$$

On the other hand, by using the rational function $f(z)=z^{-1}$, we have:
$$\sigma(a^{-1})\subset D\implies \sigma(a)\subset D^{-1}$$

Now by putting everything together we obtain, as desired:
$$\sigma(a)\subset D\cap D^{-1}=\mathbb T$$

(3) This follows by using (2), and the following rational function, with $t\in\mathbb R$:
$$f(z)=\frac{z+it}{z-it}$$

Indeed, for $t>>0$ the element $f(a)$ is well-defined, and we have:
$$\left(\frac{a+it}{a-it}\right)^*
=\frac{a-it}{a+it}
=\left(\frac{a+it}{a-it}\right)^{-1}$$

Thus $f(a)$ is a unitary, and by (2) its spectrum is contained in $\mathbb T$. We conclude that we have $f(\sigma(a))=\sigma(f(a))\subset\mathbb T$, and so $\sigma(a)\subset f^{-1}(\mathbb T)=\mathbb R$, as desired.

\medskip

(4) We have $\rho(a)\leq ||a||$ from (1). Conversely, given $\rho>\rho(a)$, we have:
$$\int_{|z|=\rho}\frac{z^n}{z-a}\,dz 
=\sum_{k=0}^\infty\left(\int_{|z|=\rho}z^{n-k-1}dz\right) a^k
=a^{n-1}$$

By applying the norm and taking $n$-th roots we obtain:
$$\rho\geq\lim_{n\to\infty}||a^n||^{1/n}$$

In the case $a=a^*$ we have $||{a^n}||=||{a}||^n$ for any exponent of the form $n=2^k$, and by taking $n$-th roots we get $\rho\geq ||{a}||$. This gives the missing inequality, namely:
$$\rho(a)\geq ||a||$$

In the general case $aa^*=a^*a$ we have $a^n(a^n)^*=(aa^*)^n$, and we get:
$$\rho(a)^2=\rho(aa^*)$$ 

Now since $aa^*$ is self-adjoint, we get $\rho(aa^*)=||{a}||^2$, and we are done.
\end{proof}

With these preliminaries in hand, we can now formulate some theorems. The basic facts about the $C^*$-algebras, that we will need here, can be summarized as:

\index{Gelfand theorem}
\index{GNS theorem}
\index{multimatrix algebra}

\begin{theorem}
The $C^*$-algebras have the following properties:
\begin{enumerate}
\item The commutative ones are those of the form $C(X)$, with $X$ compact space.

\item Any such algebra $A$ embeds as $A\subset B(H)$, for some Hilbert space $H$.

\item In finite dimensions, these are the direct sums of matrix algebras.
\end{enumerate}
\end{theorem}

\begin{proof}
All this is standard, the idea being as follows:

\medskip

(1) Given a compact space $X$, the algebra $C(X)$ of continuous functions $f:X\to\mathbb C$ is indeed a $C^*$-algebra, with norm and involution as follows:
$$||f||=\sup_{x\in X}|f(x)|\quad,\quad
f^*(x)=\overline{f(x)}$$

Observe that this algebra is indeed commutative, because:
$$f(x)g(x)=g(x)f(x)$$

Conversely, if $A$ is commutative, we can define $X=Spec(A)$ to be the space of all characters $\chi :A\to\mathbb C$, with the topology making continuous all the evaluation maps $ev_a:\chi\to\chi(a)$. We have then a morphism of algebras, as follows:
$$ev:A\to C(X)\quad,\quad 
a\to ev_a$$

Theorem 13.3 (3) shows that $ev$ is a $*$-morphism, Theorem 13.3 (4) shows that $ev$ is isometric, and finally the Stone-Weierstrass theorem shows that $ev$ is surjective.

\medskip

(2) This is standard for $A=C(X)$, where we can pick a probability measure on $X$, and set $H=L^2(X)$, and use the following embedding:
$$A\subset B(H)\quad,\quad 
f\to(g\to fg)$$

In the general case, where $A$ is no longer commutative, the proof is quite similar, by emulating basic measure theory in the abstract $C^*$-algebra setting.

\medskip

(3) Assuming that $A$ is finite dimensional, we can first decompose its unit as  follows, with $p_i\in A$ being central minimal projections:
$$1=p_1+\ldots+p_k$$

Each of the linear spaces $A_i=p_iAp_i$ is then a non-unital $*$-subalgebra of $A$, and we have a non-unital $*$-algebra sum decomposition, as follows:
$$A=A_1\oplus\ldots\oplus A_k$$

On the other hand, since each central projection $p_i$ was assumed minimal, we have unital $*$-algebra isomorphisms as follows, with $r_i=rank(p_i)$:
$$A_i\simeq M_{r_i}(\mathbb C)$$

Thus, we obtain an isomorphism $A\simeq M_{r_1}(\mathbb C)\oplus\ldots\oplus M_{r_k}(\mathbb C)$, as desired.
\end{proof}

All the above was of course quite brief, but full details on this can be found in any book on functional analysis, as for instance Lax \cite{lax}. In what concerns us, we will be mainly interested in Theorem 13.4 (1), called Gelfand theorem, which suggests formulating:

\index{compact quantum space}
\index{quantum space}

\begin{definition}
Given a $C^*$-algebra $A$, not necessarily commutative, we write
$$A=C(X)$$
and call the abstract object $X$ a compact quantum space.
\end{definition}

In other words, we define the category of the compact quantum spaces $X$ to be the category of the $C^*$-algebras $A$, with the arrows reversed. Due to the Gelfand theorem, 13.4 (1) above, the category of the usual compact spaces embeds covariantly into the category of the compact quantum spaces, and the image of this embedding consists precisely of the compact quantum spaces $X$ which are ``classical'', in the sense that the corresponding $C^*$-algebra $A=C(X)$ is commutative. Thus, what we have done here is to extend the category of the usual compact spaces, and this justifies Definition 13.5.

\bigskip

In practice now, the general compact quantum spaces $X$ do not have points, but we can perfectly study them via the associated algebras $A=C(X)$, a bit in the same way as we study a compact Lie group via its associated Lie algebra, or an algebraic manifold via the ideal of polynomials vanishing on it, and so on. In short, nothing that much abstract going on here, just another instance of the old idea ``we will use algebras, no need for points'', with the remark that for us, the use of points will be actually forbidden.

\section*{13b. Quantum groups}

We will be interested in what follows in the case where the compact quantum space $X$ is a ``compact quantum group''. The axioms for the corresponding $C^*$-algebras, found by Woronowicz in \cite{wo1}, are, in a soft form, as follows:

\index{Woronowicz algebra}

\begin{definition}
A Woronowicz algebra is a $C^*$-algebra $A$, given with a unitary matrix $u\in M_N(A)$ whose coefficients generate $A$, such that the formulae
$$\Delta(u_{ij})=\sum_ku_{ik}\otimes u_{kj}$$
$$\varepsilon(u_{ij})=\delta_{ij}$$
$$S(u_{ij})=u_{ji}^*$$
define morphisms of $C^*$-algebras $\Delta:A\to A\otimes A$, $\varepsilon:A\to\mathbb C$, $S:A\to A^{opp}$.
\end{definition}

The morphisms $\Delta,\varepsilon,S$ are called comultiplication, counit and antipode. We say that $A$ is cocommutative when $\Sigma\Delta=\Delta$, where $\Sigma(a\otimes b)=b\otimes a$ is the flip. We have the following result, which justifies the terminology and axioms:

\index{cocommutative algebra}

\begin{proposition}
The following are Woronowicz algebras:
\begin{enumerate}
\item $C(G)$, with $G\subset U_N$ compact Lie group. Here the structural maps are:
\begin{eqnarray*}
\Delta(\varphi)&=&(g,h)\to \varphi(gh)\\
\varepsilon(\varphi)&=&\varphi(1)\\
S(\varphi)&=&g\to\varphi(g^{-1})
\end{eqnarray*}

\item $C^*(\Gamma)$, with $F_N\to\Gamma$ finitely generated group. Here the structural maps are:
\begin{eqnarray*}
\Delta(g)&=&g\otimes g\\
\varepsilon(g)&=&1\\ 
S(g)&=&g^{-1}
\end{eqnarray*}
\end{enumerate}
Moreover, we obtain in this way all the commutative/cocommutative algebras.
\end{proposition}

\begin{proof}
This is something very standard, the idea being as follows:

\medskip

(1) Given $G\subset U_N$, we can set $A=C(G)$, which is a Woronowicz algebra, together with the matrix $u=(u_{ij})$ formed by coordinates of $G$, given by:
$$g=\begin{pmatrix}
u_{11}(g)&\ldots&u_{1N}(g)\\
\vdots&&\vdots\\
u_{N1}(g)&\ldots&u_{NN}(g)
\end{pmatrix}$$

Conversely, if $(A,u)$ is a commutative Woronowicz algebra, by using the Gelfand theorem we can write $A=C(X)$, with $X$ being a certain compact space. The coordinates $u_{ij}$ give then an embedding $X\subset M_N(\mathbb C)$, and since the matrix $u=(u_{ij})$ is unitary we actually obtain an embedding $X\subset U_N$, and finally by using the maps $\Delta,\varepsilon,S$ we conclude that our compact subspace $X\subset U_N$ is in fact a compact Lie group, as desired.

\medskip

(2) Consider a finitely generated group $F_N\to\Gamma$. We can set $A=C^*(\Gamma)$, which is by definition the completion of the complex group algebra $\mathbb C[\Gamma]$, with involution given by $g^*=g^{-1}$, for any $g\in\Gamma$, with respect to the biggest $C^*$-norm, and we obtain a Woronowicz algebra, together with the diagonal matrix formed by the generators of $\Gamma$:
$$u=\begin{pmatrix}
g_1&&0\\
&\ddots&\\
0&&g_N
\end{pmatrix}$$

Conversely, if $(A,u)$ is a cocommutative Woronowicz algebra, the Peter-Weyl theory of Woronowicz, to be explained below, shows that the irreducible corepresentations of $A$ are all 1-dimensional, and form a group $\Gamma$, and so we have $A=C^*(\Gamma)$, as desired.
\end{proof}

In relation with the above, we should mention that there are actually some analytic subtleties here, coming from amenability, and so our quantum spaces and groups must be divided by a certain equivalence relation, for everything to work fine. To be more precise, in the context of Definition 13.6, we write $(A,u)=(B,v)$ when there is a $*$-algebra isomorphism as follows, mapping standard coordinates to standard coordinates:
$$<u_{ij}>\simeq<v_{ij}>\quad,\quad u_{ij}\to v_{ij}$$

In general now, the structural maps $\Delta,\varepsilon,S$ have the following properties:

\index{comultiplication}
\index{counit}
\index{antipode}
\index{square of antipode}

\begin{proposition}
Let $(A,u)$ be a Woronowicz algebra.
\begin{enumerate} 
\item $\Delta,\varepsilon$ satisfy the usual axioms for a comultiplication and a counit, namely:
\begin{eqnarray*}
(\Delta\otimes id)\Delta&=&(id\otimes \Delta)\Delta\\
(\varepsilon\otimes id)\Delta&=&(id\otimes\varepsilon)\Delta=id
\end{eqnarray*}

\item $S$ satisfies the antipode axiom, on the $*$-subalgebra generated by entries of $u$: 
$$m(S\otimes id)\Delta=m(id\otimes S)\Delta=\varepsilon(.)1$$

\item In addition, the square of the antipode is the identity, $S^2=id$.
\end{enumerate}
\end{proposition}

\begin{proof}
The two comultiplication axioms follow from:
\begin{eqnarray*}
(\Delta\otimes id)\Delta(u_{ij})&=&(id\otimes \Delta)\Delta(u_{ij})=\sum_{kl}u_{ik}\otimes u_{kl}\otimes u_{lj}\\
(\varepsilon\otimes id)\Delta(u_{ij})&=&(id\otimes\varepsilon)\Delta(u_{ij})=u_{ij}
\end{eqnarray*}

As for the antipode formulae, the verification here is similar.
\end{proof}

Summarizing, the Woronowicz algebras appear to have nice properties. In view of Proposition 13.7 and Proposition 13.8, we can formulate the following definition:

\index{compact quantum group}
\index{discrete quantum group}

\begin{definition}
Given a Woronowicz algebra $A$, we formally write
$$A=C(G)=C^*(\Gamma)$$
and call $G$ compact quantum group, and $\Gamma$ discrete quantum group.
\end{definition}

When $A$ is both commutative and cocommutative, $G$ is a compact abelian group, $\Gamma$ is a discrete abelian group, and these groups are dual to each other, $G=\widehat{\Gamma},\Gamma=\widehat{G}$. In general, we still agree to write, but in a formal sense:
$$G=\widehat{\Gamma}\quad,\quad\Gamma=\widehat{G}$$

\index{corepresentation}

With this in mind, let us call now corepresentation of $A$ any unitary matrix $v\in M_n(A)$ satisfying the same conditions as those satisfied by $u$, namely:
$$\Delta(v_{ij})=\sum_kv_{ik}\otimes v_{kj}\quad,\quad
\varepsilon(v_{ij})=\delta_{ij}\quad,\quad
S(v_{ij})=v_{ji}^*$$

These corepresentations can be thought of as corresponding representations of the underlying compact quantum group $G$. Following Woronowicz \cite{wo1}, we have:

\index{Haar functional}
\index{Ces\`aro limit}

\begin{theorem}
Any Woronowicz algebra has a unique Haar integration functional, 
$$\left(\int_G\otimes id\right)\Delta=\left(id\otimes\int_G\right)\Delta=\int_G(.)1$$
which can be constructed by starting with any faithful positive form $\varphi\in A^*$, and setting
$$\int_G=\lim_{n\to\infty}\frac{1}{n}\sum_{k=1}^n\varphi^{*k}$$
where $\phi*\psi=(\phi\otimes\psi)\Delta$. Moreover, for any corepresentation $v\in M_n(\mathbb C)\otimes A$ we have
$$\left(id\otimes\int_G\right)v=P$$
where $P$ is the orthogonal projection onto $Fix(v)=\{\xi\in\mathbb C^n|v\xi=\xi\}$.
\end{theorem}

\begin{proof}
Following \cite{wo1}, this can be done in 3 steps, as follows:

\medskip

(1) Given $\varphi\in A^*$, our claim is that the following limit converges, for any $a\in A$:
$$\int_\varphi a=\lim_{n\to\infty}\frac{1}{n}\sum_{k=1}^n\varphi^{*k}(a)$$

Indeed, by linearity we can assume that $a$ is the coefficient of corepresentation, $a=(\tau\otimes id)v$. But in this case, an elementary computation shows that we have the following formula, where $P_\varphi$ is the orthogonal projection onto the $1$-eigenspace of $(id\otimes\varphi)v$:
$$\left(id\otimes\int_\varphi\right)v=P_\varphi$$

(2) Since $v\xi=\xi$ implies $[(id\otimes\varphi)v]\xi=\xi$, we have $P_\varphi\geq P$, where $P$ is the orthogonal projection onto the space $Fix(v)=\{\xi\in\mathbb C^n|v\xi=\xi\}$. The point now is that when $\varphi\in A^*$ is faithful, by using a positivity trick, one can prove that we have $P_\varphi=P$. Thus our linear form $\int_\varphi$ is independent of $\varphi$, and is given on coefficients $a=(\tau\otimes id)v$ by:
$$\left(id\otimes\int_\varphi\right)v=P$$

(3) With the above formula in hand, the left and right invariance of $\int_G=\int_\varphi$ is clear on coefficients, and so in general, and this gives all the assertions. See \cite{wo1}. 
\end{proof}

Consider the dense $*$-subalgebra $\mathcal A\subset A$ generated by the coefficients of the fundamental corepresentation $u$, and endow it with the following scalar product: 
$$<a,b>=\int_Gab^*$$

We have then the following result, also due to Woronowicz \cite{wo1}:

\index{Peter-Weyl representation}
\index{Peter-Weyl theory}

\begin{theorem}
We have the following Peter-Weyl type results:
\begin{enumerate}
\item Any corepresentation decomposes as a sum of irreducible corepresentations.

\item Each irreducible corepresentation appears inside a certain $u^{\otimes k}$.

\item $\mathcal A=\bigoplus_{v\in Irr(A)}M_{\dim(v)}(\mathbb C)$, the summands being pairwise orthogonal.

\item The characters of irreducible corepresentations form an orthonormal system.
\end{enumerate}
\end{theorem}

\begin{proof}
All these results are from \cite{wo1}, the idea being as follows:

\medskip

(1) Given $v\in M_n(A)$, its intertwiner algebra $End(v)=\{T\in M_n(\mathbb C)|Tv=vT\}$ is a finite dimensional $C^*$-algebra, and so decomposes as $End(v)=M_{n_1}(\mathbb C)\oplus\ldots\oplus M_{n_r}(\mathbb C)$. But this gives a decomposition of type $v=v_1+\ldots+v_r$, as desired.

\medskip

(2) Consider indeed the Peter-Weyl corepresentations, $u^{\otimes k}$ with $k$ colored integer, defined by $u^{\otimes\emptyset}=1$, $u^{\otimes\circ}=u$, $u^{\otimes\bullet}=\bar{u}$ and multiplicativity. The coefficients of these corepresentations span the dense algebra $\mathcal A$, and by using (1), this gives the result.

\medskip

(3) Here the direct sum decomposition, which is technically a $*$-coalgebra isomorphism, follows from (2). As for the second assertion, this follows from the fact that $(id\otimes\int_G)v$ is the orthogonal projection $P_v$ onto the space $Fix(v)$, for any corepresentation $v$.

\medskip

(4) Let us define indeed the character of $v\in M_n(A)$ to be the matrix trace, $\chi_v=Tr(v)$. Since this character is a coefficient of $v$, the orthogonality assertion follows from (3). As for the norm 1 claim, this follows once again from $(id\otimes\int_G)v=P_v$. 
\end{proof}

Observe that in the cocommutative case, we obtain from (4) that the irreducible corepresentations must be all 1-dimensional, and so that we must have $A=C^*(\Gamma)$ for some discrete group $\Gamma$, as mentioned in Proposition 13.7. 

\index{cocommutative algebra}

\section*{13c. Quantum permutations}

We will be interested here in the quantum permutation groups, and their relation with the Hadamard matrices. The following key definition is due to Wang \cite{wa1}:

\index{magic matrix}
\index{magic unitary}

\begin{definition}
A magic unitary matrix is a square matrix over a $C^*$-algebra, 
$$u\in M_N(A)$$
whose entries are projections, summing up to $1$ on each row and each column.
\end{definition}

The basic examples of such matrices come from the usual permutation groups, $G\subset S_N$. Indeed, given such subgroup, the following matrix is magic:
$$u_{ij}=\chi\left(\sigma\in G\Big|\sigma(j)=i\right)$$

The interest in these matrices comes from the following functional analytic description of the usual symmetric group, from \cite{wa1}:

\begin{proposition}
Consider the symmetric group $S_N$.
\begin{enumerate}
\item The standard coordinates $v_{ij}\in C(S_N)$, coming from the embedding $S_N\subset O_N$ given by the permutation matrices, are given by $v_{ij}=\chi(\sigma|\sigma(j)=i)$.

\item The matrix $v=(v_{ij})$ is magic, in the sense that its entries are orthogonal projections, summing up to $1$ on each row and each column.

\item The algebra $C(S_N)$ is isomorphic to the universal commutative $C^*$-algebra generated by the entries of a $N\times N$ magic matrix.
\end{enumerate}
\end{proposition}

\begin{proof}
These results are all elementary, as follows:

\medskip

(1) The canonical embedding $S_N\subset O_N$, coming from the standard permutation matrices, is given by $\sigma(e_j)=e_{\sigma(j)}$. Thus, we have $\sigma=\sum_je_{\sigma(j)j}$, so the standard coordinates on $S_N\subset O_N$ are given by $v_{ij}(\sigma)=\delta_{i,\sigma(j)}$. Thus, we must have, as claimed:
$$v_{ij}=\chi\left(\sigma\Big|\sigma(j)=i\right)$$

(2) Any characteristic function $\chi\in\{0,1\}$ being a projection in the operator algebra sense ($\chi^2=\chi^*=\chi$), we have indeed a matrix of projections. As for the sum 1 condition on rows and columns, this is clear from the formula of the elements $v_{ij}$.

\medskip

(3) Consider the universal algebra in the statement, namely:
$$A=C^*_{comm}\left((w_{ij})_{i,j=1,\ldots,N}\Big|w={\rm magic}\right)$$

We have a quotient map $A\to C(S_N)$, given by $w_{ij}\to v_{ij}$. On the other hand, by using the Gelfand theorem we can write $A=C(X)$, with $X$ being a compact space, and by using the coordinates $w_{ij}$ we have $X\subset O_N$, and then $X\subset S_N$. Thus we have as well a quotient map $C(S_N)\to A$ given by $v_{ij}\to w_{ij}$, and this gives (3). See Wang \cite{wa1}.
\end{proof}

We are led in this way to the following result:

\index{Wang theorem}
\index{quantum permutation group}
\index{Wang algebra}

\begin{theorem}
The following is a Woronowicz algebra,
$$C(S_N^+)=C^*\left((u_{ij})_{i,j=1,\ldots,N}\Big|u={\rm magic}\right)$$
and the underlying compact quantum group $S_N^+$ is called quantum permutation group.
\end{theorem}

\begin{proof}
As a first remark, the algebra $C(S_N^+)$ is indeed well-defined, because the magic condition forces $||u_{ij}||\leq1$, for any $C^*$-norm. Our claim now is that we can define maps $\Delta,\varepsilon,S$ as in Definition 13.6. Consider indeed the following matrix: 
$$U_{ij}=\sum_ku_{ik}\otimes u_{kj}$$

As a first observation, we have $U_{ij}=U_{ij}^*$. In fact the entries $U_{ij}$ are orthogonal projections, because we have as well:
$$U_{ij}^2
=\sum_{kl}u_{ik}u_{il}\otimes u_{kj}u_{lj}
=\sum_ku_{ik}\otimes u_{kj}
=U_{ij}$$

In order to prove now that the matrix $U=(U_{ij})$ is magic, it remains to verify that the sums on the rows and columns are 1. For the rows, this can be checked as follows:
$$\sum_jU_{ij}
=\sum_{jk}u_{ik}\otimes u_{kj}
=\sum_ku_{ik}\otimes1
=1\otimes1$$

For the columns the computation is similar, as follows:
$$\sum_iU_{ij}
=\sum_{ik}u_{ik}\otimes u_{kj}
=\sum_k1\otimes u_{kj}
=1\otimes1$$

Thus the matrix $U=(U_{ij})$ is magic indeed, as claimed above, and so we can define a comultiplication map, simply by setting:
$$\Delta(u_{ij})=U_{ij}$$

By using a similar reasoning, and similar elementary computations, we can define as well a counit map by $\varepsilon(u_{ij})=\delta_{ij}$, and an antipode by $S(u_{ij})=u_{ji}$. Thus the Woronowicz algebra axioms from Definition 13.6 are satisfied, and this finishes the proof.
\end{proof}

The terminology comes from the following result, also from Wang \cite{wa1}:

\index{coaction}
\index{counting measure}

\begin{proposition}
The quantum group $S_N^+$ acts on the set $X=\{1,\ldots,N\}$, the corresponding coaction map $\Phi:C(X)\to C(X)\otimes C(S_N^+)$ being given by:
$$\Phi(\delta_i)=\sum_j\delta_j\otimes u_{ji}$$
In fact, $S_N^+$ is the biggest compact quantum group acting on $X$, by leaving the counting measure invariant, in the sense that $(tr\otimes id)\Phi=tr(.)1$, where $tr(\delta_i)=\frac{1}{N},\forall i$.
\end{proposition}

\begin{proof}
Our claim is that given a compact quantum group $G$, the formula $\Phi(\delta_i)=\sum_j\delta_j\otimes u_{ji}$ defines a morphism of algebras, which is a coaction map, leaving the trace invariant, precisely when the matrix $u=(u_{ij})$ is a magic corepresentation of $C(G)$. Indeed, let us first determine when $\Phi$ is multiplicative. We have:
$$\Phi(\delta_i)\Phi(\delta_k)
=\sum_{jl}\delta_j\delta_l\otimes u_{ji}u_{lk}
=\sum_j\delta_j\otimes u_{ji}u_{jk}$$

On the other hand, we have as well:
$$\Phi(\delta_i\delta_k)
=\delta_{ik}\Phi(\delta_i)
=\delta_{ik}\sum_j\delta_j\otimes u_{ji}$$

We conclude that the multiplicativity of $\Phi$ is equivalent to the following conditions:
$$u_{ji}u_{jk}=\delta_{ik}u_{ji}\quad,\quad\forall i,j,k$$

Regarding now the unitality of $\Phi$, we have the following formula:
$$\Phi(1)
=\sum_i\Phi(\delta_i)
=\sum_{ij}\delta_j\otimes u_{ji}
=\sum_j\delta_j\otimes\left(\sum_iu_{ji}\right)$$

Thus $\Phi$ is unital when the following conditions are satisfied:
$$\sum_iu_{ji}=1\quad,\quad\forall i$$

Finally, the fact that $\Phi$ is a $*$-morphism translates into:
$$u_{ij}=u_{ij}^*\quad,\quad\forall i,j$$

Summing up, in order for $\Phi(\delta_i)=\sum_j\delta_j\otimes u_{ji}$ to be a morphism of $C^*$-algebras, the elements $u_{ij}$ must be projections, summing up to 1 on each row of $u$. Regarding now the preservation of the trace condition, observe that we have:
$$(tr\otimes id)\Phi(\delta_i)=\frac{1}{N}\sum_ju_{ji}$$

Thus the trace is preserved precisely when the elements $u_{ij}$ sum up to 1 on each of the columns of $u$. We conclude from this that $\Phi(\delta_i)=\sum_j\delta_j\otimes u_{ji}$ is a morphism of $C^*$-algebras preserving the trace precisely when $u$ is magic, and since the coaction conditions on $\Phi$ are equivalent to the fact that $u$ must be a corepresentation, this finishes the proof of our claim. But this claim proves all the assertions in the statement.
\end{proof}

As a quite surprising result now, also from Wang \cite{wa1}, we have:

\begin{theorem}
We have an embedding $S_N\subset S_N^+$, given at the algebra level by: 
$$u_{ij}\to\chi\left(\sigma\Big|\sigma(j)=i\right)$$
This is an isomorphism at $N\leq3$, but not at $N\geq4$, where $S_N^+$ is not classical, nor finite.
\end{theorem} 

\begin{proof}
The fact that we have indeed an embedding as above is clear. Regarding now the second assertion, we can prove this in four steps, as follows:

\medskip

\underline{Case $N=2$}. The fact that $S_2^+$ is indeed classical, and hence collapses to $S_2$, is trivial, because the $2\times2$ magic matrices are as follows, with $p$ being a projection:
$$U=\begin{pmatrix}p&1-p\\1-p&p\end{pmatrix}$$

\underline{Case $N=3$}. It is enough to check that $u_{11},u_{22}$ commute. But this follows from:
\begin{eqnarray*}
u_{11}u_{22}
&=&u_{11}u_{22}(u_{11}+u_{12}+u_{13})\\
&=&u_{11}u_{22}u_{11}+u_{11}u_{22}u_{13}\\
&=&u_{11}u_{22}u_{11}+u_{11}(1-u_{21}-u_{23})u_{13}\\
&=&u_{11}u_{22}u_{11}
\end{eqnarray*}

Indeed, by applying the involution to this formula, we obtain from this that we have $u_{22}u_{11}=u_{11}u_{22}u_{11}$ as well, and so we get $u_{11}u_{22}=u_{22}u_{11}$, as desired.

\medskip

\underline{Case $N=4$}. Consider the following matrix, with $p,q$ being projections:
$$U=\begin{pmatrix}
p&1-p&0&0\\
1-p&p&0&0\\
0&0&q&1-q\\
0&0&1-q&q
\end{pmatrix}$$ 

This matrix is then magic, and if we choose $p,q$ as for the algebra $<p,q>$ to be infinite dimensional, we conclude that $C(S_4^+)$ is infinite dimensional as well.

\medskip

\underline{Case $N\geq5$}. Here we can use the standard embedding $S_4^+\subset S_N^+$, obtained at the level of the corresponding magic matrices in the following way:
$$u\to\begin{pmatrix}u&0\\ 0&1_{N-4}\end{pmatrix}$$

Indeed, with this in hand, the fact that $S_4^+$ is a non-classical, infinite compact quantum group implies that $S_N^+$ with $N\geq5$ has these two properties as well. See \cite{wa1}.
\end{proof}

The above results are quite surprising, and you may wonder, okay with all this mathematics, but in practice, how to intuitively accept the fact that $\{1,2,3,4\}$ has an infinity of quantum permutations. Good point, and in answer, get to learn some quantum mechanics, say from Feynman \cite{fey} or Griffiths \cite{gri} or Weinberg \cite{wei}. You will learn many interesting things from there, and above everything, become a modest person.

\section*{13d. Partitions, easiness}

In order to study the quantum permutation group $S_N^+$, we use representation theory. Things here are quite long and advanced, and for full details on what follows, you can check my book \cite{ba1}. We will need the following version of Tannakian duality:

\index{Tannakian duality}
\index{soft Tannakian duality}

\begin{theorem}
The following operations are inverse to each other:
\begin{enumerate}
\item The construction $A\to C$, which associates to any Woronowicz algebra $A$ the tensor category formed by the intertwiner spaces $C_{kl}=Hom(u^{\otimes k},u^{\otimes l})$.

\item The construction $C\to A$, which associates to a tensor category $C$ the Woronowicz algebra $A$ presented by the relations $T\in Hom(u^{\otimes k},u^{\otimes l})$, with $T\in C_{kl}$.
\end{enumerate}
\end{theorem}

\begin{proof}
This is something quite deep, going back to Woronowicz's paper \cite{wo2} in a slightly different form, with the idea being as follows:

\medskip

(1) We have indeed a construction $A\to C$ as above, whose output is a tensor $C^*$-subcategory with duals of the tensor $C^*$-category of Hilbert spaces.

\medskip

(2) We have as well a construction $C\to A$ as above, simply by dividing the free $*$-algebra on $N^2$ variables by the relations in the statement.

\medskip

Regarding now the bijection claim, some elementary algebra shows that $C=C_{A_C}$ implies $A=A_{C_A}$, and also that $C\subset C_{A_C}$ is automatic. Thus we are left with proving $C_{A_C}\subset C$. But this latter inclusion can be proved indeed, by doing some algebra, and using von Neumann's bicommutant theorem, in finite dimensions. See \cite{ba1}.
\end{proof}

We will need as well, following the classical work of Weyl, Brauer and many others, the notion of ``easiness''. Let us start with the following definition:

\index{category of partitions}

\begin{definition}
Let $P(k,l)$ be the set of partitions between an upper row of $k$ points, and a lower row of $l$ points. A set $D=\bigsqcup_{k,l}D(k,l)$ with $D(k,l)\subset P(k,l)$ is called a category of partitions when it has the following properties:
\begin{enumerate}
\item Stability under the horizontal concatenation, $(\pi,\sigma)\to[\pi\sigma]$.

\item Stability under the vertical concatenation, $(\pi,\sigma)\to[^\sigma_\pi]$.

\item Stability under the upside-down turning, $\pi\to\pi^*$.

\item Each set $P(k,k)$ contains the identity partition $||\ldots||$.

\item The set $P(0,2)$ contains the semicircle partition $\cap$.
\end{enumerate}
\end{definition} 

As a basic example, we have the category of all partitions $P$ itself. Other basic examples include the category of pairings $P_2$, or the categories $NC,NC_2$ of noncrossing partitions, and pairings. There are many other examples, and we will be back to this. 

\bigskip

The relation with the Tannakian categories and duality comes from:

\begin{proposition}
Each $\pi\in P(k,l)$ produces a linear map $T_\pi:(\mathbb C^N)^{\otimes k}\to(\mathbb C^N)^{\otimes l}$, 
$$T_\pi(e_{i_1}\otimes\ldots\otimes e_{i_k})=\sum_{j_1\ldots j_l}\delta_\pi\begin{pmatrix}i_1&\ldots&i_k\\ j_1&\ldots&j_l\end{pmatrix}e_{j_1}\otimes\ldots\otimes e_{j_l}$$
with the Kronecker type symbols $\delta_\pi\in\{0,1\}$ depending on whether the indices fit or not. The assignement $\pi\to T_\pi$ is categorical, in the sense that we have
$$T_\pi\otimes T_\sigma=T_{[\pi\sigma]}\quad,\quad 
T_\pi T_\sigma=N^{c(\pi,\sigma)}T_{[^\sigma_\pi]}\quad,\quad
T_\pi^*=T_{\pi^*}$$
where $c(\pi,\sigma)$ are certain integers, coming from the erased components in the middle.
\end{proposition}

\begin{proof}
The concatenation axiom follows from the following computation:
\begin{eqnarray*}
&&(T_\pi\otimes T_\sigma)(e_{i_1}\otimes\ldots\otimes e_{i_p}\otimes e_{k_1}\otimes\ldots\otimes e_{k_r})\\
&=&\sum_{j_1\ldots j_q}\sum_{l_1\ldots l_s}\delta_\pi\begin{pmatrix}i_1&\ldots&i_p\\j_1&\ldots&j_q\end{pmatrix}\delta_\sigma\begin{pmatrix}k_1&\ldots&k_r\\l_1&\ldots&l_s\end{pmatrix}e_{j_1}\otimes\ldots\otimes e_{j_q}\otimes e_{l_1}\otimes\ldots\otimes e_{l_s}\\
&=&\sum_{j_1\ldots j_q}\sum_{l_1\ldots l_s}\delta_{[\pi\sigma]}\begin{pmatrix}i_1&\ldots&i_p&k_1&\ldots&k_r\\j_1&\ldots&j_q&l_1&\ldots&l_s\end{pmatrix}e_{j_1}\otimes\ldots\otimes e_{j_q}\otimes e_{l_1}\otimes\ldots\otimes e_{l_s}\\
&=&T_{[\pi\sigma]}(e_{i_1}\otimes\ldots\otimes e_{i_p}\otimes e_{k_1}\otimes\ldots\otimes e_{k_r})
\end{eqnarray*}

The composition axiom follows from the following computation:
\begin{eqnarray*}
&&T_\pi T_\sigma(e_{i_1}\otimes\ldots\otimes e_{i_p})\\
&=&\sum_{j_1\ldots j_q}\delta_\sigma\begin{pmatrix}i_1&\ldots&i_p\\j_1&\ldots&j_q\end{pmatrix}
\sum_{k_1\ldots k_r}\delta_\pi\begin{pmatrix}j_1&\ldots&j_q\\k_1&\ldots&k_r\end{pmatrix}e_{k_1}\otimes\ldots\otimes e_{k_r}\\
&=&\sum_{k_1\ldots k_r}N^{c(\pi,\sigma)}\delta_{[^\sigma_\pi]}\begin{pmatrix}i_1&\ldots&i_p\\k_1&\ldots&k_r\end{pmatrix}e_{k_1}\otimes\ldots\otimes e_{k_r}\\
&=&N^{c(\pi,\sigma)}T_{[^\sigma_\pi]}(e_{i_1}\otimes\ldots\otimes e_{i_p})
\end{eqnarray*}

Finally, the involution axiom follows from the following computation:
\begin{eqnarray*}
&&T_\pi^*(e_{j_1}\otimes\ldots\otimes e_{j_q})\\
&=&\sum_{i_1\ldots i_p}<T_\pi^*(e_{j_1}\otimes\ldots\otimes e_{j_q}),e_{i_1}\otimes\ldots\otimes e_{i_p}>e_{i_1}\otimes\ldots\otimes e_{i_p}\\
&=&\sum_{i_1\ldots i_p}\delta_\pi\begin{pmatrix}i_1&\ldots&i_p\\ j_1&\ldots& j_q\end{pmatrix}e_{i_1}\otimes\ldots\otimes e_{i_p}\\
&=&T_{\pi^*}(e_{j_1}\otimes\ldots\otimes e_{j_q})
\end{eqnarray*}

Summarizing, our correspondence is indeed categorical.
\end{proof}

In relation with the quantum groups, we have the following notion:

\index{easiness}
\index{easy quantum group}

\begin{definition}
A compact quantum matrix group $G$ is called easy when
$$Hom(u^{\otimes k},u^{\otimes l})=span\left(T_\pi\Big|\pi\in D(k,l)\right)$$
for any colored integers $k,l$, for certain sets of partitions $D(k,l)\subset P(k,l)$, where
$$T_\pi(e_{i_1}\otimes\ldots\otimes e_{i_k})=\sum_{j_1\ldots j_l}\delta_\pi\begin{pmatrix}i_1&\ldots&i_k\\ j_1&\ldots&j_l\end{pmatrix}e_{j_1}\otimes\ldots\otimes e_{j_l}$$
with the Kronecker type symbols $\delta_\pi\in\{0,1\}$ depending on whether the indices fit or not. 
\end{definition}

This is something very classical, coming from old results of Brauer, which state that the groups $O_N,U_N$ are easy, coming respectively from the categories $P_2,\mathcal P_2$ of pairings, and of matching pairings. We refer to \cite{ba1} for the story, and details. In what follows we will only need such Brauer theorems for $S_N,S_N^+$, the statements here being as follows:

\index{Brauer theorem}

\begin{theorem}
We have the following results:
\begin{enumerate}
\item $S_N$ is easy, coming from the category of all partitions $P$.

\item $S_N^+$ is easy, coming from the category of all noncrossing partitions $NC$.
\end{enumerate}
\end{theorem}

\begin{proof}
This is something quite fundamental, with the proof, using the above Tannakian results and subsequent easiness theory, being as follows:

\medskip

(1) $S_N^+$. We know that this quantum group comes from the magic condition. In order to interpret this magic condition, consider the fork partition:
$$Y\in P(2,1)$$

The linear map associated to this fork partition $Y$ is then given by:
$$T_Y(e_i\otimes e_j)=\delta_{ij}e_i$$

Thus, in usual matrix notation, this linear map is given by:
$$T_Y=(\delta_{ijk})_{i,jk}$$

Now given a corepresentation $u$, we have the following formula:
$$(T_Yu^{\otimes 2})_{i,jk}
=\sum_{lm}(T_Y)_{i,lm}(u^{\otimes 2})_{lm,jk}
=u_{ij}u_{ik}$$

On the other hand, we have as well the following formula:
$$(uT_Y)_{i,jk}
=\sum_lu_{il}(T_Y)_{l,jk}
=\delta_{jk}u_{ij}$$

We conclude that we have the following equivalence:
$$T_Y\in Hom(u^{\otimes 2},u)\iff u_{ij}u_{ik}=\delta_{jk}u_{ij},\forall i,j,k$$

The condition on the right being equivalent to the magic condition, we obtain that $S_N^+$ is indeed easy, the corresponding category of partitions being, as desired:
$$D
=<Y>
=NC$$

(2) $S_N$. Here there is no need for new computations, because we have:
$$S_N=S_N^+\cap O_N$$

At the categorical level means that $S_N$ is easy, coming from:
$$<NC,\slash\hskip-2.2mm\backslash>=P$$

Alternatively, if you prefer, we can rewrite the above proof for $S_N^+$, by adding at each step the basic crossing $\slash\hskip-2.2mm\backslash$ next to the fork partition $Y$.
\end{proof}

Let us discuss now the computation of the law of the main character. This computation is the main problem regarding any compact quantum group, as shown by the following result, which summarizes the various motivations for doing this computation:

\index{main character}
\index{Kesten measure}
\index{amenability}
\index{Cayley graph}

\begin{theorem}
Given a Woronowicz algebra $(A,u)$, the law of the main character
$$\chi=\sum_{i=1}^Nu_{ii}$$
with respect to the Haar integration has the following properties:
\begin{enumerate}
\item The moments of $\chi$ are the numbers $M_k=\dim(Fix(u^{\otimes k}))$.

\item $M_k$ counts as well the lenght $p$ loops at $1$, on the Cayley graph of $A$.

\item $law(\chi)$ is the Kesten measure of the associated discrete quantum group.

\item When $u\sim\bar{u}$ the law of $\chi$ is a usual measure, supported on $[-N,N]$.

\item The algebra $A$ is amenable precisely when $N\in supp(law(Re(\chi)))$.

\item Any morphism $f:(A,u)\to (B,v)$ must increase the numbers $M_k$.

\item Such a morphism $f$ is an isomorphism when $law(\chi_u)=law(\chi_v)$.
\end{enumerate}
\end{theorem}

\begin{proof}
All this is quite advanced, the idea being as follows:

\medskip

(1) This comes from the Peter-Weyl type theory in \cite{wo1}, which tells us the number of fixed points of $v=u^{\otimes k}$ can be recovered by integrating the character $\chi_v=\chi_u^k$.

\medskip

(2) This is something true, and well-known, for $A=C^*(\Gamma)$, with $\Gamma=<g_1,\ldots,g_N>$ being a discrete group. In general, the proof is quite similar.

\medskip

(3) This is actually the definition of the Kesten measure, in the case $A=C^*(\Gamma)$, with $\Gamma=<g_1,\ldots,g_N>$ being a discrete group. In general, this follows from (2).

\medskip

(4) The equivalence $u\sim\bar{u}$ translates into $\chi_u=\chi_u^*$, and this gives the first assertion. As for the support claim, this follows from $uu^*=1\implies||u_{ii}||\leq1$, for any $i$.

\medskip

(5) This is the Kesten amenability criterion, which can be established as in the classical case, $A=C^*(\Gamma)$, with $\Gamma=<g_1,\ldots,g_N>$ being a discrete group.

\medskip

(6) This is something elementary, which follows from (1) above, and from the fact that the morphisms of Woronowicz algebras increase the spaces of fixed points.

\medskip

(7) This follows by using (6), and the Peter-Weyl type theory from \cite{wo1}, the idea being that if $f$ is not injective, then it must strictly increase one of the spaces $Fix(u^{\otimes k})$.
\end{proof}

In the case of the symmetric group $S_N$, the character result is as follows:

\index{derangement}
\index{Poisson law}
\index{inclusion-exclusion}
\index{main character}

\begin{theorem}
For the symmetric group $S_N$ the main character counts fixed points,
$$\chi(\sigma)=\#\left\{i\in\{1,\ldots,N\}\Big|\sigma(i)=i\right\}$$
and its law becomes Poisson $(1)$, in the $N\to\infty$ limit.
\end{theorem}

\begin{proof}
This is something very classical, which can be done in 3 steps, as follows:

\medskip

(1) The trace of the permutation matrices $\sigma\in S_N\subset O_N$ being the number of 1 entries, which correspond to fixed points, we have:
$$Tr(\sigma)=\#\left\{i\in\{1,\ldots,N\}\Big|\sigma(i)=i\right\}$$

If we denote by $F_i\subset S_N$ the set of permutations satisfying $\sigma(i)=i$, the number of permutations $\sigma\in S_N$ having no fixed point at all, called derangements, is:
\begin{eqnarray*}
F_\emptyset
&=&|S_N|-\sum_i|F_i|+\sum_{i<j}|F_i\cap F_j|-\ldots\ldots+(-1)^N|F_1\cap\ldots\cap F_N|\\
&=&N!-N\cdot(N-1)!+\binom{N}{2}(N-2)!-\ldots\ldots+(-1)^N\binom{N}{N}1!\\
&=&N!-\frac{N!}{1}+\frac{N!}{2}-\frac{N!}{6}+\ldots\ldots+(-1)^N\frac{N!}{N!}
\end{eqnarray*}

(2) Thus, when dividing by $N!$, and letting $N\to\infty$, we obtain:
$$P(\chi=0)\simeq\frac{1}{e}$$

(3) In fact, the same method gives the following formula, valid for any $k\in\mathbb N$:
$$P(\chi=k)\simeq\frac{1}{ek!}$$

But this shows that $\chi$ becomes Poisson (1) with $N\to\infty$, as claimed.
\end{proof}

Summarizing, we have here some interesting results regarding the classical permutation group $S_N$. In what follows we will present some similar results regarding the quantum permutation group $S_N^+$, and we will discuss the relation between the classical results and the free results, which will complement the easiness theory developed above. In order to include as well $S_N^+$ in our discussion, we will need the following result, with $*$ being the classical convolution, and $\boxplus$ being Voiculescu's free convolution operation \cite{vdn}:

\index{Poisson limit}
\index{free Poisson limit}
\index{freeness}
\index{free Poisson law}
\index{Marchenko-Pastur law}

\begin{theorem}
The following Poisson type limits converge, for any $t>0$,
$$p_t=\lim_{n\to\infty}\left(\left(1-\frac{1}{n}\right)\delta_0+\frac{1}{n}\delta_t\right)^{*n}$$
$$\pi_t=\lim_{n\to\infty}\left(\left(1-\frac{1}{n}\right)\delta_0+\frac{1}{n}\delta_t\right)^{\boxplus n}$$
the limiting measures being the Poisson law $p_t$, and the Marchenko-Pastur law $\pi_t$, 
$$p_t=\frac{1}{e^t}\sum_{k=0}^\infty\frac{t^k\delta_k}{k!}$$
$$\pi_t=\max(1-t,0)\delta_0+\frac{\sqrt{4t-(x-1-t)^2}}{2\pi x}\,dx$$
whose moments are given by the following formulae:
$$M_k(p_t)=\sum_{\pi\in P(k)}t^{|\pi|}\quad,\quad 
M_k(\pi_t)=\sum_{\pi\in NC(k)}t^{|\pi|}$$
The Marchenko-Pastur measure $\pi_t$ is also called free Poisson law.
\end{theorem}

\begin{proof}
This is something quite advanced, related to probability theory, free probability theory, and random matrices, the idea being as follows:

\medskip

(1) The first step is that of finding suitable functional transforms, which linearize the convolution operations in the statement. In the classical case this is the logarithm of the Fourier transform $\log F$, and in the free case this is Voiculescu's $R$-transform.

\medskip

(2) With these tools in hand, the above limiting theorems can be proved in a standard way, a bit as when proving the Central Limit Theorem. The computations give the moment formulae in the statement, and the density computations are standard as well.

\medskip

(3) Finally, in order for the discussion to be complete, what still remains to be explained is the precise nature of the ``liberation'' operation $p_t\to\pi_t$, as well as the random matrix occurrence of $\pi_t$. This is more technical, and we refer here to \cite{bpa}, \cite{mpa}, \cite{vdn}.
\end{proof}

Getting back now to quantum permutations, the results here are as follows:

\index{main character}
\index{truncated character}
\index{Poisson law}
\index{free Poisson law}
\index{Marchenko-Pastur law}
\index{Weingarten formula}

\begin{theorem}
The law of the main character, given by 
$$\chi=\sum_iu_{ii}$$
for $S_N/S_N^+$ becomes $p_1/\pi_1$ with $N\to\infty$. As for the truncated character 
$$\chi_t=\sum_{i=1}^{[tN]}u_{ii}$$
for $S_N/S_N^+$, with $t\in(0,1]$, this becomes $p_t/\pi_t$ with $N\to\infty$.
\end{theorem}

\begin{proof}
This is again something quite technical, the idea being as follows:

\medskip

(1) In the classical case this is well-known, and follows by using the inclusion-exclusion principle, and then letting $N\to\infty$, as in the proof of Theorem 13.23, at $t=1$. 

\medskip

(2) In the free case there is no such simple argument, and we must use what we know about $S_N^+$, namely its easiness property. We know from easiness that we have:
$$Fix(u^{\otimes k})=span(NC(k))$$

On the other hand, a direct computation shows that the partitions in $P(k)$, and in particular those in $NC(k)$, implemented as linear maps via the operation $\pi\to T_\pi$ from Proposition 13.19, become linearly independent with $N\geq k$. Thus, we have:
\begin{eqnarray*}
\int_{S_N^+}\chi^k
&=&\dim\left(Fix(u^{\otimes k})\right)\\
&=&\dim\left(span\left(T_\pi\Big|\pi\in NC(k)\right)\right)\\
&\simeq&|NC(k)|\\
&=&\sum_{\pi\in NC(k)}1^{|\pi|}
\end{eqnarray*}

In the general case now, where our parameter is an arbitrary number $t\in(0,1]$, the above computation does not apply, but we can still get away with Peter-Weyl theory. Indeed, we know from Theorem 13.10 above how to compute the Haar integration of $S_N^+$, out of the knowledge of the fixed point spaces $Fix(u^{\otimes k})$, and in practice, by using easiness, this leads to the following formula, called Weingarten integration formula:
$$\int_{S_N^+}u_{i_1j_1}\ldots u_{i_kj_k}=\sum_{\pi,\sigma\in NC(k)}\delta_\pi(i)\delta_\sigma(j)W_{kN}(\pi,\sigma)$$

Here the $\delta$ symbols are Kronecker type symbols, checking whether the indices fit or not with the partitions, and $W_{kN}=G_{kN}^{-1}$, with $G_{kN}(\pi,\sigma)=N^{|\pi\vee\sigma|}$, where $|.|$ is the number of blocks. Now by using this formula for computing the moments of $\chi_t$, we obtain:
\begin{eqnarray*}
\int_{S_N^+}\chi_t^k
&=&\sum_{i_1=1}^{[tN]}\ldots\sum_{i_k=1}^{[tN]}\int u_{i_1i_1}\ldots u_{i_ki_k}\\
&=&\sum_{\pi,\sigma\in NC(k)}W_{kN}(\pi,\sigma)\sum_{i_1=1}^{[tN]}\ldots\sum_{i_k=1}^{[tN]}\delta_\pi(i)\delta_\sigma(i)\\
&=&\sum_{\pi,\sigma\in NC(k)}W_{kN}(\pi,\sigma)G_{k[tN]}(\sigma,\pi)\\
&=&Tr(W_{kN}G_{k[tN]})
\end{eqnarray*}

The point now is that with $N\to\infty$ the Gram matrix $G_{kN}$, and so the Weingarten matrix $W_{kN}$ too, becomes asymptotically diagonal. We therefore obtain:
$$\int_{S_N^+}\chi_t^k\simeq\sum_{\pi\in NC(k)}t^{|\pi|}$$

Thus, we are led to the conclusion in the statement. For details, see \cite{ba1}.
\end{proof}

\index{Gram matrix}
\index{Weingarten matrix}

As a conclusion to all this, the usual symmetric group $S_N$ has a free analogue $S_N^+$, which is infinite at $N\geq4$. The best way to understand the liberation operation $S_N\to S_N^+$ is via Brauer theorems and easiness. An even better way, which is more advanced, is via probability theory, for the asymptotic law of the main character. All this might seem quite heavy, but hey, we are probably into some kind of quantum mechanics here.

\section*{13e. Exercises} 

There has been a lot of theory in this chapter, and as a best exercise, we can only recommend spending some time with functional analysis, operator theory, operator algebras, Hopf algebras, quantum groups, and of course quantum permutation groups. Here is however an exercise, which would certainly help in relation with all this:

\begin{exercise}
Find an alternative, more conceptual proof for the equality
$$S_3^+=S_3$$
by considering the following morphism, called universal coaction map,
$$\Phi:\mathbb C^3\to\mathbb C^3\otimes C(S_3^+)\quad,\quad 
e_i\to\sum_je_j\otimes u_{ji}$$
then by applying the Fourier transform over the group $\mathbb Z_3$ on the $\mathbb C^3$ part, and then observing that the coefficients of $u$, in Fourier transform, must clearly commute.
\end{exercise}

This might seem a bit twisted, but the exercise hides many conceptual things, to be discovered when working hard for solving it, and once all this done, the whole thing is guaranteed to look and feel quite conceptual. In addition, there is a nice relation here with the Hadamard matrices, and more specifically with the Fourier matrix $F_3$. 

\chapter{Hadamard models}

\section*{14a. The correspondence}

We discuss here the construction of the quantum permutation group $G\subset S_N^+$ associated to a complex Hadamard matrix $H\in M_N(\mathbb C)$. Although the construction $H\to G$ is something very simple, by modern standards, there is a long story with it, as follows:

\bigskip

(1) Everything goes back to an 1983 paper by Popa \cite{pop}, who made the key remark that the pairs of maximal abelian subalgebras (MASA) in the simplest von Neumann algebra, namely the matrix algebra $M_N(\mathbb C)$, are up to conjugation the algebra of diagonal matrices $\Delta\subset M_N(\mathbb C)$ and its conjugate $H\Delta H^*$ by an Hadamard matrix $H\in M_N(\mathbb C)$.

\bigskip

(2) This remark of Popa suggests spending some time in understanding the complex Hadamard matrices $H$, and among the people involved was notably Jones \cite{jo2}, \cite{jo3}, with the far more refined statement, building on Popa's remark, that associated to $H$ is some sort of abstract ``spin model'', whose partition function must be computed.

\bigskip

(3) The Jones finding can be further refined by using quantum groups, somehow in the spirit of the Yang-Baxter equation, with the result that, as announced above, there is a construction $H\to G$, with the quantum group $G$ describing the symmetries of the spin model, and with the representation theory of $G$ computing the partition function.

\bigskip

(4) These latter things go back to work of mine from the late 90s, but took some time to be axiomatized, mainly due to various hesitations in the choice of the formalism, and including a recurrent mistake at $N=4$ too. All this axiomatization work was done in the 00s, and with several other people, like Bichon, Nicoara, Schlenker involved too. 

\bigskip

(5) So, this was for the story, and as a conclusion, we have nowadays a bright, simple construction of type $H\to G$, that we will explain below, and then all sorts of other more technical things that can be explained afterwards, in relation with the work of Jones, Popa and others, and that we will briefly explain too, in what follows.

\bigskip

Getting started now, as a first observation, the complex Hadamard matrices are related to the quantum permutation groups, via the following simple fact:

\index{magic matrix}

\begin{proposition}
If $H\in M_N(\mathbb C)$ is Hadamard, the rank one projections 
$$P_{ij}=Proj\left(\frac{H_i}{H_j}\right)$$
where $H_1,\ldots,H_N\in\mathbb T^N$ are the rows of $H$, form a magic unitary.
\end{proposition}

\begin{proof}
This is clear, the verification for the rows being as follows:
\begin{eqnarray*}
\left<\frac{H_i}{H_j},\frac{H_i}{H_k}\right>
&=&\sum_l\frac{H_{il}}{H_{jl}}\cdot\frac{H_{kl}}{H_{il}}\\
&=&\sum_l\frac{H_{kl}}{H_{jl}}\\
&=&N\delta_{jk}
\end{eqnarray*}

As for the verification for the columns, this is similar, as follows:
\begin{eqnarray*}
\left<\frac{H_i}{H_j},\frac{H_k}{H_j}\right>
&=&\sum_l\frac{H_{il}}{H_{jl}}\cdot\frac{H_{jl}}{H_{kl}}\\
&=&\sum_l\frac{H_{il}}{H_{kl}}\\
&=&N\delta_{ik}
\end{eqnarray*}

Thus, we have indeed a magic unitary, as claimed. 
\end{proof}

The above result suggests the following definition:

\begin{definition}
Associated to $H\in M_N(\mathbb C)$ is the representation
$$\pi:C(S_N^+)\to M_N(\mathbb C)\quad,\quad 
u_{ij}\to Proj\left(\frac{H_i}{H_j}\right)$$
where $H_1,\ldots,H_N\in\mathbb T^N$ are the rows of $H$.
\end{definition}

The representation $\pi$ constructed above is a ``matrix model'' for the algebra $C(S_N^+)$, in the sense that the standard generators $u_{ij}\in C(S_N^+)$, and more generally any element $a\in C(S_N^+)$, gets modelled in this way by an explicit matrix $\pi(a)\in M_N(\mathbb C)$. And the point now is that, given such a model, we have the following notions:

\index{Hopf image}
\index{inner faithfulness}

\begin{definition}
Let $G$ be a compact quantum group, and let $\pi:C(G)\to M_N(\mathbb C)$ be a matrix model for the associated Woronowicz algebra. 
\begin{enumerate}
\item The Hopf image of $\pi$ is the smallest quotient Woronowicz algebra $C(G)\to C(H)$ producing a factorization of type $\pi:C(G)\to C(H)\to M_N(\mathbb C)$.

\item When the inclusion $H\subset G$ is an isomorphism, i.e. when there is no non-trivial factorization as above, we say that $\pi$ is inner faithful.
\end{enumerate}
\end{definition}

As a first observation, in the case where the model is faithful, in the sense that we have an inclusion $\pi:C(G)\subset M_N(\mathbb C)$, the Hopf image is the algebra $C(G)$ itself, and the model is inner faithful as well. However, this situation will not appear often in practice, because the existence of an embedding $C(G)\subset M_N(\mathbb C)$ forces the algebra $C(G)$ to be finite dimensional, and so $G$ to be a finite quantum group, which is something that we cannot expect, in general. At the level of non-trivial examples now, we have:

\medskip

(1) In the case where $G=\widehat{\Gamma}$ is a group dual, the model is as follows:
$$\pi:C(G)=C^*(\Gamma)\to M_N(\mathbb C)$$

Thus, this model must come from a unitary group representation $\rho:\Gamma\to U_N$, and the minimal factorization of $\pi$ is then the one obtained by taking the image:
$$\rho:\Gamma\to\Lambda\subset U_N$$

Also, the model $\pi$ is inner faithful when $\Gamma\subset U_N$. This is the main example for Definition 14.3, which provides intuition, and justifies the terminology as well.

\medskip

(2) Dually, in the case where $G$ is a classical compact group, we have a standard construction of a matrix model for $C(G)$, obtained by taking an arbitrary family of elements $g_1,\ldots,g_N\in G$, and then constructing the following representation:
$$\pi:C(G)\to\ M_N(\mathbb C)\quad,\quad 
f\to\begin{pmatrix}
f(g_1)\\
&\ddots\\
&&f(g_N)
\end{pmatrix}$$

The minimal factorization of $\pi$ is then via the algebra $C(H)$, with:
$$H=\overline{<g_1,\ldots,g_N>}\subset G$$

Also, $\pi$ is inner faithful precisely when $G=H$, and so when:
$$G=\overline{<g_1,\ldots,g_N>}$$

This is the second main example for the construction in Definition 14.3, which provides some further intuition, and once again justifies the terminology as well.

\medskip

In general, the existence and uniqueness of the Hopf image follow by dividing $C(G)$ by a suitable ideal. We refer to \cite{ba1}, \cite{bbi} for more details regarding this construction.  In relation now with the complex Hadamard matrices, we can simply combine Definition 14.2 and Definition 14.3, and we are led in this way into the following notion:

\begin{definition}
To any Hadamard matrix $H\in M_N(\mathbb C)$ we associate the quantum permutation group $G\subset S_N^+$ given by the following Hopf image factorization,
$$\xymatrix{C(S_N^+)\ar[rr]^{\pi}\ar[rd]&&M_N(\mathbb C)\\&C(G)\ar[ur]&}$$
where $\pi(u_{ij})=Proj(H_i/H_j)$, with $H_1,\ldots,H_N\in\mathbb T^N$ being the rows of $H$.
\end{definition}

This was for the general theory. Our claim now is that the construction $H\to G$ is something really useful, with $G$ encoding the combinatorics of $H$, a bit in the same way as $\mathbb Z_N$ encodes the combinatorics of $F_N$. There are several results supporting this, and we will discuss this gradually, in what follows. As a first such result, we have:

\index{tensor product}

\begin{theorem}
The construction $H\to G$ has the following properties:
\begin{enumerate}
\item For $H=F_N$ we obtain the group $G=\mathbb Z_N$, acting on itself.

\item More generally, for $H=F_G$ we obtain the group $G$ itself, acting on itself.

\item For a tensor product $H=H'\otimes H''$ we obtain a product, $G=G'\times G''$.
\end{enumerate}
\end{theorem}

\begin{proof}
All this is standard, and elementary, as follows:

\medskip

(1) The rows of the Fourier matrix $H=F_N$ are given by $H_i=\rho^i$, where $\rho=(1,w,w^2,\ldots,w^{N-1})$, with $w=e^{2\pi i/N}$. Thus, we have the following formula:
$$\frac{H_i}{H_j}=\rho^{i-j}$$

It follows that the corresponding rank 1 projections $P_{ij}=Proj(H_i/H_j)$ form a circulant matrix, all whose entries commute. Since the entries commute, the corresponding quantum group must satisfy $G\subset S_N$. Now by taking into account the circulant property of $P=(P_{ij})$ as well, we are led to the conclusion that we have $G=\mathbb Z_N$.

\medskip

(2) In the general case now, where $H=F_G$, with $G$ being an arbitrary finite abelian group, the result can be proved either by extending the above proof, of by decomposing $G=\mathbb Z_{N_1}\times\ldots\times\mathbb Z_{N_k}$ and using (3) below, whose proof is independent from (1,2).

\medskip

(3) Assume that we have a tensor product $H=H'\otimes H''$, and let $G,G',G''$ be the associated quantum permutation groups. We have then a diagram as follows:
$$\xymatrix@R=50pt@C25pt{
C(S_{N'}^+)\otimes C(S_{N''}^+)\ar[r]&C(G')\otimes C(G'')\ar[r]&M_{N'}(\mathbb C)\otimes M_{N''}(\mathbb C)\ar[d]\\
C(S_N^+)\ar[u]\ar[r]&C(G)\ar[r]&M_N(\mathbb C)
}$$

Here all the maps are the canonical ones, with those on the left and on the right coming from $N=N'N''$. At the level of standard generators, the diagram is as follows:
$$\xymatrix@R=50pt@C65pt{
u_{ij}'\otimes u_{ab}''\ar[r]&w_{ij}'\otimes w_{ab}''\ar[r]&P_{ij}'\otimes P_{ab}''\ar[d]\\
u_{ia,jb}\ar[u]\ar[r]&w_{ia,jb}\ar[r]&P_{ia,jb}
}$$

Now observe that this diagram commutes. We conclude that the representation associated to $H$ factorizes indeed through $C(G')\otimes C(G'')$, and this gives the result.
\end{proof}

Generally speaking, going beyond Theorem 14.5 is a quite difficult question. There are several computations available here, for the most regarding the deformations of the Fourier matrices, and we will be back to this later, in chapter 16 below. At a more abstract level now, one interesting question is that of abstractly characterizing the magic matrices coming from the complex Hadamard matrices, and we have here:

\begin{proposition}
Given an Hadamard matrix $H\in M_N(\mathbb C)$, the vectors 
$$\xi_{ij}=\frac{H_i}{H_j}$$
on which the magic unitary entries $P_{ij}$ project, have the following properties:
\begin{enumerate}
\item $\xi_{ii}=\xi$ is the all-one vector.

\item $\xi_{ij}\xi_{jk}=\xi_{ik}$, for any $i,j,k$.

\item $\xi_{ij}\xi_{kl}=\xi_{il}\xi_{kj}$, for any $i,j,k,l$.
\end{enumerate}
\end{proposition}

\begin{proof}
All these assertions are trivial, by using the formula $\xi_{ij}=H_i/H_j$.
\end{proof}

Let us call now magic basis of a given Hilbert space $H$ any square array of vectors $\xi\in M_N(H)$, all whose rows and columns are orthogonal bases of $H$. With this convention, the above observations lead to the following result, at the magic basis level:

\index{magic basis}

\begin{theorem}
The magic bases $\xi\in M_N(S^{N-1}_\mathbb C)$ coming from the complex Hadamard matrices are those having the following properties:
\begin{enumerate}
\item We have $\xi_{ij}\in\mathbb T^N$, after a suitable rescaling. 

\item The conditions in Proposition 14.6 are satisfied.
\end{enumerate}
\end{theorem}

\begin{proof}
By using the multiplicativity conditions (1,2,3) in Proposition 14.6, we conclude that, up to a rescaling, we must have $\xi_{ij}=\xi_i/\xi_j$, where $\xi_1,\ldots,\xi_N$ is the first row of the magic basis. Together with our assumption $\xi_{ij}\in\mathbb T^N$, this gives the result. 
\end{proof}

\section*{14b. General theory}

Let us keep discussing what happens at the general level. We will need the following result, valid in the general context of the Hopf image construction:

\index{Tannakian category}

\begin{theorem}
Given a matrix model $\pi:C(G)\to M_N(\mathbb C)$, the fundamental corepresentation $v$ of its Hopf image is subject to the Tannakian conditions
$$Hom(v^{\otimes k},v^{\otimes l})=Hom(U^{\otimes k},U^{\otimes l})$$
where $U_{ij}=\pi(u_{ij})$, and where the spaces on the right are taken in a formal sense.
\end{theorem}

\begin{proof}
This is something which follows directly from the definition of the Hopf image, without computations needed, the idea being as follows:

\medskip

(1) Since the morphisms increase the intertwining spaces, when defined either in a representation theory sense, or just formally, we have inclusions as follows:
$$Hom(u^{\otimes k},u^{\otimes l})\subset Hom(U^{\otimes k},U^{\otimes l})$$

More generally, we have such inclusions when replacing $(G,u)$ with any pair producing a factorization of $\pi$. Thus, by Tannakian duality \cite{wo2}, the Hopf image must be given by the fact that the intertwining spaces must be the biggest, subject to these inclusions.

\medskip

(2) On the other hand, since $u$ is biunitary, so is $U$, and it follows that the spaces on the right form a Tannakian category. Thus, we have a quantum group $(H,v)$ given by:
$$Hom(v^{\otimes k},v^{\otimes l})=Hom(U^{\otimes k},U^{\otimes l})$$

By the above discussion, $C(H)$ follows to be the Hopf image of $\pi$, as claimed.
\end{proof}

With the above result in hand, we can now compute the Tannakian category of the Hopf image, in the context of our Hadamard matrix construction. We are led in this way to the following technical statement, going back to Jones \cite{jo3} in an equivalent form, and which reminds a bit the transfer matrices in statistical mechanics:

\begin{theorem}
The Tannakian category of the quantum group $G\subset S_N^+$ associated to a complex Hadamard matrix $H\in M_N(\mathbb C)$ is given by
$$T\in Hom(u^{\otimes k},u^{\otimes l})\iff T^\circ G^{k+2}=G^{l+2}T^\circ$$
where the objects on the right are constructed as follows:
\begin{enumerate}
\item $T^\circ=id\otimes T\otimes id$.

\item $G_{ia}^{jb}=\sum_kH_{ik}\bar{H}_{jk}\bar{H}_{ak}H_{bk}$.

\item $G^k_{i_1\ldots i_k,j_1\ldots j_k}=G_{i_ki_{k-1}}^{j_kj_{k-1}}\ldots G_{i_2i_1}^{j_2j_1}$.
\end{enumerate}
\end{theorem}

\begin{proof}
With the notations in Theorem 14.8, we have the following formula:
$$Hom(u^{\otimes k},u^{\otimes l})=Hom(U^{\otimes k},U^{\otimes l})$$

Here, according to our conventions, the vector space on the right consists by definition of the complex $N^l\times N^k$ matrices $T$, satisfying the following relation:
$$TU^{\otimes k}=U^{\otimes l}T$$ 

If we denote this equality by $L=R$, the left term $L$ is given by:
\begin{eqnarray*}
L_{ij}
&=&(TU^{\otimes k})_{ij}\\
&=&\sum_aT_{ia}U^{\otimes k}_{aj}\\
&=&\sum_aT_{ia}U_{a_1j_1}\ldots U_{a_kj_k}
\end{eqnarray*}

As for the right term $R$, this is given by the following formula:
\begin{eqnarray*}
R_{ij}
&=&(U^{\otimes l}T)_{ij}\\
&=&\sum_bU^{\otimes l}_{ib}T_{bj}\\
&=&\sum_bU_{i_1b_1}\ldots U_{i_lb_l}T_{bj}
\end{eqnarray*}

Consider now the vectors $\xi_{ij}=H_i/H_j$. Since these vectors span the ambient Hilbert space, the equality $L=R$ is equivalent to the following equality:
$$<L_{ij}\xi_{pq},\xi_{rs}>=<R_{ij}\xi_{pq},\xi_{rs}>$$

We use now the following well-known formula, expressing a product of rank one projections $P_1,\ldots,P_k$ in terms of the corresponding image vectors $\xi_1,\ldots,\xi_k$:
$$<P_1\ldots P_kx,y>=<x,\xi_k><\xi_k,\xi_{k-1}>\ldots\ldots<\xi_2,\xi_1><\xi_1,y>$$

This gives the following formula for $L$:
\begin{eqnarray*}
<L_{ij}\xi_{pq},\xi_{rs}>
&=&\sum_aT_{ia}<P_{a_1j_1}\ldots P_{a_kj_k}\xi_{pq},\xi_{rs}>\\
&=&\sum_aT_{ia}<\xi_{pq},\xi_{a_kj_k}>\ldots<\xi_{a_1j_1},\xi_{rs}>\\
&=&\sum_aT_{ia}G_{pa_k}^{qj_k}G_{a_ka_{k-1}}^{j_kj_{k-1}}\ldots G_{a_2a_1}^{j_2j_1}G_{a_1r}^{j_1s}\\
&=&\sum_aT_{ia}G^{k+2}_{rap,sjq}\\
&=&(T^\circ G^{k+2})_{rip,sjq}
\end{eqnarray*}

As for the right term $R$, this is given by:
\begin{eqnarray*}
<R_{ij}\xi_{pq},\xi_{rs}>
&=&\sum_b<P_{i_1b_1}\ldots P_{i_lb_l}\xi_{pq},\xi_{rs}>T_{bj}\\
&=&\sum_b<\xi_{pq},\xi_{i_lb_l}>\ldots<\xi_{i_1b_1},\xi_{rs}>T_{bj}\\
&=&\sum_bG_{pi_l}^{qb_l}G_{i_li_{l-1}}^{b_lb_{l-1}}\ldots G_{i_2i_1}^{b_2b_1}G_{i_1r}^{b_1s}T_{bj}\\
&=&\sum_bG^{l+2}_{rip,sbq}T_{bj}\\
&=&(G^{l+2}T^\circ)_{rip,sjq}
\end{eqnarray*}

Thus, we obtain the formula in the statement. See \cite{bbs}.
\end{proof}

Let us discuss now the computation of the Haar functional for the quantum permutation group $G\subset S_N^+$ associated to a complex Hadamard matrix $H\in M_N(\mathbb C)$. In the general random matrix model context, we have the following formula for the Haar integration functional of the Hopf image, coming from the work of Wang in \cite{wa2}:

\index{truncated integration}
\index{inner faithfulness}

\begin{theorem}
Given an inner faithful model $\pi:C(G)\to M_N(C(T))$, we have
$$\int_G=\lim_{k\to\infty}\frac{1}{k}\sum_{r=1}^k\int_G^r$$
with the truncated integrals on the right being given by the formula
$$\int_G^r=(\varphi\circ\pi)^{*r}$$
where $\varphi=tr\otimes\int_T$ is the random matrix trace on the target algebra.
\end{theorem}

\begin{proof}
As a first observation, there is an obvious similarity here with the Woronowicz construction of the Haar measure, explained in chapter 13. In fact, the above result holds for any model $\pi:C(G)\to B$, with $\varphi\in B^*$ being a faithful trace, and with this picture in hand, the Woronowicz construction corresponds to the case $\pi=id$, and the result itself is therefore a generalization of Woronowicz's existence result for the Haar measure. In order to prove now the result, we can proceed as in chapter 13. If we denote by $\int_G'$ the limit in the statement, we must prove that this limit converges, and that we have:
$$\int_G'=\int_G$$

It is enough to check this on the coefficients of corepresentations, and if we let $v=u^{\otimes k}$ be one of the Peter-Weyl corepresentations, we must prove that we have:
$$\left(id\otimes\int_G'\right)v=\left(id\otimes\int_G\right)v$$

We know from chapter 1 that the matrix on the right is the orthogonal projection onto $Fix(v)$. Regarding now the matrix on the left, this is the orthogonal projection onto the $1$-eigenspace of $(id\otimes\varphi\pi)v$. Now observe that, if we set $V_{ij}=\pi(v_{ij})$, we have:
$$(id\otimes\varphi\pi)v=(id\otimes\varphi)V$$

Thus, as in chapter 13, we conclude that the $1$-eigenspace that we are interested in equals $Fix(V)$. But, according to Theorem 14.8, we have:
$$Fix(V)=Fix(v)$$

Thus, we have proved that we have $\int_G'=\int_G$, as desired.
\end{proof}

In practice now, we are led to the computation of the truncated integrals $\int_G^r$ appearing in the above result, and the formula of these truncated integrals is as follows:

\begin{proposition}
The truncated integrals in Theorem 14.10, namely
$$\int_G^r=(\varphi\circ\pi)^{*r}$$
are given by the following formula, in the orthogonal case, where $u=\bar{u}$,
$$\int_G^ru_{a_1b_1}\ldots u_{a_pb_p}=(T_p^r)_{a_1\ldots a_p,b_1\ldots b_p}$$
with the matrix on the right being given by the formula
$$(T_p)_{i_1\ldots i_p,j_1\ldots j_p}=\left(tr\otimes\int_T\right)(U_{i_1j_1}\ldots U_{i_pj_p})$$
where $U_{ij}=\pi(u_{ij})$ are the images of the standard coordinates in the model.
\end{proposition}

\begin{proof}
This is something straightforward, which comes from the definition of the truncated integrals. Indeed, we have the following computation:
\begin{eqnarray*}
\int_G^ru_{a_1b_1}\ldots u_{a_pb_p}
&=&(\varphi\circ\pi)^{*r}(u_{a_1b_1}\ldots u_{a_pb_p})\\
&=&(\varphi\circ\pi)^{\otimes r}\Delta^{(r)}(u_{a_1b_1}\ldots u_{a_pb_p})\\
&=&(T_p^r)_{a_1\ldots a_p,b_1\ldots b_p}
\end{eqnarray*}

In addition to this, let us mention as well that in the general compact quantum group case, where the condition $u=\bar{u}$ does not necessarily hold, an analogue of the above result holds, by adding exponents $e_1,\ldots,e_p\in\{1,*\}$ everywhere. See \cite{bbi}.
\end{proof}

Regarding now the main character, the result here is as follows:

\index{main character}
\index{truncated main character}

\begin{theorem}
In the context of Theorem 14.10, let $\mu^r$ be the law of the main character $\chi=Tr(u)$ with respect to the truncated integration:
$$\int_G^r=(\varphi\circ\pi)^{*r}$$
\begin{enumerate}
\item The law of the main character is given by the following formula:
$$\mu=\lim_{k\to\infty}\frac{1}{k}\sum_{r=0}^k\mu^r$$

\item The moments of the truncated measure $\mu^r$ are the following numbers:
$$c_p^r=Tr(T_p^r)$$
\end{enumerate}
\end{theorem}

\begin{proof}
These results are both elementary, the proof being as follows:

\medskip

(1) This follows from the general limiting formula in Theorem 14.10.

\medskip

(2) This follows from the formula in Proposition 14.11, by summing the integrals computed there over pairs of equal indices, $a_i=b_i$.
\end{proof}

In connection with the Hadamard matrices, we can use the above technology in order to compute the law of the main character, and also discuss the behavior of the construction $H\to G$ with respect to the various operations on the Hadamard matrices, such as the transposition $H\to H^t$. Following \cite{bbi}, we have the following result, at the general level:

\begin{theorem}
Consider an inner faithful model, as follows:
$$\pi:C(G)\to M_N(\mathbb C)\quad,\quad u_{ij}\to U_{ij}$$
\begin{enumerate}
\item We set $(U'_{kl})_{ij}=(U_{ij})_{kl}$, and we define a model as follows:
$$\widetilde{\rho}:C(U_N^+)\to M_N(\mathbb C)\quad,\quad 
v_{kl}\to U_{kl}'$$

\item We perform the Hopf image construction, as to get a model as follows:
$$\rho:C(G')\to M_N(\mathbb C)$$
\end{enumerate}
The operation $A\to A'$ is then a duality, in the sense that we have $A''=A$, and in the Hadamard matrix case, this duality comes from the operation $H\to H^t$.
\end{theorem}

\begin{proof}
This is something quite technical, the idea being as follows:

\medskip

(1) First, regarding the statement, the quantum group $U_N^+$ is Wang's quantum unitary group, whose standard coordinates are subject to the condition $u^*=u^{-1},u^t=\bar{u}^{-1}$. 

\medskip

(2) Observe that $U'$ is given by $U'=\Sigma U$, where $\Sigma$ is the flip. Thus this matrix is indeed biunitary, and produces a representation $\rho$ as above. 

\medskip

(3) In what regards now the proof, the fact that $A\to A'$ is a duality is clear, and the Hadamard matrix assertion can be proved via algebraic methods. See \cite{bbi}.
\end{proof}

We denote by $D$ the dilation operation for probability measures, or for general $*$-distributions, given by the formula $D_r(law(X))=law(rX)$. Following \cite{bbi}, we have:

\begin{theorem}
Consider the rescaled measure $\eta^r=D_{1/N}(\mu^r)$.
\begin{enumerate}
\item The moments $\gamma_p^r=c_p^r/N^p$ of $\eta^r$ satisfy the following formula:
$$\gamma_p^r(G)=\gamma_r^p(G')$$

\item $\eta^r$ has the same moments as the following matrix:
$$T_r'=T_r(G')$$

\item In the orthogonal case, where $u=\bar{u}$, we have:
$$\eta^r=law(T_r')$$
\end{enumerate}
\end{theorem}

\begin{proof}
All the results follow from Theorem 14.12, as follows:

\medskip

(1) We have the following computation:
\begin{eqnarray*}
c_p^r(A)
&=&\sum_i(T_p)_{i_1^1\ldots i_p^1,i_1^2\ldots i_p^2}\ldots\ldots(T_p)_{i_1^r\ldots i_p^r,i_1^1\ldots i_p^1}\\
&=&\sum_itr(U_{i_1^1i_1^2}\ldots U_{i_p^1i_p^2})\ldots\ldots tr(U_{i_1^ri_1^1}\ldots U_{i_p^ri_p^1})\\
&=&\frac{1}{N^r}\sum_i\sum_j(U_{i_1^1i_1^2})_{j_1^1j_2^1}\ldots(U_{i_p^1i_p^2})_{j_p^1j_1^1}\ldots\ldots(U_{i_1^ri_1^1})_{j_1^rj_2^r}\ldots(U_{i_p^ri_p^1})_{j_p^rj_1^r}
\end{eqnarray*}

In terms of the matrix $(U'_{kl})_{ij}=(U_{ij})_{kl}$, then by permuting the terms in the product on the right, and finally with the changes $i_a^b\leftrightarrow i_b^a,j_a^b\leftrightarrow j_b^a$, we obtain:
\begin{eqnarray*}
c_p^r(A)
&=&\frac{1}{N^r}\sum_i\sum_j(U'_{j_1^1j_2^1})_{i_1^1i_1^2}\ldots(U'_{j_p^1j_1^1})_{i_p^1i_p^2}\ldots\ldots(U'_{j_1^rj_2^r})_{i_1^ri_1^1}\ldots(U'_{j_p^rj_1^r})_{i_p^ri_p^1}\\
&=&\frac{1}{N^r}\sum_i\sum_j(U'_{j_1^1j_2^1})_{i_1^1i_1^2}\ldots(U'_{j_1^rj_2^r})_{i_1^ri_1^1}\ldots\ldots(U'_{j_p^1j_1^1})_{i_p^1i_p^2}\ldots(U'_{j_p^rj_1^r})_{i_p^ri_p^1}\\
&=&\frac{1}{N^r}\sum_i\sum_j(U'_{j_1^1j_1^2})_{i_1^1i_2^1}\ldots(U'_{j_r^1j_r^2})_{i_r^1i_1^1}\ldots\ldots(U'_{j_1^pj_1^1})_{i_1^pi_2^p}\ldots(U'_{j_r^pj_r^1})_{i_r^pi_1^p}
\end{eqnarray*}

On the other hand, if we use again the above formula of $c_p^r(A)$, but this time for the matrix $U'$, and with the changes $r\leftrightarrow p$ and $i\leftrightarrow j$, we obtain:
$$c_r^p(A')\\
=\frac{1}{N^p}\sum_i\sum_j(U'_{j_1^1j_1^2})_{i_1^1i_2^1}\ldots(U'_{j_r^1j_r^2})_{i_r^1i_1^1}\ldots\ldots(U'_{j_1^pj_1^1})_{i_1^pi_2^p}\ldots(U'_{j_r^pj_r^1})_{i_r^pi_1^p}$$

Now by comparing this with the previous formula, we obtain:
$$N^rc_p^r(A)=N^pc_r^p(A')$$ 

Thus we have the following equalities, which give the result:
$$\frac{c_p^r(A)}{N^p}=\frac{c_r^p(A')}{N^r}$$

(2) By using (1) and the formula in Theorem 14.12, we obtain:
$$\frac{c_p^r(A)}{N^p}
=\frac{c_r^p(A')}{N^r}
=\frac{Tr((T'_r)^p)}{N^r}
=tr((T'_r)^p)$$

But this gives the equality of moments in the statement.

\medskip

(3) This follows from the moment equality in (2), and from the standard fact that for self-adjoint variables, the moments uniquely determine the distribution.
\end{proof}

\section*{14c. Von Neumann algebras}

Let us discuss now some applications of the construction $H\to G$, to questions from mathematical physics. We will need some basic von Neumann algebra theory, coming as a complement to the basic $C^*$-algebra theory from chapter 13, as follows:

\index{von Neumann algebra}
\index{factor}
\index{continuous dimension}
\index{hyperfinite factor}

\begin{theorem}
The von Neumann algebras, which are the $*$-algebras of operators
$$A\subset B(H)$$
closed under the weak operator topology, making each $T\to Tx$ continuous, are as follows:
\begin{enumerate}
\item They are exactly the $*$-algebras of operators $A\subset B(H)$ which are equal to their bicommutant, $A=A''$.

\item In the commutative case, these are the algebras $A=L^\infty(X)$, with $X$ measured space, represented on $H=L^2(X)$, up to a multiplicity.

\item If we write the center as $Z(A)=L^\infty(X)$, then we have a decomposition of type $A=\int_XA_x\,dx$, with the fibers $A_x$ having trivial center,  $Z(A_x)=\mathbb C$.

\item The factors, $Z(A)=\mathbb C$, can be fully classified in terms of ${\rm II}_1$ factors, which are those satisfying $\dim A=\infty$, and having a faithful trace $tr:A\to\mathbb C$.

\item The ${\rm II}_1$ factors enjoy the ``continuous dimension geometry'' property, in the sense that the traces of their projections can take any values in $[0,1]$.

\item Among the ${\rm II}_1$ factors, the most important one is the Murray-von Neumann hyperfinite factor $R$, obtained as an inductive limit of matrix algebras.
\end{enumerate}
\end{theorem}

\begin{proof}
This is something quite heavy, the idea being as follows:

\medskip

(1) This is von Neumann's bicommutant theorem, which is well-known in finite dimensions, and whose proof in general is not that complicated, either.

\medskip

(2) It is clear, via basic measure theory, that $L^\infty(X)$ is indeed a von Neumann algebra on $H=L^2(X)$. The converse can be proved as well, by using spectral theory.

\medskip

(3) This is von Neumann's reduction theory main result, whose statement is already quite hard to understand, and whose proof uses advanced functional analysis.

\medskip

(4) This is something heavy, due to Murray-von Neumann and Connes, the idea being that the other factors can be basically obtained via crossed product constructions.

\medskip

(5) This is a gem of functional analysis, with the rational traces being relatively easy to obtain, and with the irrational ones coming from limiting arguments.

\medskip

(6) Once again, heavy results, by Murray-von Neumann and Connes, the idea being that any finite dimensional construction always leads to the same factor, called $R$.
\end{proof}

In relation now with our questions, variations of von Neumann's reduction theory idea, basically using the abelian subalgebra $Z(A)\subset A$, include the use of maximal abelian subalgebras $B\subset A$, called MASA.  In the finite von Neumann algebra case, where we have a trace, the use of orthogonal MASA is a standard method as well, and we have: 

\index{MASA}
\index{orthogonal MASA}

\begin{definition}
A pair of orthogonal MASA inside a von Neumann algebra $A$ with a trace, $tr:A\to\mathbb C$, is a pair of maximal abelian subalgebras
$$B,C\subset A$$
which are orthogonal with respect to the trace, in the sense that we have
$$(B\ominus\mathbb C1)\perp(C\ominus\mathbb C1)$$
with the scalar product being by definition given by $<b,c>=tr(bc^*)$.
\end{definition}

Observe that, by taking into account the multiples of the identity, the orthogonality condition appearing above reformulates as follows:
$$tr(bc)=tr(b)tr(c)$$

The above notion is potentially useful in the infinite dimensional context, in relation with various structure and classification problems for the ${\rm II}_1$ factors. However, as a toy example, we can try and see what happens for the simplest factor that we know, namely the matrix algebra $M_N(\mathbb C)$, with its usual trace. In this context, we have the following surprising observation of Popa \cite{pop}, making the link with the Hadamard matrices:

\begin{theorem}
Up to a conjugation by a unitary, the pairs of orthogonal MASA in the simplest factor, namely the matrix algebra $M_N(\mathbb C)$, are as follows,
$$A=\Delta\quad,\quad 
B=H\Delta H^*$$
with $\Delta\subset M_N(\mathbb C)$ being the diagonal matrices, and with $H\in M_N(\mathbb C)$ being Hadamard.
\end{theorem}

\begin{proof}
Any MASA in $M_N(\mathbb C)$ being conjugated to the diagonal algebra $\Delta$, we can assume, up to conjugation by a unitary, that we have, for a certain $U\in U_N$:
$$A=\Delta\quad,\quad 
B=U\Delta U^*$$  

Now observe that given two diagonal matrices $D,E\in\Delta$, we have:
\begin{eqnarray*}
tr(D\cdot UEU^*)
&=&\frac{1}{N}\sum_i(DUEU^*)_{ii}\\
&=&\frac{1}{N}\sum_{ij}D_{ii}U_{ij}E_{jj}\bar{U}_{ij}\\
&=&\frac{1}{N}\sum_{ij}D_{ii}E_{jj}|U_{ij}|^2
\end{eqnarray*}

Thus, the orthogonality condition $A\perp B$ reformulates as follows:
$$\frac{1}{N}\sum_{ij}D_{ii}E_{jj}|U_{ij}|^2=\frac{1}{N^2}\sum_{ij}D_{ii}E_{jj}$$

But this tells us precisely that the entries $|U_{ij}|$ must have the same absolute value:
$$|U_{ij}|=\frac{1}{\sqrt{N}}$$

Thus the rescaled matrix $H=\sqrt{N}U$ must be Hadamard, as desired.
\end{proof}

Along the same lines, but at a more advanced level, we have the following result:

\index{commuting square}

\begin{theorem}
Given a complex Hadamard matrix $H\in M_N(\mathbb C)$, the diagram formed by the associated pair of orthogonal MASA, namely
$$\xymatrix@R=35pt@C35pt{
\Delta\ar[r]&M_N(\mathbb C)\\
\mathbb C\ar[u]\ar[r]&H\Delta H^*\ar[u] }$$ is a commuting square in the sense of subfactor theory, in the sense that the expectations onto $\Delta,H\Delta H^*$ commute, and their product is the expectation onto $\mathbb C$.
\end{theorem}

\begin{proof}
It follows from definitions that the expectation $E_\Delta:M_N(\mathbb C)\to\Delta$ is the operation which consists in keeping the diagonal, and erasing the rest:
$$M\to M_\Delta$$

Consider now the other expectation, namely:
$$E_{H\Delta H^*}:M_N(\mathbb C)\to H\Delta H^*$$

It is better to identify this with the following expectation, with $U=H/\sqrt{N}$: 
$$E_{U\Delta U^*}:M_N(\mathbb C)\to U\Delta U^*$$

This latter expectation must be given by a formula of type $M\to UX_\Delta U^*$, with $X$ satisfying the following condition:
$$<M,UDU^*>=<UX_\Delta U^*,UDU^*>\quad,\quad\forall D\in\Delta$$

The scalar products being given by $<a,b>=tr(ab^*)$, this condition reads:
$$tr(MUD^*U^*)=tr(X_\Delta D^*)\quad,\quad\forall D\in\Delta$$

Thus $X=U^*MU$, and the formulae of our two expectations are as follows:
\begin{eqnarray*}
E_\Delta(M)&=&M_\Delta\\
E_{U\Delta U^*}(M)&=&U(U^*MU)_\Delta U^*
\end{eqnarray*}

With these formulae in hand, we have the following computation:
\begin{eqnarray*}
(E_\Delta E_{U\Delta U^*}M)_{ij}
&=&\delta_{ij}(U(U^*MU)_\Delta U^*)_{ii}\\
&=&\delta_{ij}\sum_kU_{ik}(U^*MU)_{kk}\bar{U}_{ik}\\
&=&\delta_{ij}\sum_k\frac{1}{N}\cdot(U^*MU)_{kk}\\
&=&\delta_{ij}tr(U^*MU)\\
&=&\delta_{ij}tr(M)\\
&=&(E_\mathbb CM)_{ij}
\end{eqnarray*}

As for the other composition, the computation here is similar, as follows:
\begin{eqnarray*}
(E_{U\Delta U^*}E_\Delta M)_{ij}
&=&(U(U^*M_\Delta U)_\Delta U^*)_{ij}\\
&=&\sum_kU_{ik}(U^*M_\Delta U)_{kk}\bar{U}_{jk}\\
&=&\sum_{kl}U_{ik}\bar{U}_{lk}M_{ll}U_{lk}\bar{U}_{jk}\\
&=&\frac{1}{N}\sum_{kl}U_{ik}M_{ll}\bar{U}_{jk}\\
&=&\delta_{ij}tr(M)\\
&=&(E_\mathbb CM)_{ij}
\end{eqnarray*}

Thus, we have indeed a commuting square, as claimed.
\end{proof}

As a conclusion, all this leads us into commuting squares and subfactor theory. So, let us explain now the basic theory here. As a first object, which will be central in what follows, we have the Temperley-Lieb algebra \cite{tli}, constructed as follows:

\index{Temperley-Lieb algebra}

\begin{definition}
The Temperley-Lieb algebra of index $N\in[1,\infty)$ is defined as
$$TL_N(k)=span(NC_2(k,k))$$
with product given by vertical concatenation, with the rule
$$\bigcirc=N$$
for the closed circles that might appear when concatenating.
\end{definition}

In other words, the algebra $TL_N(k)$, depending on parameters $k\in\mathbb N$ and $N\in[1,\infty)$, is the formal linear span of the noncrossing pairings $\pi\in NC_2(k,k)$. The product operation is obtained by linearity, for the pairings which span $TL_N(k)$ this being the usual vertical concatenation, with the conventions that things go ``from top to bottom'', and that each circle that might appear when concatenating is replaced by a scalar factor, equal to $N$. Observe that there is a connection here with $S_N^+$, and more specifically with the category of noncrossing partitions $NC$ producing $S_N^+$, due to the following fact:

\index{noncrossing partition}

\begin{proposition}
We have bijections
$$NC(k)\simeq NC_2(2k)\simeq NC_2(k,k)$$
constructed by fattening/shrinking and rotating/flattening, as follows:
\begin{enumerate}
\item The application $NC(k)\to NC_2(2k)$ is the ``fattening'' one, obtained by doubling all the legs, and doubling all the strings as well.

\item Its inverse $NC_2(2k)\to NC(k)$ is the ``shrinking'' application, obtained by collapsing pairs of consecutive neighbors.

\item The bijection $NC_2(2k)\simeq NC_2(k,k)$ is obtained by rotating and flattening the noncrossing pairings, in the obvious way.
\end{enumerate}
\end{proposition}

\begin{proof}
The fact that the two operations in (1,2) are indeed inverse to each other is clear, by computing the corresponding two compositions, with the remark that the construction of the fattening operation requires indeed the partitions to be noncrossing. Thus, we are led to the conclusions in the statement.
\end{proof}

\index{subfactor}
\index{Jones projection}
\index{Jones tower}

Getting back now to von Neumann algebras, following Jones \cite{jo1}, consider an inclusion of ${\rm II}_1$ factors, which is actually something quite natural in quantum physics: 
$$A_0\subset A_1$$

We can consider the orthogonal projection $e_1:A_1\to A_0$, and set: 
$$A_2=<A_1,e_1>$$

This procedure, discovered by Jones and called ``basic construction'', can be iterated, and we obtain in this way a whole tower of ${\rm II}_1$ factors, as follows:
$$A_0\subset_{e_1}A_1\subset_{e_2}A_2\subset_{e_3}A_3\subset\ldots\ldots$$

The basic construction is something quite subtle, making deep connections with advanced mathematics and physics. All this was discovered by Jones in the early 80s, and his main result from \cite{jo1}, which came as a big surprise at that time, along with some supplementary fundamental work, done later, in \cite{jo2}, can be summarized as follows:

\index{Jones theorem}
\index{Temperley-Lieb algebra}
\index{planar algebra}
\index{Perron-Frobenius}
\index{ADE graph}

\begin{theorem}
Let $A_0\subset A_1$ be an inclusion of ${\rm II}_1$ factors.
\begin{enumerate}
\item The sequence of Jones projections $e_1,e_2,e_3,\ldots\in B(H)$ produces a Hilbert space representation of the Temperley-Lieb algebra 
$$TL_N\subset B(H)$$
with the parameter being the index of the subfactor, $N=[A_1,A_0]$.

\item The collection $P=(P_k)$ formed by the linear spaces 
$$P_k=A_0'\cap A_k$$
which contains the image of $TL_N$, has a planar algebra structure.

\item The index $N=[A_1,A_0]$, which is by definition a Murray-von Neumann continuous quantity $N\in[1,\infty]$, must satisfy the following condition:
$$N\in\left\{4\cos^2\left(\frac{\pi}{n}\right)\Big|n\in\mathbb N\right\}\cup[4,\infty]$$
That is, in the small index range, the index of subfactors is quantized.
\end{enumerate}
\end{theorem}

\begin{proof}
This is something quite heavy, the idea being as follows:

\medskip

(1) The idea here is that the functional analytic study of the basic construction leads to the conclusion that the sequence of projections $e_1,e_2,e_3,\ldots\in B(H)$ behaves algebrically exactly as the rescaled sequence of diagrams $\varepsilon_1,\varepsilon_2,\varepsilon_3,\ldots\in TL_N$ given by: 
$$\varepsilon_1={\ }^\cup_\cap$$
$$\varepsilon_2=|\!{\ }^\cup_\cap$$
$$\varepsilon_3=||\!{\ }^\cup_\cap$$
$$\vdots$$

But these diagrams generate $TL_N$, and so we have an embedding $TL_N\subset B(H)$, where $H$ is the Hilbert space where our subfactor $A_0\subset A_1$ lives, as claimed.

\medskip

(2) Since the orthogonal projection $e_1:A_1\to A_0$ commutes with $A_0$ we have:
$$e_1\in P_2'$$

By translation we obtain $e_1,\ldots,e_{k-1}\in P_k$ for any $k$, and so we have:
$$TL_N\subset P$$

The point now is that the planar algebra structure of $TL_N$, obtained by composing diagrams, can be shown to extend into an abstract planar algebra structure of $P$.

\medskip

(3) This is something quite surprising, which follows from (1), via some clever positivity considerations, involving the Perron-Frobenius theorem. In order to best comment on what happens, let us record the first few values of the numbers in the statement:
$$4\cos^2\left(\frac{\pi}{3}\right)=1\quad,\quad 
4\cos^2\left(\frac{\pi}{4}\right)=2$$
$$4\cos^2\left(\frac{\pi}{5}\right)=\frac{3+\sqrt{5}}{2}\quad,\quad 
4\cos^2\left(\frac{\pi}{6}\right)=3$$
$$\vdots$$

In order to prove now the result, the first observation is that, when performing a basic construction, we obtain, by trace manipulations on $e_1$:
$$N\notin(1,2)$$

With a double basic construction, we obtain, by trace manipulations on $<e_1,e_2>$:
$$N\notin\left(2,\frac{3+\sqrt{5}}{2}\right)$$

With a triple basic construction, we obtain, by trace manipulations on $<e_1,e_2,e_3>$:
$$N\notin\left(\frac{3+\sqrt{5}}{2},3\right)$$

Thus, we are led to the conclusion in the statement, by a kind of recurrence, involving certain orthogonal polynomials. In practice now, the most elegant way of proving the result is by using the fundamental fact, from (1), that that sequence of Jones projections $e_1,e_2,e_3,\ldots\subset B(H)$ generate a copy of the Temperley-Lieb algebra of index $N$:
$$TL_N\subset B(H)$$

With this result in hand, we must prove that such a representation cannot exist in index $N<4$, unless we are in the following special situation:
$$N=4\cos^2\left(\frac{\pi}{n}\right)$$

But this can be proved by using some suitable trace and positivity manipulations on $TL_N$, as above. Let us mention too that, at a more advanced level, the subfactors having index $N\in[1,4]$ can be classified by ADE diagrams, and the obstruction $N=4\cos^2(\frac{\pi}{n})$ itself comes from the fact that $N$ must be the squared norm of such a graph.
\end{proof}

As before with other advanced operator algebra topics, our explanations here were quite brief. For more on all this, we recommend Jones' original paper \cite{jo1}, then his statistical mechanics paper \cite{jo2} too, and then his planar algebra paper \cite{jo3}.

\section*{14d. Spin models}

In order to explain the connection between the Hadamard matrices and the subfactors, we will need some more subfactor theory, regarding the commuting squares. Consider a commuting square in the sense of subfactor theory, denoted as follows:
$$\xymatrix@R=35pt@C35pt{
C_{01}\ar[r]&C_{11}\\
C_{00}\ar[u]\ar[r]&C_{10}\ar[u]}$$ 

The idea is that any such square $C$ produces a subfactor of the hyperfinite ${\rm II}_1$ factor $R$. And, we will see in what follows that, when applying this construction to the commuting square $C$ associated to a complex Hadamard matrix $H$, the planar algebra of the corresponding subfactor will appear as the planar algebra $P$ of the associated quantum permutation group $G\subset S_N^+$, according to the following scheme:
$$\xymatrix@R=50pt@C=100pt{
G\ar@{-->}[r]&P\\
H\ar@{-->}[u]\ar@{-->}[r]&C\ar@{-->}[u]}$$

Let us begin with some basics. Given a commuting square $C$ as above, under suitable assumptions on the inclusions $C_{00}\subset C_{10},C_{01}\subset C_{11}$, we can perform the basic construction for them, in finite dimensions, and we obtain a whole array of commuting squares:
$$
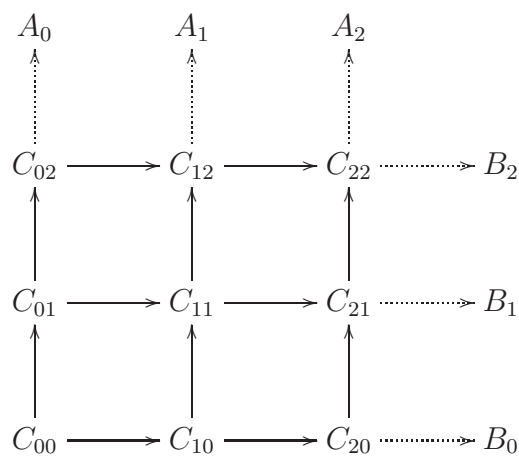
\xymatrix@R=35pt@C35pt{
A_0&A_1&A_2&\\
C_{02}\ar[r]\ar@.[u]&C_{12}\ar[r]\ar@.[u]&C_{22}\ar@.[r]\ar@.[u]&B_2\\
C_{01}\ar[r]\ar[u]&C_{11}\ar[r]\ar[u]&C_{21}\ar@.[r]\ar[u]&B_1\\
C_{00}\ar[u]\ar[r]&C_{10}\ar[u]\ar[r]&C_{20}\ar[u]\ar@.[r]&B_0}$$

Here the various $A,B$ letters stand for the von Neumann algebras obtained in the limit, which are all isomorphic to the hyperfinite ${\rm II}_1$ factor $R$. The point now is that the planar algebra of the associated subfactor can be computed explicitely, as follows:

\index{commuting square}
\index{Ocneanu compactness}

\begin{theorem}
In the context of the above diagram, the following happen:
\begin{enumerate}
\item $A_0\subset A_1$ is a subfactor, and $\{A_i\}$ is the Jones tower for it.

\item The corresponding planar algebra is given by the following formula:
$$A_0'\cap A_k=C_{01}'\cap C_{k0}$$

\item A similar result holds for the ``horizontal'' subfactor $B_0\subset B_1$.
\end{enumerate}
\end{theorem}

\begin{proof}
This is something very standard in subfactor theory, with the result itself being the starting point for various explicit constructions of subfactors, out of concrete combinatorial data, such as the construction of the ADE subfactors mentioned in the above, in the context of the Jones index theorem, the idea being as follows:

\medskip

(1) This is something quite routine, obtained by working out first the axiomatics of the Jones basic construction, and then using this result.

\medskip

(2) This is a subtle result, called Ocneanu compactness theorem \cite{ocn}, which follows by working out the linear algebra of the basic construction.

\medskip

(3) This simply follows from (1,2), by flipping the diagram.
\end{proof}

Getting back now to the Hadamard matrices, we can extend our lineup of results on the associated von Neumann algebraic aspects, namely Theorem 14.17 and Theorem 14.18, with an advanced statement, regarding subfactors, as follows:

\begin{theorem}
Given a complex Hadamard matrix $H\in M_N(\mathbb C)$, the diagram formed by the associated pair of orthogonal MASA, namely
$$\xymatrix@R=35pt@C35pt{
\Delta\ar[r]&M_N(\mathbb C)\\
\mathbb C\ar[u]\ar[r]&H\Delta H^*\ar[u] }$$ 
is a commuting square in the sense of subfactor theory, and the associated planar algebra $P=(P_k)$ is given by the following formula, in terms of $H$ itself,
$$T\in P_k\iff T^\circ G^2=G^{k+2}T^\circ$$
where the objects on the right are constructed as follows:
\begin{enumerate}
\item $T^\circ=id\otimes T\otimes id$.

\item $G_{ia}^{jb}=\sum_kH_{ik}\bar{H}_{jk}\bar{H}_{ak}H_{bk}$.

\item $G^k_{i_1\ldots i_k,j_1\ldots j_k}=G_{i_ki_{k-1}}^{j_kj_{k-1}}\ldots G_{i_2i_1}^{j_2j_1}$.
\end{enumerate}
\end{theorem}

\begin{proof}
We have two assertions here, the idea being as follows:

\medskip

(1) The fact that we have indeed a commuting square is something that we already know, coming from the orthogonal MASA result, explained in Theorem 14.18.

\medskip

(2) The computation of the associated planar algebra is possible thanks to the Ocneanu compactness theorem, corresponding to the formula in Theorem 14.22 (2). To be more precise, by doing some direct computations, which are quite similar to those in the proof of Theorem 14.9, we obtain the formula in the statement. See Jones \cite{jo3}.
\end{proof}

The point now is that all the above is very similar to Theorem 14.9. To be more precise, by comparing the above result with the formula obtained in Theorem 14.9, which is identical, we are led to the following result, clarifying the situation: 

\index{Poincar\'e series}

\begin{theorem}
Let $H\in M_N(\mathbb C)$ be a complex Hadamard matrix.
\begin{enumerate}
\item The planar algebra associated to $H$ is given by 
$$P_k=Fix(u^{\otimes k})$$
where $G\subset S_N^+$ is the associated quantum permutation group.

\item The corresponding Poincar\'e series $f(z)=\sum_k\dim(P_k)z^k$ is
$$f(z)=\int_G\frac{1}{1-z\chi}$$
which is the Stieltjes transform of the law of the main character $\chi=\sum_iu_{ii}$.
\end{enumerate}
\end{theorem}

\begin{proof}
This follows by comparing the quantum group and subfactor results:

\medskip

(1) As already mentioned above, this simply follows by comparing Theorem 14.9 with the subfactor computation in Theorem 14.23. For full details here, we refer to \cite{bbs}.

\medskip

(2) This is a consequence of (1), and of the Peter-Weyl type results from \cite{wo1}, which tell us that fixed points can be counted by integrating characters.
\end{proof}

Summarizing, we have now a clarification of the various quantum algebraic objects associated to a complex Hadamard matrix $H\in M_N(\mathbb C)$, the idea being that the central object, which best encodes the ``symmetries'' of the matrix, and which allows the computation of the other quantum algebraic objects as well, such as the associated planar algebra, is the associated quantum permutation group $G\subset S_N^+$. 

\bigskip

The above results, which are of purely algebraic nature, do not close the discussion, because we still have to understand how the subfactor itself appears from the quantum group. The result here, which is something a bit more technical, is as follows:

\begin{theorem}
The subfactor associated to $H\in M_N(\mathbb C)$ is of the form
$$A^G\subset(\mathbb C^N\otimes A)^G$$
with $A=R\rtimes\widehat{G}$, where $G\subset S_N^+$ is the associated quantum permutation group.
\end{theorem}

\begin{proof}
This is something more technical, the idea being that the basic construction procedure for the commuting squares, explained before Theorem 14.22, can be performed in an ``equivariant setting'', for commuting squares having components as follows:
$$D\otimes_GE=(D\otimes(E\rtimes\widehat{G}))^G$$

To be more precise, starting with a commuting square formed by such algebras, we obtain by basic construction a whole array of commuting squares as follows, with $\{D_i\},\{E_i\}$ being by definition Jones towers, and with $D_\infty,E_\infty$ being their inductive limits:
$$
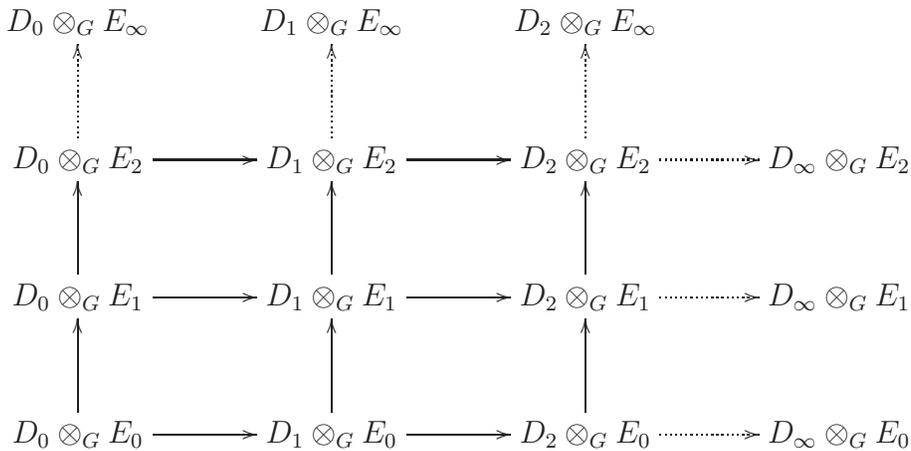
\xymatrix@R=35pt@C35pt{
D_0\otimes_GE_\infty&D_1\otimes_GE_\infty&D_2\otimes_GE_\infty\\
D_0\otimes_GE_2\ar@.[u]\ar[r]&D_1\otimes_GE_2\ar@.[u]\ar[r]&D_2\otimes_GE_2\ar@.[u]\ar@.[r]&D_\infty\otimes_GE_2\\
D_0\otimes_GE_1\ar[u]\ar[r]&D_1\otimes_GE_1\ar[u]\ar[r]&D_2\otimes_GE_1\ar[u]\ar@.[r]&D_\infty\otimes_GE_1\\
D_0\otimes_GE_0\ar[u]\ar[r]&D_1\otimes_GE_0\ar[u]\ar[r]&D_2\otimes_GE_0\ar[u]\ar@.[r]&D_\infty\otimes_GE_0}$$

The point now is that this quantum group picture works in fact for any commuting square having $\mathbb C$ in the lower left corner. In the Hadamard matrix case, that we are interested in here, the corresponding commuting square is as follows:
$$\xymatrix@R=35pt@C35pt{
\mathbb C\otimes_G\mathbb C^N\ar[r]&\mathbb C^N\otimes_G\mathbb C^N\\
\mathbb C\otimes_G\mathbb C\ar[u]\ar[r]&\mathbb C^N\otimes_G\mathbb C\ar[u] }$$ 

Thus, the subfactor obtained by vertical basic construction appears as follows:
$$\mathbb C\otimes_GE_\infty\subset\mathbb C^N\otimes_GE_\infty$$

But this gives the conclusion in the statement, with the ${\rm II}_1$ factor appearing there being by definition $A=E_\infty\rtimes\widehat{G}$, and with the remark that we have $E_\infty\simeq R$.
\end{proof}

All this is of course quite heavy, with the above results being subject to several extensions, and with all this involving several general correspondences between quantum groups, planar algebras, commuting squares and subfactors, that we will not get into.

\bigskip

As a technical comment here, it is possible to deduce Theorem 14.24 directly from Theorem 14.25, via some routine quantum group computations. However, Theorem 14.25 and its proof involve some heavy algebra and functional analysis, coming on top of the heavy algebra and functional analysis required for the general theory of the commuting squares, and this makes the whole thing quite unusable, in practice. 

\bigskip

Thus, while being technically weaker than Theorem 14.25, and dealing with pure algebra only, Theorem 14.24 above remains the main result on the subject.

\bigskip

As already mentioned in the beginning of this book, all this is conjecturally related to statistical mechanics. Indeed, the Tannakian category/planar algebra formula from Theorem 14.23 has many similarities with the transfer matrix computations for the spin models, and this is explained in Jones' paper \cite{jo3}, and known for long before that, from his 1989 paper \cite{jo2}. However, the precise significance of the Hadamard matrices in statistical mechanics, or in related areas such as link invariants, remains a bit unclear.

\bigskip

From a quantum group perspective, the same questions make sense. The idea here, which is old folklore, going back to the 1998 discovery by Wang \cite{wa1} of the quantum permutation group $S_N^+$, is that associated to any 2D spin model should be a quantum permutation group $G\subset S_N^+$, which appears by factorizing the flat representation $C(S_N^+)\to M_N(\mathbb C)$ associated to the $N\times N$ matrix of the Boltzmann weights of the model, and whose representation theory computes the partition function of the model.

\index{lattice model}
\index{spin model}
\index{transfer matrix}
\index{Boltzmann weights}

\bigskip

This is supported on one hand by Jones' theory in \cite{jo2}, \cite{jo3}, via the connecting results presented above, and on the other hand by a number of more recent results, such as those in \cite{bn2}, having similarities with the computations for the Ising and Potts models. However, the whole thing remains not axiomatized, at least for the moment, and in what regards the Hadamard matrices, their precise physical significance remains unclear.

\section*{14e. Exercises} 

As already mentioned, on several occasions, going beyond the above results is a quite difficult task, and we will partly do this in the next two chapters. There are however a few possible exercises, which are doable. Let us start with:

\begin{exercise}
Find the necessary conditions for a magic basis formed by rank $1$ projections to produce a classical quantum group, via the Hopf image construction.
\end{exercise}

Here we use the notion of magic basis, which already appeared in the above, and the application of the Hopf image construction, in order to produce a quantum permutation group, is exactly as in the context of the correspondence $H\to G$ discussed here.

\begin{exercise}
Find the necessary conditions for a magic basis formed by rank $1$ projections to produce a group dual, via the Hopf image construction.
\end{exercise}

As before with the previous exercise, after clearly formulating what precisely is to be done, this can only be a mixture of linear algebra and combinatorics.

\begin{exercise}
Prove that the generalized Fourier matrices $F_G$ are the only ones producing a classical group, or a group dual.
\end{exercise}

In relation now with operator algebras, quantum physics and more, we have the following exercise, which deals with a theme that we have not discussed yet here:

\begin{exercise}
Learn the theory of MUB, and find a relation with the quantum permutation groups.
\end{exercise}

Actually we already met the notion of MUB, in relation with the McNulty-Weigert matrices, in chapter 8 above, and the first thing is therefore to go back there, then find and read the relevant literature. And then, try to solve the exercise. 

\chapter{Generalizations}

\section*{15a. Unitary entries}

We have seen in the previous chapter that associated to any complex Hadamard matrix $H\in M_N(\mathbb C)$ is a certain quantum permutation group $G\subset S_N^+$, which describes the symmetries of the matrix. The main example for this construction $H\to G$ is, as it normally should, $F_N\to\mathbb Z_N$, and more generally, $F_G\to G$. Moreover, we have seen that all this is related to interesting questions from operator algebras, making a potential link with mathematical physics. We discuss here two extensions of the construction $H\to G$, which are both quite interesting, each having its own set of motivations, as follows:

\bigskip

(1) A first idea is that of using Hadamard matrices with noncommutative entries, $H\in M_N(A)$, with $A$ being a $C^*$-algebra. The motivation here comes from the continuous families of complex Hadamard matrices, where $A=C(X)$, and also from all sorts of other constructions involving the complex Hadamard matrices, such as the MUB.

\bigskip

(2) A second idea is that of using partial Hadamard matrices (PHM), with usual complex entries, $H\in M_{M\times N}(\mathbb C)$. Here the motivation comes from the theory of the PHM, developed at various places in this book, and also from the theory of the resulting symmetry-encoding objects $G$, which are certain interesting quantum semigroups.

\bigskip

Technically speaking now, looking at (1) and (2) above certainly suggests that there is room for some unification here, by taking about partial complex Hadamard matrices with noncommutative entries. However, this is something quite theoretical, which has not been done yet. And so again, an interesting question to be put on your to-do list. And with the warning however that, before going head-first into any kind of generalization, you should have some clear motivations, preferably coming from physics. Without clear motivation, if you just want to generalize the construction $H\to G$, you will most likely end up into some terribly complicated and abstract algebra, having 0 uses.

\bigskip

Back to work now, let us begin by discussing (1). Let $A$ be an arbitrary $C^*$-algebra. For most of the applications $A$ will be a commutative algebra, $A=C(X)$ with $X$ being a compact space, or a matrix algebra, $A=M_K(\mathbb C)$ with $K\in\mathbb N$. We will sometimes consider, as a joint generalization, the random matrix algebras $A=M_K(C(X))$. Two row or column vectors over $A$, say $a=(a_1,\ldots,a_N)$ and $b=(b_1,\ldots,b_N)$ by writing both of them horizontally, are called orthogonal when:
$$\sum_ia_ib_i^*=\sum_ia_i^*b_i=0$$

Observe that, by applying the involution, we have as well:
$$\sum_ib_ia_i^*=\sum_ib_i^*a_i=0$$

With this orthogonality notion in hand, we can formulate:

\index{generalized Hadamard matrix}

\begin{definition}
An Hadamard matrix over an arbitrary $C^*$-algebra $A$ is a square matrix $H\in M_N(A)$ such that:
\begin{enumerate}
\item All the entries of $H$ are unitaries, $H_{ij}\in U(A)$.

\item These entries commute on all rows and all columns of $H$.

\item The rows and columns of $H$ are pairwise orthogonal.
\end{enumerate}
\end{definition}

As a first remark, in the simplest case $A=\mathbb C$ the unitary group is the unit circle in the complex plane, $U(\mathbb C)=\mathbb T$, and we obtain the usual complex Hadamard matrices. In the general commutative case, $A=C(X)$ with $X$ compact space, our Hadamard matrix must be formed of ``fibers'', one for each point $x\in X$. Therefore, we obtain:

\begin{proposition}
The Hadamard matrices $H\in M_N(A)$ over a commutative algebra $A=C(X)$ are exactly the families of complex Hadamard matrices of type
$$H=\left\{H^x\Big|x\in X\right\}$$
with $H^x$ depending continuously on the parameter $x\in X$.
\end{proposition}

\begin{proof}
This follows indeed by combining the above two observations. Observe that, when we wrote $A=C(X)$ in the above statement, we used the Gelfand theorem.
\end{proof}

Let us comment now on the above axioms. For $U,V\in U(A)$ the commutation relation $UV=VU$ implies as well the following commutation relations:
$$UV^*=V^*U\quad,\quad 
U^*V=VU^*\quad,\quad 
U^*V^*=U^*V^*$$

Thus the axiom (2) tells us that the $C^*$-algebras $R_1,\ldots,R_N$ and $C_1,\ldots,C_N$ generated by the rows and the columns of $A$ must be all commutative. In view of this, we will be particulary interested in what follows in the following type of matrices:

\index{non-classical matrix}

\begin{definition}
An Hadamard matrix $H\in M_N(A)$ is called ``non-classical'' if the $C^*$-algebra generated by its coefficients is not commutative.
\end{definition}

Let us comment now on the axiom (3). According to our definition of orthogonality there are 4 sets of relations to be satisfied, namely for any $i\neq k$ we must have:
\begin{eqnarray*}
\sum_jH_{ij}H_{kj}^*
&=&\sum_jH_{ij}^*H_{kj}\\
&=&\sum_jH_{ji}H_{jk}^*\\
&=&\sum_jH_{ji}^*H_{jk}\\
&=&0
\end{eqnarray*}

Now since by axiom (1) all the entries $H_{ij}$ are known to be unitaries, we can replace this formula by the following more general equation, valid for any $i,k$:
\begin{eqnarray*}
\sum_jH_{ij}H_{kj}^*
&=&\sum_jH_{ij}^*H_{kj}\\
&=&\sum_jH_{ji}H_{jk}^*\\
&=&\sum_jH_{ji}^*H_{jk}\\
&=&N\delta_{ik}
\end{eqnarray*}

The point now is that everything simplifies in terms of the following matrices:
$$H=(H_{ij})\quad,\quad 
H^*=(H_{ji}^*)\quad,\quad 
H^t=(H_{ji})\quad,\quad 
\bar{H}=(H_{ij}^*)$$

Indeed, the above equations simply read:
$$HH^*=H^*H=H^t\bar{H}=\bar{H}H^t=N1_N$$

So, let us recall now that a square matrix $H\in M_N(A)$ is called ``biunitary'' if both $H$ and $H^t$ are unitaries. In the particular case where $A$ is commutative, $A=C(X)$, we have ``$H$ unitary $\implies$ $H^t$ unitary'', so in this case biunitary means of course unitary. In terms of this notion, we have the following reformulation of Definition 15.1:

\index{biunitary matrix}

\begin{proposition}
Assume that $H\in M_N(A)$ has unitary entries, which commute on all rows and all columns of $H$. Then the following are equivalent:
\begin{enumerate}
\item $H$ is Hadamard.

\item $H/\sqrt{N}$ is biunitary.

\item $HH^*=H^t\bar{H}=N1_N$.
\end{enumerate}
\end{proposition}

\begin{proof}
This basically follows from the above discussion, as follows:

\medskip

-- We know from definitions that the condition (1) in the statement happens if and only if the axiom (3) in Definition 15.1 is satisfied.

\medskip

--  By the above discussion, it follows that this axiom (3) in Definition 15.1 is equivalent to the condition (2) in the statement. 

\medskip

-- Regarding now the equivalence with the condition (3) in the statement, this follows from the commutation axiom (2) in Definition 15.1. 

\medskip

-- By putting now everything together, we see that all the conditions in the statement are indeed equivalent.
\end{proof}

Observe now that if $H=(H_{ij})$ is Hadamard, then so are the following matrices:
$$\bar{H}=(H_{ij}^*)\quad,\quad 
H^t=(H_{ji})\quad,\quad 
H^*=(H_{ji}^*)$$

In addition, we have the following result:

\index{matrix equivalence}

\begin{proposition}
The class of Hadamard matrices $H\in M_N(A)$ is stable under:
\begin{enumerate}
\item Permuting the rows or columns.

\item Multiplying the rows or columns by central unitaries. 
\end{enumerate}
When successively combining these two operations, we obtain an equivalence relation on the class of Hadamard matrices $H\in M_N(A)$.
\end{proposition}

\begin{proof}
This is clear from definitions, exactly as in the usual complex Hadamard matrix case. Observe that in the commutative case $A=C(X)$ any unitary is central, so we can multiply the rows or columns by any unitary. In particular in this case we can always ``dephase'' the matrix, i.e. assume that its first row and column consist of $1$ entries. Note that this operation is not allowed in the general case.
\end{proof}

Let us discuss now the tensor product operation. We have here:

\index{tensor product}

\begin{proposition}
Let $H\in M_N(A)$ and $K\in M_M(A)$ be Hadamard matrices, and assume that $<H_{ij}>$ commutes with $<K_{ab}>$. Then the ``tensor product'' 
$$H\otimes K\in M_{NM}(A)$$
given by $(H\otimes K)_{ia,jb}=H_{ij}K_{ab}$, is an Hadamard matrix.
\end{proposition}

\begin{proof}
This follows from definitions, and is as well a consequence of the more general Theorem 15.7 below, that will be proved with full details.
\end{proof}

Following Di\c t\u a \cite{dit}, the deformed tensor products can be constructed as follows:

\index{deformed tensor product}

\begin{theorem}
Let $H\in M_N(A)$ and $K\in M_M(A)$ be Hadamard matrices, and $Q\in M_{N\times M}(U_A)$. Then the ``deformed tensor product'' $H\otimes_QK\in M_{NM}(A)$, given by
$$(H\otimes_QK)_{ia,jb}=Q_{ib}H_{ij}K_{ab}$$
is an Hadamard matrix as well, provided that the entries of $Q$ commute on rows and columns, and that the algebras $<H_{ij}>$, $<K_{ab}>$, $<Q_{ib}>$ pairwise commute. 
\end{theorem}

\begin{proof}
First, the entries of $L=H\otimes_QK$ are unitaries, and its rows are orthogonal:
\begin{eqnarray*}
\sum_{jb}L_{ia,jb}L_{kc,jb}^*
&=&\sum_{jb}Q_{ib}H_{ij}K_{ab}\cdot Q_{kb}^*K_{cb}^*H_{kj}^*\\
&=&N\delta_{ik}\sum_bQ_{ib}K_{ab}\cdot Q_{kb}^*K_{cb}^*\\
&=&N\delta_{ik}\sum_jK_{ab}K_{cb}^*\\
&=&NM\cdot\delta_{ik}\delta_{ac}
\end{eqnarray*}

The orthogonality of columns can be checked as follows:
\begin{eqnarray*}
\sum_{ia}L_{ia,jb}L_{ia,kc}^*
&=&\sum_{ia}Q_{ib}H_{ij}K_{ab}\cdot Q_{ic}^*K_{ac}^*H_{ik}^*\\
&=&M\delta_{bc}\sum_iQ_{ib}H_{ij}\cdot Q_{ic}^*H_{ik}^*\\
&=&M\delta_{bc}\sum_iH_{ij}H_{ik}^*\\
&=&NM\cdot\delta_{jk}\delta_{bc}
\end{eqnarray*}

For the commutation on rows we use in addition the commutation on rows for $Q$:
\begin{eqnarray*}
L_{ia,jb}L_{kc,jb}
&=&Q_{ib}H_{ij}K_{ab}\cdot Q_{kb}H_{kj}K_{cb}\\
&=&Q_{ib}Q_{kb}\cdot H_{ij}H_{kj}\cdot K_{ab}K_{cb}\\
&=&Q_{kb}Q_{ib}\cdot H_{kj}H_{ij}\cdot K_{cb}K_{ab}\\
&=&Q_{kb}H_{kj}K_{cb}\cdot Q_{ib}H_{ij}K_{ab}\\
&=&L_{kc,jb}L_{ia,jb}
\end{eqnarray*}

The commutation on columns is similar, using the commutation on columns for $Q$:
\begin{eqnarray*}
L_{ia,jb}L_{ia,kc}
&=&Q_{ib}H_{ij}K_{ab}\cdot q_{ic}H_{ik}K_{ac}\\
&=&Q_{ib}Q_{ic}\cdot H_{ij}H_{ik}\cdot K_{ab}K_{ac}\\
&=&Q_{ic}Q_{ib}\cdot H_{ik}H_{ij}\cdot K_{ac}K_{ab}\\
&=&Q_{ic}H_{ik}K_{ac}\cdot Q_{ib}H_{ij}K_{ab}\\
&=&L_{ia,kc}L_{ia,jb}
\end{eqnarray*}

Thus all the axioms are satisfied, and $L$ is indeed Hadamard.
\end{proof}

As a basic example, we have the following construction:

\begin{proposition}
The following matrix is Hadamard,
$$M=\begin{pmatrix}x&y&x&y\\ x&-y&x&-y\\ z&t&-z&-t\\ z&-t&-z&t\end{pmatrix}$$
for any unitaries $x,y,z,t$ satisfying the following condition:
$$[x,y]=[x,z]=[y,t]=[z,t]=0$$
\end{proposition}

\begin{proof}
This follows indeed from Theorem 15.7, because we have:
$$\begin{pmatrix}1&1\\ 1&-1\end{pmatrix}\otimes_{\begin{pmatrix}x&y\\ z&t\end{pmatrix}}\begin{pmatrix}1&1\\ 1&-1\end{pmatrix}
=\begin{pmatrix}x&y&x&y\\ x&-y&x&-y\\ z&t&-z&-t\\ z&-t&-z&t\end{pmatrix}$$

In addition, the commutation relations in Theorem 15.7 are satisfied indeed.
\end{proof}

The usual complex Hadamard matrices were classified by Haagerup in \cite{ha1} at $N=2,3,4,5$. In what follows we investigate the case of the general Hadamard matrices. We use the equivalence relation constructed in Proposition 15.5. We first have:

\begin{proposition}
The $2\times 2$ Hadamard matrices are all classical, and are all equivalent to the Fourier matrix $F_2$.
\end{proposition}

\begin{proof}
Consider indeed an arbitrary $2\times 2$ Hadamard matrix:
$$H=\begin{pmatrix}A&B\\ C&D\end{pmatrix}$$

We already know that $A,D$ each commute with $B,C$. Also, we have:
$$AB^*+CD^*=0$$

We deduce that $A=-CD^*B$ commutes with $D$, and that $C=-AB^*D$ commutes with $B$. Thus our matrix is classical, any since all unitaries are now central, we can dephase our matrix, which follows therefore to be the Fourier matrix $F_2$.
\end{proof}

Let us discuss now the case $N=3$. Here the classification in the classical case uses the key fact that any formula of type $a+b+c=0$, with $|a|=|b|=|c|=1$, must be, up to a permutation of terms, a ``trivial'' formula of the following type, with $j=e^{2\pi i/3}$:
$$a+ja+j^2a=0$$

Here is the noncommutative analogue of this simple fact:

\begin{proposition}
Assume that we have a vanishing sum of unitaries:
$$a+b+c=0$$ 
Then this sum must be of the following special type, 
$$a+wa+w^2a=0$$
with $w$ being a unitary satisfying $1+w+w^2=0$.
\end{proposition}

\begin{proof}
Since $-c=a+b$ is unitary we have the following formula:
$$(a+b)(a+b)^*=1$$

Thus we have $ab^*+ba^*=-1$, and so we obtain:
$$ab^*ba^*+(ba^*)^2=-ba^*$$

But with $w=ba^*$ we obtain from this equality that we have:
$$1+w^2=-w$$

Thus, we are led to the conclusion in the statement. 
\end{proof}

With the above result in hand, we can start the $N=3$ classification. We first have the following technical result, that we will improve later on:

\begin{proposition}
Any $3\times 3$ Hadamard matrix must be of the form
$$H=\begin{pmatrix}a&b&c\\ ua&uv^*w^2vb&uv^*wvc\\ va&wvb&w^2vc\end{pmatrix}$$
with $w$ being subject to the equation $1+w+w^2=0$.
\end{proposition}

\begin{proof}
Consider an arbitrary Hadamard matrix $H\in M_3(A)$. We define $a,b,c,u,v,w$ as for that part of the matrix to be exactly as in the statement, as follows:
$$H=\begin{pmatrix}a&b&c\\ ua&x&y\\ va&wvb&z\end{pmatrix}$$

Let us look first at the scalar product between the first and third row:
$$vaa^*+wvbb^*+zc^*=0$$

By simplifying we obtain $v+wv+zc^*=0$, and by using Proposition 15.10 we conclude that we have $1+w+w^2=0$, and that $zc^*=w^2v$, and so $z=w^2vc$, as claimed. The scalar products of the first column with the second and third ones are:
$$a^*b+a^*u^*x+a^*v^*wvb=0$$
$$a^*c+a^*u^*y+a^*v^*w^2vc=0$$

By multiplying to the left by $va$, and to the right by $b^*v^*$ and $c^*v^*$, we obtain:
$$1+vu^*xb^*v^*+w=0$$
$$1+vu^*yc^*v^*+w^2=0$$

Now by using Proposition 15.10 again, we obtain:
$$vu^*xb^*v^*=w^2$$
$$vu^*yc^*v^*=w$$

Thus $x=uv^*w^2vb$ and $y=uv^*wvc$, and we are done.
\end{proof}

We can already deduce now a first classification result, as follows:

\begin{proposition}
There is no Hadamard matrix $H\in M_3(A)$ with self-adjoint entries.
\end{proposition}

\begin{proof}
We use Proposition 15.11. Since the entries are idempotents, we have:
$$a^2=b^2=c^2=u^2=v^2=(uw)^2=(vw)^2=1$$

It follows that our matrix is in fact of the following form:
$$H=\begin{pmatrix}a&b&c\\ ua&uwb&uw^2c\\ va&wvb&w^2vc\end{pmatrix}$$

The commutation between $H_{22},H_{23}$ reads:
\begin{eqnarray*}
[uwb,wvb]=0
&\implies&[uw,wv]=0\\
&\implies&uwwv=wvuw\\
&\implies&uvw=vuw^2\\
&\implies&w=1
\end{eqnarray*}

Thus we have reached to a contradiction, and we are done.
\end{proof}

Let us go back now to the general case. We have the following technical result, which refines Proposition 15.11, and which will be in turn further refined, later on:

\begin{proposition}
Any $3\times 3$ Hadamard matrix must be of the form
$$H=\begin{pmatrix}a&b&c\\ ua&w^2ub&wuc\\ va&wvb&w^2vc\end{pmatrix}$$
where $(a,b,c)$ and $(u,v,w)$ are triples of commuting unitaries, and $1+w+w^2=0$.
\end{proposition}

\begin{proof}
We use Proposition 15.11. With $e=uv^*$, the matrix there becomes:
$$H=\begin{pmatrix}a&b&c\\ eva&ew^2vb&ewvc\\ va&wvb&w^2vc\end{pmatrix}$$

The commutation relation between $H_{22},H_{32}$ reads:
\begin{eqnarray*}
[ew^2vb,wvb]=0
&\implies&[ew^2v,wv]=0\\
&\implies&ew^2vwv=wvew^2v\\
&\implies&ew^2v=wvew\\
&\implies&[ew,wv]=0
\end{eqnarray*}

Similarly, the commutation between $H_{23},H_{33}$ reads:
\begin{eqnarray*}
[ewvc,w^2vc]=0
&\implies&[ewv,w^2v]=0\\
&\implies&ewvw^2v=w^2vewv\\
&\implies&ewv=w^2vew^2\\
&\implies&[ew^2,w^2v]=0
\end{eqnarray*}

We can rewrite this latter relation by using the formula $w^2=-1-w$, and then, by further processing it by using the first relation, we obtain:
\begin{eqnarray*}
[e(1+w),(1+w)v]=0
&\implies&[e,wv]+[ew,v]=0\\
&\implies&2ewv-wve-vew=0\\
&\implies&ewv=\frac{1}{2}(wve+vew)
\end{eqnarray*}

We use now the key fact that when an average of two unitaries is unitary, then the three unitaries involved are in fact all equal. This gives:
$$ewv=wve=vew$$

Thus we obtain $[w,e]=[w,v]=0$, so $w,e,v$ commute. Our matrix becomes:
$$H=\begin{pmatrix}a&b&c\\ eva&w^2evb&wevc\\ va&wvb&w^2vc\end{pmatrix}$$

Now by remembering that $u=ev$, this gives the formula in the statement.
\end{proof}

We can now formulate our main classification result, as follows:

\begin{theorem}
The $3\times 3$ Hadamard matrices are all classical, and are all equivalent to the Fourier matrix $F_3$.
\end{theorem}

\begin{proof}
We know from Proposition 15.13 that we can write our matrix in the following way, where $(a,b,c)$ and $(u,v,w)$ pairwise commute, and where $1+w+w^2=0$:
$$H=\begin{pmatrix}a&b&c\\ au&buw&cuw^*\\ av&bvw^*&cvw\end{pmatrix}$$

We also know that $(a,u,v)$, $(b,uw,vw^*)$, $(c,uw^*,vw)$ and $(ab,ac,bc,w)$ have entries which pairwise commute. We first show that $uv$ is central. Indeed, we have:
\begin{eqnarray*}
buv
&=&buvww^*\\
&=&b(uw)(vw^*)\\
&=&(uw)(vw^*)b\\
&=&uvb
\end{eqnarray*}

Similarly, $cuv=uvc$. It follows that we may in fact suppose that $uv$ is a scalar. But since our relations are homogeneous, we may assume in fact that $u=v^*$. Let us prove now that we have $[abc,vw^*]=0$. Indeed, we have the following computation:
\begin{eqnarray*} 
abc
&=&a(bc)ww^*\\
&=&aw(bc)w^*\\
&=&av(wv^*)bcw^*\\
&=&avb(wv^*)cw^*\\
&=&v(ab)wv^*cw^*\\
&=&vw(ab)v^*cw^*\\
&=&vw(ab)w(w^*v^*)cw^*\\
&=&vw^2(ab)c(w^*v^*)w^*\\
&=&vw^*abcv^*w
\end{eqnarray*}

We know also that $[b,vw^*]=0$. Hence $[ac,vw^*]=0$. But $[ac,w^*]=0$. Hence $[ac,v]=0$. But $[a,v]=0$. Hence $[c,v]=0$. But $[c,vw]=0$. So $[c,w]=0$. But $[bc,w]=0$. So $[b,w]=0$. But $[b,v^*w]=0$ and $[ab,w]=0$, so respectively $[b,v]=0$ and $[a,w]=0$. Thus all operators $a,b,c,v,w$ pairwise commute, and we are done.
\end{proof}

At $N=4$ now, the classification work for the usual complex Hadamard matrices uses the fact that an equation of type $a+b+c+d=0$ with $|a|=|b|=|c|=|d|=1$ must be, up to a permutation of the terms, a ``trivial'' equation of the following form:
$$a-a+b-b=0$$

In our setting, however, we have for instance:
$$\begin{pmatrix}a&0\\ 0&x\end{pmatrix}+\begin{pmatrix}-a&0\\ 0&y\end{pmatrix}+\begin{pmatrix}b&0\\ 0&-x\end{pmatrix}+\begin{pmatrix}-b&0\\ 0&-y\end{pmatrix}=0$$

It is probably possible to further complicate this kind of identity, and this makes the $N=4$ classification a quite difficult task. As for the case $N=5$ or higher, things here are most likely very complicated, and we will stop our classification work here.

\section*{15b. Quantum groups}

With the above basic theory developed, let us get now to the point where we wanted to get. The generalized Hadamard matrices produce quantum groups, as follows:

\index{magic matrix}

\begin{theorem}
If $H\in M_N(A)$ is Hadamard, the following matrices $P_{ij}\in M_N(A)$ form altogether a magic matrix $P=(P_{ij})$, over the algebra $M_N(A)$:
$$(P_{ij})_{ab}=\frac{1}{N}H_{ia}H_{ja}^*H_{jb}H_{ib}^*$$
Thus, we can let $\pi:C(S_N^+)\to M_N(A)$ be the representation associated to $P$, mapping $u_{ij}\to P_{ij}$, and then factorize this representation as follows,
$$\pi:C(S_N^+)\to C(G)\to M_N(A)$$
with the closed subgroup $G\subset S_N^+$ chosen minimal.
\end{theorem}

\begin{proof}
The magic condition can be checked in three steps, as follows:

\medskip

(1) Let us first check that each $P_{ij}$ is a projection, i.e. that we have $P_{ij}=P_{ij}^*=P_{ij}^2$. Regarding the first condition, namely $P_{ij}=P_{ij}^*$, this simply follows from:
\begin{eqnarray*}
(P_{ij})_{ba}^*
&=&\frac{1}{N}(H_{ib}H_{jb}^*H_{ja}H_{ia}^*)^*\\
&=&\frac{1}{N}H_{ia}H_{ja}^*H_{jb}H_{ib}^*\\
&=&(P_{ij})_{ab}
\end{eqnarray*}

As for the second condition, $P_{ij}=P_{ij}^2$, this follows from the fact that all the entries $H_{ij}$ are assumed to be unitaries, i.e. follows from axiom (1) in Definition 15.1:
\begin{eqnarray*}
(P_{ij}^2)_{ab}
&=&\sum_c(P_{ij})_{ac}(P_{ij})_{cb}\\
&=&\frac{1}{N^2}\sum_cH_{ia}H_{ja}^*H_{jc}H_{ic}^*H_{ic}H_{jc}^*H_{jb}H_{ib}^*\\
&=&\frac{1}{N}H_{ia}H_{ja}^*H_{jb}H_{ib}^*\\
&=&(P_{ij})_{ab}
\end{eqnarray*}

(2) Let us check now that fact that the entries of $P$ sum up to 1 on each row. For this purpose we use the equality $H^*H=N1_N$, coming from the axiom (3), which gives:
\begin{eqnarray*}
(\sum_jP_{ij})_{ab}
&=&\frac{1}{N}\sum_jH_{ia}H_{ja}^*H_{jb}H_{ib}^*\\
&=&\frac{1}{N}H_{ia}(H^*H)_{ab}H_{ib}^*\\
&=&\delta_{ab}H_{ia}H_{ib}^*\\
&=&\delta_{ab}
\end{eqnarray*}

(3) Finally, let us check that the entries of $P$ sum up to 1 on each column. This is the tricky check, because it involves, besides axiom (1) and the formula $H^t\bar{H}=N1_N$ coming from axiom (3), the commutation on the columns of $H$, coming from axiom (2):
\begin{eqnarray*}
(\sum_iP_{ij})_{ab}
&=&\frac{1}{N}\sum_iH_{ia}H_{ja}^*H_{jb}H_{ib}^*\\
&=&\frac{1}{N}\sum_iH_{ja}^*H_{ia}H_{ib}^*H_{jb}\\
&=&\frac{1}{N}H_{ja}^*(H^t\bar{H})_{ab}H_{jb}\\
&=&\delta_{ab}H_{ja}^*H_{jb}\\
&=&\delta_{ab}
\end{eqnarray*}

Thus $P$ is indeed a magic matrix in the above sense, and we are done.
\end{proof}

As an illustration, consider a usual Hadamard matrix $H\in M_N(\mathbb C)$. If we denote its rows by $H_1,\ldots,H_N$ and we consider the vectors $\xi_{ij}=H_i/H_j$, then we have:
$$\xi_{ij}=\left(\frac{H_{i1}}{H_{j1}},\ldots,\frac{H_{iN}}{H_{jN}}\right)$$

Thus the orthogonal projection on this vector $\xi_{ij}$ is given by:
\begin{eqnarray*}
(P_{\xi_{ij}})_{ab}
&=&\frac{1}{||\xi_{ij}||^2}(\xi_{ij})_a\overline{(\xi_{ij})_b}\\
&=&\frac{1}{N}H_{ia}H_{ja}^*H_{jb}H_{ib}^*\\
&=&(P_{ij})_{ab}
\end{eqnarray*}

We conclude that we have $P_{ij}=P_{\xi_{ij}}$ for any $i,j$, so our construction from Theorem 15.15 is compatible with the construction for the usual complex Hadamard matrices.

\bigskip

Let us discuss now the computation of the quantum permutation groups associated to the deformed tensor products of Hadamard matrices. This is actually something that we have not discussed in chapter 14, when talking about the usual Hadamard models, so the results below are relevant even in the case of these usual models. Let us begin with a study of the associated magic unitary. We have here the following result:

\index{deformed tensor product}

\begin{proposition}
The magic unitary associated to $H\otimes_QK$ is given by
$$P_{ia,jb}=R_{ij}\otimes\frac{1}{N}(Q_{ic}Q_{jc}^*Q_{jd}Q_{id}^*\cdot K_{ac}K_{bc}^*K_{bd}K_{ad}^*)_{cd} $$
where $R_{ij}$ is the magic unitary matrix associated to $H$.
\end{proposition}

\begin{proof}
With standard conventions for deformed tensor products and for double indices, the entries of $L=H\otimes_QK$ are by definition the following elements: 
$$L_{ia,jb}=Q_{ib}H_{ij}K_{ab}$$

Thus the projections $P_{ia,jb}$ constructed in Theorem 15.15 are given by:
\begin{eqnarray*}
(P_{ia,jb})_{kc,ld}
&=&\frac{1}{MN}L_{ia,kc}L_{jb,kc}^*L_{jb,ld}L_{ia,ld}^*\\
&=&\frac{1}{MN}(Q_{ic}H_{ik}K_{ac})(Q_{jc}H_{jk}K_{bc})^*(Q_{jd}H_{jl}K_{bd})(Q_{id}H_{il}K_{ad})^*\\
&=&\frac{1}{MN}(Q_{ic}Q_{jc}^*Q_{jd}Q_{id}^*)(H_{ik}H_{jk}^*H_{jl}H_{il}^*)(K_{ac}K_{bc}^*K_{bd}K_{ad}^*)
\end{eqnarray*}

In terms now of the standard matrix units $e_{kl},e_{cd}$, we have:
\begin{eqnarray*}
&&P_{ia,jb}\\
&=&\frac{1}{MN}\sum_{kcld}e_{kl}\otimes e_{cd}\otimes(Q_{ic}Q_{jc}^*Q_{jd}Q_{id}^*)(H_{ik}H_{jk}^*H_{jl}H_{il}^*)(K_{ac}K_{bc}^*K_{bd}K_{ad}^*)\\
&=&\frac{1}{MN}\sum_{kcld}\left(e_{kl}\otimes 1\otimes H_{ik}H_{jk}^*H_{jl}H_{il}^*\right)(1\otimes e_{cd}\otimes Q_{ic}Q_{jc}^*Q_{jd}Q_{id}^*\cdot K_{ac}K_{bc}^*K_{bd}K_{ad}^*)
\end{eqnarray*}

Since the quantities on the right commute, this gives the formula in the statement.
\end{proof}

In order to investigate the Di\c t\u a deformations, we use:

\index{free wreath product}

\begin{definition}
Let $C(S_M^+)\to A$ and $C(S_N^+)\to B$ be Hopf algebra quotients, with fundamental corepresentations denoted $u,v$. We let
$$A*_wB=A^{*N}*B/<[u_{ab}^{(i)},v_{ij}]=0>$$
with the Hopf algebra structure making $w_{ia,jb}=u_{ab}^{(i)}v_{ij}$ a corepresentation.
\end{definition}

The fact that we have indeed a Hopf algebra follows from the fact that $w$ is magic. In terms of quantum groups, if $A=C(G)$, $B=C(H)$, we write $A*_wB=C(G\wr_*H)$:
$$C(G)*_wC(H)=C(G\wr_*H)$$

The $\wr_*$ operation is the free analogue of $\wr$, the usual wreath product, and we refer for instance to \cite{ba1} for more on this. With this convention, we have the following result:

\begin{theorem}
The representation associated to $L=H\otimes_QK$ factorizes as
$$\xymatrix{C(S_{NM}^+)\ar[rr]^{\pi_L}\ar[rd]&&M_{NM}(\mathbb C)\\&C(S_M^+\wr_*G_H)\ar[ur]&}$$
and so the quantum group associated to $L$ appears as a subgroup $G_L\subset S_M^+\wr_*G_H$.
\end{theorem}

\begin{proof}
We use the formula in Proposition 15.16. For simplifying writing we agree to use instead of expressions of type $H_{ia}H_{ja}^*H_{jb}H_{ib}^*$, fractions as follows, by keeping in mind that the variables are only subject to the commutation relations in Definition 15.1:
$$\frac{H_{ia}H_{jb}}{H_{ja}H_{ib}}$$

Our claim is that the factorization can be indeed constructed, as follows:
$$U_{ab}^{(i)}=\sum_jP_{ia,jb}\quad,\quad 
V_{ij}=\sum_aP_{ia,jb}$$

Indeed, we have three verifications to be made, as follows:

\medskip

(1) We must prove that the elements $V_{ij}=\sum_aP_{ia,jb}$ do not depend on $b$, and generate a copy of $C(G_H)$. But if we denote by $(R_{ij})$ the magic matrix for $H$, we have indeed:
\begin{eqnarray*}
V_{ij}
&=&\frac{1}{N}\left(\frac{Q_{ic}Q_{jd}}{Q_{id}Q_{jc}}\cdot\frac{H_{ik}H_{jl}}{H_{il}H_{jk}}\cdot\delta_{cd}\right)_{kc,ld}\\
&=&((R_{ij})_{kl}\delta_{cd})_{kc,ld}\\
&=&R_{ij}\otimes 1
\end{eqnarray*}

(2) We prove now that for any $i$, the elements $U_{ab}^{(i)}=\sum_jP_{ia,jb}$ form a magic matrix. Since $P=(P_{ia,jb})$ is magic, the elements $U_{ab}^{(i)}=\sum_jP_{ia,jb}$ are self-adjoint, and we have $\sum_bU_{ab}^{(i)}=\sum_{bj}P_{ia,jb}=1$. The fact that each $U_{ab}^{(i)}$ is an idempotent follows from:
\begin{eqnarray*}
&&((U_{ab}^{(i)})^2)_{kc,ld}\\
&=&\frac{1}{N^2M^2}\sum_{mejn}\frac{Q_{ic}Q_{je}}{Q_{ie}Q_{jc}}\cdot\frac{H_{ik}H_{jm}}{H_{im}H_{jk}}\cdot\frac{K_{ac}K_{be}}{K_{ae}K_{bc}}\cdot\frac{Q_{ie}Q_{nd}}{Q_{id}Q_{ne}}\cdot\frac{H_{im}H_{nl}}{H_{il}H_{nm}}\cdot\frac{K_{ae}K_{bd}}{K_{ad}K_{be}}\\
&=&\frac{1}{NM^2}\sum_{ejn}\frac{Q_{ic}Q_{je}Q_{nd}}{Q_{jc}Q_{id}Q_{ne}}\cdot\frac{H_{ik}H_{nl}}{H_{jk}H_{il}}\delta_{jn}\cdot\frac{K_{ac}K_{bd}}{K_{bc}K_{ad}}\\
&=&\frac{1}{NM^2}\sum_{ej}\frac{Q_{ic}Q_{je}Q_{jd}}{Q_{jc}Q_{id}Q_{je}}\cdot\frac{H_{ik}H_{jl}}{H_{jk}H_{il}}\cdot\frac{K_{ac}K_{bd}}{K_{bc}K_{ad}}\\
&=&\frac{1}{NM}\sum_j\frac{Q_{ic}Q_{jd}}{Q_{jc}Q_{id}}\cdot\frac{H_{ik}H_{jl}}{H_{jk}H_{il}}\cdot\frac{K_{ac}K_{bd}}{K_{bc}K_{ad}}\\
&=&(U_{ab}^{(i)})_{kc,ld}
\end{eqnarray*}

Finally, the condition $\sum_aU_{ab}^{(i)}=1$ can be checked as follows:
\begin{eqnarray*}
\sum_aU_{ab}^{(i)}
&=&\frac{1}{N}\left(\sum_j\frac{Q_{ic}Q_{jd}}{Q_{id}Q_{jc}}\cdot\frac{H_{ik}H_{jl}}{H_{il}H_{jk}}\cdot\delta_{cd}\right)_{kc,ld}\\
&=&\frac{1}{N}\left(\sum_j\frac{H_{ik}H_{jl}}{H_{il}H_{jk}}\cdot\delta_{cd}\right)_{kc,ld}\\
&=&1
\end{eqnarray*}

(3) It remains to prove that we have $U_{ab}^{(i)}V_{ij}=V_{ij}U_{ab}^{(i)}=P_{ia,jb}$. First, we have:
\begin{eqnarray*}
(U_{ab}^{(i)}V_{ij})_{kc,ld}
&=&\frac{1}{N^2M}\sum_{mn}\frac{Q_{ic}Q_{nd}}{Q_{id}Q_{nc}}\cdot\frac{H_{ik}H_{nm}}{H_{im}H_{nk}}\cdot\frac{K_{ac}K_{bd}}{K_{ad}K_{bc}}\cdot\frac{H_{im}H_{jl}}{H_{il}H_{jm}}\\
&=&\frac{1}{NM}\sum_n\frac{Q_{ic}Q_{nd}}{Q_{id}Q_{nc}}\cdot\frac{H_{ik}H_{jl}}{H_{nk}H_{il}}\delta_{nj}\cdot\frac{K_{ac}K_{bd}}{K_{ad}K_{bc}}\\
&=&\frac{1}{NM}\cdot\frac{Q_{ic}Q_{jd}}{Q_{id}Q_{jc}}\cdot\frac{H_{ik}H_{jl}}{H_{jk}H_{il}}\cdot\frac{K_{ac}K_{bd}}{K_{ad}K_{bc}}\\
&=&(P_{ia,jb})_{kc,ld}
\end{eqnarray*}

The remaining computation is similar, as follows:
\begin{eqnarray*}
(V_{ij}U_{ab}^{(i)})_{kc,ld}
&=&\frac{1}{N^2M}\sum_{mn}\frac{H_{ik}H_{jm}}{H_{im}H_{jk}}\cdot\frac{Q_{ic}Q_{nd}}{Q_{id}Q_{nc}}\cdot\frac{H_{im}H_{nl}}{H_{il}H_{nm}}\cdot\frac{K_{ac}K_{bd}}{K_{ad}K_{bc}}\\
&=&\frac{1}{NM}\sum_n\frac{Q_{ic}Q_{nd}}{Q_{id}Q_{nc}}\cdot\frac{H_{ik}H_{nl}}{H_{jk}H_{il}}\delta_{jn}\cdot\frac{K_{ac}K_{bd}}{K_{ad}K_{bc}}\\
&=&\frac{1}{NM}\cdot\frac{Q_{ic}Q_{jd}}{Q_{id}Q_{jc}}\cdot\frac{H_{ik}H_{jl}}{H_{jk}H_{il}}\cdot\frac{K_{ac}K_{bd}}{K_{ad}K_{bc}}\\
&=&(P_{ia,jb})_{kc,ld}
\end{eqnarray*}

Thus we have checked all the relations, and we are done.
\end{proof}

In general, the problem of further factorizing the above representation is a quite difficult one, and this even in the case of the usual Hadamard matrices. For a number of results here, which are however quite specialized, we refer to \cite{bbi} and related papers.

\section*{15c. Partial permutations}

Let us discuss now another generalization of the construction $H\to G$, which is independent from the one above. The idea, following \cite{bsk}, will be that of looking at the partial Hadamard matrices (PHM), and their connection with the partial permutations. Let us start with the following standard definition:

\index{partial permutation}

\begin{definition}
A partial permutation of $\{1\,\ldots,N\}$ is a bijection 
$$\sigma:X\simeq Y$$
between two subsets of the index set, as follows:
$$X,Y\subset\{1,\ldots,N\}$$
We denote by $\widetilde{S}_N$ the set formed by such partial permutations. 
\end{definition}

We have $S_N\subset\widetilde{S}_N$, and the embedding $u:S_N\subset M_N(0,1)$ given by the standard permutation matrices can be extended to an embedding $u:\widetilde{S}_N\subset M_N(0,1)$, as follows:
$$u_{ij}(\sigma)=
\begin{cases}
1&{\rm if}\ \sigma(j)=i\\
0&{\rm otherwise}
\end{cases}$$

By looking at the image of this embedding, we see that $\widetilde{S}_N$ is in bijection with the matrices $M\in M_N(0,1)$ having at most one 1 entry on each row and column. In analogy now with Wang's theory in \cite{wa1}, we have the following definition:

\index{submagic matrix}
\index{free partial permutation}

\begin{definition}
A submagic matrix is a matrix $u\in M_N(A)$ whose entries are projections, which are pairwise orthogonal on rows and columns. We let $C(\widetilde{S}_N^+)$ be the universal $C^*$-algebra generated by the entries of a $N\times N$ submagic matrix. 
\end{definition}

Here the fact that the algebra $C(\widetilde{S}_N^+)$ is indeed well-defined is clear. As a first observation, this algebra has a comultiplication, given by the following formula:
$$\Delta(u_{ij})=\sum_ku_{ik}\otimes u_{kj}$$

This algebra has as well a counit, given by the following formula:
$$\varepsilon(u_{ij})=\delta_{ij}$$

Thus $\widetilde{S}_N^+$ is a quantum semigroup, and we have maps as follows, with the bialgebras at left corresponding to the quantum semigroups at right:
$$\begin{matrix}
C(\widetilde{S}_N^+)&\to&C(S_N^+)\\
\\
\downarrow&&\downarrow\\
\\
C(\widetilde{S}_N)&\to&C(S_N)
\end{matrix}
\quad \quad \quad:\quad \quad\quad
\begin{matrix}
\widetilde{S}_N^+&\supset&S_N^+\\
\\
\cup&&\cup\\
\\
\widetilde{S}_N&\supset&S_N
\end{matrix}$$

The relation of all this with the PHM is immediate, appearing as follows:

\index{PHM}
\index{complex PHM}
\index{quantum semigroup}

\begin{theorem}
If $H\in M_{M\times N}(\mathbb T)$ is a PHM, with rows denoted $H_1,\ldots,H_M\in\mathbb T^N$, then the following matrix of rank one projections is submagic:
$$P_{ij}=Proj\left(\frac{H_i}{H_j}\right)$$ 
Thus $H$ produces a representation $\pi_H:C(\widetilde{S}_M^+)\to M_N(\mathbb C)$, given by $u_{ij}\to P_{ij}$, that we can factorize through $C(G)$, with the quantum semigroup $G\subset\widetilde{S}_M^+$ chosen minimal.
\end{theorem}

\begin{proof}
We have indeed the following computation, for the rows:
\begin{eqnarray*}
\Big\langle\frac{H_i}{H_j},\frac{H_i}{H_k}\Big\rangle
&=&\sum_l\frac{H_{il}}{H_{jl}}\cdot\frac{H_{kl}}{H_{il}}\\
&=&\sum_l\frac{H_{kl}}{H_{jl}}\\
&=&<H_k,H_j>\\
&=&\delta_{jk}
\end{eqnarray*}

The verification for the columns is similar, as follows:
\begin{eqnarray*}
\left<\frac{H_i}{H_j},\frac{H_k}{H_j}\right>
&=&\sum_l\frac{H_{il}}{H_{jl}}\cdot\frac{H_{jl}}{H_{kl}}\\
&=&\sum_l\frac{H_{il}}{H_{kl}}\\
&=&N\delta_{ik}
\end{eqnarray*}

Regarding now the last assertion, we can indeed factorize our representation as indicated, with the existence and uniqueness of the bialgebra $C(G)$, with the minimality property as above, being obtained by dividing $C(\widetilde{S}_M^+)$ by a suitable ideal. See \cite{bsk}.
\end{proof}

Summarizing, we have a generalization of the $H\to G$ construction from chapter 14.  The very first problem is that of deciding under which exact assumptions our construction is in fact ``classical''. In order to explain the answer here, we will need:

\index{pre-Latin square}

\begin{definition}
A pre-Latin square is a square matrix 
$$L\in M_M(1,\ldots,N)$$
having the property that its entries are distinct, on each row and each column. 
\end{definition}

Given such a pre-Latin square $L$, to any $x\in\{1,\ldots,N\}$ we can associate the partial permutation $\sigma_x\in\widetilde{S}_M$ given by the following formula:
$$\sigma_x(j)=i\iff L_{ij}=x$$

With this construction in hand, we denote by $G\subset\widetilde{S}_M$ the semigroup generated by these partial permutations $\sigma_1,\ldots,\sigma_N$, and call it semigroup associated to $L$. Also, given an orthogonal basis $\xi=(\xi_1,\ldots,\xi_N)$ of $\mathbb C^N$, we can construct a submagic matrix $P\in M_M(M_N(\mathbb C))$, according to the following formula:
$$P_{ij}=Proj(\xi_{L_{ij}})$$

With these notations, we have the following result, from \cite{bsk}:

\begin{theorem}
If $H\in M_{N\times M}(\mathbb C)$ is a PHM, the following are equivalent:
\begin{enumerate}
\item The semigroup $G\subset\widetilde{S}_M^+$ is classical, i.e. $G\subset\widetilde{S}_M$.

\item The projections $P_{ij}=Proj(H_i/H_j)$ pairwise commute.

\item The vectors $H_i/H_j\in\mathbb T^N$ are pairwise proportional, or orthogonal. 

\item The submagic matrix $P=(P_{ij})$ comes for a pre-Latin square $L$.
\end{enumerate}
In addition, if so is the case, $G$ is the semigroup associated to $L$.
\end{theorem}

\begin{proof}
This is something standard, as follows:

\medskip

$(1)\iff(2)$ is clear from definitions.

\medskip

$(2)\iff(3)$ comes from the fact that two rank 1 projections commute precisely when their images coincide, or are orthogonal.

\medskip

$(3)\iff(4)$ is clear again.

\medskip

As for the last assertion, this is something standard, coming from Gelfand duality, which allows us to compute the Hopf image, in combinatorial terms. See \cite{bsk}. 
\end{proof}

We call ``classical'' the matrices in Theorem 15.23, that we will study now. Let us begin with a study at $M=2$. We make the following convention, where $\tau$ is the transposition, $ij$ is the partial permutation $i\to j$, and $\emptyset$ is the null map:
$$\widetilde{S}_2=\{id,\tau,11,12,21,22,\emptyset\}$$

With this convention, we have the following result:

\begin{proposition}
A partial Hadamard matrix $H\in M_{2\times N}(\mathbb T)$, in dephased form
$$H=\begin{pmatrix}1&\ldots&1\\ \lambda_1&\ldots&\lambda_N\end{pmatrix}$$
is of classical type when one of the following happens:
\begin{enumerate}
\item Either $\lambda_i=\pm w$, for some $w\in\mathbb T$, in which case $G=\{id,\tau\}$.

\item Or $\sum_i\lambda_i^2=0$, in which case $G=\{id,11,12,21,22,\emptyset\}$
\end{enumerate}
\end{proposition}

\begin{proof}
With $1=(1,\ldots,1)$ and $\lambda=(\lambda_1,\ldots,\lambda_N)$, the matrix formed by the vectors $H_i/H_j$ is $(^1_{\bar{\lambda}}{\ }^\lambda_1)$. Since $1\perp\lambda,\bar{\lambda}$ we just have to compare $\lambda,\bar{\lambda}$, and we have two cases:

\medskip

(1) Case $\lambda\sim\bar{\lambda}$. This means that we have $\lambda^2\sim1$, and so $\lambda_i=\pm w$, for some complex number $w\in\mathbb T$. In this case the associated pre-Latin square is $L=(^1_2{\ }^2_1)$, and the partial permutations $\sigma_x$ associated to $L$, as above, are as follows:
$$\sigma_1=id\quad,\quad 
\sigma_2=\tau$$

We obtain from this that we have, as claimed:
$$G=<id,\tau>=\{id,\tau\}$$

(2) Case $\lambda\perp\bar{\lambda}$. This means $\sum_i\lambda_i^2=0$. In this case the associated pre-Latin square is $L=(^1_3{\ }^2_1)$, the associated partial permutations $\sigma_x$ are given by:
$$\sigma_1=id\quad,\quad 
\sigma_2=21\quad,\quad 
\sigma_3=12$$

The semigroup generated by these partial permutations is:
$$G=<id,21,12>=\{id,11,12,21,22,\emptyset\}$$

Thus, we are led to the conclusion in the statement.
\end{proof}

The matrices in (1) are, modulo equivalence, those which are real. As for the matrices in (2), these are parametrized by the solutions $\lambda\in\mathbb T^N$ of the following equations:
$$\sum_i\lambda_i=\sum_i\lambda_i^2=0$$

In general, it is quite unclear on how to deal with these equations. Observe however that, as a basic example here, we have the upper $2\times N$ submatrix of $F_N$, with $N\geq3$. We refer to \cite{bop}, \cite{bsk} and related papers, for more on these questions.

\section*{15d. Fourier matrices}

Let us discuss now in detail the truncated Fourier matrix case. First, we have the following result, that we already know from chapter 14, but that we will present here with a complete proof, as an illustration for Theorem 15.23:

\begin{proposition}
The Fourier matrix, which is as follows, with $w=e^{2\pi i/N}$, 
$$F_N=(w^{ij})$$
is of classical type, and the associated group $G\subset S_N$ is the cyclic group $\mathbb Z_N$.
\end{proposition}

\begin{proof}
Since $H=F_N$ is a square matrix, the associated semigroup $G\subset\widetilde{S}_N^+$ must be a quantum group, $G\subset S_N^+$. We must prove that we have  $G=\mathbb Z_N$. Let us set:
$$\rho=(1,w,w^2,\ldots,w^{N-1})$$

The rows of $H$ are then given by $H_i=\rho^i$, and so we have:
$$\frac{H_i}{H_j}=\rho^{i-j}$$

We conclude that $H$ is indeed of classical type, coming from the Latin square $L_{ij}=j-i$ and from the following orthogonal basis:
$$\xi=(1,\rho^{-1},\rho^{-2},\ldots,\rho^{1-N})$$

We have $G=<\sigma_1,\ldots,\sigma_N>$, where $\sigma_x\in S_N$ is given by:
$$\sigma_x(j)=i\iff L_{ij}=x$$

Now from $L_{ij}=j-i$ we obtain $\sigma_x(j)=j-x$, and so:
$$G=\{\sigma_1,\ldots,\sigma_N\}\simeq\mathbb Z_N$$

Thus, we are led to the conclusion in the statement.
\end{proof}

\index{truncated Fourier matrix}

We will be interested in what follows in the truncated Fourier matrices. Let $F_{M,N}$ be the upper $M\times N$ submatrix of $F_N$, and $G_{M,N}\subset\widetilde{S}_M$ be the associated semigroup. The simplest case is that when $M$ is small, and we have here the following result:

\begin{theorem}
In the $N>2M-2$ regime, $G_{M,N}\subset\widetilde{S}_M$ is formed by the maps 
\vskip-10mm$$\begin{matrix}\\ \\ \\ \sigma=\ \ \\ \end{matrix}
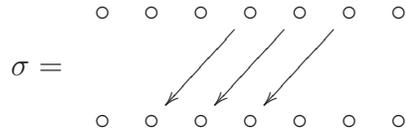
\xymatrix@R=10mm@C=2mm{
\circ&\circ&\circ&\circ\ar[dll]&\circ\ar[dll]&\circ\ar[dll]&\circ\\
\circ&\circ&\circ&\circ&\circ&\circ&\circ}$$
that is, $\sigma:I\simeq J$, $\sigma(j)=j-x$, with $I,J\subset\{1,\ldots,M\}$ intervals, independently of $N$.
\end{theorem}

\begin{proof}
For $\widetilde{H}=F_N$ the associated Latin square is circulant, given by:
$$\widetilde{L}_{ij}=j-i$$

Thus, the pre-Latin square that we are interested in is given by:
$$L=\begin{pmatrix}
0&1&2&\ldots&M-1\\
N-1&0&1&\ldots&M-2\\
N-2&N-1&0&\ldots&M-3\\
\ldots\\
N-M+1&N-M+2&N-M+3&\ldots&0
\end{pmatrix}$$

Observe that, due to our $N>2M-2$ assumption, we have $N-M+1>M-1$, and so the entries above the diagonal are distinct from those below the diagonal. Let us compute now the partial permutations $\sigma_x\in\widetilde{S}_M$ given by:
$$\sigma_x(j)=i\iff L_{ij}=x$$

We have $\sigma_0=id$, and then $\sigma_1,\sigma_2,\ldots,\sigma_{M-1}$ are as follows:
\vskip-7mm
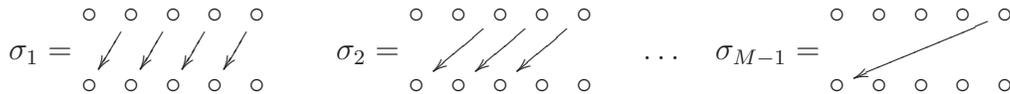
$$\begin{matrix}\\ \\ \sigma_1=\\ \end{matrix}
\xymatrix@R=5mm@C=1mm{
\circ&\circ\ar[dl]&\circ\ar[dl]&\circ\ar[dl]&\circ\ar[dl]\\
\circ&\circ&\circ&\circ&\circ}$$
$$\begin{matrix}\\ \\ \sigma_2=\\ \end{matrix}
\xymatrix@R=5mm@C=1mm{
\circ&\circ&\circ\ar[dll]&\circ\ar[dll]&\circ\ar[dll]\\
\circ&\circ&\circ&\circ&\circ}$$
$$\vdots$$
\vskip-7mm
$$\begin{matrix}\\ \\ \sigma_{M-1}=\\ \end{matrix}
\xymatrix@R=5mm@C=1mm{
\circ&\circ&\circ&\circ&\circ\ar[dllll]\\
\circ&\circ&\circ&\circ&\circ}$$

Observe that we have the following formulae, for these maps:
$$\sigma_2=\sigma_1^2$$
$$\sigma_3=\sigma_1^3$$
$$\vdots$$
$$\sigma_{M-1}=\sigma_1^{M-1}$$

As for the remaining partial permutations, these are given by:
$$\sigma_{N-1}=\sigma_1^{-1}$$
$$\sigma_{N-2}=\sigma_2^{-1}$$
$$\vdots$$
$$\sigma_{N-M+1}=\sigma_{M-1}^{-1}$$

The corresponding diagrams are as follows:
\vskip-7mm$$\begin{matrix}\\ \\ \sigma_{N-1}=\\ \end{matrix}
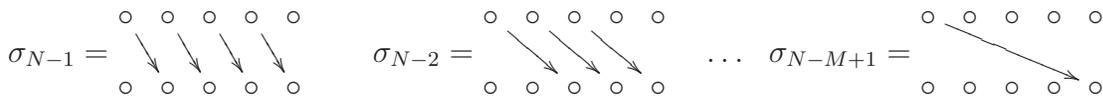
\xymatrix@R=5mm@C=1mm{
\circ\ar[dr]&\circ\ar[dr]&\circ\ar[dr]&\circ\ar[dr]&\circ\\
\circ&\circ&\circ&\circ&\circ}$$
$$\begin{matrix}\\ \\ \sigma_{N-2}=\\ \end{matrix}
\xymatrix@R=5mm@C=1mm{
\circ\ar[drr]&\circ\ar[drr]&\circ\ar[drr]&\circ&\circ\\
\circ&\circ&\circ&\circ&\circ}$$
$$\vdots$$
\vskip-7mm
$$\begin{matrix}\\ \\ \sigma_{N-M+1}=\\ \end{matrix}
\xymatrix@R=5mm@C=1mm{
\circ\ar[drrrr]&\circ&\circ&\circ&\circ\\
\circ&\circ&\circ&\circ&\circ}$$

We conclude that we have the following generation result:
$$G_{M,N}=<\sigma_1>$$

Now if we denote by $G_{M,N}'$ the semigroup in the statement, we have $\sigma_1\in G_{M,N}'$, and so we have an inclusion as follows:
$$G_{M,N}\subset G_{M,N}'$$

The reverse inclusion can be established as follows:

\medskip

(1) Assume first that $\sigma\in G_{M,N}'$, $\sigma:I\simeq J$ has the property $M\in I,J$:
\vskip-10mm$$\begin{matrix}\\ \\ \\ \sigma=\ \ \\ \end{matrix}
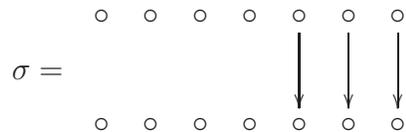
\xymatrix@R=10mm@C=2mm{
\circ&\circ&\circ&\circ&\circ\ar[d]&\circ\ar[d]&\circ\ar[d]\\
\circ&\circ&\circ&\circ&\circ&\circ&\circ}$$

Then we can write $\sigma=\sigma_{N-k}\sigma_k$, with $k=M-|I|$, so we have $\sigma\in G_{M,N}$.

\medskip

(2) Assume now that $\sigma\in G_{M,N}'$, $\sigma:I\simeq J$ has just the property $M\in I$ or $M\in J$:
\vskip-10mm$$\begin{matrix}\\ \\ \\ \sigma'=\ \ \\ \end{matrix}
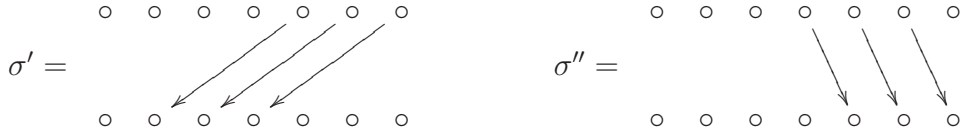
\xymatrix@R=10mm@C=2mm{
\circ&\circ&\circ&\circ&\circ\ar[dlll]&\circ\ar[dlll]&\circ\ar[dlll]\\
\circ&\circ&\circ&\circ&\circ&\circ&\circ}$$
$$\begin{matrix}\\ \\ \\ \sigma''=\ \ \\ \end{matrix}\xymatrix@R=10mm@C=2mm{
\circ&\circ&\circ&\circ\ar[dr]&\circ\ar[dr]&\circ\ar[dr]&\circ\\
\circ&\circ&\circ&\circ&\circ&\circ&\circ}$$

In this case we have as well $\sigma\in G_{M,N}$, because $\sigma$ appears from one of the maps in (1) by adding a ``slope'', which can be obtained by composing with a suitable map $\sigma_k$.

\medskip

(3) Assume now that $\sigma\in G_{M,N}'$, $\sigma:I\simeq J$ is arbitrary:
\vskip-10mm$$\begin{matrix}\\ \\ \\ \sigma=\ \ \\ \end{matrix}
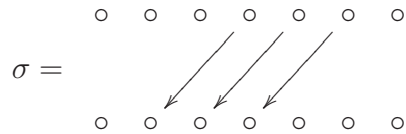
\xymatrix@R=10mm@C=2mm{
\circ&\circ&\circ&\circ\ar[dll]&\circ\ar[dll]&\circ\ar[dll]&\circ\\
\circ&\circ&\circ&\circ&\circ&\circ&\circ}$$

Then we can write $\sigma=\sigma'\sigma''$ with $\sigma':L\simeq J$, $\sigma'':I\simeq L$, where $L$ is an interval satisfying $|L|=|I|=|J|$ and $M\in L$, and since $\sigma',\sigma''\in G_{M,N}$ by (2), we are done.
\end{proof}

Summarizing, we have so far complete results at $N=M$, and at $N>2M-2$. In the remaining regime, $M<N\leq2M-2$, the semigroup $G_{M,N}\subset\widetilde{S}_M$ looks quite hard to compute, and for the moment there are only partial results regarding it. For a partial permutation $\sigma:I\simeq J$ with $|I|=|J|=k$, set $\kappa(\sigma)=k$. We have:

\begin{theorem}
The following semigroup components, with $k>2M-N$,
$$G_{M,N}^{(k)}=\left\{\sigma\in G_{M,N}\Big|\kappa(\sigma)=k\right\}$$ 
are in the $M<N\leq2M-2$ regime the same as those in the $N>2M-2$ regime.
\end{theorem}

\begin{proof}
In the $M<N\leq2M-2$ regime the pre-Latin square that we are interested in has as usual 0 on the diagonal, and then takes its entries from the following set, in a uniform way from each of the 3 components:
$$S=\{1,\ldots,N-M\}\cup\{N-M+1,\ldots,M-1\}\cup\{M,\ldots,N-1\}$$

Here is an illustrating example, at $M=6,N=8$: 
$$L=\begin{pmatrix}
{\bf 0}&1&2&{\bf 3}&{\bf 4}&{\bf 5}\\
7&{\bf 0}&1&2&{\bf 3}&{\bf 4}\\
6&7&{\bf 0}&1&2&{\bf 3}\\
{\bf 5}&6&7&{\bf 0}&1&2\\
{\bf 4}&{\bf 5}&6&7&{\bf 0}&1\\
{\bf 3}&{\bf 4}&{\bf 5}&6&7&{\bf 0}
\end{pmatrix}$$

The point now is that $\sigma_1,\ldots,\sigma_{N-M}$ are given by the same formulae as those in the proof of Theorem 15.26, then $\sigma_{N-M+1},\ldots,\sigma_{M-1}$ all satisfy $\kappa(\sigma)=2M-N$, and finally $\sigma_M,\ldots,\sigma_{N-1}$ are once again given by the formulae in the proof of Theorem 15.26. Now since we have $\kappa(\sigma\rho)\leq\min(\kappa(\sigma),\kappa(\rho))$, adding the maps $\sigma_{N-M+1},\ldots,\sigma_{M-1}$ to the semigroup $G_{M,N}\subset\widetilde{S}_M$ computed in the proof of Theorem 15.26 won't change the $G_{M,N}^{(k)}$ components of this semigroup at $k>2M-N$, and this gives the result.
\end{proof}

\section*{15e. Exercises} 

We have seen in this chapter two recent generalizations of the construction $H\to G$ from chapter 14, and going beyond the results presented here, even with some simple exercises, is no easy task. As a first exercise, however, we have:

\begin{exercise}
Write down a complete, simplified proof for the factorization
$$\xymatrix{C(S_{NM}^+)\ar[rr]^{\pi_L}\ar[rd]&&M_{NM}(\mathbb C)\\&C(S_M^+\wr_*G_H)\ar[ur]&}$$
found above, for $L=H\otimes_QK$, in the scalar matrix case.
\end{exercise}

To be more precise, the problem is that of reviewing the proof of the above factorization, checking what simplifies in the scalar matrix case, and writing this down.

\begin{exercise}
Prove that the number of partial permutations is given by
$$|\widetilde{S}_N|=\sum_{k=0}^Nk!\binom{N}{k}^2$$
that is, $1,2,7,34,209,\ldots\,$, and that we have the estimate
$$|\widetilde{S}_N|\simeq N!\sqrt{\frac{\exp(4\sqrt{N}-1)}{4\pi\sqrt{N}}}$$
in the $N\to\infty$ limit.
\end{exercise}

Here the first assertion is easy, and the second one is difficult.

\begin{exercise}
Prove that we have an isomorphism
$$C(\widetilde{S}_2^+)\simeq\left\{(x,y)\in C^*(D_\infty)\oplus C^*(D_\infty)\Big|\varepsilon(x)=\varepsilon(y)\right\}$$
where $\varepsilon:C^*(D_\infty)\to\mathbb C1$ the usual counit map.
\end{exercise}

As a first step here, we would need a structure result for the $2\times2$ submagic matrices.

\begin{exercise}
Develop a theory of partial Hadamard matrices with noncommutative entries, and of the associated quantum permutation semigroups.
\end{exercise}

The statement here is of course quite loose, as is always the case with research-grade exercises, and anything is welcome, the more the better.

\chapter{Fourier models}

\section*{16a. Deformations}

In this chapter we go back to the usual complex Hadamard matrices, $H\in M_N(\mathbb C)$. We know that associated to any such matrix is a certain quantum permutation group $G\subset S_N^+$, which describes the symmetries of the matrix. The main example for this construction $H\to G$ is, as expected, $F_N\to\mathbb Z_N$, and more generally, $F_G\to G$, for any finite abelian group $G$. There are of course many things that can be said about the correspondence $H\to G$, but the main question remains the explicit computation of $G$, in terms of $H$. Here we discuss this question for the deformed Fourier matrices.

\bigskip

Contrary to many other things discussed in this book, this is something that has been intensively studied, and not that the known results are fully satisfactory, but at least they lie at the level of what the experts can do. The story of the subject is as follows:

\bigskip

(1) The origins of the question go back to some discussions, and even papers, written by Bichon, Nicoara, Schlenker and myself in the mid 00s, containing a few mistakes, which ruined the thing, initially. Be said in passing, regarding wrong papers, never ever do that, if possible, and for good reason. Not with respect to mathematics and the community, who are legendary slow anyway in digesting new things, but with respect to yourself, and your business. Believe me, with any wrong paper, you dig your own grave.

\bigskip

(2) Towards the end of the 00s, some computations by Nicoara and his students on one hand, and some computations of Burstein, a student of Jones, on the other \cite{bur}, done in the commuting square and subfactor context, showed that the problem for the deformed Fourier matrices is very interesting, and far more complicated than previously thought. In the context of the correspondence $H\to G$, as above, the study was done short after, in a joint paper by Bichon and myself \cite{bbi}, that we will explain in what follows.

\bigskip

(3) Finally, and as a third piece of the story, the paper \cite{bbi}, which contains several exciting things, had several follow-ups, both by Bichon and by myself, which are extremely technical, and barely readable, and that you will certainly be able to find on the internet, if interested, just by following citations, as usual. These papers are, needless to say, correct, but really tough, and the problem for younger generations is that of going beyond that. In my opinion, and Bichon's too, this is certainly possible, and very interesting.

\bigskip

Getting to work now, following \cite{bbi}, we would like to discuss the computation of the quantum groups associated to the Di\c t\u a deformations of the tensor products of Fourier matrices. Let us begin by recalling the construction of the Fourier matrix models:

\index{Fourier model}

\begin{definition}
Associated to a finite abelian group $G$ is the matrix model 
$$\pi:C(G)\to M_G(\mathbb C)$$
coming from the following magic matrix, 
$$(U_{ij})_{kl}=\frac{1}{N}F_{i-j,k-l}$$
where $F=F_G$ is the Fourier matrix of $G$.
\end{definition}

Let us recall as well the construction of the deformed Fourier models:

\index{deformed Fourier matrix}

\begin{definition}
Given two finite abelian groups $G,H$, we consider the corresponding deformed Fourier matrix, given by the formula
$$(F_G\otimes_Q F_H)_{ia,jb}=Q_{ib}(F_G)_{ij}(F_H)_{ab}$$
and we factorize the associated representation $\pi_Q$ of the algebra $C(S_{G\times H}^+)$,
$$\xymatrix@R=40pt@C=40pt
{C(S_{G\times H}^+)\ar[rr]^{\pi_Q}\ar[rd]&&M_{G\times H}(\mathbb C)\\&C(G_Q)\ar[ur]_\pi&}$$
with $C(G_Q)$ being the Hopf image of this representation $\pi_Q$.
\end{definition}

Explicitely computing the above quantum permutation group $G_Q\subset S_{G\times H}^+$, as function of the parameter matrix $Q\in M_{G\times H}(\mathbb T)$, will be our main purpose, in what follows. In order to do so, we will need the following elementary result:

\begin{proposition}
If $G$ is a finite abelian group then
$$C(G)=C(S_G^+)\Big/\left<u_{ij}=u_{kl}\Big|\forall i-j=k-l\right>$$
with all the indices taken inside $G$.
\end{proposition}

\begin{proof}
As a first observation, the quotient algebra in the statement is commutative, because we have the following relations:
$$u_{ij}u_{kl}=u_{ij}u_{i,l-k+i}=\delta_{j,l-k+i}u_{ij}$$
$$u_{kl}u_{ij}=u_{i,l-k+i}u_{ij}=\delta_{j,l-k+i}u_{ij}$$

Thus if we denote the algebra in the statement by $C(H)$, we have $H\subset S_G$. Now since $u_{ij}(\sigma)=\delta_{i\sigma(j)}$ for any $\sigma\in H$, we obtain:
$$i-j=k-l\implies(\sigma(j)=i\iff\sigma(l)=k)$$

But this condition tells us precisely that $\sigma(i)-i$ must be independent on $i$, and so, for some $g\in G$, we have $\sigma(i)=i+g$. Thus we have $\sigma\in G$, as desired.
\end{proof}

In order to factorize the representation in Definition 16.2, we will need:

\index{free wreath product}
\index{wreath product}

\begin{definition}
Gives two Hopf algebra quotients, as follows,
$$C(S_M^+)\to A\quad,\quad 
C(S_N^+)\to B$$
with fundamental corepresentations denoted $u,v$, we let
$$A*_wB=A^{*N}*B/<[u_{ab}^{(i)},v_{ij}]=0>$$
with the Hopf algebra structure making $w_{ia,jb}=u_{ab}^{(i)}v_{ij}$ a corepresentation.
\end{definition}

The fact that we have indeed a Hopf algebra follows from the fact that $w$ is magic. In terms of quantum groups, let us write:
$$A=C(G)\quad,\quad 
B=C(H)$$

We can write then the Hopf algebra constructed above as follows: 
$$A*_wB=C(G\wr_*H)$$

In other words, we make the following convention:
$$C(G)*_wC(H)=C(G\wr_*H)$$

The $\wr_*$ operation is then the free analogue of $\wr$, the usual wreath product. For details regarding this construction, we refer to \cite{bbi}, or to the book \cite{ba1}. Now with this notion in hand, we can factorize representation $\pi_Q$ in Definition 16.2, as follows:

\begin{theorem}
We have a factorization as follows,
$$\xymatrix@R=40pt@C=40pt
{C(S_{G\times H}^+)\ar[rr]^{\pi_Q}\ar[rd]&&M_{G\times H}(\mathbb C)\\&C(H\wr_*G)\ar[ur]_\pi&}$$
given on the standard generators by the formulae
$$U_{ab}^{(i)}=\sum_jW_{ia,jb}\quad,\quad 
V_{ij}=\sum_aW_{ia,jb}$$
independently of $b$, where $W$ is the magic matrix producing $\pi_Q$.
\end{theorem}

\begin{proof}
With $K=F_G,L=F_H$ and $M=|G|,N=|H|$, the formula of the magic matrix $W\in M_{G\times H}(M_{G\times H}(\mathbb C))$ associated to $H=K\otimes_QL$ is as follows:
\begin{eqnarray*}
(W_{ia,jb})_{kc,ld}
&=&\frac{1}{MN}\cdot\frac{Q_{ic}Q_{jd}}{Q_{id}Q_{jc}}\cdot\frac{K_{ik}K_{jl}}{K_{il}K_{jk}}\cdot\frac{L_{ac}L_{bd}}{L_{ad}L_{bc}}\\
&=&\frac{1}{MN}\cdot\frac{Q_{ic}Q_{jd}}{Q_{id}Q_{jc}}\cdot K_{i-j,k-l}L_{a-b,c-d}
\end{eqnarray*}

Our claim now is that the representation $\pi_Q$ constructed in Definition 16.2 can be factorized in three steps, up to the factorization in the statement, as follows:
$$\xymatrix@R=70pt@C=60pt
{C(S_{G\times H}^+)\ar[rr]^{\pi_Q}\ar[d]&&M_{G\times H}(\mathbb C)\\
C(S_H^+\wr_*S_G^+)\ar[r]\ar@{.>}[rru]&C(S_H^+\wr_*G)\ar[r]\ar@{.>}[ur]&C(H\wr_*G)\ar@{.>}[u]}$$

Indeed, these factorizations can be constructed as follows:

\medskip

(1) The construction of the map on the left is standard, by checking the relations for the free wreath product, and this produces the first factorization. 

\medskip

(2) Regarding the second factorization, the one in the middle, this comes from the fact that since the elements $V_{ij}$ depend on $i-j$, they satisfy the defining relations for the quotient algebra $C(S_G^+)\to C(G)$, coming from Proposition 16.3. 

\medskip

(3) Finally, regarding the third factorization, the one on the right, observe that the above matrix $W_{ia,jb}$ depends only on $i,j$ and on $a-b$. By summing over $j$ we obtain that the elements $U_{ab}^{(i)}$ depend only on $a-b$, and we are done.
\end{proof}

Summarizing, we already have some advances on our problem, the quantum group that we want to compute appearing as a subgroup of a certain free wreath product. In order to further factorize the above representation, we use:

\index{crossed product}
\index{crossed coproduct}

\begin{definition}
If $H\curvearrowright\Gamma$ is a finite group acting by automorphisms on a discrete group, the corresponding crossed coproduct Hopf algebra is
$$C^*(\Gamma)\rtimes C(H)=C^*(\Gamma)\otimes C(H)$$
with comultiplication given by the following formula,
$$\Delta(r\otimes\delta_k)=\sum_{h\in H}(r\otimes\delta_h)\otimes(h^{-1}\cdot r\otimes\delta_{h^{-1}k})$$
for $r\in\Gamma$ and $k\in H$. The corresponding quantum group is denoted $\widehat{\Gamma}\rtimes H$.
\end{definition}

Observe that $C(H)$ is a subcoalgebra, and that $C^*(\Gamma)$ is not a subcoalgebra. Now back to the factorization in Theorem 16.5, the point is that we have:

\begin{proposition}
With $L=F_H,N=|H|$ we have an isomorphism
$$C(H\wr_*G)\simeq C^*(H)^{*G}\rtimes C(G)$$
given by $v_{ij}\to1\otimes v_{ij}$ and by
$$u_{ab}^{(i)}=\frac{1}{N}\sum_cL_{b-a,c}c^{(i)}\otimes 1$$
on the standard generators.
\end{proposition}

\begin{proof}
We know that the algebra $C(H\wr_*G)$, constructed according to our above conventions, is the quotient of $C(H)^{*G}*C(G)$ by the following relations:
$$[u_{ab}^{(i)},v_{ij}]=0$$

Now since the variable $v_{ij}$ depends only on $j-i$, we obtain:
$$[u_{ab}^{(i)},v_{kl}]
=[u_{ab}^{(i)},v_{i,l-k+i}]
=0$$

Thus, we are in a usual tensor product situation, and we have:
$$C(H\wr_*G)=C(H)^{*G}\otimes C(G)$$

Consider now the Fourier transform over $H$, which is a map as follows:
$$\Phi:C(H)\to C^*(H)$$

We can compose the above identification with the following map:
$$\Psi=\Phi^{*G}\otimes id$$

Thus, we obtain an isomorphism as in the statement. Now observe that we have:
$$\Phi(u_{ab})=\frac{1}{N}\sum_cL_{b-a,c}c$$

Thus the formula for the image of $u_{ab}^{(i)}$ is indeed the one in the statement.
\end{proof}

Here is now our key result, which will lead to further factorizations:

\begin{proposition}
With $c^{(i)}=\sum_aL_{ac}u_{a0}^{(i)}$ and $\varepsilon_{ke}=\sum_iK_{ik}e_{ie}$ we have:
$$\pi(c^{(i)})(\varepsilon_{ke})=\frac{Q_{i,e-c}Q_{i-k,e}}{Q_{ie}Q_{i-k,e-c}}\varepsilon_{k,e-c}$$
In particular if $c_1+\ldots+c_s=0$ then the matrix
$$\pi(c_1^{(i_1)}\ldots c_s^{(i_s)})$$
is diagonal, for any choice of the indices $i_1,\ldots,i_s$.
\end{proposition}

\begin{proof}
With $c^{(i)}$ as in the statement, we have the following formula:
\begin{eqnarray*}
\pi(c^{(i)})
&=&\sum_aL_{ac}\pi(u_{a0}^{(i)})\\
&=&\sum_{aj}L_{ac}W_{ia,j0}
\end{eqnarray*}

On the other hand, in terms of the basis in the statement, we have:
$$W_{ia,jb}(\varepsilon_{ke})=\frac{1}{N}\delta_{i-j,k}\sum_d\frac{Q_{id}Q_{je}}{Q_{ie}Q_{jd}}L_{a-b,d-e}\varepsilon_{kd}$$

We therefore obtain, as desired:
\begin{eqnarray*}
\pi(c^{(i)})(\varepsilon_{ke})
&=&\frac{1}{N}\sum_{ad}L_{ac}\frac{Q_{id}Q_{i-k,e}}{Q_{ie}Q_{i-k,d}}L_{a,d-e}\varepsilon_{kd}\\
&=&\frac{1}{N}\sum_d\frac{Q_{id}Q_{i-k,e}}{Q_{ie}Q_{i-k,d}}\varepsilon_{kd}\sum_aL_{a,d-e+c}\\
&=&\sum_d\frac{Q_{id}Q_{i-k,e}}{Q_{ie}Q_{i-k,d}}\varepsilon_{kd}\delta_{d,e-c}\\
&=&\frac{Q_{i,e-c}Q_{i-k,e}}{Q_{ie}Q_{i-k,e-c}}\varepsilon_{k,e-c}
\end{eqnarray*}

Regarding now the last assertion, this follows from the fact that each matrix of type $\pi(c_r^{(i_r)})$ acts on the standard basis elements $\varepsilon_{ke}$ by preserving the left index $k$, and by rotating by $c_r$ the right index $e$. Thus when we assume $c_1+\ldots+c_s=0$ all these rotations compose up to the identity, and we obtain indeed a diagonal matrix.
\end{proof}

We have now all needed ingredients for refining Theorem 16.5, as follows:

\begin{theorem}
We have a factorization as follows,
$$\xymatrix@R=40pt@C=30pt
{C(S_{G\times H}^+)\ar[rr]^{\pi_Q}\ar[rd]&&M_{G\times H}(\mathbb C)\\&C^*(\Gamma_{G,H})\rtimes C(G)\ar[ur]_\rho&}$$
where the group on the bottom is given by
$$\Gamma_{G,H}=H^{*G}\Big/\left<[c_1^{(i_1)}\ldots c_s^{(i_s)},d_1^{(j_1)}\ldots d_s^{(j_s)}]=1\Big|\sum_rc_r=\sum_rd_r=0\right>$$
with the above conventions and notations.
\end{theorem}

\begin{proof}
Assume that we have a representation, as follows:
$$\pi:C^*(\Gamma)\rtimes C(G)\to M_L(\mathbb C)$$

Let $\Lambda$ be a $G$-stable normal subgroup of $\Gamma$, so that $G$ acts on $\Gamma/\Lambda$, and we can form the product $C^*(\Gamma/\Lambda)\rtimes C(G)$, and assume that $\pi$ is trivial on $\Lambda$. Then $\pi$ factorizes as:
$$\xymatrix@R=40pt@C=30pt
{C^*(\Gamma)\rtimes C(G)\ar[rr]^\pi\ar[rd]&&M_L(\mathbb C)\\&C^*(\Gamma/\Lambda)\rtimes C(G)\ar[ur]_\rho}$$

With $\Gamma=H^{*G}$, and by using the above results, this gives the result.
\end{proof}

In what follows we will restrict attention to the case where the parameter matrix $Q$ is generic, and we prove that, in this case, the representation in Theorem 16.9 is the minimal one. Our starting point is the group $\Gamma_{G,H}$ found above. Let us formulate:

\begin{definition}
Associated to two finite abelian groups $G,H$ is the discrete group
$$\Gamma_{G,H}=H^{*G}\Big/\left<[c_1^{(i_1)}\ldots c_s^{(i_s)},d_1^{(j_1)}\ldots d_s^{(j_s)}]=1\Big|\sum_rc_r=\sum_rd_r=0\right>$$
where the superscripts refer to the $G$ copies of $H$, inside the free product.
\end{definition}

We will need a more convenient description of this group. The idea here is that the above commutation relations can be realized inside a suitable semidirect product. Given a group acting on another group, $H\curvearrowright G$, we denote as usual by $G\rtimes H$ the semidirect product of $G$ by $H$, which is the set $G\times H$, with multiplication as follows:
$$(a,s)(b,t)=(as(b),st)$$

Now given a group $G$, and a finite abelian group $H$, we can make $H$ act on $G^H$, in the obvious way, and then form the following crossed product:
$$K=G^H\rtimes H$$

Since the elements of type $(g,\ldots,g)$ are invariant under the action of $H$, we can form as well the following crossed product:
$$K'=(G^H/G)\rtimes H$$

We can identify $G^H/G\simeq G^{|H|-1}$ via the following map:
$$(1,g_1,\ldots,g_{|H|-1})\to(g_1,\ldots,g_{|H|-1})$$

Thus, we obtain a crossed product $G^{|H|-1}\rtimes H$. With these notations, we have the following result, regarding the group from Definition 16.10:

\begin{proposition}
The group $\Gamma_{G,H}$ has the following properties:
\begin{enumerate}
\item We have an isomorphism as follows:
$$\Gamma_{G,H}\simeq\mathbb Z^{(|G|-1)(|H|-1)}\rtimes H$$

\item We have as well an isomorphism as follows,
$$\Gamma_{G,H}\subset\mathbb Z^{(|G|-1)|H|}\rtimes H$$
given on the standard generators by the formulae
$$c^{(0)}\to(0,c)\quad,\quad 
c^{(i)}\to(b_{i0}-b_{ic},c)$$
where $b_{ic}$ are the standard generators of $\mathbb Z^{(|G|-1)|H|}$.
\end{enumerate}
\end{proposition}

\begin{proof}
We prove these assertions at the same time. We must prove that we have group morphisms, given by the formulae in the statement, as follows:
\begin{eqnarray*}
\Gamma_{G,H}
&\simeq&\mathbb Z^{(|G|-1)(|H|-1)}\rtimes H\\
&\subset&\mathbb Z^{(|G|-1)|H|}\rtimes H
\end{eqnarray*}

Our first claim is that the formula in (2) defines a morphism as follows:
$$\Gamma_{G,H}\to\mathbb Z^{(|G|-1)|H|}\rtimes H$$

Indeed, we know that the elements $(0,c)$ produce a copy of $H$. Also, we have a group embedding as follows:
$$H\subset\mathbb Z^{|H|}\rtimes H\quad,\quad 
c\to(b_0-b_c,c)$$

Thus the elements $C^{(i)}=(b_{i0}-b_{ic},c)$ produce a copy of $H$, for any $i\neq 0$. In order to check now the commutation relations, observe that we have:
$$C_1^{(i_1)}\ldots C_s^{(i_s)}
=\left(b_{i_10}-b_{i_1c_1}+b_{i_2c_1}-b_{i_2,c_1+c_2}+\ldots+b_{i_s,c_1+\ldots+c_{s-1}}-b_{i_s,c_1+\ldots+c_s},\sum_rc_r\right)$$

Thus $\sum_rc_r=0$ implies the following condition:
$$C_1^{(i_1)}\ldots C_s^{(i_s)}\in\mathbb Z^{(|G|-1)|H|}$$

Since we are now inside an abelian group, we have the commutation relations, and our claim is proved. By using the general crossed product considerations before the statement, it is routine to construct an embedding as follows:
$$\mathbb Z^{(|G|-1)(|H|-1)}\rtimes H\subset \mathbb Z^{(|G|-1)|H|}\rtimes H$$

To be more precise, we would like this embedding to be such that we have group morphisms whose composition is the group morphism just constructed, as follows:
\begin{eqnarray*}
\Gamma_{G,H}
&\to&\mathbb Z^{(|G|-1)(|H|-1)}\rtimes H\\
&\subset&\mathbb Z^{(|G|-1)|H|}\rtimes H
\end{eqnarray*}

It remains to prove that the map on the left is injective. For this purpose, consider the following morphism:
$$\Gamma_{G,H}\to H\quad,\quad 
c^{(i)}\to c$$

The kernel $T$ of this morphism is formed by the elements of type $c_1^{(i_1)} \ldots c_s^{(i_s)}$, with $\sum_rc_r=0$. We therefore obtain an exact sequence, as follows:
$$1\to T\to\Gamma_{G,H}\to H\to1$$

This sequence splits by $c\to c^{(0)}$, so we have:
$$\Gamma_{G,H}\simeq T\rtimes H$$

Now by the definition of $\Gamma_{G,H}$, the subgroup $T$ constructed above is abelian, and is moreover generated by the following elements:
$$(-c)^{(0)}c^{(i)}\quad,\quad c\neq0$$

Finally, the fact that $T$ is freely generated by these elements follows from the computation in the proof of Proposition 16.13 below.
\end{proof}

\section*{16b. Generic parameters}

As already mentioned, we will be interested in what follows in the case where the deformation matrix $Q$ is generic. Our genericity assumptions are as follows:

\index{root independence}
\index{generic deformation}

\begin{definition}
We use the following notions:
\begin{enumerate}
\item We call $p_1,\ldots,p_m\in\mathbb T$ root independent if for any $r_1,\ldots, r_m\in\mathbb Z$ we have:
$$p_1^{r_1}\ldots p_m^{r_m}=1\implies r_1=\ldots=r_m=0$$

\item A matrix $Q\in M_{G\times H}(\mathbb T)$, taken to be dephased,
$$Q_{0c}=Q_{i0}=1$$ 
is called generic if the elements $Q_{ic}$, with $i,c\neq0$, are root independent.
\end{enumerate}
\end{definition}

In what follows we will do the computation for such matrices. Our main result will show that the associated quantum group does not depend in fact of the matrix. In order to do the computation, we will need the following technical result:

\begin{proposition}
Assume that $Q\in M_{G\times H}(\mathbb T)$ is generic, and set:
$$\theta_{ic}^{ke}=\frac{Q_{i,e-c}Q_{i-k,e}}{Q_{ie}Q_{i-k,e-c}}$$
For every $k \in G$, we have a representation $\pi^k:\Gamma_{G,H}\to U_{|H|}$ given by:
$$\pi^k(c^{(i)})\epsilon_e=\theta_{ic}^{ke}\epsilon_{e-c}$$
The family of representations $(\pi^k)_{k \in G}$ is projectively faithful, in the sense that if for some $t \in \Gamma_{G,H}$ we have that $\pi^k(t)$ is a scalar matrix for any $k$, then $t=1$.
\end{proposition}

\begin{proof}
The representations $\pi^k$ arise as above. With $\Gamma_{G,H}=T\rtimes H$, as in the proof of Proposition 16.11, we see that for $t\in\Gamma_{G,H}$ such that $\pi^k(t)$ is a scalar matrix for any $k$, then $t\in T$, since the elements of $T$ are the only ones having their image by $\pi^k$ formed by diagonal matrices. Now write $t$ as follows, with the generators of $T$ being as in the proof of Proposition 16.11, and with $R_{ic}\in\mathbb Z$ being certain integers:
$$t=\prod_{i \not=0, c\not=0} ((-c)^{(0)}(c)^{(i)})^{R_{ic}}$$

Consider now the following quantities:
\begin{eqnarray*}
A(k,e)
&=&\prod_{i\neq0}\prod_{c\neq0}(\theta_{ic}^{ke}(\theta_{0c}^{ke})^{^{-1}})^{R_{ic}}\\
&=&\prod_{i\neq0}\prod_{c\neq0} (\theta_{ic}^{ke})^{R_{ic}}(\theta_{0c}^{ke})^{-R_{ic}}\\
&=&\prod_{i\neq0}\prod_{c\neq0}(\theta_{ic}^{ke})^{R_{ic}}
\cdot\prod_{c\neq0}(\theta_{0c}^{ke})^{-\sum_{i\neq0}R_{ic}}\\
&=&\prod_{j\neq0}\prod_{c\neq0} (\theta_{jc}^{ke})^{R_{jc}}
\cdot\prod_{c\neq0}\prod_{j\neq0}(\theta_{jc}^{ke})^{\sum_{i\neq0}R_{ic}}\\
&=&\prod_{j\neq0}\prod_{c\neq0}(\theta_{jc}^{ke})^{R_{jc}+\sum_{i\neq0}R_{ic}}
\end{eqnarray*}

We have then the following formula, valid for any $k,e$:
$$\pi^k(t)(\epsilon_e)= A(k,e)\epsilon_e$$

Our assumption is that for any $k$, and for any $e,f$, we have: 
$$A(k,e)=A(k,f)$$

By using now the root independence of the elements $Q_{ic}$, with $i,c\neq0$, we see that this implies $R_{ic}=0$ for any $i,c$, and this proves our assertion.  
\end{proof}

We will need as well the following technical result:

\begin{proposition}
Consider a surjective Hopf algebra map
$$\pi:C^*(\Gamma)\rtimes C(H)\to L$$
such that $\pi_{|C(H)}$ is injective, and such that for $r\in\Gamma$ and $f\in C(H)$, we have:
$$\pi(r\otimes 1)=\pi(1\otimes f)\implies r=1$$ 
Then $\pi$ is an isomorphism.
\end{proposition}

\begin{proof} 
We use here various Hopf algebra tools. Consider the following algebra:
$$A=C^*(\Gamma)\rtimes C(H)$$

In order to prove the result, we start with the following standard Hopf algebra exact sequence, where $i(f)=1\otimes f$, and where $p=\varepsilon\otimes 1$:
$$\mathbb C\to C(H)\overset{i}\to A \overset{p}\to C^*(\Gamma)\to
\mathbb C$$ 

Since $\pi\circ i$ is injective, and the Hopf subalgebra $\pi\circ i(C(H))$ is central in $L$, we can form the following quotient Hopf algebra:
$$\overline{L} = L/ (\pi\circ i(C(H))^+L$$

We obtain in this way another exact sequence, as follows:
$$\mathbb C\longrightarrow C(H)\overset{\pi\circ i}\longrightarrow L\overset{q}\longrightarrow\overline{L}\longrightarrow\mathbb C$$ 

Note that this sequence is indeed exact, e.g. by centrality. Thus, we get the following diagram with exact rows, with the Hopf algebra map on the right being surjective:
$$\xymatrix@R=50pt@C=50pt
{\mathbb C\ar[r]&C(H)\ar@2@{-}[d]\ar[r]^i&A\ar[d]^\pi\ar[r]^p&C^*(\Gamma)\ar[r]\ar[d]&\mathbb C\\
\mathbb C\ar[r]&C(H)\ar[r]^{\pi\circ i}&L\ar[r]^q&\overline{L}\ar[r]&\mathbb C}$$

Since a quotient of a group algebra is still a group algebra, we get a commutative diagram with exact rows as follows:
$$\xymatrix@R=50pt@C=50pt
{\mathbb C\ar[r]&C(H)\ar@2@{-}[d]\ar[r]^i&A\ar[d]^\pi\ar[r]^p&C^*(\Gamma)\ar[r]\ar[d]&\mathbb C\\
\mathbb C\ar[r]&C(H)\ar[r]^{\pi\circ i}&L\ar[r]^{q'}&C^*(\overline{\Gamma})\ar[r]&\mathbb C}$$

Here the map on the right is induced by a surjective group morphism, as follows:
$$u:\Gamma\to\overline{\Gamma}\quad,\quad 
g\to\overline{g}$$

By the five lemma, which is something very classical in algebra, we just have to show that $u$ is injective. So, let $g \in \Gamma$ be such that $u(g)=1$. We have then:
$$q' \pi(g \otimes 1)
=u p(g\otimes 1)
=u(g)
=\overline{g}
=1$$

For $g\in\Gamma$, let us set:
$$_gA=\left\{a \in A \ \Big| \ p(a_1) \otimes a_2= g \otimes a\right\}$$
$$_{\overline{g}}L= \left\{l \in L \ \Big| \ q'(l_1) \otimes l_2= \overline{g} \otimes l\right\}$$

The commutativity of the square on the right ensures that we have:
$$\pi(_gA) \subset {_{\overline{g}}L}$$

Then with the previous $g$, we have, by exactness of the sequence:
$$\pi(g \otimes 1) \in {_{\overline{1}}L} = \pi i (C(H))$$

Thus, for some $f \in C(H)$, we must have:
$$\pi(g \otimes 1)= \pi(1 \otimes f)$$

We conclude by our assumption that $g=1$.
\end{proof}

We have now all the needed ingredients for proving a main result, as follows:

\index{generic deformation}

\begin{theorem}
When $Q$ is generic, the minimal factorization for $\pi_Q$ is 
$$\xymatrix@R=45pt@C=40pt
{C(S_{G\times H}^+)\ar[rr]^{\pi_Q}\ar[rd]&&M_{G\times H}(\mathbb C)\\&C^*(\Gamma_{G,H})\rtimes C(G)\ar[ur]_\pi&}$$
where on the bottom
$$\Gamma_{G,H}\simeq\mathbb Z^{(|G|-1)(|H|-1)}\rtimes H$$
is the discrete group constructed above.
\end{theorem}

\begin{proof}
We want to apply Proposition 16.13 to the following morphism, arising from the factorization in Theorem 16.9, where $L$ denotes the Hopf image of $\pi_Q$:
$$\theta : C^*(\Gamma_{G,H})\rtimes C(G)\to L$$

To be more precise, this morphism produces the following commutative diagram:
$$\xymatrix@R=40pt@C=40pt
{C(S_{G\times H}^+) \ar[rr]^{\pi_Q} \ar[dr]_{} \ar@/_/[ddr]_{}& & M_{G\times H}(\mathbb C) \\
& L \ar[ur]_{} & \\
& C^*(\Gamma_{G,H})\rtimes C(G) \ar@{-->}[u]_\theta \ar@/_/[uur]_{\pi}& 
}$$

The first observation is that the injectivity assumption on $C(G)$ holds by construction, and that for $f \in C(G)$, the matrix $\pi(f)$ is ``block scalar'', the blocks corresponding to the indices $k$ in the basis $\varepsilon_{ke}$ in the basis from Proposition 16.13. Now for $r \in \Gamma_{G,H}$ with $\theta(r\otimes 1)=\theta(1 \otimes f)$ for some $f \in C(G)$, we see, using the commutative diagram, that we will have that $\pi(r \otimes 1)$ is block scalar. By Proposition 16.11, the family of representations $(\pi^k)$ of $\Gamma_{G,H}$, corresponding to the blocks $k$, is projectively faithful, so $r=1$. We can apply indeed Proposition 16.13, and we are done.
\end{proof}

Summarizing, we have computed the quantum permutation groups associated to the Di\c t\u a deformations of the tensor products of Fourier matrices, in the case where the deformation matrix $Q$ is generic. For some further computations, in the case where the deformation matrix $Q$ is no longer generic, we refer to \cite{bbi} and follow-up papers.

\section*{16c. Kesten measures}

Let us compute now the Kesten measure $\mu=law(\chi)$, in the case where the deformation matrix is generic, as before. Our results here will be a combinatorial moment formula, a geometric interpretation of it, and an asymptotic result. We first have:

\begin{theorem}
We have the moment formula
$$\int\chi^p
=\frac{1}{|G|\cdot|H|}\#\left\{\begin{matrix}i_1,\ldots,i_p\in G\\ d_1,\ldots,d_p\in H\end{matrix}\Big|\begin{matrix}[(i_1,d_1),(i_2,d_2),\ldots,(i_p,d_p)]\ \ \ \ \\=[(i_1,d_p),(i_2,d_1),\ldots,(i_p,d_{p-1})]\end{matrix}\right\}$$
where the sets between square brackets are by definition sets with repetition.
\end{theorem}

\begin{proof}
According to the various formulae above, the factorization found in Theorem 16.15 is, at the level of standard generators, as follows:
$$\begin{matrix}
C(S_{G\times H}^+)&\to&C^*(\Gamma_{G,H})\otimes C(G)&\to&M_{G\times H}(\mathbb C)\\
u_{ia,jb}&\to&\frac{1}{|H|}\sum_cF_{b-a,c}c^{(i)}\otimes v_{ij}&\to&W_{ia,jb}
\end{matrix}$$

Thus, the main character of the quantum permutation group that we found in Theorem 16.15 is given by the following formula:
\begin{eqnarray*}
\chi
&=&\frac{1}{|H|}\sum_{iac}c^{(i)}\otimes v_{ii}\\
&=&\sum_{ic}c^{(i)}\otimes v_{ii}\\
&=&\left(\sum_{ic}c^{(i)}\right)\otimes\delta_1
\end{eqnarray*}

Now since the Haar functional of $C^*(\Gamma)\rtimes C(H)$ is the tensor product of the Haar functionals of $C^*(\Gamma),C(H)$, this gives the following formula, valid for any $p\geq1$:
$$\int\chi^p=\frac{1}{|G|}\int_{\widehat{\Gamma}_{G,H}}\left(\sum_{ic}c^{(i)}\right)^p$$

Consider the elements $S_i=\sum_cc^{(i)}$. By using the embedding in Proposition 16.11 (2), with the notations there we have:
$$S_i=\sum_c(b_{i0}-b_{ic},c)$$

Now observe that these elements multiply as follows:
$$S_{i_1}\ldots S_{i_p}=\sum_{c_1\ldots c_p}
\begin{pmatrix}
b_{i_10}-b_{i_1c_1}+b_{i_2c_1}-b_{i_2,c_1+c_2}&&\\
+b_{i_3,c_1+c_2}-b_{i_3,c_1+c_2+c_3}+\ldots\ldots&,&c_1+\ldots+c_p&\\
\ldots\ldots+b_{i_p,c_1+\ldots+c_{p-1}}-b_{i_p,c_1+\ldots+c_p}&&
\end{pmatrix}$$

In terms of the new indices $d_r=c_1+\ldots+c_r$, this formula becomes:
$$S_{i_1}\ldots S_{i_p}=\sum_{d_1\ldots d_p}
\begin{pmatrix}
b_{i_10}-b_{i_1d_1}+b_{i_2d_1}-b_{i_2d_2}&&\\
+b_{i_3d_2}-b_{i_3d_3}+\ldots\ldots&,&d_p&\\
\ldots\ldots+b_{i_pd_{p-1}}-b_{i_pd_p}&&
\end{pmatrix}$$

Now by integrating, we must have $d_p=0$ on one hand, and on the other hand:
$$[(i_1,0),(i_2,d_1),\ldots,(i_p,d_{p-1})]=[(i_1,d_1),(i_2,d_2),\ldots,(i_p,d_p)]$$

Equivalently, we must have $d_p=0$ on one hand, and on the other hand:
$$[(i_1,d_p),(i_2,d_1),\ldots,(i_p,d_{p-1})]=[(i_1,d_1),(i_2,d_2),\ldots,(i_p,d_p)]$$

Thus, by translation invariance with respect to $d_p$, we obtain:
$$\int_{\widehat{\Gamma}_{G,H}}S_{i_1}\ldots S_{i_p}
=\frac{1}{|H|}\#\left\{d_1,\ldots,d_p\in H\Big|\begin{matrix}[(i_1,d_1),(i_2,d_2),\ldots,(i_p,d_p)]\ \ \ \ \\=[(i_1,d_p),(i_2,d_1),\ldots,(i_p,d_{p-1})]\end{matrix}\right\}$$

It follows that we have the following moment formula:
$$\int_{\widehat{\Gamma}_{G,H}}\left(\sum_iS_i\right)^p
=\frac{1}{|H|}\#\left\{\begin{matrix}i_1,\ldots,i_p\in G\\ d_1,\ldots,d_p\in H\end{matrix}\Big|\begin{matrix}[(i_1,d_1),(i_2,d_2),\ldots,(i_p,d_p)]\ \ \ \ \\=[(i_1,d_p),(i_2,d_1),\ldots,(i_p,d_{p-1})]\end{matrix}\right\}$$

Now by dividing by $|G|$, we obtain the formula in the statement.
\end{proof} 

The formula in Theorem 16.16 can be interpreted as follows:

\index{Gram matrix}

\begin{theorem}
With $M=|G|,N=|H|$ we have the formula
$$law(\chi)=\left(1-\frac{1}{N}\right)\delta_0+\frac{1}{N}law(A)$$
where the matrix on the right,
$$A\in C(\mathbb T^{MN},M_M(\mathbb C))$$
is given by $A(q)=$ Gram matrix of the rows of $q$.
\end{theorem}

\begin{proof}
According to Theorem 16.16, we have the following formula:
\begin{eqnarray*}
\int\chi^p
&=&\frac{1}{MN}\sum_{i_1\ldots i_p}\sum_{d_1\ldots d_p}\delta_{[i_1d_1,\ldots,i_pd_p],[i_1d_p,\ldots,i_pd_{p-1}]}\\
&=&\frac{1}{MN}\int_{\mathbb T^{MN}}\sum_{i_1\ldots i_p}\sum_{d_1\ldots d_p}\frac{q_{i_1d_1}\ldots q_{i_pd_p}}{q_{i_1d_p}\ldots q_{i_pd_{p-1}}}\,dq\\
&=&\frac{1}{MN}\int_{\mathbb T^{MN}}\sum_{i_1\ldots i_p}\left(\sum_{d_1}\frac{q_{i_1d_1}}{q_{i_2d_1}}\right)\left(\sum_{d_2}\frac{q_{i_2d_2}}{q_{i_3d_2}}\right)\ldots\left(\sum_{d_p}\frac{q_{i_pd_p}}{q_{i_1d_p}}\right)dq
\end{eqnarray*}

Consider now the Gram matrix in the statement, namely:
$$A(q)_{ij}=<R_i,R_j>$$

Here $R_1,\ldots,R_M$ are the rows of the following matrix:
$$q\in \mathbb T^{MN}\simeq M_{M\times N}(\mathbb T)$$

We have then the following computation:
\begin{eqnarray*}
\int\chi^p
&=&\frac{1}{MN}\int_{\mathbb T^{MN}}<R_{i_1},R_{i_2}><R_{i_2},R_{i_3}>\ldots<R_{i_p},R_{i_1}>\\
&=&\frac{1}{MN}\int_{\mathbb T^{MN}}A(q)_{i_1i_2}A(q)_{i_2i_3}\ldots A(q)_{i_pi_1}\\
&=&\frac{1}{MN}\int_{\mathbb T^{MN}}Tr(A(q)^p)dq\\
&=&\frac{1}{N}\int_{\mathbb T^{MN}}tr(A(q)^p)dq
\end{eqnarray*}

But this gives the formula in the statement, and we are done.
\end{proof}

In general, the moments of the Gram matrix $A$ are given by a quite complicated formula, and we cannot expect to have a refinement of Theorem 16.17, with $A$ replaced by a plain, non-matricial random variable, say over a compact abelian group. However, this kind of simplification appears at $M=2$, and since phenomenon this is quite interesting, we will explain this now. As a first remark, at $M=2$ we have:

\begin{proposition}
For $F_2\otimes_QF_H$, with $Q\in M_{2\times N}(\mathbb T)$ generic, we have
$$N\int\left(\frac{\chi}{N}\right)^p=\int_{\mathbb T^N}\sum_{k\geq0}\binom{p}{2k}\left|\frac{a_1+\ldots+a_N}{N}\right|^{2k}da$$
where the integral on the right is with respect to the uniform measure on $\mathbb T^N$.
\end{proposition}

\begin{proof}
In order to prove the result, consider the following quantity, which appeared in the proof of Theorem 16.17:
$$\Phi(q)=\sum_{i_1\ldots i_p}\sum_{d_1\ldots d_p}\frac{q_{i_1d_1}\ldots q_{i_pd_p}}{q_{i_1d_p}\ldots q_{i_pd_{p-1}}}$$

We can ``half-dephase'' the matrix $q\in M_{2\times N}(\mathbb T)$ if we want to, as follows:
$$q=\begin{pmatrix}1&\ldots&1\\ a_1&\ldots&a_N\end{pmatrix}$$

Let us compute now the above quantity $\Phi(q)$, in terms of the numbers $a_1,\ldots,a_N$. Our claim is that we have the following formula:
$$\Phi(q)=2\sum_{k\geq0}N^{p-2k}\binom{p}{2k}\left|\sum_ia_i\right|^{2k}$$

Indeed, the idea is that:

\medskip

-- The $2N^k$ contribution will come from $i=(1\ldots1)$ and $i=(2\ldots2)$.

\medskip

-- Then we will have a $p(p-1)N^{k-2}|\sum_ia_i|^2$ contribution coming from indices of type $i=(2\ldots 21\ldots1)$, up to cyclic permutations.

\medskip

-- Then we will have a $2\binom{p}{4}N^{p-4}|\sum_ia_i|^4$ contribution coming from indices of type $i=(2\ldots 21\ldots12\ldots21\ldots1)$.

\medskip

-- And so on. 

\medskip

In practice now, this gives the result. Indeed, in order to prove our claim, in order to find the $N^{p-2k}|\sum_ia_i|^{2k}$ contribution, we have to count the circular configurations consisting of $p$ numbers $1,2$, such that the $1$ values are arranged into $k$ non-empty intervals, and the $2$ values are arranged into $k$ non-empty intervals as well. Now by looking at the endpoints of these $2k$ intervals, we have $2\binom{p}{2k}$ choices, and this gives the above formula. Now by integrating, this gives the formula in the statement.
\end{proof}

Observe now that the integrals in Proposition 16.18 can be computed as follows:
\begin{eqnarray*}
\int_{\mathbb T^N}|a_1+\ldots+a_N|^{2k}da
&=&\int_{\mathbb T^N}\sum_{i_1\ldots i_k}\sum_{j_1\ldots j_k}\frac{a_{i_1}\ldots a_{i_k}}{a_{j_1}\ldots a_{j_k}}da\\
&=&\#\left\{i_1\ldots i_k,j_1\ldots j_k\Big|[i_1,\ldots,i_k]=[j_1,\ldots,j_k]\right\}\\
&=&\sum_{k=\sum r_i}\binom{k}{r_1,\ldots,r_N}^2
\end{eqnarray*}

We obtain in this way the following ``blowup'' result, for our measure:

\index{blowup}

\begin{proposition}
For $F_2\otimes_QF_H$, with $Q\in M_{2\times N}(\mathbb T)$ generic, we have
$$\mu=\left(1-\frac{1}{N}\right)\delta_0+\frac{1}{2N}\left(\Psi^+_*\varepsilon+\Psi^-_*\varepsilon\right)$$
where $\varepsilon$ is the uniform measure on $\mathbb T^N$, and where the blowup function is:
$$\Psi^\pm(a)=N\pm\left|\sum_ia_i\right|$$
\end{proposition}

\begin{proof}
We use the formula found in Proposition 16.18, along with the following standard identity, coming from the Taylor formula:
$$\sum_{k\geq0}\binom{p}{2k}x^{2k}=\frac{(1+x)^p+(1-x)^p}{2}$$

By using this identity, Proposition 16.18 reformulates as follows:
$$N\int\left(\frac{\chi}{N}\right)^p=\frac{1}{2}\int_{\mathbb T^N}\left(1+\left|\frac{\sum_ia_i}{N}\right|\right)^p+\left(1-\left|\frac{\sum_ia_i}{N}\right|\right)^p\,da$$

Now by multiplying by $N^{p-1}$, we obtain the following formula:
$$\int\chi^k=\frac{1}{2N}\int_{\mathbb T^N}\left(N+\left|\sum_ia_i\right|\right)^p+\left(N-\left|\sum_ia_i\right|\right)^p\,da$$

But this gives the formula in the statement, and we are done.
\end{proof}

We can further improve the above result, by reducing the maps $\Psi^\pm$ appearing there to a single one, and we are led to the following statement: 

\begin{theorem}
For $F_2\otimes_QF_H$, with $Q\in M_{2\times N}(\mathbb T)$ generic, we have
$$\mu=\left(1-\frac{1}{N}\right)\delta_0+\frac{1}{N}\Phi_*\varepsilon$$
where $\varepsilon$ is the uniform measure on $\mathbb Z_2\times\mathbb T^N$, and where the blowup map is:
$$\Phi(e,a)=N+e\left|\sum_ia_i\right|$$
\end{theorem}

\begin{proof}
This is clear indeed from Proposition 16.19.
\end{proof}

As already mentioned, the above results at $M=2$ are something quite special. In the general case, $M\in\mathbb N$, it is not clear how to construct a nice blowup of the measure. All the above results are quite interesting in the general context of subfactor theory, where the blowup question is one of the main open questions, related to the continuations of Jones' planar algebra work in \cite{jo3}, and to many other things, mainly coming from advanced quantum physics. For more on all this, we refer to \cite{bbi} and its previous versions, which were more subfactor-centered, and which can be found on the internet.

\section*{16d. Poisson laws}

Let us go back now to the general case, where $M,N\in\mathbb N$ are arbitrary. The problem that we would like to solve is that of finding the good regime, of the following type, where the measure in Theorem 16.16 converges, after some suitable manipulations:
$$M=f(K)\quad,\quad 
N=g(K)\quad,\quad 
K\to\infty$$

As before by following \cite{bbi}, we will see that this is indeed possible, and that as limiting laws we have some very interesting objects, namely some versions of the Marchenko-Pastur laws, or free Poisson laws, that we met at the end of chapter 13. Let us first recall from there the definition and main properties of these laws, in the general context:

\begin{theorem}
The following Poisson limits converge, for any $t>0$,
$$p_t=\lim_{n\to\infty}\left(\left(1-\frac{t}{n}\right)\delta_0+\frac{t}{n}\delta_1\right)^{*n}\quad,\quad 
\pi_t=\lim_{n\to\infty}\left(\left(1-\frac{t}{n}\right)\delta_0+\frac{t}{n}\delta_1\right)^{\boxplus n}$$
the limiting measures being the Poisson law $p_t$, and the Marchenko-Pastur law $\pi_t$,
$$p_t=\frac{1}{e^t}\sum_{k=0}^\infty\frac{t^k\delta_k}{k!}\quad,\quad 
\pi_t=\max(1-t,0)\delta_0+\frac{\sqrt{4t-(x-1-t)^2}}{2\pi x}\,dx$$
with at $t=1$, the Marchenko-Pastur law being given by the following formula:
$$\pi_1=\frac{1}{2\pi}\sqrt{4x^{-1}-1}\,dx$$
Moreover, the moments of these laws are given by the formulae
$$M_k(p_t)=\sum_{\pi\in P(k)}t^{|\pi|}\quad,\quad 
M_k(\pi_t)=\sum_{\pi\in NC(k)}t^{|\pi|}$$
where $|.|$ is the number of blocks. 
\end{theorem}

\begin{proof}
All this is standard probability and free probability theory:

\medskip

(1) In what regards the classical results, concerning $p_t$, the standard way of viewing them is by defining the Poisson law $p_t$ by the formula in the statement, then by establishing the Poisson Limiting Theorem (PLT) via Fourier transform, and finally by working out the moment formula either by recurrence, or from Fourier via cumulants.

\medskip

(2) In the free case now, in relation with $\pi_t$, pretty much the same procedure can be used, with however the change that the study of free PLT comes first, by using Voiculescu's $R$-transform, which produces then via Stieltjes inversion the formula of $\pi_t$ in the statement. We refer here to \cite{vdn}, or to any other free probability book. 
\end{proof}

In order to establish our results, we have to do some combinatorics. We denote by $NC(p)$ the set of noncrossing partitions of $\{1,\ldots,p\}$, and for $\pi\in P(p)$ we denote by $|\pi|\in\{1,\ldots,p\}$ the number of blocks. We will also use some standard tools from combinatorics, such as the Kreweras complementation, which are well-known in free probability \cite{vdn}. With these conventions, we have the following result from \cite{bbi}, regarding the moments $c_p$ of the measure that we are interested in, computed in Theorem 16.16:

\begin{proposition}
With $M=\alpha K,N=\beta K$, $K\to\infty$ we have:
$$\frac{c_p}{K^{p-1}}\simeq\sum_{r=1}^p\#\left\{\pi\in NC(p)\Big||\pi|=r\right\}\alpha^{r-1}\beta^{p-r}$$
In particular, with $\alpha=\beta$ we have:
$$c_p\simeq\frac{1}{p+1}\binom{2p}{p}(\alpha K)^{p-1}$$
\end{proposition}

\begin{proof}
We use the combinatorial formula in Theorem 16.16. Our claim is that, with $\pi=\ker(i_1,\ldots,i_p)$, the corresponding contribution to $c_p$ is:
$$C_\pi\simeq
\begin{cases}
\alpha^{|\pi|-1}\beta^{p-|\pi|}K^{p-1}&{\rm if}\ \pi\in NC(p)\\
O(K^{p-2})&{\rm if}\ \pi\notin NC(p)
\end{cases}$$

As a first observation, the number of choices for a multi-index $(i_1,\ldots,i_p)\in X^p$ satisfying the condition $\ker i=\pi$ is:
$$M(M-1)\ldots (M-|\pi|+1)\simeq M^{|\pi|}$$

Thus, we have the following estimate:
$$C_\pi\simeq M^{|\pi|-1}N^{-1}\#\left\{d_1,\ldots,d_p\in Y\Big|[d_\alpha|\alpha\in b]=[d_{\alpha-1}|\alpha\in b],\forall b\in\pi\right\}$$

Consider now the following partition:
$$\sigma=\ker d$$

The contribution of $\sigma$ to the above quantity $C_\pi$ is then given by:
$$\Delta(\pi,\sigma)N(N-1)\ldots(N-|\sigma|+1)\simeq\Delta(\pi,\sigma)N^{|\sigma|}$$

Here the quantities on the right are as follows:
$$\Delta(\pi,\sigma)=\begin{cases}
1&{\rm if}\ |b\cap c|=|(b-1)\cap c|,\forall b\in\pi,\forall c\in\sigma\\
0&{\rm otherwise}
\end{cases}$$

We use now the standard fact that for $\pi,\sigma\in P(p)$ satisfying $\Delta(\pi,\sigma)=1$ we have:
$$|\pi|+|\sigma|\leq p+1$$

\index{Kreweras complementation}

In addition, the equality case is known to happen when $\pi,\sigma\in NC(p)$ are inverse to each other, via Kreweras complementation. This shows that for $\pi\notin NC(p)$ we have:
$$C_\pi=O(K^{p-2})$$

Also, this shows that for $\pi\in NC(p)$ we have:
\begin{eqnarray*}
C_\pi
&\simeq&M^{|\pi|-1}N^{-1}N^{p-|\pi|-1}\\
&=&\alpha^{|\pi|-1}\beta^{p-|\pi|}K^{p-1}
\end{eqnarray*}

Thus, we have obtained the result.
\end{proof}

We denote by $D$ the dilation operation for probability measures, given by:
$$D_r(law(X))=law(rX)$$

With this convention, we have the following result, based on Proposition 16.22:

\begin{theorem}
With $M=\alpha K,N=\beta K$, $K\to\infty$ we have:
$$\mu=\left(1-\frac{1}{\alpha\beta K^2}\right)\delta_0+\frac{1}{\alpha\beta K^2}D_{\frac{1}{\beta K}}(\pi_{\alpha/\beta})$$
In particular with $\alpha=\beta$ we have:
$$\mu=\left(1-\frac{1}{\alpha^2K^2}\right)\delta_0+\frac{1}{\alpha^2K^2}D_{\frac{1}{\alpha K}}(\pi_1)$$
\end{theorem}

\begin{proof}
At $\alpha=\beta$, this follows from Proposition 16.22. In general now, we have:
\begin{eqnarray*}
\frac{c_p}{K^{p-1}}
&\simeq&\sum_{\pi\in NC(p)}\alpha^{|\pi|-1}\beta^{p-|\pi|}\\
&=&\frac{\beta^p}{\alpha}\sum_{\pi\in NC(p)}\left(\frac{\alpha}{\beta}\right)^{|\pi|}\\
&=&\frac{\beta^p}{\alpha}\int x^pd\pi_{\alpha/\beta}(x)
\end{eqnarray*}

When $\alpha\geq\beta$, where $d\pi_{\alpha/\beta}(x)=\varphi_{\alpha/\beta}(x)dx$ is continuous, we obtain:
\begin{eqnarray*}
c_p
&=&\frac{1}{\alpha K}\int(\beta Kx)^p\varphi_{\alpha/\beta}(x)dx\\
&=&\frac{1}{\alpha\beta K^2}\int x^p\varphi_{\alpha/\beta}\left(\frac{x}{\beta K}\right)dx
\end{eqnarray*}

But this gives the formula in the statement. When $\alpha\leq\beta$ the computation is similar, with a Dirac mass as 0 dissapearing and reappearing, and gives the same result.
\end{proof}

Let us state as well an explicit result, regarding densities:

\begin{theorem}
With $M=\alpha K,N=\beta K$, $K\to\infty$ we have:
$$\mu=\left(1-\frac{1}{\alpha\beta K^2}\right)\delta_0+\frac{1}{\alpha\beta K^2}\cdot\frac{\sqrt{4\alpha\beta K^2-(x-\alpha K-\beta K)^2}}{2\pi x}\,dx$$
In particular with $\alpha=\beta$ we have:
$$\mu=\left(1-\frac{1}{\alpha^2K^2}\right)\delta_0+\frac{1}{\alpha^2K^2}\cdot\frac{\sqrt{\frac{4\alpha K}{x}-1}}{2\pi}$$
\end{theorem}

\begin{proof}
According to the formula for the density of the free Poisson law, the density of the continuous part $D_{\frac{1}{\beta K}}(\pi_{\alpha/\beta})$ is indeed given by:
$$\frac{\sqrt{4\frac{\alpha}{\beta}-(\frac{x}{\beta K}-1-\frac{\alpha}{\beta})^2}}
{2\pi\cdot\frac{x}{\beta K}}=\frac{\sqrt{4\alpha\beta K^2-(x-\alpha K-\beta K)^2}}{2\pi x}$$

With $\alpha=\beta$ now, we obtain the second formula in the statement, and we are done.
\end{proof}

Observe that at $\alpha=\beta=1$, where $M=N=K\to\infty$, the above measure is:
$$\mu=\left(1-\frac{1}{K^2}\right)\delta_0+\frac{1}{K^2}D_{\frac{1}{K}}(\pi_1)$$

This measure is supported by $[0,4K]$. On the other hand, since the groups $\Gamma_{M,N}$ are all amenable, the corresponding measures are supported on $[0,MN]$, and so on $[0,K^2]$ in the $M=N=K$ situation. The fact that we do not have a convergence of supports is not surprising, because our convergence is in moments.

\bigskip

The above results are of course not the end of the story, because we have now to understand what happens in the case of non-generic parameters. There has been some technical work here, by Bichon and by myself, and as a sample result here, we have:

\index{free Poisson law}
\index{Marchenko-Pastur law}

\begin{theorem}
Given two finite abelian groups $G,H$, having cardinalities 
$$|G|=M\quad,\quad 
|H|=N$$ 
consider the main character $\chi$ of the quantum group associated to $\mathcal F_{G\times H}$. We have then
$$law\left(\frac{\chi}{N}\right)=\left(1-\frac{1}{M}\right)\delta_0+\frac{1}{M}\,\pi_t$$
in moments, with $M=tN\to\infty$, where $\pi_t$ is the free Poisson law of parameter $t>0$. In addition, this formula holds for any generic fiber of $\mathcal F_{G\times H}$.
\end{theorem}

\begin{proof}
We already know that the second assertion holds, as explained above. Regarding now the first assertion, our first claim is that for the representation coming from the parametric matrix $\mathcal F_{G\times H}$ we have the following formula, where $M=|G|,N=|H|$, and the sets between brackets are sets with repetitions:
$$c_p^r=\frac{1}{M^{r+1}N}\#\left\{\begin{matrix}i_1,\ldots,i_r,a_1,\ldots,a_p\in\{0,\ldots,M-1\},\\
b_1,\ldots,b_p\in\{0,\ldots,N-1\},\\
[(i_x+a_y,b_y),(i_{x+1}+a_y,b_{y+1})|y=1,\ldots,p]\\
=[(i_x+a_y,b_{y+1}),(i_{x+1}+a_y,b_y)|y=1,\ldots,p], \forall x
\end{matrix}\right\}$$

Indeed, by using the general moment formula with $K=F_G$, $L=F_H$, we have the following formula for the above numbers:
\begin{eqnarray*}
&&c_p^r\\
&=&\frac{1}{(MN)^r}\int_{T^r}\sum_{i_1^1\ldots i_p^r}\sum_{b_1^1\ldots b_p^r}\frac{Q^1_{i_1^1b_1^1}Q^1_{i_1^2b_2^1}}{Q^1_{i_1^1b_2^1}Q^1_{i_1^2b_1^1}}\ldots\frac{Q^1_{i_p^1b_p^1}Q^1_{i_p^2b_1^1}}{Q^1_{i_p^1b_1^1}Q^1_{i_p^2b_p^1}}\ldots\ldots\frac{Q^r_{i_1^rb_1^r}Q^r_{i_1^1b_2^r}}{Q^r_{i_1^rb_2^r}Q^r_{i_1^1b_1^r}}\ldots\frac{Q^r_{i_p^rb_p^r}Q^r_{i_p^1b_1^r}}{Q^r_{i_p^rb_1^r}Q^r_{i_p^1b_p^r}}\\
&&\hskip15mm\frac{1}{M^{pr}}\sum_{j_1^1\ldots j_p^r}\frac{K_{i_1^1j_1^1}K_{i_1^2j_2^1}}{K_{i_1^1j_2^1}K_{i_1^2j_1^1}}\ldots\frac{K_{i_p^1j_p^1}K_{i_p^2j_1^1}}{K_{i_p^1j_1^1}K_{i_p^2j_p^1}}\ldots\ldots\frac{K_{i_1^rj_1^r}K_{i_1^1j_2^r}}{K_{i_1^rj_2^r}K_{i_1^1j_1^r}}\ldots\frac{K_{i_p^rj_p^r}K_{i_p^1j_1^r}}{K_{i_p^rj_1^r}K_{i_p^1j_p^r}}\\
&&\hskip15mm\frac{1}{N^{pr}}\sum_{a_1^1\ldots a_p^r}\frac{L_{a_1^1b_1^1}L_{a_1^2b_2^1}}{L_{a_1^1b_2^1}L_{a_1^2b_1^1}}\ldots\frac{L_{a_p^1b_p^1}L_{a_p^2b_1^1}}{L_{a_p^1b_1^1}L_{a_p^2b_p^1}}\ldots\ldots\frac{L_{a_1^rb_1^r}L_{a_1^1b_2^r}}{L_{a_1^rb_2^r}L_{a_1^1b_1^r}}\ldots\frac{L_{a_p^rb_p^r}L_{a_p^1b_1^r}}{L_{a_p^rb_1^r}L_{a_p^1b_p^r}}\,dQ\\
\end{eqnarray*}

Since we are in the Fourier matrix case, $K=F_G,L=F_H$, we can perform the sums over $j,a$. To be more precise, the last two averages appearing above are respectively:
\begin{eqnarray*}
\Delta(i)&=&\prod_x\prod_y\delta(i^x_y+i^{x+1}_{y-1},i^{x+1}_y+i^x_{y-1})\\
\Delta(b)&=&\prod_x\prod_y\delta(b^x_y+b^{x+1}_{y-1},b^{x+1}_y+b^x_{y-1})
\end{eqnarray*}

We therefore obtain the following formula for the truncated moments of the main character, where $\Delta$ is the product of Kronecker symbols constructed above:
\begin{eqnarray*}
&&c_p^r\\
&=&\frac{1}{(MN)^r}\int_{T^r}\sum_{\Delta(i)=\Delta(b)=1}\frac{Q^1_{i_1^1b_1^1}Q^1_{i_1^2b_2^1}}{Q^1_{i_1^1b_2^1}Q^1_{i_1^2b_1^1}}\ldots\frac{Q^1_{i_p^1b_p^1}Q^1_{i_p^2b_1^1}}{Q^1_{i_p^1b_1^1}Q^1_{i_p^2b_p^1}}\ldots\ldots\frac{Q^r_{i_1^rb_1^r}Q^r_{i_1^1b_2^r}}{Q^r_{i_1^rb_2^r}Q^r_{i_1^1b_1^r}}\ldots\frac{Q^r_{i_p^rb_p^r}Q^r_{i_p^1b_1^r}}{Q^r_{i_p^rb_1^r}Q^r_{i_p^1b_p^r}}\,dQ
\end{eqnarray*}

Now by integrating with respect to $Q\in(\mathbb T^{G\times H})^r$, we are led to counting the multi-indices $i,b$ satisfying several conditions. First, we have the following condition:
$$\Delta(i)=\Delta(b)=1$$

We have as well the following conditions, where the sets between brackets are by definition sets with repetitions:
$$\begin{bmatrix}
i_1^1b_1^1&\ldots&i_p^1b_p^1&i_1^2b_2^1&\ldots&i_p^2b_1^1
\end{bmatrix}=
\begin{bmatrix}
i_1^1b_2^1&\ldots&i_p^1b_1^1&i_1^2b_1^1&\ldots&i_p^2b_p^1
\end{bmatrix}$$
$$\vdots$$
$$\begin{bmatrix}
i_1^rb_1^r&\ldots&i_p^rb_p^r&i_1^1b_2^r&\ldots&i_p^1b_1^r
\end{bmatrix}
=\begin{bmatrix}
i_1^rb_2^r&\ldots&i_p^rb_1^r&i_1^1b_1^r&\ldots&i_p^1b_p^r
\end{bmatrix}$$

In a more compact notation, the moment formula that we obtain in this way is therefore as follows:
$$c_p^r=\frac{1}{(MN)^r}\#\left\{i,b\Big|\Delta(i)=\Delta(b)=1,\ [i^x_yb^x_y,i^{x+1}_yb^x_{y+1}]=[i^x_yb^x_{y+1},i^{x+1}_yb^x_y],\forall x\right\}$$

Now observe that the above Kronecker type conditions $\Delta(i)=\Delta(b)=1$ tell us that the arrays of indices $i=(i^x_y),b=(b^x_y)$ must be of the following special form:
$$\begin{pmatrix}i^1_1&\ldots&i^1_p\\&\ldots\\ i^1_r&\ldots&i^r_p\end{pmatrix}=\begin{pmatrix}i_1+a_1&\ldots&i_1+a_p\\&\ldots\\ i_r+a_1&\ldots&i_r+a_p\end{pmatrix}$$

$$\begin{pmatrix}b^1_1&\ldots&b^1_p\\&\ldots\\ b^1_r&\ldots&b^r_p\end{pmatrix}=\begin{pmatrix}j_1+b_1&\ldots&j_1+b_p\\&\ldots\\ j_r+b_1&\ldots&j_r+b_p\end{pmatrix}$$

Here all the new indices $i_x,j_x,a_y,b_y$ are uniquely determined, up to a choice of $i_1,j_1$. Now by replacing $i^x_y,b^x_y$ with these new indices $i_x,j_x,a_y,b_y$, with a $MN$ factor added, which accounts for the choice of $i_1,j_1$, we obtain the following formula:
$$c_p^r=\frac{1}{(MN)^{r+1}}\#\left\{i,j,a,b\Big|\begin{matrix}[(i_x+a_y,j_x+b_y),(i_{x+1}+a_y,j_x+b_{y+1})]\\
=[(i_x+a_y,j_x+b_{y+1}),(i_{x+1}+a_y,j_x+b_y)],\forall x\end{matrix}\right\}$$

Now observe that we can delete if we want the $j_x$ indices, which are irrelevant. Thus, we obtain the announced formula. The continuation is via combinatorics.
\end{proof}

There are many interesting questions that are still open, regarding the computation of the spectral measure in the case where the parameter matrix $Q$ is not generic, and also regarding the computation for the deformations of the generalized Fourier matrices, which are not necessarily of Di\c t\u a type. We refer here to \cite{bbi} and related papers.

\section*{16e. Exercises} 

To start with, we have the following exercise from the previous chapter, which is related to the above, and that we reproduce here, in case you have not solved it yet:

\begin{exercise}
Write down a complete, simplified proof for the factorization
$$\xymatrix{C(S_{NM}^+)\ar[rr]^{\pi_L}\ar[rd]&&M_{NM}(\mathbb C)\\&C(S_M^+\wr_*G_H)\ar[ur]&}$$
found in the previous chapter, for $L=H\otimes_QK$, in the scalar matrix case.
\end{exercise}

This exercise is important, because it is related to the first factorization performed in this chapter, in the context of the Fourier models.

\begin{exercise}
Work out the combinatorial details of the computation for deformed Fourier models with formal parameters, outlined in the proof of Theorem 16.25.
\end{exercise}

This is actually quite unobvious, but finding the relevant literature and writing up a concise account of what is done there would do.

\begin{exercise}
Do some computations for the deformations of $F_4$, at non generic values of the parameter, and write down what you found.
\end{exercise}

And that is all. In the hope that you liked the present book, and that we will hear from you soon, with interesting results about the Hadamard matrices. There are just so many things to be done, all interesting. You can't go wrong with these matrices.

\baselineskip=14pt

\printindex

\end{document}